\documentclass[a4paper]{article}

\usepackage{amsmath,amssymb,amsthm,epsfig}
\usepackage[utf8]{inputenc}
\usepackage{color}
\usepackage[pdftex,colorlinks]{hyperref}

\usepackage[all]{xy}

\makeatother

\newtheorem{theor}{theorem}[section]
\newtheorem{theorem}[theor]{Theorem}
\newtheorem{defi}[theor]{Definition}
\newtheorem{defis}[theor]{Definitions}
\newtheorem{prop}[theor]{Proposition}
\newtheorem{corollary}[theor]{Corollary}
\newtheorem{remarks}[theor]{Remarks}
\newtheorem{remark}[theor]{Remark}
\newtheorem{lemma}[theor]{Lemma}
\newtheorem{compari}[theor]{Comparison}
\newtheorem{exas}[theor]{Examples}
\newtheorem{exa}[theor]{Example}

\newcommand{\Sys}{\operatorname{Sys}}
\newcommand{\Sec}{\operatorname{Sect}}

\newcommand{\N}{\mathbb{N}}
\newcommand{\Z}{\mathbb{Z}}
\newcommand{\R}{\mathbb{R}}
\newcommand{\C}{\mathbb{C}}
\newcommand{\f}{\rightarrow}

\newcommand{\scal}{\operatorname{scal}}
\newcommand{\CAT}{\operatorname{CAT}}

\newcommand{\id}{\operatorname{id}}

\newcommand{\Ker}{\operatorname{Ker}}

\newcommand{\Vol}{\operatorname{Vol}}

\newcommand{\Ric}{\operatorname{Ric}}

\newcommand{\Ent}{\operatorname{Ent}}
\newcommand{\Min}{\operatorname{Min}}
\newcommand{\Max}{\operatorname{Max}}
\newcommand{\sys}{\operatorname{sys}}
\newcommand{\dias}{\operatorname{Dias}}
\newcommand{\diam}{\operatorname{diam}}
\newcommand{\inj}{\operatorname{inj}}

\newcommand{\g}{\operatorname{\gamma}}
\newcommand{\e}{\operatorname{\varepsilon}}

\title{Curvature-Free Margulis Lemma\\
for Gromov-Hyperbolic Spaces}
\author{G.~Besson, G.~Courtois, S.~Gallot, A.~Sambusetti}

\begin{document}

\maketitle

\begin{abstract}
We prove curvature-free versions of the celebrated Margulis Lemma. We are interested by both the algebraic aspects and the geometric ones, with however an emphasis on the second and we aim at giving quantitative (computable) estimates of some important invariants. Our goal is to get rid of the pointwise curvature assumptions in order to extend the results to more general spaces such as certain metric spaces. Essentially the upper bound on the curvature is replaced by the assumption that the space is hyperbolic in the sense of Gromov and the lower bound of the curvature by an upper bound on the entropy of which we recall the definition.
\end{abstract}

\tableofcontents

\section{Introduction}

The celebrated Margulis Lemma is the keystone of a beautiful theory of the structure of complete Riemannian manifolds with bounded sectional curvature. It has two main aspects: the first one is algebraic and concerns the fundamental group of the manifold, the second one is more geometric and yields a thin-thick decomposition of the manifold. To be more precise let us state a weak version of this lemma pertaining to the first aspect (see \cite{Mar}, \cite{Ba-Gr-Sc}, \cite{Bu-Za} Section 37.3).

\begin{theorem}\label{Margulisclassic}
There exist constants $ \varepsilon (n)>0 $ and $ C(n)>0 $ such that, for every complete Riemannian manifold $ M $ whose sectional curvature satisfies $-1\leq \Sec (M)\leq 0 $, every point $ p\in M $ and every $ \varepsilon \leq \varepsilon (n) $, the subgroup $ \Gamma_{\varepsilon}(p) $ of $ \pi_1 (M) $ generated by the loops at $ p $ of length less than $ \varepsilon $ is virtually nilpotent. Furthermore, the index of the nilpotent subgroup is bounded above by $ C(n) $.
\end{theorem}

This statement is a weak version, indeed in the strong one the upper bound on the sectional curvature could be positive, with extra assumption though.
A version of this theorem for manifolds of strictly negative curvature was simultaneously proved by E. Heintze in his habilitationsschrift of 1976 (see \cite{zbMATH01825238}).

\smallskip
The history of this result goes back to Bieberbach Theorem  (\cite{Bie})  which describes the discrete subgroups of the isometry group of the Euclidean spaces and consequently gives a structure theorem for the flat manifolds and orbifolds. Later, this result was extended to the study of discrete subgroups of Lie groups
by H. Zassenhaus (\cite{Zass}) and to locally symmetric spaces by D. Kazhdan and G. Margulis (\cite{KM}, using Zassenhaus' Lemma). Recently there has been progresses on the question of extending this result to different spaces or curvature conditions: for example, after a short sketch of proof by J. Cheeger and T. Colding (see \cite{CC96}) 
under the hypothesis \lq \lq Ricci curvature bounded from below" and a first complete proof of V. Kapovitch, A. Petrunin and W. 
Tuschmann (\cite{KPT}) under the hypothesis \lq \lq sectional curvature bounded from below", V. Kapovitch and B. Wilking (\cite{KW}) recently established a Margulis-like Lemma under the hypothesis \lq \lq Ricci curvature bounded from below by $- (n-1)$" 
instead of \lq \lq $-1\leq \Sec (M)\leq 0 $" (see also \cite{Cou} for references and a detailed exposition).

\smallskip
This paper is the first of a series of articles devoted to this theme. Here we are interested by both aspects, algebraic and geometric, with however an emphasis on the second  and we aim at giving quantitative (computable) estimates of some important invariants. Our goal is to get rid of the pointwise curvature assumptions, as mentioned in the title, in order to extend the results to more general spaces such as certain metric spaces.  Essentially the upper bound on the curvature is replaced by the assumption that the space is $ \delta $-hyperbolic (in the sense of Gromov, see \cite{Gr4} and Section \ref{coherence} of Appendix \ref{outils} for precise definitions) and the lower bound by an upper bound on the entropy which we define below. Notice that $ \delta $ behaves like a distance and it is rather $ \delta^{-2} $ which is curvature-like.\\
The starting point of the ideas developed in this paper is the prepublication \cite{BCG}, which initially concerned the isometric actions of a more limited class of 
groups on less general types of spaces, 
namely fundamental groups of manifolds with sectional curvature $\sigma \le -1$ and with injectivity radius $\ge i_0 > 0$ (and groups such that any non 
abelian subgroup with two generators admits an injective homomorphism into such a fundamental group). Several developments of the ideas contained in 
\cite{BCG} were established by G. Reviron (\cite{Re}), F. Zuddas (\cite{Zu}, \cite{Zu2}), F. Cerocchi (\cite{Ce}), F. Cerocchi and A. Sambusetti 
(\cite{CS2}, \cite{CS3} and \cite{CS}, this last paper being devoted to prove Margulis' properties in the abelian setting)\ .

\medskip
Let $ (X,d) $ be a (non-elementary) metric space  which we assume to be \textsl{proper}, \textsl{i.e} the closed metric balls are compact. We only consider metric spaces which are geodesic. More precisely, a geodesic segment is the image of an interval of $ I\subset\textbf{R} $ by an isometric map from $ I $ into $ X $. The space $ (X,d) $ is said to be \textsl{geodesic} if any two points of $ X $ are joined by at least one geodesic segment. Let $ \mu $ be a positive (non identically zero) Borel measure. We call $ (X,d, \mu) $ a \textsl{measured metric space}.

\begin{defi}
Let $ (X,d,\mu ) $ be a metric measured space we define its entropy by
 $$ \Ent (X,d, \mu)=\liminf_{R\to+\infty}\frac{1}{R}\ln\big (\mu (B_X(x,R))\big ) $$ 
where $ B_X(x,R) $ is the open ball of radius $ R $ around $ x\in X $. Furthermore, the entropy is independent of $ x $.
\end{defi}
In the sequel we will consider a group $ \Gamma $ acting by isometries on $ (X,d) $ \textsl{properly} and, often, co-compactly. We recall that the action is said to be \textsl{proper} if for $ x\in X $ and for all $ R>0 $, the number of $ \gamma\in\Gamma $ satisfying $ d(x,\gamma x)\leq R $ is finite; this does not depend on $ x $. In that case, for any measure $ \mu $ \textsl{invariant} by $ \Gamma $ the above definition yields the same number which we call the entropy of $ (X,d) $ and denote by $ \Ent (X,d) $. If $ (X,d) $ is a $ \delta $-hyperbolic geodesic metric space and $ \gamma $ a torsion-free isometry, we define the {\textsl{asymptotic displacement} of $ \gamma $ (sometimes also called \lq \lq stable displacement") 
by
$$ \ell (\gamma )=\lim_{k\to +\infty} \frac{1}{k}d(x, \gamma^kx)  ,$$ 
this definition does not depend on the choice of $ x\in X $.
One of our results is the following (see Proposition \ref{minorentropie}).

\begin{prop}
For every non-elementary $ \delta $-hyperbolic metric space $ (X,d) $ and every group $ \Gamma $ acting properly by isometries on $ (X,d) $, if $ \diam (\Gamma \backslash X) \leq D<+\infty $, then 
 $$ \Ent (X,d) \geq \frac{\ln 2}{L+17\delta +2 D}\geq \frac{\ln 2}{27\delta +10 D}, $$ 
where $ L=\inf\big\{\ell (\gamma): \gamma \textrm{ hyperbolic element of } \Gamma\setminus\{e\}\big\} $.
\end{prop}

Note that in the proof we show that, in the above situation, there always exists an hyperbolic element $ \gamma\in\Gamma\setminus\{e\}$
 which satisfies $ \ell (\gamma )\leq 8 D + 10 \delta$.

Now, let $ \Gamma $ be a group which is non-elementary (\textsl{i.e.} whose boundary has at least three points) and $ \Sigma$ be a finite generating set for $\Gamma$. The Cayley graph of $ \Gamma $ defined by $ \Sigma $ is a metric space when endowed with the distance such that the edges have length $ 1 $. We say that $ (\Gamma , \Sigma ) $ is a $ \delta $-hyperbolic group if this metric space is $ \delta $-hyperbolic. The group $ \Gamma $ acts by isometries on this metric space and the quotient is compact and has diameter $ 1 $. The entropy of this metric space is denoted by $ \Ent(\Gamma, \Sigma ) $ and is called the algebraic entropy of $ \Gamma $ with respect to $ \Sigma $. We also define the algebraic entropy of $ \Gamma $ by
 $$ \Ent (\Gamma )=\inf_\Sigma \{\Ent(\Gamma, \Sigma )\} , $$ 
the infimum being taken among all finite generating sets $ \Sigma $.

The study of the algebraic entropy of groups with exponential growth has recently made progresses. When the group acts on a Hadamard manifold the three first authors have proved a quantitative version of the Tits alternative, see \cite{BCG2} and the references herein. 

A corollary of the above Proposition is the following statement (see \ref{minorentalg}),

\begin{corollary}\label{cor:intro-entalg} 
Let $ \Gamma $ be a non-elementary group and $ \Sigma $ a finite generating set such that $ (\Gamma , \Sigma ) $ is $ \delta $-hyperbolic, then
 $$ \Ent (\Gamma , \Sigma )\geq \frac{\ln 2}{27\delta +10}\,..$$ 
\end{corollary}

\begin{remarks}
Once the first version of this article was completed E.~Breuillard mentioned to us his joint work with K.~Fujiwara (see, \cite{BF}) which contains an improvement of \cite{BCG} and \cite{BCG2}; Corollary \ref{cor:intro-entalg} is then similar to their result.
\end{remarks}

These two estimates do not provide any lower bound of the algebraic entropy of the group $\Gamma$. Indeed, when the generating system varies 
in order to minimize the Entropy, its hyperbolicity constant may go to infinity. We obtain such a lower bound of the algebraic entropy as a consequence of 
our Bishop-Gromov inequality (see Theorems \ref{intro:cocompact2} and \ref{cocompact2}), namely:

\begin{corollary}\label{intro-entalg3}
Let $\Gamma$ be a non virtually cyclic Gromov-hyperbolic group then, for every positive constant $M$, if there exists a finite system $S_0$ of generators of 
$\Gamma$ such that $(\Gamma, S_0)$ is $\delta_0$-hyperbolic and $ \Ent (\Gamma, S_0) \, (\delta_0+ 1) \le M$, then the algebraic entropy of 
$\Gamma$ and of any finitely generated and non virtually cyclic subgroup $\Gamma'$ of $\Gamma$  is bounded from below by 
$\dfrac{\ln 2}{42\, N \left(\left[3^{12} \, e^{ 490 \,M }\right] + 1 \right) + 2}$, where $N(\cdot )$ is the function which appears in Theorem \ref{BGT}.
\end{corollary}

This study starts with the simple remark that if two elements, $a$ and $b$, of a discrete subgroup $ \Gamma $ of the isometry group of a Hadamard manifold $ X $ generate a free group and if their displacements at  $ x\in X $, that is $ d(x, ax) $ and $ d(x, bx) $, are small, then the entropy of $ X $ is big. Hence an upper bound on the entropy prevents the subgroup of $ \Gamma $ generated by the elements with small displacement at $ x $ to be algebraically "big". Nevertheless, even in the case of controlled entropy, free-subgroups or free-semigroups do exist but their generators must have large displacements.\\ 
Conversely, if the asymptotic displacements of two independent elements $a$ and $b$ of $ \Gamma $ are bounded from below, then there 
exist bounded powers of $a$ and $b$ which generate a free semi-group (see Proposition \ref{granddeplacementx}). This underlines the importance of computing a universal lower bound of the asymptotic displacements of all the torsion-free elements of the group (see Theorems \ref{intro:minorsystglobale0}
and \ref{Cat2}).

This can be made effective and the next theorem is in this spirit. 

A metric space $ (X,d) $ is said to be \textsl{geodesically complete} if all geodesic segments can be defined on $ \textbf{R} $. It is called a \textsl{Busemann space} (see \cite{papadop}, p. 187) if the distance $ d $ is convex, that is if the function $ d(c(t), c'(t)) $ is a convex function of $ t\in [0,1] $ for two geodesic segments $ c $ and $ c' $, affinely reparametrized. We have (see \ref{freesemigroup}),

\begin{theorem}\label{free-semi-group}
Let $ (X,d) $ be a connected, geodesically complete, non-elementary $ \delta $-hyperbolic  Busemann metric space. Let $ \Gamma $ be a torsion-free discrete subgroup of the isometry group of $ (X,d) $. We assume that $ \diam (\Gamma \backslash X) \leq D<+\infty $ and $ \Ent (X,d)\leq H $. For all pairs of elements $ a, b $ of $ \Gamma $, if the subgroup generated by $a$ and $b$ is not cyclic, then, for all integers $ p,q\geq S(\delta, H, D) $ one of the two semi-groups generated by $ \{a^p, b^q\} $ or by  $ \{a^p, b^{-q}\} $ is free.

Here $ S(\delta, H, D) $ is a function of $ \delta $, $ H $ and $ D $ which we describe precisely. 
\end{theorem}
In the same spirit we can minimise the normalised volume entropy on certain closed manifolds. Let us recall that, for a closed Riemannian manifold $(M,g)$, by abuse of language, we denote by ${\rm Ent}( M, g)$ the entropy of the metric space ($\widetilde M, d_{\tilde g}, dv_{\tilde g}$) where $\tilde g$ is the pulled back metric on the universal cover $\widetilde M$ of $M$. We then prove the following theorem (see Subsection \ref{appli:manifolds} for the necessary definitions and Theorem \ref{theo}).

\begin{theorem}\label{intro:theo} Let $M$ be a n-dimensional essential closed manifold, $n\geq 2$. Assume that the fundamental group $\Gamma$ of $M$ is non elementary, torsion free and admits a generating set $S$ such that $\Gamma$ is $\delta$-hyperbolic with respect to $S$ and satisfies 
${\rm Ent}(\Gamma, S) \leq H$ then, for every Riemannian metric $g$ on $M$,
$${\rm Ent}( M, g) ^n {\rm Vol}( M, g)\geq C(n,\delta, H) >0.$$
\end{theorem}

We extend this result to polyhedrons in Theorem \ref{theo:poly} of Subsection \ref{subsect:poly} (see the begining of Subsection \ref{subsect:poly} for the notion of Riemannian polyhedrons). This opens the applications to a wide range of metric spaces some of which are described right after Theorem \ref{theo:poly}, including $\CAT(0)$-square complexes with hyperbolic fundamental groups as well as higher dimensional constructions related to cube complexes.

One of the key tools used in proving the main results of this article is a Bishop-Gromov-like theorem which yields an explicit link between the algebraic and geometric aspects. The next theorem proves such 
a Bishop-Gromov inequality in the case of Gromov-hyperbolic metric spaces, where the hypothesis \lq \lq Ricci curvature bounded from below" is replaced by the much weaker (see subsection
\ref{comparaison}) hypothesis \lq \lq Entropy bounded from above". The following statement is a simplification of Theorem \ref{cocompact2} for the purpose of this introduction.

  \begin{theorem}\label{intro:cocompact2}
Let $(X , d) $ be a $ \delta$-hyperbolic metric space, for every proper action by isometries of a group $\Gamma$ on $(X,d)$ such that the diameter of $\Gamma \backslash X$ and the
entropy of $(X , d)$ are respectively bounded by $D$ and $H$, then, for every $x \in X$
  \begin{itemize}
    \item[(i)] for every $\Gamma$-invariant measure $\mu$ on $X$, for every $R\geq r\ge\frac{5}{2} (7 D + 4 \delta)$ one has,
$$\dfrac{\mu \left( \overline B_X \big( x , R \big)\right)}{\mu \big( B_X ( x , r )\big)} \le 3 e^{H D} \left( \frac{R}{r} \right)^{25/4}  \left( \frac{R}{r} \right)^{6 H D}  e^{6 H  (R - \frac{4}{5} r)} \,,$$
    \item[(ii)] for every $R\geq r\ge 10\,(D + \delta )$, the counting measure $ \mu^\Gamma_x = \sum_{\gamma \in \Gamma} \delta_{\gamma x}$ of the orbit $\Gamma x$ verifies the inequalities:
$$\dfrac{\mu_{x}^{\Gamma} \left(\overline B_X \big( x , R \big)\right)}{\mu_{x}^{\Gamma} \big( B_X ( x , r )\big)} < 3 \left( \frac{R}{r} \right)^{25/4} e^{6 H (R- \frac{4}{5} r)} \,.$$
  \end{itemize}
\end{theorem}

It is possible to reinterpret these Bishop-Gromov-like inequalities (i) and (ii) in terms of the doubling properties in the sense of Definitions \ref{doublefaible} (see the comments following the statement of Theorem \ref{cocompact2}). The second inequality is interesting since the counting measure concerns the algebraic properties of the group and its geometric action. It thus make the link between the two aspects of our study and it is somehow a good surprise that, despite the fact that the counting measure is the most primitive one in this context, strong results could be obtained. 

We also remind the reader that the classical Bishop-Gromov inequality for manifolds has been a revolutionary tool which, in particular, led to compactness as well as finiteness results. In a forthcoming paper (\cite{BCGS2}), we shall prove finiteness and compactness results for compact quotients of metric measure spaces satisfying a weak Bishop-Gromov inequality similar to the above theorem \ref{intro:cocompact2}, these results will be applied in particular to compact 
quotients of $\delta$-hyperbolic metric spaces with bounded entropy (see the chapter 3 of \cite{Be-Ga} for a report presenting these results).

Finally, we mention a result related to the thick-thin decomposition. Let us first give the definition of an interesting family of groups (see \ref{systoleaction}).

\begin{defi}\label{intro:systoleaction}
Given parameters $\delta_0\geq 0$  and  $\e'_0 > 0$, we denote by $\text{\rm Hyp}_{\rm thick} (\delta_0, \e'_0)$ the set of non virtually cyclic groups $\Gamma$ 
which admit a proper, possibly non co-compact, action by isometries on some $\delta_0$-hyperbolic metric space $ (X_0, d_0)$ such that every 
torsion-free $\g\in \Gamma\setminus\{e \}$ verifies $\ell(\g) \ge \e'_0$.
\end{defi}
Notice that the space $ (X_0, d_0) $ may depend on $ \Gamma$ and that non trivial examples of such groups are given in the present article. For the sake of simplicity we shall assume, in this introduction, that $\Gamma$ is torsion-free.

In subsection \ref{comparacyl}, we compare the class of groups $\text{\rm Hyp}_{\rm thick} (\delta_0, \e'_0)$ with the class of acylindrically hyperbolic groups
(with some normalization), proving that this last class is included in $\text{\rm Hyp}_{\rm thick} (\delta_0, \e'_0)$.

Now, if $ \Gamma $ acts on a metric space $ (X,d) $ by isometries, we define the \textsl{pointwise and global systoles} by,
 $$ \sys_{\Gamma} (x)=\inf_{\gamma\in\Gamma\setminus\{e\}} \{d(x,\gamma x)\} ,\quad \Sys_\Gamma (X) = \inf_{x\in X}\{\sys_{\Gamma} (x)\} . $$ 
The next theorem is a curvature-free version of the celebrated \lq \lq collar lemma" for metric spaces. It shows that, under the hypotheses, if one has a small loop at some point any other independent loop at this point is long. More precisely (see Theorem \ref{minorsystglobale}),

\begin{theorem}\label{intro:minorsystglobale}
Let $\delta_0\geq 0$, and $\e'_0,\,  H > 0$, there exists an integer $ n'_0$ depending on $\delta_0$ and $\e'_0$ only such that, for any (torsion-free) element $\Gamma$ of 
$\text{\rm Hyp}_{\rm thick} (\delta_0, \e'_0)$, for any proper 
action, by isometries preserving the measure, of $\Gamma$ on a connected metric measured space $ (Y , d, \mu) $ whose entropy is bounded from above by $H$ we have,  
 \begin{itemize}
    \item[(i)] (Collar Lemma)  if $y \in Y$ and $\sigma \in \Gamma^*$ verify $d(y, \sigma y)=\sys_\Gamma (y) \le \frac{1}{2 n'_0 H}$, then every 
$\gamma \in \Gamma$ which does not commute with $\sigma$ satisfies $d(y, \gamma y) \ge \frac{1}{2 H}\   \ln  \left(\dfrac{1 }{n'_0 H \sys_\Gamma (y)}\right) - \frac{1}{2} n'_0  \sys_\Gamma (y)\,,$
    \item[(ii)] if moreover $(Y,d)$ is path-connected, then
$$ \text{\rm Sys}_\Gamma (Y) \ge \dfrac{1}{2 n'_0  H}\, e^{- 4 H \diam (\Gamma \backslash Y)}\,.$$
\end{itemize}
\end{theorem}
Notice that if the quotient of $Y$ by $\Gamma$ is non compact then Inequality $(ii)$ is trivial. Theorem \ref{intro:minorsystglobale} is stated in a weak form for the sake of simplicity and the reader is referred to Theorem \ref{minorsystglobale} for a more general statement.

Let us summarize this result by the following sentence: \textsl{a thin tube is long and has simple topology}. More precise statements on the topology and the structure of thin tubes will be given in Subsection \ref{structuremince}. Notice that in the inequality $ i) $ the right hand side goes to $ +\infty $ logarithmically when $\sys_{\Gamma} (y)$ goes to $ 0 $, exactly like in the standard collar theorem. Furthermore, in Theorem \ref{intro:minorsystglobale}, 
the metric space  $(Y,d)$ is not assumed to be $ \delta $-hyperbolic. A striking case is when $Y$ is a manifold which carries a Riemannian metric $ g_0 $ of sectional curvature less than $-1 $, Theorem \ref{intro:minorsystglobale} then applies to any other Riemannian metric on $Y$ whose entropy is less than $ H $.

In order to make Theorem \ref{intro:minorsystglobale} effective we have to provide a lower bound of  $ \e_0' $ in terms of the data. Such lower estimates are given in this article in various situations. One particular case is when we consider a proper co-compact action by isometries on a $ \delta $-hyperbolic space for which $ \e_0' $ can be taken to be the global systole of the action. The next theorem gives such a bound (see \ref{Cat2}).

\begin{theorem}\label{intro:minorsystglobale0} 
Let $ (X,d) $ be a $ \delta $-hyperbolic, non elementary, geodesically complete, Busemann space whose entropy satisfies $ \Ent (X,d)\leq H $. Let $ \Gamma $ acting properly, co-compactly by isometries on this space such that $ \diam (\Gamma\backslash X)\leq D $. Then, for any torsion-free element $ \gamma $ of $ \Gamma\setminus\{e\} $ we have,
 $$ \ell (\gamma) > s_0(\delta, H, D) , $$ 
for $ s_0(\delta, H, D) $ a function of $ \delta $, $ H $ and $ D $ which we describe. If furthermore we assume that $ \Gamma $ is torsion-free, we get,
 $$ \Sys_\Gamma (X) > s_0(\delta, H, D) . $$ 
\end{theorem}
This last result is the main step in the proof of Theorem \ref{free-semi-group}  and, together with Theorem \ref{free-semi-group}, they are the key tools which lead to effectiveness in Theorem \ref{intro:minorsystglobale}.

Let us finish this introduction by stating a finiteness/compactness result for Riemannian manifolds whose Ricci curvature
is bounded from below (Theorems \ref{fini0} and \ref{compacite0}).

\begin{theorem}\label{intro:fini0}
Given $n \ge 2$, $D, K , i_0 > 0$ and $\delta_0 , \e'_0 > 0$, let us consider the set of Riemannian $n$-dimensional manifolds $(M,g)$ which verify the
following hypotheses:
  \begin{itemize}
    \item[(i)] the fundamental group $\Gamma_M$ of $M$ is torsion-free and belongs to $\text{\rm Hyp}_{\rm thick} (\delta_0, \e'_0)$,
    \item[(ii)] $\Ric_g \ge -(n-1) K^2 \cdot g$ and $\diam (M,g) \le D$,
    \item[(iii)] the injectivity radius of its Riemannian universal cover $(\widetilde M , \tilde g)$ is bounded from below by $i_0$.
  \end{itemize}
 Then, this set contains only finitely many differentiable structures and is a finite union of compacts for the $C^{0,s}$-topology (see Definition 
\ref{Cconvergence} and Theorem \ref{compacite0} for clarifications).
\end{theorem}
Applications of this last result are finiteness, and rigidity results for Einstein manifolds and existence of upper bounds for the possible values 
of their scalar curvature (see Section \ref{Einsteindiscrete}).

\medskip
We now describe the plan of this article. Section 2 contains the basic notations. Section 3 contains the definition of the various entropies and several versions of the doubling property which are compared according to their generality.  Section 4 is devoted to describing the techniques used to produce discrete free groups and free semi-groups of isometries of a $ \delta $-hyperbolic space, in particular an adapted version of 
the Ping-Pong method and the precise study of two kinds of Margulis constants which (when bounded from below by $16\delta$) guarantee the existence 
of free subgroups or free semigroups generated by elements with bounded displacements. 
Section 5 begins with one of the main tools running all over the paper: a Bishop-Gromov and doubling property for any $ \delta $-hyperbolic space which admits  a co-compact action of a group of isometries, see Subsection \ref{doublehyp}. This tool allows to prove lower bounds of the exponential growth of spaces and groups and several Margulis properties, which yield a lower bound of the asymptotic displacements of all the
torsion-free elements of the group.
In Section 6 we develop the idea which we call \textsl{transplantation of  Margulis properties}. \textsl{Grosso modo}, the underlying philosophy is that, if a discrete group acts properly by isometries on a Gromov-hyperbolic  space, then it inherits, from this action, algebraic properties which in turn translate into Margulis type properties when it acts on another metric measured space whose entropy is bounded. This is the section in which the interplay between algebraic and topological properties is the most enlightening. Section 7 is devoted to applications. We first show in Theorem \ref{theo} that closed manifolds which are essential and whose fundamental groups are hyperbolic (see the precise statement) have a minimal entropy, for a given volume, bounded away from zero by an explicit constant. We extend this result to polyhedrons in Subsection \ref{subsect:poly}. The end of this section concerns Einstein metrics on manifolds for which we prove compactness and finiteness results in our context. Finally, in Section 8, we recall the basic facts about Gromov-hyperbolic spaces and their isometries.

\section*{Acknowledgements} The authors thanks Misha Gromov for enlightening comments. We also wish to thank  Marc Bourdon, Frédéric Haglund, Damian Osajda, and St\'ephane Sabourau for detailed discussions and Emmanuel Breuillard for interesting exchanges. G\'erard Besson is supported by the ERC Advanced Grant 320939, GETOM and ANR CCEM. Sylvestre Gallot wish to thank the Dipartimento SBAI and the Istituto Guido Castelnuovo dell'Universit\`a  di Roma \lq \lq Sapienza" for their support.

\section{Basic definitions and notations}\label{notations}

\begin{defis}\label{proprediscr0} 
Let $ (X,d) $ be a metric space and $ \Gamma $ be a group acting by isometries on 
 $ (X,d) $,
\begin{itemize}
  \item[(i)] the space $ (X,d) $, is said to be \emph{\lq \lq proper"} if every closed ball is compact.
  \item[(ii)] the action of $ \Gamma $ is said to be  \emph{\lq \lq proper"} if, for at least one $ x \in X $ (and then  $ \forall x \in X $), and $ \forall R > 0 $, the set of $ \gamma \in \Gamma $ such that 
 $ d(x, \gamma  x) \le R $ is finite. 
  \item[(iii)] the action is said to be  \emph{\lq \lq discrete"} if it is faithful and if the image of $ \Gamma $ (by this action) in the isometry group of $ (X,d) $ is a discrete subgroup (for the topology of uniform convergence on compact sets).
\end{itemize}
\end{defis}

\begin{remark}\label{proprediscr} 
On a proper space $ (X,d) $, every faithful action is proper if and only if it is discrete.
\end{remark}

The proof of Remark \ref{proprediscr} will be given in Proposition \ref{discret1}.

Except for the results described in Section \ref{entropy&doubling}, 
the metric spaces and the actions that we study will be assumed to be proper. 

Classically, one defines the systole of a Riemannian manifold $ (Y,h) $ to be the infimum of the length of non-homotopically trivial loops.
If $ (Y,h) $ is viewed as the quotient of its universal cover $ (\widetilde Y , \tilde h) $ by the action of its fundamental group $ G $, the systole coincides with the invariant
 $ \inf_{x \in \widetilde Y}\left(\min_{\gamma \in G \setminus \{e\}}\  d(x, \gamma   x ) \right) $ of the action of $ G $ on $ (\widetilde Y , \tilde h) $.

We generalize this notion in the following definitions:
\begin{defi}\label{systoles} 
Let $ (X,d) $ be a metric space, let $ \Gamma $ be a group acting by isometries $ (X,d) $,
\begin{itemize}
  \item $ \forall x \in X $, the pointwise systole of this action is defined by
 $ \text{\rm sys}_\Gamma (x) := \inf_{g \in \Gamma \setminus \{e\}} d(x, g  x ) $,
  \item  the global systole of this action is defined by $ \text{\rm Sys}_\Gamma (X) := \inf_{x \in X} \text{\rm sys}_\Gamma (x) $,
  \item  the diastole of this action is defined by $ \dias_\Gamma (X) :=
\sup_{x \in X} \text{\rm sys}_\Gamma (x) $.
\end{itemize}
\end{defi}

The systole, in the classical sense, of a Riemannian manifold $ (Y, h) $ is then the global systole of the action of its fundamental group on the universal cover $ (\widetilde Y , \tilde h) $ of $ (Y, h) $, endowed with the Riemannian distance associated to the pulled-back Riemannian metric $ \tilde h $. 

Notice that in our situation, contrarily to the Riemannian case, the pointwise systole and hence the global systole of the action of a group $ \Gamma $ on a metric space $ (X,d) $ could vanish. One of our goals will be to find large classes of metric spaces and group actions such that the pointwise and global systoles of these spaces and actions are bounded by the same positive constant.

In the sequel, except in section \ref{transport}, for any $ \delta > 0 $,  we shall consider metric spaces which are $ \delta$-hyperbolic in the sense of Gromov. A definition will be given at the beginning of Subsection \ref{coherence}. Notice that this definition implies that the space is automatically geodesic and proper. On such a space we shall also consider any proper action by isometries of a group $ \Gamma $ (see Definitions \ref{proprediscr0}) which is assumed to be non-virtually cyclic. 

\begin{defi}\label{Gammaepsilon}
We note $ \Gamma^* $ the set $ \Gamma \setminus \{e\} $. We also note $ \Sigma_r (x) $ (resp. $ \widehat{\Sigma}_r (x) $) the set (which is finite when the action is proper)  
of $ \g \in \Gamma^*$ such that $ d(x , \g x) \le  r $ (resp.  $ d(x , \g x) < r $), and we note $ \Gamma_r (y) $ 
(resp. $ \widehat{\Gamma}_r (x) $) the subgroup of $ \Gamma $ generated by $ \Sigma_r (x) $ (resp. by 
 $ \widehat{\Sigma}_r(x) $).
\end{defi}

In the sequel we work in the general frame of  {\em metric measured spaces} $ (X,d,\mu) $.
We denote by $ B_X (x,r) $ (resp. $ \overline{B}_X (x,r) $) the \emph{open (resp. closed) ball} with centre $ x $ and radius $ r> 0 $ of the space $ (X,d) $.

A metric space $ (X,d) $ is said to be a \emph{length space} if, $ \forall (x,y) \in X \times X $, there exists a continuous rectifiable path from $ x $ to $ y $ and  if $ d(x,y) $ is the infimum of the length the paths joining $ x $ and $ y $.

In a length space, we call {\it geodesic}  the image of an interval $ I\subset \R $ by an isometric embedding $ c $ (i. e. satisfying $ d(c(t),c(s)) = |t-s|  ,\forall t,s \in I $). A geodesic is then, by definition, {\em minimizing}.
When we restrict $ I $ to be $ [a,b] $, or  $ [a, +\infty [ $, or $ ] - \infty , +\infty [ $, this geodesic will be respectively called a {\em geodesic segment}, a 
{\em geodesic ray}, or a {\em  geodesic line}.

We call {\em local geodesic} a map $ c $ from an interval $ I\subset \R $ in $ X $ which is locally minimizing\footnote{In the case of Riemannian manifolds the Riemannian geodesics are local ones.}, that is $ \forall t \in I $, there exists an $ \e > 0 $ such that $ c $ is a geodesic on the interval $ ]  t-\e , t+ \e [  \cap   I $.

We will often consider {\emph{geodesic spaces}. A metric space is said to be geodesic if any two points can be joined by at least one geodesic segment. Several geodesic segments may join two distinct points $x, y$: by abuse of notation we will denote by $ [x,y] $ anyone of these geodesic segments.

A geodesic metric space $ (X,d) $ is said to have the {\em property of geodesic extension} if, for all local geodesic $ c : [a,b] \to X $ ($ a \ne b $), there exists $ \e > 0 $ and a local geodesic $ c' : 
[a,b+ \e ] \to X $ which extends $ c $ (i.e. $ c'_{|_{[a,b]}} = c $). This space will be said to be {\em geodesically complete} if every local geodesic $ c : [a , b ] \to X $ ($ a \ne b $)  can be extended in a local geodesic $ \bar c :   ]- \infty , +\infty [ \to X $.
It is worth recalling that a \emph{complete} metric space $ (X,d) $, which furthermore is geodesic, has the geodesic extension property if and only if it is geodesically complete (cf. \cite{BH} Lemme II.5.8 (1) p. 208).

On these metric spaces a group $ \Gamma $ will act. We will only consider proper actions by isometries (see the definitions \ref{proprediscr0}), which implies, in particular, that $ \Gamma $ acts via a representation $ \varrho : \Gamma \f \text{Isom} (X,d) $ whose kernel is a finite and normal subgroup of $ \Gamma $ and whose image $ \varrho (\Gamma) $ is a discrete subgroup of the group $ \textrm{Isom} (X,d) $ of isometries of $ (X,d) $ (see Lemma \ref{reductisom}, whose results (i) and (ii) are valid on general metric spaces); this also implies that the stabiliser $ \text{\rm Stab}_\Gamma (x) $ in $ \Gamma $ of any point $ x $ is finite. Finally, the quotient space $ \Gamma \backslash X $ will be endowed with the {\em quotient distance} $ \bar d $ (see definition in Lemma \ref{autofidele} (i)). 

For such an action of a group $ \Gamma $ on a metric space $ X $, a  certain number of results that we prove will be called \emph{\lq \lq Margulis properties"}. This means that for each of the problems $ {\rm M}_i \ (i= 1, 2, 3, 4) $ stated below, we will try to compute explicitly universal constants (i. e. valid on the largest set $ {\cal M}_i $ of metric spaces $ X $, of groups $ \Gamma $ and of actions of  $ \Gamma $).

\begin{itemize}
   \item {\bf Problem $ {\rm M}_1 $:} Find a set $ {\cal M}_1 $ and a constant $ \varepsilon_1 >0 $ such that, $ \forall X \in {\cal M}_1 ,\ \forall x \in X$, 
$\Gamma_{\varepsilon_1} (x)$ is virtually nilpotent.
  \item {\bf Problem $ {\rm M}_2 $:} Find a set $ {\cal M}_2 $ and a constant $ \varepsilon_2 >0 $ 
such that, $\forall X \in {\cal M}_2,\ \text{\rm Dias}_\Gamma (X) \ge\varepsilon_2 $.
  \item {\bf Problem $ {\rm M}_3 $:} Find a set $ {\cal M}_3 $ and a constant $ \varepsilon_3 >0 $ 
such that $\forall X \in {\cal M}_3,\  \text{\rm Sys}_\Gamma (X) \ge\varepsilon_3 . $ 
  \item {\bf Problem $ {\rm M}_4 $:} Find a set $ {\cal M}_4 $ and constants $ \e_4  , C_0  , C_1 >0 $ satisfying the following property: $ \forall \varepsilon
\in \, ]0, \e_4 ] $, if $R(\e) = C_1  \ln \left(\dfrac{C_0}{\varepsilon} \right) $, then $\forall X \in {\cal M}_4$, $ \forall x \in X $, if there exists a torsion-free element 
$ \sigma\in\Gamma^* $ such that $ d(x, \sigma x) \le \varepsilon $, then $\Gamma_{R(\e)} (y)$  is a virtually cyclic subgroup containing $\sigma$ (notice that
$R(\e) \f +\infty$ when $\e \f 0$).
\end{itemize}

Among the group actions to which this work apply, we will consider the action of a discrete group $ \Gamma $, generated by a finite set $ \Sigma $, on its {\em Cayley graph} $ {\cal G}_\Sigma (\Gamma) $.
We denote by $ | \gamma |_\Sigma $ the {\em word metric} related to $ \Sigma $ and $ d_\Sigma $ the  associated \emph{algebraic distance} on $ \Gamma $ (i. e. $ d_\Sigma (\g, g) := |\g^{-1} g   |_\Sigma $) as well as the length distance on the graph $ {\cal G}_\Sigma (\Gamma) $. By abuse of language we use the same notation for these two distances. 

The measures $ \mu $ considered in this article are {\em Borel}, non-negative and non identically vanishing. The main examples are the following:
\begin{itemize}
    \item the {\em counting measure}  $ \# $, for discrete sets,
    \item the {\em orbital counting measure}  $ \mu_{x}^{\Gamma} = \sum_{\gamma \in \Gamma}   \delta_{\gamma   x} $ on the orbit $ \Gamma \cdot x $ of a point  $ x $, associated to a proper action of a group 
    $ \Gamma $ on a space $ X $, where $ \delta_y $ denotes  the Dirac measure at $ y $,
    \item the {\em $ 1$-dimensional measure induced on the Cayley graph} $ {\cal G}_\Sigma (\Gamma) $, given by the length of the edges. 
    \item the {\em  Riemannian measure} $ dv_g $ on a Riemannian manifold $ (X,g) $. 
\end{itemize}

The metric spaces which are our main concern are (except in section \ref{transport}) the {\em $ \delta $-hyperbolic spaces in Gromov sense}. We recall their definition and basic properties in the Appendix (section \ref{outils}). We only consider $ \delta $-hyperbolic spaces which are geodesic and proper without necessarily recalling it. A group $ \Gamma $ with a finite generating set $ \Sigma $ is said to be a {\em $ \delta $-hyperbolic group} if its Cayley graph $ {\cal G}_\Sigma (\Gamma) $ endowed with the distance $ d_\Sigma $ is $ \delta $-hyperbolic in Gromov sense.

For a $ \delta $-hyperbolic space (in Gromov sense) $ (X,d) $, we use the symbols $ \partial X $ to denote its  {\em ideal boundary} and $ L\Gamma $ to denote the {\em limit set} of a discrete group $ \Gamma $ acting by isometries on $ X $ (i. e. the set of accumulation points in $ \partial X $ of any orbit $ \Gamma\cdot x $). Finally an hyperbolic space or group is called {\em elementary} if its ideal boundary has at most two points; similarly, any action of a group $ \Gamma $ on an hyperbolic space $ (X ,d) $ whose limit set has at most two points is called {\em elementary}; see subsection \ref{nilpotents} for more informations about elementary groups or actions.

In Section \ref{transport}, we leave the realm of $ \delta $-hyperbolic spaces to study general metric spaces and focus on the type of groups whose actions on metric measured spaces still satisfy the Margulis properties which were proved (in sections \ref{sectionlibres} and \ref{deltahyp}) to be valid on $ \delta $-hyperbolic spaces. On all these spaces the key invariant, which is a guideline all along this article and replaces the curvature, is the {\em entropy of a metric measured space}. It will be defined, discussed and compared to other invariants in the next section.

\section{Entropy, Doubling and Packing Properties}\label{entropy&doubling}

\subsection{Entropies}\label{entropies}

\begin{defi}\label{recouvrant}
Let $ \Gamma $ be a group acting on $ (X ,d) $, we call
\emph{covering domain} a subset $ K \subset X $ such that 
 $ \bigcup_{\g \in \Gamma}\  \g K = X $ and  \emph{fundamental domain} a covering domain such that 
 $ \g K^o \cap K^o = \emptyset $ for all $ \g \in \Gamma^* $ (where  $ K^o $ denote the interior of $ K $).
The action of $ \Gamma $ is said to be \emph{co-compact} if there exists a compact covering domain.
\end{defi}
\begin{defi}\label{Entropies0}
The entropy of a metric measured space $ (X ,d, \mu) $ (denoted by $ \Ent (X ,d, \mu) $) is the lower limit (when $ R \to +\infty $)
of $ \dfrac{1}{R}   \ln \Big( \mu \big( B_X (x , R )\big) \Big) $. It does not depend of the choice of $ x $.
\end{defi}
This invariant, possibly infinite, gives an asymptotic hence weak information on the geometry of the metric space (see subsection \ref{comparaison}), nevertheless it becomes interesting
when there exists a group $ \Gamma $ acting properly by isometries on $ (X ,d) $ (and possibly co-compactly) and when we restrict ourselves to Borel measures $ \mu $ which are
invariant by this action. A particular role will be played by the \emph{counting measure $ \mu_{x}^{\Gamma} $ on the orbit $ \Gamma \cdot x $ of a point $ x $ }. Notice that in the co-compact case the entropy does not depend on the chosen $ \Gamma$-invariant measure, has shown in the

\begin{prop}\label{Entropies1} Let $ (X ,d) $ be a non compact metric space and $ \Gamma $ be a group acting properly and co-compactly on $ (X ,d) $ 
by isometries. For every non trivial measure $ \mu $ on $ X $ which is invariant by this action, if there exists some compact covering domain of finite $ \mu$-measure, 
then  $\Ent (X ,d, \mu) = \Ent (X ,d, \mu_{x}^{\Gamma}) $ for every $x \in X$.\\
If, furthermore, $ (X ,d) $ is a length space, then $ \Ent (X ,d, \mu) $ is the limit (when $ R \to +\infty $) 
of $ \dfrac{1}{R}   \ln \Big( \mu \big( B_X (x , R )\big) \Big) $.
\end{prop}
Following this proposition \emph{we shall use the notation $ \Ent (X ,d) $ instead of
 $ \Ent (X ,d, \mu) $ for this type of measures}.\\
The proof of Proposition \ref{Entropies1} relies on the two following lemmas.

\begin{lemma}\label{recouvrement} Let $ \Gamma $ acting properly on a metric space $ (X ,d) $, if the quotient $ \Gamma \backslash X $ is compact, then the space $ (X ,d) $ is proper and every closed ball of radius at least equal to the diameter of $ \Gamma \backslash X $ is a compact covering domain for this action.
\end{lemma}

\begin{proof} If $ \Gamma \backslash X $ is compact of diameter $ D $, for all $ R \ge D $, for all $ x \in X $ and for all sequence $ (y_n)_{n \in \N} $ of points in the closed ball $ \overline B_X (x , R) $,
 there exists a subsequence (denoted by $ (y_n)_{n \in \N} $) whose image by $ \pi $ converges in $ (\Gamma \backslash X , \bar d) $. By definition of the quotient distance $ \bar d $ 
(see Lemma \ref{autofidele} (i)) and since $ \overline B_X (x , R) $ is a covering domain, there exists a point $ y_\infty \in \overline B_X (x , R) $ and a sequence $ (\g_n)_{n \in \N} $ 
of elements of $ \Gamma $ such that $ d(y_n , \g_n   y_\infty ) \f 0 $ when $ n \to +\infty $. Consequently there exists $ N \in \N $ such that, for all $ n \ge N $, we have
 $$d(x , \g _n x) \le d(x, y_n) + d(y_n, \g_n   y_\infty ) + d(\g_n   y_\infty , \g_n   x) \le  2   R + 1 . $$
The action being proper, this implies that the sequence $ (\g_n)_{n \ge N} $ take finitely many values, and hence admits a constant subsequence equal to $ \g \in \Gamma $. There thus exists a subsequence of 
 $ (y_n)_{n \in \N} $ which converges towards $ \g y_\infty $. This shows that the closed ball of radius
 $ R \ge D $ are compact covering domains, and hence that every closed ball is compact. 
\end{proof}

\begin{lemma}\label{Entropies2} Let $ (X ,d) $ be a non compact metric space and $ \Gamma $ be a group acting properly and co-compactly on $ (X ,d) $ by isometries. 
Let  $ \mu $ be a $\Gamma$-invariant measure and $ K $ be a compact covering domain, then, $ \forall x \in X $, $ \forall R, R' $ such that $ 0 < R' < R $, we have
 $$\dfrac{ \mu \big(B_X(x, R') \big)}{\mu_{x}^{\Gamma}\big(B_X(x, R' + \diam (K))\big) } \cdot \mu_{x}^{\Gamma} \big( B_X(x, R-R')\big)  \le 
\mu \big( B_X(x, R)\big) \le \mu (K)  \cdot  \mu_{x}^{\Gamma} \big( B_X(x, R + \diam ( K) )\big) $$
\end{lemma}

\begin{proof} The proof which follows is a variation on Proposition 2.3 of \cite{Re}.
By the definition of a covering domain, there exists $ g\in \Gamma $ such that
 $ x \in g\cdot K $. For the sake of simplicity let $ K' := g\cdot K $, denoting by $ \widehat{\Sigma}_{r} (x) $ the set of 
 $ \gamma \in \Gamma $ such that $ d(x, \g x ) < r $ and $ D := \diam (K) $ the diameter of $ K $, 
the triangle inequality gives:
 $$\mu (B_X(x, R) \le \mu \left( \cup_{\g \in  \widehat{\Sigma}_{R+D} (x)}\  \g K'\right) \le \sum_{\g \in  
\widehat{\Sigma}_{R+D} (x)}
\ \mu ( \g K') = \mu ( K)\cdot \mu_{x}^{\Gamma} \big( B_X(x, R + D)\big)  . $$
To prove the first inequality of Lemma \ref{Entropies2}, recall that, if $ \mu_1 $ and $ \mu_2 $ are two Borel measures on $ X $, 
for all Borel set $ U \subset X $, we have
\begin{equation}\label{Fubini}
\int_U   \mu_2 \big( B_X (x , r)\big) \ d \mu_1(x) = \int_X   \mu_1 \big( B_X (x , r) \cap U \big) \ d \mu_2(x) \ ,
\end{equation}
indeed, the two members of this inequality are equal to 
 $ \int_X \int_X \mathbf 1_{[0 , r[} \big(d(x,y) \big)   \mathbf 1_{U} (x) \ d \mu_2(y)  \ d \mu_1(x) $.\\
Replacing $ \mu_1 $ by $ \mu_{x}^{\Gamma} $ and $ \mu_2 $ by $ \mu $, equality \eqref{Fubini}, 
and the fact that $ \mu \big( B_X ( \g  x , R') \big) = \mu \big( B_X (x , R')\big) $ (thanks to the $ \Gamma$-invariance 
of $ d $ and $ \mu $) yields
 $$\mu_{x}^{\Gamma}\big(B_X(x,R-R') \big)\cdot \mu \big( B_X (x , R') \big) =
\int_{B_X(x,R-R')}    \mu \big( B_X (z , R') \big)\ d \mu_{x}^{\Gamma}(z) = $$
 $$\int_{B_X(x,R)} \mu_{x}^{\Gamma}\big(B_X(y, R') \cap B_X(x,R-R') \big) \ d\mu (y) \le \mu_{x}^{\Gamma}\big(B_X(x,R' + D) \big)
\cdot \mu \big(B_X(x,R) \big) \ , $$
where the last equality follows from the fact that, by definition of a covering domain, for all $ y \in X $, there exists
 $ g \in \Gamma $ such that $ B_X(y, R') \subset B_X (g  x, R' + D) $. 
\end{proof}

\begin{proof}[End of the proof of Proposition \ref{Entropies1}:] Let us chose a compact covering domain $ K $ such that
 $ \mu (K) < + \infty $ and a point $ x\in X $, there exists a point $ x' \in K $ and an isometry $ g \in \Gamma $ such that $ x = g  x' $; we set $ D := \diam (K) $ and chose $ R' = D + 1 $, 
so that $ K \subset B_X(x', R') $. The fact that $ \mu_{x'}^{\Gamma}  \big( B_X(x', R + D)\big) < +\infty $ for all
 $ R> 0 $ (since the action is proper), input in the second inequality of Lemma \ref{Entropies2}, 
has two consequences: on the one hand $ \mu \big( B_X(x', R )\big) < +\infty $ for all $ R> 0 $, on the other hand,
if $ \mu (K) $ vanishes, we would have $ \mu \big( B_X(x', R)\big) = 0 $ for all $ R> 0 $ and, as increasing union of $ B_X(x',R) $, $ \mu (X) $ would vanish too, which would imply that $ \mu $ is trivial, 
this contradicts the hypothesis. We then have $ \mu (K) > 0 $, which implies that $ \mu (B_X(x', R')) > 0 $. Inputing these two positivity results and the finiteness of $ \mu (K) $ and of $ \mu \big( B_X(x', R' )\big) $ 
in the inequalities of  Lemma \ref{Entropies2}, we obtain that the lower limits 
(when $ R \to +\infty $) of $ \dfrac{1}{R}   \ln \Big( \mu \big( B_X (x' , R )\big) \Big) $, of
 $ \dfrac{1}{R}   \ln \Big( \mu_{x'}^{\Gamma} \big( B_X (x' , R )\big) \Big) $ and of 
 $ \dfrac{1}{R}   \ln \Big( \mu_{x}^{\Gamma}\big( B_X (x , R )\big) \Big) $ coincide. The last coincidence is a consequence of the fact that $ \mu_{ g  x'}^{\Gamma} = \mu_{ x'}^{\Gamma} $ 
and that $ \mu_{ x'}^{\Gamma} \big( B_X (  g  x' , R )\big) = 
\mu_{ x'}^{\Gamma} \big( B_X (  x' , R )\big) $. This proves that $ \Ent (X ,d, \mu) = \Ent (X ,d, \mu_{x}^{\Gamma}) $.

Let us now assume that $ (X ,d) $ is a length space, then the above properties and Lemma \ref{Entropies2} imply that, if $ \dfrac{1}{R}   \ln \Big( \mu_{x}^{\Gamma}
\big( B_X (x , R )\big) \Big) $ has a limit when $ R \to +\infty $, it is the same for
 $ \dfrac{1}{R}   \ln \Big( \mu \big( B_X (x , R )\big) \Big) $. The existence of such a limit is deduced  from the fact that the function $ R \mapsto 
\ln \Big( \mu_{x}^{\Gamma}\big( B_X (x , R + 2  D)\big) \Big) $ is non decreasing and sub-additive, see Property 2.5 of \cite{Re} for a complete proof. 
\end{proof}
Two classical examples:

\smallskip
-- The classical notion of \emph{Volume entropy of a closed Riemannian manifold} $ (M , g) $, 
is defined, following Definition \ref{Entropies0}, as $ \Ent (\widetilde M , d_{\tilde g} , dv_{\tilde g}) $, 
where $ (\widetilde M ,\tilde g) $ is the Riemannian universal cover of $ (M , g) $, $ d_{\tilde g} $ its Riemannian distance and $ dv_{\tilde g} $ its Riemannian measure.
In this case $ \Gamma $ is the fundamental group of $ M $ acting by isometric deck transformations on $ (\widetilde M ,\tilde g) $. By Proposition \ref{Entropies1}, one can replace, in this definition, $ dv_{\tilde g} $ by the counting measure
 $ \mu_{x}^{\Gamma} $ on any orbit, or by any other $\Gamma$-invariant measure $ \mu $.

\smallskip
-- The notion of \emph{algebraic entropy of a finitely generated group $ \Gamma $ with a finite generating set $ \Sigma $ }. In this case, the algebraic entropy of $ \Gamma $ related to $ \Sigma $ (also called {\em rate of exponential growth} or {\em critical exponent} of $ \Gamma $ related to $ \Sigma $) is denoted by $ \Ent  (\Gamma, \Sigma) $) and can be defined in two different ways as follows~: either as the entropy of the metric measured space $ (\Gamma , d_\Sigma , \# ) $, where $ d_\Sigma $ is the algebraic distance defined in Section \ref{notations}, or as the entropy of the metric measured space $ \left( {\cal G}_{\Sigma }(\Gamma),   d_\Sigma  , \mu \right) $,  where $ {\cal G}_{\Sigma }(\Gamma) $ is the Cayley graph of $ \Gamma $ associated to $ \Sigma $, and $ d_\Sigma $ and $ \mu $ are respectively the length distance and the $ 1$-dimensional measure induced on the graph (see Section \ref{notations}).

\smallskip
A link between these two notions is given by the following classical result. Denote by $ d_x $ the \emph{geometric pseudo-distance} 
defined on $ \Gamma $ by  $ d_x (\g , \g') := d(\g x, \g'   x) $, the balls of $ d_x $ are well defined on $ \Gamma $, one can then define the entropy of the pseudo-metric measured space $ (\Gamma , d_x, \#) $ 
as the lower limit, as $ R \f +\infty $, of $ \frac{1}{R}\; \ln \left(\# \{\g  :  d_x (e, \g) < R\}\right) $. Denoting by $ \text{\rm Stab}_\Gamma (x) $ the stabilizer of $ x $ in $ \Gamma $, the pseudo-distance induces a distance on $ \Gamma / \text{\rm Stab}_\Gamma (x) $ and the entropy of the metric measured space $ (\Gamma/\text{\rm Stab}_\Gamma (x) , d_x, \#) $ coincides with the entropy of $ (\Gamma , d_x, \#) $ when $ \text{\rm Stab}_\Gamma (x) $ is finite.

\begin{lemma}\label{comparentropi} {\rm  -- } Let $ \Gamma $ be any group which admits a finite generating set $ \Sigma $ and which acts properly (by isometries) on 
a metric space $ ( X , d )$. For all $\Gamma$-invariant measure $ \mu $ on $ X $, for all $ x \in X $, we have
 $$\Ent (X,d,\mu) \ge \Ent (X,d,\mu_{x}^{\Gamma}) = \Ent (\Gamma, d_x , \# ) \ge  \dfrac{1}{\Max_{\sigma \in \Sigma} 
  d(x, \sigma   x )} \, \Ent (\Gamma, \Sigma) . $$
\end{lemma}

\begin{proof} The definitions of $ \mu_{x}^{\Gamma} $ and of the pseudo-distance $ d_x $ implying that $ \mu_{x}^{\Gamma}\big( B_X (x , R )\big) = 
\# \left\{ \g  :  d_x (e , \g) < R \right\} $, the equality $ \Ent (\Gamma, d_x , \# ) = \Ent (X,d,\mu_{x}^{\Gamma}) $ follows. 
The triangle inequality implies that $ d(x, \g  x) \le M\cdot d_\Sigma (e , \g) $, where $ M := \Max_{\sigma \in \Sigma}   d(x, \sigma   x ) $, hence that
 $$\dfrac{1}{R}   \ln \Big( \mu_{x}^{\Gamma}\big( B_X (x , R )\big) \Big) \ge \dfrac{1}{R}   \ln \left( \# 
\left\{ \g  :  d_\Sigma (e , \g) < \frac{R}{M}\right\} \right) , $$
which, by taking the lower limit when $ R \f +\infty $, proves the last inequality of Lemma
\ref{comparentropi} (notice that this result is still valid when $ \Ent (X,d,\mu_{x}^{\Gamma}) = +\infty $).

The total measure of $ X $ being strictly positive, there exists $ R' > 0 $ such that $ \mu \big( B(x , R')\big) > 0 $. If there exists $R$ such that 
$ \mu \big( B(x , R)\big) = + \infty $, then $\Ent (X,d,\mu) = +\infty$ and the inequality $\Ent (X,d,\mu) \ge \Ent (X,d,\mu_{x}^{\Gamma}) $ is trivially verified, hence
we shall suppose that $\forall R ,\  \mu \big( B(x , R)\big) < + \infty $.
Formula \eqref{Fubini} (where $ \mu_1 $ is replaced by $ \mu_{x}^{\Gamma} $ and $ \mu_2 $ by $ \mu $) then yields, for all $ R>0 $,
 $$\mu_{x}^{\Gamma}\big(B_X(x,R) \big)\cdot \mu \big( B_X (x , R') \big) =
\int_{B_X(x,R)}    \mu \big( B_X (z , R') \big)\ d \mu_{x}^{\Gamma}(z) $$
\begin{equation}\label{mugammavsmu}
= \int_{B_X(x,R+R')} \mu_{x}^{\Gamma}\big(B_X(y, R') \cap B_X(x,R) \big) \ d\mu (y) 
\le \int_{B_X(x,R+R')} \mu_{x}^{\Gamma}\big(B_X(y, R')\big) \ d\mu (y)  ,
\end{equation}
If $ d(y , \Gamma x ) \ge R' $, then $ B(y,R') \cap \Gamma x = \emptyset $ and $ \mu_{x}^{\Gamma}
\big(B_X(y, R')\big) = 0 $; if  $ d(y , \Gamma x ) <  R' $, there exists $ g\in \Gamma $ such that 
 $ d(y , g x ) <  R' $ and the triangle inequality ensures that $ B_X(y, R') \subset B_X(g x , 2  R') $, and hence that
 $$\mu_{x}^{\Gamma} \big(B_X(y, R')\big) \le  \mu_{x}^{\Gamma} \big(B_X(g x , 2  R')\big) = 
\mu_{x}^{\Gamma} \big(B_X( x , 2  R')\big) . $$
The last equality follows from the $ \Gamma- $ invariance of the measure $ \mu_{x}^{\Gamma} $ and the fact that
 $ B_X(g x , 2  R') = g\big( B_X( x , 2  R')\big) $. Plugging these two estimates in Inequality
\eqref{mugammavsmu}, we obtain
 $$\mu_{x}^{\Gamma}\big(B_X(x,R) \big)\cdot \mu \big( B_X (x , R') \big) \le \mu_{x}^{\Gamma} \big(B_X( x , 2  R')\big)
\cdot \mu \big( B_X(x,R+R')\big). $$
Taking the logarithm of both sides, dividing by $ R+R' $ and taking the lower limit of both sides when $ R \f +\infty $, we deduce that $ \Ent (X,d,\mu) \ge \Ent (X,d,\mu_{x}^{\Gamma}) $, which ends the proof. 
\end{proof}

\subsection{Doubling and Packing Properties}

\begin{defis}\label{doublefaible}
Let $ C_0 > 1 $ and $ I \subset    ] 0 , + \infty [ $ be an interval. We consider a metric measured space $ (X,d , \mu) $ and a point $ x\in X $.
\begin{itemize}
    \item[(i)] $ (X, d, \mu) $ is said to satisfy the \emph{$ C_0 $-doubling for all balls centred at $ x $ and of radius $ r \in I $} if
\begin{equation}\label{defdouble}
\forall r \in I, \ \ \  0 < \mu  \big(B_X(x , r )\big) < +\infty 
 \ \ \ \text{and} \ \ \  \dfrac{\mu \big(B_X(x , 2  r )\big)}{\mu \big(B_X(x , r )\big)} \le C_0 .
\end{equation}
    \item[(ii)] $ (X, d, \mu) $ is said to satisfy the \emph{weak $ C_0 $-doubling around $ x $ at the scale $ r_0 $} if it satisfies the $ C_0 $-doubling for all balls centred 
at $ x $ and of radius $ r \in \left[ \frac{r_0}{2}   ,   2  r_0\right] $.
    \item[(iii)] $ (X, d, \mu) $ is said to satisfy the \emph{strong $ C_0 $-doubling around $ x $ at the scale $ r_0 $} if it satisfies the $ C_0 $-doubling for all balls 
centred at $ x $ and of radius $ r \in \left] 0   ,   2  r_0\right] $.
\end{itemize}
If Condition \eqref{defdouble} is satisfied for all $ x \in X $, we say, in case i), that $ (X,d, \mu) $ satisfies the \emph{ $ C_0 $-doubling for all balls of radius $ r \in I $}, in case ii), that it satisfies the \emph{weak $ C_0 $-doubling at scale $ r_0 $} and, in case iii), that it satisfies the \emph{strong $ C_0 $-doubling at scale $ r_0 $}.

In all cases $ C_0 $ is called the (doubling) \emph{amplitude} and $ r_0 $ the (doubling) \emph{scale}\footnote{Notice 
that for $ \lambda> 0 $, neither the scale nor the amplitude are changed when the measure is multiplied by $ \lambda $. When the distance is multiplied by  $ \lambda $ the amplitude remains unchanged but the scale is multiplied by $ \lambda $.}.
\end{defis}

Given a proper action (by isometries) of a group $\Gamma$ on $ (X,d ) $, the application of these definitions to the counting measure $ \mu_{x}^{\Gamma} $ of the orbit of $x$ (introduced after Definition \ref{Entropies0}) will be important in the sequel.
In Lemma  \ref{croissanceN} we shall show that, if there exists some $ \Gamma$-invariant measure which satisfies a doubling condition, then 
$ \mu_{x}^{\Gamma} $ satisfies the same doubling condition (thus this last condition is weaker). One difficulty though comes from the fact that, for all $ y \in X $ such that $ d(y , \Gamma \cdot x) \ge r_0$, $ \mu_{x}^{\Gamma} \big( B_X(y , r_0)\big) = 0 $. Consequently, the doubling condition \eqref{defdouble} does not make 
sense for balls centred at $ y $ and of radius $ r < r_0 $; consequently the weak $ C_0 $-doubling at a scale $ r_1 \le r_0 $ around 
$ y $ does not make sense. For this reason, in the sequel, doubling conditions for the measure $ \mu_{x}^{\Gamma} $ will only concern balls centred at $ x $ 
(or at any point of its orbit). 

Let us remark that it is not necessary to assume that the action of $ \Gamma $ is fixed point-free. Indeed, if the stabilizer or any point of the orbit of $ x $ is not trivial, it is finite by the properness assumption on the action. We could think of the counting measure as defined by $ \mu'_x (A) = \# \big( A \cap \Gamma  x \big) $ for any $ A \subset X $; we then have $ \mu_{x}^{\Gamma}= \# \big( \text{\rm Stab}_\Gamma (x) \big)\cdot \mu'_x $; this implies that $ \mu_{x}^{\Gamma} $ satisfies the $ C_0 $-doubling condition for all balls of radius $ r \in I $ centred at $ x $ if and only if $ \mu'_x $ does.

\emph{For this reason, in the sequel, when we shall consider doubling properties satisfied by the counting measure this will indifferently refer to $ \mu'_x $ or to 
$ \mu_{x}^{\Gamma} $.}

Finally, related to doubling conditions, the packing condition can be stated as follows,
\begin{defi}\label{packing} 
For all $ N_0 \in \N^* $ and $ r_0 > 0 $, a metric space $ (X,d ) $ is said to satisfy the \emph{packing condition with bound $ N_0 $ at scale $ r_0 $ }
if, for all $ x \in X $, the maximal number of disjoint balls of radius $ r_0 /2 $ included in the ball $ B_X(x, 9  r_0 ) $ is not greater than $ N_0 $.
\end{defi}

\subsection{Comparison between the various possible hypotheses~: bound on the entropy, doubling and packing conditions, Ricci curvature}\label{comparaison}

Let us start by a general comment. An upper bound on the entropy as well as the weak doubling condition for the counting measure on an orbit of the action of a group $ \Gamma $ makes sense on general metric spaces. On the other hand, a lower bound on the Ricci curvature concerns, in our context, mainly Riemannian manifolds. However, even if we restrict ourselves to Riemannian manifolds, where $ \Gamma $ is the fundamental group of a compact manifold acting by deck transformations on its universal cover, the comparison between all these conditions is roughly summarized as follows~:

\begin{compari}\label{comparaison1}
An upper bound on the entropy is a condition that is strictly weaker than the weak doubling (around $x$) of the counting measure of an orbit $ \Gamma x $ (cf. Lemma \ref{entropiebis} (ii)), which is itself strictly weaker than a packing condition at a similar scale (cf. Lemma \ref{packing1}), itself weaker than the weak 
doubling condition on the Riemannian measure (cf. Lemma \ref{packing2}), itself strictly weaker than the strong doubling at the same scale (cf. Lemma \ref{exemple} and its Corollary \ref{exemple1}), itself strictly weaker than a lower bound on the Ricci curvature (cf. Lemma \ref{Riccidouble}).\\
If we furthermore restrict ourselves to comparing the various weak doubling conditions, the smaller the scale, the stronger the condition (cf. Lemma \ref{entropiebis} (i)).
\end{compari}

It is worth noticing that, in the co-compact case, on a Riemannian manifold $ (X,g) $, for any scale $ r_0 $, there always exists a constant $ C:= C (r_0, X , g) $ such that the Riemannian measure satisfies the weak $ C$-doubling at the scale $ r_0 := r_0 (X,g)$. Similarly, there always trivially exists a constant $ H := H (X,g)$ which bounds from above the entropy of $ (X,g) $. However, the idea in this article is to compute these constants $C, r_0$ and $H$ in such a way that they are 
simultaneously valid on all the elements $(X,d)$ of a family of metric measured spaces (i. e. independently of $(X,d)$); a first condition will be said to be strictly 
weaker than a second one if the family of spaces satisfying the second one is strictly included in the family satisfying the first one and if the constants appearing in 
the definition of the first condition can be (explicitly) computed in terms of those appearing in the definition of the second condition. The next lemmas makes this philosophy precise.

\begin{lemma}[Entropy bounded from above versus  weak doubling of the counting measure]\label{entropiebis}
For every proper action (by isometries) of a group $\Gamma$ on a length space $ (X,d ) $ such that $ \Gamma \backslash X $ is compact with diameter $ \le D $, if 
there exists $ x_0\in X $ such that the counting measure $ \mu_{x_0}^{\Gamma} $ of its orbit satisfies the weak $ C_0 $-doubling around $ x_0 $ for at 
least one scale $ R_0 > D $, then
\begin{itemize}
    \item[(i)] $ \mu_{x_0}^{\Gamma} $ satisfies the weak $ C_0^{3 (1 +\left[\frac{2 R_1}{R_0}\right])} $-doubling around $ x_0 $ 
at every scale $ R_1 \ge  R_0 $,
    \item[(ii)] the entropy of $ (X,d ) $ is bounded from above by $ \dfrac{3}{R_0} \ln C_0 $.
\end{itemize}
\end{lemma}

In (ii) the condition \lq \lq entropy bounded from above" is strictly weaker than the weak doubling one, as can be seen when considering any lattice $ \Gamma $ 
acting on $ \R^n $. Indeed, $ \R^n $ has zero entropy however the $ C_0 $-doubling at scale $ R_0 $ is not any more satisfied when $ n $ is large enough.

\begin{proof} Let us recall that $ \widehat{\Sigma}_{r} (x_0) $ denotes the set of $ \gamma \in \Gamma $ such that $ d(x_0, \g   x_0 ) < r $ and  let $ R\geq R_0 $. 
 From equality \eqref{Fubini} (applied for two measures equal to $ \mu_{x_0}^{\Gamma} $), we deduce that
 $$\int_{B_X(x_0 , R)}   \mu_{x_0}^{\Gamma} \big( B_X (x , 4  R_0)\big) \ d \mu_{x_0}^{\Gamma} (x) =
 \int_{X} \ \mu_{x_0}^{\Gamma} \big(  B_X (x , 4   R_0) \cap B_X(x_0 , R) \big) \ d \mu_{x_0}^{\Gamma} (x) ; $$
this and the fact that $ x \mapsto\mu_{x_0}^{\Gamma} \big( B_X (x , 4  R_0)\big) $ is constant on the support of $ \mu_{x_0}^{\Gamma} $ yields
\begin{equation}\label{Fubini1}
 \mu_{x_0}^{\Gamma} \big( B_X(x_0 , R)\big) \cdot \mu_{x_0}^{\Gamma} \big(  B_X (x_0 , 4  R_0)\big) \ge \sum_{\g \in  \widehat{\Sigma}_{R+ R_0} (x_0)}   \mu_{x_0}^{\Gamma} 
\big(  B_X ( \g x_0 , 4   R_0) \cap B_X(x_0 , R) \big) .
\end{equation}

\begin{itemize}
  \item If $ d(x_0 , \g   x_0) \le |R - 4  R_0| $, then $ B_X (  \g   x_0 , 4   R_0) \cap B_X(x_0 , R) $ contains one of the balls $ B_X(x_0 , R) $ or 
$ B_X (  \g   x_0 , 4   R_0) $, thus $ \mu_{x_0}^{\Gamma} \big(  B_X ( \g   x_0 , 4   R_0) \cap B_X(x_0 , R) \big) \ge  \mu_{x_0}^{\Gamma} \big( B_X(x_0 , R_0) \big)$.
  \item If $ |R - 4  R_0| < d(x_0 , \g   x_0) \le R + R_0 $, since we are on a length space, for all  $ \e $ such that $ 0< \e < R_0 - D $, there exists a point $ y(\g) \in X $ (chosen for example on an almost minimizing path from $ x_0 $ to $ \g x_0 $, whose length is less than $ d(x_0 , \g x_0)+\e $) 
such that,
 $$d \big(x_0 , y(\g) \big) = \frac{1}{2} \left[ d(x_0 , \g x_0)+R-4R_0\right] \ \text{ and }\ 
d \big(\g x_0 , y(\g) \big) < \frac{1}{2} \left[ d(x_0 , \g x_0)-R+4R_0 + 2  \e\right] . $$
We easily check that $B_X \left( y(\g ),\frac{3  R_0}{2} - \e \right) $ is included in $ B_X \left( y(\g ), \frac{1}{2} \left[ R + 4  R_0 - d(x_0 , \g  x_0) \right] -\e \right)$,
itself included in $B_X (  \g x_0, 4 R_0) \cap B_X (x_0 , R)$.
Since there exists $ g \in \Gamma $ such that $ d \big(g  x_0 , y(\g ) \big) \le D < R_0 - \e $, the intersection
 $ B_X (  \g x_0 , 4   R_0) \cap B_X(x_0 , R) $ contains the ball $ B_X (  g x_0 ,  \frac{R_0}{2}) $.
\end{itemize}
For all $ \g \in  \widehat{\Sigma}_{R+ R_0} (x_0) $, we thus have $ \mu_{x_0}^{\Gamma} \big(  B_X ( \g   x_0 , 4   R_0) \cap B_X(x_0 , R) \big) 
\ge  \mu_{x_0}^{\Gamma} \big( B_X(x_0 ,  \frac{R_0}{2}) \big) $; plugging this inequality in \eqref{Fubini1} yields
$$\mu_{x_0}^{\Gamma} \big( B_X(x_0 , R)\big) \cdot \mu_{x_0}^{\Gamma} \big(  B_X (x_0 , 4  R_0)\big) \ge  
\mu_{x_0}^{\Gamma} \big(  B_X (x_0 , R+  R_0)\big) \cdot \mu_{x_0}^{\Gamma} \big( B_X(x_0 ,  R_0/2) \big) $$
and thus $ \dfrac{ \mu_{x_0}^{\Gamma} \big( B_X (x_0 , R +  R_0)\big)}{ \mu_{x_0}^{\Gamma} \big( B_X (x_0 , R)\big)} 
\le \dfrac{ \mu_{x_0}^{\Gamma} \big( B_X (x_0 , 4 R_0)\big)}{\mu_{x_0}^{\Gamma} \left[ B \left ( x_0   ,  \frac{1}{2}  R_0 \right)\right]} \le C_0^3 $,
and consequently, for all $ k \in \N^* $, $ \dfrac{ \mu_{x_0}^{\Gamma} \big( B_X (x_0 , R + k   R_0)\big)}{ \mu_{x_0}^{\Gamma} \big( B_X (x_0 , R)\big)} 
\le C_0^{3 k} $, which implies on the one hand that 
$ \dfrac{ \mu_{x_0}^{\Gamma} \big( B_X (x_0 , 2  R)\big)}{ \mu_{x_0}^{\Gamma} \big( B_X (x_0 , R)\big)} \le C_0^{3 (1 + \left[\frac{R}{R_0}\right])} $ 
and ends the proof of (i), and on the other hand that the entropy is bounded above by $ \dfrac{3}{ R_0} \ln (C_0) $, which ends the proof of (ii) (recall that the entropy does not depend on the choice of the measure by Proposition \ref{Entropies1}). 
\end{proof}

The lemma that follows describes the basic tools to compare doubling and packing.

\begin{lemma}\label{compardoubling}
For every proper action (by isometries) of a group $\Gamma$ on a metric space $ (X,d ) $, for every $ x \in X $, every integer $ k \ge 5 $ and every $ r > 0 $, we 
denote by $ N_k^\Gamma (r) $ the maximal number of disjoint balls of radius $ \frac{r}{2} $, centred on points of the orbit $ \Gamma \cdot x $ that can be included 
in the ball $ B_X (x , k \frac{r}{2})$. Then, for every $ \Gamma $-invariant measure $ \mu $, we have
 $$\dfrac{\mu_{x}^{\Gamma} \left(\overline B_X \big( x ,  (k-1)  \frac{r}{2}\big)\right)}{\mu_{x}^{\Gamma} 
\big( B_X ( x ,  r)\big)} \le N_k^\Gamma (r) \le
\dfrac{\mu \left(B_X \big( x ,  k  \frac{r}{2}  \big)\right)}{\mu \left(B_X \big( x ,   \frac{r}{2} \big)\right)} \ . $$
\end{lemma}

\begin{proof} For sake of simplicity we set $ N := N_k^\Gamma (r) $; denote by $ \gamma_1 x, \ldots ,\gamma_N  x $ the centres of a maximal packing of the ball $ B_X (x , k \frac{r}{2} ) $ by
disjoint balls of radius $ \frac{r}{2} $ centred on points of $ \Gamma \cdot x $. 
By $ \Gamma $-invariance of $ \mu $ and of the distance, $ \mu \left( B_X(\gamma_i \, x , \frac{ r}{2}) \right) = \mu \left( B_X(x , \frac{r}{2}) \right) $ for all $ i $, and consequently $ N \cdot  \mu \left(B_X(x , \frac{r}{2}) \right) \le \mu \left( B_X(x , k\frac{r}{2}) \right) $; this proves the right inequality of Lemma \ref{compardoubling}.\\
If there exists $ \gamma \in \Gamma $ such that $ \gamma   x \in \overline B_X \left( x   ,   (k-1) \frac{r}{2}\right) \cap \Gamma \cdot x $ and that 
$\g x \notin \cup_{i = 1}^N B_X \left( \gamma_i x, r\right)$, then the balls of radius $ \frac{r}{2} $, centred on $ \gamma_1 x, \ldots , \gamma_N (x) , \gamma x $ 
are disjoint and included in $ B_X (x , k  \frac{r}{2}) $. This contradicts the maximality of the packing by the balls 
$\left( B_X(\gamma_i \, x , \frac{ r}{2}) \right)_{1 \le i \le N}$; it follows that 
$$\overline B_X \left( x   ,   (k-1) \frac{r}{2}\right)  \cap \Gamma \cdot x \  \subset  \ \cup_{i=1}^N   
\Big(B_X \left( \gamma_i  x   ,   r \right) \cap \Gamma \cdot x \Big) , $$
and thus that $ \# \left( \overline B_X \left( x   ,   (k-1) \frac{r}{2}\right)  \cap \Gamma \cdot x   \right) \le 
N\ \# \left( B_X \left( x   ,   r \right) \cap \Gamma \cdot x   \right) $.
This gives 
 $$\dfrac{\mu_{x}^{\Gamma} \left(\overline B_X \left( x   ,   (k-1) \frac{r}{2}\right)\right)}{\mu_{x}^{\Gamma} 
\left(B_X ( x ,  r )\right)} =
\dfrac{\# \left(\overline B_X \left( x   ,   (k-1) \frac{r}{2}\right) \cap \Gamma \cdot x \right) }{\# 
\left( B_X \left( x   ,   r \right) \cap \Gamma \cdot x   \right) }\le N \ ,$$
and this ends the proof.
\end{proof} 

\begin{lemma}[Weak doubling of the counting measure versus packing]\label{packing1} 
Let $ (X,d ) $ be a metric space satisfying the packing condition with bound $ N_0 $ at scale $ r_0 $ and  $ \Gamma $ a group acting on  $ (X,d ) $ properly and by isometries. For all $ x\in X $, the counting measure
 $ \mu_{x}^{\Gamma} $ on the orbit $ \Gamma \cdot x $ satisfies the weak $ N_0 $-doubling at scale  $ 2 r_0 $ around $ x $.
\end{lemma}

\begin{proof} Recalling that $ N_{18}^\Gamma (r_0) $ is the maximal number of disjoints balls of radius $ \frac{r_0}{2} $ (centred at points of 
$ \Gamma \cdot x $) included in $ B_X (x , 9 r_0 ) $, Definition \ref{packing} gives $ N_{18}^\Gamma (r_0)  \le N_0 $. Then the first inequality of Lemma \ref{compardoubling} (with $ k = 18 $) implies that, for all $ r \in \left[ r_0 , 4 r_0 \right] $,
 $$\dfrac{\mu_{x}^{\Gamma} \left[B_X \left( x ,  2  r \right)\right]}{\mu_{x}^{\Gamma} \left[B_X \left( x ,  r \right)\right]} \le
\dfrac{\mu_{x}^{\Gamma} \left[B_X \left( x ,  \frac{17}{2} r_0  \right)\right]}{\mu_{x}^{\Gamma} \left[B_X \left( x ,  r_0 \right)\right]} \le N_{18}^\Gamma (r_0)  \le N_0 . $$
\end{proof} 

The series of examples that follow shows that these various conditions are strictly different.

\small

\begin{exas}\label{exemplestrict}
\emph{We choose a scale $ r_0 > 0$. In each of the following examples we construct a sequence of Riemannian manifolds
 $ (X_k , g_k) $ and of groups  $ \Gamma_k $ acting properly and isometrically on them, such that  the maximal number $ N^k (r_0) $ of disjoint balls
 of radius $ r_0/2 $ which are included in a ball of $ (X_k , g_k) $ of radius  $ 9 r_0 $ goes to infinity with $ k $. This shows that, for a given $ N_0 $ and for $ k $ 
sufficiently large,} $ (X_k , g_k) $ does not satisfy the packing condition with bound $ N_0 $ at scale $ r_0 $. On the other hand we show that, for every $k$, 
the counting measure $ \mu_{x_k}^{\Gamma_k} $ of the orbit of this action satisfies the weak $ C_0 $-doubling at scale $ 2r_0 $ (where $ C_0 $ and $ r_0 $ 
are independent of $ k $).
{\em 
\begin{itemize}
    \item [(1)] Let $ \Gamma_k := \frac{1}{k}\cdot \Z \times k \cdot \Z^{k-1} $ be a lattice acting on 
 $ (\R^k, can.) $ and $ x \in \R^k $, then the measure $ \mu_{x}^{\Gamma_k} $ satisfies the weak $ 2 $-doubling condition at scale $ 2r_0 $ (for any $ r_0>0 $) around $ x $ as soon as $ k > 8r_0 $. On the other hand, $ N^k (r_0) $ goes to infinity with the dimension $ k $.

\item [(2)] This example generalizes to the case where each $ (X_k,g_k) $ is the universal Riemannian cover of a closed Riemannian manifold 
$ (\overline X_k, \bar g_k) $ whose sectional curvature is between $-1 $ and $-K^2 $ and $ \Gamma_k $ is its fundamental group. A variation on Margulis Lemma 
shows the existence of $ \e_0 > 0 $ (independent of $k$) such that, in each connected component of the set of points where the injectivity radius is $ < \e_0 $,
there exists a periodic geodesic $ \bar c_k $ of length less than $ 2 \e_0 $ and whose homotopy class generates $ \widehat{\Gamma}_{2 \e_0} (x_k) $, where $ x_k $ is the lift of a point $ \bar x_k\in\bar c_k $. Since
 $ \widehat{\Gamma}_{2 \e_0} (x_k) $ acts by translation on the lift $ c_k $ of $ \bar c_k $ passing through $ x_k $ and $ d(x_k, \g x_k) \ge 2 \e_0 $ for every $ \g \in \Gamma_k \setminus \widehat{\Gamma}_{2\e_0} (x_k) $, the measure
 $ \mu_{x_k}^{\Gamma_k} $ satisfies the weak $ 2 $-doubling condition at scale $ r_0 $ for each  $ r_0 \le \e_0/2$. On the contrary, choosing any $x_k \in X_k$, when 
$k \f + \infty$, if $\dim (X_k) \f + \infty$, then the packing parameter $ N^k (r_0)$ goes to $ + \infty $.

\item [(3)] In this example, the sequence $ (X_k , g_k)_{k \in \N^*} $ is constructed so that  the dimension of $ X_k $ is fixed equal to $ n $ and that the topology of 
$ X_k $ becomes more and more complicated when $ k \f +\infty $; this complexity can be estimated (for example) by the Euler characteristic; we fix the scale $ r_0 $ 
to be any positive number (independent of $k$).\\
Let $ (\overline X, \bar g) $ be a Riemannian manifold, $ (X,g) $ its Riemannian universal cover and $ \Gamma $ its fundamental group. For $ \e > 0 $ small enough,
let $\overline N(\e)$ be the maximal number of disjoint balls of radius $ \e $ that can be included in $ (\overline X, \bar g) $ and let 
$\bar x_1 , \ldots \bar x_{\overline N(\e)}$ be the centres of these balls, then these centers are a $2 \e$-lattice in $\overline X$:
indeed, if there exists $ \bar x \in \overline X $ such that $ \bar x \notin \cup_{i = 1}^{N (\e) } B_{\overline X} \left( \bar x_i, 2 \e\right)$, then 
$B_{\overline X} \left( \bar x, \e \right) \cap B_{\overline X} \left( \bar x_i, \e \right) = \emptyset$ for every $i$, in contradiction with the maximality of the 
packing $\left( B_{\overline X} \left( \bar x_i, \e\right) \right)_{i = 1}^{N (\e) } $, thus $ \overline X \subset  \cup_{i = 1}^{N (\e) } B_{\overline X} 
\left( \bar x_i, 2 \e\right)$.
We modify the metric $ \bar g $ on each ball $ B_{\overline X} (\bar x_i , 2 \e^2) $ so that the new metric, still denoted by $ \bar g $, is flat on 
$ B_{\overline X} (\bar x_i , \e^2) $. The lattice $ \overline R_\varepsilon $ lifts as a $ 2\e $-lattice $ R_\varepsilon = \{x_i\}_{i\in I} $ of $ (X, g) $, globally 
$ \Gamma $-invariant.

Let $ (Y,h) $ be any closed Riemannian manifold, whose diameter is equal to $ 2 r_0 $ and two points $ y, y' $ such that $ d_Y (y, y') = 2r_0 $. We modify the metric, as above, so that it becomes flat in a neighbourhood of $ y $; the ball $ B_Y (y ,\e^3) $ is then isometric to each ball $ B_{X} (x_i , \e^3) $, which ensures that the boundaries of the balls $ B_{X} (x_i , \e^3) $ are $ B_Y (y ,\e^3) $ are both isometric to the Euclidean sphere $ \mathbb S^{n-1} \left(\e^3\right) $ of radius $ \e^3 $. 
We denote by $ (Y',  h') $ the Riemannian manifold with boundary obtained\footnote{The Riemannian metric $ h' $  is only piecewise $ C^1 $ but it could easily be smoothed out.} by gluing a cylinder $ C_\varepsilon := [0 , \e^{2}] \times \mathbb S^{n-1} \left(\e^3\right) $ to $ Y \setminus B_Y (y ,\e^3) $, identifying 
$ \partial  B_Y (y ,\e^3) $ with $ \{\e^2\} \times \mathbb S^{n-1} \left(\e^3\right) $. Let us consider a family $ \left(Y'_i \right)_{i \in I} $ of copies of $ (Y', h') $ and 
let us glue each $ Y'_i $ to $ X \setminus \left( \bigcup_{i \in I} B_{X} (x_i , \e^3)\right) $ identifying $ \partial Y'_i = \{0\} \times \mathbb S^{n-1} \left(\e^3\right) $ 
with the connected component $ \partial B_{X} (x_i , \e^3) $ of the boundary of $ X \setminus \left( \bigcup_{i \in I} B_{X} (x_i , \e^3)\right) $.
If we choose $ \e = \e_k \to 0 $ we then obtain Riemannian manifolds called $ (X_k, g_k) $.

We now consider maps $ f_k : X_k \to X $, which  send $ Y'_i $ onto the ball $ B_{X} (x_i , \e_k^3) $, contracting $ Y\subset Y'_i $ on $ x_i $, sending the generatrices of the cylinder $ C_{\varepsilon_k}\subset Y'_i $ onto the rays of the ball and such that $ f_k $ restricted to $ X \setminus \left( \bigcup_{i \in I} B_{X} (x_i , \e_k^3)\right) $ is the identity.  
The map $ f_k :( X_k , g_k)\to (X,g) $ is then contracting and, furthermore, there exists a sequence $ \eta_k $ going to zero with $ k $ such that 
\begin{equation}\label{Hausdorff}
\forall  x , z \in X \setminus \left( \cup_{i \in I} B_{X} (x_i , \e_k^3)\right),
\ \ d_X (f_k(x),f_k(z) ) = d_X (x,z ) \ge (1+ \eta_k)^{-1}  d_{X_k}(x,z) - \eta_k \  .
\end{equation}
Let $ x \in X \setminus \left( \bigcup_{i \in I} B_{X} (x_i , \e_k^3)\right) $ and $J := \{i \in I : x_i \in B_X (x, 6  r_0)\}$. The number of elements of $ J $, denoted by 
$\widetilde N_{\e_k}(r_0) $, is the number of elements of $ 2\e_k $-lattice in a fixed ball, thus it goes to infinity when $ \e_k \f 0 $, hence when $ k \f +\infty $.\\
For every $ i \in J $, we have $ \ Y'_i \subset B_{X_k} (x, 9  r_0) $: indeed, for any $ z \in  Y'_i $, there exists $ z' \in \partial B_{X} (x_i , \e^3) $ such that $ d_{X_k} (z,z') \le 2  r_0 + \e_k^2 $ and Inequality \eqref{Hausdorff} implies that $d_{X_k} (x,z') \le (1+ \eta_k) \big( d_{X} (x , z') + \eta_k \big)$, which yields
$$d_{X_k} (x,z') \le(1+ \eta_k) ( d_{X} (x,x_i) + \e_k^3 + \eta_k ) \le (1+ \eta_k) (6 r_0+ \e_k^3 + \eta_k )< \frac{13}{2} r_0 \, .$$
From the three last inequalities and the triangle inequality, we deduce that $ d_{X_k} (x,z) <  \frac{13}{2} r_0 + 2  r_0 + \e_k^2< 9  r_0 $, hence that 
$ Y'_i \subset B_{X_k} (x, 9  r_0) $. As each $ Y'_i $ contains a ball of radius $ r_0 $, and as $ Y'_i \cap Y'_j = \emptyset $ when $i\ne j$, the maximal number of 
disjoint balls of radius $ r_0 $ included in $ B_{X_k} (x, 9 r_0) $ is at least equal to $ \widetilde N_{\e_k}(r_0) $), hence it goes to infinity when $ k \f +\infty $.
Hence, \emph{for any choice of the constants $N_0$ and $r_0$, for $ k $ sufficiently large, $ (X_k , g_k) $ does not satisfy the packing condition with bound $ N_0 $ at 
scale $ r_0 $.} 

On the other hand, the above construction being $ \Gamma $-invariant, $ \Gamma $ also acts by isometries on $ (X_k,g_k) $ and we have $ f_k (\g x) = \g f_k(x) $ for all $ \gamma \in \Gamma $ and $ x \in X_k $;
the application $ f_k $ being contracting, we deduce that $ \{ \g :   \g x\in B_{X_k} (x , r)\} \subset \{ \g :  \g f_k(x) \in B_{X} (f_k (x), r)\} $. The manifold 
$ (X,g) $ (and the action of $\Gamma$ on it) being fixed, there exists 
 $ C_0 > 1 $ such that the counting measure $ \mu_{f_k(x)}^{\Gamma} $ on the orbit of $ f_k(x) $ under the action of $ \Gamma $ satisfies $ \mu_{f_k(x)}^{\Gamma} \left(B_{X} (f_k (x), 8  r_0) \right) \le C_0 $. 
The previous inclusion shows that the counting measure  $ \mu_{x}^{\Gamma} $ of the orbit of $ x\in X_k $ 
for the action 
satisfies $ \mu_{x}^{\Gamma} \left( B_{X_k} (x , 8  r_0)\right) \le C_0 $, which gives,
for all $ r \in [r_0 , 4r_0] $,
 $$\dfrac{\mu_{x}^{\Gamma} \left( B_{X_k} (x , 2  r)\right)}{\mu_{x}^{\Gamma} \left( B_{X_k} (x , r)\right)} \le
\mu_{x}^{\Gamma} \left( B_{X_k} (x , 2  r)\right) \le C_0  . $$
\end{itemize}
}
\end{exas}

\normalsize

\begin{lemma}[Packing versus weak doubling]\label{packing2}
Let $ (X,d ) $ be a metric space, if there exists a measure $ \mu $ which satisfies the weak $ C_0 $-doubling for every ball of radius $ r \in \left[ \frac{r_0}{2} , 
9r_0 \right] $ centred at \emph{every point}\footnote{The fact that the measure must satisfy the $ C_0 $-doubling condition for \emph{all} balls of radius 
$ r \in \left[ \frac{r_0}{2} , 9 r_0 \right] $, imposes that for all $ x \in X $, the ball $ B_X(x, \frac{r_0}{2}) $ has positive measure. This, in particular, excludes the 
counting measure of an orbit of a group $ \Gamma $ acting properly by isometries when the diameter of $ \Gamma \backslash X $ is greater than $ r_0 $.} $ x\in X $, 
then $ (X,d ) $ satisfies the packing condition with bound $ C_0^6 $ at scale $ r_0 $.
\end{lemma}

\begin{proof} Let $ x \in X $ and consider a packing of $ B_X(x,  9r_0 ) $ by $ N $ disjoint balls $ B_X \left(x_{i}, \frac{r_0}{2} \right) $. Amongst these balls, let us denote by
 $ B_X \left(x_{i_0}, \frac{r_0}{2} \right) $ one which has minimal measure, we then have 
 $$N\cdot \mu \Big( B_X \left(x_{i_0}, \frac{r_0}{2} \right)\Big) \le \mu \big( B_X(x,  9r_0 ) \big) \le 
\mu \big( B_X \left(x_{i_0}, 18  r_0 \right)\big) . $$
We thus get $ N \le \dfrac{\mu \big( B_X \left(x_{i_0}, 18r_0 \right)\big) }{\mu \Big( B_X \left(x_{i_0}, 
\frac{r_0}{2} \right)\Big)} \le C_0^6 . $ 
\end{proof}

Lemmas \ref{packing1} and \ref{packing2} asserts that a weak doubling condition on the counting measure of the orbit of a point is weaker than the weak doubling condition for any other measure around every point. It is however interesting to give a pointwise and precise result, this is the

\begin{lemma}[Counting measure vs any measure]\label{croissanceN} 
For any proper action (by isometries) of a group $ \Gamma $ on a proper metric space $ (X,d) $, for any $ x\in X $,  if there exists a $ \Gamma $-invariant 
measure which satisfies the $ C_0 $-doubling for all balls centred at $ x $, of radius $ r \in \left[\frac{1}{2} r_0 , r_1\right] $ (where $ r_1 \ge 
\frac{5}{4} r_0 $), then the counting measure $ \mu_{x}^{\Gamma} $ of $ \Gamma \cdot x $ satisfies 
$\dfrac{\mu_{x}^{\Gamma} ( B_X(x , 2r )) }{ \mu_{x}^{\Gamma} ( B_X(x , r ) )} \le C_0^3 $ for every $ r \in \big[r_0  , \frac{4}{5} r_1\big] $.
\end{lemma}

\begin{proof} 
Lemma \ref{compardoubling}, where we choose $ k = 5 $, ensures that, for all 
 $ r \in \big[r_0 , \frac{4}{5}r_1\big] $, we have
$$\dfrac{\mu_{x}^{\Gamma} \left[B_X \left( x ,  2  r \right)\right]}{\mu_{x}^{\Gamma} \left[B_X \left( x ,  r \right)\right]} 
\le \dfrac{\mu ( B_X(x , \frac{5}{2}   r )) }{ \mu ( B_X(x , \frac{r}{2} ) )} \le \dfrac{\mu ( B_X(x ,  r )) }{ \mu ( B_X(x , \frac{r}{2} ) )}\cdot 
\dfrac{\mu ( B_X(x , \frac{5}{4}   r )) }{ \mu ( B_X(x , r ) )} \cdot  \dfrac{\mu ( B_X(x , \frac{5}{2}   r )) }{\mu ( B_X(x , \frac{5}{4}   r )) }\le C_0^3 \, ,$$
since $ \frac{r}{2},  r $ and $ \frac{5}{4}r $ belong to the interval $ \left[\frac{1}{2}r_0 , r_1\right] $. 
\end{proof}

\begin{lemma}[Strong doubling versus Ricci curvature bounded below]\label{Riccidouble}
For all $ K \ge 0 $ and $ r_0 > 0 $, we set $ C_0 = 2^n \big(\cosh (K r_0) \big)^{n-1} $. Let $ (X,g) $ be a complete Riemannian manifold such that its Ricci curvature satisfies, $ \text{\rm Ric}_g \ge - (n-1)  K^2 \cdot g $, then $ (X,g) $ satisfies the $ C_0 $-doubling condition for all balls of radius $ r \in \left] 0 , r_0 \right] $ (independently of their centre).
\end{lemma}

\begin{proof} Bishop-Gromov's comparison theorem implies that, for all $ x \in X $ and all $ r \le r_0 $,
 $$\dfrac{\Vol_g  B_X(x , 2r )}{\Vol_g B_X(x , r )} \le \dfrac{\int_0^{2r} \left(\frac{1}{K} \sinh (K t)\right)^{n-1} dt}
{\int_0^{r} \left(\frac{1}{K}\sinh (K t)\right)^{n-1} dt} = 2 \cdot \dfrac{\int_0^{K r} \left(\sinh (2 t)\right)^{n-1} dt }{\int_0^{K r} \left(\sinh (t)\right)^{n-1} dt } \le 2^n (\cosh K r)^{n-1} . $$
\end{proof}

It is easy to find examples of sequences $ \big( X , g_k \big)_{k \in \N} $ of Riemannian manifolds which all satisfy the strong $ C_0 $-doubling condition at scale $ r_0 $ (where $ C_0 $ and $ r_0 $ do not depend on $ k $) 
and, nevertheless, such that the minimum of their Ricci curvature goes to $- \infty $ when $ k $ goes to $ +\infty $. One example can be obtained by gluing two copies of $ \R^n \setminus B^n $ on their boundary $ \mathbb S^{n-1} $ by constructing, on the resulting space, a sequence of smooth Riemannian metrics converging towards the singular metric obtained by endowing each copy with its Euclidean metric.

\medskip
It is clear (from  Definitions \ref{doublefaible} (ii) and (iii)) that the strong $ C_0 $-doubling condition is stronger than the weak one 
(at the same scale, around the same point, with the same amplitude). In fact, it is strictly stronger by the following results:

\begin{lemma}[Weak vs strong doubling]\label{exemple}
Let $ (X, g) $ be a complete Riemannian manifold, of dimension $ n \ge 2 $, whose Ricci curvature satisfies $ \text{\rm Ric}_g \ge - (n-1)  K^2 \cdot g $;
fix any $ r_0>0 $ and set $ C_0 = 1 + 2^{2 n} \big(\cosh (K r_0) \big)^{2n-2} $.
Let $ (Y_k , h_k) $ be any sequence of closed $n$-dimensional Riemannian manifolds, whose diameters and volumes go to zero when $ k \f +\infty $.
Let $ X_k=X\# Y_k $ be the connected sum of $ X $ with $ Y_k $, then, the metric $ g_k $ on each $ X_k $ being obtained by gluing $g$ and $h_k$, the corresponding 
 Riemannian measure on $ (X_k, g_k) $ satisfies the weak $ C_0 $-doubling condition at scale $ r_0 $ (for $ k $ large enough).\\
Furthermore, for any $ r_0 > 0$ and $ C_0 \ge  1 + 2^{2 n} \big(\cosh (K r_0) \big)^{2n-2} $, there exists a choice of the sequence $ (Y_k , h_k) $ such that 
all the $ (X_k , g_k) $ satisfy the weak $ C_0 $-doubling condition at scale $ r_0 $ and none of them satisfies the strong $ C_0 $-doubling condition at the same scale.
\end{lemma}

\begin{corollary}\label{exemple1} 
For any $ C_0 \ge 2^{2 n} + 1 $ and $ r_0 > 0 $, the weak  $ C_0 $-doubling condition at scale $ r_0 $ does not imply any restriction to the topology and the geometry of the balls of radius less than or equal to $ \frac{r_0}{20} $. Hence it does not imply any restriction to the local topology and geometry.
\end{corollary}

\begin{proof}[Proof of Lemma \ref{exemple}]
As the Ricci curvature of $ (X,g) $ is $\,\ge - (n-1)  K^2 \cdot g $, it follows from Lemma \ref{Riccidouble} that, for every $ r_0>0 $, the Riemannian measure of 
$ (X,g) $ satisfies the $ C'_0 $-doubling condition for all balls of radius $ r \in \left] 0 , 2r_0\right] $, where $ C'_0 = 2^n \big(\cosh (2 K r_0) \big)^{n-1} $.\\
Let us make the construction of  $ (X_k , g_k )= X \# Y_k $ more precise: we cut out balls $ B_k := B_X (x_0 , r_k) \subset X $ and $ B'_k := B_{Y_k} (y_k , r_k) \subset Y_k $, of radius $ r_k < \frac{1}{100} \Min \big( \inj (X,g) ; \inj (Y_k , h_k)\big) $. We then glue, as in Example \ref{exemplestrict} (3), $ X \setminus B_k $ and $ Y_k \setminus B'_k $ at the two ends of a cylinder $ \C_k := [0 , r_k^{3/4}] \times \mathbb S^{n-1} \left(r_k\right) $ of radius $ r_k $. Again, similarly to Example \ref{exemplestrict} (3), we construct a contracting map
 $ f_k : X_k \to X $ such that there exists a sequence $ \left(\eta_k \right)_{k \in \N} $, going to zero when $ k $ goes to $ +\infty $, such that $ d_X (f_k(x) , f_k(y) )  \ge d_{X_k}(x,y) -  \eta_k $.
Since $ f_k $ is contracting and onto, we have
\begin{equation}\label{volumapprox}
\Vol_{g_k} \big( B_{X_k} (x, r) \big) \ge \Vol_{g} \left(  f_k \big( B_{X_k} (x, r) \big) \right) \ge \Vol_{g} \left(  B_{X} \big( f_k(x)  , 
r - \eta_k  \big)\right)  .
\end{equation}
Now, from the definition of $ f_k $, there exists a sequence $ \left(V_k \right)_{k \in \N} $, going to zero when $ k $ goes to $ +\infty $, such that, for all domain 
$ A \subset X_k $, 
 $$\Vol_{g} \left(  f_k \left( A \right)\right) \ge \Vol_{g_k} \left( A \cap (X \setminus B_k) \right) \ge \Vol_{g_k} \left(A \right) - 
\Vol_{g_k} \left( Y_k \setminus B'_k\right) - \Vol_{g_k} \left( C_k \right)  \ge \Vol_{g_k} \left(A \right) - V_k  . $$
From this and the fact that $ f_k $ is contracting (which implies that $ f_k \left( B_{X_k} (x, 2 r) \right)  \subset
B_{X} \left( f_k (x) , 2r  \right) $) we deduce that
 $$\Vol_{g_k} \left( B_{X_k}(x, 2r) \right) \le \Vol_{g} \left(   f_k \left( B_{X_k}(x, 2r) \right)\right) + V_k
\le \Vol_{g} \left( B_X (f_k (x) , 2r) \right)+ V_k  . $$
As $ \frac{r}{2} $ and $ r $ belong to $ \left]0, 2r_0\right] $ for every $ r \in \left[ \frac{r_0}{2}, 2r_0\right] $,
this last inequality and \eqref{volumapprox} imply that there exists $ N_0 \in \N $ such that, for all 
 $ k \ge N_0 $, for all $ x\in X $ and for all $ r \in \left[ \frac{r_0}{2} , 2r_0\right] $, one has
 $$\dfrac{\Vol_{g_k}  B_{X_k}  (x , 2r )}{\Vol_{g_k} B_{X_k} (x , r )} \le \dfrac{\Vol_{g} \left( B_X (f_k (x) , 2r) \right) 
+ V_k}{\Vol_{g} \left( B_X \big( f_k (x) , r - \eta_k  \big) \right) } \le 
\dfrac{\Vol_{g} \left( B_X (f_k (x) , 2r) \right)}{\Vol_{g} \left( B_X \big( f_k (x) , r/2 \big) \right) }+ 1 \le \big(C'_0\big)^2 
+ 1 = C_0 \,. $$
Hence, for all $ k \ge N_0 $, the Riemannian measure of $ (X_k , g_k ) $ satisfies the weak $ C_0 $-doubling condition at scale $ r_0 $ around every point $ x \in X_k $. This shows the first part of the lemma.

Consider any sequence of $n$-manifolds $(Y_k)_{k \in \N^*}$ such that each $Y_k$ admits\footnote{Up to taking a subsequence, one can verify that the existence 
of such a sequence 
of metrics $  h'_k $ is equivalent to the existence of a sequence of metrics $  h''_k $ whose sectional curvature is $ \le - 1 $ and whose injectivity radius goes to $+\infty$.} a Riemannian metric $h'_k$ whose sectional curvature is $ \, \le - k^2 $ and whose injectivity radius is $ \, \ge 4 r_0$. Choosing a real sequence
$ \left(\varepsilon_k \right)_{k \in \N} $, going to zero when $ k $ goes to $ +\infty $ and verifying $\forall k,\ \e_k \le \diam (Y_k , h'_k)^{-2}$, we define the metric
$h_k$ on $Y_k$ by $ h_k := \varepsilon_k^2 \cdot h'_k $. Constructing $(X_k , g_k)$ as a connected sum of $ (X , g, x_0) $ and $ (Y_k , h_k , y_k) $ as above, we
choose a point $ x_k \in Y_k $ at $ h'_k $-distance of $ y_k $ at least $ 5 r_0 $. For all $ r \le 2r_0 $, we have $B_{(Y_k , h_k)} (x_k , 4 \e_k r_0 ) = B_{(Y_k , h'_k)} (x_k , 4 r_0) \subset Y \setminus B'_k$, hence $B_{X_k}  (x_k , \varepsilon_k r' ) = B_{(Y_k , h'_k)} (x_k , r')$ for every $r' \le 4 r_0$; this yields
$$\dfrac{\Vol_{g_k}  B_{X_k}  (x_k , 2\varepsilon_k r )}{\Vol_{g_k}  B_{X_k} (x_k ,\varepsilon_k r )} =
\dfrac{\Vol_{h'_k} B_{(Y_k, h'_k)}  (x_k , 2r )}{\Vol_{h'_k} B_{(Y_k , h'_k)} (x_k , r )} \ge \dfrac{\int_0^{2r} \sinh^{n-1} (kt) dt}{\int_0^{r} \sinh^{n-1} (kt) dt} 
\ge e^{(n-1) k r} . $$
The doubling amplitude (for balls whose radius is close to zero) then goes to $ + \infty $ with $ k $. Hence, for all $ C_0> 1 $ and $ r_0 > 0 $, and for all $ k $ large enough,
none of the manifolds $ (X_k, g_k) $ satisfy the strong $ C_0 $-doubling as scale $ r_0 $ around $ x_k $. 
\end{proof}
\begin{proof}[Proof of Corollary \ref{exemple1}]
Let $ (X,g) $ be a $ n $-dimensional flat torus. As in the proof of Lemma  \ref{exemple}, we construct $(X_k , g_k)$ as a connected sum of $ (X , g, x_0) $ and 
$ (Y_k , h_k , y_k) $, where the diameters and the volumes of the $ (Y_k , h_k )$'s go to zero when $ k \f +\infty $.
Since there is no restriction on the topology of  $ Y_k $ and on the geometry (modulo homotheties) of $ h_k $ and as the gluing is made within a ball of radius 
much smaller than $ \frac{r_0}{20} $, there is no restriction on the topology and on the geometry of the balls of  $ X_k $ of radius less than $ \frac{r_0}{20} $, 
centred around a point  $ x_k \in Y_k \setminus B'_k $. 
Nevertheless Lemma \ref{exemple} (where we make $K = 0$) shows that, for any $ r_0 > 0 $ and every $ C_0 \ge 2^{2 n} + 1 $, there exists $ N_0 \in \N $ 
such that (for every $ k \ge N_0 $) the Riemannian measure of $ ( X_k , g_k) $ satisfies the weak $ C_0 $-doubling condition at scale $ r_0 $. 
\end{proof}

\subsection{Doubling property induced on subgroups}

An important feature of the counting measure on a group, that we will extensively use in the sequel, is that the doubling condition transfers to subgroups. The next proposition makes this more precise.

\begin{prop}\label{sousgroupes}
Let $ \Gamma $ be a group acting properly by isometries on a metric space $ (X,d) $. Let $ \Gamma '\subset\Gamma $ be a subgroup and $ x_0 \in X $ a point. We then have,
\begin{itemize} 
    \item[(i)] If the counting measure $ \mu_{x_0}^{\Gamma} $ satisfies the $ C $-doubling condition for all balls of radius $ r \in \left[ \frac{1}{2} r_0  , \frac{5}{4}r_1 \right] $ centred 
around $ x_0 $ (where $ 0 < r_0 \le r_1 $) the measure $ \mu_{x_0}^{\Gamma'} $ satisfies:
 $$\forall r \in \left[ r_0 , r_1 \right]  ,\quad \dfrac{\mu_{x_0}^{\Gamma'}  \big( \overline {B}_X( x_0 , 2 r ) \big)}{\mu_{x_0}^{\Gamma'} \big( B_X (x_0 , r ) \big)} \le C^3 . $$
    \item[(ii)] If the counting measure $ \mu_{x_0}^{\Gamma} $ satisfies $ \dfrac{\mu_{x_0}^{\Gamma} 
\big( B_X (x_0 , 2r ) \big)}{\mu_{x_0}^{\Gamma} \big( B_X (x_0 , r ) \big)}
\le C    e^{\alpha  r} $ for all $ r \in [ \frac{1}{2} r_0 , +\infty [ $, then the measure $ \mu_{x_0}^{\Gamma'} $ satisfies
 $$\forall r \in \left[ r_0 , +\infty \right[  ,\quad \dfrac{\mu_{x_0}^{\Gamma'} \big(\overline {B}_X (x_0 , 2r ) \big)}{\mu_{x_0}^{\Gamma'} \big( B_X (x_0 , r ) \big)}
\le C^3  e^{\frac{19}{8}\alpha  r}  . $$
\end{itemize} 
\end{prop}
 
\begin{proof}
Let us recall that, for all $ r > 0 $, we have denoted by $ N_5^\Gamma (r) $ (resp. $ N_5^{\Gamma'} (r) $) the maximal number of disjoint balls of radius $ \frac{r}{2} $ contained in the ball $ B_X (x_0 , \frac{5}{2}r) $, and centred around points of the orbit $ \Gamma \cdot x_0 $ (resp. of the orbit $ \Gamma' \cdot x_0 $). 
Being in  $ \Gamma' \cdot x_0 $ is more restrictive than being in $ \Gamma \cdot x_0 $, hence we have $ N_5^{\Gamma'} (r) \le N_5^{\Gamma} (r) $.
Lemma \ref{compardoubling} applied twice, firstly to the counting measure $ \mu_{x_0}^{\Gamma'} $ on the orbit $ \Gamma' \cdot x_0 $, and secondly to the counting measure $ \mu_{x_0}^{\Gamma} $ on the orbit $ \Gamma \cdot x_0 $, yields
\begin{equation}\label{sousgroupe1}
  \begin{split}
\dfrac{\mu_{x_0}^{\Gamma'} \left(\overline B_X \big( x_0 ,  2r\big)\right)}{\mu_{x_0}^{\Gamma'} 
\big( B_X ( x_0 ,  r)\big)} & \le N_5^{\Gamma'} (r) \le  N_5^{\Gamma} (r) \le
\dfrac{\mu_{x_0}^{\Gamma} \left(B_X \big( x_0 , \frac{5}{2}   r  \big)\right)}{\mu_{x_0}^{\Gamma} 
\left(B_X \big( x_0 ,   \frac{r}{2} \big)\right)}\\
& \le  \dfrac{\mu_{x_0}^{\Gamma} \big( B_X( x_0 , \frac{5}{2} r ) \big)}{ \mu_{x_0}^{\Gamma} \big( B_X( x_0 , \frac{5}{4} r) \big)} \cdot   \dfrac{\mu_{x_0}^{\Gamma} \big( B_X( x_0 , \frac{5}{4} r ) \big)}{ \mu_{x_0}^{\Gamma} 
\big( B_X( x_0 , \frac{5}{8} r) \big)}  \cdot   \dfrac{\mu_{x_0}^{\Gamma} 
\big( B_X( x_0 , r ) \big)}{ \mu_{x_0}^{\Gamma} \big( B_X( x_0 , \frac{1}{2} r) \big)}\ .
  \end{split}
\end{equation}
$\bullet$ For all $ r \in\left[ r_0 ,  r_1 \right] $, since $ \frac{1}{2}r $, $ \frac{5}{8}r $ and $ \frac{5}{4}r $ belong to $ \left[ \frac{1}{2} r_0  , \frac{5}{4} r_1 \right] $, the doubling condition, assumed in $ (i) $, satisfied by $ \mu_{x_0}^{\Gamma} $ allows to bound from above by $ C $ each term in the last inequality of \eqref{sousgroupe1}, which proves that 
   $ \dfrac{ \mu_{x_0}^{\Gamma'} \left[ \overline B_X(x_0 , 2 r)\right]}{\mu_{x_0}^{\Gamma'} \left[ B_X (x_0 , r)\right]} \le C^3 $ and ends the proof of $ (i) $.

$\bullet$ For all $ r \in \left[ r_0 , +\infty \right[ $, the hypothesis assumed in $ (ii) $ allows to bound from above each term in the last inequality of \eqref{sousgroupe1} by, respectively, $ C e^{\frac{5}{4} \alpha  r} $, $ C e^{\frac{5}{8} \alpha  r} $ and  $ C e^{\frac{1}{2} \alpha  r} $. This proves $ (ii) $.
\end{proof}

\section{Free Subgroups}\label{sectionlibres} 

\small
\emph{In this  section, we only consider torsion-free isometries of a $ \delta $-hyperbolic space $ (X, d) $ (hence geodesic and proper by Definition \ref{hypdefinition}). 
To each isometry  $ \gamma $  of  $ (X, d) $, we associate its minimal displacement  $ s(\g) = \inf_{\g \in \Gamma^*} d(x, \g  x) $  and its asymptotic displacement  $ \ell(\g) = \lim_{n \to +\infty} \frac{d(x,\g^nx)}{n} $  (see Definitions \ref{deplacements}). If  $ \gamma $  is hyperbolic, $ {\cal G} (\g) $  is the set of oriented geodesics joining its fixed points  $ \gamma^- $  and  $ \gamma^+ $, and $ M(\g ) $  is the union of these geodesic lines (see Definitions \ref{faisceau}). For the basic properties concerning projections in hyperbolic spaces the reader is referred to Section \ref{sectionprojections}.}

\normalsize

The following remark will be useful in all this paper:
\begin{remark}\label{Dirichlet}
Let $ (X,d) $ be a $ \delta $-hyperbolic space and $ \g $  an isometry, then the three conditions \lq \lq $\g$ is non elliptic", \lq \lq $\g$ is torsion-free and the action 
of the group $ \langle \gamma \rangle $ generated by $\g$ is proper" and \lq \lq $\g$ is torsion-free and the action of the group $ \langle \gamma \rangle $ 
is discrete" are equivalent; in the sequel, we shall thus indifferently use one or the others of these three equivalent hypothesis.
\end{remark}

\begin{proof}
The fact that the assumptions \lq \lq proper" and \lq \lq discrete" (concerning the action of $ \langle \gamma \rangle $) are equivalent is due to Lemma \ref{discret1}.
If $\g$ is torsion-free and if the group $ \langle \gamma \rangle $ acts properly on $ (X,d)$, then $\g$ is non elliptic by Remark \ref{kpointsfixes} (i).\\
Conversely, if $\g$ is supposed non elliptic, then it is hyperbolic or parabolic by Theorem \ref{ellparahyp} and (by the definition of these notions given before
Theorem \ref{ellparahyp}) it is torsion-free and, for every $x \in X$, $ d(x , \g^k x)$ goes to $+\infty$ when $k \f \pm \infty$; hence the 
action of $ \langle \gamma \rangle $ on $ (X,d)$ is proper.
\end{proof}

\subsection{Ping-pong Lemma}\label{pingpong}
Let  $ (X,d) $  be a $ \delta $-hyperbolic space and $ \g $  a non elliptic isometry. For all  $ x \in X $  we denote by  $ D_\gamma (x) $  and  $ D'_\gamma (x) $ the 
\emph{Dirichlet domains (respectively open and closed) of the point  $ x $  for the action of}  $ \langle \gamma \rangle $, i. e. the subsets defined by\footnote{The minima in question in the sequel are achieved since the torsion-free and discreteness hypotheses ensure (by Remark \ref{kpointsfixes} (i)) that  $ \g $  is hyperbolic or parabolic, and hence that  $ g^k x $  goes to infinity when  $ k \f \pm \infty $.}:
\begin{equation}
D_\gamma (x) = \{ y :  \Min_{k \in  \mathbb Z^*} d(y , \gamma^k x ) > d(y,x) \} \ , \ 
D'_\gamma (x) = \{ y : \Min_{k \in  \mathbb Z}  d(y , \gamma^k x ) = d(y,x) \}.
\end{equation}

Let us remark that  $D'_{\g^{-1}} (x) = D'_{\g} (x)$, and that (for all  $ k \in \Z $)  $ D'_\gamma (\gamma^k   x) = 
\gamma^k \left(D'_\gamma (x) \right)$  and let us define the \emph{attraction (resp. repulsion)} domains
 $U^+_{\g} (x)$  (resp. $U^-_{\g} (x)$) of  $\g$  associated to $x$  by:
\begin{equation}\label{attraction}
U_\gamma^+ (x) := \cup_{k \in \mathbb N^*} \, D'_\gamma (\gamma^k   x) =
\cup_{k \in \mathbb N^*}\  \gamma^k \left(D'_\gamma (x) \right)\,,\quad  U_\gamma^- (x)  := U_{\g^{-1}}^+ (x) \ .
\end{equation}
Let us remark that, by Definition \eqref{attraction},
 $U^-_{\gamma} (x)  \cup  U^+_\gamma (x) = \cup_{p \in \Z^*}\  D'_\gamma (\gamma^p   x)$ 
and, by taking complements, one has:
$$X \setminus \left(U^-_{\gamma} (x)  \cup  U^+_\gamma (x)\right) =   \{ y \ : \ \forall p \in \Z^* \ \  
\Min_{k \in \Z} d(y , \g^k   x ) < d(y, \g^p x) \} $$
 $$= \{ y\  :\  \Min_{k \in \Z} d(y , \g^k   x ) < \Min_{p \in \Z^*}d(y, \g^p x) \} = \{ y\   : \ d(y ,  x ) < 
\Min_{p \in \Z^*}d(y, \g^p x) \} \,.$$
A consequence of this equality and of Definition \eqref{attraction} are the properties:
\begin{equation}\label{attraction1}
X \setminus \left(U^-_{\gamma} (x)  \cup  U^+_\gamma (x)\right) =  D_\gamma (x)\,, \quad \forall k \in \Z^* \ 
\g^k (D'_\gamma (x)) \subset U^-_{\gamma} (x)  \cup  U^+_\gamma (x)\ .
\end{equation}
Remark also that, as $\g$ is a parabolic or hyperbolic isometry by Theorem \ref{ellparahyp}, then $ \forall k \in \Z^*, \ \g^k x \ne x$ and $x \in D_\gamma (x) $.

\begin{defis}\label{Schottky} Let us consider two non elliptic isometries $a$ and $b$ of $(X,d)$ and a point $x \in X$,
\begin{itemize}
\item[(i)] $a$ and $b$ are said to be in Schottky position with respect to $x$ if the subsets $U^-_a (x)  \cup  U^+_a (x)$ and $U^-_b (x) \cup  U^+_b (x)$ are disjoints.
\item[(ii)] $a$ and $b$ are said to be in half-Schottky position with respect to $x \in X$ if  $U^+_a (x) \subset  X\setminus\big( U^-_b (x) \cup  U^+_b (x) \big)$ 
and $U^+_b (x) \subset  X \setminus \big( U^-_a (x) \cup  U^+_a (x) \big)$.
\end{itemize}
\end{defis}

\begin{prop}\label{Schottky1}
Let $\Sigma$ be a finite set of non elliptic isometries of the Gromov-hyperbolic space $(X,d)$. 
If every pair of elements of  $\Sigma$ is in Schottky position with respect to the same point $x \in X$, then $\Sigma$ generates a free subgroup of the group of isometries of $(X,d)$.
\end{prop}

\begin{proof}
Every non-trivial relation involving elements of $\Sigma$ is written as $s_1^{p_1}\cdot s_2^{p_2} \ldots s_{m}^{p_m} = e$, where $s_1,  s_2,  \ldots , s_{m}$ are elements of $\Sigma$ 
such that $s_{i} \ne s_{i+1}$ for all $i \in \{1 , \ldots , m-1\}$, and where  $p_1,  p_2,  \ldots , p_{m}$ are elements of $\Z^*$.

Let us prove by induction that, for all $i \in \{1 , \ldots , m\}$, $s_i^{p_i}\cdot s_{i+1}^{p_{i+1}} \ldots s_{m}^{p_m} (x) \in U^-_{s_i} (x) 
\cup  U^+_{s_i} (x)$. This is true for $i = m$, since $x \in D_{s_m} (x)$, and hence $s_{m}^{p_m} (x)\in  s_{m}^{p_m}\left(D_{s_m} (x)\right)\subset  U^-_{s_m} (x)\cup  U^+_{s_m} (x)$;  
this last inclusion is a consequence of \eqref{attraction1}.
We now assume that $s_i^{p_i}\cdot s_{i+1}^{p_{i+1}} \ldots s_{m}^{p_m} (x) \in U^-_{s_i} (x)\cup  U^+_{s_i} (x)$; since $s_i$ and $s_{i-1}$ are distinct elements of $\Sigma$, they are in Schottky position with respect to $x$, we then have $U^-_{s_i} (x) \cup  U^+_{s_i} (x) \subset X \setminus \left(U^-_{s_{i-1}} (x) \cup  U^+_{s_{i-1}} (x) \right) =  D_{s_{i-1}} (x) $, where this last equality follows from \eqref{attraction1}.
From this and from \eqref{attraction1}, we deduce that
$$s_{i-1}^{p_{i-1}}\cdot  s_i^{p_i} \ldots s_{m}^{p_m} (x) \in s_{i-1}^{p_{i-1}} \left(U^-_{s_i} (x) 
\cup  U^+_{s_i} (x)\right) \subset s_{i-1}^{p_{i-1}} \left(  D_{s_{i-1}} (x) \right) \subset U^-_{s_{i-1}} (x) \cup  U^+_{s_{i-1}} (x)\,,$$
which proves the induction. It then follows that $x = s_1^{p_1}\cdot s_2^{p_2} \ldots s_{m}^{p_m} (x)$ is an element of $U^-_{s_{1}} (x) \cup  U^+_{s_{1}} (x)$, which contradicts the fact that $x \in  D_{s_{1}} (x)$ (see above) since, by \eqref{attraction1}, $D_{s_{1}} (x) \cap \left( U^-_{s_{1}} (x) \cup  
U^+_{s_{1}} (x) \right) = \emptyset$. Consequently, there is no non-trivial relation between elements of $\Sigma$. 
\end{proof}

\begin{prop}\label{Schottky2} Let $ \{a , b\} $ be a pair of non elliptic isometries of $(X,d)$, if $a$ and $b$ 
are in half-Schottky position with respect to some point $x \in X$, then $\{a , b\}$ generates a free semi-group.
\end{prop}

Before proving Proposition \ref{Schottky2}, we establish the following lemma.

\begin{lemma}\label{pingpongball} Let $\{a , b\}$ be a pair of non elliptic isometries of $(X,d)$, if $a$ and $b$ 
are in half-Schottky position with respect to some point $x \in X$, then any product $R(a,b)$ of positive powers of $a$ and  $b$ satisfies 
$R(a,b)x \in U^+_{a} (x)$ (resp. $R(a,b)x\in U^+_{b} (x)$) if the first factor of the product
$R(a,b)$ is a power of $a$ (resp. of $b$).
\end{lemma}

\begin{proof}
It is sufficient to give the proof when the first factor of $R(a,b)$ is a power of $a$. Let us then assume that 
$$R(a,b) =  a^{p_0} b^{q_0}\,\ldots\,a^{p_k} b^{q_k}\quad \text{or}\quad R(a,b) =  a^{p_0} b^{q_0}\,\ldots\,a^{p_k} b^{q_k} a^{p_{k+1}}\,,$$
where all $ p_i $'s and $ q_i $'s belong to $\mathbb N^*$.

We note that $b^{q_k} (x)\in b^{q_k}\big(D_b (x)\big ) \subset U^+_{b} (x)$, where the inclusion follows for the first equality in \eqref{attraction}. Similarly we show that $a^{p_{k+1}} (x) \in U^+_{a} (x)$. As $a$ and $b$ are in half-Schottky position with respect to $x$, we also have
$$a^{p_i} \big( U^+_{b} (x) \big) \subset a^{p_i} \left( X \setminus \big( U^-_a (x) \cup  U^+_a (x) \big) \right) 
= a^{p_i} \big( D_a (x)\big) \subset U^+_{a} (x)$$
where the equality and  the last inclusion follows from the properties \eqref{attraction1} and \eqref{attraction}. Similarly we show that $b^{q_i} \big( U^+_{a} (x) \big) \subset  U^+_{b} (x)$. Iterating and alternating these two properties, we deduce that $a^{p_0} b^{q_0}\,\ldots\, a^{p_k} b^{q_k} (x)$ and $a^{p_0} b^{q_0}\, \ldots\,  a^{p_k} b^{q_k} a^{p_{k+1}} (x)$ belong to $ U^+_{a} (x)  $ and this ends the proof. 
\end{proof}

\begin{proof}[End of the proof of Proposition \ref{Schottky2}]
We show that, if the semi-group generated by $a$ and $b$ is not free, one of the two following alternatives is satisfied:
  \begin{itemize}
      \item[\bf Case 1]: There exists a relation involving $a$ and $b$ which has, after suitable simplifications, the form
      $R(a,b)=e= id_X$, where $R(a,b)$ is a product of strictly positive powers of $a$ and $b$.
      \item[\bf Case 2]: There exists a relation involving $a$ and $b$ which has, after suitable simplifications, the form
      $R_1(a,b)=R_2(a,b)$, where $R_1(a,b)$ and $R_2(a,b)$ are products of strictly positive powers of $a$ and $b$. We leave to the reader to check that we may furthermore assume that the first letter appearing in $R_1(a,b)$ is different from the first letter appearing in $R_2(a,b)$.
  \end{itemize}

Lemma \ref{pingpongball} implies that every product $R(a,b)$ of strictly positive powers of $a$ and $b$ satisfies $R(a,b)x \in  U^+_a (x) \cup  U^+_b (x)$. On the other hand the first equality of \eqref{attraction1} (and the fact that $x \in D_\gamma (x) $ when $\g$ is non elliptic) gives
$x \in D_a (x) \cap D_b (x) = X \setminus \left( U^-_a (x) \cup  U^+_a (x) \cup  U^-_b (x) \cup  U^+_b (x) \right)$, it is not possible to have $R(a,b)  x = x$, hence the relation $R(a,b) = e$. This takes care of Case 1.\\
If now $R_1(a,b)$ and $R_2(a,b)$ are non-trivial products of strictly positive powers of $a$ and $b$ with the first letter in $R_1(a,b)$ different from the first letter in $R_2(a,b)$, Lemma \ref{pingpongball} implies that 
either $R_1(a,b)x \in  U^+_a (x)$ and $R_2(a,b)x \in  U^+_b (x)$, if the first letter of $R_1(a,b)$ is $a$, or $R_1(a,b)x\in  U^+_b (x)$ and $R_2(a,b)x \in  U^+_a (x)$, if the first letter of $R_1(a,b)$ is $b$. This shows that we cannot have $R_1(a,b)x = R_2(a,b)x$. This takes care of Case 2 and ends the proof of Proposition \ref{Schottky2}.
\end{proof} 

The next proposition generalizes, in the context of isometric actions on a $\delta$-hyperbolic metric space, a result proved by Th.~Delzant for hyperbolic groups (see \cite{De}, Lemma 1.2, p. 179).

\begin{prop}\label{libre1}
Let $ (X , d)$ be a  $\delta$-hyperbolic space and $a$ and $b$ two isometries of $(X,d)$.
  \begin{itemize}
      \item[(i)] If there exists a point $x\in X$ such that $d(a^px, b^q x) > \Max \left[d(x, a^px), d(x, b^qx) \right] + 2\delta$ for all $(p,q)\in \mathbb Z^*\times 
\mathbb Z^*$, then $a$ and $b$ are in Schottky position with respect to $ x $ and the group generated by $a$ and $b$ is free.
      \item[(ii)]  If there exists a point $x\in X$ such that $d(a^p x, b^q x) > \Max \left[d(x, a^p x), d(x, b^q x) \right] + 2\delta$ for all $(p,q) \in\left(\mathbb Z^* \times\mathbb Z^*\right)\setminus\left(\mathbb Z^- \times \mathbb Z^-\right)$, then $a$ and $b$ are in half-Schottky position with respect to $x \in X$ and the semi-group generated by $a$ and $b$ is free.
  \end{itemize}
\end{prop}

\begin{remark}\label{implicites}
Any of the hypotheses made in points (i) or (ii) of Proposition \ref{libre1} automatically implies that $a$ and $b$ are non elliptic and that the point $x$ appearing in 
these two points is never fixed by any element of $\langle a \rangle\setminus \{\id_X\}$ and  by any element of $\langle b \rangle \setminus \{\id_X\}$.
\end{remark}

\begin{proof}[Proof of Remark \ref{implicites}] 
Any of the hypotheses made in \ref{libre1} (i) or (ii) implies that, for all $(p,q)\in \mathbb (Z^*\times \N^*)$ and all $(p,q)\in \mathbb (N^*\times \Z^*$, we have
$$\Max \left[d(x, a^p x), d(x, b^q  x) \right]  +  \Min \left[d(x, a^px), d(x, b^qx) \right] = d(x, a^px ) + d(x, b^qx )\ge d(a^px, b^q x) $$
$$> \Max \left[d(x, a^px), d(x, b^qx) \right] + 2\delta\,.$$
We deduce that $\Min \left[d(x, a^p x), d(x, b^q x) \right] > 2 \delta$ for all $(p,q)\in\mathbb N^*\times \mathbb N^*$, hence for all $(p,q)\in
\mathbb Z^*\times \mathbb Z^*$ and thus $a^px\ne x$ and $b^qx\ne x$: this proves the second part of the remark.
We also deduce that, for all $(k,p)\in \mathbb Z\times\mathbb Z$ such that $p\ne k$, we have $d(a^kx,  a^px) > 2 \delta$ and $d(b^k x,  b^p x ) > 2\delta$; this
shows that the sequence $\left( a^k x \right)_{k \in\N}$ does not admit any Cauchy subsequence, and is thus unbounded (by the properness of $(X,d)$); hence
$a$ is non elliptic. We similarly show that $b$ is non elliptic.
\end{proof} 

\begin{proof}[Proof of Proposition \ref{libre1}] 
Remark \ref{implicites} ensures that $a$  and $b$ are non elliptic, we can then use the definitions and results of Section \ref{pingpong}.\\
\emph{Proof of (i)}: Let $y\in U^-_b (x) \cup  U^+_b (x)$, Definitions \eqref{attraction} shows that there exists $k \in \mathbb Z^*$ such that $y \in D'_b (b^k   x)$, which implies that
$d(y, b^k x)\le d(y,x)$.  Let us choose this number $k$. For all $p \in \mathbb Z^*$, the Quadrangle Lemma \ref{proprietes} (ii), applied to the points $y$, $a^p  x$, $x$ and $b^k x$, yields
 $$d(y, x) + d(a^px, b^kx) \le \Max \left[d(y, b^kx ) + d(x , a^px)\, ;\, d( y, a^px ) + d(x, b^kx) \right]+2\delta \,.$$
The hypothesis allows then to deduce $d(y , x)  < \Max \left[d(y , b^k   x )\,;\, d( y , a^p   x ) \right]  $, hence that  $d(y , x) < d(y, a^p x) $ for all $p \in \mathbb Z^*$, which in turn implies
 $y\in D_a (x) =  X \setminus \big( U^-_a (x) \cup  U^+_a (x) \big) $ for all $y \in  U^-_b (x) \cup  U^+_b (x)$, hence that $U^-_b (x) \cup  U^+_b (x) \subset  X \setminus \big( U^-_a (x) \cup  U^+_a (x) \big)$. 
Definition \ref{Schottky} (i), then shows that $a$ and $b$ are in Schottky position with respect to $x$. Proposition \ref{Schottky1} then ensures that the subgroup of the isometry group of $(X,d)$ generated by $a$ and $b$ is free.\\
\emph{Proof of (ii)}: The proof that $U^+_b (x) \subset  X\setminus\big( U^-_a (x)\cup  U^+_a (x) \big)$ is identical to the proof of (i), except that we need to choose a random point $y\in U^+_b (x)$, and consequently assume that $k \in \mathbb N^*$. The proof that $U^+_a (x)\subset  X\setminus\big( U^-_b (x)\cup  U^+_b (x) \big)$ is done along the same lines, exchanging the roles of $a$ and $b$. 
From Definition \ref{Schottky} (ii),it follows that $a$ and $b$ are in half-Schottky position with respect to $x$. Proposition \ref{Schottky2} then shows that  the semi-group generated by $a$ and $b$ is free. 
\end{proof} 

\begin{corollary}\label{librebis}
Let $ (X , d)$ be a  $\delta$- hyperbolic space and $\Gamma$ a group acting on $ (X , d)$, for every pair $a,\, b$ of elements of $\Gamma$, if there 
exists some point $x\in X$ such that $d(a^p x, b^q x) > \Max \left[d(x, a^p x), d(x, b^q x) \right] + 2\delta$ for all $(p,q) \in\left(\mathbb Z^* \times\mathbb Z^*\right)\setminus\left(\mathbb Z^- \times \mathbb Z^-\right)$, then the semi-group generated by $a$ and $b$ is free.
\end{corollary}

\begin{proof} Let $\varrho$ be the representation from $\Gamma$ to $\text{Isom} (X,d)$ associated to the action under consideration; as (by definition of 
$\varrho$) $\g x := \varrho(\g) x$ for every $\g \in \Gamma$, the above hypothesis may be rewritten as $d \big(\varrho(a)^p x , \varrho(b)^q x \big) > 
\Max \left[d(x, \varrho(a)^p x), d(x, \varrho(b)^q x) \right] + 2\delta$ for all 
$(p,q) \in\left(\mathbb Z^* \times\mathbb Z^*\right)\setminus\left(\mathbb Z^- \times \mathbb Z^-\right)$. Applying Proposition \ref{libre1}, 
$\varrho(a)$ and $\varrho(b)$ generate
a free semi-group. If $a$ and $b$ do not generate a free semi-group then any non trivial relation between positive powers of $a$ and $b$ maps to a similar
non trivial relation between positive powers of $\varrho(a)$ and $\varrho(b)$, in contradiction with the fact that $\varrho(a)$ and $\varrho(b)$ generate a free 
semi-group. Hence $a$ and $b$ generate a free semi-group.
\end{proof} 

\subsection{When the asymptotic displacement is bounded below}

\begin{prop}\label{granddeplacementx} 
Let $(X , d)$ be a  $\delta$-hyperbolic space and $a$ and $b$ two isometries 
such that the group generated by $a$ and $b$ is a discrete non virtually cyclic subgroup of the isometry group, then
  \begin{itemize} 
      \item[(i)] if $s(a),  s(b) > 13\delta$,  one of the two semi-groups generated by $\{ a  ,  b\}$ or by $\{ a  ,  b^{-1}\}$ is free,
      \item[(ii)] if $\ell (a),  \ell (b) > 0$ then, for every integers $p,  q >  \dfrac{13 \delta}{\Min (\ell(a) , \ell (b))}$, one of the two semi-groups generated by 
$\{ a^p  ,  b^q\}$ or by 
      $\{ a^p  ,  b^{-q}\}$ is free.
  \end{itemize} 
\end{prop}

\begin{corollary}\label{granddeplacement} 
For any proper action (by isometries) of a group $\Gamma$ on a $\delta$-hyperbolic space $(X , d)$, for every hyperbolic elements $ a$ and $b $ of $\Gamma$
which generate a non virtually cyclic subgroup of $\Gamma$, for every integers $p,  q > \dfrac{13 \delta}{\Min (\ell(a) , \ell (b))}$,
one of the two semi-groups generated by $\{ a^p  ,  b^q\}$ or by $\{ a^p  ,  b^{-q}\}$ is free.
\end{corollary}

Before proving Proposition \ref{granddeplacementx} and Corollary \ref{granddeplacement}, we state and prove the following lemma.

\begin{lemma}\label{quasigeodesic}
Let $ \g $ be an hyperbolic isometry acting on a $\delta$-hyperbolic space $(X, d)$ and verifying $s (\g) > \frac{11}{2}\delta$, let $c \in {\cal G} (\g)$ be any oriented geodesic from $\gamma^- $ to $\gamma^+ $, and let $x\in M(\g) $, we then have
\begin{itemize}
\item[(i)] $ d(x, \g x)\le\ell (\g) + 4\delta$,
\item[(ii)] if $c(t_k)$ denotes a projection of $\gamma^k x $ on $c$, then the sequence $\left( t_k\right)_{k \in \Z}$ is strictly increasing.
\end{itemize}
\end{lemma}

\begin{proof}[Proof of (ii)] By contradiction, assume that the sequence $\left( t_k\right)_{k \in \Z}$ is not strictly increasing. As $t_k \f \pm\infty$ when $ k \f \pm\infty $, there exists $ p \in \Z $ such that 
$t_{p} \le\Min\big( t_{p-1} , t_{p+1} \big)$. Let us set $\delta_k = d\big( \g^k x , c(t_k) \big)$, Proposition \ref{geodasympt} (i) implies that $\delta_k \le 2 \delta$ for all $k \in \Z $ (since $c$ and $\g^k \circ c$ are two geodesics of $\cal{G} (\g)$). 
We then have the following.
  \begin{itemize}
      \item If $ t_{p}\le t_{p-1}\le t_{p+1}$, then the triangle inequality, Lemma \ref{quasigeod} (ii), and the hypothesis made on $s(\g)$ imply that:
 $$t_{p+1}-t_{p}  + 2 \delta_{p-1} + \delta_{p+1} +\delta_{p} \ge d(\g^{p} x, \g^{p-1} x) + d ( \g^{p-1} x, \g^{p+1} x)$$
        \begin{equation}\label{quasigeodesic1}
\ge  d (x, \g   x) + d (x, \g^{2} x) \ge 2 d(x, \g x) + s(\g) - \delta\ge 3 \,s (\g) - \delta > \frac{25}{2}  \delta\,.
        \end{equation}
As $\delta_k \le 2 \delta$ for all $k\in \Z$, we deduce from this inequality that $d \big( c(t_{p}) , c(t_{p + 1}) \big) > 3 \delta$ and we can apply Lemma \ref{ecartement} to get:
 $$d (x, \g x) + 6\delta = d (\g^{p} x, \g^{p+1} x)  + 6\delta \ge  t_{p+1} - t_{p} +\delta_{p+1} + \delta_{p} \ge 2d(x, \g x) + s (\g) - \delta -  2 \delta_{p-1} \,,$$
where the last inequality comes from \eqref{quasigeodesic1}; we deduce that 
 $2s(\g) \le d (x, \g x) +  s(\g) \le 7\delta +  2 \delta_{p-1} \le 11\delta $, which contradicts the hypothesis on $s(\g)$.
      \item If $t_{p}\le t_{p+1}\le t_{p-1}$, exchanging the roles of $p-1$ and $p + 1$ and using arguments similar to the previous case we get: 
      $2s(\g) \le d(x, \g x) +  s(\g) \le 7\delta +  2 \delta_{p+1}\le 11\delta$, in contradiction to the hypothesis on $s(\g)$.
\end{itemize}
This shows that the sequence $\left( t_k\right)_{k \in \Z}$ is strictly increasing.
\end{proof}

\begin{proof}[Proof of (i)] As we saw before, for all $k\in\N^*$, $\delta_0,\delta_k \le 2\delta$, the triangle inequality gives:
$$d(x, \g^k x) - 4\delta\le d(x, \g^k x) - \delta_0 - \delta_k \le |t_{k}  - t_{0}| \le d(x, \g^k x) + \delta_0 + \delta_k \le d(x, \g^k x) + 4\delta\,.$$
From (ii) we know that the sequence $\left( t_i\right)_{i \in \N}$ is increasing, we deduce that
$$\lim_{k \f +\infty} \left(\frac{1}{k}\sum_{i = 0}^{k - 1} |t_{i+1}  - t_{i}|\big) \right) = \lim_{k \f +\infty} \left( \frac{1}{k} |t_{k}  - t_{0}| \big) \right) = \lim_{k \f +\infty} \left(  \frac{1}{k}   d \big( x , \g^k x \big)\right) = \ell(\g)\,,$$
which implies that $\inf_{i \in \N} \  |t_{i+1}  - t_{i}| \le  \ell(\g)$. We deduce that
 $$d(x, \g x) = \inf_{i \in \N} d (\g^i x, \g^{i+1} x) \le \inf_{i \in \N} \left(|t_{i+1}  - t_{i}| + \delta_i + \delta_{i+1} \right)\le \ell (\g) + 4\delta\,.$$
\end{proof}
\begin{proof}[Proof of Proposition \ref{granddeplacementx}]
Notice that, in the assertion (i) as in the assertion (ii) of Proposition \ref{granddeplacementx}, the hypotheses imply that
$a$ and $b$ are hyperbolic (by Lemma \ref{ellpositive}) because $\ell (a),  \ell (b) > 0$ in both cases (in the case of assertion (i), this is deduced
from the inequality $\ell (\g) \ge s(\g) - \delta > 12 \delta$ proved in Lemma \ref{quasigeod} (i)).
To simplify, let us denote by $N_1$ the smallest integer strictly greater than $\dfrac{13 \delta}{\Min (\ell(a) , \ell (b))}$.
We first show that $(i) \implies (ii)$. Indeed, for every $p, q \ge N_1$, the group $\langle a^p, b^q\rangle$ is not virtually cyclic (by the proposition \ref{actionelementaire} (vi)), and 
$s(a^p) \ge \ell (a^p) \ge N_1  \ell (a) > 13 \delta $,
while $s(b^q) \ge \ell (b^q) \ge N_1  \ell (b) > 13 \delta $. Property (i) then implies that one of the two semi-group 
generated by $\{ a^p  ,  b^q\}$ or by $\{ a^p  ,  b^{-q}\}$ is free.

{\it Proof of (i) :}  Under the hypotheses of (i), the lemma \ref{quasigeod} (i) proves that 
$\ell(a)\ge s(a) - \delta > 12\delta$ and that $\ell(b) \ge s(b) - \delta > 12\delta $; a consequence is that 
  \begin{equation}\label{sminore}
\forall k\in \Z^*\,,\quad s(a^k) > \max \big(13,  12|k|\big)\cdot\delta\quad\text{and}\quad s(b^k)> \max\big(13, 12  |k| \big)\cdot \delta\,;
  \end{equation}
indeed, this property is satisfied (by hypothesis) when $k = \pm 1$. When $|k| \ge 2$, we have $s(a^k) \ge \ell(a^k) = |k|\ell(a) >  |k|\cdot 12\delta\ge 24\delta$ and we show in the same way that
$s(b^k) > |k|\cdot 12\delta \ge 24\delta $.

$a^-$ and $a^+$ (resp. $b^-$ and $b^+$) being the fixed points of $a$ (resp. of $b$), we denote by $c_a$ (resp. $c_b$) one of the oriented geodesics from 
$a^-$ to $a^+$ (resp. from $b^-$ to $b^+$). If the limit sets $\{ a^+  , a^- \}$ and $\{ b^+  , b^- \}$ have a common point, they are identical by Proposition 
\ref{actionelementaire} (i) and are the fixed points of all the elements of $\langle a , b \rangle$, which in turn implies, from Proposition \ref{actionelementaire} (ii), that $\langle a , b \rangle$ is virtually cyclic which is excluded.
Hence $\{ a^+  , a^- \}$ and $\{ b^+  , b^- \}$ have no common point.\\
Proposition \ref{projgeod} (ii) then shows that there exist points $x_0 = c_a (s_0)$  and $x_0' = c_a (s_0')$ on the geodesic line $c_a$ such that, for every sequence 
$\left(t_n\right)_{n \in \mathbb N}$ going to $+\infty$, there exists a sequence 
$\left(\varepsilon_n\right)_{n \in \mathbb N}$ of strictly positive real numbers, going to zero when $n \to + \infty$, such that, for all $t \in \mathbb R$,
 \begin{equation}\label{granddeplace1}
d\left(c_a (t), c_b (t_n)\right) \ge  d\left(c_b (t_n), x_0\right) +  d\left(x_0, c_a(t)\right)-5 \delta - \varepsilon_n\, .
  \end{equation}
\begin{equation}\label{granddeplace2}
d\left(c_a (t), c_b (-t_n)\right) \ge  d\left(c_b (-t_n), x_0'\right) +  d\left(x_0', c_a(t)\right) - 5 \delta -\varepsilon_n \,.
  \end{equation}
{\it By eventually changing $b$ in $ b^{-1} $ (which changes the orientation of $c_b$) hence exchanging $x_0$ and $x_0'$, we may assume in the sequel that $s_0' \le s_0$.}

Let us note $y_0 = c_b (r_0)$ a projection of $ x_0  $ on the geodesic $c_b$. Denote by $c_a (s_k)$ a projection of $a^k x_0$ on $c_a$ and by
 $c_b (r_q)$ a projection of $b^q y_0$ on $c_b$. Proposition \ref{geodasympt} (i), applied to the geodesics $c_a$ and $a^k \circ c_a$ (resp. to the geodesics $c_b$ and $b^q 
\circ c_b)$, gives:
  \begin{equation}\label{granddeplace3}
\forall k, q \in \Z\,,\quad  d\big(a^k x_0,  c_a (s_k)\big) \le 2\delta\quad\text{ and resp.}\quad d \big(b^q y_0,  c_b (r_q)\big)\le 2  \delta\,.
  \end{equation}
For all $q \in \N^*$ and all $k \in \Z^*$, we consider a sequence $\left(t_n\right)_{n \in \mathbb N}$ which goes to $+\infty$. 
By hypothesis $s (b) > 13 \delta$, we can then apply Lemma \ref{quasigeodesic}, which shows that the sequence
 $n \to r_n $ is strictly increasing. In particular, we have $r_0 < r_q$ and hence, when $n$ is large enough,
 $d \big( y_0, c_b (t_n) \big) = d \big( y_0, c_b (r_q) \big) +  d \big(c_b (r_q), c_b (t_n) \big)$.
Applying Lemma \ref{projection} and the triangle inequality, gives, when $n$ is large enough,
$$d \big( x_0, c_b (t_n) \big) + 2 \delta \ge  d \big( x_0, y_0 \big)+ d \big( y_0, c_b (t_n) \big) = d \big( b^q x_0,  b^q y_0 \big)+ d \big( y_0, c_b (r_q) \big) + d\big(c_b (r_q), c_b (t_n) \big)$$
 $$\ge d \big( b^q x_0, c_b (r_q) \big) -  2  d\big( c_b (r_q), b^q y_0 \big) + d \big( y_0,  b^q y_0\big) + d\big(c_b (r_q), c_b (t_n) \big)\,.$$
From the second inequality \eqref{granddeplace3} we deduce:
  \begin{equation}\label{granddeplace4}
d\big( x_0, c_b (t_n) \big)\ge  d\big( b^q x_0, c_b (t_n) \big) -  6\delta + d\big( y_0, b^q y_0\big)\,. 
  \end{equation}
Applying inequality \eqref{granddeplace1} and the quadrangle inequality (Lemma \ref{proprietes} (ii)), we get, when $\varepsilon_n$ is small enough,
 $$d \big( x_0, b^q x_0 \big) + d \big( x_0, c_a (s_k) \big) + d \big( x_0, c_b (t_n) \big) - 7\delta - \varepsilon_n \le d\big( x_0, b^q x_0 \big) + d \big( c_a (s_k),  c_b (t_n) \big) - 2\delta $$
 $$\le \Max \left[ d \big( x_0, c_a (s_k) \big) + d \big(  b^q  x_0, c_b (t_n)  \big)\, ;\, d \big(  b^q  x_0, c_a (s_k) \big) + d \big( x_0,  c_b (t_n)  \big)\right]$$
\begin{equation}\label{granddeplace5}
\le d \big(  b^q  x_0, c_a (s_k) \big) + d \big( x_0,  c_b (t_n)  \big)\,;
\end{equation}
where the last inequality follows from Inequality \eqref{granddeplace4} and from the fact that $d\big( x_0, b^q x_0 \big) + d \big( y_0, b^q y_0 \big)\ge 2s(b^q) 
> 26 \delta $ by Property \eqref{sminore}. 
When $ n $ is large enough, for all $q\in\N^*$ and all $k \in \Z^*$, applying the inequality $\Min\left[ d\big( x_0, b^q x_0 \big)\,,\, d\big( x_0, a^k x_0 \big) \right] 
\ge \Min (s(a^k)\,,\, s(b^p) )> 13 \delta + \e_n$ (which follows from Property \eqref{sminore} and from the fact that $n$ is large), applying Inequalities \eqref{granddeplace5}, \eqref{granddeplace3} and finally the triangle inequality, we deduce that 
$$\Max \left[ d\big( x_0, b^q x_0 \big)\,,\, d \big( x_0, a^k x_0 \big)  \right] + 2\delta < d \big( x_0, b^q x_0 \big) + d \big( x_0, a^k x_0 \big) - 11\delta - \varepsilon_n$$
\begin{equation}\label{granddeplace5b}
\le d \big( x_0, b^q x_0 \big) + d \big( x_0, c_a (s_k) \big) - 9\delta - \varepsilon_n \le d \big(  b^q  x_0, c_a (s_k) \big)- 2\delta \le d \big( b^q  x_0, a^k x_0 \big)\,.
\end{equation}

Let us now consider the case where $k \in \mathbb N^*$ and $q  \in - \mathbb N^*$ and where $\left(t_n\right)_{n \in \mathbb N} $ is
a real valued sequence which goes to $ +\infty $. Lemma \ref{quasigeodesic} and
the hypothesis on $ s (a) $ ensure that the sequence $ k \mapsto s_k $ is strictly increasing. We have  $ s_0' \le s_0 < s_k $ and then 
$d \big(x_0', c_a (s_k)  \big) = d \big( x_0', x_0 \big) +  d \big(x_0, c_a (s_k)  \big) $.
Inequality (\ref{granddeplace2}) and the triangle inequality yields, when $ n $ is large enough,
 $$d\left( c_a (s_k),  c_b (-t_n) \right) + 5 \delta + \varepsilon_n  \ge  d\left( c_b (-t_n) ,  x_0'\right) +  d\left(x_0' ,   c_a (s_k) \right) $$
 $$\ge d\left( c_b (-t_n) ,  x_0'\right) +  d\left(x_0',  x_0) + d( x_0 ,   c_a (s_k) \right) \ge  d\left(  c_b (-t_n) ,  x_0\right) +  d\left( x_0 ,   c_a (s_k) \right)\,.$$
From this and from the Quadrangle Lemma \ref{proprietes} (ii), when $ n \in \mathbb N^* $ is large enough, we get
 $$d ( x_0, b^{q}  x_0 ) + d\left(  c_b (-t_n),  x_0\right) +  d\left( x_0,   c_a (s_k) \right) - 7\delta - \varepsilon_n \le d ( x_0, b^{q}  x_0 ) + d\left(  c_a (s_k),  c_b (-t_n) \right) - 2 \delta$$
\begin{equation}\label{granddeplace6}
\le \Max \left[ d ( x_0,   c_a (s_k)  ) + d (  b^{q}  x_0,  c_b (-t_n) ); d ( b^{q}  x_0,  c_a (s_k  )  ) + d ( x_0,   c_b (-t_n) \right]\,.
\end{equation}
As Lemma \ref{quasigeodesic} and the hypothesis on $s (b)$ ensures that $n\mapsto r_n$ is strictly increasing, we have $-t_n < r_q < r_0$ (when 
$ n $ is large), which implies that
$d \big( y_0, c_b (- t_n) \big) = d \big( y_0, c_b (r_q) \big) +  d \big(c_b (r_q), c_b (-t_n) \big)$.
With this equality, Lemma \ref{projection} and the triangle inequality, we get
$$d ( x_0, c_b (-t_n) ) + 2\delta \ge d ( x_0, y_0 ) + d( y_0, c_b (-t_n) ) = d ( x_0, y_0 ) + d( y_0,  c_b (r_{q} ) ) + d( c_b (r_{q} ),  c_b (-t_n)  )$$
$$\ge d ( x_0, y_0 ) + d( y_0, b^q y_0  ) - 2\delta + d( c_b (r_{q} ),  c_b (-t_n)  )\,,$$
where the last inequality follows from \eqref{granddeplace3}. From this last inequality, using the triangle inequality, the fact that $ d( y_0, b^{q} y_0) > 13\delta $
and \eqref{granddeplace3}, we deduce
$$d (  b^{q}  x_0, c_b (-t_n) ) < d (  b^{q}  x_0,  b^{q}  y_0 ) + d ( b^{q}  y_0, c_b (r_q) ) + d( c_b (r_{q} ),  c_b (-t_n)  ) +  d( y_0, b^q y_0)  - 13 \delta$$
$$\le d ( x_0,  y_0 ) + d( c_b (r_{q} ),  c_b (-t_n)  )  +  d( y_0, b^q y_0) - 11  \delta\le d ( x_0, c_b (-t_n) ) - 7 \delta\,,$$
which implies 
$$d ( x_0, c_a (s_k)  ) + d (  b^{q}  x_0, c_b (-t_n) )< d ( x_0, c_a (s_k)  ) + d ( x_0, c_b (-t_n) ) - 7\delta $$
$$< d ( x_0,  b^{q}  x_0 ) + d\left(  c_b (-t_n),  x_0\right) +  d\left( x_0,   c_a (s_k) \right) - 7\delta - \varepsilon_n \,.$$
Plugging this last inequality in \eqref{granddeplace6}, we deduce that
 $$d (x_0, b^{q}x_0) +  d\left( x_0,   c_a (s_k) \right) -  7\delta - \varepsilon_n \le d ( b^{q}  x_0 ,  c_a (s_k)    ) $$
which yields, using the triangle inequality and estimates \eqref{granddeplace3}, 
$$d ( x_0, b^{q}  x_0 ) +  d\left( x_0,    a^k  x_0  \right) -  9  \delta - \varepsilon_n \le d ( b^{q}  x_0,  a^k  x_0 ) + 2\delta\,.$$
Now, applying the inequality $\Min\left[ d \big( x_0, b^q x_0 \big), d \big( x_0, a^k x_0 \big) \right] \ge \Min (s(a^k), s(b^p) )> 13   \delta + \e_n $
(which follows from property \eqref{sminore} and is valid when $ \varepsilon_n$ is small enough), we deduce that
$$\forall k, p \in \N^*\,,\quad\Max  \left[d ( x_0, b^{q}  x_0 ),   d( x_0 ,  a^k  x_0  \right] + 2\delta < d ( b^{q}  x_0,  a^k  x_0 )\,.$$
This last inequality and \eqref{granddeplace5b} show that, for all $(p,q) \in\left(\mathbb Z^* \times \mathbb Z^*\right)\setminus  \left(\mathbb Z^- \times \mathbb Z^-\right)$, 
we have 
\begin{equation}\label{granddeplt}
d ( b^{q}  x_0,  a^k  x_0 )  > \Max  \left[ d ( x_0, b^{q}  x_0 ),   d( x_0, a^k  x_0  \right] + 2 \delta\,.
\end{equation}
Proposition \ref{libre1} (ii) then shows that the semi-group generated by $a$ and $b$ is free. This ends the proof of  (i) and hence of Proposition \ref{granddeplacementx}. 
\end{proof}

\begin{proof}[Proof of Corollary \ref{granddeplacement}]
If $\varrho$ is the representation $\Gamma \f \text{Isom} (X,d)$ associated to the action of $\Gamma$ on $ (X,d)$, Lemma \ref{reductisom} (v) and (vi)
proves that $\varrho (a)$ and $\varrho (b)$ are hyperbolic isometries of $(X,d)$ satisfying $\ell (\varrho (a)) = \ell (a) > 0 $ and $\ell (\varrho (b)) = 
\ell (b) > 0$. As $\langle a, b\rangle$ is non virtually cyclic and as the action is proper, Lemma \ref{reductisom} (ii) and (vii) guarantees that 
$\varrho (a)$ and $\varrho (b)$ generate a non virtually cyclic discrete subgroup of isometries of $ (X,d)$.
For every $ p,  q > \dfrac{13 \delta}{\Min (\ell(a) , \ell (b))} = \dfrac{13 \delta}{\Min \big(\ell (\varrho (a)) , \ell (\varrho (b))\big)} $, Proposition 
\ref{granddeplacementx} (ii) (applied to $\varrho (a)$ and $\varrho (b)$) implies that $\varrho (a^p)$ and $\varrho (b^q)$ (or $\varrho (a^p)$ and 
$\varrho (b^{-q})$) generate a free semi-group; eventually changing $b$ in $ b^{-1} $, suppose that $\varrho (a^p)$ and $\varrho (b^q)$ generate a free 
semi-group.
If $a^p$ and $b^q$ do not generate a free semi-group then any non trivial relation between positive powers of $a^p$ and $b^q$ maps to a similar
non trivial relation between positive powers of $\varrho(a^p)$ and $\varrho(b^q)$, in contradiction with the fact that $\varrho(a^p)$ and $\varrho(b^q)$ generate 
a free semi-group. Hence $a^p$ and $b^q$ generate a free semi-group.
\end{proof}

\subsection{When some Margulis constant is bounded below:}\label{free1}

\small

For every group $G$ and every $A, B , S \subset G$, we denote by $A\cdot B$ the image of $A\times B$ by the map $(\g, g) \mapsto \g \cdot g$ and define (by 
induction) $S^k$ as $S^{k-1}\cdot S$.\\
The aim of this subsection is to prove that, for any finite set $S$ of isometries of any $\delta$-hyperbolic space $(X , d) $,which generates a non virtually cyclic
group $\langle S \rangle$ of isometries, if the Margulis constant of $S$ (see Definitions \ref{cteMargulis}) is $\ge C\, \delta$, then $S^{16}$ contains two elements 
which generate a free semi-group. In the case where $\# S = 2$, modifying a little the definition of the Margulis constant, we get a sufficient condition for the group
$\langle S \rangle$ to be free. These two \lq \lq Margulis constants" are defined as follows:

\normalsize

\begin{defis}\label{cteMargulis}
For any metric space $(X,d)$, denote by $\text{\rm Isom} (X,d)$ the group of the isometries of $(X,d)$ and define
\begin{itemize}
\item[(i)] the Margulis constant $L(a,b)$ of a pair $\{a , b\} \subset \text{\rm Isom} (X,d)$ as the infimum of the function 
$(x, p ,q) \mapsto \Max \, [\, d(x, a^p x)\, ; \, d(x , b^q x)\, ]$ for all the $(x, p ,q)  \in X \times \Z^* \times \Z^*$,
\item[(ii)] the  Margulis constant $L^*(S)$ of a finite set $ S \subset \text{\rm Isom} (X,d)$ as the infimum of the function 
$ x \mapsto \Max_{\g \in S}  d(x, \g x)$ for all the $x \in X$.
\end{itemize}
\end{defis}

Remark that, as $ L(a,b) \le L^*(\{a , b\}) $, assuming a lower bound for $L^*(\{a , b\})$  is a weaker hypothesis than assuming a lower bound for $L(a,b)$.

\subsubsection{When the  Margulis constant  $L^*$ is bounded below:}\label{freestar}

{\bf a) When $L^*$ is applied to a pair of isometries}

\begin{theorem}\label{minorLstar}
On any $ \delta$-hyperbolic space $(X , d) $, for every pair $\{a , b\} \subset \text{\rm Isom} (X,d)$ such that  $L^*(\{a , b\}) \ge \frac{23}{2} \,\delta$ and
which generates a discrete and non virtually nilpotent subgroup, one has
\begin{itemize}
\item[(i)] either $\Max [ \ell(a) , \ell(b)] \ge \delta$ and one of the semi-groups generated by $\{a^{14}, ba^{14} b^{-1}\}$
or by $\{a^{14}, ba^{-14} b^{-1}\}$ is free,
\item[(ii)] or all the elements $\g \in \{ ab , a b^{-1},  ba , b a^{-1}, a^{-1} b,  a^{-1} b^{-1} , b^{-1} a , b^{-1} a^{-1} \}$
verifies $\ell (\g) \ge 2\, \delta$ and one of the semi-groups generated by $\{(ab)^{7}, (ba)^{7}\}$ or by $\{(ab)^{7}, (ba)^{-7}\}$ is free,
\end{itemize}
\end{theorem}

Before proving this Theorem, let us first establish the two following preliminary Lemmas:
\begin{lemma}\label{minorLstar1}
Let $(X , d) $ be a $ \delta$-hyperbolic space, for any pair $\{a , b\} \subset \text{\rm Isom} (X,d)$ such that $L^*(\{a , b\}) > 3\, \delta + 
\Max \big( \ell(a) , \ell(b)\big)$, for every $x \in X$, the middle points $m_a$ (resp. $m_b$) of any geodesic segment from $x$ to $a x$ (resp. from $x$ to $b x$) satisfy
$ d(m_a , m_b) \ge L^*(\{a , b\}) -\left( \frac{1}{2} \big(\ell (a) + \ell (b)\big)+ 3\,\delta \right)$
\end{lemma}

\begin{proof} For sake of simplicity, we set $ L^* :=  L^*(\{a , b\})$ and $D := d(m_a , m_b) $; consider any geodesic $ c : [0 , D ] \to X $ from $m_a$ to $m_b$. 
By Lemma \ref{milieux}, we have
$$ d(m_a , a \,m_a)  \le 3\, \delta +  \ell(a) \ \ \ \ ;  \ \ \ \ d(m_b , b \, m_b)  \le 3\, \delta +  \ell(b)$$
From this and from the triangle inequality, it follows that $ d \big(c(t) , a\circ c(t)\big) \le 2 \, t + 3\, \delta +  \ell(a)$ and $d \big(c(t) , b\circ c(t)\big) \le 
2 (D-t) + 3\, \delta +  \ell(b)$; by the definition of $ L^*(\{a , b\})$, this gives, $\forall t \in [0 , D]$:
\begin{equation}\label{minorLstar2}
L^* \le \Max \big[ \,d \big(c(t) , a\circ c(t)\big)\,,\,  d \big(c(t) , b\circ c(t)\big) \,\big] \le \Max \big[ 2 \, t +  \ell(a)\,,\, 2 (D-t)  +  \ell(b) \,\big]+ 3\, \delta \ .
\end{equation}
Applying this inequality for $t= 0$ and $t=D$ (and the hypothesis $\ell(a) , \ell(b) < L^* - 3\delta$) yields 
$$ \ell (b) <  L^* -  3\, \delta \le 2 \, D  +  \ell(a)\ \ \ \ \text{and}  \ \ \ \  \ell (a) <  L^* -  3\, \delta \le 2 \, D  +  
\ell(b) \ ;$$
these inequalities imply that $|\ell(b) - \ell(a)| < 2\, D$, thus that $t_0 := \frac{1}{4} \big( 2\, D + \ell (b) - 
\ell (a)\big)$ verifies $ 0 < t_0 < D$; we may then apply the inequality \eqref{minorLstar2} to $t=t_0$ and obtain:
$$ L^* -  3\, \delta \le \Max \big[ 2 \, t _0+  \ell(a)\,,\, 2 (D-t_0)  +  \ell(b) \,\big] = d(m_a , m_b)  + 
\frac{1}{2} \big( \ell (a) + \ell (b)\big) \ ,$$
which concludes.
\end{proof}

\begin{lemma}\label{minorLstar0}
Let $(X , d) $ be a $ \delta$-hyperbolic space, for any finite set $S\subset \text{\rm Isom} (X,d) $ such that $ L^* (S) > 4\, \delta$, for any positive value 
$\e < L^* (S) - 4\, \delta $ and any point $x$ such that $ \Max_{\g \in S}  d(x, \g x) < L^* (S) + \e$, if $a \in S$ verifies $d(x,a\,x) = \Max_{\g \in S}  d(x, \g x) $, 
then either $\ell (a) > \delta$ or there exists $b \in S \setminus \{a, a^{-1}\}$ such that $d(x, b\,x) \ge L^* (S) - \delta -\e$.
\end{lemma}

\begin{proof} If $\ell (a) > \delta$, the lemma is proved, we shall therefore suppose that $\ell (a) \le \delta$.\\
Let us fix geodesic segments $[x , a\,x]$, $[x , a^2 x]$ and $[a\, x , a^2 x] := a ([x , a\,x])$, 
and a geodesic triangle 
$\Delta = [x , a\,x ,  a^2 x ]$, whose sides are these geodesic segments. Let us consider
its approximation by a tripod $ f_\Delta : (\Delta , d) \to  (T_\Delta , d_T) $ (the construction of this approximation
is described before Lemma \ref{prodist}); let $c$ be the branching point of this tripod and let us denote 
by $\alpha , \beta , \gamma$ the respective lengths of the branches $ [c , f_\Delta ( a \,x)],\  [c , f_\Delta ( a^2 x)], 
\  [c , f_\Delta ( x)]$ of this tripod. 
By Lemma \ref{prodist}, $d(x, a\,x)$, $d(x, a^2 x)$ and $d(a\,x, a^2 x)$ are respectively equal to 
$d_T \left( f_\Delta ( x) , f_\Delta ( a \,x)\right) $, $d_T \left( f_\Delta ( x) , f_\Delta ( a^2 x)\right) $ and
$d_T \left( f_\Delta ( a\,x) , f_\Delta ( a^2 x)\right) $ and we deduce that 
$\alpha + \g = d(x, a\,x) =  d(a\, x, a^2x) = \alpha + \beta $ (which means that $\beta = \g$), and that
$d(x, a^2 x) = \g + \beta = 2 \beta$. A consequence of this and of Lemma \ref{puissances} (i) is that
$2\, \beta \le \alpha + \beta + \ell (a) + 2 \delta$ and so that
$$\alpha \ge \frac{1}{2} \big(d(x, a\,x) - \ell (a) \big) - \delta \ge  \frac{1}{2} \big( L^* (S) - 3 \delta \big)
> \frac{1}{2} (\delta +\e) \ .$$
Let $y$ be the point of $[x , a\,x]$ such that $d(x,y) =  \frac{1}{2} (\delta + \e) $, then $a\,y$ is the point
of $[a\, x , a^2 x]$ such that  $d(a\,x, a\,y) =  \frac{1}{2} (\delta + \e) $; it follows from Lemma \ref{prodist} 
that $d_T \left( f_\Delta ( a\,x) , f_\Delta ( a\,y)\right) = \frac{1}{2} (\delta + \e) < \alpha$ and that the point
$f_\Delta ( a\,y)$ is on the branch $ [c , f_\Delta ( a \,x)]$ of the tripod. It comes from this (and from the fact that
$d_T \left( f_\Delta ( x) , f_\Delta ( a\,x) \right) \ge  L^* (S) > 4\,\delta +  \e > 
d_T \left( f_\Delta ( a\,x) , f_\Delta ( a\,y)\right) + d_T \left( f_\Delta ( x) , f_\Delta ( y)\right) $) that 
$$d_T \left( f_\Delta ( y) , f_\Delta ( a\,y)\right) = d_T \left( f_\Delta ( x) , f_\Delta ( a\,x) \right) -
d_T \left( f_\Delta ( a\,x) , f_\Delta ( a\,y)\right) - d_T \left( f_\Delta ( x) , f_\Delta ( y)\right) $$
$$ <  L^* (S) + \e - \delta - \e \ .$$
It follows from Lemma \ref{proprietes} (i) that $ d(y, a\,y ) = d(y, a^{-1}y) \le 
d_T \left( f_\Delta ( y) , f_\Delta ( a\,y)\right) + \delta <  L^* (S) $, there therefore 
exists some $b \in S$ such that $d(y, b\,y) \ge  L^* (S)$, and we automatically get that $b \notin \{a, a^{-1}\}$,
the triangle inequality concludes that  $d(x, b\,x) \ge  L^* (S) - \delta - \e$.
\end{proof}

\begin{proof}[Proof of Theorem \ref{minorLstar}]
We may assume that $\Max [ \ell(a) , \ell(b)] < \delta$ (if not, the 
Theorem \ref{minorLstar} is automatically verified); hence
we have to prove that, under this assumption, every $\g \in \{ ab , a b^{-1},  ba , a^{-1} b \}$ verifies
$\ell (\g) \ge 2\, \delta$. As $L^*(\{a , b\}) = L^*(\{a , b^{-1}\}) = L^*(\{b , a\}) = L^*(\{a^{-1} , b\})$ (because 
$d(x, \g x) $ is always equal to $d(x, \g^{-1} x) $), if the pair $\{a , b\}$ satisfies the assumptions of Theorem \ref{minorLstar}, then these assumptions are also satisfied by the pairs $\{a , b^{-1}\}$, $\{b , a\}$ and $\{a^{-1} , b\}$;
it is thus sufficient to prove that, under the assumptions of Theorem \ref{minorLstar}, if 
$\Max [ \ell(a) , \ell(b)] < \delta$, then $\ell (ab) \ge 2\,\delta$. By contradiction, let us suppose that $\ell (ab) < 
2\,\delta$.\\
For any $\e> 0$ such that $\e < \delta - \Max [ \ell(a) , \ell(b)] $, let us fix some point $x$ such that 
$$ L^*(\{a , b\}) \le \Max \,[\, d(x , a \,x) \,;\,  d(x , b \,x) \,] <  L^*(\{a , b\}) + \e \ .$$
The above assumptions and the lemma \ref{minorell} then give:
\begin{equation}\label{minorLstar5}
d(a\,x , b\, x) \le  \Max \,[d(x,a\, x)\,;\, d( x, b x) ] + 6\,\delta\ .
\end{equation}
On the other hand, as the hypothesis and the choice of $\e$ imply that $ L^*(\{a , b\})  \ge \frac{23}{2} \,\delta  
> 4\, \delta + 2 \e$, the aforementioned choice of $x$ and the lemma \ref{minorLstar0} prove that 
\begin{equation}\label{minorLstar6}
 |d(x, a\, x) - d(x, b\, x)| = \Max \,[d(x,a\, x)\,;\, d( x, b\, x) ] - \Min \,[d(x,a\, x)\,;\, d( x, b\, x) ]  < \delta + 2 \, \e \ .
\end{equation}
For the sake of simplicity, let us denote by $ L^*$ the Margulis constant $ L^*(\{a , b\})$ and by $m_a$ 
(resp. by $m_b$) the middle point 
of some geodesic segment from $x$ to $a x$ (resp. from $x$ to $b x$). From the assumptions, it comes that 
$L^* > 4 \, \delta >  3\, \delta + \Max \big( \ell(a) , \ell(b)\big)$ and, applying Lemma \ref{minorLstar1}, that
\begin{equation}\label{minorLstar3}
 d(m_a , m_b) \ge L^* - \frac{1}{2} \big(\ell (a) + \ell (b)\big) - 3\,\delta > \frac{1}{2} \left(  L^* + \frac{9}{2}\, 
\delta - \Min [\ell (a) , \ell (b)]\right) \ ,
\end{equation}
where the last inequality follows from the assumption $L^*(\{a , b\}) \ge \frac{23}{2} \,\delta >  
 \frac{21}{2}  \,\delta + \Max [\ell (a) , \ell (b)]$.

Let us now consider any geodesic triangle $\Delta = [x , a\,x ,  b\,x]$ (with vertices $x , a\,x ,  b\,x$)
and its approximation by a tripod $ f_\Delta : (\Delta , d) \to  (T_\Delta , d_T) $ (the construction of this approximation
is described in the beginning of the subsection \ref{coherence}); let $c$ be the branching point of this tripod and let us denote 
by $\alpha , \beta , \gamma$ the respective lengths of the branches $ [c , f_\Delta ( a \,x)],\  [c , f_\Delta ( b \,x)], 
\  [c , f_\Delta ( x)], \ $ of this tripod. A consequence of Lemma \ref{prodist} is that $f_\Delta ( m_a)$ and 
$f_\Delta ( m_b)$ are the middle points of the sides $[f_\Delta ( x) , f_\Delta (a\, x)]$ and 
$[f_\Delta ( x) , f_\Delta (b\, x)]$ of this tripod; from this, from Lemma \ref{proprietes} (i) and from 
\eqref{minorLstar3} we deduce that
\begin{equation}\label{minorLstar4}
d_T \left( f_\Delta ( m_a) , f_\Delta ( m_b) \right) \ge d(m_a , m_b) - \delta> \frac{1}{2} \left(  L^* + \frac{5}{2}\, 
\delta - \Min [\ell (a) , \ell (b)]\right) \ .
\end{equation}
A consequence of Lemma \ref{prodist} and of \eqref{minorLstar6} and \eqref{minorLstar4} is that 
$$ |d_T \left( f_\Delta ( x) , f_\Delta ( m_a) \right) - d_T \left( f_\Delta ( x) , f_\Delta ( m_b) \right)| = 
|d(x , m_a) - d(x , m_b)| = \frac{1}{2} |d(x, a\, x) - d(x, b\, x)|$$
$$< \frac{\delta}{2} + \e< d_T \left( f_\Delta ( m_a) , f_\Delta ( m_b) \right) \ ; $$
it follows that the side $[f_\Delta ( x) , f_\Delta (a\, x)] := [f_\Delta ( x) , c] \cup [c , f_\Delta (a\, x)]$ of the tripod
contains $f_\Delta ( m_a) $, but not $f_\Delta ( m_b) $ and that the side $[f_\Delta ( x) , f_\Delta (b\, x)] := 
[f_\Delta ( x) , c] \cup [c , f_\Delta (b\, x)]$ of the tripod contains $f_\Delta ( m_b) $, but not $f_\Delta ( m_a) $, 
which proves that $ f_\Delta ( m_a)$ (resp. $ f_\Delta ( m_b)$) belongs to the branch 
$ ] c , f_\Delta (a\, x)]$ (resp. to the branch $ ] c , f_\Delta (b\, x)]$).
From this we deduce that
$$ d_T ( c ,  f_\Delta ( m_a)) = d_T \left( f_\Delta ( x) , f_\Delta ( m_a) \right) - d_T \left( c ,  f_\Delta ( x)  \right) = \frac{\alpha - \gamma}{2} \ \ ; \ \   d_T ( c ,  f_\Delta ( m_b)) = \frac{\beta - \gamma}{2} $$
(where the proof of the second equality is similar to the proof of the first one), proving by the way that 
$\g \le \Min (\beta ,\alpha)$. A consequence of these last equalities, 
of the fact that $ f_\Delta ( m_a)$ (resp. $ f_\Delta ( m_b)$) belongs to two different branches of the tripod, and of 
the first inequalities of \eqref{minorLstar4} and \eqref{minorLstar3} is that
\begin{equation}\label{minorLstar7}
\frac{\alpha - \gamma}{2}  + \frac{\beta - \gamma}{2}  =  d_T \left( f_\Delta ( m_a) , f_\Delta ( m_b) \right) 
\ge  d(m_a , m_b) - \delta  \ge L^* - \frac{1}{2} \big(\ell (a) + \ell (b)\big)  - 4\,\delta \,.
\end{equation}
As, by Lemma \ref{prodist}, $d(x, a\,x)$, $d(x, b\,x)$ and $d(a\,x, b\,x)$ coincide respectively with 
$d_T \left( f_\Delta ( x) , f_\Delta ( a \,x)\right) $, $d_T \left( f_\Delta ( x) , f_\Delta ( b \,x)\right) $ and
$d_T \left( f_\Delta ( a\,x) , f_\Delta ( b \,x)\right) $, and then with $\g + \alpha$, $\g + \beta$ and $ \alpha + \beta$
respectively, we deduce from \eqref{minorLstar5} and  \eqref{minorLstar6} that $\Min (\alpha, \beta) -  6\, \delta \le \g \le \Min (\alpha, \beta) $ and 
$\Min (\alpha, \beta)  > \Max (\alpha, \beta)  - \delta  - 2 \,\e $;
plugging these two estimates in \eqref{minorLstar7}, and recalling that $ \Max [ \ell(a) , \ell(b)] < \delta - \e $ by 
the choice of $\e$, we get
$$\frac{\delta}{2} + \e > \frac{1}{2} \big(\Max (\alpha, \beta)  - \Min (\alpha, \beta)\big) = \frac{\alpha + \beta}{2} 
- \Min (\alpha, \beta) \ge  \frac{\alpha + \beta}{2} -\g  - 6 \, \delta  = \frac{\alpha - \gamma}{2}  + 
\frac{\beta - \gamma}{2} - 6 \, \delta $$
$$\ge L^* -\frac{1}{2} \big(\ell (a) + \ell (b)\big)  - 10\,\delta > L^* - 11\, \delta + \e\ ,$$
in contradiction with the assumption 
$L^* \ge \frac{23}{2}\,\delta$. Hence the assumption $\ell(ab) < 2\, \delta$ is false, and this 
proves the first part of Theorem \ref{minorLstar}.

\emph{Proof of the second part of Theorem \ref{minorLstar}:}
if the subgroup $\langle a , b \rangle$ generated by $a$ and $b$ is discrete and not virtually cyclic,
then firstly $\langle a , b a b^{-1} \rangle$ is not virtually cyclic by Proposition \ref{actionelementaire} (v) 
and secondly, as $\langle a ,b\rangle = \langle ab ,b\rangle$ is not virtually cyclic, $\langle ab , ba\rangle = \langle ab , b (ab) b^{-1}\rangle $ is not virtually cyclic 
too by Proposition \ref{actionelementaire} (v). From the first part of Theorem \ref{minorLstar}, we know that

$\bullet$ either $\ell ( b a b^{-1}) = \ell(a) \ge \delta$, and the proposition \ref{granddeplacementx} (ii) then implies 
that one of the two semi-group generated by $\{a^{14}, ba^{14} b^{-1}\}$ or by $\{a^{14}, ba^{-14} b^{-1}\}$ is 
free, 

$\bullet$ either $\ell(b) \ge \delta$, and the proposition \ref{granddeplacementx} (ii) then implies 
that one of the two semi-group generated by $\{b^{14}, a b^{14} a^{-1}\}$ or by $\{b^{14}, a b^{-14} a^{-1}\}$ is 
free,

$\bullet$ or $\ell (ba) = \ell(ab) \ge 2\, \delta$, and the proposition \ref{granddeplacementx} (ii) then implies 
that one of the two semi-group generated by $\{(ab)^{7}, (ba)^{7}\}$ or by $\{(ab)^{7}, (ba)^{-7}\}$ is free.
\end{proof}

\begin{lemma}\label{minorLetoile1}
Let $(X , d) $ be a $ \delta$-hyperbolic space, for any pair $\{a , b\}$ of isometries of $(X , d) $,
and for any point $x \in X$ such that $d(x,a\,x) \ge d(x,b\,x)$, if $\ell (a)$, $\ell (b)$, $\ell (ab)$, $\ell (b^{-1} a)$ 
and $\ell (b^{-2} a)$ 
are (strictly) bounded above by $\delta$, then the middle point $m$ of any geodesic segment from $x$ to $ax$ satisfies
$ \Max [d(m, a\, m) \,; \, d(m,b\,m)] < \frac{31}{2}\,\delta$.
\end{lemma}

Though it looks similar to Theorem \ref{minorLstar} (assuming however stronger hypotheses), this Lemma is original for the point $m$ only depends on $a$ and 
not on $b$. This point will be important in the proof of  Theorem \ref{minorLetoile}, which extends Theorem \ref{minorLstar} to the case where the pair of 
isometries $\{a , b\}$ is replaced by any finite set $S$ of isometries.

\begin{proof} Lemma \ref{milieux} and the assumptions prove that
\begin{equation}\label{minorLetoile2}
d(m, a\, m) \le 3\,\delta + \ell (a) < 4\, \delta
\end{equation}
It is thus sufficient to prove that $d(m , b\, m) < \frac{31}{2}\, \delta$.
As $\ell (a), \,\ell (b),\, \ell (ab) < \delta$, Lemma \ref{minorell} yields
\begin{equation}\label{minorLetoile3}
d(a x , b x) \le  \Max \,[d(x,a x)\,;\, d(x, b x) ] + \dfrac{11}{2}\,\delta = d(x, a x) + \dfrac{11}{2}\,\delta \ .
\end{equation}
Similarly, as  $\ell (b^{-1}), \, \ell (b^{-1} a) ,\, \ell (b^{-2} a) < \delta$, Lemma \ref{minorell} also gives
$$ d (b^{-1} x , b^{-1} a x) \le  \Max \,[d(x, b^{-1} x)\,;\, d(x ,  b^{-1} a  x) ] + 
\dfrac{11}{2}\,\delta \ ,$$
hence
\begin{equation}\label{minorLetoile4}
d(x , ax) \le  \Max \,[d(x,b x)\,;\, d(a x, b x) ] + \dfrac{11}{2}\,\delta \ .
\end{equation}
Let us fix geodesic segments $[x , a\,x]$, $[x , b\,x]$ and $[a\, x , b\,x]$, and a geodesic triangle 
$\Delta = [x , a\,x ,  b\,x ]$, whose sides are these geodesic segments. Let us consider
its approximation by a tripod $ f_\Delta : (\Delta , d) \to  (T_\Delta , d_T) $ (the construction of this approximation
is described in the beginning of the subsection \ref{coherence}); let $c$ be the branching point of this tripod and let us denote 
by $\alpha , \beta , \gamma$ the respective lengths of the branches $ [c , f_\Delta ( a \,x)],\  [c , f_\Delta ( b \,x)], 
\  [c , f_\Delta ( x)], \ $ of this tripod. 
By Lemma \ref{prodist} $d(x, a\,x)$, $d(x, b\,x)$ and $d(a\,x, b\,x)$ are respectively equal to 
$d_T \left( f_\Delta ( x) , f_\Delta ( a \,x)\right) $, $d_T \left( f_\Delta ( x) , f_\Delta ( b \,x)\right) $ and
$d_T \left( f_\Delta ( a\,x) , f_\Delta ( b \,x)\right) $ and we get:
\begin{equation}\label{minorLetoile5}
d(x, a\,x) = \g + \alpha \ \ \ ,\ \ \ d(x, b\,x) = \g + \beta \ \ \ ,\ \ \ d(a\, x, b\,x) = \alpha + \beta \ .
\end{equation}
The assumption $d(x,b\,x) \le d(x,a\,x)$ means that $\beta \le \alpha$. A consequence of \eqref{minorLetoile5}
and \eqref{minorLetoile3} is that $\beta \le  \gamma + \frac{11}{2}\,\delta $, while \eqref{minorLetoile5} and
\eqref{minorLetoile4} imply that $\Min(\g , \alpha) \le \beta + \frac{11}{2}\,\delta $. We summarize all these estimates
in the inequalities:
\begin{equation}\label{minorLetoile6}
\min (\gamma , \alpha) - \frac{11}{2}\,\delta\le \beta \le \min \left(\g+ \frac{11}{2}\,\delta \,,\, \alpha \right)
\end{equation}
$\bullet$ {\bf First case: if $\alpha \le \g$:} Let us denote by $m_1$ the middle point of $[x , b\,x]$, as $\beta \le \alpha$,
one has $ 0 < \dfrac{\g + \beta}{2} \le \dfrac{\g + \alpha}{2} \le \g$, which means that $ f_\Delta ( m) $ and 
$ f_\Delta ( m_1) $ both belong to the branch $ [c , f_\Delta ( x)]$ of the tripod $T$ and satisfy
$$d_T \left(  f_\Delta ( m) \,,\,  f_\Delta ( m_1)  \right) = \dfrac{\g + \alpha}{2} - \dfrac{\g + \beta}{2} \le 
\frac{11}{4}\,\delta \ ,$$
where the last of these inequalities is derived from \eqref{minorLetoile6}. A consequence of this and of the 
Lemma \ref{proprietes} (i) is that $d(m, m_1) \le \frac{15}{4}\,\delta$. 
On the other hand, applying the assumption
and Lemma \ref{milieux}, it comes: $d(m_1 , b\, m_1) \le 3 \, \delta + \ell (b) < 4 \, \delta $. These two last estimates 
and the triangle inequality give: 
$$ d(m , b\,m) \le d(m_1 , b\, m_1) + 2\, d(m, m_1) < \frac{23}{2}\,\delta \ .$$
$\bullet$ {\bf Second case: if $\alpha > \g$:}  Let us denote now by $m_2$ the middle point of $[a\, x , b\,x]$; 
as $ \dfrac{\g + \alpha}{2}$ and $\dfrac{\alpha + \beta}{2}$ are both smaller than $\alpha$, the points $ f_\Delta ( m) $ and 
$ f_\Delta ( m_2) $ both belong to the branch $ [c , f_\Delta ( a\, x)]$ of the tripod $T$ and satisfy
$$d_T \left(  f_\Delta ( m) \,,\,  f_\Delta ( m_2)  \right) = \left|\dfrac{\g + \alpha}{2} - \dfrac{\alpha+ \beta}{2} \right| \le 
\frac{11}{4}\,\delta \ ,$$
because $ |\beta - \g| \le  \frac{11}{2}\, \delta$ by \eqref{minorLetoile6}. This inequality and Lemma \ref{proprietes} (i) yields $d(m, m_2) \le \frac{15}{4}\,\delta$.
The image by $a^{-1}$ of $[a\, x , b\,x]$ being denoted by $[x , a^{-1} b x]$, whose middle point is 
$m'_2 := a^{-1} m_2$, Lemma \ref{milieux} and the fact that $\ell ( a^{-1} b) = \ell ( b^{-1} a) < \delta$ by
hypothesis imply that
$$ d(m_2 , b a^{-1 }m_2) = d(m'_2 , a^{-1 } b \,m'_2) \le 3\,\delta + \ell ( a^{-1} b)<4\,\delta \ .$$
Using this last inequality and the aforementioned inequality $d(m, m_2) \le \frac{15}{4}\,\delta$, we get:
$$ d(m , b a^{-1 }m) \le  d(m_2 , b a^{-1 }m_2) + 2 \, d (m, m_2) < \frac{23}{2}\,\delta \ ,$$
This last estimate, together with \eqref{minorLetoile2} and with the triangle inequality gives:
$$d(m, b\,m) \le d(m , b a^{-1 }m) + d( b a^{-1 }m , b\, m) < \frac{23}{2}\,\delta  + d( m , a\, m)<  
\frac{31}{2}\,\delta \ .$$
This upper bound of $d(m, b\,m)$, together with \eqref{minorLetoile2}, ends the proof.
\end{proof}

{\bf b) When $L^*$ is applied to any finite set of isometries}

For any finite set $\Sigma$ of elements of a given group $G$, we denote by $\Sigma^{-1}$ the set 
$\{ \g : \g{-1} \in \Sigma\}$; $\Sigma$ is said to be \emph{\lq \lq symmetric"} if $\Sigma^{-1} = \Sigma$.
To any such finite set $\Sigma$, one associates its "symmetrized set" $S := \Sigma \cup \Sigma^{-1} $;
notice that the groups generated by $ \Sigma $ and by $ S = \Sigma \cup \Sigma^{-1} $ coincide, that
their Cayley graphs verify $ {\cal G}_S = {\cal G}_\Sigma$, and thus that their algebraic word distances 
$ d_S$ and $d_\Sigma$ (associated respectively to the sets of generators $ S$ and $\Sigma$, see their definition 
in section \ref{notations}) coincide.

\begin{theorem}\label{minorLetoile}
Let $(X , d) $ be a $ \delta$-hyperbolic space, for any finite symmetric set $S$ of isometries, if $L^* (S) \ge 
\frac{31}{2} \delta$, then 
\begin{itemize}
\item[(i)] there exists $ \g_0 \in S^3$ such that $\ell (\g_0) \ge \delta$
\item[(ii)] moreover, if the subgroup generated by $S$ is discrete and not virtually cyclic, there exists $\sigma \in S$ 
such that 
one of the two semi-groups generated by $ \{\g_0^{14} \, , \,  \sigma \g_0^{14} \sigma^{-1}\}$ or by 
$ \{\g_0^{14}\, , \, \sigma \g_0^{-14} \sigma^{-1}\}$ is free.
\end{itemize}
\end{theorem}

\begin{proof} Arguing by contradiction, let us suppose that $\ell (\g) < \delta$ for every $ \g \in S^3$, let us fix any point
$x \in X$ and denote by $a$ an element of $S$ such that $d(x,a\,x) = \Max_{g \in S}  d(x, g x) $. For any $g \in S$,
we have $d(x,a\,x) \ge d(x, g x) $, and  (by assumption) $\ell (a)$, $\ell (g)$, $\ell (ag)$, $\ell (g^{-1} a)$ and 
$\ell (g^{-2} a)$ are (strictly) bounded above by $\delta$, Lemma \ref{minorLetoile1} then proves that the middle 
point $m$ of any geodesic segment from $x$ to $ax$ satisfies $ d(m,g\,m)< \frac{31}{2}\,\delta$; as this inequality 
is valid for every $g \in S$, we deduce that  $L^* (S) < \frac{31}{2} \delta$, a contradiction with the hypothesis 
which proves (i).\\
From (i), there exists $ \g_0 \in S^3$ such that $\ell (\g_0) \ge \delta$, on the other hand, there exists $\sigma \in S$
such that the subgroup generated by $\{\g_0 , \sigma\}$ is not virtually cyclic, otherwise (by the 
Proposition \ref{actionelementaire} (vii)) the subgroup generated by $S$ would be virtually cyclic. 
By the proposition \ref{actionelementaire} (v), we deduce that
$ \{\g_0 \, , \,  \sigma \g_0 \sigma^{-1}\}$ generates a non virtually cyclic discrete subgroup.
Hence we may apply the proposition \ref{granddeplacementx} (ii) to the pair $ \{\g_0 \, , \,  \sigma \g_0 \sigma^{-1}\}$,
which proves that one of the two semi-groups generated by $ \{\g_0^{14} \, , \,  \sigma \g_0^{14} \sigma^{-1}$ or by 
$ \{\g_0^{14}\, , \, \sigma \g_0^{-14} \sigma^{-1}\}$ is free.
\end{proof}



\subsubsection{When the  Margulis constant  $L$ is bounded below:}\label{freeL}

Let us now assume that $a$ and $b$ are hyperbolic isometries, we then have

\begin{lemma}\label{perpendiculaire} If the group $\langle a, b \rangle$ generated by $a$ and $b$ is a non virtually cyclic discrete subgroup of the isometry group of $(X,d)$, then,  
for all $R> 0$ such that $M_R (a)$ and  $M_R (b)$ are non empty, there exists points $x_0 \in M_R (a) $ and $ y_0 \in M_R (b) $ such that $ d(x_0, y_0 )  $ is the minimum of $d(x, y )$ 
among $(x,y)\in M_R (a) \times M_R (b)$.
\end{lemma}

\begin{proof} 
Let us denote by $a^+$ and $a^-$ (resp. $b^+$ and $b^-$) the points in the limit set of $a$ (resp. of $b$).
Let us assume that there exists a sequence $(x_n, y_n) $ of elements of $M_R (a) \times M_R (b)$ which goes to infinity and such that $d(x_n, y_n)$ is bounded.
In the following we assume that $x_n$ goes to infinity, the argument would be the same with $y_n$. In this case, Lemma \ref{tube} (ii) shows the existence of a subsequence $x_{n_k}$ 
which converges to $a^+$ or $a^-$ and, as $ d(x_{n_k}, y_{n_k})$ is bounded, it follows that the sequence $ y_{n_k} $ converges towards $ a^+ $ or $ a^- $. Lemma
\ref{tube} (ii) shows that the limit of  $ y_{n_k} $ can only be $ b^+ $ or $ b^- $, hence we have 
 $ \{a^- , a^+\} \cap \{b^- , b^+\} \ne \emptyset $,  which implies (by Proposition 
\ref{actionelementaire} (i)) that $ \{a^- , a^+\} = \{b^- , b^+\}  $ and consequently that all  elements of the group
 $ \langle a , b \rangle $ have $ \{a^- , a^+\} $ as fixed points set. Then, from Proposition
\ref{actionelementaire} (ii) we deduce that 
 $ \langle a , b \rangle $ is virtually cyclic. This contradicts the  hypothesis, then $ d(x,y) $ goes 
to $ +\infty $ when $ (x,y)  $ does. This ensures the existence of a point $(x_0 , y_0 ) \in M_R (a) \times M_R (b) $ where the function $ (x,y) \mapsto d(x , y )$ achieves its minimum. 
 \end{proof}

Clearly, Lemma \ref{perpendiculaire} is trivial when $ M_R (a) \cap M_R (b) \ne \emptyset $ since then $(x_0 , x_0 )$  with $x_0 \in M_R (a) \cap M_R (b)$ is a solution.
It is thus important to have a criterion characterizing those values of $ R $ such that $ M_R (a) \cap M_R (b)= \emptyset $. This is the goal of the next result. 

\begin{lemma}\label{disjoints} Let $ R>0$ be such that $M_R (a)$ and $M_R (b)$ are non-empty. If $L(a,b) > R $, the Margulis domains $ M_R (a) $ and $ M_R (b) $ are disjoints.
Conversely, if the Margulis domains $ M_R (a) $ and $ M_R (b) $ are  disjoints, then $ L(a,b) \ge R $.
\end{lemma}

\begin{proof} If there exists  $x\in M_R (a)\cap M_R (b) $, there exists $ (p,q) \in \mathbb Z^* \times \mathbb Z^* $ such that
 $ d(x , a^p  x) \le R $ and $ d(x , b^q  x) \le R$. From the definition de $ L(a,b) $, this implies that $ L(a,b) \le R$.
 
 Conversely, if the Margulis domains $ M_R (a) $ and $ M_R (b) $ are disjoints, then
every point $ x \in X $ satisfies $x\in X\setminus M_R (a)$ or $x\in X\setminus M_R (b)$,
which implies that, for all $p, q \in \mathbb Z^*$, $d(x , a^p x ) > R$ or 
 $d(x, b^q x ) > R$. Taking the infimum, we deduce that $L(a,b)\ge R$. 
 \end{proof}

\begin{theorem}\label{margulislibre}
Let $ (X , d) $ be a  $ \delta$-hyperbolic space. Let $a$ and $b$ be two hyperbolic isometries acting on $ (X,d) $ such that the subgroup $\langle a, b \rangle$ is discrete and non-virtually cyclic. If the Margulis constant satisfies $L(a,b) > 23 \delta$, we then have the following alternatives:
   \begin{itemize}
      \item [(i)] if $ \ell(a), \ell (b)\le 13\delta$, the subgroup generated by $\{a, b\}$ is free,
      \item [(ii)] if $\ell(a), \ell (b) > 13 \delta$, one of the two semi-groups generated by $ \{ a,  b\} $ or by $ \{ a,  b^{-1}\} $ is free,
      \item [(iii)] if $\ell(a)\le 13 \delta$ and $\ell (b) > 13\delta$, one of the two semi-groups generated by $ \{b,   a b a^{-1}\} $ or by $ \{b,   a  b^{-1} a^{-1}\} $ is free,
      \item [(iv)] if $\ell(a)> 13\delta$ and $\ell (b)\le 13 \delta$, one of the two semi-groups generated by $\{a,   b a b^{-1}\} $ or by $\{a,   b  a^{-1} b^{-1}\}$ is free.
  \end{itemize}
\end{theorem}

The proof relies of the following proposition:

\begin{prop}\label{libre3}
Let $ (X, d) $ be a  $\delta$-hyperbolic space.
Let $a$ and $b$ be two hyperbolic isometries acting on $(X,d)$ such that the subgroup $\langle a, b \rangle$ is discrete and non-virtually cyclic. If \footnote{In this Proposition,  the hypothesis $L(a, b)-\Max [\ell(a),   \ell (b) ]> 10 \delta$, can be replaced by: there does not exist points $ x \in X $ and integers $p, q \in \N^*$ which satisfies simultaneously $d(x, a^p x) \le 10 \delta + \Max   [ \ell(a),   \ell(b) ]$ and $d(x, b^q x) \le 10 \delta + \Max   [ \ell(a),   \ell(b) ]$. The equivalence between this new statement and Proposition \ref{libre3} 
follows from Lemma \ref{disjoints} and the beginning of the proof of Proposition \ref{libre3}.} the Margulis constant $L(a,b)$ satisfies
 $L(a,b)-\Max [\ell(a),   \ell(b) ]>10 \delta$, then $a$ and $b$ generate a free subgroup  of the isometry group of $(X, d)$.
\end{prop}

\begin{proof}[Proof of Proposition \ref{libre3}] 
For the sake of simplicity let  $ \ell_0=Max [ \ell  (a), \ell (b) ]$. The hypothesis $L(a,b) > 10 \delta + \ell_0$ allows to choose $R_0$ such that 
$10 \delta + \ell_0 < R_0 < L(a,b)$. Let us denote by $\e$ any real number such that
 $ 0 < \e < R_0 - (10 \delta + \ell_0) $ and define $ \e' >0$ by $ R_0 = 10 \delta + \ell_0 +\e + 2  \e' $. We also set $ r_0 := \ell_0 + \delta + \e$,
Lemmas \ref{MRferme} (ii) and \ref{disjoints} implies that $ M_{r_0} (a) $ and $ M_{r_0} (b) $ are closed disjoints and non empty sets. Furthermore, from Lemma \ref{perpendiculaire}, 
we can choose two points  
 $ x_0 \in M_{r_0} (a) $ and $ y_0 \in M_{r_0} (b) $ such that $ d(x_0, y_0 )  $ is the minimum of
 $d(x , y )$ when $ (x, y) $ runs through $ M_{r_0} (a) \times M_{r_0} (b) $. We then fix a geodesic $ [x_0 , y_0] $ between these two points. Let us first prove the following property:
\begin{equation}\label{xecarte}
\exists  x\in  [x_0 , y_0] \;\text{ such that }\;\forall (p,q) \in \mathbb Z^* \times \mathbb Z^*,\quad d(x, a^p  x ) > R_0\; {\rm and }\; d(x, b^q  x ) >  R_0\,.
\end{equation}
Indeed, one has $L(a,b) >  R_0$, Lemma \ref{disjoints} then implies that $ M_{R_0} (a) $ and 
 $ M_{R_0} (b) $ are closed disjoints sets  containing respectively $ x_0 $ and $ y_0 $. Their intersections with $ [x_0 , y_0] $ are then closed disjoints and non empty whose union 
 cannot be equal to $ [x_0 , y_0] $. Consequently, there exists a point $ x \in  [x_0, y_0] $ which is not in 
 $ M_{R_0} (a) \cup M_{R_0} (b) $ and hence satisfies Property \eqref{xecarte}.

Let us now fix such a point $ x \in  [x_0, y_0]  $ (given by Property \eqref{xecarte}).
The fact that $x \notin M_{R_0} (a) $, that $ x_0 $ is a projection of $ x $ on $ M_{r_0} (a) $, that
 $ a^p x_0 \in M_{r_0} (a) $, and Lemma \ref{distants} imply that
\begin{equation}\label{xeloigne}
\forall p \in \mathbb Z^* \ \ d(x, a^p  x_0) \ge d( x, x_0)  > \frac{1}{2}  (R_0 - r_0)
= \frac{9}{2} \delta +  \varepsilon' \,.
\end{equation}

We choose $ p $ and $ q $ in $ \mathbb Z^* $.
We denote by $ [x, x_0] $ the segment of the geodesic $ [x_0, y_0] $ between $ x $ and $ x_0 $ 
and we call $ [a^p  x, a^p  x_0] $ the image of this geodesic by $ a^p $.
Let us choose arbitrary geodesics $ [x, a^p  x] $, $ [x_0, a^p  x_0] $ and $ [x , a^p   x_0] $.
Inequalities (\ref{xecarte}) and (\ref{xeloigne}), show that there exist points $ u $, $ u' $ and 
 $ u'' $, respectively on the geodesics $ [x, a^p  x] $, $ [x, a^p   x_0] $ 
and $ [x,  x_0] $, such that $ d(x, u) =  d(x, u') =  d(x, u'') =  \frac{5}{2} \delta + \varepsilon' $.

We now consider the triangles $ \Delta = [x , a^p  x,  a^p  x_0] $ and $ \Delta' = [x,  a^p x_0, x_0] $ 
defined by the above geodesics and their approximations by the associated tripods $ f_\Delta : (\Delta, d) \to  (T_\Delta, d_T) $ and 
 $ f_{\Delta'}: (\Delta', d) \to (T_{\Delta'}, d_{T'}) $ (see Sub-section \ref{coherence}); the branching points of these tripods are respectively denoted by $ c $ and $ c' $. The length of the sides of the tripod 
 $ (T_\Delta, d_T) $ (resp. $ (T_{\Delta'}, d_{T'}) $) whose ends are  
 $ f_\Delta (x), f_\Delta ( a^p  x) $ and $ f_\Delta ( a^p  x_0) $ (resp. $ f_{\Delta'} (x), f_{\Delta'} ( a^p  x_0) $ and $ f_{\Delta'} (x_0) $) are denoted by $ \alpha ,\,\beta $ and $ \gamma $ 
(resp. $ \alpha' ,\,\beta'  $ and $ \gamma' $). 
Finally we denote by $ c'_2  $ the inverse image in $ [ x_0 , a^p (x_0) ] $ by $ f_{\Delta'} $ 
of the branching point $ c' $ of the tripod $ T_{\Delta'} $. Lemma \ref{quasiconvexe} (iii) ensures that 
 $ c'_2 \in  M_{r_0 + 2 \delta} (a) $ and, as $ x  $ belongs to the closure of $ X \setminus M_{R_0} (a) $, one has (from Lemma \ref{distants})
 $$\alpha' = d_{T'}\big(f_{\Delta'} (x),  c' \big) \ge d(x, c'_2) - \delta \ge  
d \big( x, M_{r_0 + 2 \delta} (a)\big) - \delta \ge \frac{1}{2}  (R_0 - r_0 -  2 \delta ) - \delta= \frac{5}{2}  
\delta + \e',$$
where the last equality follows from the choice of $ R_0 $ and $ r_0 $. It then follows that the points $ f_{\Delta'} (u') $ and
 $ f_{\Delta'} (u'') $ both belong to the side $ [c',  f_{\Delta'} (x)] $ of the tripod 
 $ (T_{\Delta'}, d_{T'}) $ and hence that $ f_{\Delta'} (u') = f_{\Delta'} (u'') $. The approximation lemma \ref{proprietes} (i) then yields 
$ d(u', u'') \le \delta $.
 
On the other hand, the map $ f_\Delta $ being an isometry in restriction to each side of $ \Delta $ (cf. Lemma \ref{prodist}), one has
 $$\alpha + \gamma = d(x, a^p   x_0) \ge d(x, x_0) = d( a^p   x, a^p   x_0) = \beta + \gamma , $$
which ensures that $ \alpha \ge \beta$; since we also have $ \alpha + \beta = 
 d(x, a^p   x ) > R_0 $, Inequality \eqref{xecarte} implies that 
 $$\alpha > \frac{R_0}{2} > 5  \delta + \varepsilon' > d(x, u) = d(x, u') =  d_{T} ( f_{\Delta} (x), f_{\Delta} (u))
=  d_{T} ( f_{\Delta} (x), f_{\Delta} (u'))\,.$$
This shows that the points $ f_{\Delta} (u) $ and $ f_{\Delta} (u') $ are on the side 
 $ [c,  f_{\Delta} (x)] $ of the tripod $ T_{\Delta} $ and satisfy $ f_{\Delta} (u) = f_{\Delta} (u')$; 
from Lemma \ref{proprietes} (i) we get that $ d(u, u') \le \delta$.
This last inequality together with $ d(u', u'') \le \delta$ previously proved show that $d(u, u'') \le 2 \delta$, by the triangle inequality.
We denote by $ [x, y_0] $ the segment of the geodesic $ [x_0 , y_0] $ between $ x $ and $ y_0$. Replacing
$a$ by $b$ and $ p $ by $ q $ in the previous argument, we construct geodesics
$ [x, b^q  x] $ and $ [x, b^q   y_0] $ and points $ v $, $ v' $ and $ v'' $, respectively on the geodesics
 $ [x, b^q  x] $, $ [x, b^q   y_0]$ and $ [x,  y_0] $, such that
 $ d(x, v) =  d(x, v') =  d(x, v'') =  \frac{5}{2} \delta + \varepsilon'$.
As before we prove that $ d(v, v'') \le 2 \delta$.

This last estimate, the inequality $ d(u, u'') \le 2 \delta$, the triangle inequality and the fact that $ u''$,
 $ x $ and $ v'' $ belong to the same minimizing geodesic $[x_0, y_0] $ yield
\begin{equation}\label{ecart1}
d(u, v) \ge d(u'', x) + d( x, v'') - d(u, u'') - d(v, v'') \ge 5\delta + 2\varepsilon' - 4\delta = \delta +
2  \varepsilon'\,.
\end{equation}
Let us now consider the triangle $ \Delta = [x, a^p  x,  b^q  x] $ and the associated tripod  
 $ f_\Delta : (\Delta , d) \to  (T_\Delta , d_T) $, whose sides with ends $ \ f_\Delta ( a^p  x)$, $ f_\Delta ( b^q  x)$ and  $ f_\Delta (x)  $ 
have length denoted by $ \alpha, \beta $ and  $ \gamma$ respectively, and we denote by $c $ the branching point of the tripod $ T_\Delta$.
From Inequality (\ref{ecart1}) and Lemma \ref{proprietes} (i), we deduce that
\begin{equation}\label{ecart2}
d_T\left(f_\Delta (u), f_\Delta ( v)\right) \ge 2 \varepsilon' > 0\,.
\end{equation}
The application $f_\Delta$ restricted to each side of $ \Delta $ is an isometry \
(cf. Lemma \ref{prodist}), we then obtain:
 $$d_T\left(f_\Delta (x), f_\Delta (u)\right) = d(x, u) = \frac{5}{2} \delta + \varepsilon' 
= d( x, v ) = d_T\left(f_\Delta (x), f_\Delta (v)\right)\,.$$

 \begin{itemize}
    \item If $ \gamma \ge \frac{5}{2} \delta + \varepsilon' $, then $ f_\Delta (u)$ and $f_\Delta (v)$ both belong to the same side $ [c, f_\Delta (x) ] $ of the tripod $ T_\Delta $, and hence satisfy 
 $ f_\Delta (u) = f_\Delta (v)$, which contradicts Inequality (\ref{ecart2}).
    \item If $ \frac{5}{2} \delta < \gamma < \frac{5}{2} \delta + \varepsilon'$, then $ f_\Delta (u) $ and $ f_\Delta (v) $ respectively belong to the sides 
 $ [c , f_\Delta (a^p  x) ] $ and $ [c , f_\Delta (b^q  x) ] $ of the tripod $ T_\Delta $, and then
 $$d_T\left(c, f_\Delta (v)\right) =  d_T\left(c, f_\Delta (u)\right) = d_T\left(f_\Delta (x), f_\Delta (u)\right) - d_T\left(f_\Delta (x), c\right) = \frac{5}{2} \delta 
+ \varepsilon' - \gamma < \varepsilon'\,,$$
which implies that $ d_T\left(f_\Delta (u), f_\Delta ( v)\right) < 2 \varepsilon' $, in contradiction with Inequality (\ref{ecart2}).
 \end{itemize}
 
The only possibility is then that 
\begin{equation}\label{ecart3}
 \gamma \le \frac{5}{2} \delta\,.
\end{equation}
On the other hand we have
 $ \Min (\alpha , \beta ) = \Min \left[ d(x, a^p  x);  d(x, b^q  x)\right] - \gamma > R_0 - \gamma$, which yields
 $$d(a^p  x, b^q  x) = \alpha + \beta = \Max (\alpha, \beta ) + \Min (\alpha, \beta ) >  \Max (\alpha + \gamma, \beta +\gamma) + R_0 - 2 \gamma\,, $$
which, using the definition of $ R_0 $ and Inequality (\ref{ecart3}), shows
 $$d(a^p  x, b^q  x) - \Max \left[d(x, a^p  x);  d(x, b^q  x)\right] > R_0 - 2 \gamma \ge 10 \delta + \varepsilon + 2  \varepsilon' - 5 \delta > 5  \delta\,.$$
This last inequality being valid for all $(p, q) \in \mathbb Z^* \times \mathbb Z^*$,
we conclude by applying  Proposition \ref{libre1} (i). 
\end{proof}
 
\begin{proof}[End of the proof of Theorem \ref{margulislibre}]
  We now have the following alternatives:
  \begin{itemize}
    \item If $\ell(a), \ell (b) \le  13 \delta$, as $ L(a,b) > 23  \delta $, we have $ L(a,b) - \Max   [ \ell(a),   \ell (b) ]> 10 \delta$ and Proposition \ref{libre3} implies that the subgroup generated by $ \{a, b\} $ is free.
    \item If $\ell(a), \ell (b) >  13 \delta$, Lemma \ref{quasigeod} (i) ensures that $s(a),  s(b) >  13 \delta$  and Proposition \ref{granddeplacementx} (i) then implies that one of the two semi-groups generated by 
 $ \{ a,  b\} $ or by $ \{ a,  b^{-1}\} $ is free.
    \item If $ \ell(a) \le  13 \delta $ and $ \ell (b) > 13  \delta $, then $ \ell (a  b a^{-1}) =  \ell (b) > 13 \delta$ and Lemma \ref{quasigeod} (i) ensures that $ s(b),  s(a  b a^{-1}) >  13 \delta$. Proposition \ref{granddeplacementx} (i) then implies that one of the two semi-groups generated by $ \{b,   a b a^{-1}\} $ or by $ \{b,   a  b^{-1} a^{-1}\} $ is free.
    \item If $\ell(a)  >  13 \delta$ and $\ell (b) \le 13  \delta $, exchanging the roles of $a$ and $b$, the same proof shows that one of the two semi-groups generated by
 $\{a,   b a b^{-1}\} $ or by $ \{a,   b  a^{-1} b^{-1}\} $ is free. 
  \end{itemize}
 \end{proof} 

\subsection{Free semi-groups for convex distances}

In this subsection we fix $ \delta >0$, $H>0$ and $D>0$ and we study the $ \delta$-hyperbolic spaces $(X, d )$
which are \emph{Busemann spaces} (i. e. the distance $d$ is convex in the sense of  Definition \ref{dconvexe}) endowed with a proper action by isometries of a group $ \Gamma $. We assume that the entropy  
of $(X, d)$ and the diameter of $ \Gamma \backslash  X$ are bounded above by $ H $ and $ D $ respectively. To these parameters $ \delta $, $ H $ and $ D $,
we associate a function $s_0( \delta, H, D) > 0 $ defined by Equality \eqref{defminorant}.

Contrarily to Subsection \ref{free1}, the goal of the present subsection is to state a result which apply to \emph{every pairs} of elements of $ \Gamma $ without any restriction on $ \ell(a)$,
$ \ell(b) $ or on the Margulis constants $L(a,b)$ or $L^*(a,b)$. The price to pay is a slightly stronger 
hypothesis on the geometry of the metric space under consideration.

\begin{theorem}\label{freesemigroup} 
 Let  $ (X, d )$ be any connected, geodesically complete, Busemann, non elementary $\delta$-hyperbolic space.  
Let $ \Gamma $ be any group acting properly by isometries on (X, d) such that the diameter of $ \Gamma \backslash X $ and the entropy of $(X, d)$ are respectively bounded above by $ D $ and $ H $. For every pair $a, b$ of torsion-free elements of $\Gamma $ which generates a non virtually cyclic subgroup, 
for every integers $ p,  q  \ge   \dfrac{13 \delta}{s_0(\delta , H, D)} $ one of the two semi-groups generated by $ \{ a^p,  b^q\}$ or by $ \{ a^p,  b^{-q}\}$ 
is free.
\end{theorem}

The proof of this theorem follows from a uniform lower bound of the asymptotic displacement,  whose proof will be given in Theorem \ref{Cat2}, that we assume here:

\begin{proof} For sake of simplicity, define $s_0 :=  s_0 (\delta , H, D)$.
Theorem \ref{Cat2} (i) shows that every pair of torsion-free elements $a, b \in \Gamma^* $ verifies $ \ell (a), \ell (b) > s_0 $; for every 
$ p,  q  \ge   \dfrac{13 \delta}{s_0} $, Corollary \ref{granddeplacement} then allows to deduce that 
one of the two semi-groups generated by $\{ a^p ,  b^q \}$ 
or by $\{ a^p , b^{-q}\}$  is free.
\end{proof}

\section{Co-compact actions on Gromov-hyperbolic spaces}\label{deltahyp}

\small
\emph{In this section, for any fixed positive constants $\delta$, $H$ and $D$, we are concerned by all the
$\delta$-hyperbolic spaces $(X,d )$ (which are geodesic and proper spaces by Definition \ref{hypdefinition}) 
and by every proper action (by isometries) of a group $\Gamma$ on $(X,d )$ such that the entropy of $(X,d )$ 
and the diameter of $\Gamma \backslash  X$ are respectively bounded above by $H$ and $D$.
Let us recall that $\Gamma \backslash  X$ is then compact (by Lemma \ref{autofidele} (ii)) and that, if $\Gamma$ is torsion-free, its action on $X$ is faithful and fixed point free (by Lemma \ref{autofidele} (iv)).\\
Let us recall that  $\Sigma_r (x) := \{ \g \in \Gamma^*  :  d( x, \gamma\, x) \le r\}$ and that $\Gamma_{r}(x)$ is the 
subgroup generated by $\Sigma_r (x) $ (see Definitions \ref{Gammaepsilon}).
When the action of $\Gamma$ on $(X,d )$ is co-compact, we recall that the entropy of $(X,d )$ can be computed
for any Borel $\Gamma$-invariant measure and is independent of this measure.}

\normalsize

\subsection{A Bishop-Gromov inequality for Gromov-hyperbolic spaces}\label{doublehyp}

By Definitions \ref{doublefaible}, a doubling property concerns balls whose radius lies in a given interval and whose doubling amplitude
$ C_0 $ is constant. On a general $n$-dimensional Riemannian manifold, the classical Bishop-Gromov inequality is a doubling property (for balls whose 
radius lies in 
$] 0 , +\infty[$) whose doubling amplitude depends on the radius of the balls and on a lower bound $\text{Ric}_{\rm min}$ of the Ricci curvature: if one 
is not aiming for a sharp inequality, when the Ricci curvature is not supposed to be nonnegative,  it can be rewritten:
$$\dfrac{\Vol \, B (x , 2R )}{\Vol \, B (x , R ) } \le 2^n \,\exp \left(\sqrt{(n-1) |\text{Ric}_{\rm min}|}\, R\right) \ .$$
For this reason, in the sequel, for any nonnegative measure $\mu$, a doubling property of the kind 
$\dfrac{\mu \big(B (x , 2R )\big)}{\mu \big(\, B (x , R )\big) } \le C_1 \, e^{C_2 \,H\, R}$, where $H$ is a parameter which replaces the lower bound
of the curvature, will be called a (generalized) \lq \lq Bishop-Gromov inequality" (instead of \lq \lq doubling property"). The following Theorem proves such 
a Bishop-Gromov inequality in the
case of Gromov-hyperbolic metric spaces, where the hypothesis \lq \lq Ricci curvature bounded from below" is replaced by the much weaker (see subsection
\ref{comparaison}) hypothesis \lq \lq Entropy bounded from above".

\begin{theorem}\label{cocompact2}
Let $(X , d) $ be any $ \delta$-hyperbolic metric space, for every proper action 
(by isometries) of a group $\Gamma$ on $(X,d)$ such that the diameter of $\Gamma \backslash X$ and the
entropy of $(X , d) $ are respectively bounded by $D$ and $H$, then, for every $x \in X$
\begin{itemize}
\item[(i)] for every $\Gamma$-invariant measure $\mu$ on $X$, for every $r \ge \frac{5}{2} (7 D + 4 \delta)$,
$$ \dfrac{\mu \big( B_X (x , \frac{6}{5} \,r ) \big) }{ \mu \big( B_X (x , r)\big) } \le 1 + 2 \, e^{H D} \, e^{\frac{6}{5} \, H r}\ \ \ \ \ \ \text{ \rm and }\ \ \ \ \ \ 
\dfrac{\mu \big( B_X (x , 2 \,r ) \big) }{ \mu \big( B_X (x , r)\big) } \le  3^{4}\, e^{4 H D}\, e^{\frac{13}{2} H r}$$
$$\text{and }\ \ \  \  \ \ \ \ \ \  \ \ \   \forall R \ge r  \ \ \ \     \dfrac{\mu \left( \overline B_X \big( x , R \big)\right)}{\mu \big( B_X ( x , r )\big)} \le 3 e^{H D} \left( \frac{R}{r} \right)^{25/4}  \left( \frac{R}{r} \right)^{6 H D}  e^{6 H  (R - \frac{4}{5} r)} \, .$$

\item[(ii)] for every $r \ge 10\,(D + \delta )$, the counting measure $ \mu^\Gamma_x = \sum_{\gamma \in \Gamma} \delta_{\gamma x}$ of the orbit $\Gamma x$ verifies the inequalities :
$$ \dfrac{\mu_{x}^{\Gamma} \big( B_X (x , \frac{6}{5} \,r ) \big) }{ \mu_{x}^{\Gamma} \big( B_X (x , r)\big) } \le 1 + 2 e^{\frac{6}{5} \,H\, r}\ \ \ 
\text{ \rm , }\ \ \ \dfrac{ \#  \big( B_X (x , 2 r ) \cap \Gamma x\big)}{\#  \big( B_X (x , r ) \cap \Gamma x\big)} = \dfrac{\mu_{x}^{\Gamma} 
\big( B_X (x , 2 \,r ) \big) }{ \mu_{x}^{\Gamma} \big( B_X (x , r)\big) } \le  3^{4} e^{\frac{13}{2} H r} \,. $$
$$\text{and }\ \ \  \ \ \  \ \ \ \ \ \  \ \ \   \forall R \ge r  \ \ \ \     \dfrac{\mu_{x}^{\Gamma} \left(\overline B_X \big( x , R \big)\right)}{\mu_{x}^{\Gamma} \big( B_X ( x , r )\big)} 
< 3 \left( \frac{R}{r} \right)^{25/4} e^{6 H (R- \frac{4}{5} r)} \, .$$
\end{itemize}
\end{theorem}

If one wants to revisit this Bishop-Gromov-like inequality (i) (resp. (ii)) in terms of doubling (see Definitions \ref{doublefaible}), the second inequality of (i)
(resp. of (ii)) says that the measure $\mu$ (resp. the counting measure of any orbit of the action of $\Gamma$) satisfies the $ C_R$-doubling property, where $ C_R := 3^{4} \, e^{4 H D}\, e^{\frac{13}{2} H R} $ (resp. where $ C_R := 3^{4} e^{\frac{13}{2} H R}$) 
for all the balls of radius $r \in \left[\frac{5}{2} (7 D + 8 \delta) , R\right]$ (resp. of radius $r \in \left[10\,(D+ \delta ) , R\right]$).

\smallskip
The first step of the proof of Theorem \ref{cocompact2} is the following

\begin{lemma}\label{cocompact3}
On any $\delta$-hyperbolic space $(X,d )$, for every $R, \, R' \in \,] 0 , +\infty[$, and for any pair of points $x , \, y$ 
such that $d(x,y) < R+R'$, there exists a point  $y' \in X $ such that
$$B_X(x, R ) \cap B_X(y,R') \subset B_X(y', r)\ , \  \ \ \text{where}\ \ \ \  r = \Min \Big( R\,, \, R'\, ,\, 
\frac{1}{2} \,\big( R+R' - d(x,y) \big)+ \delta \Big) \ .$$
\end{lemma}

\begin{proof}
As $B_X(x, R ) \cap B_X(y,R')  = \emptyset$ when $d(x,y) \ge R+R'$, we shall only study the case 
where $d(x,y) < R+R'$. Even if this means exchanging the names of the points and of the radii, we
may suppose that $ R' \le R$.

\smallskip
-- If $d(x,y) \le  R - R' + 2 \, \delta$, the lemma is trivially verified when 
choosing $ y' = y$, because $ r = \Min \Big( R\,, \, R'\, ,\, \frac{1}{2} \,\big( R+R' - d(x,y) \big)+ \delta \Big) $ 
is then equal to $R'$ and $B_X(x, R ) \cap B_X(y,  R') \subset  B_X(y,  R') = B_X(y', r)$.

\smallskip
-- If $ R- R' + 2 \delta\le d(x,y) < R+R' $, let us denote by $c = [ x , y ]$ any geodesic such that $c(0) = x $
and $ c ( d(x,y) ) = y$. Let $y' := c \big( \frac{1}{2} \, (R - R' + d(x,y) )\big)$. For any point $z \in B_X(x, R ) 
\cap B_X(y,R') $, Lemma \ref{rectangletriangle} gives
$$ d(y' , z ) + d(x,y)   \le  \Max \left[ d(x,z) + d( y , y' ) \,,\,d(y,z) + d( x , y' ) \right] +  \delta $$
$$< \Max \left[ R+ d(x,y) - \frac{1}{2} \, \big(R - R' + d(x,y) \big)\ ,\  R' +  
\frac{1}{2} \, \big(R - R' + d(x,y) \big)\right] +  \delta\, ;$$
this implies that $d(y' ,z ) <  \frac{1}{2} \, \big(R + R' - d(x,y) \big) +  \delta = \Min \Big( R\,, \, R'\, ,\, 
\frac{1}{2} \,\big( R+R' - d(x,y) \big)+ \delta \Big) $.
\end{proof}

The second step of the proof of Theorem \ref{cocompact2} is the
\begin{lemma}\label{righthand}
Under the hypotheses of Theorem \ref{cocompact2}, for every $\Gamma$-invariant measure $\mu$ on $X$, for every $R' , R$ such that
$ 4 \delta \le R' \le R$, we have 
$$\dfrac{\int_{B_X (x , R) \setminus B_X(x,R- \frac{1}{2}R')} \mu \big( B_X(y, R') \big)\  d\mu (y)}{\mu \Big( B_X \big( x\,,\, \frac{3}{4} \, R'+ \delta + 
D\big)\Big) } \le \mu \Big( B_X (x , R+R') \setminus B_X(x,R- \frac{1}{2}R')\Big)$$
\end{lemma}

\begin{proof} In this proof we shall write $B (x,r)$ instead of $B_X (x,r)$ for sake of simplicity.
Equation \eqref{Fubini} gives
$$ \int_{B(x,R)}\ \mu \big( B(y, R') \big)\  
d\mu (y) =  \int_{X}\ \mu \big( B(y, R') \cap B(x, R)\big)\  d\mu (y)
\le \int_{B(x,R- \frac{1}{2}R')}\ \mu \big( B(y, R') \big)\  d\mu (y) +  \ $$
\begin{equation}\label{cococompact}
 + \int_{B (x , R + R') \setminus B(x,R- \frac{1}{2}R')}\ \mu \big( B(y, R') \cap B(x, R)\big)\  d\mu (y) \ .
\end{equation}
For every $y \in B (x , R + R') \setminus B (x,R- \frac{1}{2}R')$, we have $ R- R' + 2 \delta\le d(x,y) < R+R' $ and $\Min \Big( R\,, \, R'\, ,\, 
\frac{1}{2} \,\big( R+R' - d(x,y) \big)+ \delta \Big) = \frac{1}{2} \,\big( R+R' - d(x,y) \big)+  \delta \le \frac{3}{4} \, R'+ \, \delta $; then, applying Lemma \ref{cocompact3}, there exists a point $y'$ such that
$$B (y, R') \cap B(x, R) \subset B \Big( y' \,,\, \frac{1}{2} \, \big(R + R' - d(x,y) \big) + \delta \Big)
\subset B \Big( y', \frac{3}{4} \, R'+ \delta \Big) \ ;$$
as there exists $\gamma' \in \Gamma$ such that $d(y', \gamma' x) \le D$, the triangle inequality gives:
\small
$$\mu \big( B(y, R') \cap B (x, R)\big) \le \mu \Big( B \big( y'\,,\, \frac{3}{4} \, R'+ \delta \big)\Big)
\le \mu \Big( B \big( \gamma' \, x\,,\, \frac{3}{4} \, R'+  \delta + D\big)\Big) =
\mu \Big( B \big( x\,,\, \frac{3}{4} \, R'+  \delta + D\big)\Big)$$
\normalsize
where the last equality follows from the invariance of the measure and of the distance under the action of $\Gamma$. 
Plugging this estimate in \eqref{cococompact}, we obtain
$$\int_{B (x , R) \setminus B(x,R- \frac{1}{2}R')} \mu \big( B(y, R')\big)\,  d\mu (y) \le \mu \Big( B \big( x , \frac{3}{4} \, R'+ \delta + D\big)\Big)
\int_{B (x , R + R') \setminus B(x,R- \frac{1}{2}R')} d\mu (y) $$
and this ends the proof.
\end{proof}

The third step of the proof of Theorem \ref{cocompact2} is the
\begin{lemma}\label{lefthand}
Under the hypotheses of Theorem \ref{cocompact2}, for every $x \in X$, one has:
\begin{itemize}
\item[(i)] for every $R' , R$ such that $ 12 (D + \delta)\le R' \le R$, we have 
$$\dfrac{ \mu_{x}^{\Gamma}\big(B_X(x, R')\ \big) }{ \mu_{x}^{\Gamma}\Big( B_X \big( x\,,\, \frac{5}{6} \, R'\big)\Big) }  \le
\dfrac{ \mu_{x}^{\Gamma}\Big( B_X (x , R+R') \Big) -  \mu_{x}^{\Gamma}\Big( B_X(x,R- \frac{R'}{2})\Big) }{\mu_{x}^{\Gamma}\Big( B_X (x , R) \Big) -  
\mu_{x}^{\Gamma}\Big( B_X(x,R- \frac{R'}{2})\Big) } $$

\item[(ii)] for every $R' , R$ such that $ 3 (7 D + 4 \delta)\le R' \le R - D$, we have 
$$\dfrac{\mu\big(B_X(x, R')\ \big) }{ \mu \Big( B_X \big( x\,,\, \frac{5}{6} \, R'\big)\Big) }  \le
\dfrac{ \mu \Big( B_X (x , R+R'+ D) \Big) -  \mu \Big( B_X(x,R- \frac{R' + D}{2})\Big) }{\mu \Big( B_X (x , R) \Big) -  
\mu\Big( B_X(x,R- \frac{R' + D}{2})\Big) }$$
\end{itemize}
\end{lemma}

\begin{proof} In this proof we shall write $B (x,r)$ instead of $B_X (x,r)$ for sake of simplicity.\\
The measure $\mu_{x}^{\Gamma}$ and the distance being $\Gamma$-invariant, we have $\mu_{x}^{\Gamma} \big( B (y, R')\big) = 
\mu_{x}^{\Gamma} \big( B (x, R')\big) $ for every $ y \in \Gamma \, x$; from this and from the fact that the support of $\mu_{x}^{\Gamma}$ is  $\Gamma \, x$, 
we deduce that
$$\mu_{x}^{\Gamma} \big( B(x, R')\big) \cdot \mu_{x}^{\Gamma} \big(B (x , R) \setminus B(x,R- \frac{1}{2}R') \big) =
\int_{B (x , R) \setminus B(x,R- \frac{1}{2}R')} \mu_{x}^{\Gamma} \big( B(x, R')\big)\,  d\mu_{x}^{\Gamma} (y) $$
$$ = \int_{B (x , R) \setminus B(x,R- \frac{1}{2}R')} \mu_{x}^{\Gamma} \big( B(y, R')\big)\,  d\mu_{x}^{\Gamma} (y) \le
\mu_{x}^{\Gamma}  \Big( B (x , R+R') \setminus B(x,R- \frac{1}{2}R')\Big) \cdot$$
$$\mu_{x}^{\Gamma}  \Big( B \big( x\,,\, \frac{3}{4} \, R'+ \delta + D\big)\Big) \le \mu_{x}^{\Gamma}  \Big( B (x , R+R') \setminus B(x,R- \frac{1}{2}R')\Big) \cdot \mu_{x}^{\Gamma}  \Big( B \big( x, \frac{5}{6} R'\big)\Big)\, ,$$
and this proves (i).

\smallskip
Let $R'' := R' + D$ and let $y$ be any point of $B_X (x , R) \setminus B_X(x,R- \frac{1}{2}R'')$; as $d(y , \Gamma x) \le D$, there exists some $\g \in \Gamma$ 
such that $\g \big(B(x , R'' -D) \big) = B(\g x , R'' -D) \subset B(y' , R'')$, and thus (as the measure $\mu$ is $\Gamma$-invariant), it yields
$ \mu  \big[B(x , R') \big] = \mu  \big[B(x , R'' -D) \big] =  \mu  \big(\g \big(B(x , R'' -D)\big) \big)\le \mu \big( B (y, R'') \big)$. From this and from Lemma \ref{righthand},
we get:
$$\mu \big( B(x, R')\big) \cdot \mu \big(B (x , R) \setminus B(x,R- \frac{1}{2}R'') \big) \le \int_{B (x , R) \setminus B(x,R- \frac{1}{2}R'')} \mu \big( B_X(y, R'') \big)\,
 d\mu (y)$$
$$ \le \mu \Big( B (x , R+R'') \setminus B (x,R- \frac{1}{2}R'')\Big) \cdot \mu \Big( B \big( x\,,\, \frac{3}{4} \, R''+  \delta + 
D\big)\Big)\, ;$$
as $\frac{3}{4} \, R''+  \delta + D = \frac{3}{4} \, R'+  \delta + \frac{7}{4} D \le \frac{5}{6} \, R'$ (because $\frac{1}{12} \, R' \ge \frac{7}{4} D + \delta$), 
it comes: 
$$\mu \big( B(x, R')\big) \cdot \mu \big(B (x , R) \setminus B(x,R- \frac{1}{2}\, R'') \big) \le \mu \Big( B (x , R+R'') \setminus B (x,R- \frac{1}{2}\, R'')\Big) \cdot 
\mu \Big( B \big( x ,  \frac{5}{6} \, R'\big)\Big)\, ,$$
and this proves (ii).
\end{proof}

\begin{proof}[End of the proof of Theorem \ref{cocompact2}] 
 In this proof we shall write $B (x,r)$ instead of $B_X (x,r)$ for sake of simplicity; when $\mu$ is any $\Gamma$-invariant measure (resp. when 
$\mu = \mu_{x}^{\Gamma}$), for every $R' \ge 3 (7 D + 4 \delta)$ (resp. for every $R' \ge 12 (D +  \delta)$), we shall define $R''$ as $R' + D$ (resp. $R''$ as $R'$). 
In both cases, the result of Lemma \ref{lefthand} then writes:
\begin{equation}\label{resume}
\dfrac{\mu\big(B (x, R')\ \big) }{ \mu \Big( B \big( x\,,\, \frac{5}{6} \, R'\big)\Big) }  \le
\dfrac{ \mu \Big( B (x , R+ R'') \Big) -  \mu \Big( B (x,R- \frac{R''}{2})\Big) }{\mu \Big( B (x , R) \Big) -  \mu\Big( B (x,R- \frac{R''}{2})\Big) }
\end{equation}
Let us put $C := \dfrac{ \mu \big(B(x, R')\ \big) }{ \mu \Big( B \big( x\,,\, \frac{5}{6} \, R'\big)\Big) } - 1$
and $ a_k =  \mu \big(B (x, k\,\frac{ R''}{2})\ \big)$, a consequence of \eqref{resume} (where we replace $R$ by $k\,\frac{ R''}{2}$) is that,
for every $ k \ge 2$, one has $ a_{k+2} -  a_{k-1}  \ge (C+ 1) \left(  a_{k} -  a_{k-1} \right)$, and thus
$ a_{k+2} -  a_{k}  \ge  C \left(  a_{k} -  a_{k-1} \right)$. 
Making the sum of this last inequality for all the integers $k \in [2 , n ]$, we obtain
$$ 2\, (a_{n+2} - a_1 ) \ge a_{n+2} - a_3 +  a_{n + 1} - a_2  \ge C\left(  a_{n} -  a_{1} \right) \ ,$$
from which comes that $ a_{2\, n}  \ge \left( \frac{C}{2}\right)^{n-1} \left(  a_{2} -  a_{1} \right) + a_1$ and, as $a_{2} -  a_{1} > 0$ 
(for the annulus $ B (x, k\, R'')\setminus B (x, k\,\frac{ R''}{2})$ contains a ball $B(z ,D)$ of radius $D$ and $\mu \big( B(z ,D)\big)$ cannot vanish because 
the measure of $X = \cup_{\g \in \Gamma} \g \big( B(z ,D)\big)$ is not trivial), a consequence of this is
$$H \ge \Ent (X,d) = \lim_{n \to +\infty} \left(\frac{1}{n R'' } \  \ln \left(\mu \big(B (x, n\, R'') \big)\right)\right)
=  \lim_{n \to +\infty} \left(\frac{1}{n R'' } \  \ln \left( a_{2 n} \right)\right) \ge \frac{1}{R'' } \ \ln \left( \frac{C}{2}\right) \ .$$
For every $R' \ge 3 (7 D + 4 \delta)$ (resp. for every $R' \ge 12 (D + \delta)$), this yields $\dfrac{ \mu \big(B (x, R')\ \big) }{ \mu \Big( B \big( x\,,\, \frac{5}{6} \, R'\big)\Big) } \le 1 + 2 \, e^{H R''}$; making $r := \frac{5}{6} \, R'$, it follows that 
$\dfrac{ \mu \big(B (x, \frac{6}{5} \,r)\ \big) }{ \mu \Big( B \big( x ,  r\big)\Big) } \le 1 + 2 \, e^{H R''}$ for every $r \ge 3\cdot \dfrac{5}{6} \, (7 D + 4 \delta)$ (resp. for every $R' \ge 12 \cdot \dfrac{5}{6} \,(D + \delta)$).\\
In the case of the measure $\mu^\Gamma_x$, as $R'' = R' = \frac{6}{5} \, r$, we get:
$$ \forall r \ge 10 (D + \delta) \ \ \ \ \ \ \dfrac{\mu^\Gamma_x  \big(B (x, \frac{6}{5} \, r)\ \big) }{ \mu^\Gamma_x  \Big( B \big( x , r\big)\Big) } 
\le 1 + 2 \, e^{\frac{6}{5} \, H r} \le 3 \, e^{\frac{6}{5} \, H r} \,.$$
In the case of any other $\Gamma$-invariant measure $\mu$, as $R'' = R' + D = \frac{6}{5} \, r + D$, we obtain:
$$\forall r \ge \frac{5}{2} (7 D + 4 \delta) \ \ \ \ \ \ \dfrac{\mu \big(B (x, \frac{6}{5} \, r)\ \big) }{ \mu \Big( B \big( x , r\big)\Big) } 
\le 1 + 2 \, e^{H D} \, e^{\frac{6}{5} \, H r} \le 3 \, e^{H D} \, e^{\frac{6}{5} \, H r} \  .$$
This proves the first inequalities in (i) and (ii).\\
By iteration, replacing in succession $r$ by $r\,,\, \frac{6}{5} \,r\, , \ldots ,\,\left(\frac{6}{5}\right)^3 r$ in these two last inequalities we finish by proving that
the measure $\mu^\Gamma_x$ verifies
$$ \forall r \ge 10 (D + \delta) \ \ \ \ \ \ \dfrac{ \mu^\Gamma_x \big(B (x, 2 r) \big) }{ \mu^\Gamma_x \big( B (x, r )\big) }  \le 3^4 \, e^{\left(\frac{6}{5} + \ldots +
\left(\frac{6}{5}\right)^4\right) H r} \le 3^4 \, e^{\frac{13}{2}\,H r} $$
and that any other $\Gamma$-invariant measure $\mu$ satisfies:
$$ \forall r \ge \frac{5}{2} (7 D + 4 \delta)  \ \ \ \ \ \ \dfrac{ \mu \big(B (x, 2 r) \big) }{ \mu \big( B (x, r )\big) }  \le 3^4 \, e^{4 H D}\, e^{\left(\frac{6}{5} + \ldots +
\left(\frac{6}{5}\right)^4\right) H r} \le 3^4 \, e^{4 H D}  \, e^{\frac{13}{2}\,H r} \ .$$
This proves the second inequalities in (i) and (ii).

\smallskip
Let us now prove the third inequality of (ii): given $r , R $ such that $10\,(D+ \delta ) \le r \le R$, choose $p \in \N$ such that 
$ \left(\frac{6}{5}\right)^{p} \le \frac{R}{r}< \left(\frac{6}{5}\right)^{p+1}$; if $p = 0$, then the third inequality of (ii) is a trivial consequence of the first 
inequality of (i) and 
of the fact that, in this case, $\overline B_X \big( x , R \big) \subset B_X \big(x , \frac{6}{5} r\big) $, let us now suppose 
that $p \ge 1$ and choose $\e> 0$ sufficiently small in order that $\frac{R + \e}{r}< \left(\frac{6}{5}\right)^{p+1}$, applying the first of the inequalities 
(ii) to each term of the following product, we get:
$$  \dfrac{\mu_{x}^{\Gamma} \left(\overline B_X \big( x , R\big)\right)}{\mu_{x}^{\Gamma} \big( B_X ( x , r )\big)} \le
\dfrac{\mu_{x}^{\Gamma} \left(B_X \big( x , R+ \e \big)\right)}{\mu_{x}^{\Gamma} \big( B_X ( x , r )\big)}= 
 \prod_{i = 0}^{p-1} \dfrac{\mu_{x}^{\Gamma} \left(B_X \big( x ,  \big(\frac{5}{6}\big)^{i} (R+ \e ) \big)\right)}{\mu_{x}^{\Gamma} \left(B_X \big( x ,  \big(\frac{5}{6}\big)^{i+1}  (R+ \e ) \big)\right)}\cdot \dfrac{\mu_{x}^{\Gamma} \left(B_X \big( x ,  \big(\frac{5}{6}\big)^{p}  (R+ \e ) \big)\right)}{\mu_{x}^{\Gamma} \left(B_X \big( x ,  r \big)\right)}$$
$$< 3 \exp \left( \frac{6}{5} H r \right)\cdot \prod_{i = 0}^{p-1} 3 \exp \left( \left(\frac{5}{6}\right)^{i} H  (R+ \e )\right) < 3 \cdot 3^{p} 
\exp \left(6 H  (R+ \e )\left(1 - \left(\frac{5}{6}\right)^{p}\right) + \frac{6}{5} H r \right) \ .$$
$$<3\cdot \left( \frac{6}{5}\right)^{p \frac{\ln 3}{\ln (6/5)}}  \exp \left(6 H  (R+ \e )\left(1 -  \frac{r}{ (R+ \e )}\right) + \frac{6}{5} H r \right) < 3 \left( \frac{R}{r} \right)^{25/4} e^{6 H  (R+ \e - \frac{4}{5} r)}\, .$$
Now, making $\e \f 0$, we get
$$  \dfrac{\mu_{x}^{\Gamma} \left( \overline B_X \big( x , R \big)\right)}{\mu_{x}^{\Gamma} \big( B_X ( x , r )\big)} \le 3 \left( \frac{R}{r} \right)^{25/4} e^{6 H (R- \frac{4}{5} r)}\, .$$

\smallskip
Let us now prove the third inequality of (i): given $r , R $ such that $\frac{5}{2} (7 D + 4 \delta)\le r \le R$, choose $p \in \N$ such that 
$ \left(\frac{6}{5}\right)^{p} \le \frac{R}{r}< \left(\frac{6}{5}\right)^{p+1}$; if $p = 0$, then (ii) is a trivial consequence of the first inequality of (i) and 
of the fact that, in this case, $\overline B_X \big( x , R \big) \subset B_X \big(x , \frac{6}{5} r\big) $, let us now suppose 
that $p \ge 1$ and choose $\e> 0$ sufficiently small in order that $\frac{R + \e}{r}< \left(\frac{6}{5}\right)^{p+1}$, applying the first of the inequalities 
(i) to each term of the following product, we get:
$$  \dfrac{\mu \left(\overline B_X \big( x , R\big)\right)}{\mu \big( B_X ( x , r )\big)} \le
\dfrac{\mu \left(B_X \big( x , R+ \e \big)\right)}{\mu \big( B_X ( x , r )\big)}= 
 \prod_{i = 0}^{p-1} \dfrac{\mu \left(B_X \big( x ,  \big(\frac{5}{6}\big)^{i} (R+ \e ) \big)\right)}{\mu \left(B_X \big( x ,  \big(\frac{5}{6}\big)^{i+1}  (R+ \e ) \big)\right)}\cdot \dfrac{\mu \left(B_X \big( x ,  \big(\frac{5}{6}\big)^{p}  (R+ \e ) \big)\right)}{\mu \left(B_X \big( x ,  r \big)\right)}$$
$$< 3 e^{H D} e^{\frac{6}{5} H r}\cdot \prod_{i = 0}^{p-1} 3 e^{H D} \exp \left( \left(\frac{5}{6}\right)^{i} H  (R+ \e )\right) < 3 \cdot 3^{p} e^{(p+1) H D}
\exp \left(6 H  (R+ \e )\left(1 - \left(\frac{5}{6}\right)^{p}\right) + \frac{6}{5} H r \right) \ .$$
$$<3 e^{H D} \,\left( \frac{6}{5}\right)^{p \frac{\ln 3}{\ln (6/5)}} \left( \frac{6}{5}\right)^{p \frac{H D}{\ln (6/5)}} \exp \left(6 H  (R+ \e )\left(1 -  \frac{r}{ (R+ \e )}\right) + \frac{6}{5} H r \right) $$
$$ < 3 e^{H D} \left( \frac{R}{r} \right)^{25/4}  \left( \frac{R}{r} \right)^{6 H D}  e^{6 H  (R+ \e - \frac{4}{5} r)}\, .$$
Now, making $\e \f 0$, we get
$$  \dfrac{\mu \left( \overline B_X \big( x , R \big)\right)}{\mu \big( B_X ( x , r )\big)} \le 3 e^{H D} \left( \frac{R}{r} \right)^{25/4}  \left( \frac{R}{r} \right)^{6 H D}  e^{6 H  (R - \frac{4}{5} r)}\, .$$
\end{proof}

Considering any given co-compact action of a group $\Gamma$ on a $\delta$-hyperbolic space $(X,d)$ and 
any cyclic subgroup $\Gamma'$ of $\Gamma$, a consequence of revisiting Theorem \ref{cocompact2} (ii) as a doubling property (see the few lines after 
Theorem \ref{cocompact2}) is a universal upper bound $N(R)$ for the number of points of the orbit $\Gamma'\,x$ 
which are contained in the ball $B_X(x,R)$ (see Lemma \ref{majorationorbitale}). 
The fact that the bound $N(R)$ grows exponentially with $R$ 
may look strange for cyclic groups, but the
 Lemma \ref{puissances} (iii) proves that, for intermediate values of $R$, $N(R)$ must depend exponentially on $R$,
and this will be a fundamental tool in the proofs of the quantitative 
Tits' alternative (Theorem \ref{Tits3}) and of the lower bound of the systole given by Theorem \ref{Cat2}.\\
Let us recall that, to each hyperbolic isometry $\g$ and to each point $x$ of a $\delta$-hyperbolic space $(X,d)$,
one associates its \lq \lq displacement radius at $x$", i. e. $ R_{\g} (x) := \Min_{k \in \N^*} d(x, \g^{k} x) $
(see Definition \ref{rayondeplact}).

\begin{lemma}\label{majorationorbitale}
Let $(X , d) $ be any $ \delta$-hyperbolic (non elementary) space, for any proper action 
(by isometries) of a group $\Gamma$ on $(X,d)$ such that the diameter of $\Gamma \backslash X$ 
is bounded above by $D$, for every torsion-free element $\g \in \Gamma^*$ and for every $x_0 \in X$ 
such that $ R_{\g} (x_0) \ge 20 ( D + \delta ) $, then

\begin{itemize}
\item[(i)]  for every $p \in \N^*$
$$\# \left\{ k \in \Z : d \left(x_0 , \g^k x_0\right) < 2^p\, R_{\g}  (x_0) \right\} \le 3^{12 p} \, e^{16\, \Ent(X,d) 
\left( 2^p - 1\right)  R_{\g} (x_0)} \ .$$

\item[(ii)] at any point  $x \in X$ and for every $R \ge R_{\g} (x)$
$$ \# \{ k \in \Z : d \left(x , \g^k x\right) < R \} \le 3^{12} \left( \dfrac{R + 2\, d(x_0 , x)}{R_{\g}  (x_0)}\right)^{12 \, \frac{\ln 3}{\ln 2} } e^{ \frac{65}{2}\,  \Ent(X,d) \, (R + 2\, d(x_0 , x) )} \  .$$
\end{itemize}
\end{lemma}

\begin{proof}
First notice that, by Lemma \ref{reductisom} (i), $\g$ being torsion-free, the action of $\langle \g \rangle$ on $X$ is faithful.
Let $R_0:=20(D+ \delta)$. For sake of simplicity, let us write $ H$ instead of $\Ent(X,d)$. Theorem \ref{cocompact2} (ii) proves that,
for every $ R \ge \frac{1}{2}\, R_0$,
$\,\dfrac{\mu_{x_0}^{\Gamma} \big( B_X (x_0 , 2 \,R ) \big) }{ \mu_{x_0}^{\Gamma} \big( B_X (x_0 , R )\big) }
\le  3^{4} \, e^{\frac{13}{2} \, H  \, R}\,   \le C \, e^{\,\alpha  R}\, $, where $C=3^4$ and $\alpha= \frac{13}{2}\,H$.

 \emph{Proof of (i)}. As the inequality (i) is trivially verified when $ p= 0$, we shall suppose that $ p \ge 1$.
For every $i \in \N$, one has $ 2^i R_{\g} (x_0) \ge R_0$, we may thus apply the proposition \ref{sousgroupes} (ii)
in the case where $\Gamma' $ is the cyclic subgroup generated by $\g$, which implies that, for every $i \in \N$,
$$\dfrac{ \# \{ k \in \Z : d \left(x_0 , \g^k x_0\right) < 2^{i + 1} R_{\g}  (x_0)\}}{ \# \{ k \in \Z : d \left(x_0 , \g^k x_0\right) < 2^{i} R_{\g}  (x_0) \}} \le  C^3 \, e^{\frac{19}{8}\,\alpha \, 2^{i} R_{\g}  (x_0)}\ .$$
Making the product for all the integers $i  \in [ 0 , p-1]$, we infer that
$$\# \left\{ k \in \Z : d \left(x_0 , \g^k x_0\right) < 2^p R_{\g} (x_0) \right\}\le
\dfrac{ \# \{ k \in \Z : d \left(x_0 , \g^k x_0\right) < 2^{p} R_{\g}  (x_0)\}}{ \# \{ k \in \Z : d \left(x_0 , \g^k x_0\right) 
<  R_{\g}  (x_0) \}} $$
\begin{equation}\label{borneorbite}
\le C^{3 p} \, e^{\frac{19}{8}\, \alpha \left( 2^p - 1\right)  R_{\g} (x_0)} = 3^{12 p} \, e^{\frac{247}{16}\, H
\left( 2^p - 1\right)  R_{\g} (x_0)} \ ,
\end{equation}
where the first inequality is deduced from the fact that
$ \left\{ k \in \Z : d \left(x_0 , \g^k x_0\right) < R_{\g} (x_0) \right\} = \{ 0\}$.

\smallskip
\emph{Proof of (ii)}. For every $ k \in \Z^* $, the triangle inequality and the invariance by $\g$ of the distance give
$ d \left(x_0 , \g^k x_0 \right) \le  d \left(x , \g^k x\right) + 2\, d(x_0 ,x)$, this has two consequences: the first one is
the inequality $ R + 2\, d(x_0 ,x) \ge R_{\g} (x) + 2\, d(x_0 ,x) \ge R_{\g} (x_0)$, the second one is the fact that
$ \{ k \in \Z : d \left(x , \g^k x\right) < R \}$ is included in $\{ k \in \Z : d \left(x_0 , \g^k x_0\right) < R + 2\, d(x_0 ,x) \}$; 
let $p$ be the integer such that $2^{p-1} R_{\g} (x_0) <  R + 2\, d(x_0 ,x) \le 2^{p} R_{\g}  (x_0)$, from what 
precede and from \eqref{borneorbite}, we deduce:
$$ \# \left\{ k \in \Z : d \left(x , \g^k x\right) < R \right\} \le \# \left\{ k \in \Z : d \left(x_0 , \g^k x_0\right) < 
2^p R_{\g} (x_0) \right\} 
\le \left(2^p\right)^{12 \frac{\ln 3}{\ln 2}} \, e^{\frac{247}{16}\, H
\left( 2^p - 1\right)  R_{\g} (x_0)} .$$
The inequality (ii) then comes from this and from the choice of $p$ which implies that $ 2^p < 2 \,
\dfrac{R + 2\, d(x_0 , x)}{R_{\g}  (x_0)}$.
\end{proof}

\subsection{A Tits alternative and lower bounds for the entropies}\label{Tits}

\small

\emph{To each isometry $\g$ and to each point $x$ of a $\delta$-hyperbolic space $(X,d)$, 
one associates its displacement 
radius at $x$, i. e. $ R_{\g} (x) := \Min_{k \in \N^*} d(x, \g^{k} x)$ (see Definition \ref{rayondeplact}), its minimal 
displacement $s(\g) := \inf_{x \in X} d(x, \g \,x)$ and its asymptotic displacement $ \ell(\g) := \lim_{n \to +\infty} 
d(x, \g^n x)/n$ (see Definitions \ref{deplacements}). Let us also recall that the Margulis domain $M_R(\g)$ is the set
$\{ x \in X  :  R_{\g} (x) \le R\}$  (see Definition \ref{Margudomain} and Remark \ref{Margudomain1}), that 
$\g^-$ and $\g^+$ denote the fixed points of $\g$, that $ {\cal G} (\g)$ is the set of the geodesics $c$ such that $c(-\infty) = \g^-$ and $c(+\infty) = \g^+$, that $ M(\g)$ is the union of the images of these geodesics, and that $M_{\rm{min}} (\g)$ 
is the set of the points where the function $ x \mapsto d(x, \g x)$ attains its minimum (see Definitions
\ref{faisceau}).}

\normalsize

\subsubsection{The alternative: entropy vs asymptotic displacement}

\begin{theorem}\label{Tits3}
For any proper action (by isometries) of any group $\Gamma$ on any (non elementary) $\delta$-hyperbolic space
 $(X,d )$, if the diameter of the quotient $\Gamma \backslash X$ is bounded above by $D$, let $K:=\Min \left(\dfrac{1}{\delta} \,, \,1 \right)$, then

\begin{itemize}
  \item[(i)] either $\Ent (X,d) > \dfrac{K}{750} $,
  \item[(ii)] or $\Ent (X,d) \le \dfrac{K}{750}$, and then
$ \ell(\g) > 3^{-34}\, R_0 \, e^{- \,\frac{4}{29}\, K R_0} $
for every torsion-free $\g \in \Gamma^*$ (where $R_0 := \Max  \left[ 20\, \left( D+ \delta\right)  \, , 
\, 720\, \delta  \right] $).
\end{itemize}
\end{theorem}

\begin{exa}\label{sharpness}
\emph{Let $X$ be the Cayley graph of the free group $\Gamma$ with $2$ generators $a, b$, endowed 
with the following length-structure~: any path in the  Cayley graph being a concatenation of edges of types $[g , ga^{\pm 1}]$ or 
$[g , gb^{\pm 1}]$, its length is computed by deciding that the length of the edges $[g , ga^{\pm 1}]$ is 
$\alpha$ and that the length of the edges $[g , gb^{\pm 1}]$ is $\beta$. Let us suppose that  $0< \alpha \le \beta$. The corresponding length-distance $d_{\alpha, \beta}(x,y)$ between $x$ and $y$ is the minimal length of the paths joining $x$ to $y$. As the Cayley graph is a tree, it is clear that $(X, d_{\alpha, \beta})$ is $0$-hyperbolic.\\
The action of $\Gamma$ on the Cayley graph $X$ is the canonical one, by left-translations.
The quotient $( \Gamma \backslash X , \bar d_{\alpha, \beta} )$ being a union of two circles, of respective 
lengths $\alpha$ and $\beta$, its diameter is $D = \frac{1}{2} (\alpha + \beta)$.\\
When $\alpha$ goes to $0$ and $\beta$ is fixed, then $\Ent (X, d_{\alpha, \beta})$ goes to $+\infty$ and we are in the case $(i)$ of Theorem \ref{Tits3}.\\
When $\alpha$ goes to $0$ and $\beta = \dfrac{1}{\alpha}$, then $\Ent (X, d_{\alpha, \beta})$ goes to $0$ and we are in the 
case $(ii)$ of Theorem \ref{Tits3}. Moreover, as $\ell(a) = \alpha$ goes to $0$, it is coherent that the lower bound
given by $(ii)$ should go to zero when $2D= \alpha + \dfrac{1}{\alpha}$ goes to infinity. This proves that any
lower bound of the $\ell (\g)$'s must depend on a geometric invariant such as the diameter of the quotient.
}
\end{exa}

The following Lemma, though trivial, is often used in the sequel, it is thus necessary for the sake of clarity and to avoid
repetitions. The aim is to prove that, when a group acts properly on a Gromov-hyperbolic space, most of the properties
of this action deduce from the analogous properties of its image by the representation into the isometry group of $ (X,d)$.

\begin{lemma}\label{reductisom}
On every Gromov-hyperbolic space $(X,d)$, every proper action of a group $\Gamma$ is via a representation
$\varrho : \Gamma \f  \text{Isom} (X,d)$ which enjoys the following properties:

\begin{itemize}
\item[(i)] $\Ker \varrho$ is a finite normal subgroup of $\Gamma$, consequently all its elements have torsion,

\item[(ii)] $\varrho (\Gamma)$ is a discrete subgroup of $\text{Isom} (X,d)$,

\item[(iii)] if the diameter of $\Gamma \backslash X$ is bounded, as $\varrho (\Gamma)\backslash X = 
\Gamma \backslash X$, for the quotient metric, 
$\varrho (\Gamma)\backslash X $ is compact with the same diameter as $ \Gamma \backslash X$,

\item[(iv)] when $\Gamma \backslash X$ is compact, a measure $\mu$ on $X$ is $\varrho (\Gamma)$-invariant 
if and only if it is $\Gamma$-invariant; as a 
consequence, its entropy is unmoved when the action of $\Gamma$ via $\varrho$ is replaced by the 
canonical action of $\varrho (\Gamma)$,

\item[(v)] for every $\g \in \Gamma^*$, $\varrho (\g)$ is torsion-free if and only if $\g$ is torsion-free;
$\g$ is hyperbolic if and only if $\varrho (\g)$ is an hyperbolic isometry,

\item[(vi)] for every $\g \in \Gamma^*$, $\ell \big ( \varrho (\g) \big) = \ell \big ( \g \big) $.

\item[(vii)] for every subgroup $G$ of $\Gamma$, if $\varrho (G)$ is virtually cyclic, then 
$G$ is virtually cyclic.

\end{itemize}
\end{lemma}

\begin{proof} 
(i) is an immediate consequence of the definition \ref{proprediscr0} (ii) of a proper action, which implies that
$ \{ \g \in \Gamma : d(x,\g x) = 0\}$ is a finite set. (ii) is directly deduced from Lemma \ref{discret1} and from the
fact that the canonical action of $\varrho (\Gamma)$ is faithful and proper.\\
If  $ \Gamma \backslash X$ has bounded diameter, it is compact by Lemma \ref{autofidele} (ii), properties
(iii) and (iv) then follow immediately. (vi) is derived from the fact that every $\g \in \Gamma$ acts as the isometry 
$ \varrho (\g) $. If $ \g^n = e  $ then $\varrho (\g)^n = \id_X$;
conversely $\varrho (\g)^n = \id_X$ if and only if $ \g^n \in \Ker \varrho $, consequently $\g^n$ (and thus $\g$)
has torsion by (i). As $\ell \big ( \varrho (\g) \big) > 0$ iff $ \ell ( \g ) > 0 $ by (vi), $\g$ is hyperbolic if and 
only if $\varrho (\g)$ is an hyperbolic isometry by Lemma \ref{ellpositive}; this proves (v).\\
Proof of (vii): If $\varrho (G)$ is finite, as $\Ker \varrho $ is a finite normal subgroup of $G$ by (i), then $G$ is finite and (vii) is proved in this case.
If $\varrho (G)$ is infinite virtually cyclic, there exists an infinite cyclic subgroup $Z$ with finite index in $\varrho (G)$ (i. e. $ \varrho (G)/Z$ is finite);
then the subgroup $\varrho^{-1} (Z)$ is a subgroup of finite index in $G$ (because $\varrho$ induces a one to one map between the quotients
$G/ \varrho^{-1} (Z)$ and $\varrho (G)/Z$). By restriction $\varrho : \varrho^{-1} (Z) \f Z$ is a homomorphism whose kernel is the finite normal
subgroup $\varrho^{-1} (Z) \cap \Ker \varrho = \Ker \varrho$, in other words, $\varrho^{-1} (Z) $ is a finite-by-cyclic group, and it is a classical result that
\lq \lq finite-by-cyclic" implies \lq \lq cyclic-by-finite", i. e. that $\varrho^{-1} (Z) $ is virtually cyclic\footnote{The proof of the fact that $\varrho^{-1} (Z) $ is virtually cyclic is direct: indeed, there exists a splitting, i. e. a map $ \sigma : Z \f \varrho^{-1} (Z)$ such that $\varrho \circ \sigma = \id_Z$ (you just have 
to map the unique generator $\tau$ of $Z$ on any $\tau' \in \varrho^{-1} (\{\tau\})$. For every $g \in \varrho^{-1} (Z)$, as 
$\big(g.\Ker \varrho \big) \cap \sigma (Z) = \{\sigma \big(\varrho (g)\big)\}$, then $g \in \sigma \big(\varrho (g)\big).\Ker \varrho $ and 
$\varrho^{-1} (Z) = \sigma (Z).\Ker \varrho = \Ker \varrho. \sigma (Z)$, thus the quotient $\varrho^{-1} (Z) / \sigma (Z)$is finite, for it is in one to one correspondence with the finite set $\Ker \varrho$, consequently $\varrho^{-1} (Z)$ is virtually cyclic.}.
As $G/\varrho^{-1} (Z)$ is finite, it follows that $G$ is virtually cyclic.
\end{proof}

\begin{lemma}\label{minordist3}
Under the assumptions and notations of Theorem \ref{Tits3}, for every torsion-free  $\g \in \Gamma^*$, for every
$x_0 \in X$ such that $ R_{\g} (x_0) \ge R_0 $,

\begin{itemize}
    \item[(i)] there exists a positive integer $k_0$ such that $R_{\g} (x_0) = d(x_0, \g^{k_0} x_0) $ and then
$$ \Ent (X,d) \le \dfrac{1}{750\, \delta} \implies  k_0\, \ell(\g) \ge  3^{-12} R_{\g}  (x_0)  \   e^{- 16\,  \Ent (X,d) \,\,R_{\g} (x_0)}\ ,$$

    \item[(ii)] for every $ z \in X$, if $\Ent (X,d) \le \dfrac{1}{750\, \delta}$, then
$$ \dfrac{\ell(\g)}{R_{\g}  (x_0)} > \dfrac{3^{-12}}{2} 
\left( \dfrac{5 \,d(x_0,z) + 3\, d(z,\g \, z)}{R_{\g}  (x_0)}\right)^{- 20 } e^{- \frac{65}{2}\,  \Ent (X,d) \, [5 \,d(x_0,z) + 3\, d(z,\g \, z)]}\ . $$
\end{itemize}
\end{lemma}

\begin{proof} Let $\varrho$ be the representation from $\Gamma$ to $\text{Isom} (X,d)$ associated to the action
under consideration, Lemma \ref{reductisom} proves that the action of $\varrho (\Gamma)$ on $(X,d)$ also
verifies the assumptions of Theorem \ref{Tits3}, in particular
$\varrho (\g)$ is torsion-free when $\g$ is torsion-free, $\ell \big ( \varrho (\g) \big) = \ell \big ( \g \big) $ and
$ R_{\varrho (\g)} (x)  = R_\gamma (x) $ for every $x \in X$; therefore, if the conclusions of Lemma \ref{minordist3}
are satisfied when the group is $\varrho (\Gamma)$, they are also satisfied when the group is $\Gamma$.
To prove Lemma \ref{minordist3}, it is thus
sufficient to prove it when $\Gamma$ is a subgroup of $\text{Isom} (X,d)$; this is what we shall suppose in the 
sequel of this proof.\\
From Lemma \ref{autofidele} (ii), $\Gamma \backslash X$ is then compact and $\g$ cannot be a parabolic isometry
by the proposition \ref{actioncocompacte} (ii); as $\g$ is torsion-free, it is not elliptic (by Remark \ref{kpointsfixes} (i))
and therefore is an hyperbolic isometry by Theorem \ref{ellparahyp}.\\
By Lemma \ref{MRnontout}, there exists points $x \in X$ such that $R_{\g} (x) \ge R_0$. As the action is proper,
there exists an integer $k_0 > 0$ such that $d(x_0, \g^{k_0} x_0) =  R_{\g} (x_0)$; as the lemma \ref{puissances} (iii) 
gives, for any $  n \in \N^*$,
\begin{equation}\label{majordome}
d(x_0 , \g^{k_0 n} x_0) \le R_{\g} (x_0) + (n-1)\, \ell(\g^{k_0}) + 4\ \delta\, \dfrac{\ln (n)}{\ln (2)} \  ,
\end{equation}
we get
\small
 $$\# \{ n \in \N^* : d(x_0 , \g^{k_0 n} x_0) < 2 \, R_{\g}  (x_0)\} \ge \# \left\{\, n \in \N^* : (n-1)\,k_0\, \ell(\g) \le 
\frac{ R_{\g}  (x_0)}{2} \  \text{ and } 4\,\delta \,\dfrac{\ln (n)}{\ln (2)} < \frac{ R_{\g}  (x_0)}{2} \, \right\}$$
\begin{equation}\label{majormin}
\ge \Min \left( 1 + \left[\dfrac{ R_{\g}  (x_0) }{2\,k_0\, \ell(\g)} \right] \ \ ;\ \  e^{\dfrac{ R_{\g}  (x_0) \ln 2}{8 \,\delta}} -1 \right)
>  \Min \left( \dfrac{ R_{\g}  (x_0) }{2\,k_0\, \ell(\g)} \,\ ; \,\  e^{\dfrac{ R_{\g}  (x_0)}{12 \,\delta}} \right) \ ,
\end{equation}
\normalsize
where the last inequality follows from the fact that, as $ R_{\g} (x_0) \ge R_0 \ge 720 \,\delta$, one has
$$ e^{ R_{\g}  (x_0) \ln 2 / (8 \,\delta)}  - e^{ R_{\g}  (x_0)/ (12 \,\delta)}\  = \ e^{ R_{\g}  (x_0)/ (12 \,\delta)}\Big( 
e^{\frac{ R_{\g}  (x_0)}{\delta} \left( \frac{\ln 2}{8} \,- \, \frac{1}{12} \right)} -1\Big) \ge e^{60}\left( e^{720 \,
\left( \frac{\ln 2}{8} \,- \, \frac{1}{12} \right)} -1\right) > 1 \ .$$
On the other hand, property (i) of Lemma \ref{majorationorbitale} gives
$$\# \left\{ k \in \Z : d \left(x_0 , \g^k x_0\right) < 2\, R_{\g} (x_0) \right\} \le 3^{12 } \, e^{16\,  \Ent (X,d) \,  R_{\g} (x_0)} \ .$$
Frow this last inequality and from \eqref{majormin} follows:
\begin{equation}\label{proprietei}
\Min \left( \frac{ R_{\g}  (x_0) }{k_0\, \ell(\g)} \,\ ; \,\  2\, \exp{\left(\dfrac{ R_{\g}  (x_0)}{12 \,\delta}\right)} \right) +1<
3^{12 } \, e^{16\, \Ent(X,d)  \,R_{\g} (x_0)}\ , 
\end{equation}
\emph{Proof of (i) :} If $ \Ent (X,d) \le \dfrac{1}{750 \, \delta}$, recalling the existence of a point $ x_0 \in X$
such that $R_{\g} (x_0) \ge R_0$ and that the inequalities $ \frac{R_{\g}  (x_0)}{24 \,\delta} 
\ge \frac{R_0}{24 \,\delta} \ge 30 > 12 \, \ln 3 $ are valid in this case, the inequality \eqref{proprietei} gives:
$$ \Min \left( \dfrac{ R_{\g}  (x_0) }{k_0\, \ell(\g)} \,\ ; \,\  2\, e^{\dfrac{ R_{\g}  (x_0)}{12 \,\delta}} \right) +1<
3^{12 } \, e^{16\, \Ent (X,d) R_{\g} (x_0)} <  3^{12 } \,e^{\dfrac{2 \,R_{\g} (x_0)}{93 \,\delta}} <  e^{\dfrac{R_{\g} (x_0)}{12 \,\delta}}\ . $$
This proves that $ \dfrac{R_{\g}  (x_0) }{k_0\, \ell(\g)} < 3^{12 } \, e^{16\, \Ent (X,d) R_{\g} (x_0)}  $ and ends
the proof of (i). 

\medskip

\emph{Proof of (ii) :} For the sake of simplicity, in this proof, we shall set $ H := \Ent(X,d)$, 
$ R' := d(x_0,z)$, $A := d(z, \g\, z) $ and $ R := 3(R' + A) = 3 \big( \,d(x_0,z) + d(z,\g \, z) \big)$.\\ 
By Lemma \ref{puissances} (iii), for every $k \in \N^*$, we have:
$$d(z, \g^{k} z) \le  d(z, \g z)+ (k - 1)\,\ell(\g) + 4\,\delta \  \dfrac{\ln (k)}{\ln (2)} \le  A + (k - 1)\,\ell(\g) + 
4\,\delta \  \dfrac{\ln (k)}{\ln (2)}\ ,\ \ \ \text{ and consequently:}$$
$$\# \{ k \in \N^* \, :\  d \left(z , \g^k z\right) < R \} \ge \# \left\{ k \in \N \, :\, k < \dfrac{R - 2\,R'-2\,A}{2 \, \ell(\g)} +1 \ 
\text{and } \dfrac{\ln (k)}{\ln (2)} \le \dfrac{R + 2\,R'}{8 \, \delta}\right\} \ ,$$
it follows that
$$\# \{ k \in \Z \, :\  d \left(z , \g^k z\right) < R \}  > \Min \left( \dfrac{R - 2\,R' -2\,A}{\ell(\g)} \ , \ 
2\cdot 2^{(R  + 2\,R')/ 8 \, \delta} - 1\right) \ .$$
By the definition of $R$, one has $R \ge d(z, \g z) \ge R_\gamma (z)$, and we are authorized to apply the inequality (ii) 
of Lemma \ref{majorationorbitale} which (related to the previous inequality, noticing that $2^{\frac{R  + 2\,R'}{8  
\, \delta}} > 1 $), provides:
\begin{equation}\label{majorentropie}
\Min \left( \dfrac{R - 2\,R' -2\,A}{\ell(\g)} \ , \ 2^{( R  + 2\,R')/ 8 \, \delta}\right) < 3^{12} \left( \dfrac{R + 2\,R'}{R_{\g}  (x_0)}\right)^{12 \, \frac{\ln 3}{\ln 2} } e^{ \frac{65}{2}\, H \, (R + 2\,R' )}\ . 
\end{equation}
By the triangle inequality $R_{\g} (x_0)  \le 2\,R' + A$ thus, by the definition of $R$, one has
$ 2 (R + 2\, R') \ge 10 \, R' + 6\, A > 5 \, R_{\g}(x_0)$, and consequently
$  \dfrac{R_{\g} (x_0)}{2 (R + 2\,R')} \ \ln \left(  \dfrac{2 (R + 2\,R')}{R_{\g} (x_0)} \right) < \dfrac{\ln 5}{5}$.
From this and from the fact that $ H \le \dfrac{1}{750\, \delta} $ and $R_{\g}  (x_0) \ge R_0 \ge 720 \delta$
by assumption, we infer the two inequalities:
$$ \dfrac{65}{128}\, \dfrac{R  + 2\,R'}{8 \, \delta} \,\ln 2\   > \dfrac{65}{2}\, H\, (R + 2\,R' )
\ \ \ ,\ \ \  \dfrac{32}{128}\, \dfrac{R  + 2\,R'}{8 \, \delta}\, \ln 2 > 12 \ \frac{\ln 3}{\ln 2}\ 
\ln \left(  \dfrac{2 (R + 2\,R')}{R_{\g}  (x_0)} \right) \ ,$$
which imply that
$$2^{( R  + 2\,R')/ 8 \, \delta} > 3^{12} \left( \dfrac{R + 2\,R'}{R_{\g}  (x_0)}\right)^{12 \, \frac{\ln 3}{\ln 2} } 
e^{ \frac{65}{2}\, H \, (R + 2\,R' )}  $$
and, putting this estimate in \eqref{majorentropie}, we get that
$$\dfrac{R - 2\,R' -2\,A}{\ell(\g)} < 3^{12} \left( \dfrac{R + 2\,R'}{R_{\g}  (x_0)}\right)^{12 \, \frac{\ln 3}{\ln 2} } 
e^{ \frac{65}{2}\, H \, (R + 2\,R' )}\  ;$$
as we have seen that $R_{\g} (x_0) \le 2\,R' + A$, this proves (ii) because it proves that
$$ \dfrac{\ell(\g)}{R_{\g}  (x_0)} \ge \dfrac{\ell(\g)}{2 (R - 2\,R' -2\,A)}> \dfrac{3^{-12}}{2} \left( \dfrac{R + 2\,R'}{R_{\g}  (x_0)}\right)^{-20 } e^{- \frac{65}{2}\, H \, (R + 2\,R' )} \ . $$
\end{proof}

\begin{proof}[Proof of Theorem \ref{Tits3}]
 Let $\varrho$ be the representation from $\Gamma$ to $\text{Isom} (X,d)$ associated to the action
under consideration, we have already seen (see the beginning of the proof of Lemma \ref{minordist3}) 
that, if the assumptions of Theorem \ref{Tits3} are verified when the group is $\Gamma$, they are also
verified when the group is $\varrho (\Gamma)$ and that, if the conclusions of Theorem \ref{Tits3} are 
valid when the group is $\varrho (\Gamma)$, they are also valid when the group is $\Gamma$.
In order to prove Theorem \ref{Tits3}, it is thus
sufficient to prove it when $\Gamma$ is a subgroup of $\text{Isom} (X,d)$; this is what we shall suppose in the 
sequel of this proof.\\
$\Gamma \backslash X$ is then compact by Lemma \ref{autofidele} (ii). For the sake of simplicity, let 
$ H := \Ent (X,d)$.\\ 
If $ H > \dfrac{1}{750} \,\,\Min \left(\dfrac{1}{\delta} \,, \,1 \right)$ or if $\ell(\g) > \Max \left(\dfrac{\delta}{2500} \,,\, 10^{-5} D\right)  >  \dfrac{3^{-12}}{5}\ R_0 $ for every torsion-free element $\g \in \Gamma^*$, then the theorem
\ref{Tits3} is trivially proved, this is the reason why, from now on, (arguing by contradiction) we shall suppose that
$ H \le \dfrac{1}{750} \,\,\Min \left(\dfrac{1}{\delta} \,, \,1 \right)$ and that there exist torsion-free elements
$\g \in \Gamma^*$ such that $\ell(\g) \le  \Max \left(\dfrac{\delta}{2500} \,,\, 10^{-5} D\right) $; let us fix any
of these elements, denoted by $\g$. By Proposition  \ref{actioncocompacte} (ii), $\g$ cannot be parabolic 
and, as $\g$ is torsion-free, it is not elliptic by Remark \ref{kpointsfixes} (i), thus $\g$ is an hyperbolic isometry (by Theorem \ref{ellparahyp})
satisfying $\ell (\g) > 0$ (by Lemma \ref{ellpositive}).\\
As $R_0 := \Max \left[ 20\, ( D+ \delta) \, , \, 720\, \delta  \right]$,
the lemma \ref{quasigeod} (i) and the upper bound of $\ell(\g) $ which has just been assumed give
\begin{equation}\label{majores}
s(\g) := \inf_{x \in X}d(x, \g x) \le  \ell(\g) + \delta \le \delta + \Max \left(  \dfrac{\delta}{2500}\,,\, 10^{-5} D\right) 
< \dfrac{R_0}{719} \ ,
\end{equation}
Applying Lemmas \ref{MRferme} (ii) and \ref{MRnontout} to the Margulis domain $ M_{R_0} (\g)$, we know that
$ M_{R_0} (\g) \ne \emptyset $ and $X \setminus M_{R_0} (\g) \ne \emptyset $, a consequence (using the intermediate
value Theorem) is the existence of some point $x_1 \in X$ and of some $k_0 \in \N^*$ such that 
$ R_{\g} (x_1) := \Min_{k \in \N^*} d(x_1, \g^k x_1) = R_0$ and $ d(x_1, \g^{k_0} x_1) = R_0$; Lemma
\ref{minordist3} (i) then implies that
\begin{equation}\label{Hmajore}
 k_0\, \ell(\g) \ge \varepsilon'_0 \ \ \ \ \text{where}\quad \varepsilon'_0 := 3^{-12 }\, R_0 \, 
e^{- 16  H\, R_0} \ .
\end{equation}
Let us denote by $k_1$ the smallest integer such that $  k_1  k_0 \, \ell(\g) >  3 \delta $, \eqref{Hmajore} implies that
\begin{equation}\label{k1}
1 \le k_1 \le \left[ \dfrac{3\, \delta}{\e_0'}\right] + 1 \ \ \ \ \ \text{ and }\ \ \ \ \   (k_1 -1)\, k_0 \, \ell(\g) 
\le 3\, \delta \ .
\end{equation}
Let $ g := \g^{k_1 k_0}$, the inequality \eqref{majordome} guarantees that
\begin{equation}\label{puissancek0}
d(x_1 , g\, x_1) \le  R_0+ (k_1 - 1)\,k_0 \,\ell(\g) + 4\,\delta \  \dfrac{\ln (k_1)}{\ln (2)} \le
 \dfrac{9\,R_0}{8 }+ (k_1 - 1)\,k_0 \,\ell(\g) + 45 \, \delta \ ,
\end{equation}
where the last inequality is deduced from the fact that, by a direct computation using \eqref{k1} and the inequalities 
$R_0 \ge 720 \delta$ and $\ln (1+ t) \le \ln t + 1/t$ when $t > 0$, one has:
$$\ln (k_1) \le \ln \left(1 + 3^{13}\, \frac{\delta}{R_0}\,e^{16\, H\, R_0} \right) \le
\ln \left(1 + 3^{13}\, \frac{\delta}{R_0}\,\, e^{\dfrac{2 \,R_0}{93 \,\delta}} \right) $$
$$\le  \ln \left(3^{13}\, \frac{\delta}{R_0}\, e^{\dfrac{2 \,R_0}{93 \,\delta}} \right) + 
3^{-13}\, \frac{R_0}{\delta} \le 13\, \ln 3 - \ln 720 + 
\dfrac{R_0}{\delta}  \left( \frac{2}{93} + 3^{-13}\right)\le \dfrac{\ln 2}{4 \delta} \left( 45 \delta  + 
\dfrac{R_0}{8} \right).$$
As $\g$ is hyperbolic, the action of $\g$ on $X \cup \partial X$ has exactly two fixed points $\g^-\,, \g^+ \in \partial X$;
as $M(g) =M(\g) $ (because $\g^-$ and $\g^+$ are also the fixed points of $g$), and as $\ell(g) =  k_1\, k_0 \, \ell(\g) 
> 3\,\delta$, Lemmas \ref{distgeod} and \eqref{puissancek0} give:
\begin{equation}\label{distCg}
d(x_1, M(\g)) = d(x_1,M(g)) \le \frac{1}{2} \big(d(x_1, g \,x_1) - \ell(g)\big) + 3\,\delta \le
\frac{1}{2} \left( \dfrac{9\, R_0}{8 } - k_0 \, \ell(\g) \right) + \frac{51}{2}\, \delta \ .
\end{equation}
For any $\e > 0$, let us fix any geodesic $c \in {\cal G}(\g)$ such that the distance from $x_1 $ to the image 
of $c$ is smaller than
$d(x_1, M(\g)) + \e$, let us denote by $c(t_k)$ a projection of $ g^k x_1$ onto the geodesic line $c$; for every 
$ k \in \Z$ such that $d\big(c(t_k) \,,\,c(t_{k+1})\big) > 3\,\delta$, a consequence of the inequality 
\eqref{puissancek0}, of the lemma \ref{ecartement}, and of the fact that $g^{- k}\circ c$ and $g^{- (k+1)}\circ c$ 
are both geodesics of ${\cal G}(\g)$ is that
$$ \dfrac{9 \,R_0}{8 } + (k_1 - 1)\,k_0 \,\ell(\g) + 51\,\delta \ge d(x_1 , g\, x_1) + 6\,\delta =
d\big( g^k x_1 , g^{k+1}x_1 \big) + 6\,\delta  $$
$$\ge d\big( g^k x_1 , c(t_k)\big) + d\big(c(t_k) \,,\,c(t_{k+1})\big) + d\big( c(t_{k+1}) , g^{k+1}x_1 \big) \ge
2\, d(x_1, M(\g)) +  d\big(c(t_k) \,,\,c(t_{k+1})\big) \ .$$
From this, from \eqref{distCg} and from the definition of $  k_1$, we infer that, for every $k \in \Z$,
$$d\big(c(t_k) \,,\,c(t_{k+1})\big) \le \Max \left( \dfrac{9 \,R_0}{8 } + (k_1 - 1)\,k_0 \,\ell(\g) + 51\,\delta 
- 2\, d(x_1, M(\g)) \ , \ 3 \delta \right) =$$
$$ =  \dfrac{9 \,R_0}{8 } + (k_1 - 1)\,k_0 \,\ell(\g) + 51\,\delta - 2\, d(x_1, M(\g)) \le  \dfrac{9 \,R_0}{8 } +
 54\,\delta - 2\, d(x_1, M(\g)) \ ,$$
the last inequality being deduced from \eqref{k1}. This implies:
\begin{equation}\label{distc}
\forall t \in \R \ \ \exists k \in \Z \ \ \text{ such that } \  \  d\big(c(t) \,,\,c(t_{k})\big) \le \frac{9 }{16 } \,R_0  + 27\,\delta -  d(x_1, M(\g))\ .
\end{equation}
Let us denote by $c(t'_k )$ a projection of $g^k \circ c (t_0)$ onto the geodesic line $c$, as the geodesics
$c$ and $ g^k \circ \, c$ both admit $\gamma^-$ and $\gamma^+$ as endpoints, the proposition \ref{geodasympt} (i)
implies that $ d\big(g^k \circ c (t_0) ,  c(t'_k )\big) \le 2\, \delta$ and thus that
$$ \forall k \in \Z \ \ \ \ d(g^k  x_1 , c(t_k)\big) \le d(g^k  x_1 , c(t_k')\big) \le d(g^k  x_1 , g^k \circ c (t_0)\big) + 
d\big(g^k \circ c (t_0) ,  c(t'_k )\big) $$
$$\le d(x_1 ,  c (t_0)\big) + 2\, \delta   \le d(x_1, M(\g)) + 2\, \delta+ \e \ .$$
From this and from \eqref{distc}, we deduce that, for every $t\in \R$,
\begin{equation}\label{ecartgeo}
\min_{k \in \Z}  d\big(c(t) \,,\, g^k  x_1\big) \le \min_{k \in \Z} d\big(c(t) , c(t_k)\big) + \max_{k \in \Z}  d\big(c(t_k) \,,\, g^k  x_1\big)\le \frac{9 }{16 } \,R_0 + 29\,\delta + \e <  \dfrac{8}{13}\,R_0  + \e \ .
\end{equation}
Let us now notice that, as $\g$ is  hyperbolic, $M_{\rm{min}} (\g)$ is closed and non empty set by Lemma \ref{tube} (iv). Let us fix any point $y \in M_{\rm{min}} (\g)$ and projections $ c(t)$ and $ c(t')$ of $ y$ and $\g y$ (respectively)
onto the geodesic line $c$; applying \eqref{ecartgeo}, let us fix $k_2 \in \Z$ such that $d\big(c(t) \,,\, g^{k_2} \, x_1\big) 
\le \frac{9 }{16 } \,R_0 + 29\,\delta + \e$; we have only two possible cases:

\begin{itemize}

\smallskip
\item {\bf Case 1 : } If $ d (c(t) ,  c(t')  ) > 3 \,\delta$, the triangle inequality and the lemma \ref{ecartement} give:
$$ d(y, c(t)) + 3\delta + d(y , c(t')) -d(y, \g y) < d(y, c(t)) + d(c(t) , c(t')) + d(\g y , c(t'))  \le d(y, \g y) +  6 \, \delta\ ;$$
noticing that $ d(y , c(t')) \ge d(y, c(t) ) + d(c(t) , c(t') ) - 2 \, \delta > d(y, c(t) ) +   \delta$ (by Lemma \ref{projection})
and that $d(y,\g y) = s(\g) < \frac{R_0}{719}  $ (by \eqref{majores}) and plugging these two estimates into the previous
inequality, we obtain that $d(y,c(t)) < s(\g) + \delta < \frac{R_0}{719} +  \delta$.
From this, from the definition of $k_2$ and from the fact that  $29\,\delta + \delta \le \frac{R_0}{24}$ (by definition),
we deduce, when $\e$ is sufficiently small, that
$$ d( x_1 , g^{- k_2}\, y) = d( g^{k_2} \,x_1 , y) \le d( g^{k_2} \,x_1 , c(t)) + d(c(t) , y) <\dfrac{25}{41} \, R_0 \ .$$
As $M_{\rm{min}} (\g)$ is invariant by $\g^k$ for any $k\in \Z$, it is invariant by $ g^{-k_2} $ and thus 
$ g^{-k_2}\, y \in M_{\rm{min}} (\g)$; a consequence of the last inequality is that $d \big(x_1 ,  M_{\rm{min}} (\g) \big) < 
\dfrac{25}{41} \, R_0$ and thus
\begin{equation}\label{firststep}
\exists x \in M_{\rm{min}} (\g) \ \ \ \text{ such that }\ \ \ d( x_1 , x) < \dfrac{25}{41} \, R_0 
\end{equation}
Let us fix such a point $x \in M_{\rm{min}} (\g)$ and let us define $R:= 5 \,d(x_1, x ) + 3\, d(x,\g \,x) $; the 
estimates \eqref{firststep} and \eqref{majores} lead to
$$ R = 5 \,d(x_1, x ) + 3\, d(x,\g \,x) = 5 \,d(x_1,x) + 3\, s(\g) \le \left(\frac{125}{41} + \frac{3}{719}\right) R_0
<  \dfrac{77}{25}\, R_0 \ .$$
Noticing that $R_{\g} (x_1) = R_0$ and that (by assumption) $ H := \Ent (X,d) \le \dfrac{1}{750\, \delta}$, we may 
apply Lemma \ref{minordist3} (ii) to the points $x_1$ and $x$, which gives:
\begin{equation}\label{firstcase}
\dfrac{\ell(\g)}{R_0} > \dfrac{3^{-12}}{2} \left( \dfrac{R}{R_0}\right)^{- 20 } e^{- \frac{65}{2}\,  H\, R}
>  \dfrac{3^{-12}}{2}  \left( \dfrac{25}{77}\right)^{20}  e^{- 101\,  H R_0}\ .
\end{equation}

\medskip
\item {\bf Case 2: } If $ d (c(t) , c(t') ) \le 3 \,\delta$, let us denote by $ c\,(t")$ a projection of $\g \circ\, c\, (t)$
onto the geodesic line $c$, the triangle inequality then gives:
\begin{equation}\label{deplaceCg}
d \big( c(t) , \g \circ\, c\, (t) \big) \le | t- t'| + | t' - t" | + d \big( c(t") , \g \circ\, c\, (t) \big)
\end{equation}

\begin{itemize}
\item If $| t' - t" | \le 3 \, \delta$,  as the geodesics $c$ and $\g \circ \, c$ both admit $\gamma^-$ and $\gamma^+$ as 
endpoints, the proposition \ref{geodasympt} (i) implies that $ d \big(  \g \circ\, c\, (t) , c(t") \big) \le 2 \, \delta$; from this,
from \eqref{deplaceCg}, and from the assumptions on  $ d (c(t) , c(t') ) $ and on $| t' - t" | $, we get that 
$d \big( c(t) , \g \circ\, c\, (t) \big) \le 8 \, \delta$.

\item If $| t' - t" | > 3 \, \delta$, the lemma \ref{ecartement} gives:
$$ d \big( y \,, c\, (t) \big) +  6 \,\delta = d \big( \g y \,, \g \circ\, c\, (t) \big) +  6 \,\delta \ge d \big( \g y \,, c\, (t') \big)
+ | t' - t" | + d \big( c(t") , \g \circ\, c\, (t) \big) $$
$$\ge d \big( y \,, c\, (t') \big) - d(y, \g y) + | t' - t" | + d \big( c(t") , \g \circ\, c\, (t) \big) $$ 
$$\ge d \big( y \,, c\, (t) \big) + | t- t'| - 2\delta - s(\g) + | t' - t" | + d \big( c(t") , \g \circ\, c\, (t) \big) \ ,  $$
the last inequality being a consequence of Lemma \ref{projection}. Transferring this last estimate in \eqref{deplaceCg},
we obtain $d \big( c(t) , \g \circ\, c\, (t) \big) \le 8 \, \delta + s (\g) < \dfrac{R_0}{719} + 8 \, \delta$.
\end{itemize}
In case 2, we therefore always have
\begin{equation}\label{longitude}
d \big( c(t) , \g \circ\, c\, (t) \big) < \dfrac{R_0}{719} + 8 \, \delta <  \dfrac{R_0}{79} \ .
\end{equation}
Let us now redefine $R:= 5 \,d( g^{k_2} x_1 , c(t) ) + 3\, d(c(t) ,\g \circ\, c(t)) $; from the definition of $k_2$  (justified by
\eqref{ecartgeo}) and from \eqref{longitude}, if $\e$ is sufficiently small, we infer the estimate:
\begin{equation}\label{estimeR}
R := 5 \,d( g^{k_2}  x_1 , c(t) ) + 3\, d(c(t) ,\g \circ\, c(t)) < \dfrac{40}{13}\, R_0+ \dfrac{3}{79} \, R_0
< \dfrac{78}{25}\, R_0 \ .
\end{equation}
As  $ g := \g^{k_1 k_0}$, for every $ p \in \Z^*$ we have $ d \big(\g^p (g^{k_2} \,x_1)\, ,\, g^{k_2} \,x_1 \big) =  
d \big(\g^p \,x_1 , x_1 \big)$, and consequently $  R_{\g}  (g^{k_2} \,x_1 ) = R_{\g} (x_1) = R_0$; as moreover 
$  H = \Ent (X,d) \le \dfrac{1}{750\, \delta}$ by assumption, we may apply Lemma \ref{minordist3} (ii) to the points
$ g^{k_2} \,x_1$ and $ c(t) $ and get:
\begin{equation}\label{secondcase}
\dfrac{\ell(\g)}{R_0} > \dfrac{3^{-12}}{2} \left( \dfrac{R }{R_0}\right)^{- 20 } e^{- \frac{65}{2} H \,R} 
> \dfrac{3^{-12}}{2}  \left( \dfrac{78}{25} \right)^{- 20}  e^{- 102 \,H R_0 } > 3^{-34} \, e^{- 102 \,H R_0}\ ,
\end{equation}
where the second inequality comes from \eqref{estimeR} and (as $R$ depends on $\e$) is valid because $\e$ has 
been chosen sufficiently small.
\end{itemize}

We can now summarize all the cases that we have considered in the whole of this proof :
if $ H \le \dfrac{1}{750} \,\,\Min \left(\dfrac{1}{\delta} \,, \,1 \right)$, either $\ell(\g) > 
\Max \left(\dfrac{\delta}{2500} \,,\, 10^{-5} D\right)  >  \dfrac{3^{-12}}{5}\ R_0 $ or
$$\dfrac{\ell(\g)}{R_0} > 3^{-34} \, e^{- 102 \,H R_0 } \ge 3^{-34}\, \exp \left( - \dfrac{4}{29}\cdot  
\min \left( \dfrac{1}{\delta}, 1 \right)\cdot  R_0\right) \ ,$$
this last inequality being a consequence of the union of \eqref{secondcase} and \eqref{firstcase}. This ends the proof. 
\end{proof}


\subsubsection{Explicit universal lower bounds for the exponential growth}\label{quantitative}

In the case where the quotient space $\Gamma \backslash X$ has bounded diameter, we have the:

\begin{prop}\label{minorentropie}
For every non elementary $\delta$-hyperbolic metric space $(X,d )$, for every proper action (by isometries)
of a group $\Gamma$ on $(X,d)$, if the quotient $\Gamma \backslash X$ has diameter bounded above
by $D$, then
$$\Ent (X,d) \ge  \dfrac{\ln 2}{27 \delta + 10 D }\,,$$
\end{prop}

On Riemannian manifolds, we noticed (in the subsection \ref{comparaison}) that the assumption
\emph{Ricci curvature bounded from below} is much stronger than the assumption \emph{entropy bounded from above};
hence Proposition \ref{minorentropie}, being an obstruction to the smallness of the entropy, can be viewed as
a generalised version of the classical obstructions for the Ricci curvature to be almost nonnegative.

Example \ref{sharpness} proves that, in Proposition \ref{minorentropie}, the assumption on the diameter is 
necessary~: in fact let us consider the free group $\Gamma$ with $2$ generators $a, b$ and its
Cayley graph $X$ endowed with the length distance $d_{\alpha, \beta}$ defined (in Example \ref{sharpness}) 
by deciding that the length of the edges $[g , ga^{\pm 1}]$ (resp. $[g , gb^{\pm 1}]$)  is $\alpha$  
(resp.$\beta$). If $\beta$ goes to $+\infty$ and if $\alpha$ is fixed, we have seen in Example \ref{sharpness} 
that $\Ent (X, d_{\alpha, \beta})$ goes to $0$ and that $\diam (\Gamma \backslash X )$ goes to $+\infty$, while 
$ (X, d_{\alpha, \beta})$ is $0$-hyperbolic.

The following lower bound of the Entropy is valid in the case where $\Gamma \backslash X$ is non compact and 
in the case where it is compact with unbounded diameter. However, Example \ref{sharpness} proves that 
one needs to replace the bound on the diameter by the bound on another geometric invariant, as can be seen in the
next proposition.

\begin{prop}\label{entropienoncpct1}
For every $\delta$-hyperbolic, non elementary, metric space $(X,d )$, for every proper action (by isometries)
of a group $\Gamma$ on $(X,d)$, if there exists $L,C > 0$ and two hyperbolic elements $a , b \in 
\Gamma^*$ such that
\begin{itemize}
  \item $\langle a , b \rangle$ is not virtually cyclic,
  \item $\ell(a) , \ell(b) \le L$,
  \item $d \big( M(a) , M(b) \big) \le C$,
\end{itemize}
then, for any $\Gamma$- invariant Borel measure $\mu$ on $X$, 
$\Ent (X,d, \mu) \ge \dfrac{\ln 2}{17 \delta + L+ C}$.
\end{prop}

\begin{corollary}\label{entropienoncpct2}
For every $\delta$-hyperbolic, non elementary, metric space $(X,d )$, for every proper action (by isometries)
of a group $\Gamma$ on $(X,d)$, if there exists $L, R > 0$ and an hyperbolic element 
$a \in \Gamma^*$ such that
\begin{itemize}
  \item $\ell(a) \le L $,
  \item there exists $x \in M(a)$ such that $\Gamma_R (x)$ is not virtually cyclic,
\end{itemize}
then, for any $\Gamma$- invariant Borel measure $\mu$ on $X$, 
$\Ent (X,d, \mu) \ge \dfrac{\ln 2}{17 \delta + L+ R}$.
\end{corollary}

Theorem \ref{tubelong} (ii) allows to re-interpret the parameter $R$ which appears in
the assumptions of Corollary \ref{entropienoncpct2} as an upper bound of the length of a 
\emph{thin} tubular neighbourhood (in $\Gamma \backslash X$) of the image of $M(a)$ by 
the quotient map $X \f \Gamma \backslash X$. Thus the geometric meaning of the assumptions of  
Corollary \ref{entropienoncpct2} is that there exists at least one of these thin tubes whose length is bounded 
from above.

\emph{The lower bound of $\Ent (X,d) $ given by Corollary \ref{entropienoncpct2} must 
depend on the infimum $R$ of the $r$'s such that $\Gamma_r (x)$ is not virtually cyclic and must 
go to $0$ when $R \f +\infty$.} 

In fact, revisiting the example \ref{sharpness}, let us consider the free group $\Gamma$ with $2$ generators 
$a, b$ and its Cayley graph $X$ endowed with the length distance $d_{\alpha, \beta}$ defined 
(in the example \ref{sharpness}) by deciding that the length of the edges $[g , g\,a^{\pm 1}]$ (resp. 
$[g , g\,b^{\pm 1}]$)  is $\alpha$  (resp.$\beta$). Let us recall that $ (X, d_{\alpha, \beta})$ is $0$-hyperbolic,
and moreover $\CAT (K)$ for any $K \le 0$; the $a$-invariant geodesic $c_a$ is the concatenation of the 
edges $ \ldots \cup [a^{-k-1}, a^{-k}] \cup [a^{-k}, a^{-k +1}] \ldots \cup [a^{k-1}, a^k] \cup [a^{k}, a^{k+1}] \cup \ldots$ (and the same for the 
$b$-invariant geodesic $c_b$), and these two invariant geodesics $c_a$ and 
$c_b$ both contain the identity element $ e $. We consequently have $\ell(a) = \alpha$, $\ell(b) = \beta$.
Notice that $ e \in M(a) = c_a$ and let us fix $x_0 = e$. If $\beta \to +\infty$ and if $\alpha$ is fixed then, 
on $ (X, d_{\alpha, \beta})$, the infimum $R$ of the $r$'s such that $\Gamma_r (x_0) $ is not virtually cyclic is 
equal to $\beta$: in fact $\Gamma_\beta (x_0) $ is equal to $\Gamma$ and is thus not virtually cyclic and 
$\Gamma_r (x_0) = \langle a \rangle $ is cyclic for every $  r \in [\alpha, \beta[$. 
It follows that, for this value of $ R = \beta$, all the assumptions of Corollary \ref{entropienoncpct2} 
are fulfilled on $ (X, d_{\alpha, \beta})$ (except the fact that $R$ is bounded), in spite of this the entropy
is not bounded from below because $\Ent (X, d_{\alpha, \beta}) \to 0$ when $R := \beta \to +\infty$.

The entropy of a group with respect to a given system of generators being defined in Subsection
\ref{entropies}, we obtain the

\begin{corollary}\label{minorentalg}
For every non elementary group $\Gamma$, endowed with a finite generator system $\Sigma$ with respect to which
$\Gamma$ is $\delta$-hyperbolic, we have
$$\Ent (\Gamma,\Sigma)\ge \dfrac{\ln 2}{27 \delta + 10 } \ .$$
\end{corollary}

Previous estimates in the case of groups acting on Cartan-Hadamard spaces were given in \cite{BCG} and \cite{BCG2}.
In \cite{BF}, E. Breuillard and K. Fujiwara independently obtained, under the same assumptions and when $\delta \ge 1$,
the inequality $\Ent (\Gamma,\Sigma)\ge \dfrac{\ln 2}{400150000\, \delta}$ (see \cite{BF}, Corollary 13.2). However, this estimate and our Corollary \ref{minorentalg} do not provide any lower bound of the algebraic 
entropy of the group $\Gamma$. Indeed, if $(\Sigma_i)_{i \in \N}$ is an entropy-minimizing sequence of generating systems, the corresponding sequence of 
hyperbolicity constants of $(\Gamma, \Sigma_i)$ may go to infinity. See Corollary \ref{entalg3} and Theorem \ref{entalg2} (iii) for lower bounds of the 
algebraic entropy.

\begin{proof}[Proof of Proposition \ref{entropienoncpct1}]
As $ \ell(a) , \ell(b) > 0$, let $N, N'$ be the smallest integers such that 
$$ (N-1)  \ell (a) \le 13  \delta < N \ell(a) ,\quad  (N'-1)  \ell (b) \le 13  \delta < N' \ell(b) ,$$
we have:
\begin{equation}\label{ellaN}
13  \delta < \ell(a^N) \le 13  \delta + \ell(a) \le 13  \delta + L  \ ,\quad 
13  \delta < \ell(b^{N'}) \le 13  \delta + \ell(b) \le 13  \delta + L\,.
\end{equation}
Let us denote by $F$ the subgroup generated by $ \{ a^N ,b^{N'}\}$, it is not virtually cyclic by the proposition
\ref{actionelementaire} (vi); applying the corollary \ref{granddeplacement}, we get that one of the two semi-groups generated by $ \{ a^N , b^{N'}\}$ 
or by $ \{ a^N ,  b^{-N'}\}$ is a free one, and thus that the entropy of $F$ with respect to the 
complete system of generators $ S := \{ a^N , a^{-N},b^{N'}, b^{-N'}\}$ is
at least $\ln 2$.\\
For any $ \e > 0 $, let $x_0 \in M(a) = M(a^{N})$ and $x_1 \in M(b)= M(b^{N'}) \big) $ be two points 
which satisfy 
$d( x_0 , x_1) < C + \e$; let us fix some geodesic $[x_0 , x_1]$
connecting $x_0$ and $x_1$ and the middle-point $m$ of this geodesic. 
As $d(x_0 , a^N x_0 ) \le  \ell (a^N ) + 4\delta \le 17 \delta + L $ and $d(x_1 , b^{N'} x_1 ) \le \ell (b^{N'}) + 4\delta 
\le  17 \delta + L$ by Lemma  \ref{quasigeodesic} (i) and the estimates \eqref{ellaN}, the triangle inequality gives:
$$d (m , a^{-N}  m) = d (m , a^N  m) \le d(m, x_0) +  d(x_0 , a^N x_0 ) + d(a^N x_0 , a^N  m) \le C + 17 \delta + L + 
\e \ ;$$
by a similar proof we get $ d (m ,  b^{-N'}  m) = d (m ,  b^{N'}  m) \le  C + 17 \delta + L + \e$. Using Lemma 
\ref{comparentropi}, we obtain, for any $\Gamma$- invariant Borel measure $\mu$ on $X$,
$$\Ent (X,d,\mu) \cdot (17 \delta + L + C + \e) \ge \Ent (X,d,\mu) \cdot \Max \left[  d (m ,  a^{\pm N}  m)  ;\, 
d (m ,  b^{\pm N'} m) \right]$$
$$ \ge \Ent (F , S) \ge \ln 2 $$ 
and, when $\e \to 0$, this proves Proposition \ref{entropienoncpct1}. 
\end{proof}

\begin{proof}[Proof of Corollary \ref{entropienoncpct2}]
Choosing $ x_0 \in M(a)$ such that 
$\Gamma_R (x_0) = \langle \Sigma_R (x_0) \rangle$ is not virtually cyclic, Corollary \ref{elementaryaction} (iii) then implies the existence of some
$b \in  \Sigma_R (x_0)$ such that $ \langle a , b \rangle$ is not virtually cyclic, and Corollary \ref{elementaryaction} (i) then guarantees that
$ \langle a , b a b^{-1} \rangle$ is not virtually cyclic.
Moreover one has $d \big( M(a) , M(b a b^{-1})\big) =d \left( M(a) , b\big( M(a) \big)\right) \le 
d(x_0, b(x_0) \le R$.
Because of this and of the fact that $\ell(b a b^{-1}) = \ell(a) \le L$, we may apply Proposition 
\ref{entropienoncpct1} (where we replace $C$ by $R$)
to the two hyperbolic isometries $a$ and $b a b^{-1}$, which concludes.
\end{proof}

Before proving Proposition \ref{minorentropie}, let us first prove the

\begin{lemma}\label{loxodromique}
For every Gromov-hyperbolic space $(X,d)$, every infinite, discrete and co-compact subgroup $\Gamma$ of the group of 
isometries of $(X,d)$ contains at least one hyperbolic element $\g$ such that there exists $x \in X$ satisfying $ d (x, \g  x) \le  8 D + 10 \delta $,
where $D$ is an upper bound of the diameter of $\Gamma \backslash X$. A consequence is that the ideal boundary of $(X,d)$ contains
at least $2$ points.
\end{lemma}

\begin{proof} The action of $\Gamma$ on $(X,d)$ being proper (by Proposition \ref{discret1}), every sequence $\left( \g_n\right)_{n \in \N}$ of
distinct elements of $\Gamma$ verifies $d(x, \g_n x) \f +\infty$ when $n \to +\infty$ (for every $x \in X$). This proves that $(X,d)$ is a
unbounded, geodesic and proper space, thus that its ideal boundary $\partial X$ is non empty. 
Let $x_0$ be any fixed origin, $\theta$ be any point of $\partial X$, and $c$ a geodesic ray such that $ c(0) = x_0$ and $ c(+\infty) = \theta$. 
Let us fix $\e$ such that $3 \e$ is the infimum of $d(x, g x) - (8 D + 10 \delta)$ when $g$ runs in the set $\{g : d(x, g x) > 8 D + 10 \delta\}$; as the action of 
$\Gamma$ is proper (by Proposition \ref{discret1}), this infimum is attained and thus $\e >0$; moreover, by the choice of $\e$, every $ \g \in \Gamma$ such that $ d(x , \g x) \le 8 D + 10 \delta + 2 \e$ satisfies automatically $ d(x , \g x) \le 8 D + 10 \delta$.
Let us define $R = 3D + 5 \delta + \e$ and denote by $x$ and $y$ (respectively) the points $c(R)$ and $c(2 R)$ of the geodesic ray $c$. 
As $ d(y, \Gamma  x)$ and $ d( x_0 , \Gamma  x)$ are bounded above by $D$, there exists $g , h \in \Gamma$  such that $d(x_0, gx) \le D$ and 
$d(y , hx) \le D$. The triangle inequality then implies that $  R- D \le d(x, g x) \le R + D$, $  R- D \le d(x, h x) \le R + D$ and $ 2R - 2D \le d(gx, h x) \le 2R + 2D$; 
we then have
$$ d\big(g x,h x\big) \ge 2 R - 2 D = R + D +  5\,\delta  + \e  \ge \Max \,[d(x, g x) \,;\, d(x, h x) ] + 5\,\delta + \e \ .$$
Lemma \ref{minorell} (which is a variation of the lemma 9.2.3 page 98 of \cite{CDP})
then implies that either $\ell (g) \ge 3\delta + \e> 0$, or $\ell (h) \ge 3\delta +\e > 0$, 
or $\ell (g h) = \ell (h g) \ge 2\, \e > 0$, which implies (by Lemma \ref{ellpositive}) that (at least) one of the three isometries 
$g$, $h$ or $g h$ is hyperbolic; hence there exists an hyperbolic isometry $\g$ (equal to $g$, $h$ or $g h$) such that 
$ d(x , \g x) \le 2R + 2D = 8 D + 10 \delta + 2 \e$, and then we have $ d(x , \g x) \le 8 D + 10 \delta$
as proved above. A consequence is that the ideal boundary contains the two fixed point of this hyperbolic isometry and this ends the proof
\end{proof}

\begin{proof}[Proof of Proposition \ref{minorentropie}] Lemma \ref{loxodromique} guarantees that one can fix an element $ a \in \Gamma^*$ 
such that $a$ acts as an hyperbolic isometry and that $ 0 < \ell(a) \le 8 D + 10 \delta $. Choose some point $x_0 \in M(a)$. A result of M. Gromov 
(\cite{Gr1} Proposition 3.22, whose proof is written for Riemannian manifolds, but is still valid on path-connected metric spaces\footnote{in fact, 
the path-connectedness implies that, for every 
$\g \in \Gamma$ and every $\e > 0$, there exists a finite set $ \{y_0, y_1 , \ldots , y_N\} \subset X$
such that $y_0 = y$, $y_N = \g y$ and $d(y_{i-1}, y_i) < \e $ for every $i \in \{1 , \ldots , N \}$.
Let us choose $\g_0, \g_1 , \ldots , \g_N \in \Gamma$ such that $ \g_0 = e$, $\g_N = \g $ and 
$d (y_i , \g_i y ) \le D$, we then get $\g = \sigma_1 \cdot \ldots  \cdot \sigma_N$, where 
$\sigma_i = \g_{i-1}^{-1}\cdot \g_i \in \Sigma_{2 D+ \e} (y)$ and the finiteness of $\Sigma_{3 D} (y)$ 
proves that $\Sigma_{2 D+ \e}   (y) = \Sigma_{2 D} (y)$ when $\e$ is sufficiently small.}) proves that 
$\Sigma_{2 D} (x_0) := \{ \sigma \in \Gamma^*  : d(x_0, \sigma x_0) \le 2 D\}$ is a complete system of 
generators of $\Gamma $, the properness of the action implying the finiteness of $\Sigma_{2 D} (x_0) $.
This proves the existence of some $ b \in \Sigma_{2 D} (x_0) $ (possibly elliptic !) such that $\langle a,b \rangle$ is 
not virtually cyclic (otherwise $\Gamma$ would be virtually cyclic by the corollary \ref{elementaryaction} (iii) and $(X,d)$ would then be elementary 
by Propositions \ref{actioncocompacte}  (iii) and \ref{actioncocompacte}  (iv)). By Corollary \ref{elementaryaction} (i) $ \langle  a , b a b^{-1})\rangle$  
is not virtually cyclic too; as moreover $ \ell(b a b^{-1}) = \ell(a) \le 8 D + 10 \delta $ and 
$d \big( M(a) , M(b a b^{-1})\big) = d \left( M(a) , b\big( M(a) \big) \right) \le d(x_0, b(x_0) \le 2D$, we may apply the proposition \ref{entropienoncpct1} 
(where $L$ and $C$ are replaced by $8 D + 10 \delta  $ and $2D$ respectively) to the two hyperbolic isometries
$a$ and $b a b^{-1}$, which ends the proof.
\end{proof}

\begin{remark}
\emph{Proposition \ref{minorentropie} could also be deduced from Theorem \ref{Tits3}, but the proof given
above is much more simple. This stresses the fact that Theorem \ref{Tits3} is
a more powerful result: in fact, if we suppose that the entropy of $(X,d)$
is smaller than $\dfrac{1}{750} \, \Min \left( \dfrac{1}{\delta}, 1\right)$, Proposition \ref{minorentropie} is a consequence of the fact that
there exists} a torsion-free $\g \in \Gamma^*$ \emph{such that $ \ell (\g) \le C(\delta, D)$ while the 
Theorem \ref{Tits3} implies that} every torsion-free $\g \in \Gamma^*$ \emph{verifies 
$ \ell (\g) \ge C'(\delta, D)$.}
\end{remark}

\begin{proof}[Proof of Corollary \ref{minorentalg}]
As seen in the sections \ref{notations} and \ref{entropies}, the $\delta$-hyperbolicity of the group $\Gamma$ 
(when endowed with the generator system $\Sigma$) means that its 
Cayley graph $ {\cal G}_{\Sigma} (\Gamma) $, endowed with the algebraic distance $ d_\Sigma$ (defined in 
the section \ref{notations}) is a $\delta$-hyperbolic metric space. The action (by left-translations) of $ \Gamma$ 
on $ {\cal G}_{\Sigma} (\Gamma) $ is isometric and proper, because the balls of $(G , d_\Sigma) $ are finite sets.
As $\Sigma$ is finite, $\left( {\cal G}_{\Sigma} (\Gamma) ,  d_\Sigma\right)$ is a proper space and the quotient
$ \Gamma \backslash {\cal G}_{\Sigma} (\Gamma)$ is a union of a finite number of circles whose diameter 
is bounded by $ D = 1$.
As $ \Gamma$ is non elementary, its Cayley graph is non elementary too. We can thus apply the  Proposition 
\ref{minorentropie}, which implies that
$$\Ent (\Gamma ,\Sigma) := \Ent \left({\cal G}_{\Sigma} (\Gamma) ,   d_\Sigma\right) \ge 
\dfrac{\ln 2}{10  D + 27 \delta }  =  \dfrac{\ln 2}{10 + 27 \delta } \ .$$
\end{proof}


\subsubsection{Implicit universal lower bounds for the exponential growth}

\small

\emph{The results of this subsection are in some sense stronger than the results of the previous subsection 
\ref{quantitative} because
they provide lower bounds for the algebraic entropy of any non virtually cyclic (eventually non cocompact) 
group acting properly on a $\delta$-hyperbolic space, in particular they bound from below the algebraic entropy of any
non virtually cyclic subgroup of a $\delta$-hyperbolic group (independently of its system of generators).
On the other hand, the results of this subsection are in some sense weaker: in fact, the lower bounds given in 
the previous subsection \ref{quantitative} are explicit while the bounds which will be established in the present 
subsection, though universal, cannot be specified, because they all depend on the non explicit function 
$ p \mapsto N(p)$ which appears in the theorem \ref{BGT} of E. Breuillard, B. Green and T. Tao. From this 
function, for every values of $ C_0 > 1$ and of $H, D, \delta > 0$, we define the two following universal constants:
\begin{equation}\label{N1N0}
N_1 (C_0) = 3 N \left( \left[ C_0^3\right] + 1\right) \ \ \ ;\ \ \ N_0 (\delta , H , D) := N \left(\left[3^{12} \, 
e^{ 490 \,H ( D + \delta )   }\right] + 1 \right) \ .
\end{equation}
Let us furthermore recall that, for every $A, B \subset G$, we denote by $A\cdot B$ the image of $A\times B$ by the 
map $(\g, g) \mapsto \g \cdot g$ and define (by induction) $S^k$ as $S^{k-1}\cdot S$.}

\normalsize

\begin{theorem}\label{entalg2}
Let $(X , d) $ be any $ \delta$-hyperbolic (non elementary) metric space, for every proper action 
(by isometries) of a group $\Gamma$ on $(X,d)$ such that the diameter of $\Gamma \backslash X$ and the
entropy of $(X , d) $ are respectively bounded by $D$ and $H$, then (denoting by $N_0$ the above universal 
constant  $N_0 (\delta, H, D)$)
\begin{itemize}
\item[(i)] for any finite symmetric subset $S$ of $\Gamma$ which generates a non virtually cyclic subgroup,
there exists at least one $ \g_0 \in S^{3 N_0}$ such that $\ell (\g_0) \ge \delta$ and there exists 
$\sigma \in S$ such that one of the two semi-groups generated by $ \{\g_0^{14} \, , \,  \sigma \g_0^{14} 
\sigma^{-1} \}$ or by $ \{\g_0^{14}\, , \, \sigma \g_0^{-14} \sigma^{-1}\}$ is free.

\item[(ii)] The algebraic entropy\footnote{Recall that, in the absence of additional specification, the algebraic entropy of a group $G$ is the 
infimum (with respect to every  finite 
system $S$ of generators of $G$) of $\Ent (G,S)$.} of any finitely generated and non virtually cyclic subgroup $\Gamma'$ of $\Gamma$ 
 is bounded from below by $\dfrac{\ln 2}{42\, N_0 + 2}$.

\item[(iii)]  The algebraic entropy of $\Gamma$ is bounded from below by $\dfrac{\ln 2}{42\, N_0 + 2}$.
\end{itemize}
\end{theorem}

Given a finitely generated group $\Gamma$, let us recall that, to any choice of a finite system $S$ of generators of $\Gamma$ corresponds a Cayley graph,
which is a metric space when endowed with the algebraic distance $d_S$. When $\Gamma$ is Gromov-hyperbolic, this metric space is Gromov-hyperbolic too 
and we shall denote by $\delta (\Gamma , S)$ its hyperbolicity constant\footnote{When $\Gamma$ is Gromov-hyperbolic, for every finite system $S$ of generators, $\delta (\Gamma , S) < +\infty$.}. We then have

\begin{corollary}\label{entalg3}
Let $\Gamma$ be a non virtually cyclic Gromov-hyperbolic group then, for every positive constant $M$, if there exists a finite system $S_0$ of generators of 
$\Gamma$ such that $\Ent (\Gamma, S_0) \, (\delta (\Gamma , S_0) + 1) \le M$, then the algebraic entropy of $\Gamma$ and of any finitely 
generated and non virtually cyclic subgroup $\Gamma'$ of $\Gamma$  is bounded from below by 
$\dfrac{\ln 2}{42\, N \left(\left[3^{12} \, e^{ 490 \,M }\right] + 1 \right) + 2}$ where $ p \mapsto N(p)$ is the function which appears in Theorem \ref{BGT}.
\end{corollary}

Note that E. Breuillard and K. Fujiwara also give a lower bound of the algebraic entropy of an hyperbolic group and of its non virtually 
cyclic subgroups, 
their estimate also depends on the existence of a generating system $S_0$ of $\Gamma$, however they suppose that $\delta (\Gamma , S_0)$ and the cardinality $\vert S_0\vert$ of this generating 
system are bounded above. They can then replace our use of the Bishop-Gromov inequality (Theorem \ref{cocompact2}) by the fact that 
$ \dfrac{\# B_{S_0} (e , N +1)}{\# B_{S_0} (e , N)} \le 2 \,\# (S_0) - 1$, which clearly implies that $\Ent (\Gamma, S_0) \le \ln \big(2\, \# (S_0) - 1 \big)$ 
(see \cite{BF}, Corollary 13.4).

\medskip
Theorem \ref{entalg2} is in fact a corollary of the following proposition, which is also valid in the non cocompact case:

\begin{prop}\label{entalg1}
Let $(X , d) $ be any $\delta$-hyperbolic (non elementary) metric space, for every proper action 
(by isometries) of a group $\Gamma$ on $(X,d)$ if, for every point $x \in X$, the counting measure $
\mu^{\Gamma}_x$ of the orbit $\Gamma x$ satisfies the $ C_0 $-doubling condition for all the balls 
(centered at $ x $) of radius 
$ r \in \left[ \frac{1}{2}\, r_0 \,,\, \frac{5}{4}\, r_0  \right]$  (where $r_0 > \frac{31}{2} \, \delta$ and $C_0 > 1$ are 
arbitrary constants) then (denoting by $N_1$ the above universal constant $ N_1 (C_0)$)
\begin{itemize}

\item[(i)] for any finite symmetric subset $S$ of $\Gamma$ which generates a non virtually nilpotent subgroup,
there exists at least one $ \g_0 \in S^{N_1}$ such that $\ell (\g_0) \ge \delta$ and there exists 
$\sigma \in S$ such that one of the two semi-groups generated by $ \{\g_0^{14} \, , \,  \sigma \g_0^{14} 
\sigma^{-1} \}$ or by $ \{\g_0^{14}\, , \, \sigma \g_0^{-14} \sigma^{-1}\}$ is free.

\item[(ii)] The algebraic entropy of any finitely generated and non virtually nilpotent subgroup $\Gamma'$ of $\Gamma$, 
with respect to any finite system of generators of $\Gamma'$, is bounded from below by $\dfrac{\ln 2}{14\, N_1 + 2}$.
\end{itemize}
\end{prop}

\begin{proof}
Let $N'_1 = \frac{1}{3}\, N_1 = N \left( \left[ C_0^3\right] + 1\right) $; let us denote by $G = \langle S\rangle$ 
the subgroup generated by $S$ and by $A$ the set $\{\g \in G : d(x, \g x) \le r_0\}$. By the proposition \ref{sousgroupes},
the $ C_0 $-doubling condition assumed for the counting measure $\mu^{\Gamma}_x$ of the orbit $\Gamma x$ implies
that the counting measure $\mu^{G}_x$ of the orbit of $G$ satisfies the condition
$$\dfrac{\# (A \cdot A)}{\# A} \le \dfrac{\mu_{x}^{G}  \big( \overline {B}_X( x , 2 r_0 ) \big)}{\mu_{x}^{G} \big( B_X (x , r_0 ) \big)} \le C_0^3 \ .$$
As $G := \langle S \rangle$ is not virtually nilpotent, the theorem \ref{BGT} then proves that $S^{N'_1}$ is not contained in $A$; there consequently 
exists some $\g \in S^{N'_1}$ such that $d(x, \g x) > r_0 > \frac{31}{2}\, \delta$.\\ 
If $\varrho$ is the representation $\Gamma \f \text{Isom} (X,d)$ associated to the action of $\Gamma$ on $ (X,d)$, there exists 
$\varrho (\g) \in \varrho (S^{N'_1}) = \varrho (S)^{N'_1}$ such that $d(x, \varrho (\g) x) > r_0 > \frac{31}{2}\, \delta$. As this is valid for every $x \in X$, we obtain 
the estimate: $L^* \big(\varrho (S)^{N'_1}\big) \ge  r_0 > \frac{31}{2}\, \delta$, and the theorem \ref{minorLetoile} (where we replace 
$S$ by $\varrho (S)^{N'_1}$) then implies that there exists $ g \in \varrho (S)^{3\, N'_1} =  \varrho (S^{N_1})  $ such that $\ell (g) \ge \delta$;
there thus exists $ \g_0 \in S^{N_1} $ such that $\varrho (\g_0) = g$ and thus (by Lemma \ref{reductisom} (vi)) $\ell (\g_0) = \ell (g) \ge \delta$.\\
Now, as $\langle S \rangle$ is not virtually cyclic, by Corollary \ref{elementaryaction} (iii), there exists $\sigma \in S$ such that the subgroup generated by 
$\{\g_0 , \sigma \}$ is not virtually cyclic and, by Corollary \ref{elementaryaction} (i), $ \{\g_0 \, , \,  \sigma \g_0 \sigma^{-1}\}$ generates a non 
virtually cyclic subgroup of $\Gamma$. Hence we may apply Corollary \ref{granddeplacement} to the pair $ \{\g_0 \, , \,  \sigma \g_0 \sigma^{-1}\}$,
which proves that one of the two semi-groups generated by $ \{\g_0^{14} \, , \,  \sigma \g_0^{14} \sigma^{-1}\}$ 
or by $ \{\g_0^{14}\, , \, \sigma \g_0^{-14} \sigma^{-1}\}$ is free. This proves (i).\\
This also implies that the (algebraic) entropy of the subgroup $G'$ generated by $ S' :=\{\g_0^{14} \, , \,  
\sigma \g_0^{14} \sigma^{-1}\}$ (with respect to the generator system $ S'$) is bounded from below by $\ln 2$.
Let us consider the metric space $G$, endowed with the algebraic distance $d_S$ associated to the generator system 
$S$, and the faithful action (by left translations) of $G'$ on $G$, as this action is proper and $d_S$-isometric, we can apply the lemma  \ref{comparentropi}, which gives
$$\Ent (G, S) := \Ent (G, d_S , \# ) \ge  \dfrac{ \Ent (G', S')}{\Max \left[d_S \left(e,\g_0^{14}\right)\,;\, 
d_S \left(e, \sigma \g_0^{14} \sigma^{-1} \right)\right]} \ge \dfrac{\ln 2}{14\, N_1 + 2} . $$
Let us now consider any finitely generated and non virtually cyclic subgroup $\Gamma'$ of $\Gamma$, and any finite 
system of generators $\Sigma$ of $\Gamma'$, let us denote by $S$ the symmetrized of $\Sigma$,
as $\Gamma' = \langle S \rangle$, we can apply the previous inequality
(replacing $G$ by $\Gamma'$), which ends the proof of the proposition \ref{entalg1}.
\end{proof}

\begin{proof}[End of the proof of Theorem \ref{entalg2}]
Theorem \ref{cocompact2} (ii) proves that,
for every point $x \in X$, the counting measure $\mu^{\Gamma}_x$ of the orbit $\Gamma x$ satisfies the 
$ C_0 $-doubling condition (with $C_0 := 3^{4} \, e^{ \frac{490}{3} \,H ( D + \delta )}$)
for all the balls (centered at $ x $) of radius 
$ r \in \left[ \frac{1}{2}\, R_0 \,,\, \frac{5}{4}\, R_0  \right]$  (where $R_0 = 20\, (D + \delta) > 
\frac{31}{2} \, \delta$). We may therefore apply Proposition \ref{entalg1}, while replacing $N_1$ by $3\, N_0$ because
$$N_1 = N_1 (C_0) = 3\, N \left( \left[ C_0^3\right] + 1\right) = 3 \, N \left(\left[3^{12} \, 
e^{ 490 \,H ( D + \delta )   }\right] + 1 \right) = 3\, N_0 \ .$$
This proves the properties (i) and (ii) of Theorem \ref{entalg2}. Moreover, if $\varrho$ is the representation $\Gamma \f \text{Isom} (X,d)$ associated 
to the action of $\Gamma$ on $ (X,d)$, $\varrho (\Gamma)$ is not virtually cyclic, otherwise $(X,d)$ would be elementary by Propositions 
\ref{actioncocompacte}  (iii) and (iv); this implies that $\Gamma$ is not virtually cyclic and thus that $(ii) \implies (iii)$.
\end{proof}

\begin{proof}[End of the proof of Corollary \ref{entalg3}]
We apply Theorem \ref{entalg2} to the $\delta$-hyperbolic space $(X,d)$, where $X$ is the Cayley graph of $\Gamma$
(associated to the system of generators $S_0$) and where $d$ is the canonical extension to this Cayley graph of the 
word-distance $d_{S_0}$ on $\Gamma$ (hence on the vertices of the Cayley graph). Notice that, in this case, 
the diameter of $\Gamma \backslash X$ is equal to $1$, the hyperbolicity constant of $(X,d)$ is $\delta (\Gamma , S_0)$, the 
entropy of $(X,d)$ is equal to the entropy of $ (\Gamma, S_0)$ and is thus bounded above by $\dfrac{M}{1+ \delta (\Gamma , S_0)}$. As $\Gamma$ is not
virtually cyclic, $(X,d)$ is non elementary by Proposition \ref{actioncocompacte} (iv) and (iii).
Let us consider the action (by left translations) of $\Gamma$ on its Cayley graph $(X,d)$; as this action is proper, co-compact and $d$-isometric, we can apply Theorem \ref{entalg2} (ii) and (iii), replacing (in this Theorem) $D$ by $1$ and 
$N_0 (\delta, H, D)$ by $ N_0 \left(\delta (\Gamma , S_0) \,,\, \dfrac{M}{1+ \delta (\Gamma , S_0)}\, ,1 \right)$. The corollary \ref{entalg3} then follows.
\end{proof}

\subsection{Margulis Lemmas for group actions on Gromov-hyperbolic spaces}\label{Margulis}

\small

\emph{For every group $G$ and every $A, B \subset G$, we shall still denote by $A\cdot B$ the image of $A\times B$ 
by the map $(\g, g) \mapsto \g \cdot g$ and define (by induction) $S^k$ as $S^{k-1}\cdot S$.\\
Let us recall that we consider any proper action (by isometries) of any group $\Gamma$ on any
$\delta$-hyperbolic (and thus geodesic and proper) space $(X,d )$ such that the Entropy of $(X,d )$ (with respect to 
at least one $\Gamma$-invariant measure) and the diameter of $\Gamma \backslash  X$ are respectively bounded
above by $H$ and $D$. For any $x \in X$ and any $r>0$, let us also recall that the subgroup generated by 
$\Sigma_r (x) := \{ \gamma \in \Gamma^*  : d(x, \gamma x) \leq r \}$ is denoted by $\Gamma_r(x)$.In the following subsections, making use of the function $N( \bullet )$ which will be defined in 
the following Theorem \ref{BGT}, given any $\delta, \, H, \, D > 0$,  we shall consider the universal constants
\begin{equation}\label{defepsilon0}
R_0 := 20 ( D + \delta )\ \ \ ,\ \ \  N_0 := N \left(\left[3^{12} \, e^{ 490 \,H ( D + \delta ) }\right] + 1 \right) 
\ \ \ ,\ \ \  \e_0 (\delta, H , D) := \dfrac{R_0}{N_0}\ ,
\end{equation}
\begin{equation}\label{defminorant}
 s_0 (\delta, H,D) :=  2\cdot 3^{-12} R_0 \,e^{ - \frac{1}{2} \, (N_0 + 10 )\, (13\, H\,R_0 + 6) } \ .
\end{equation}
}

\normalsize

\subsubsection{A first Margulis Lemma for $\delta$-hyperbolic spaces}\label{Margulis1}

\begin{theorem}\label{BGT} \emph{(E. Breuillard, B. Green, T. Tao, \cite{BGT} Corollary 11.2)}
For every $ p \in \N^*$, there exists $ N = N(p) \in \N^*$ (only depending on $p$) such that the following 
property holds for every group $ G$ and any finite symmetric system $S$ of generators of $G$: if there exists 
some $A \subset G$ which contains $S^{N(p)}$ and satisfies $\# (A\cdot A) \le p \,\# (A)$, then
$G$ is virtually nilpotent.
\end{theorem}

The following corollary settles a first Margulis property under a weak doubling assumption at some scale; 
it improves a previous corollary of
Breuillard, Green and Tao in \cite{BGT}: in fact the \lq \lq packing" hypothesis of \cite{BGT} (see
Definition \ref{packing}) is replaced here by a weak doubling condition on the counting measure of an orbit of the action 
of a given group $\Gamma$, this last hypothesis being weaker than the packing one, as proved in Lemma \ref{packing1}, 
because of the fact that, when counting the maximal number of disjoint balls $B_X(x_i , r)$ included in a bigger one, this 
number is greater when we put no condition on the centers $x_i$ than when we compel the $x_i$'s to be located 
on the same orbit of $\Gamma$.

\begin{corollary}\label{Marg1} For every data $r_0 >0$ and $C_0 > 1$, on every proper
metric space $(X,d )$, for every proper action (by isometries) of any group $\Gamma$ on this space, and for every
$x_0 \in X$, if the counting measure $\mu_{x_0}^{\Gamma}$ of the orbit $ \Gamma \cdot x_0$ verifies the
$C_0$-doubling condition for all the balls of radius $r \in \left[ \frac{1}{2} r_0  ,   \frac{5}{4}  r_0\right]$ centered
at the same point $x_0 $, then the subgroup $\Gamma_{\varrho} (x_0)$ is virtually nilpotent for every 
$\varrho \le  \dfrac{r_0}{N \left(\left[C_0^3 \right] + 1\right)}$ (where $N(\cdot )$ is defined in Theorem \ref{BGT}).
\end{corollary}

In the $\delta$-hyperbolic case, we deduce the following Margulis Lemma

\begin{theorem}\label{cocompact1}
On any $\delta$-hyperbolic space $(X,d )$, for every proper action (by isometries) of any group $\Gamma$ on this 
space such that the diameter of $\Gamma \backslash X$ and the Entropy of $(X,d )$ are respectively bounded above 
by $D$ and $H$, for every $x\in X$ and every $ r \le  \e_0 (\delta, H , D) $, the subgroup $\Gamma_{r} (x) $ 
is virtually cyclic.
\end{theorem}

\small{Notice that the celebrated original Margulis Lemma was settled for Riemannian manifolds whose sectional curvature and dimension are bounded; 
a more recent celebrated result was proved by V. Kapovitch and B. Wilking (\cite{KW}), in which the above universal constant $\e_0$ of Theorem 
\ref{cocompact1} is replaced by another universal constant depending on a lower bound of Ricci curvature, on an upper 
bound of the diameter and on the dimension). With respect to this last result, in Theorem \ref{cocompact1}, we have replaced the bounds on Ricci 
curvature and on the dimension by an upper bound of the Entropy (see section \ref{comparaison} for a comparison between all these hypotheses).}

\normalsize

\medskip
In the proof of these two last results, the key role is played first by Theorems \ref{cocompact2} and \ref{BGT}, 
and secondly by Proposition \ref{sousgroupes} (i).

\begin{proof}[Proof of Corollary \ref{Marg1}]
Once admitted Proposition \ref{sousgroupes} (i), the proof follows the one of the corollary 11.17 of 
\cite{BGT}. In fact, let $ G :=\Gamma_{\varrho} (x_0)$ and $S := \Sigma_{\varrho}  (x_0)$, according to 
Proposition \ref{sousgroupes} (i), the doubling condition assumed on the counting measure of the orbit 
$ \Gamma \cdot x_0$ of $\Gamma$ implies that the counting measure $\mu_{x_0}^{G}$ of the orbit $ G \cdot x_0$ 
of the subgroup $ G$ verifies
\begin{equation}\label{ssgroupe1} 
\dfrac{\mu_{x_0}^{G} \left[ \overline B_X (x_0 , 2 \,r_0)\right]}{ \mu_{x_0}^{G} \left[ \overline B_X (x_0 , r_0)\right]} \le
\dfrac{ \mu_{x_0}^{G} \left[ \overline B_X (x_0 , 2 \,r_0)\right]}{ \mu_{x_0}^{G} \left[ B_X (x_0 , r_0)\right]} \le 
C_0^3 \ .
\end{equation}
We then apply Theorem \ref{BGT}, where we put $p = \left[C_0^3\right] + 1 $ and
$ A := \{ \gamma \in G : d( x_0 , \gamma x_0 ) \le r_0\}$; in fact, taking $\e_0' :=  \frac{r_0}{N (p)}$, the triangle inequality 
and the hypothesis $ \varrho \le \e_0'$ guarantees that $  S^{N(p)} \subset  \Sigma_{\e_0'} (x_0)^{N(p)} \subset A$; 
as $A \cdot A \subset  \{ \gamma \in G  : d( x_0 , \gamma x_0 ) \le 2\, r_0\}$, inequality \eqref{ssgroupe1} implies
$$ \dfrac{\# (A\cdot A)}{\# (A)} \le  \dfrac{\# \left(\{ \gamma \in G: d( x_0 , \gamma x_0 ) 
\le 2 \, r_0\}\right)}{\# \left(\{ \gamma \in G: d( x_0 , \gamma x_0 ) \le r_0\}\right)}
=\dfrac{ \mu_{x_0}^{G} \left[ \overline B_X (x_0 , 2 \,r_0)\right]}{ \mu_{x_0}^{G} 
\left[ \overline B_X (x_0 , r_0)\right]}\le  C_0^3 \le p \ ;$$
Theorem \ref{BGT} then implies that the subgroup $ G =\Gamma_{\varrho} (x_0)$ generated by 
$S := \Sigma_{\varrho}  (x_0)$ is virtually nilpotent.
\end{proof}

\begin{proof}[End of the proof of Theorem \ref{cocompact1}]
$ \Gamma\backslash X$ is compact by Lemma \ref{autofidele} (ii). Let $R_0 := 20 ( D + \delta)$ and 
$C_0 = 3^{4} \, e^{\frac{490}{3} \,H \,  ( D + \delta)} $; under the hypotheses of Theorem \ref{cocompact1}, a consequence of revisiting 
Theorem \ref{cocompact2} (ii) as a doubling property (see comments after this Theorem) is that, for every $x \in X$ and every $ R \ge \frac{1}{2} R_0$, 
one has $ \dfrac{ \mu_{x}^{\Gamma}\big( B_X (x , 2 \,R ) \big) }{ \mu_{x}^{\Gamma}\big( B_X (x ,R )\big) }\le  3^{4} \, e^{\frac{13}{2} \,H \, R}$, 
hence that the counting measure $\mu_{x}^{\Gamma}$ of the orbit $\Gamma\, x$ verifies 
the $ C_0$-doubling condition for all the balls centered at $x$, of radius $R \in \left[ \frac{1}{2} R_0 ,
\frac{5}{4}  R_0\right]$. As $N_0 = N \left(\left[C_0^3 \right] + 1\right) $ and $\e_0 (\delta, H , D) = R_0/N_0$ by 
\eqref{defepsilon0}, the corollary \ref{Marg1} then implies that, for every $ r \le \e_0 (\delta, H , D) $,
$\Gamma_{r} (x)$ is virtually nilpotent. If $\varrho$ is the representation $\Gamma \f \text{Isom} (X,d)$ associated to the action of $\Gamma$ on $ (X,d)$, 
it follows that $\varrho \big(\Gamma_{r} (x) \big)$ is virtually nilpotent too, and thus virtually cyclic according to Proposition \ref{actioncocompacte} (v);
Lemma \ref{reductisom} (vii) then guarantees that $\Gamma_{r} (x)$ is virtually cyclic.
\end{proof}

This leads to the following lower bound of the Margulis constant $L(a,b)$ (see Definition \ref{cteMargulis}):

\begin{corollary}\label{MargL}
On any $\delta$-hyperbolic space $(X,d )$, for every proper action (by isometries) of any group $\Gamma$ on this 
space such that the diameter of $\Gamma \backslash X$ and the Entropy of $(X,d )$ are respectively bounded above 
by $D$ and $H$, every pair $\{a , b \}$ of torsion-free elements of $\Gamma^*$ which does not generate a 
virtually cyclic subgroup verifies $ L(a,b) \ge   \e_0 (\delta, H , D) $.
\end{corollary}

\begin{proof} 
$a$ and $b$ cannot act as parabolic isometries by the proposition \ref{actioncocompacte} (ii); as $a$ and $b$ are torsion-free, they are not elliptic (by 
Remark \ref{kpointsfixes} (i)), they therefore act as hyperbolic isometries by Theorem \ref{ellparahyp}.\\
For sake of simplicity, let $ \e_0 := \e_0 (\delta, H , D)$; if $ L(a,b) < \varepsilon_0$, there
exists $x \in X$ and $ (p,q) \in \Z^*  \times \Z^*$ such that  $ \Max  [d(x, a^p x)\, ; \, d(x, b^q x)] < \varepsilon_0$, hence such that $ \langle a^p , b^q  \rangle \subset \Gamma_{\varepsilon_0} (x)$.
As $ \Gamma_{\varepsilon_0} (x)$ is virtually cyclic by Theorem \ref{cocompact1}, $ \langle a^p, b^q  \rangle$ is virtually cyclic too, and Corollary 
\ref{elementaryaction} (ii) then implies that $ \langle a, b \rangle$ is virtually cyclic, in contradiction with the hypothesis.
\end{proof}

\subsubsection{A lower bound of the diastole}\label{Margulis2}

\begin{defis}\label{thinthick}
In any metric space $(X,d )$, for every proper action (by isometries) of any group $\Gamma$ on this space,

\begin{itemize}

\item  at any point $x \in X$,  $\sys_\Gamma (x)$ is the minimum of $d(x, \g x)$ when $\g$ runs through the elements of $\Gamma^*$
($\sys_\Gamma (x) > 0$ when no element $\g \in \Gamma^*$ fixes $x$),

\item the $ r$-thin subset\footnote{The use of the word \lq \lq thin" is justified by the results of the section
\ref{structuremince}.} of $X$ is the open set $ X_r := \{x \in X : \sys_\Gamma (x)  < r\}$,

\item at any point $x \in X$, $\sys^{\diamond}_\Gamma (x)$ is the minimum of $d(x, \g x)$ when $\g$ runs through the torsion-free
elements of $\Gamma^*$,

\item the torsion-free $ r$-thin subset of $X$ is the open set $ X^{\diamond}_r := \{x \in X : \sys^{\diamond}_\Gamma  (x)  < r\}$,

\item the {\em $ r$-thick subset of $X$} is the complement of $X_{r}$ in $X$.

\end{itemize}
\end{defis}

\begin{prop}\label{Marg2bis}
On any connected non elementary $\delta$-hyperbolic space $(X,d )$, for every proper action (by isometries) of any group 
$\Gamma$ on this 
space such that the diameter of $\Gamma \backslash X$ and the Entropy of $(X,d )$ are respectively bounded above 
by $D$ and $H$, for every $ r \le  \e_0 (\delta, H , D) $ (where $\e_0 (\delta, H , D)$ is defined at \eqref{defepsilon0}), 
the (torsion-free) $ r$-thin subset $X^{\diamond}_{r}$ of $X$ is either empty or not connected; in 
particular there exists a point $x \in X $ such that $\sys^{\diamond}_\Gamma (x) \ge \e_0 (\delta, H , D)$.\\
If moreover $\Gamma$ is torsion-free, then $\dias_{\Gamma}(X) \ge \e_0 (\delta, H , D)$.
\end{prop}

This Proposition will be generalized to actions on metric measured spaces by Theorem \ref{transyst}. However the present Proposition
is more direct and the lower bound it provides is greater than the one given in Theorem \ref{transyst}. This two results are both based on 
the same Proposition \ref{transsystprep}, whose statement and proof are given in section \ref{transport}.

\begin{remark}\label{withtorsion}
\emph{When torsion elements are admitted in $\Gamma$, there is no possible universal lower bound of the diastole
(and a fortiori of the pointwise or global systole) under the hypotheses of Proposition \ref{Marg2bis}, as proved by 
the following example: given any non elementary $\delta$-hyperbolic space $(X_0,d_0 )$ and any proper action 
(by isometries) of any group $\Gamma_0$ on this space such that the diameter of $\Gamma_0 \backslash X_0$ and the 
Entropy of $(X_0,d_0 )$ are respectively bounded above by $D$ and $H$, we construct the metric space $(X,d )$ as 
the product of $(X_0,d_0 )$ with the circle $(\mathbb T\, , \,\e . \,can)$ and the group $\Gamma$ as the product of 
$\Gamma_0$ with the group $R_n$ of the rotations whose angles are $\dfrac{2 k \pi}{n}$ (where $k \in \{0 , 1, \ldots , n-1\}$). For sufficiently small values of $\e$, $(X,d )$ is a non elementary $(\delta +2 \pi \e)$-hyperbolic space and the
canonical product-action of $\Gamma$ on $(X,d )$ verifies all the assumptions of Proposition \ref{Marg2bis}. However,
we have $\forall x \ \sys_\Gamma (x) = 2 \pi \e/n$, which goes to zero when $\e$ goes to zero or when $n$ goes to
$+\infty$.
}
\end{remark}

\begin{proof}[Proof of Proposition \ref{Marg2bis}] 
By Lemma \ref{autofidele} (ii), as the diameter of $\Gamma \backslash X$ is finite, the action of $\Gamma$ on $ (X, d_0)$ is co-compact. 
Under the assumptions of Theorem \ref{Marg2bis}, we may apply Theorem \ref{cocompact1}, which proves that the subgroup $\Gamma_{r} (x) $ 
is virtually cyclic for every $x\in X$ and every $ r \le  \e_0 (\delta, H , D) $. We may thus apply Proposition \ref{transsystprep} (ii), where
we identify the two spaces $(Y,d)$ and $ (X, d_0)$ with $ (X, d)$ and the two actions of $\Gamma$ on these two spaces with the action
of $\Gamma$ on $ (X, d)$ considered here; Proposition \ref{transsystprep} (ii) then guarantees that 
$ X^{\diamond}_{r}$ is disconnected or empty for every $ r \le  \e_0 (\delta, H , D) $. 
Hence $X^{\diamond}_{\e_0}$ is disconnected or empty and $ X \setminus X^{\diamond}_{\e_0} $ is not empty, consequently every $x \in  
X \setminus X^{\diamond}_{\e_0}$ verifies $\sys^{\diamond}_\Gamma (x) \ge \e_0 $.\\
If $\Gamma$ is torsion-free, every $x \in  X \setminus X^{\diamond}_{\e_0} $ satisfies $\sys_\Gamma (x) 
= \sys^{\diamond}_\Gamma (x) \ge \e_0 $, so $\dias_{\Gamma}(X) \ge \e_0$.
\end{proof}

\subsubsection{A lower bound of the global systole for Busemann spaces}

\begin{theorem}\label{Cat2}
On any non elementary, geodesically complete, $\delta$-hyperbolic, Busemann space (in the sense of Definition 
\ref{dconvexe}) $(X,d )$, for every proper action (by isometries) of any group 
$\Gamma$ on this 
space such that the diameter of $\Gamma \backslash X$ and the Entropy of $(X,d )$ are respectively bounded above 
by $D$ and $H$, one has
\begin{itemize}
\item[(i)] $\ell (\g) > s_0 (\delta , H, D)$ for every torsion-free $\g \in \Gamma^*$,

\item[(ii)] $\text{\rm Sys}_\Gamma (X)  >  s_0 (\delta , H, D) $ if $\Gamma$ is torsion-free,
\end{itemize}
where $ s_0 ( \delta , H,D)$ is the universal constant defined at \eqref{defminorant}.
\end{theorem}

In this Theorem \emph{the hypotheses \lq \lq torsion-free" are necessary}, as proved by Remark 
\ref{withtorsion}.

\smallskip
\emph{The hypothesis \lq \lq $\delta$-hyperbolic" is also necessary}: in fact, let us consider the Riemannian product of 
a fixed compact Riemannian manifold $(Y,g)$ of sectional curvature $\sigma \le -1$ by a circle of length $2\, \pi \, \e$,
the systole of the action of its fundamental group $\Gamma$ on its Riemannian universal covering $ \big(\widetilde Y 
\times \R\, ,\, \tilde g \oplus (dt)^2 \big)$ goes to zero with $\e$, though $ \big(\widetilde Y \times \R\, , \, \tilde g 
\oplus (dt)^2 \big)$ verifies all the assumptions of Theorem \ref{Cat2} (except the $\delta$-hyperbolicity) and 
though $\Gamma$ is a torsion-free 
subgroup of the isometry group of $ \big(\widetilde Y \times \R\, , \, \tilde g \oplus (dt)^2 \big)$.\\
\emph{The hypothesis \lq \lq $\diam (\Gamma \backslash X) \le D$" is also necessary}: in fact let $(\Sigma, g_n)$ be 
a sequence of  hyperbolic surfaces whose diameter goes to $+\infty$ with $n$, then their injectivity radius goes to zero
and consequently the systole of the action of their fundamental group $\Gamma$ on their Riemannian cover 
$(\mathbb H^2 , can.)$ goes to zero, though $(\mathbb H^2 , can.)$ verifies all the assumptions of Theorem \ref{Cat2} 
(except the upper bound of the diameter) and though $\Gamma$ is torsion-free.

\smallskip
\emph{The hypothesis \lq \lq $\,\Ent (X,d) \le H$" is also necessary}: in fact, let us consider any surface 
$\Sigma$ obtained by connected sum from two compact pointed hyperbolic surfaces $(\Sigma_1 , g_1, x_1)$ and 
$(\Sigma_2, g_2, x_2)$ whose injectivity radii at $x_1$ and $x_2$ are larger than some fixed $\e_1 > 0$: more 
precisely, for every positive $\e << \e_1$ we glue $\Sigma_1 \setminus B_{\Sigma_1} (x_1, \e)$ and $\Sigma_2 
\setminus B_{\Sigma_2} (x_2, \e)$ at the two ends of the cylinder $ [-1 , 1 ]\times \mathbb S^1 $, identifying 
respectively $ \partial B_{\Sigma_1} (x_1, \e)$ and $ \partial B_{\Sigma_2} (x_2, \e)$ with $\{-1\} \times  \mathbb S^1 $ 
and $\{1\} \times  \mathbb S^1 $; we endow $\Sigma_1 \setminus B_{\Sigma_1} (x_1, \e)$ and $\Sigma_2 
\setminus B_{\Sigma_2} (x_2, \e)$ with their hyperbolic metrics $g_1$ and $g_2$, and $ [-1 , 1 ]\times \mathbb S^1 $ 
with the metric $ h_{\e} := (dt)^2 + b_{\e}(t)^2 (d\theta)^2$, where $  b_{\e}(t) := \dfrac{\sinh \e \,}{\cosh K_{\e}} 
\cdot \cosh (K_{\e}\, t)$ and where $ K_{\e} $ is chosen in order that $K_{\e} \tanh K_{\e} = \dfrac{1}{\tanh \e}$; 
the metric $g_{\e}$ on $\Sigma$ is obtained by gluing these three metrics (this gluing is $C^1$ because 
$b_{\e}(-1) = b_{\e}(1) = \sinh \e $ and $- b_{\e}'(-1) = b_{\e}'(1) = \cosh \e $). Then the surfaces 
$(\Sigma, g_{\e})$ have bounded diameters, sectional curvature $\sigma \le -1$, and injectivity radius 
$\pi \,  b_{\e}(0) < \pi\,\dfrac{\sinh (\e) \,}{\cosh (1/\e)}$\footnote{In fact, as 
$ K_{\e} > \dfrac{1}{\tanh \e} > \dfrac{1}{\e}$, we have $ b_{\e}(0) = \dfrac{\sinh \e \,}{\cosh K_{\e}} < 
\dfrac{\sinh (\e) \,}{\cosh (1/\e)}$, which goes to zero with $\e$.}, the systole of the action of their fundamental
group $\Gamma$ on their Riemannian universal coverings $(\widetilde \Sigma, \tilde g_{\e})$ is equal to $\pi \,  
b_{\e}(0)$) and goes to zero with $\e$, despite the facts that $\Gamma$ is torsion-free, that the diameters of 
$\Gamma \backslash (\widetilde \Sigma, \tilde g_{\e}) = (\Sigma, g_{\e})$ are 
uniformly bounded, that $(\widetilde \Sigma, \tilde g_{\e})$ is a $\CAT (-1)$ space and is therefore a Busemann space
 and a $\delta_0$-hyperbolic space (with $\delta_0 = \ln 3$, by Corollary 1.4.2 page 12 of \cite{CDP}). Hence all the
hypotheses of Theorem \ref{Cat2} are verified by the actions of $\Gamma$ on the spaces $(\widetilde \Sigma, 
\tilde g_{\e})$, except the \lq \lq Entropy bounded" one. In fact $\Ent (\widetilde \Sigma, \tilde g_{\e}) \f +\infty$
when $\e \to 0$ by Theorem \ref{minorsystglobale} (ii)).

\smallskip
However, deciding if the hypothesis \lq \lq $ (X,d)$ is a Busemann space" is necessary in Theorem \ref{Cat2}
is an open problem.

\medskip
Before proving Theorem \ref{Cat2}, we shall provide (in Lemmas \ref{minordist2} and \ref{minordist4}) lower bounds of the distance between the boundary of a  Margulis domain $M_{R}(\gamma)$  and every Margulis domain $M_{\e}(\gamma)$ such that $\ell(\g) \le \e < R$. The first Lemma investigate the case where $R = R_0 := 20 \, (D + \delta)$, the 
second Lemma, using the convexity of the distance, is concerned by the case where $\e = \ell(\g) $ and $ \ell(\g) < R \le R_0$. These two Lemmas are the keys of the proof of Theorem \ref{Cat2}: in fact, if $\e = \ell (\g)$ and if 
$\e_0$ is the lower bound of 
the Margulis constant given by the corollary \ref{MargL}, these Lemmas provide a lower bound of the distance 
between $M_{\e}(\gamma)$ and the boundary of $M_{\e_0}(\gamma)$ which goes to $+\infty$ when $\e \to 0$ and, 
as this distance is bounded from above in terms of the diameter, it comes that $\e = \ell(g)$ cannot be small.

\begin{lemma}\label{minordist2}
On any non elementary $\delta$-hyperbolic space $(X,d )$, for every proper action (by isometries) of any group 
$\Gamma$ on this space such that the diameter of $\Gamma \backslash X$ is bounded above by $D$, for every
torsion-free $\g \in\Gamma^*$, for every $ (x_0 , x) \in X\times X$ such that $ R_{\g} (x_0) \ge 20 ( D + \delta )$ and 
$R_{\g} (x) \le R_{\g} (x_0)$, one has
$$ \dfrac{d(x_0 ,x)}{R_{\g} (x_0)}  \ge \dfrac{\ln \left( 3^{-12} \left( 2 
\left[\dfrac{R_{\g} (x_0)}{R_{\g} (x)}\right] + 1  \right) \right)}{65 \,  \Ent(X,d) \cdot  R_{\g} (x_0) + 14} - \frac{1}{2} \ . $$
\end{lemma}

\begin{proof}
For the sake of simplicity, let $N' := \left[\dfrac{R_{\g}  (x_0)}{R_{\g} (x)}\right]$ and $H := \Ent(X,d)$.
By definition of $R_{\g} (x)$, there exists $p \in \N^*$ such that $ d \left(x , \g^p  x\right)  =  R_{\g} (x) $ 
and the triangle inequality then gives $d(x, \g^{k p} x) < R_{\g}  (x_0) + \varepsilon $ for every $ k \in \Z$ 
such that $ |k| \le N'$ and for any $ \e > 0 $; applying Lemma \ref{majorationorbitale} (ii) (which is a corollary of the Bishop-Gromov-like inequality of 
Theorem \ref{cocompact2} (ii)) in the case where $R = R_{\g}  (x_0) + \e$, we obtain
\small
$$ 2 N' + 1 \le \# \left\{ k \in \Z : d(x , \g^k x) <  R_{\g}  (x_0)+ \varepsilon \right\}
\le 3^{12} \left( 1+ \dfrac{ 2 \, d(x_0 ,x) + \varepsilon}{R_{\g}  (x_0)}\right)^{\!\!\!\! 12 \frac{\ln 3}{\ln 2}}
\hspace{-3mm} \cdot  e^{\frac{65}{2} H \left(R_{\g}  (x_0) + 2\, d(x_0 ,x) + \varepsilon \right) }\ , $$
\normalsize
when $\e \to 0$, noticing that $\Max_{t \in \R^+} \left(t^\beta \, e^{- (\beta /e)\, t} \right) = 1$, we infer the inequalities
\small
$$ 2 N'+ 1 \le  3^{12} \left( 1+ \dfrac{ 2 \, d(x_0 ,x) }{R_{\g} (x_0)}\right)^{12\frac{\ln 3}{\ln 2}} 
 e^{\frac{65}{2} H \left(R_{\g}  (x_0) + 2\, d(x_0 ,x) \right) }
\le 3^{12}  \, e^{\left( \frac{65}{2} H R_{\g}  (x_0) + \frac{12 \ln  3 }{e\,\ln 2}  \right) \big(1 + 2\, \frac{d(x_0 ,x)}{R_{\g}  (x_0)} \big) } \ ,$$
\normalsize
which implies that
$$ \frac{d(x_0 ,x)}{R_{\g}  (x_0)}  \ge \dfrac{\ln \left( 3^{-12} \left( 2 N'+ 1 \right) \right)}{65 H\, R_{\g}  (x_0) +   \dfrac{24 \ln 3}{e\,\ln 2} } - \dfrac{1}{2} \ .$$
\end{proof}

Let us recall that the universal constants $\e_0 (\delta , H, D)$, $N_0$ and $R_0$ are defined in \eqref{defepsilon0} and
that $s_0 (\delta , H, D)$ is defined in \eqref{defminorant}. For the sake of simplicity, we shall use the notations
$s_0 $ for $s_0 (\delta , H, D)$ and  $\e_0 $ for $\e_0 (\delta , H, D)$ in the sequel. Though the final aim of 
Theorem \ref{Cat2} is to prove that 
the hypothesis $\ell (\g) \le s_0 (\delta , H, D)$ of the following Lemma is false, it is interesting to investigate what 
would be the consequences of such an hypothesis when arguing by contradiction.

\begin{lemma}\label{minordist4} 
Under the hypotheses of Theorem \ref{Cat2}, every hyperbolic element $\g \in\Gamma^*$ such that 
$\ell (\g) \le s_0 (\delta , H, D)$ admits a $\g$-invariant geodesic line $c_{\g}$ on which $\g$ acts by 
translation of length $\ell (\g)$ and, for every $r$ such that  $ \e_0/2 < r \le R_0$ 
and for every $x \in X$ such that $ R_{\g} (x) \ge r$, the distance from 
$x$ to the geodesic line $c_{\g}$ is bounded from below by $R_0 \dfrac{r - \e_0/2}{R_0 - \e_0/2}
\left(\dfrac{3}{2} + \dfrac{N_0}{5} \right)$.
\end{lemma}

\begin{proof}
Let $\g^-$ and $\g^+$ be the two points of the ideal boundary which are fixed by $\g$.
From Lemma \ref{invariantgeod} follows the existence of a geodesic line $c_\gamma$ joining $\g^-$ to 
$\g^+$ such that, for every $t \in \R$, $ \g \big( c_{\g} (t) \big) = c_{\g} \big(t + \ell (\g) \big)$. A first consequence 
is that the Margulis domain $M_{\rho} (\g)$ is non empty for every $\rho \ge \ell(\g)$, a second one is that 
$M_{\ell(\g)} (\g)$ coincides with $ M_{\rm{min}} (\g)$ and contains $c_\gamma$.\\
Let us consider any $x$ such that $ R_{\g} (x) \ge r$ and denote by $\bar x$ a projection of $x$ on the 
geodesic line $c_\gamma$, the geodesic completeness and the lemma \ref{localglobal} then prove that 
the geodesic segment $ [ \bar x , x ]$ can be extended as a geodesic ray (denoted by $c$) whose origin is $\bar x$ and which contains $x$. Lemma \ref{convex1} asserts that the function $f : t \mapsto  
d \big(c (t) , c_\gamma \big)$ is convex; as $f(0) = 0$ and as there exists $t_0 = d(\bar x , x) >0$ such that 
$ f(t) = t$ for every $ t \in [0,t_0]$, we have $f(t) \ge t$ for every 
$ t \in [t_0 , +\infty[$, therefore $d \big(c (t) , c_\gamma \big) = t$ for every $t \ge 0$, and thus 
 $c (+\infty)$ is a point of the ideal boundary which is different from $\g^- = c_\gamma (-\infty)$ and from 
$\g^+ = c_\gamma (+\infty)$; it then follows from this and from Lemma \ref{tube} (ii) that 
$c(t) \notin M_{R_0} (\g )$ when $t$ is great enough and, by the Intermediate Value Theorem, there exists
a point $x_0$ on the geodesic ray $c$ such that $ R_{\g } (x_0) = R_0$. As $R_\gamma (x_0) = R_0$ and 
$R_\gamma (\bar x) = \ell (\g)$, Lemma \ref{minordist2} gives:
\begin{equation}\label{distminore}
\dfrac{d(x_0 , \bar x)}{R_0}  \ge \dfrac{\ln \left( 3^{-12} \left( 2 
\left[\dfrac{R_0}{\ell(\g)}\right] + 1  \right) \right)}{65 \, H \, R_0+ 14} - \frac{1}{2}> \frac{3}{2} + \frac{N_0}{5} \ ,
\end{equation} 
where the last inequality is a corollary of the following inequality, which is a consequence of
the assumption $ \ell (\g) \le  s_0 $, of the definition \eqref{defminorant} 
of $ s_0 $ and of the following direct computation:
\begin{equation*}
\dfrac{ \ln  \left[3^{-12} \left( 2 \left[\dfrac{R_0}{\ell(\g)}\right] + 1 \right)\right] }{65 \, H \,R_0 + 14} >
 \frac{1}{2} \dfrac {(N_0 +10 ) (13 H\,R_0 + 28/5)}{65 \, H \,R_0 + 14} = 2 + \frac{N_0}{5} \ .
\end{equation*}
When $ x_0 \in  [ \bar x , x ] $, the inequality \eqref{distminore} implies that 
$d(x , \bar x ) \ge d(x_0 , \bar x) > R_0 \left(\frac{3}{2} + \frac{N_0}{5}\right)$ and proves the lemma in this case.\\
\emph{From now on we shall therefore suppose that $ x \in  [ \bar x , x_0 ] $.}\\
By the properness of the action and the definition of $R_\gamma (x_0)$, there exists $k_0 \in \N^*$ such that
$d(x_0, \g^{k_0} x_0) =  \Min_{k \in \Z^*} d(x_0, \g^{k} x_0) = R_0$. We shall now prove the inequality:
\begin{equation}\label{ellminore}
k_0\, \ell(\g) \le \frac{1}{2} \e_0 \  .
\end{equation}
Arguing by contradiction, let us suppose that $k_0\, \ell(\g) > \frac{1}{2} \e_0$ and let us denote by $k_1$ the smallest 
element of $\N^*$ such that $  k_1  k_0 \, \ell(\g) >  3 \delta $, then
\begin{equation}\label{k1bis}
1 \le k_1 \le \left[ \dfrac{6\, \delta}{\e_0}\right] + 1 \ \ \ \ \ \text{ and }\ \ \ \ \   (k_1 -1)\, k_0 \, \ell(\g) 
\le 3\, \delta \ ,
\end{equation}
where the first property comes from the assumption $ k_0\, \ell(\g) > \frac{1}{2} \e_0$ and the second one 
(which remains valid when $ k_0\, \ell(\g) > 3 \delta$) from the definition of $k_1$. Let $ h := \g^{k_1 k_0}$, 
according to \eqref{majordome} we have
\begin{equation}\label{puissancek0bis}
d(x_0 , h\, x_0) \le  R_0+ (k_1 - 1)\,\ell(\g^{k_0}) + 4\,\delta \ \dfrac{\ln k_1}{\ln 2}\le 
R_0 + (k_1 - 1)\,k_0 \,
\ell(\g) + \dfrac{4 \,\delta}{\ln 2}\, \ln\left( 1 + \frac{6 \, \delta}{\e_0}\right)    \, ,
\end{equation}
where the last inequality is deduced from \eqref{k1bis}. As $ d( \bar x, h\,\bar x ) =  k_1\, k_0 \, \ell(\g)> 3 \delta$, Lemma \ref{ecartement} gives:
$$d(x_0 , h\, x_0) \ge d( x_0 , \bar x) + d( \bar x, h\,\bar x ) + d( h x_0, h\,\bar x ) - 6\,\delta \ge 
2\,  d( x_0 , \bar x)  + k_1\, k_0 \, \ell(\g)  - 6\,\delta \, ;$$
joining this last inequality to \eqref{puissancek0bis} and \eqref{distminore}, it comes
\begin{equation*}
\frac{1}{2} \left( R_0 + \dfrac{4 \,\delta}{\ln 2}\, \ln\left( 1 + \frac{6 \, \delta}{\e_0}\right) + 6\,\delta \right)\ge
d( x_0 , \bar x) >  \left( \frac{3}{2} + \frac{N_0}{5}\right) \, R_0\ .
\end{equation*}

As $ \delta / R_0 \le 1/20$, as $ \delta / \e_0 = N_0\, \delta/R_0$ and as $R_0 \ge 20\,\delta$, we deduce 
that
$$ \frac{3}{20} + \dfrac{1}{10 \ln 2}\, \ln\left( 1 + \frac{3}{10} \, N_0\right) \ge  \dfrac{2}{\ln 2} \,\dfrac{\delta}{R_0}\, \ln\left( 1 + \frac{6 \, \delta}{\e_0}\right) + 3\, \dfrac{\delta}{R_0} 
> 1 + \frac{N_0}{5} ;$$
this last inequality being false because $N_0 \ge 1$, the assumption $k_0\, \ell(\g) > \frac{1}{2} \e_0$ is false and this proves 
\eqref{ellminore}.

\smallskip
The distance between the two geodesic segments $[\bar x , x_0]$ and $[\g^{k_0}\bar x , \g^{k_0} x_0]$ being convex 
(see Definition \ref{dconvexe}), as $ x \in  [ \bar x , x_0 ] $, we get:
$$r \le R_{\g} (x)  \le d(x , \g^{k_0} x ) \le \dfrac{d(\bar x , x)}{d(\bar x , x_0)}\, d(x_0 , \g^{k_0} x_0) 
+  \left( 1 - \dfrac{d(\bar x , x)}{d(\bar x , x_0)}\right)\, d(\bar x , \g^{k_0} \bar x) \ , $$
and, as $ d(x_0 , \g^{k_0} x_0) = R_0$ and $d(\bar x , \g^{k_0} \bar x) = k_0 \,\ell(\g) \le \frac{1}{2} \e_0 $ by 
\eqref{ellminore}, we infer that
\begin{equation*}
d(\bar x, x) \ge \dfrac{r - \e_0/2}{R_0 - \e_0/2}\  d(\bar x , x_0) \ge R_0 \left(\dfrac{r - \e_0/2}{R_0 - \e_0/2}\right)
\left(  \frac{3}{2} + \frac{N_0}{5} \right)\, .
\end{equation*}
\end{proof}

\begin{proof}[Proof of Theorem \ref{Cat2}] 
As (i) evidently implies (ii), we only prove (i). If $\varrho$ is the representation $\Gamma \f \text{Isom} (X,d)$ 
associated to the action of $\Gamma$ on $ (X,d)$, Lemma \ref{reductisom} proves that the canonical action of 
 $\varrho (\Gamma)$ on $(X,d)$ also verifies the hypotheses of Theorem \ref{Cat2}, in particular (see Lemma \ref{reductisom} (v) and (vi)) $\g$ is torsion-free if and only if si $\varrho (\g)$  is torsion-free and then $\ell (\varrho (\g)) 
= \ell (\g)$. Therefore, in order to prove Theorem \ref{Cat2}, it is sufficient to prove it when $\Gamma$ is a subgroup of 
$\text{Isom} (X,d)$ acting properly, this is what we shall suppose in the whole of this proof.\\
According to Lemma \ref{autofidele} (ii),
$ \Gamma\backslash X$ is compact and Proposition \ref{actioncocompacte} (ii) then implies that every torsion-free 
element of $\Gamma$ is an hyperbolic isometry.
\emph {Arguing by contradiction, let us suppose that there exists at least one hyperbolic element $ \g \in \Gamma^*$ 
such that $\ell (\g) \le  s_0$; let us fix this element.}
Proposition \ref{actionelementaire} (ii) proves that, if $\g^-$ and $\g^+$ are the points of the ideal boundary which 
are fixed by $\g$, the subgroup $\Gamma_{\gamma} \, := \, \Big \{ g \in \Gamma :  
g(\{ \g^- , \g^+\}) = \{ \g^- , \g^+\}\Big \}$ is the maximal virtually cyclic subgroup containing $\g$.
By Lemma \ref{minordist4}, we can fix a $\g$-invariant geodesic line $c_{\g}$, which joins $\g^-$ and $\g^+$, 
on which $\g$ acts by translation of length $\ell (\g)$.

For every $\e \in \, ]0 , \frac{\e_0}{1000} [$, let us fix any point $x$ such that $R_{\g} (x) = \varepsilon_0 - \varepsilon$ 
(such a point exists because $\ell (\g) \le s_0< \varepsilon_0 - \varepsilon$ and thus $M_{\varepsilon_0 - \varepsilon} (\g)$ and 
$X \setminus M_{\varepsilon_0 - \varepsilon} (\g)$ 
are both non empty closed subsets of the connected space $X$). Let $\bar x$ be a projection of $x$ onto $c_{\g}$.

For every $g \in \Gamma_{\gamma}$, let us denote 
by $ \bar x_g$ a projection of $g \,\bar x$ on the geodesic line $c_\gamma$; as $g \,\bar x$ is located on the geodesic
line $ g \circ c_\gamma$ which also joins $\g^-$ and $\g^+$, we have $ d (g \, \bar x \, ,\, \bar x_g ) \le 2\,\delta $
according to Proposition \ref{geodasympt} (i), which leads to
$d( x, g\, \bar x )  \ge d( x,  \bar x_g  ) - d (g \, \bar x \, ,\, \bar x_g ) \ge  d( x,  \bar x) -  2\,\delta$. Using this last 
inequality and the lower 
bound of $d( x,  \bar x)$ (in terms of $r := \varepsilon_0 - \varepsilon$) given by Lemma \ref{minordist4}, we get
\begin{equation}\label{compardist} 
\forall g \in  \Gamma_{\gamma} \ \ \ \ \ \ 
d( x, g\, \bar x )  \ge  d( x,  \bar x) -  2\,\delta \ge\frac{1}{2}\, \left(\e_0 - 2\, \e \right) \left(  \frac{3}{2} + \frac{N_0}{5} \right) - 2\,\delta \  ,
\end{equation} 
For any $g \in \Gamma \setminus \Gamma_{\gamma}$, Proposition \ref{actionelementaire} (ii) implies that
the subgroup generated by $\g$ and $g$ is not virtually cyclic and Proposition \ref{actionelementaire} (v) asserts 
that the subgroup generated by $\g$ and $g \g g^{-1}$ is not virtually cyclic. Corollary \ref{MargL} then implies that
$ L(\g , g \g g^{-1}) > \varepsilon_0 - \varepsilon $ hence, by Lemma \ref{disjoints}, the Margulis domains verify $M_{\varepsilon_0 - \varepsilon} (\g) \cap M_{\varepsilon_0 - \varepsilon} (g \g g^{-1}) = \emptyset$. As 
$x \in M_{\varepsilon_0 - \varepsilon} (\g)$, then 
$x \notin M_{\varepsilon_0 - \varepsilon} (g \g g^{-1})$ and we thus have $R_{g \g g^{-1}} (x) > \e_0 - \e$.
As $\g$ acts by translation of length $\ell (\g)$ on the geodesic line $c_\gamma$, then $g \g g^{-1} $ 
acts by translation of length $\ell (\g)=\ell (g \g g^{-1}) \le s_0$ on the geodesic line $g \circ c_\gamma$,
we can thus apply Lemma \ref{minordist4} (where we replace $r$ by $\varepsilon_0 - \varepsilon$) to the 
hyperbolic isometry $g \g g^{-1} $ and we get:
\begin{equation}\label{compardist1} 
\forall g \in \Gamma \setminus \Gamma_{\gamma} \ \ \ \ \ \ 
d( x, g\, \bar x )  \ge  d( x,  g \circ c_\gamma)  \ge \frac{1}{2}\, \left(\e_0 - 2\, \e \right) \left(  \frac{3}{2} + \frac{N_0}{5} \right)  \  ,
\end{equation} 
By definition of the diameter of $\Gamma \backslash X$, one has $D \ge \min_{g \in \Gamma}\  d(x , g\, \bar x)$;
using inequalities \eqref{compardist} and \eqref{compardist1} (when $\e \to 0$) this gives
$$ D \ge \dfrac{\e_0}{2}\left( \frac{3}{2} + \frac{N_0}{5}\right) - 2\, \delta =  
\dfrac{R_0}{2\,N_0}\left( \frac{3}{2} + \frac{N_0}{5}\right) - 2\, \delta > \dfrac{R_0}{10}  - 2\, \delta = 2 \,D\ ,$$
this contradiction proves that the hypothesis $\ell (\g) \le  s_0 (\delta , H, D)$ is never verified, which ends the proof. 
\end{proof}

\section{Transplantation of Margulis' Properties}\label{transport}

\small
In contrast with the statements of section \ref{deltahyp}, the results of the present section concern the actions of groups $\Gamma$ on measured 
metric spaces $(Y, d ,\mu)$ \emph{which are no longer supposed to be Gromov-hyperbolic}; the only geometric hypothesis which will be assumed on $(Y, d ,\mu)$ is the 
\lq \lq bounded entropy" one\footnote{This upper bound on the entropy is a way to fix a limit to the scale: indeed, in the absence of such a rescaling, 
it would be impossible to bound from below distances, displacements or systole, for these invariants go to zero when multiplying
the distance $d$ by a factor $\e$ going to zero. Among all the hypotheses limiting the scale, the \lq \lq bounded entropy 
one" is the weakest possible, as proved in subsection \ref{comparaison} (devoted to comparing the possible hypotheses).}.
The other hypotheses are algebraic intrinsic conditions on the groups $ \Gamma$, these conditions being inherited from the existence 
of a system of generators which endows
$\Gamma$ with a structure of hyperbolic group or from the existence of an action of $\Gamma$ on an hyperbolic metric space admitting some geometric
bounds. The aim of this section is to prove that the actions of such groups on any measured metric space $(Y, d ,\mu)$ (with bounded entropy)
verify several of the aforementioned Margulis'properties: see Theorems \ref{transportnil}, \ref{transyst}, \ref{minorsystglobale}, 
\ref{topotubes} and \ref{tubelong} for precise statements: roughly speaking the Margulis'properties that we already proved (in section \ref{deltahyp}) for 
the actions of these groups on Gromov-hyperbolic spaces are still valid for the actions of these groups on any  measured metric space.
A first version of these ideas was introduced in \cite{BCG} but, in \cite{BCG} the class of groups $\Gamma$ under consideration was more limited: they
were fundamental groups of manifolds with sectional curvature $\sigma \le -1$ and with injectivity radius $\ge i_0 > 0$ and of groups such that any non 
abelian subgroup with two generators admits an injective homomorphism into such a fundamental group.

\normalsize

\smallskip
{\bf Notations: }Let us denote by $\Gamma^{\diamond} $ the subset of $ \Gamma^* = \Gamma \setminus \{e\}$ whose elements are the torsion-free elements  
of $ \Gamma$ ($\Gamma^{\diamond}$ is not a group).
For any $y \in Y$ and any $r>0$, recall that $\Sigma_r (y) := \{ \gamma \in \Gamma^*  : d(y, \gamma y) \le r \}$,
($\Gamma_r(y)$ being the subgroup generated by $\Sigma_r (y) $), that $\sys_\Gamma (y)$  (resp. $\sys^{\diamond}_\Gamma (y)$) is the minimum of $d(y, \g y)$ when $\g$ runs through the elements of $\Gamma^*$ (resp. through the elements of $\Gamma^{\diamond}$); denote by $Y_r$ 
(resp. by  $ Y^{\diamond}_{r}$) the set of the $y$'s in $Y$ satisfying $\sys_{\Gamma} (y) < r$ (resp. 
$\sys^{\diamond}_{\Gamma} (y) < r$).
recall also that $\dias_\Gamma (Y)$ (resp. $\Sys_\Gamma (Y) $) is the supremum (resp. the infimum) on $Y$ of the function 
$y \mapsto \sys_{\Gamma} (y)$.\\
For every subsets $A$ and $B$ of a given group $\Gamma$, recall that $A\cdot B$ is the set of the products 
$a b$, where $(a,b) \in A \times B$ and that $A^n$ is defined (from $A^2 = A \cdot A$) by iteration of the equality 
$A^n = A^{n-1} \cdot A$.

\smallskip
For sake of simplicity, given a co-compact action of a group $\Gamma$ on a metric space $ (X, d_0)$, we call \lq \lq co-diameter of this action" the
diameter of $\Gamma \backslash X$ for the quotient-metric $\bar d_0$ (see definition in Lemma \ref{autofidele}).

\normalsize
Given the positive constants $\delta_0, H_0 , D_0, \e'_0$, recalling the definition \eqref{N1N0} of the universal constant $N_0 := 
N_0 (\delta_0, H_0 , D_0)$, let us define
\begin{equation}\label{universalconstants1}
M_0 = M_0 (\delta_0, H_0 , D_0) := 42 N_0 + 3\ \ \ \ \ ;\ \ \ \ \ \alpha_0 = \alpha_0 (\delta_0, H_0 , D_0) := \dfrac{\ln 2}{M_0 (\delta_0, H_0 , D_0)}
\ \ \ ,\ \ \ 
\end{equation}
\begin{equation}\label{universalconstants2}
M'_0 = M'_0 (\delta_0 , H_0) := N_0 (\delta_0, H_0 , 1) \ \ \ ;\ \ \ \alpha'_0 = \alpha'_0 (\delta_0, H_0) := \alpha_0 (\delta_0, H_0 , 1)\ \ \ ;\ \ \ 
\end{equation}

\begin{equation}\label{universalconstants3}
r_0 = r_0 (\delta_0 , \e'_0) := \dfrac{\e'_0 \ln 2}{13 \delta_0 + 4 \e'_0} \ \ \ ;\ \ \ n_0 = n_0 (\delta_0 , \e'_0) := 
\left[\dfrac{13 \delta_0 + \e'_0}{\e'_0}\right]
\end{equation}

\subsection{The classes of groups which are considered here}\label{classesalg}

When a metric space $(X,d)$ is endowed with a proper, co-compact action by isometries of a group $\Gamma$, recall that its entropy may be computed with 
respect to any $\Gamma$-invariant measure (see Proposition \ref{Entropies1}) and is independent of the choice of this measure.

\subsubsection{Definitions and first properties}

\begin{defi}\label{hyperbolicaction}
Given any real parameters $\delta_0, H_0, D_0 > 0$, we denote by 
\begin{itemize}
\item $\text{\rm Hyp}^{\star}_{\rm action} (\delta_0, H_0, D_0)$ the set of groups $G$ which admit a proper action by isometries on some connected, non elementary, $\delta_0$-hyperbolic metric space whose entropy and co-diameter are bounded from above by $H_0$ and $D_0$ respectively,
\item $\text{\rm Hyp}_{\rm action} (\delta_0, H_0, D_0)$ the set of all non 
virtually cyclic subgroups $\Gamma$ of groups $G$ belonging to the set $\text{\rm Hyp}^{\star}_{\rm action} (\delta_0, H_0, D_0)$.
\end{itemize}
\end{defi}
Notice that a group $G \in \text{\rm Hyp}^{\star}_{\rm action} (\delta_0, H_0, D_0)$ acts co-compactly on some $\delta_0$-hyperbolic metric space
but that the induced action of a subgroup $\Gamma$ of $G$ is in general not cocompact.

\begin{defi}\label{Busemannaction}
Given any real parameters $\delta_0, H_0, D_0 > 0$, we denote by 
\begin{itemize}
\item $ \text{\rm Hyp}^{\star}_{\rm convex} (\delta_0, H_0, D_0)$ the set of groups $G$ which admit
a proper action by isometries on some connected, non elementary, geodesically complete, Busemann, $\delta_0$-hyperbolic metric space 
whose entropy and co-diameter are bounded from above by $H_0$ and $D_0$ respectively,
\item $\text{\rm Hyp}_{\rm convex} (\delta_0, H_0, D_0)$ the set of all non virtually cyclic subgroups $\Gamma$ of groups $G$ belonging to the set
$\text{\rm Hyp}^{\star}_{\rm convex} (\delta_0, H_0, D_0)$.
\end{itemize}
\end{defi}

\begin{defi}\label{hyperbolicgroup}
Given any real parameters $\delta_0, H_0 > 0$, we denote by 
\begin{itemize}
\item $\text{\rm Hyp} (\delta_0, H_0)$ the set of non virtually cyclic groups $G$ which admit a finite system of generators $S_0$ such that 
$G$  is $\delta_0$-hyperbolic (with respect the associated algebraic distance $d_{S_0}$) and such that $\Ent (G, S_0) \le H_0 $,
\item $\text{\rm Hyp}_{\rm sub} (\delta_0, H_0)$ the set of all non virtually cyclic subgroups $\Gamma$ of groups $G$ belonging to the set $\text{\rm Hyp}(\delta_0, H_0)$.
\end{itemize}
\end{defi}

Notice that subgroups of hyperbolic groups are in general not hyperbolic.

\begin{defi}\label{systoleaction}
Given any parameters $\delta_0 ,\e'_0 > 0$, we denote by $\text{\rm Hyp}_{\rm thick} (\delta_0, \e'_0)$ the set of non virtually cyclic groups $\Gamma$ 
which admit a proper (possibly non co-compact) action by isometries on some $\delta_0$-hyperbolic metric space $ (X, d_0)$ such that every 
torsion-free $\g \in \Gamma^*$ verifies $\ell(\g) \ge \e'_0$.
\end{defi}

To these definitions, for convenience, we add the following one.
\begin{defi}\label{systoleaction0}
We define $$\text{\rm Hyp}_{\rm thick} = \bigcup_{\delta_0\geq 0,\, \e'_0>0} \text{\rm Hyp}_{\rm thick} (\delta_0, \e'_0)\,.$$ 
\end{defi}
The idea for this definition is to consider groups $\Gamma$ which belong to some $\text{\rm Hyp}_{\rm thick} (\delta_0, \e'_0)$ without specifying the parameters $\delta_0$ and $\e'_0$, see, for example, Corollary \ref{cor:positiveentropy}, Theorem \ref{Einstein2} and its Corollaries \ref{Einstein3}, 
\ref{cor:S(Einstein)}, \ref{cor:Einstein4}  and Proposition \ref{hypisole}.

Let us remark that the set of all non virtually cyclic subgroups $\Gamma$ of groups $G$ belonging to $\text{\rm Hyp}_{\rm thick} (\delta_0, \e'_0)$ (resp.
to $\text{\rm Hyp}_{\rm thick}$) coincides\footnote{The proof is as follows: let $G$ be any non virtually cyclic group which admits a proper action 
by isometries on some $\delta_0$-hyperbolic metric space $ (X, d_0)$ such that every torsion-free $g \in G^*$ verifies $\ell(g) \ge \e'_0$ then, for every 
(non virtually cyclic) subgroup $\Gamma$ of $G$, the induced action of $\Gamma$ on $ (X, d_0)$ satisfies the same properties.} with $\text{\rm Hyp}_{\rm thick} (\delta_0, \e'_0)$ (resp. with $\text{\rm Hyp}_{\rm thick}$).

Let us now list some properties of the above defined sets of groups:

\begin{lemma}\label{freesg}
Given any parameters $\delta_0 ,\e'_0 > 0$, introduce $ n_0 = n_0 (\delta_0 , \e'_0)$ as in \eqref{universalconstants3}.
In every group $ \Gamma \in \text{\rm Hyp}_{\rm thick}  (\delta_0, \e'_0)$, for every pair of torsion-free elements $a, b \in \Gamma$ which generates a 
non virtually cyclic subgroup, for every integers $p,q \ge n_0$,  one of the two semi-groups generated by $\{ a^{p}  ,  b^{q}\}$ or 
by $\{ a^{p}  ,  b^{- q}\}$ is free.
\end{lemma}

\begin{proof} By definition of $\text{\rm Hyp}_{\rm thick}  (\delta_0, \e'_0)$, there exists some $\delta_0$-hyperbolic metric space $ (X, d_0)$ and a 
proper action (by isometries) of $\Gamma$ on $ (X, d_0)$ such that every pair of torsion-free elements $a, b \in \Gamma$ verifies
$\ell(a), \ell (b) \ge \e'_0$. Corollary \ref{granddeplacement} then  implies that, for every integers $p,q \ge n_0$,  one of the two semi-groups 
generated by $\{ a^{p}  ,  b^{q}\}$ or by $\{ a^{p}  ,  b^{- q}\}$ is free.
\end{proof}

\begin{lemma}\label{HypinHyp}
For every $\delta_0 ,H_0 > 0$,  $\text{\rm Hyp}_{\rm sub}(\delta_0, H_0)$ is a subset of $\text{\rm Hyp}_{\rm action} (\delta_0, H_0, 1)$.
\end{lemma}

\begin{proof} Let $\Gamma$ be any non virtually cyclic subgroup of some group $G \in  \text{\rm Hyp} (\delta_0, H_0)$
By definition of a $\delta_0$-hyperbolic group, there exists a finite system $S_0$ of generators of $G$ such that the Cayley graph 
$ {\cal G}_{S_0}(G) $ of $G$ (endowed with the canonical extension of the algebraic distance $d_{S_0}$, see definitions in Section \ref{notations}) is 
a connected $\delta_0$-hyperbolic metric space whose entropy is smaller than $H_0$ 
(see the examples before Lemma \ref{comparentropi}). As $G$ acts properly (by isometries) by left translations on this Cayley graph and as 
the diameter of ${G \backslash \cal G}_{S_0}(G) $ is smaller than $1$, Propositions \ref{actioncocompacte} (iii) and (iv) then imply that 
$\big({\cal G}_{S_0}(G) , d_{S_0}\big)$ is a non elementary $\delta_0$-hyperbolic space. This proves that $G \in \text{\rm Hyp}^{\star}_{\rm action} (\delta_0, H_0, 1)$ and thus that $\Gamma \in \text{\rm Hyp}_{\rm action} (\delta_0, H_0, 1)$.
\end{proof}

It is clear that $\text{\rm Hyp}(\delta_0, H_0)$ is a subset of $\text{\rm Hyp}_{\rm sub} (\delta_0, H_0)$, we also get

\begin{lemma}\label{nonvirtcyclic}
For every $\delta_0, H_0, D_0 > 0$, every $\Gamma \in \text{\rm Hyp}^{\star}_{\rm action} (\delta_0, H_0, D_0)$ (and consequently every $\Gamma \in 
\text{\rm Hyp}^{\star}_{\rm convex} (\delta_0, H_0, D_0)$) is a non virtually cyclic group. Hence 
$\text{\rm Hyp}^{\star}_{\rm action} (\delta_0, H_0, D_0)$ and $\text{\rm Hyp}^{\star}_{\rm convex}(\delta_0, H_0, D_0)$ are respectively included in 
$\text{\rm Hyp}_{\rm action} (\delta_0, H_0, D_0)$ and in $\text{\rm Hyp}_{\rm convex} (\delta_0, H_0, D_0)$.
\end{lemma}

\begin{proof} 
As $\Gamma \in \text{\rm Hyp}^\star_{\rm action} (\delta_0, H_0, D_0)$, there exists a proper action (by isometries) on some connected, non elementary, $\delta_0$-hyperbolic 
metric space $(X,d_0)$ whose diameter is bounded by $D_0$.
Let $\varrho_0 : \Gamma \f  \text{Isom} (X,d_0)$ be the representation associated to this action; as this action is co-compact (by Lemma \ref{autofidele} (ii))
on a non elementary hyperbolic 
space, $\varrho_0 (\Gamma)$ is not virtually cyclic by Propositions \ref{actioncocompacte} (iii) and (iv), and thus $\Gamma$ is not virtually cyclic.
\end{proof}

\begin{lemma}\label{torsionfreehyperbolic}
For every $\delta_0, H_0, D_0, \e'_0 > 0$, for any $\Gamma \in \text{\rm Hyp}^\star_{\rm action} (\delta_0, H_0, D_0)\cup \text{\rm Hyp} (\delta_0, H_0 )
\cup \text{\rm Hyp}_{\rm thick} (\delta_0, \e'_0)$, the action of $\Gamma$ on the $\delta_0$-hyperbolic metric space $ (X, d_0)$ which is mentioned in the 
definitions of these three sets of groups satisfies the two following properties:
\begin{itemize}
\item[(i)] every element $g \in \Gamma$ acts on $ (X, d_0)$ as an elliptic or hyperbolic isometry (i. e. none of these elements acts as a parabolic isometry),
\item[(ii)] every torsion-free $\g \in \Gamma$ acts as an hyperbolic isometry (i. e. $\ell (\g) > 0$ by Lemma \ref{ellpositive}).
\end{itemize}
Consequenly these two properties remain valid for every subgroup $\Gamma'$ of $\Gamma$.
\end{lemma}

\begin{proof} 
Let $\varrho_0 : \Gamma \f  \text{Isom} (X,d_0)$ be the representation associated to the action of $\Gamma$ on $ (X, d_0)$; for every $\g \in \Gamma$, 
$\varrho_0 (\gamma)$ is elliptic if and only if $\g$ has torsion (by Lemma \ref{reductisom} (v) and Remark \ref{kpointsfixes} (i)).
Consequently, if $\Gamma \in \text{\rm Hyp}_{\rm thick} (\delta_0, \e'_0)$, for every $\g \in \Gamma$, $\varrho_0 (\gamma)$ is either elliptic or verifies 
$\ell \big(\varrho_0 (\gamma)\big) = \ell (\g) > 0$ (by Definition \ref{systoleaction} and Theorem \ref{ellparahyp}), thus $\varrho_0 (\gamma)$ is either elliptic or hyperbolic (by Lemma \ref{ellpositive}) and this proves properties (i) and (ii) in this case.\\
Let us now suppose that $\Gamma \in \text{\rm Hyp}^\star_{\rm action} (\delta_0, H_0, D_0)\cup \text{\rm Hyp} (\delta_0, H_0)$: then, as 
$ \varrho_0 ( \Gamma)$ acts co-compactly on $ (X, d_0)$ (by Lemma \ref{autofidele} (ii)), for every $\g \in \Gamma$, $\varrho_0 (\gamma)$ is a non 
parabolic isometry by Proposition \ref{actioncocompacte} (ii); this implies that, for every torsion-free $\g \in \Gamma$, $\varrho_0 (\gamma)$ is non
elliptic and non parabolic, thus it is an hyperbolic isometry and it satisfies $\ell (\g) > 0$ by Lemma \ref{ellpositive}.
\end{proof}

\begin{lemma}\label{convexinthick}
For every $\delta_0, H_0, D_0 > 0$, $\text{\rm Hyp}_{\rm convex} (\delta_0, H_0, D_0)$ is included in $\text{\rm Hyp}_{\rm thick}  (\delta_0, s_0)$, 
where $ s_0 := s_0 (\delta_0 , H_0, D_0) $ is the universal constant defined at \eqref{defminorant}. 
\end{lemma}

\begin{proof} Let $\Gamma$ be any non virtually cyclic subgroup of any group $G \in\text{\rm Hyp}^{\star}_{\rm convex} (\delta_0, H_0, D_0)$, then 
$G$ is non virtually cyclic (by Lemma \ref{nonvirtcyclic}) and there exists a proper action (by isometries) of $G$ on some (connected, 
non elementary) geodesically complete, Busemann, $\delta_0$-hyperbolic metric space $(X,d_0)$, whose entropy and co-diameter are bounded from above by $H_0$ and $D_0$ respectively.
As the action of $G$ on $(X,d_0)$ satisfies all the hypotheses of Theorem \ref{Cat2}, it follows from this theorem that the action of every 
torsion-free $\g \in G^*$ on $(X,d_0)$ verifies $\ell (\g) > s_0 := s_0 (\delta_0 , H_0, D_0)$; hence 
$G \in \text{\rm Hyp}_{\rm thick}  (\delta_0, s_0)$, and the induced action of $\Gamma$ on $ (X, d_0)$ is proper, by isometries and still verifies $\ell(\g) > s_0$ for every torsion-free $\g \in \Gamma^*$. This proves that $\Gamma \in \text{\rm Hyp}_{\rm thick} (\delta_0, s_0)$.
\end{proof}

\subsubsection{Comparison with acylindrical hyperbolicity}\label{comparacyl}

The notion of acylindrical action was introduced on trees by Z. Sela and on metric spaces by B. Bowditch; the following definitions were explicited by D. Osin in 
\cite{Osi}:

\begin{defis}\label{acylindrique}
The action by isometries of a group $\Gamma$ on a metric space $(X,d)$ is acylindrical if, for every $R>0$, there exist positive numbers $A = A(\Gamma, R)$ and 
$N = N(\Gamma, R)$ such that, for all $x,y \in X$ verifying $d(x,y) \ge A$, one has
$$\# \left\{ \g \in \Gamma : \, d(x, \g x)\le R \text{ and } d(y,\g y)\le R \right\} \le N \,.$$
A group $\Gamma$ is called acylindrically hyperbolic if it admits a non elementary acylindrical action on some Gromov-hyperbolic 
space\footnote{By the corollary 14.4 of \cite{BHS}, the set of acylindrically hyperbolic groups contains the set of hierarchically hyperbolic groups 
(the definition of this last set is given, for instance, in \cite{BHS}).} .
\end{defis}

The hyperbolic spaces involved in Definitions \ref{acylindrique} are not supposed to be proper, and there are many qualitative applications in 
the non proper case. In the following definition, we shall suppose these hyperbolic spaces to be proper, because we aim at quantitative estimates:

\begin{defis}\label{acylindrique0}
A group $\Gamma$ is called properly acylindrically hyperbolic if it admits a non elementary acylindrical action on some proper 
Gromov-hyperbolic space.
The set of properly acylindrically hyperbolic groups will be denoted $\text{\rm Hyp}_{\rm acyl}$ in the sequel.
\end{defis}

Every properly acylindrically hyperbolic group admits a  non elementary, proper, acylindrical action on some proper Gromov-hyperbolic space by the

\begin{remark}\label{acylindricproper}
Every acylindrical action of a group $G$ on a proper metric space $(X,d)$ is proper.
\end{remark}

\begin{proof} Let $\varrho : G \f \text{Isom} (X,d)$ be the representation associated to the action (i.e. $g x := \varrho (g) (x)$ for every $g \in G$ and every
$x \in X$; we define $\Gamma := \varrho (G)$.\\
On the space $C(X,X)$ of continuous maps from $X$ to $X$, the compact-open topology is metrisable\footnote{Let $d_n (f,g) := 
\sup_{x \in K_n} d_X \big( f(x) , g(x)\big)$, the distance on $C(X,X)$ which induces the compact-open topology may be defined as 
$d (f,g) := \sum_{n = 1}^{+ \infty} 2^{-n}\frac{d_n (f,g)}{1 + d_n (f,g)}$.} because $(X,d)$ is a metric space
which is the increasing union (for all $n \in \N^*$) of the compact balls $K_n := \overline B_X(x_0,n)$, where $x_0 \in X$. Hence 
$\text{Isom} (X,d)$ is a metric space for the induced metric.
To prove that $\Gamma$ is a closed discrete subset of $\text{Isom} (X,d)$, endowed with the compact-open topology, it is thus sufficient to prove that 
every converging sequence of elements of $\Gamma$ is stationary.\\ 
Arguing by contradiction, suppose that there exists a converging sequence $ \left(\g_n \right)_{n \in \N}$ of distinct elements of $\Gamma$. 
Let us fix a pair of points $x,y \in X$ verifying $d(x,y) \ge A(G , 1)$ and define $H_{x,y} (1) := \left\{ g \in G : \, d(x, g x)\le 1 \text{ and } d(y,g y)\le 
1 \right\}$. As $ \left(\g_n \right)_{n \in \N}$ uniformly converges on $\{ x , y \}$, there exists $n\in \N$ such that 
$$ \forall p \in \N \  \ \ \ \ d(x, \g_n^{-1} \g_{n+p} x) =  d(\g_n x , \g_{n+p} x) < 1\ \ \ \text{and} \ \ \ 
d(y, \g_n^{-1} \g_{n+p} y) < 1 \, ,$$
and $H_{x,y} (1)$ is thus infinite, in contradiction with the definition of acylindricity, which implies that $H_{x,y} (1)$ is finite.
We conclude that $\Gamma$ is a closed discrete subgroup of $\text{Isom} (X,d)$ and Proposition \ref{discret1} then implies that the action of 
$\Gamma := \varrho (G)$ on $X$ is proper. Now, as $\Ker \varrho $ is included in $H_{x,y} (1)$, it is finite and the action of $G$ is also proper.
Moreover every element of $\Ker \varrho$ is a torsion one.
\end{proof}

Let us first notice that every group $\Gamma$ which admits a proper co-compact action on some proper Gromov-hyperbolic space $(X, d)$ is automatically 
a properly acylindrically hyperbolic group, for the number of $\g \in \Gamma$ such that $d(x , \g x) \le R$ is then bounded above independently of $x \in X$, this 
bound depending on $R$, on the group, on its action and on the diameter of $\Gamma \backslash X$. 
Hence, being verified by all the groups that we shall consider, i.e. by all the groups which admit a proper co-compact action on some proper 
Gromov-hyperbolic space, the condition \lq \lq $\Gamma$ is a properly acylindrically hyperbolic group" alone is not a discriminating criterion, and it is thus 
not conclusive 
when one aims at quantitative estimates of geometric constants (as Margulis' ones) which measure the action of such groups on any metric space, especially 
since we want these quantitative estimates (as those given by Theorems \ref{transportnil}, \ref{transyst}, \ref{minorsystglobale}, \ref{topotubes} and \ref{tubelong}) to be independent of this space and to be uniformly bounded on a whole class of these groups, .\\
Hence, in order to compute such universal estimates, it is necessary to bound some of the parameters involved in Definitions 
\ref{acylindrique} and \ref{acylindrique0} independently 
on $\Gamma$. In this spirit, there exists a previous work by K. Fujiwara (\cite{Fuj}, Theorem 14) who established
a version of Lemma \ref{freesg}, under stronger hypotheses (see the sequel) and with stronger conclusion, it is the

\begin{theorem}\label{Fujiwara} {\rm (K. Fujiwara)} Suppose that $\Gamma$ acts acylindrically on a $\delta$-hyperbolic graph $X$, with 
parameters $A = A(\Gamma, R)$ and $N = N(\Gamma, R)$. Then there exists a constant $M$, depending only on $\delta$, $N(\Gamma, 20 \delta)$, $N(\Gamma, 200 \delta)$, $A(\Gamma, 20 \delta)$ and $A(\Gamma, 200 \delta)$ with the following property: suppose that $a, b \in \Gamma$ act
hyperbolically and assume that, for any $p, q \in \Z\setminus \{0\}$, $\left[a^p, b^q\right]$ is not trivial then, for every $p, q \ge M$, the subgroup
generated by $a^p$ and $ b^q$ is free.
\end{theorem}

In the same spirit, but aiming to make the estimates independent of $\Gamma$, we define the following larger class of groups:

\begin{defis}\label{HTactions} Given a constant $\delta_0 > 0$, a function $ N_0 : \ ]\, 0 , +\infty\,[\  \f \N^*$, and any action by isometries of a group 
$\Gamma$ on a Gromov-hyperbolic space $(X,d)$, this action is said to be $(\delta_0 , N_0 (\cdot) )$-properly acylindrically hyperbolic if it is non elementary, 
if $(X,d)$ is
a proper $\delta_0$-hyperbolic space and if, for every $R>0$, there exists $A \ge 0$ such that, for all $x,y \in X$ verifying $d(x,y) \ge A$, one has
$$ \#  \left\{ \g \in \Gamma : d(x, \g x)<R \text{ and } d(y, \g y)<R \right\} \le N_0 (R)\, .$$
A group is said to be $(\delta_0 , N_0 (\cdot) )$-properly acylindrically hyperbolic if it admits a $(\delta_0 , N_0 (\cdot) )$-properly acylindrically hyperbolic
action on some proper $\delta_0$-hyperbolic space.\\
We moreover denote by $\text{\rm Hyp}_{\rm acyl} \big(\delta_0, N_0 (\cdot)\big)$ the set of $(\delta_0 , N_0 (\cdot) )$-properly acylindrically hyperbolic groups.
\end{defis}

The following result proves that the hypothesis \lq \lq $\Gamma \in \text{\rm Hyp}_{\rm thick} (\delta_0, \e'_0)$" is weaker than the hypothesis 
\lq \lq $\Gamma$ is $(\delta_0 , N_0 (\cdot) )$-properly acylindrically hyperbolic":

\begin{prop}\label{AHimpliesHT}
For every $(\delta_0 , N_0 (\cdot) )$-properly acylindrically hyperbolic action of a group $\Gamma$ on some proper $\delta_0$-hyperbolic space 
$(X,d)$, every torsion-free element $\g \in \Gamma$ verifies $\ell (\g) \ge \e'_0$, where 
$\e'_0 := \dfrac{21\delta_0}{N_0(20 \,\delta_0) +2}$. Consequently $\text{\rm Hyp}_{\rm acyl} \big(\delta_0, N_0 (\cdot)\big) \subset
\text{\rm Hyp}_{\rm thick} (\delta_0, \e'_0)$ and $\text{\rm Hyp}_{\rm acyl} \subset \text{\rm Hyp}_{\rm thick}$.
\end{prop}

For every choice of the function $N_0 (\cdot)$ there exist actions of groups $\Gamma$ on proper $2$-hyperbolic spaces $(X,d)$ which are not 
$(2 , N_0 (\cdot) )$-properly acylindrically hyperbolic, though they verify $\ell (\g) \ge 1$ for every torsion-free element $\g \in \Gamma$: just consider a sequence 
$(\Gamma_n)_{n \in \N}$ of groups such that $\Gamma_n := F \times G_n$, where $F$ is a free group and $G_n$ a finite group such that
$\# G_n \ge n$. Endow each $\Gamma_n $ with the system of generators $ \Sigma_n := S \cup G_n$, where $S$ is the canonical system of generators of $F$.
As the diameter of $G_n$, with respect to the associated algebraic distance $d_{G_n}$ is equal to $1$, the Cayley graph of $(\Gamma_n,  \Sigma_n)$ 
is a proper $\delta_0$-hyperbolic space, with ($\delta_0:= 2$) and, for every pair $x, y$ of points of this graph, every element $\g \in \{1\} \times G_n$ verifies 
$d(x,g.x)\le 1 \text{ and } d(y,g.y)\le 1$. Hence, for every $n > N_0 (1)$, the action is not $(\delta_0 , N_0 (\cdot) )$-properly acylindrically hyperbolic.
On the contrary, for every $n \in \N$, every torsion-free element $\g \in \Gamma_n$ is a non trivial element of $F \times  \{1\}$ and thus verifies $\ell (\g) \ge 1$.

\begin{proof}[Proof of Proposition \ref{AHimpliesHT}] Let $\Gamma$ be any $(\delta_0 , N_0 (\cdot) )$-properly acylindrically hyperbolic group, there 
then exist a proper $\delta_0$-hyperbolic space $(X,d)$ 
and a non elementary action of $\Gamma$ on $(X,d)$ which verifies, for the value $A_0 := A(\Gamma, 20 \delta_0) $ of the parameter and for every
$x,y $ such that $d(x,y) \ge A_0$,

\begin{equation}\label{AHimpliesHT1}
\# \left\{ g \in \Gamma :  \   d(x,g.x)\le 20\, \delta_0 \text{ and } d(y,g.y)\le 20\, \delta_0 \right\} \le N_0 (20 \delta_0)\,.
\end{equation}

By Proposition \ref{actionelementaire} (iii), $\Gamma$ is non virtually cyclic; moreover the action of $\Gamma$  on $(X,d)$  is proper by Remark \ref{acylindricproper}.\\
We first prove, by contradiction, that none of the elements of $\Gamma$ acts by parabolic isometry. Suppose that $\g \in \Gamma$ acts on $(X,d)$ by parabolic 
isometry, it would then admit only one fixed point $ \g^\infty$ on the ideal boundary $ \partial X$; let us fix a geodesic ray $c$ such that $c (+\infty) =  \g^\infty$. For the sake of simplicity, define $N := N_0 (20 \delta_0)$ and the constant $T_N := \sup_{n \in \{1, \ldots , N\}} d \big(c(0) \, ,\, \g^n c(0) \big) $.\\
For every $t , t' \in [T_N , +\infty[ $ such that $|t'-t| \ge A_0$, for all $n \in \Z$ such that $|n| \le N$, Lemma \ref{Busemann} (applied to the parabolic 
isometry $\g^n$) guarantees that 
$$ d\big( c(t) , \g^n c(t)\big) \le 7 \delta_0 < 20 \delta_0  \ \ \ \ \ \text{and} \ \ \ \ \ d\big( c(t') , \g^n c(t')\big) \le 7 \delta_0 
< 20 \delta_0\ .$$
A consequence of this and of \eqref{AHimpliesHT1} is that
$$2\, N +1 = \#  \left\{ \g^n : n \in \Z\ \  \text{s.t. } |n| \le N \right\} \le \#  \left\{ \g^n : d\big( c(t) , \g^n  \,c(t)\big) < 20 \delta_0 \ \ \text{and}  \ \  d\big( c(t') , \g^n  \,c(t')\big) < 20 \delta_0  \right\} $$
$$\le \#  \left\{ g \in \Gamma    :     d\big( c(t) , g \,c(t)\big) < 20 \delta_0 \ \ \text{and}  \ \  d\big( c(t') , g \,c(t')\big) < 20 \delta_0  \right\} \le  N_0 (20 \delta_0) = N \, ,$$
a contradiction which proves that none of the elements of $\Gamma$ acts by parabolic isometry.

\smallskip
A consequence of this, of Theorem \ref{ellparahyp} and of Remark \ref{kpointsfixes} (i) is that every torsion-free $\g \in \Gamma^*$ acts on $(X,d)$ by 
hyperbolic isometry. As $\ell(\g) > 0$ 
by Lemma \ref{ellpositive}, we can define $I := \left\{p \in \Z  \  :  \  \dfrac{11 \, \delta_0}{2 \,\ell (\g) } < |p| \le \dfrac{16 \, \delta_0}{\ell (\g)} \right\}$.
For every $p \in I$, one has $s(\g^p) \ge \ell (\g^p) = |p|\,\ell(\g) > \frac{11 }{2  } \, \delta_0$, where the minimal displacement $s(\g)$ is introduced in 
Definitions \ref{deplacements}, we thus can apply Lemma \ref{quasigeodesic} (i), which implies that, for every $x \in M(\g)$ and every $p \in I$, 
$ d(x , \g^p x) \le \ell(\g^p) + 4 \, \delta_0 \le 20\, \delta_0$. There exist two points $x, y$ of $ M(\g)$ such that  $d(x,y) \ge A_0$ (one can for instance 
choose them on the same geodesic $ c \in \cal{G} (\g)$), which thus satisfy the property:
$$\forall p \in I \ \ \ \ \ \  d(x , \g^p x) \le 20\, \delta_0 \ \ \text{and}\ \ d(y , \g^p y) \le 20\, \delta_0 \, .$$
A consequence is that, for every $p \in I$, $\g^p \in \left\{ g \in \Gamma :  d(x,g.x) \le 20\, \delta_0 \text{ and } d(y,g.y) \le 20\, \delta_0 \right\} $.
We deduce, from this and from \eqref{AHimpliesHT1}, that
$$N_0(20 \delta_0) \ge \# \left\{ g \in \Gamma : d(x,g.x) \le 20\, \delta_0 \text{ and } d(y,g.y) \le 20\, \delta_0 \right\} \ge  \# (I) \ge 2\left( \dfrac{21 \, \delta_0}{2 \,\ell (\g)} - 1\right),$$
which gives $\ell (\g) \ge \e'_0$ and, comparing to Definition \ref{systoleaction}, implies that $\Gamma \in \text{\rm Hyp}_{\rm thick} (\delta_0, \e'_0)$.\\
Now, as every $0$-hyperbolic space is $\delta$-hyperbolic for every $\delta > 0$, it is clear, from Definitions \ref{acylindrique0} that, for every $\Gamma \in \text{\rm Hyp}_{\rm acyl}$, there exist $\delta > 0$,
a proper $\delta$-hyperbolic space $(X,d)$, and an isometric action of $\Gamma$ on $(X,d)$ verifying the following property: for every $R>0$, there exists $A = A(\Gamma, R)$ and $N = N(\Gamma, R)$ such that, for all $x,y \in X$ such that $d(x,y) \ge A$, one has
$$\# \left\{ \g \in \Gamma : \, d(x, \g x)\le R \text{ and } d(y,\g y)\le R \right\} \le N_\Gamma (R) \,,$$
where $N_\Gamma (R) := N(\Gamma, R)$. It follows, by Definition \ref{HTactions}, 
that $\Gamma \in \text{\rm Hyp}_{\rm acyl} \big(\delta, N_\Gamma (\cdot)\big)$ and we just proved that 
$\Gamma \in \text{\rm Hyp}_{\rm thick} (\delta , \e')$, where $\e'  := \dfrac{21\delta}{N_\Gamma(20 \,\delta) +2} > 0$. Hence $ \Gamma \in 
\text{\rm Hyp}_{\rm thick}$.
\end{proof}

\subsection{A first Margulis Property}

\begin{defi}\label{Margulisconstant}
Let $\Gamma$ be a group acting properly, by isometries, on a metric space $ (Y , d) $, the \lq \lq Margulis invariant\rq \rq\, of this action, 
denoted by $\text{\rm Marg}_\Gamma (Y, d) $, is defined by
$$\text{\rm Marg}_\Gamma (Y, d)  := \sup \{r \, | \,  \Gamma _r (x) {\rm\hbox{ is virtually cyclic for  every x}\in Y} \} \, ,$$
where $\Gamma _r (x)$ is the subgroup of $\Gamma$ generated by all the $\g \in \Gamma$ such that $d(x, \g x) \le r$.
\end{defi}

The Margulis invariant is generally defined as $\sup \{r \, | \,  \Gamma _r (x) {\rm\hbox{ is virtually nilpotent for  every x}\in Y} \} $, the above definition
looks different. In fact, the two definitions are equivalent in the cases that we are considering, as proved by the following

\begin{remark}\label{Margulisconstantbis}
Let $\Gamma$ be a group acting properly, by isometries, on a metric space $ (Y , d) $, if $\Gamma$ is torsion-free, then 
$\text{\rm Marg}_\Gamma (Y, d) $ is equal to $ \sup \{r \, | \,  \Gamma _r (x) {\rm\hbox{ is cyclic for  every x}\in Y} \}$,
if $\Gamma$ is Gromov-hyperbolic, then $\text{\rm Marg}_\Gamma (Y, d) $ is equal to  $\sup \{r \, | \,  \Gamma _r (x) {\rm\hbox{ is virtually nilpotent for  every x}\in Y} \} $.
\end{remark}

\begin{proof} If $\Gamma$ is torsion-free then, by Lemma \ref{virtuelcyclique}, $\Gamma _r (x)$ is virtually cyclic if and only if it is cyclic; this proves the 
first part of this remark. Now, if there exists a finite system of generators of $\Gamma$, denoted by $S$, such that the Cayley graph of $\Gamma$, 
endowed with the algebraic accessibility distance $d_S$ associated to $S$ (defined at the last page of section \ref{notations}), is Gromov-hyperbolic,
as the canonical action of $\Gamma$ on its Cayley graph is proper, co-compact and by isometries, we can apply Proposition \ref{actioncocompacte} (v)
which shows that $\Gamma _r (x)$ is virtually cyclic if and only if it is virtually nilpotent; this proves the second part of this remark.
\end{proof}

Given any parameters $\delta_0, H_0, D_0, \e'_0 > 0$, we recall the definitions of  $\alpha_0 = \alpha_0 (\delta_0, H_0 , D_0)$ (in \eqref{universalconstants1}), of $\alpha'_0 = \alpha'_0 (\delta_0 , H_0)$ (in \eqref{universalconstants2}) and of $ r_0 = r_0 (\delta_0 , \e'_0)$ 
(in \eqref{universalconstants3}). The following results concern the sets of groups $\text{\rm Hyp}_{\rm action} (\delta_0, H_0, D_0)$, 
$\text{\rm Hyp}_{\rm sub} (\delta_0, H_0)$ and $\text{\rm Hyp}_{\rm thick}  (\delta_0, \e'_0)$ introduced in Definitions \ref{hyperbolicaction}, 
\ref{hyperbolicgroup} and \ref{systoleaction} respectively.

\begin{theorem}\label{transportnil}
Given any parameters $\delta_0, H_0, D_0, \e'_0 , H> 0$, every proper action (by isometries preserving the measure) of every group $\Gamma$ on 
any metric measured space $ (Y , d, \mu) $ whose entropy is bounded from above by $H$ has the following properties:
\begin{itemize}
\item[(i)] If $\Gamma \in \text{\rm Hyp}_{\rm action} (\delta_0, H_0, D_0)$, for every $y \in Y$, the subgroup $\Gamma_{\alpha_0/H} (y)$ is virtually cyclic; in other words $\text{\rm Marg}_\Gamma (Y, d)  \cdot \Ent (Y,d, \mu ) \ge \alpha_0 (\delta_0, H_0, D_0) > 0$.
\item[(ii)] If $\Gamma \in \text{\rm Hyp}_{\rm sub} (\delta_0, H_0)$, for every $y \in Y$, the subgroup $\Gamma_{\alpha'_0/H} (y)$ is virtually cyclic; in other words $\text{\rm Marg}_\Gamma (Y, d)  \cdot \Ent (Y,d, \mu ) \ge \alpha'_0 (\delta_0, H_0) > 0$.
\item[(iii)] If $ \Gamma \in \text{\rm Hyp}_{\rm thick}  (\delta_0, \e'_0)$, for every $y \in Y$ such that $\Sigma_{r_0/H} (y)$ contains a torsion-free element
(or is empty), the subgroup $\Gamma_{r_0/H} (y)$ is virtually cyclic; in particular, when $\Gamma$ is torsion-free, 
$\text{\rm Marg}_\Gamma (Y, d)  \cdot \Ent (Y,d, \mu ) \ge r_0 (\delta_0, \e'_0) > 0$.
\end{itemize}
\end{theorem}

\begin{remark}
 The fact that we can consider only subgroups $\Gamma$ of groups  $G\in \text{\rm Hyp} (\delta_0, H_0)$, that is groups 
$\Gamma\in\text{\rm Hyp}_{\rm sub} (\delta_0, H_0)$, will be important in the next section, in Theorems \ref{theo} and \ref{theo:poly} and the application of \ref{theo:poly} concerning the example of $\CAT (0)$-cube complexes. 
\end{remark}

\begin{proof}[Proof of (i)] 
As $\Gamma$ is an element of $\text{\rm Hyp}_{\rm action} (\delta_0, H_0, D_0)$, it is a non virtually cyclic subgroup of some group
$G \in \text{\rm Hyp}^{\star}_{\rm action} (\delta_0, H_0, D_0)$. The difficult point here is that the existence of an action of $\Gamma $ on $ (Y , d, \mu) $ 
does not imply the existence of an action of $G$ on $ (Y , d, \mu) $. 
Considering the action of $\Gamma$ on $ (Y , d, \mu) $, let $S$ be the finite subset of $\Gamma$ defined by $S := \Sigma_{\alpha_0 /H} (y)$, it is a finite set for this action is proper. 
As (by definition) the subgroup generated by $S$ in $\Gamma$ is $\langle S\rangle = \Gamma_{\alpha_0 /H} (y)$, and as $S \subset \Gamma \subset G$,
$ \Gamma_{\alpha_0 /H} (y)$ is also the subgroup generated by $S$ in $G$ and, if $\Gamma_{\alpha_0 /H} (y)$ is not virtually cyclic,
then, as $G$  acts (properly, by isometries) on some connected $\delta_0$-hyperbolic metric space $ (X, d_0)$, whose entropy and 
co-diameter are bounded (from above) by $H_0$ and $D_0$ respectively, we may apply Theorem \ref{entalg2} (ii) to this action of $G$ and
conclude that $\Ent (\langle S\rangle , S) > \alpha_0$.\\
From this and from Lemma \ref{comparentropi} one deduces that
$$ H \cdot \frac{\alpha_0}{H} \ge \Ent (Y,d, \mu) \cdot \Max_{\sigma \in S} d(y, \sigma y) \ge \Ent (\langle S\rangle , S) 
> \alpha_0 \,.$$
As this is false, $\Gamma_{\alpha_0 /H} (y)$ is virtually cyclic by contradiction.
\end{proof}
\begin{proof}[Proof of (ii)] 
As $\Gamma \in \text{\rm Hyp}_{\rm action} (\delta_0, H_0, 1)$ by Lemma \ref{HypinHyp}, (ii) is an immediate corollary of (i).
\end{proof}
\begin{proof}[Proof of (iii)]
If $\Sigma_{r_0/H} (y)$ is empty, then $\Gamma_{r_0/H} (y)$ is trivial and the proposition is proved; let us thus suppose that 
$\Sigma_{r_0/H} (y)$ contains a torsion-free element, denoted by $\g_0$. Introduce $ n_0 = n_0 (\delta_0 , \e'_0)$ as in \eqref{universalconstants3}.
By definition of $\text{\rm Hyp}_{\rm thick} (\delta_0, \e'_0)$, there exists some $\delta_0$-hyperbolic metric space $ (X, d_0)$ and a 
proper action (by isometries) of $\Gamma$ on $ (X, d_0)$ such that $\ell(\g_0) \ge \e'_0$. 
Arguing by contradiction, suppose that $\Sigma_{r_0/H} (y)$ generates a non virtually cyclic subgroup of $\Gamma$, then Corollary \ref{elementaryaction} 
(iii) guarantees the existence of some $ \sigma \in \Sigma_{r_0/H} (y)$ such that $\langle \g_0 , \sigma \rangle$ is non virtually cyclic and one deduces 
(by Corollary \ref{elementaryaction} (i)) that $\langle \g_0 , \sigma \g_0 \sigma^{-1} \rangle$ is non virtually cyclic; by Lemma \ref{freesg}, one 
of the two semi-groups generated by $\{ \g_0^{n_0}  ,  \sigma \g_0^{n_0} \sigma^{-1}\}$ or by $\{ \g_0^{n_0}  ,  \sigma \g_0^{- n_0} \sigma^{-1}\}$ 
is free.\\
Now, applying Lemma \ref{comparentropi} to the action of $\Gamma$ on $ (Y , d, \mu) $ and using the fact that \linebreak 
$\Max \big( d(y, \g_0^{n_0} y)\, ;\,  d(y,   \sigma \g_0^{\pm n_0} \sigma^{-1} y)\big) \le (n_0 + 2)\frac{r_0}{H}$, we get:
$$\ln 2 > (n_0 + 2) r_0 =  H (n_0 + 2) \frac{r_0}{H} \ge \Ent (Y,d,\mu)\, \Max \big( d(y, \g_0^{n_0} y)\, ;\,  d(y,   \sigma \g_0^{\pm n_0} \sigma^{-1} y)\big) $$
$$ \ge \Ent \big(\langle  \g_0^{n_0}  ,  \sigma \g_0^{\pm n_0} \sigma^{-1}\rangle , \{ \g_0^{n_0}  ,  \sigma \g_0^{ \pm n_0} \sigma^{-1}\} \big) \ge 
\ln 2 \ ,$$
a contradiction which can be avoided only if $\Sigma_{r_0/H} (y)$ generates a virtually cyclic subgroup of $\Gamma$.
\end{proof}

\subsection{A lower bound of the diastole}\label{diastolebound}

Given any parameters $\delta_0, H_0, D_0, \e'_0 > 0$, we recall the definitions of  $\alpha_0 = \alpha_0 (\delta_0, H_0 , D_0)$ (in \eqref{universalconstants1}), of $\alpha'_0 = \alpha'_0 (\delta_0 , H_0)$ (in \eqref{universalconstants2}) and of $ r_0 = r_0 (\delta_0 , \e'_0)$ 
(in \eqref{universalconstants3}). The following Theorem concerns the sets of groups $\text{\rm Hyp}_{\rm action} (\delta_0, H_0, D_0)$, 
$\text{\rm Hyp}_{\rm sub} (\delta_0, H_0)$ and $\text{\rm Hyp}_{\rm thick}  (\delta_0, \e'_0)$ introduced in Definitions \ref{hyperbolicaction}, \ref{hyperbolicgroup} 
and \ref{systoleaction} respectively.

\begin{theorem}\label{transyst}
Given any parameters $\delta_0, H_0, D_0, \e'_0 , H> 0$, every proper action (by isometries preserving the measure) of a group $\Gamma$ on 
a connected metric measured space $ (Y , d, \mu) $ whose entropy is bounded from above by $H$ has the following properties:
\begin{itemize}
\item[(i)] If $\Gamma \in \text{\rm Hyp}_{\rm action} (\delta_0, H_0, D_0)$, there exists $y\in Y$ such that $ \sys^{\diamond}_\Gamma (y)   \ge \frac{\alpha_0}{H} $; 
in addition, for every $r \in \left] 0 ,\frac{\alpha_0}{H}\right]$, $ Y^{\diamond}_{r}$ is disconnected or empty.
\item[(ii)] If $\Gamma \in \text{\rm Hyp}_{\rm sub} (\delta_0, H_0)$, there exists $y\in Y$ such that $ \sys^{\diamond}_\Gamma (y) \ge \frac{\alpha'_0}{H} $;
 in addition, for every $r \in \left] 0 ,\frac{\alpha'_0}{H}\right]$, $ Y^{\diamond}_{r}$ is disconnected or empty.
\item[(iii)] If $ \Gamma \in \text{\rm Hyp}_{\rm thick}  (\delta_0, \e'_0)$, there exists $y\in Y$ such that $ \sys^{\diamond}_\Gamma (y) \ge \frac{r_0}{H} $; 
in addition, for every $r \in \left] 0 ,\frac{r_0}{H}\right]$, $ Y^{\diamond}_{r}$ is disconnected or empty.
\end{itemize}
If moreover $\Gamma$ is torsion-free, these three results remain valid if we replace $ \sys^{\diamond}_\Gamma (y)$ by $ \sys_\Gamma (y)$
and the disconnectedness of $ Y^{\diamond}_{r}$ by the disconnectedness of $ Y_{r}$.
\end{theorem}

Though $Y^{\diamond}_{r}$ is disconnected (by this Theorem) for the values of $r$ in the interval mentioned above, its projection $\overline Y^{\diamond}_{r}$ on $\Gamma \backslash Y$ may sometimes be connected, see however 
Theorem \ref{topotubes}.

\smallskip
When torsion elements are admitted in $\Gamma$, there is no possible universal lower bound of the diastole
(and a fortiori of the pointwise or global systole) under the hypotheses of Theorem \ref{transyst}, 
as proved by the example given in Remark \ref{withtorsion}; in the same vein, this example also proves that, $\sys^{\diamond}_\Gamma (y) $ and 
$ Y^{\diamond}_{r}$ cannot be replaced by $\sys_\Gamma (y) $ and $ Y_{r}$ in the conclusions of this theorem concerning the case with torsion.

The proof of this Theorem is based on the following Lemma and Proposition:

\begin{lemma}\label{transsystprep0}
For every proper action (by isometries) of any group $\Gamma$ on any Gromov-hyperbolic metric space $ (X, d_0)$, every
hyperbolic element $\g$ of $\Gamma$ is contained in a unique maximal virtually cyclic subgroup.
\end{lemma}

\begin{proof}
Let $\varrho_0 : \Gamma \f  \text{Isom} (X,d_0)$ be the representation associated to the action of $\Gamma$ on $ (X, d_0)$. As $\varrho_0 (\g)$ is a
hyperbolic isometry, Proposition \ref{actionelementaire} (ii) implies that $\varrho_0 (\gamma)$ is contained in 
a unique maximal virtually cyclic subgroup $G_\gamma$ of $\varrho_0 (\Gamma)$ (namely 
$G_\gamma = \{g \in \Gamma : g \big(\{\g^- , \g^+\} \big) = \{\g^- , \g^+\}\}$). A consequence is that $\varrho_0^{-1} (G_\gamma)$ is a virtually cyclic
subgroup of $\Gamma$ (by Lemma \ref{reductisom} (vii)) which contains $\g$; moreover any other virtually cyclic subgroup $G$ of $\Gamma$ which contains 
$\g$ verifies $\varrho_0 (G) \subset G_\gamma$ (because of the maximality of $G_\gamma$) and thus $G \subset \varrho_0^{-1} (G_\gamma)$. 
We conclude that $\varrho_0^{-1} (G_\gamma)$ is the unique maximal virtually cyclic subgroup of $\Gamma$ which contains $\g$. 
\end{proof}

\begin{prop}\label{transsystprep}
Let us consider any proper action (by isometries) of any non virtually cyclic group $\Gamma$ on any connected metric space $ (Y , d)$ and every $r> 0$ such that the subgroup 
$\Gamma_{r} (y)$ is virtually cyclic at every point $y $ of the \lq \lq thin\rq\rq\, subset $ Y^{\diamond}_{r}$; if the same group $\Gamma$ admits another 
proper action (by isometries) on some Gromov-hyperbolic metric space $ (X, d_0)$ such that $\Gamma\backslash X$ is 
compact (or, more generally, such that none of the elements of $\Gamma$ acts by parabolic isometry of $ (X, d_0)$), then the action on 
$Y$ has the following properties:
\begin{itemize}
\item[(1)] for every $y \in Y^{\diamond}_{r}$, there exists a maximal virtually cyclic subgroup $G(y)$ of $\Gamma$ which contains $\Sigma_{r} (y)$;
moreover the map $y \mapsto G (y)$ is constant on each connected component of $Y^{\diamond}_{r}$.
\item[(2)] the \lq \lq thin\rq\rq\, subset $ Y^{\diamond}_{r}$ is empty or disconnected.
\end{itemize}
\end{prop}

\begin{proof}
Let $\varrho_0 : \Gamma \f  \text{Isom} (X,d_0)$ be the representation associated to the action of
$\Gamma$ on $ (X, d_0)$; for every torsion-free $\g \in \Gamma$, $\varrho_0 (\gamma)$ is also torsion-free (by Lemma 
\ref{reductisom} (v)), thus it is an isometry of $ (X, d_0)$ which is non elliptic (by Remark \ref{kpointsfixes} (i)) and non parabolic (by hypothesis in the
non co-compact case, by Proposition \ref{actioncocompacte} (ii) in the co-compact case); hence $\varrho_0 (\gamma)$ is an hyperbolic isometry of 
$ (X, d_0)$ by Theorem \ref{ellparahyp}; this implies that $\ell(\g) > 0$ by Lemma \ref{ellpositive}.

\emph{Proof of (1)} : For every $y \in Y^\diamond_{r}$, as $\Gamma_{r} (y) = \langle \Sigma_{r} (y) \rangle$ is (by hypothesis) a 
virtually cyclic subgroup of $\Gamma$ which contains a torsion-free (and thus hyperbolic) element $\g$, Lemma \ref{transsystprep0} guarantees the existence of a unique 
maximal virtually cyclic subgroup $G (y)$ of $\Gamma$ which contains $\g$, and thus contains $\Sigma_{r} (y) $ by the maximality of $G (y)$. 
By continuity of the distance and 
properness of the action\footnote{In fact, the properness of the action implies the existence of some $\beta = \beta (y) >0 $ such that
$ \Min_{\g \in \Gamma \setminus \Sigma_r (y)} d(y, \g y) >  r + \beta$; it follows that, for every $y' \in B_Y(y,\beta / 2 )$,
one has $\Sigma_r (y')  \subset \Sigma_r (y) $.}, for every $y \in Y^{\diamond}_{r}$ and for every $y'$ in a sufficiently small neighbourhood 
of $y$ (thus contained in the open set $Y^{\diamond}_{r}$), one has $ \Sigma_{r} (y') \subset \Sigma_{r} (y)$, and consequently $G (y') = G (y)$ by the
maximality of $G(y')$.
The map $y \mapsto G (y)$, being locally constant, is constant on each connected component of $Y^{\diamond}_{r}$.

\emph{Proof of (2)} :  Suppose that $Y^{\diamond}_{r}$ is non empty. For every $y \in Y^{\diamond}_{r}$ and every
$g \in \Gamma$, as $\Sigma_{r} (g y) = g \cdot \Sigma_{r} (y) \cdot g^{-1}$ and as  $g \cdot G (y) \cdot g^{-1}$ is the unique maximal virtually cyclic 
subgroup of $\Gamma$ which contains $g \cdot \Sigma_{r} (y) \cdot g^{-1}$, one gets $G (g y) = g \cdot G (y)\cdot g^{-1}$.\\
If $Y^{\diamond}_{r}$ is connected, $G (y)$ is constant on $Y^{\diamond}_{r}$ and equal to a fixed maximal virtually cyclic subgroup $G$, which 
contains all the subgroups $\Gamma_{r} (y)$ for all the  $y$'s in $Y^{\diamond}_{r}$.
For every $y \in Y^\diamond_r$ and every $g \in \Gamma$, the equality $\sys^\diamond_{\Gamma} (g y) = 
\sys^\diamond_{\Gamma} (y)$ implies that $g y $ also belongs to $Y^\diamond_{r}$, and it follows that $\, g \cdot G \cdot g^{-1} = 
g \cdot G (y)\cdot g^{-1} = G (g y) = G  \,$ and thus $ G$ is a normal virtually cyclic subgroup of $\Gamma$.\\
As $G$ is normal in $\Gamma$, and as (for every $ y \in Y^\diamond_r $) $G$ contains $\Sigma_{r} (y)$, which contains some torsion-free element 
$\g$ such that $\ell (\gamma) > 0$ (as proved above) then, for every $g \in \Gamma$, $g \g g^{-1}  \in G$ and $\langle \g , g \g g^{-1} \rangle$ is 
virtually cyclic, thus 
(using Corollary \ref{elementaryaction} (i)) $\langle \g , g\rangle$ is virtually cyclic, and $\langle \g , g\rangle \subset G$ by the maximality of $G$. 
As this is valid for every $g \in \Gamma$, it implies that $\Gamma = G$, in contradiction with the fact that $\Gamma$ is not virtually 
cyclic by hypothesis. Hence $Y^{\diamond}_{r}$ is disconnected.
\end{proof}

\begin{proof}[Proof of Theorem \ref{transyst} (i)] As $\Gamma$ is an element of $\text{\rm Hyp}_{\rm action} (\delta_0, H_0, D_0)$, 
it is a non virtually cyclic subgroup of some group $G \in\text{\rm Hyp}^{\star}_{\rm action} (\delta_0, H_0, D_0)$; as $G$ admits a proper action (by isometries)
on some $\delta_0$-hyperbolic metric space $ (X, d_0)$, and as this action is co-compact by Lemma \ref{autofidele} (ii) because $G \backslash X$ has
bounded diameter, Proposition \ref{actioncocompacte} (ii) implies that none of the elements of $G$ (and thus none of the elements of $\Gamma$) acts 
by parabolic isometry of $ (X, d_0)$. We may apply Theorem \ref{transportnil} (i) which proves that $\Gamma_r (y)$ is virtually cyclic for 
every $ r \le \frac{\alpha_0}{H}$ and every $y \in Y$. 
As all its hypotheses are verified, we may thus apply Proposition \ref{transsystprep} (2), 
which guarantees that $ Y^{\diamond}_{r}$ is disconnected or empty.\\
Choosing now $r = \frac{\alpha_0}{H}$, we obtain that $Y^{\diamond}_{\alpha_0/H}$ is disconnected or empty and, as
$Y$ is connected, it comes that $Y \setminus Y^{\diamond}_{\alpha_0/H} \ne \emptyset$ and any point $y \in Y \setminus Y^{\diamond}_{\alpha_0/H}$
verifies $ \sys^{\diamond}_\Gamma (y)   \ge \frac{\alpha_0}{H} $.\\
In the case where $\Gamma$ is torsion-free, one has $ \sys^{\diamond}_\Gamma (\cdot )  = \sys_\Gamma (\cdot )$ and $Y^\diamond_{r} = Y_r$
is disconnected or empty for every $r \in \left] 0 ,\frac{\alpha_0}{H}\right]$; in
particular there exists $y \in Y$ such that $ \sys_\Gamma (y)   \ge \frac{\alpha_0}{H} $.
\end{proof}

\begin{proof}[Proof of Theorem \ref{transyst} (ii)]
The hypotheses of Theorem \ref{transyst} (ii) imply that one may apply Lemma \ref{HypinHyp}, which proves that $\Gamma \in 
\text{\rm Hyp}_{\rm action} (\delta_0, H_0, 1)$, and we end the proof by applying Theorem \ref{transyst} (i).
\end{proof}

\begin{proof}[Proof of Theorem \ref{transyst} (iii)]
As $ \Gamma \in \text{\rm Hyp}_{\rm thick}  (\delta_0, \e'_0)$, it is non virtually cyclic by definition and 
it admits a proper action (by isometries) on 
some $\delta_0$-hyperbolic space $ (X, d_0)$ such that no element acts as a parabolic isometry (by Lemma \ref{torsionfreehyperbolic}).\\
Considering now the action of $\Gamma$ on $(Y,d)$, notice that, for every $ y \in Y^\diamond_{r_0/H}$, one has $\sys^{\diamond}_\Gamma (y) <r_0/H$, and there exists a torsion-free element $\g \in \Sigma_{r_0/H}(y)$; we may thus apply Theorem \ref{transportnil} (iii) which proves that 
$\Gamma_{r_0/H} (y)$ 
is virtually cyclic for every  $ y \in Y^\diamond_{r_0/H}$; a consequence of this and of the inclusions $ Y^{\diamond}_{r} \subset  Y^\diamond_{r_0/H}$ 
and $\Gamma_{r} (y) \subset \Gamma_{r_0/H} (y)$ (for every $ r \le \frac{r_0}{H}$) is that $\Gamma_r (y)$ is virtually cyclic for every 
$ r \le \frac{r_0}{H}$ and every $ y \in Y^\diamond_{r}$. As all its hypotheses are verified, we may thus apply Proposition \ref{transsystprep} (2), which 
guarantees that, for every $ r \le \frac{r_0}{H}$, $ Y^{\diamond}_{r}$ is disconnected or empty and, as $Y$ is connected, it comes that $Y \setminus Y^{\diamond}_{r_0/H} \ne \emptyset$ and 
any point $y \in Y \setminus Y^{\diamond}_{r_0/H}$ verifies $ \sys^{\diamond}_\Gamma (y)   \ge \frac{r_0}{H} $.\\
In the case where $\Gamma$ is torsion-free, one has $ \sys^{\diamond}_\Gamma (\cdot )  = \sys_\Gamma (\cdot )$ and $Y^\diamond_{r} = Y_r$
is disconnected or empty for every $r \in \left] 0 ,\frac{r_0}{H}\right]$; in particular there exists $y \in Y$ such that $ \sys_\Gamma (y)   \ge \frac{r_0}{H} $.
\end{proof}

An important consequence is the following Corollary \ref{cor:positiveentropy}, whose assumption (iii) is valid if one replaces the class $\text{\rm Hyp}_{\rm thick}$ by the
slightly smaller\footnote{Indeed, instead of considering non virtually cyclic groups, in Definitions \ref{systoleaction2}, we consider groups 
whose action on the ad hoc $\delta_0$-hyperbolic space is non elementary. In the case where this action is not co-compact, this new hypothesis is slightly stronger than the 
previous one (see Propositions \ref{actionelementaire} (iii) and \ref{actioncocompacte} (iv)).} class of groups 
$ \widetilde{\text{\rm Hyp}}_{\rm thick}$, defined as follows:

\begin{defis}\label{systoleaction2}
Given any parameters $\delta_0 ,\e'_0 > 0$, we denote by $\widetilde{\text{\rm Hyp}}_{\rm thick} (\delta_0, \e'_0)$ the set of groups $\Gamma$ 
which admit a proper non elementary (possibly non co-compact) action by isometries on some $\delta_0$-hyperbolic metric space $ (X, d_0)$ such that 
every torsion-free $\g \in \Gamma^*$ verifies $\ell(\g) \ge \e'_0$.\\
We then define $ \widetilde{\text{\rm Hyp}}_{\rm thick}$ as $\bigcup_{\delta_0\geq 0,\, \e'_0>0} \widetilde{\text{\rm Hyp}}_{\rm thick} 
(\delta_0, \e'_0)$.
\end{defis}

\begin{corollary}\label{cor:positiveentropy} Every metric measured space $(Y,d,\mu )$ verifies $\mathrm{Ent}(Y,d,\mu )>0 $ if it admits a proper action (by isometries preserving the measure) of some group $\Gamma$ which verifies one of the three following hypotheses:
\begin{itemize}
\item[(i)] $\Gamma$ is finitely generated and belongs to $\text{\rm Hyp}_{\rm action}(\delta_0, H_0, D_0) 
\cup \text{\rm Hyp}_{\rm sub} (\delta_0, H_0)$,

\item[(ii)] $\Gamma$ contains at least one torsion-free element and belongs to $\text{\rm Hyp}_{\rm action} (\delta_0, H_0, D_0)$, or to
$\text{\rm Hyp}_{\rm sub}(\delta_0, H_0)$, or to $\text{\rm Hyp}_{\rm thick}$,

\item[(iii)] $\Gamma \in \widetilde{\text{\rm Hyp}}_{\rm thick}$.
\end{itemize}
\end{corollary}

\begin{proof} As $ \text{\rm Hyp}_{\rm sub} (\delta_0, H_0)  \subset \text{\rm Hyp}_{\rm action} (\delta_0, H_0, 1) $ by Lemma
\ref{HypinHyp}, it is sufficient
to make the proof in the case where $ \Gamma \in \text{\rm Hyp}_{\rm action} (\delta_0, H_0, D_0) \cup \text{\rm Hyp}_{\rm thick}$.\\
If we had $\mathrm{Ent}(Y,d,\mu ) = 0 $, we could apply Theorem \ref{transyst} for $H$ arbitrarily close to $0$ and infer that 
the supremum (for all $y\in Y$) of $ \sys^{\diamond}_\Gamma (y)$ is infinite. This would imply that $\Gamma$ only has torsion elements, in contradiction
with hypothesis (ii), proving that $\mathrm{Ent}(Y,d,\mu ) > 0 $ in this case.

Hypothesis (iii) is stronger than hypothesis (ii): indeed, if $\Gamma \in \widetilde{\text{\rm Hyp}}_{\rm thick}$, Gromov's classification (see the 
beginning of subsection \ref{nilpotents}) proves that, as the action of $\Gamma$ on some Gromov-hyperbolic space $(X,d_0)$ is non elementary, then it 
contains hyperbolic (thus torsion-free) elements. Hence $\mathrm{Ent}(Y,d,\mu ) > 0 $ under the hypothesis (iii).\\
Hypothesis (i) is stronger than hypothesis (ii) because, if $\Gamma \in \text{\rm Hyp}_{\rm action} (\delta_0, H_0, D_0) $ is finitely generated 
it is a finitely generated and non virtually cyclic subgroup of a group $G$ which admits a proper action (by isometries) on some (connected, non elementary) 
$\delta_0$-hyperbolic metric space $(X,d_0)$ whose entropy and co-diameter are bounded (from above) by $H_0$ and $D_0$ respectively. We can thus apply 
Theorem \ref{entalg2} (i), which proves the existence, in $\Gamma$, of a torsion-free element $\g_0 \in \Gamma$.
This proves that $\mathrm{Ent}(Y,d,\mu ) > 0 $ under the hypotheses (i).
\end{proof}

\subsection{A lower bound of the global systole}\label{subsection:global}

For every proper action (by isometries) of every group $\Gamma$ on any metric space $ (Y , d) $, let us define the \lq \lq topological radius\rq\rq\, $\text{\rm  Toprad} (y)$ as the supremum of the values $r \in \R^+$ such that $\Gamma_{2 r} (y)$
is virtually cyclic.

For any $y\in Y$ the \lq \lq topological radius" is a local version of the \lq \lq Margulis invariant" of Definition \ref{Margulisconstant} and it verifies
$$ \inf_{y \in Y}\text{\rm  Toprad} (y) = \frac{1}{2}\, \text{\rm Marg}_\Gamma (Y, d)\,.$$

\smallskip
Using the term \lq \lq topological radius" for this invariant can be justified by the following observation: when $Y$ is simply connected and when
$\pi : Y \f \overline Y = \Gamma \backslash Y$ is the quotient map (and universal covering of $\overline Y$), $\Gamma_{2 r} (y)$ coincides with the 
image (by the homomorphism associated to the canonical injection) of the fundamental group of the ball $B_{\overline Y} (\pi(y), r)$ of $\overline Y$ in 
the fundamental group of $\overline Y$ (see Theorem \ref{tubelong} (v)), the topological radius is then the supremum of the $r$'s such that the image 
of the fundamental group of the ball $B_{\overline Y} (\pi(y), r)$ of $\overline Y$ in the fundamental group of $\overline Y$ is virtually cyclic.

Recall that $ s_0 := s_0 (\delta_0 , H_0, D_0) $ is the universal constant defined at \eqref{defminorant} and
define the function $ N' : \R_+^* \times \R_+^* \f  \N^*$ by $N' ( \delta, \e) := \left[\dfrac{13 \delta+  \e}{\e}\right]$.
Given any parameters $\delta_0, H_0, D_0, \e'_0 > 0$, The following Theorem concerns the sets of groups $\text{\rm Hyp}_{\rm convex} 
(\delta_0, H_0, D_0)$ and $\text{\rm Hyp}_{\rm thick}  (\delta_0, \e'_0)$ introduced in Definitions \ref{Busemannaction} and \ref{systoleaction} 
respectively, and may be viewed as a generalization of the celebrated \lq \lq Collar Lemma":

\begin{theorem}\label{minorsystglobale}
Given any parameters $\delta_0,\, H_0, \, D_0 ,\, \e'_0,\,  H > 0$, for every element $\Gamma$ of \linebreak
$\text{\rm Hyp}_{\rm convex} (\delta_0, H_0, D_0)$ (resp. of 
$\text{\rm Hyp}_{\rm thick} (\delta_0, \e'_0)$), defining the integer $ n'_0 := N' ( \delta_0 , s_0)$ (resp. $ n'_0 := N' ( \delta_0 , \e'_0)$), for any proper 
action (by isometries preserving the measure) of $\Gamma$ on any connected metric measured space $ (Y , d, \mu) $ whose entropy is bounded from above by $H$, then  

\begin{itemize}
\item[(i)] for every $\e \le \frac{1}{2 n'_0 H}$, at any point $ y \in Y$ such that $\text{\rm sys}_\Gamma^{\diamond} (y) \le \e$, $\Gamma_R (y)$ is 
virtually cyclic for every $R \le \frac{1}{2 H}\   \ln  \left(\dfrac{1 }{n'_0 H \e}\right) - \frac{1}{2} n'_0  \e$, and $\text{\rm  Toprad} (y) \ge \frac{1}{4 H}\   \ln  \left(\dfrac{1 }{n'_0  H \e}\right) - \frac{1}{4} n'_0  \e $.

\item[(ii)] if moreover $(Y,d)$ is path-connected, then\footnote{In the following inequalities the first one is trivial when $\Gamma$ is virtually cyclic
(because $\text{\rm  Toprad} (y) $ is then infinite) and the second one is trivial when $\Gamma \backslash Y$ is non compact (because the diameter 
of $\Gamma \backslash Y$ is then infinite by Lemma \ref{autofidele} (i) and (ii)).} 
$$  \inf_{y \in Y} \text{\rm sys}_\Gamma^{\diamond} (y) > \dfrac{1}{2 n'_0  H}\, \exp \big(- 4 H \, \sup_{y \in Y} [\text{\rm  Toprad} (y)] \big)
\ge \dfrac{1}{2 n'_0  H}\, e^{- 4 H \diam (\Gamma \backslash Y)}\ ,$$
\item[(iii)] (Collar Lemma) if $\Gamma$ is torsion-free, if $y \in Y$ and $\sigma \in \Gamma^*$ verify $d(y, \sigma y) \le \e \le \frac{1}{2 n'_0 H}$, then every 
$\gamma \in \Gamma^*$ which does not commute with $\sigma$ satisfies $d(y, \gamma y) \ge \frac{1}{2 H}\   \ln  \left(\dfrac{1 }{n'_0 H \e}\right) - \frac{1}{2} n'_0  \e $.
\end{itemize}
\end{theorem}

Before proving this Theorem, let us establish two preliminary Lemmas

\begin{lemma}\label{entropielongueurs}
Let $L$ be a free semi-group with $2$ generators $\gamma_1$ and $\gamma _2$ endowed with any distance
$d$ invariant by left translations (by any $\g \in L$). For every $(l_1, l_2) \in \, ]0 , +\infty[^2$ such that 
$ l_1\ge d(e,\gamma_1)$ and $ l_2\ge d(e,\gamma _2)$, the entropy of $(L,d)$ (for the counting measure)
verifies
$$\Ent(L,d) \ge \sup_{a \in ]0,+\infty[} 
\left[ \Max \left(\frac{1}{l_1+a l_2} , \frac{1}{ l_2+a l_1}\right) \cdot 
\Bigl( (1+a) \ln (1+a) - a  \ln a\Bigr)\right]\, .$$
\end{lemma}

\begin{proof}
For every $R>0$, define $L_R := \{ \g \in L : d(e,\gamma ) \le R\}$.
For every $(p_1,p_2) \in (\mathbb {N}^*)^2$ such that $p_1 l_1 + p_2 l_2 \le R$, denote by $\Lambda_{p_1,p_2}$
the set of the elements of $L$ which are products of $p_1$ times $\gamma _1$ and $p_2$ times $\gamma _2$ (in any
order). As $L$ is a free semi-group, one has $\# \big(\Lambda_{p_1,p_2}\big)= C^{p_1}_{p_1+p_2}$. The triangle
inequality and the invariance of $d$ by left translations imply that $\Lambda_{p_1,p_2} \subset L_R$ and one thus gets
$$\# L_R \ge C^{p_1}_{p_1+p_2} \ge  e^{-(1 + \frac{1}{2p_1})} 
\frac{(p_1+p_2)^{p_1+p_2+{1/2}}}{(p_1)^{p_1+1/2} (p_2)^{p_2 + 1/2}}~,$$
where the last inequality follows from the property $\int^{i+1}_i f(x)\ dx \le \frac{1}{2} \left[ f(i+1) + f(i)\right]$ 
when $f(x) = \ln \left(\frac{p+x}{x}\right)$. For every $a \in ]0,+\infty[$, define $p_1 = \left[\frac{R}{ l_1 + a l_2}\right]$ 
and $p_2 = \left[a  p_1\right]$. From the above, as the respective limits of $\frac{p_1}{R}$ and $\frac{p_2}{R}$ (when 
$R \to +\infty$) are $\frac{1}{ l_1 + a l_2}$ and $\frac{a}{ l_1 + a l_2}$, we infer that
$$\liminf_{R \to +\infty}\left[ \frac{1}{R}\,\,\ln (\# L_R)\right] \ge \frac{1}{ l_1 + a\ell_2} \left[ (1+a) \ln (1+a) - a \ln (a)\right]~.$$
The same inequality holds when exchanging $l_1$ and $l_2$ and this ends the proof.
\end{proof}

\begin{lemma}\label{entropieaction}
Let $L$ be a group with $2$ generators $\gamma_1$ and $\gamma _2$ such that one of the two semi-groups 
generated by $\{\gamma_1 , \gamma _2\}$ or by $\{\gamma_1 ,\gamma _2^{-1}\}$ is free and let us denote 
by $L^+$ this free semi-group. For every proper action of $L$ on any metric space $( Y , d ) \, $, for any $L$-invariant
measure $\mu $ on $Y$ and for any $y \in Y$, one has:

\begin{itemize}

\item[(i)] $\ \ \ \ \  \ \ \ \Ent (Y, d ,\mu) . \Max \left[ d( y , \gamma_1 \,y ) \, , \,  d( y , \gamma_2 \,y ) \right] \ 
\ge \ \ln 2 $;

\item[(ii)] $\ \ \ \  \ \ \ \Min \left[d(y,\gamma_1 \, y)\, ,\, d(y,\gamma_2 \, y)\right] > \dfrac{1}{\Ent (Y, d ,\mu)}\cdot
 e^{- \Ent (Y, d ,\mu)\cdot  \Max \left[d(y,\gamma_1 \, y)\,,\, d(y,\gamma_2 \, y)\right]} \ ;.$

\item[(iii)] $\ \ \ $ for every $ \g \in  L^+ \setminus \{e\}$, the set of fixed points of $\g$ is empty.

\end{itemize}
\end{lemma}

\begin{proof}
Replacing eventually $\gamma _2$ by $\gamma _2^{-1}$, we may suppose that the semi-group
generated by $\gamma_1$ and $\gamma _2$ is a free one. For sake of simplicity, we define $H := \Ent (Y, d ,\mu) $.\\
The action being proper, the stabilizer $\text{Stab}_L (y)$ of any point $y$ is finite and thus it contains only 
elements with torsion. As every element $\g \ne e$ of the free semi-group $ L^+$ is torsion-free
(the relation $\g^p = e$ being prohibited), it comes that $L^+ \cap \text{Stab}_L  (y) = \{e\}$ and this proves (iii).\\
Denote by $d_\Sigma $ the algebraic distance on $L$ associated to the system of generators $ \Sigma = 
\{ \gamma_1 , \gamma_2\}$ (for a definition of the algebraic distance, see section \ref{entropies}). As the
number of elements $\g \in L^+$ such that $d(e, \g) \le n$ is greater than $2^n$, we deduce that 
$ \Ent  (L, \Sigma ) \ge \ln 2$. The lemma \ref{comparentropi} ends the proof of (i), for it proves that
$$\Ent (Y, d, \mu) \cdot \Max \left[ d( y , \gamma_1 \,y ) \, , \,  d( y , \gamma_2 \,y ) \right] \ge  \Ent  (L, \Sigma ) 
\ge \ln 2 \ .$$
On $L$, we consider the pseudo-distance $d_y$ defined by $d_y(\gamma ,\gamma ') = d(\gamma \, y,\gamma ' \, y)$; 
as the elements of $L^+ \setminus \{e\}$ do not have any fixed point, $d_y$ is actually a distance when restricted to 
$L^+$.\\
Let $ l_1 = d(y,\gamma _1 \,y) = d_y (e , \gamma_1)$ and $\ell_2 = d(y,\gamma_2 \,y) = d_y (e , \gamma_2)$; as 
the pseudo-distance $d_y$ is invariant by left-translations, we may apply first Lemma \ref{comparentropi}, and 
afterwards Lemma \ref{entropielongueurs}, which give
$$ H \ge \Ent(Y, d, \mu) = \Ent(Y,d , \mu^L_y) \ge \Ent(L^+,d_y, \#) $$
$$ \ge \sup_{a \in ]0,+\infty[} 
\left[ \Max \left(\frac{1}{ l_1 + a l_2}\, , \, \frac{1}{ l_2+a l_1}\right) 
\cdot \Bigl( (1+a) \ln (1+a) - a \ln a \Bigr)\right]~,$$
where $\#$ is the counting measure of $L^+$ and $\mu^L_y$ is (as usual) the counting measure of the orbit of the 
action of $L$ on $Y$. Choosing $a =  H\, l_1 $ in this last inequality and using the fact that $\ln (1+a)>\frac{a}{1+a}$ 
for every $a > 0$, we obtain: $ H \,\ell_2 >  - \ln \left( H \, l_1\right)$, and thus $ H\, l_1 > e^{-H\, l_2}$.\\
By the same proof, choosing $ a = H \, \ell_2 $ in the above inequality, we get $ H\, l_2 > e^{-H\, l_1}$.
It follows that $ \Min \left( l_1, l_2\right) > \dfrac{1}{H}\cdot  e^{- H  \Max \left( l_1, l_2\right)} $ and this ends the 
proof of the part {\it (ii)} of the lemma.
\end{proof}

\begin{proof}[End of the proof of Theorem \ref{minorsystglobale}] 
For sake of simplicity, for any $\e \le \frac{1}{2 n'_0 H}$, we define $R'_\varepsilon := \frac{1}{2 H}\   \ln  \left(\frac{1 }{n'_0 H \e}\right) - 
\frac{1}{2} n'_0 \e$. 
\begin{itemize}
\item {\bf In the case where $\Gamma \in \text{\rm Hyp}_{\rm thick} (\delta_0, \e'_0)$:} there then exists a proper (eventually non co-compact) action 
by isometries on some $\delta_0$-hyperbolic metric space $ (X, d_0)$ such that every torsion-free $g \in \Gamma^*$ verifies $\ell(g) \ge \e'_0$.\\
We first prove (i): as $\text{\rm sys}_\Gamma^{\diamond} (y) \le \e$, there exists a torsion-free element $\sigma \in \Sigma_{\e} (y)$ and the action of
$\sigma $ on $ (X, d_0)$ verifies $\ell(\sigma) \ge \e'_0 > 0$.
Coming back to the action of $\Gamma$ on $ (Y , d, \mu) $ and arguing by contradiction, suppose that the subgroup 
$\Gamma_{R'_\varepsilon} (y)$ generated by $\Sigma_{R'_\varepsilon} (y) $ is 
non virtually cyclic, then (by Corollary \ref{elementaryaction} (iii)) there exists $\g \in \Sigma_{R'_\varepsilon} (y) $ such that $ \langle \sigma, \g \rangle $
is non virtually cyclic and (by Corollary \ref{elementaryaction} (i)) $\langle\sigma , \g \sigma \g^{-1} \rangle$ is then non virtually cyclic.
Applying Lemma \ref{freesg} to the torsion-free pair $\{\sigma , \g \sigma \g^{-1} \}$, we deduce that one of the two semi-groups generated by 
$\{\sigma^{n'_0} , \g \sigma^{n'_0} \g^{-1} \}$ or by $\{\sigma^{n'_0} , \g \sigma^{- n'_0} \g^{-1} \}$ is free (here $n'_0 := N' (\delta_0 , \e'_0) =
\left[\frac{13 \delta+  \e'_0}{\e'_0}\right]$). From this and from Lemma \ref{entropieaction} (ii), as the triangle inequality guarantees that 
$ d(y, \sigma^{n'_0} y) \le n'_0\e$ and $d(y, \g \sigma^{n'_0} \g^{-1}  y) \le 2 R'_\varepsilon + n'_0 \e$, we deduce that
$$ H n'_0 \e \ge \Ent (Y,d, \mu) \,\Min \left[d(y, \sigma^{n'_0} y)\, ,\, d(y, \g \sigma^{n'_0} \g^{-1}  y)\right] $$
$$ > \exp \big(- \Ent (Y,d, \mu) \, \Max \left[d(y, \sigma^{n'_0} y)\,,
\, d(y, \g \sigma^{n'_0} \g^{-1}  y)\right] \big) \ge  e^{- H (2 R'_\varepsilon + n'_0 \e)} = n'_0 H \e\ .$$
This contradiction proves that $ \Gamma_{R'_\varepsilon} (y)$ is virtually cyclic, thus that $\Gamma_{R} (y)$  is virtually cyclic for every 
$R \le R'_\varepsilon$ (for $\Gamma_{R} (y) \subset \Gamma_{R'_\varepsilon} (y)$). A consequence of this and of the definition of the topological 
radius is that $\text{\rm  Toprad} (y) \ge \frac{1}{2} \,R'_\varepsilon$. This ends the proof of (i).

\smallskip
Let us now prove (ii): if $\inf_{y \in Y} \text{\rm sys}_\Gamma^{\diamond} (y) \ge \frac{1}{2 n'_0 H}$, the two inequalities of (ii) are trivially verified,
we may thus suppose in the sequel that $\inf_{y \in Y} \text{\rm sys}_\Gamma^{\diamond} (y)  < \frac{1}{2 n'_0 H}$ and choose $\e$ and $y$ such that
$ \text{\rm sys}_\Gamma^{\diamond} (y)  < \e < \frac{1}{2 n'_0 H}$, the first inequality of (ii) is then a direct consequence of the last inequality of (i).\\
For sake of simplicity, define $D := \diam (\Gamma \backslash Y)$; if $D = +\infty$ the last inequality of (ii) is trivially verified; in order to prove the 
last inequality of (ii), it is thus sufficient to prove that $\sup_{y \in Y} [\text{\rm  Toprad} (y)]  \le D $ 
when $D < +\infty$. Arguing by contradiction, suppose that there exists $y \in Y$ such that $\text{\rm  Toprad} (y) > D $ then $\Gamma_{2D} (y)$ would be
virtually cyclic. We  have seen (see the proof of Proposition \ref{minorentropie}) that Proposition 3.22 of \cite{Gr1} applies to path-connected metric spaces and proves that $\Sigma_{2D} (y)$ is a (finite) system of generators of $\Gamma$, as a consequence $\Gamma = \Gamma_{2D} (y)$ would be
virtually cyclic, in contradiction with the hypothesis. This proves that $\text{\rm  Toprad} (y) \le D $ for every $y \in Y$ and ends the proof of (ii).

\smallskip
We now prove (iii): for every $\g \in \Gamma^*$, if $d(y,\g y) \le R'_\varepsilon $, then the group generated by $\sigma$ and $\g$ is 
virtually cyclic by (i), thus it is cyclic by Remark \ref{virtuelcyclique} and then $\g$ commutes with $\sigma$.

\item {\bf In the case where $\Gamma \in \text{\rm Hyp}_{\rm convex} (\delta_0, H_0, D_0)$:} as $\text{\rm Hyp}_{\rm convex} 
(\delta_0, H_0, D_0)$ is included in $\text{\rm Hyp}_{\rm thick} (\delta_0, s_0)$ by Lemma \ref{convexinthick}, the validity of Theorem 
\ref{minorsystglobale} in the case where $\Gamma \in \text{\rm Hyp}_{\rm convex} (\delta_0, H_0, D_0)$ is a consequence of its validity for every 
$\Gamma \in \text{\rm Hyp}_{\rm thick} (\delta_0, s_0)$.
\end{itemize}
\end{proof}

\subsection{Structure of thin subsets in quotients of metric measured spaces}\label{structuremince}

\subsubsection{Generic topological Lemma}\label{topogeneral}

Let $ (Y , d) $ be any metric space and $\Gamma$ any group acting properly (by isometries) and without fixed point on $ (Y , d ) $. 
Call $ \pi : Y \to \overline Y= \,\Gamma \backslash Y$ the 
quotient map. For every connected open subset $\overline V$ of $\overline Y$, let $\bar j$ be the inclusion mapping
$\overline V \hookrightarrow \overline Y$. The set of connected components of 
$ V:=\pi^{-1} (\overline V)$ being $\left\{ V ^i :  i\in I \right\}$, every $\g \in \Gamma$ maps every connected
component of $\pi^{-1} (\overline V)$  onto some connected component of $\pi^{-1} (\overline V)$; we thus denote by
$ \widehat\Gamma^i_V $ the subgroup of those $\g$ such that $\g (V^i) = V^i $. This allows to define the quotient 
spaces $\widehat Y ^i := \widehat\Gamma^i_V \backslash Y\, $ and $\widehat V ^i := \widehat\Gamma^i_V \backslash V^i\, $ 
and the quotient mapping $ \widehat \pi_i : Y \to \widehat \Gamma_V^i \backslash Y $; the map $\pi$ then goes down to the 
quotient and provides the mapping $\bar\pi_i \,:\,  \widehat \Gamma_V^i \backslash Y \f \Gamma \backslash Y $ which
satisfies $ \pi = \bar\pi_i \circ  \widehat \pi_i $. As the inclusion mapping $ j : V^i \hookrightarrow Y $ trivially commutes 
with the two actions of $ \widehat\Gamma^i_V$ on $ V^i $ and on $Y$, it gives (going down to the quotients) the 
canonical inclusion mapping
$ j' : \widehat\Gamma_V^i \backslash V^i \hookrightarrow \widehat\Gamma_V^i \backslash Y $;
we then define $ \,p_i $ (resp.$\widehat p_i $) as the restriction of $\pi$ (resp. of 
$\widehat \pi_i$) to $V^i $ at the origin and to $\overline  V$ (resp. to $\widehat\Gamma_V^i \backslash V^i $) at the aim, equivalently $ \bar j \circ p_i = \pi \circ j $ (resp. $ j' \circ \widehat p_i = \widehat \pi_i \circ j$).
Consider now the map $\bar\pi_i  \circ j' : \widehat\Gamma_V^i \backslash V^i  
\f \Gamma \backslash Y = \overline Y$, as its image is included in the image of $\pi \circ j$, thus in $\overline  V$ 
(because $\pi \circ j =  \bar \pi_i  \circ  \widehat \pi_i \circ j = \bar \pi_i  \circ j' \circ \widehat p_i $ and $\widehat p_i $ is 
surjective), it gives (by restriction at the aim) a map $\bar p_i :  \widehat\Gamma_V^i \backslash V^i \rightarrow \overline V $ 
such that $\bar\pi_i  \circ j'  =\bar j \circ \bar p_i $. Moreover, one gets $ \bar p_i  \circ \widehat p_i  = p_i$ because
$$\bar j \circ (\bar p_i  \circ \widehat p_i) = \bar\pi_i \circ ( j' \circ \widehat p_i) =  \bar\pi_i \circ (\widehat \pi_i \circ j)
= \pi \circ j = \bar j \circ p_i \ .$$
All these results are summarized in the following diagram:

\vspace{-5mm}

\large
\begin{equation}\label{diagramme}
\xymatrix{
\Gamma \curvearrowright  \hspace{-12mm}            & \hspace{2mm} Y \hspace{2mm}  \ar[r]^{\hspace{-10mm} \widehat \pi_i}   \ar@/^2pc/[rr]^{\pi}                               & \hspace{4mm} \widehat Y^i =  \widehat\Gamma_V^i \backslash Y  \hspace{1mm}  \ar[r]^{\bar \pi_i} 
                                 &  \hspace{5mm} \overline Y  = \Gamma \backslash Y  \hspace{2mm}  \\
\widehat\Gamma^i_V \curvearrowright  \hspace{-10mm}    &  \hspace{2mm} V^i  \hspace{2mm}  \ar@{^{(}->}[u]  \ar[r]^{\hspace{-10mm}  \widehat  p_i}     \ar@/_2pc/[rr]^{p_i}        & \hspace{4mm} \widehat V^i  =  \widehat\Gamma_V^i \backslash V^i  \hspace{1mm}  \ar[r]^{\bar p_i}  \ar@{^{(}->}[u]    
                                 & \hspace{2mm}  \overline V \hspace{4mm}   \ar@{^{(}->}[u]  
 }
\end{equation}
\normalsize

At last, we shall denote by $\bar d$ the distance on $\overline Y = \Gamma \backslash Y$ induced (by quotient) from the distance $d$ on $Y$ (this distance $\bar d$ is defined at Lemma \ref{autofidele} (i)).  

\smallskip
The following Lemma is semi-classical and will be used several times in this section:

\begin{lemma}\label{lemmetopologique} \emph{(generic topological Lemma)}
With the above notations, for every proper action (without fixed point and by isometries) of any group $\Gamma$ on any metric space 
$ (Y , d) $, 
\begin{itemize}
\item[(i)]  for every $i \in I$, $\pi$ and $\widehat \pi_i$ are open mappings, which are locally isometric coverings;
\item[(ii)] for every $i \in I$, $\pi (V^i) = \overline V$ and $p_i$ and $\widehat p_i\,:  V^i  \to \widehat  V^i $ are locally isometric coverings;
\item[(iii)] for every pair $ V^i , \, V^j$ of connected components of $ \pi^{-1} (\overline V)$, there exists 
$ \g_{i,j} \in \Gamma $ such that $ \g_{i,j} (V^i) = V^j$ and $  \widehat\Gamma^j_V = \g_{i,j} \widehat\Gamma^i_V 
\g_{i,j}^{-1} $;
\item[(iv)] for every $i \in I$, $\bar p_i$ is a locally isometric homeomorphism which preserves the path lengths;
\item[(v)] if moreover $Y$ is locally path-connected and simply connected and if, for every connected component $ V^i $ of 
$ V :=\pi^{-1} (\overline V)$ and for every $y \in V^i$, we denote by $ i_*\,: \pi_1 (\overline V ,
\pi (y) ) \to \pi_1 \big(\overline Y , \pi (y) \big)$ the morphism induced by the canonical injection $ i : \overline V 
\hookrightarrow \overline Y$ then the 
subgroup $\widehat\Gamma^i_V$ is identified with $i_*  \pi_1 (\overline V , \pi (y) ) $ via the canonical 
isomorphism\footnote{Given any $y \in Y$, the canonical isomorphism $\Gamma \f \pi_1 \big(\overline Y , \pi (y) \big)$
maps every $\g \in \Gamma$ onto the homotopy class of $\pi \circ c$, where $c$ is any continuous path
with endpoints $y$ and $\g y$.} $\Gamma \f \pi_1 \big(\overline Y , \pi (y) \big) $.
\end{itemize}
\end{lemma}

\begin{proof}
By Lemma \ref{autofidele} (iii), the action of $\Gamma$ on $ (Y , d ) $, being proper and without fixed point, 
is also faithful and discrete, we may thus consider $\Gamma$ as a discrete subgroup of the group of isometries of
$ (Y , d) $ which satisfies $d(y , \g y ) >0$ for every $y \in Y$ and every $\g \in \Gamma^*$. This and the properness 
of the action imply the existence of some $\sigma \in \Gamma^*$ such that $\sys_\Gamma (y) = d(y , \sigma\, y) >0$.

\begin{itemize}

\item \emph{Proof of (i)}: Choose any $y \in Y$ and any $\e \le \frac{1}{2} \, \sys_\Gamma (y) $;
as, for any $(y,z) \in Y^2$, one has $ \bar d \big(\pi (y) , \pi (z) \big) = \inf_{\g \in \Gamma} d(\g y , z)$ then $ \pi^{-1} 
\big(B_{\overline Y} (\pi (y) , \e) \big)$ is the disjoint union of the balls $ B_Y  ( \g y , \e )$ (for all the $\g$'s in $\Gamma$).
Consequently $\pi$ is a covering and, for any $\g \in \Gamma$, $\pi : \big( B_Y ( \g y , \e/2 ) \,,\, d \big) \to 
\big( B_{\overline Y} ( \pi (y) , \e/2 ) \,,\, \bar d \big)$ is an isometry. Moreover, for every open subset $U$ of $Y$,
$\pi^{-1} \big( \pi (U)\big) = \cup_{\g \in \Gamma} \, \g (U)$ is also an open subset of $Y$ and, by definition of the
quotient-topology, $\pi (U)$ is an open subset of $\overline Y$ and $\pi$ is an open map.\\
The proof of the fact that $\widehat \pi_i$ is a locally isometric covering and an open map is similar (as the action of 
$\widehat\Gamma^i_V$ is the restriction of the action of $\Gamma$, notice that 
$\sys_{\widehat\Gamma^i_V} (y) \ge \sys_\Gamma (y) >0$). 

\item \emph{Proof of (ii)}: As $ V^i $ is open, for every $y \in V^i$, there exists $\e > 0 $ such that 
$B_Y(y , \e) \subset V^i$; from this and from the definition of $ \widehat\Gamma^i_V $, it follows that, for every 
$ \g \in \widehat\Gamma^i_V$, one has $B_Y ( \g y , \e) = \g \big( B_Y(y, \e) \big) \subset \g (V^i) = V^i$. 
Replacing eventually $\e$ by a smaller value (still called $\e$), (i) implies that, for every $ \g \in \widehat\Gamma^i_V$,
$\widehat \pi_i$ is an isometry from the ball $ B_Y ( \g  y , \e) \subset V^i$ onto the ball $B_{\widehat Y_i} (\widehat \pi_i(y) , \e) 
\subset \widehat \pi_i (V^i) = \widehat p_i (V^i)$ of $ \widehat Y_i \,$; as $ \widehat \pi_i^{-1} \big( B_{\widehat Y_i} (\widehat \pi_i(y) , \e) \big) = \cup_{\g \in \widehat\Gamma^i_V} B_Y ( \g y , \e)= \widehat p_i^{-1} \big( B_{\widehat Y_i} (\widehat p_i(y) , \e) \big)$, this proves that $\widehat p_i$ is still a locally isometric covering.\\
In order to prove that $p_i$ is also a covering, first prove that $\pi (V^i) = \overline V$: as (by construction)
$ \pi^{-1} (\overline V)$ is (globally) $\Gamma$-invariant, every $\g \in \Gamma$ only exchanges the connected 
components of $ \pi^{-1} (\overline V)$ and, as two of these connected components are either identical or disjoint,
for every pair $V_i, \, V_j$ of these connected components and every $\g , \, g \in \Gamma$, one has
\begin{equation}\label{intersectioncomp}
\g (V^i) \cap g(V^j) \ne \emptyset \implies  \g (V^i) = g (V^j) \ ,
\end{equation}
which implies that
\begin{equation}\label{intersectionproj}
\pi (V^i) \cap \pi (V^j) \ne \emptyset \implies \pi (V^i) = \pi (V^j) \ :
\end{equation}
indeed, if there existed some $ y \in Y$ such that $\pi (y) \in \pi (V^i) \cap \pi (V^j)$, there would exist $\g , \, g 
\in \Gamma$ such that $ y \in \g (V^i)$ and $ y \in g (V^j)$ and this, by \eqref{intersectioncomp}, would imply that
$ \g (V^i) = g (V^j)$, and thus that $\pi (V^i) = \pi (V^j)$.\\
As $ \overline V = \cup_{i \in I} \, \pi (V^i)$, where (by (i) and \eqref{intersectionproj}) the subsets 
$\big(\pi (V^i)\big)_{i \in I}$ are open and pairwise disjoint or identical,
the connectedness of $ \overline V $ guarantees that all the $\pi (V^i)$'s coincide, and
we thus conclude that $p_i (V^i) = \pi (V^i) = \overline V$ for every $ i \in I$. This also proves that $p_i$ and $\bar p_i$ are surjective maps.\\
Mimicking the above proof, for every $ y \in V^i$, there exists $\e > 0 $ (small enough) such that one has simultaneously 
$\e \le \frac{1}{4} \, \sys_\Gamma (y)$ and $B_Y(y , \e) \subset V^i$, this implies that $B_Y ( \g y , \e) = \g 
\big( B_Y(y, \e) \big) \subset \g (V^i)$ for every $ \g \in \Gamma$; from this we deduce that\\
--  if $\g \in \widehat\Gamma^i_V$, then $B_Y ( \g y , \e) \subset V^i  = \g (V^i) $,\\
--  if $\g \notin \widehat\Gamma^i_V$, then $ V^i \ne \g (V^i) $, thus $  \g (V^i)  \cap V^i = \emptyset$ by 
\eqref{intersectioncomp}, and $ B_Y ( \g y , \e) \cap V^i = \emptyset$;\\
from this, we obtain that
$$ p_i^{-1} \big( B_{ \overline Y }(\pi(y) , \e )\big) = V^i \cap \pi^{-1} \big(  B_{ \overline Y } (\pi(y) , \e )\big) = \bigcup_{\g \in \Gamma} 
 \big(B_Y ( \g y , \e) \cap V^i \big) = \bigcup_{\g \in \widehat\Gamma^i_V} \, B_Y ( \g y , \e) \ ,$$
where, in this series of equalities, the last union is a disjoint one. This proves that $ p_i$ is a covering, and moreover 
a local isometry because $\pi$ is an isometry from the ball $ B_Y ( \g  y , \e) \subset V^i$ onto the ball 
$ B_{\overline Y}(\pi(y) , \e ) \subset \overline V$ of $(\overline Y , \bar d)$.

\item \emph{Proof of (iii)} : A consequence of (ii) is that $ \pi (V^i) = \overline V = \pi (V^j)$; this implies that, for every 
$y \in V^i$, the set $\pi^{-1} \big( \{\pi (y)\} \big) \cap V^j$ is not empty, and there thus exists some 
$ \g_{i,j} \in \Gamma $ such that $ \g_{i,j} (y) \in V^j$; as $ \g_{i,j} (V^i) \cap V^j \ne \emptyset$, using 
\eqref{intersectioncomp}, we therefore get $ \g_{i,j} (V^i) = V^j $. Another consequence of this is the global invariance
of $V^j$ (resp. of $V^i$) under the action of $\g_{i,j}  \widehat\Gamma^i_V  \g_{i,j}^{-1} $ (resp. of 
$\g_{i,j}^{-1}  \widehat\Gamma^j_V  \g_{i,j} $); this proves that $ \widehat\Gamma^j_V = \g_{i,j} 
\widehat\Gamma^i_V \g_{i,j}^{-1} $.

\item \emph{Proof of (iv)} :  The open set $V^i$ being globally invariant by $\widehat\Gamma^i_V$, one has
$ \widehat \pi_i^{-1} \left( \widehat \pi_i (V^i)\right) = V^i$, which implies that 
$ j' \circ \widehat p_i (V^i) = \widehat \pi_i \circ j (V_i) = \widehat \pi_i (V_i)$ 
is an open subset of $\widehat Y^i $. As $\widehat p_i $ and $p_i = \bar p_i \circ \widehat p_i  $ are surjective, the map $\bar p_i$ is well defined and bijective
from $ \widehat V^i $ onto $\overline V$ because, for every $y, \,y' \in V^i$, one has:
$$\bar p_i \circ \widehat p_i (y) = \bar p_i \circ \widehat p_i (y') \iff \bar j \circ  p_i (y) = \bar j \circ  p_i  (y') \iff
\pi \circ j(y) = \pi \circ j(y') \iff \pi  (y) = \pi (y') $$ 
$$ \iff \exists \g \in \Gamma \text{ such that } y' = \g y \iff  \exists \g \in \widehat\Gamma^i_V \text{ such that } 
y' = \g y   \iff   \widehat p_i (y) = \widehat p_i (y') \ ,$$
where the fifth equivalence comes from the fact that, if $y$ and $ \g y$ both lie in $V^i$, then $V^i \cap \g (V^i) \ne \emptyset$ and, using 
\eqref{intersectioncomp}, $ \g \in \widehat\Gamma^i_V$. As $\widehat p_i$ and $ p_i$ are moreover local isometries, then $\bar p_i ,$ is also a local 
isometry, and consequently a homeomorphism. 

\item \emph{Proof of (v)} : Notice that, by (i) and the fact that $Y$ is simply connected, $\pi : Y \f \bar Y$ is the universal covering of $\bar Y$.
Being simultaneously connected and locally path-connected, the connected components $V^i$ of $ \pi^{-1} (\overline V)$ 
are path-connected. Consider any of these connected components, denoted by $V^i$, any point $y \in V^i$ and the point  $\bar y = \pi (y) $. The 
canonical isomorphism $\psi_y : \Gamma \f \pi_1 \big(\overline Y , \bar y\big)$ maps every $\g \in \widehat\Gamma^i_V$, onto the homotopy class 
$[ i \circ \pi \circ c ] \in \pi_1 (\overline Y , \bar y ) $ of the loop $ i \circ \pi \circ c $, where $c$ is any continuous path with endpoints $y$ and $\g y$ which lies in $V^i$ (such a path exists for $V^i$ is path-connected); as $[ i \circ \pi \circ c ]  =  i_* ([\pi \circ c ]) \in  i_* \big(\pi_1 (\overline V , \bar y ) \big)$, it proves 
that $\psi_y ( \widehat\Gamma^i_V) \subset  i_* \big(\pi_1 (\overline V , \bar y ) \big)$.\\
Conversely, every element $ \beta \in  i_* \big(\pi_1 (\overline V , \bar y ) \big)$ can be represented by a loop $\bar c$ lying in $\overline V $, with basepoint 
$\bar y$; this loop can be lifted as a path $c$, with origin at $y$ and lying in $\pi^{-1} (\overline V)$; the endpoint of $c$ is thus an element of 
$\pi^{-1} (\bar y) = \Gamma y$ which does not depend on the choice of the loop $\bar c$ in the homotopy class $\beta$, because 
two loops which are homotopic (in $\overline Y$) lift as two paths with the same endpoints in the universal cover $Y$. Denoting by $\g y$ the endpoint of
$c$, it follows that $\g y$ (and thus $\g$, for $\Gamma$ acts without fixed point) does not depend on the choice of the loop $\bar c$ in the homotopy class 
$\beta$. Moreover, as $\bar c = \pi \circ c$, then $c$ lies in $\pi^{-1} (\overline V)$, thus its endpoints $y$ and $\g y$ belong to the same connected 
component $V^i$ of $ \pi^{-1} (\overline V)$, proving that $ \g (V^i) \cap V^i \ne \emptyset$, thus (by \eqref{intersectioncomp}) that $ \g(V^i) = V^i $
and then that $\g$ is an element of $\widehat\Gamma^i_V$.
Mapping $\beta$ onto this $\g$, we obtain a well-defined map $\phi_y :  i_* \big(\pi_1 (\overline V , \bar y ) \big) \f \widehat\Gamma^i_V$.
One verifies easily that, by construction, for every $\g \in \widehat\Gamma^i_V$ and every $ \beta \in  i_* \big(\pi_1 (\overline V , \bar y ) \big)$,
one has $\phi_y \circ \psi_y (\g) = \g$ and $\psi_y \circ \phi_y (\beta) = \beta$, and this ends the proof.
\end{itemize}
\end{proof}

\subsubsection{On the topology of thin subsets}\label{topominces}

On any metric space $(Y,d )$, for any proper action, by isometries, without fixed point, of any group $\Gamma$ on $(Y,d)$, for any
given $ r \in \, ] 0 , +\infty [$, we recall that the $ r$-thin subset of $Y$ (resp. of $\overline Y = \Gamma \backslash Y$) is the open set $Y_{r} =
\{ y \in Y : \text{\rm sys}_\Gamma (y)  < r\}$ (resp. its image $\overline Y_r $ by the quotient-map $\pi : Y \to \overline{Y}= \,\Gamma \backslash Y$)
and that the (torsion-free) $ r$-thin subset of $Y$ (resp. of $\overline Y = \Gamma \backslash Y$) is the open set $Y^{\diamond}_{r} =
\{ y \in Y : \text{\rm sys}^{\diamond}_\Gamma (y)  < r\}$ (resp. its image $\overline Y_r^{\diamond}$ by the quotient-map $\pi : Y \to \overline{Y}= \
\Gamma \backslash Y$).\\
Recalling the Definitions \ref{Gammaepsilon}, we denote by $\Sigma_r (x)$ (resp. by $\widehat{\Sigma}_r (x)$) the finite set of $\g \in \Gamma^*$ 
satisfying $d(x , \g x) \le  r$ (resp.  $d(x , \g x) < r$), and $\Gamma_r (y)$ (resp. $\widehat{\Gamma}_r (x)$) the subgroup of $\Gamma$ generated by 
$\Sigma_r (x)$  (resp. by $ \,\widehat{\Sigma}_r(x)$).\\
For the definitions of geodesics and local geodesics, see section \ref{notations}.

\begin{remark}\label{banal} The functions $ y \mapsto \sys_\Gamma (y) $ and $ y \mapsto \sys^{\diamond}_\Gamma (y) $ are both invariant by the action 
of $\Gamma$, there thus exists two functions $\overline{\sys}_\Gamma $ and $\overline{\sys}^{\diamond}_\Gamma $ from $\overline{Y} =
\Gamma \backslash  Y $ to $\R^+$ such that $ \sys_\Gamma = \overline{\sys}_\Gamma \circ \pi$ and $ \sys^{\diamond}_\Gamma = 
\overline{\sys}^{\diamond}_\Gamma \circ \pi$. Consequently, for every $r > 0$, the sets $Y_r$ and $Y^{\diamond}_r$ are stable under the action of 
$\Gamma$, which means on the one hand that $ Y_r = \pi^{-1} (\overline Y_r )$ and $ Y^{\diamond}_r = \pi^{-1} (\overline Y^{\diamond}_r )$, thus that 
$\overline Y_r = \Gamma \backslash   Y_r$ and $\overline Y^{\diamond}_r = \Gamma \backslash   Y^{\diamond}_r$, and on the other hand that
$\overline Y_r $ (resp. $\overline Y_r^{\diamond}$) coincide with the subset of $\bar y \in \overline Y$ such that $ \overline{\sys}_\Gamma (\bar y)  < r$
(resp. such that $ \overline{\sys}^{\diamond}_\Gamma (\bar y)  < r$)
\end{remark} 

The proof of this Remark is straightforward, the fact that $g (Y_r^{\diamond}) = Y_r^{\diamond}$ (for every $g \in \Gamma$) coming from the fact
that $\Gamma^{\diamond} = g \Gamma^{\diamond} g^{-1}$ because $\g$ is torsion-free iff $g \g g^{-1}$ is torsion-free.

\medskip
Fix arbitrarily some $ r > 0$, denote by $\left\{ Y_r^i :  i\in I \right\}$ (resp. by $\{ \overline{Y}_{r}^j :  j\in J \}$) the set
of connected components of $ Y^\diamond_r = \pi^{-1} (\overline Y^\diamond_r)$ (resp. of $\overline{Y}^\diamond_{r}$). The continuous mapping $\pi : Y \to \overline Y$ maps each 
connected component $Y_r^i $ of $Y^\diamond_r $ into some connected component of $\overline{Y}^\diamond_{r} $, that we shall denote 
by $\overline{Y}_{r}^{k(i)} $; 
this defines a map $ i \mapsto k(i)$ and moreover proves that, for every $(i,j) \in I \times J$,
\begin{equation}\label{intersect}
Y_r^i \cap \pi ^{-1}(\overline{Y}_{r}^{j}) \ne \emptyset \iff \pi (Y_r^i) \cap \overline{Y}_{r}^{j} \ne \emptyset \iff 
\pi (Y_r^i) \subset \overline{Y}_{r}^{j} \iff  i \in k^{-1}(\{j\}) \ .
\end{equation}
A consequence is that, if $i \in k^{-1}(\{j\}) $ (resp. if $i \notin k^{-1}(\{j\}) $) then $Y_{r}^i \subset \pi ^{-1}(\overline{Y}_{r}^j) $ (resp. $Y_{r}^i 
\cap \pi ^{-1}(\overline{Y}_{r}^j) = \emptyset $), which implies that
$$\pi ^{-1}(\overline{Y}_{r}^j) = Y_r \cap  \pi ^{-1}(\overline{Y}_{r}^j) = \bigcup_{i \in I} Y_{r}^i \cap \pi ^{-1}(\overline{Y}_{r}^j) 
= \bigcup_{i \in k^{-1}(\{j\})} Y_{r}^i \cap \pi ^{-1}(\overline{Y}_{r}^j) =  \bigcup_{i \in k^{-1}(\{j\})} Y_{r}^i \ .$$
It follows that $\{Y_r^i : i \in k^{-1}(\{j\})\}$ is the set of connected components of $\pi ^{-1}(\overline{Y}_{r}^j)$.
Hence, one may apply the results of section \ref{topogeneral} where we replace $\overline V $ by $ \overline{Y}_{r}^j$ and $V$ by 
$ \bigcup_{i \in k^{-1}(\{j\})} Y_{r}^i $ and where we denote by $ \widehat\Gamma^i_r $ the subgroup of the elements $\g \in \Gamma$ which
verifies $\g (Y_{r}^i) = Y_{r}^i$, these results allow to define the quotient spaces $\widehat Y^i:= \widehat\Gamma^i_r \backslash Y$ and 
$\widehat Y^i_r :=    \widehat\Gamma^i_r \backslash Y^i_r $ and the quotient maps $ \widehat \pi_i : Y \to  
\widehat Y^i $ and $ \widehat p_i: Y^i_r \to  \widehat Y^i_r $.
As in subsection \ref{topogeneral}, let $j_i$ and $\bar j_i$ be the respective inclusions $ Y^i_r \hookrightarrow Y$ and $ \overline{Y}_{r}^{k(i)} 
\hookrightarrow \overline{Y}$; the map $j_i$ induces (by quotients) a map $ j'_i : \widehat\Gamma_r^i \backslash Y^i_r \hookrightarrow \
\widehat\Gamma_r^i \backslash Y $, that we also consider as an inclusion; then we have (by construction) $\widehat \pi_i \circ j_i = j'_i \circ \widehat p_i$
and $\widehat p_i$ may be considered as a restriction (on both sides) of $\widehat \pi_i $.
As $\pi (Y_{r}^i) \subset \overline{Y}_{r}^{k(i)}$, one may define $p_i : Y^i_r \f \overline{Y}_{r}^{k(i)} $ as the restriction (on both sides) of $\pi$ 
(namely $\pi \circ j_i = \bar j_i \circ p_i$).
As $\pi$ is constant on the orbits of $ \widehat\Gamma_r^i $, the map $p_i$ is also constant on the orbits of $ \widehat\Gamma_r^i $, and thus
the maps $\pi $ and $p_i$ respectively go down as maps $ \bar \pi_i  \,: \,  \widehat\Gamma^i_r \backslash Y \f \Gamma \backslash Y$ and $ \bar p_i  \,: \,  \widehat\Gamma^i_r \backslash Y^i_r  \f \overline{Y}_{r}^{k(i)} $ satisfying $ \bar \pi_i \circ \widehat \pi_i = \pi$ and $ \bar p_i \circ \widehat p_i = p_i $.
We shall now prove that $\bar p_i$ may be considered as the restriction (on both sides) of $\bar \pi_i $, in fact
$$ \forall y \in Y^i_r \ \ \ \  \bar j_i \circ \bar p_i  \big(\widehat p_i (y) \big) = \bar j_i \circ  p_i (y) =  \pi \circ j_i (y) = \bar \pi_i \circ \widehat \pi_i \circ j_i (y) = 
\bar \pi_i \circ j'_i \circ \widehat p_i (y)   = \bar \pi_i \circ j'_i \big( \widehat p_i (y) \big) \ ,$$
and (as $\widehat p_i$ is surjective) this proves that $\bar j_i \circ \bar p_i = \bar \pi_i \circ j'_i$.\\
All these definitions and properties are summarized in the diagram below:

 \vspace{-5mm}

\large
\begin{equation}\label{diagramme2} 
\xymatrix{
\Gamma \curvearrowright  \hspace{-12mm}            & \hspace{2mm} Y \hspace{2mm}  \ar[r]^{\hspace{-10mm} \widehat \pi_i}   \ar@/^2pc/[rr]^{\pi}                               & \hspace{4mm}   \widehat Y^i =   \widehat\Gamma_r^i \backslash Y  \hspace{1mm}  \ar[r]^{\bar \pi_i} 
                                 &  \hspace{5mm} \overline{Y}  = \Gamma \backslash Y  \hspace{2mm}  \\
\widehat\Gamma^i_r \curvearrowright  \hspace{-10mm}    &  \hspace{2mm} Y^i_r  \hspace{2mm}  \ar@{^{(}->}[u]  \ar[r]^{\hspace{-10mm}  \widehat p_i}     \ar@/_2pc/[rr]^{p_i}        & \hspace{4mm} \widehat Y^i_r  =   \widehat\Gamma_r^i \backslash Y^i_r  \hspace{1mm}  \ar[r]^{\bar p_i}  \ar@{^{(}->}[u]    
                                 &    \hspace{4mm} \overline{Y}_r^{k(i)}     \hspace{5mm} \ar@{^{(}->}[u]  
 }
\end{equation}

\normalsize

The following Proposition summarizes and completes the above prologue. Notice that, when $r \le \inf_{y \in Y}\sys^\diamond_\Gamma (y) $, then 
$ Y^\diamond_r$, $\overline{Y}^\diamond_{r}$ and $I$ are empty and the following Proposition is then trivial.

\begin{prop}\label{topotubes1} 
With the above notations, for every proper action (without fixed point and by isometries) of any group $\Gamma$ on any metric space $ (Y , d) $, 
\begin{itemize}
\item[(o)] for every connected component $\overline{Y}_{r}^j$ of $\overline{Y}^\diamond_{r}$, the set of connected components of $\pi ^{-1}(\overline{Y}_{r}^j)$ 
is $\left\{Y_r^i : i \in k^{-1} \big(\{j\}\big)\right\}$;

\item[(i)] for every pair $ Y_{r}^i  , \  Y_{r}^{i'}$ of connected components of $ Y^\diamond_r$, the indices $ k(i)$ and $ k(i')$ are equal if and only if there exists 
$\g \in \Gamma$ such that $\g (Y_{r}^i ) = Y_{r}^{i'}$, and then $ \widehat{\Gamma}_{r}^{i'} = \g  \widehat{\Gamma}_{r}^{i} \g^{-1} $;

\item[(ii)] for every $i \in I$, $\pi (Y_{r}^i) = \overline{Y}_{r}^{k(i)}$ and the maps $\pi$, $\widehat \pi_i$, $p_i$ and $\widehat p_i $ are locally 
isometric coverings (and thus they are open mappings);

\item[(iii)] for every $i \in I$, the map $\bar p_i$  is a homeomorphism from $\widehat Y^i_r$ onto $\overline{Y}_{r}^{k(i)}$, which is a local isometry 
(and thus preserves the path-lengths);

\item[(iv)] Il $(Y,d)$ is a length space, for every $ i \in I$, $\widehat{\Gamma}_r^i$ contains the group generated by all the torsion-free elements of the
union of all the sets $\widehat \Sigma_{r} (y)$ for all the $y$'s in $Y^i_{r}$;

\item[(v)] If $(Y, d)$ is geodesic, each connected component $\overline{Y}_{r}^{j}$ of $\overline{Y}^\diamond_{r}$ (whose closure is compact or which satisfies 
$\liminf \overline{\sys}^\diamond_\Gamma (\bar y) \ge r$ when $\bar y$ goes to infinity in $\overline{Y}_{r}^{j}$)
contains a non homotopically trivial closed local geodesic, whose length (or period) is equal to the minimum (on $\overline{Y}_{r}^{j}$)
of the function $ \bar y \mapsto \overline{\sys}^\diamond_\Gamma (\bar y)$.
\end{itemize}
\end{prop}

\begin{proof} The point (o) has been proved in the beginning of the subsection \ref{topominces}. The property (o) and the fact that $\Gamma $ acts properly 
without fixed point prove that the hypotheses of Lemma \ref{lemmetopologique} are verified, where we replace the open subset $\overline V$ of the 
Lemma \ref{lemmetopologique} by any connected component $\overline{Y}_{r}^{j}$ of $\overline{Y}_r$ and the set of connected components of 
$\pi ^{-1}(\overline V)$ by the set $\left\{Y_r^i : i \in k^{-1} \big(\{j\}\big)\right\}$ of the connected components of $\pi ^{-1} \big(\overline{Y}_{r}^j \big)$;
Lemma \ref{lemmetopologique} (iii) then implies that
$$ i, i' \in k^{-1}(\{j\}) \implies  \exists \g \in \Gamma \  \text{ such that } \   \g (Y_{r}^i ) = Y_{r}^{i'} 
\  \text{ and } \  \widehat{\Gamma}_{r}^{i'} = \g  \widehat{\Gamma}_{r}^{i} \g^{-1}  \ ,$$
as the converse implication is trivially true, this proves the point (i).\\
The point (ii) is an immediate application of Lemma \ref{lemmetopologique} (i) and (ii);
a consequence of the lemma \ref{lemmetopologique} (ii) is moreover that, for every $ j \in J$ and every $ i \in k^{-1}(\{j\})$, the restrictions
$p_i$ and $\widehat p_i$ to $ Y_{r}^i  $  are locally isometric coverings from $ Y_{r}^i  $ onto $\overline{Y}_{r}^{j}$  and $\widehat Y^i_r$ respectively.\\
The lemma \ref{lemmetopologique} (iv) implies that the map $\bar p_i$ is a homeomorphism from $\widehat Y^i_r$ onto $\overline{Y}_{r}^{k(i)}$, which is a local isometry and preserves the path-lengths. This ends the proof of point (iii).\\
\emph{Proof of (iv) :}
Suppose now that $(Y,d)$ is a length space, for any $ i \in I$ and any point $y$ of the corresponding connected component $ Y_r^i$ of $ Y^\diamond_r$, 
let us consider any torsion-free $\g \in \widehat{\Sigma}_r (y) $; it verifies (by definition) $ d(y, \g y) <  r$ and as $(Y,d)$ is a length space, for every 
$ 0 < \eta < r - d(y, \g y) $, there exists a path $c$ which joins the points $y $ and $ \g y $ such that any point $u$ of this path satisfies
$ d(y ,  u) + d(u, \g y)  < d(y, \g y) + \eta < r$, consequently, we get
$$d (u, \g u ) \le d(u, \g y)  + d(\g y , \g u) =  d(y ,  u) + d(u, \g y) \le \text{length of } c  <   r \ ;$$
from this and from the fact that $\gamma$ is torsion-free, we deduce that $\sys^\diamond_\Gamma (u) < r$, and thus that the image of $c$ is entirely contained in $Y^\diamond_r$; it follows that $\g y$ belongs to the same connected component of $Y^\diamond_r$ as $y$, namely it belongs to $Y_r^i$.
From this, we infer that $\g (Y_r^i) \cap  Y_r^i \ne \emptyset$ and (by \eqref{intersectioncomp}) that $\g (Y_r^i) = Y_r^i$, hence that every torsion-free 
$\g \in \widehat{\Sigma}_r (y) $ belongs to $\widehat \Gamma_{r} (y)$ and (iv) is proved.\\

\emph{Proof of (v) :} Suppose now that $(Y,d)$ is geodesic:
by definition of $ \overline Y^\diamond_r$, for any of its connected components $ \overline Y_r^j$ which satisfies the hypotheses of (v),
the function $\overline{\text{\rm sys}}^\diamond_\Gamma$ (defined at Remark \ref{banal}) is strictly smaller 
than $r$ on $ \overline Y_r^j$, and at least equal to $r$ on the boundary of $ \overline Y_r^j$ or at infinity; 
thus this function (when restricted to $ \overline Y_r^j$) attains its minimum at some point $\bar y_j \in \overline Y_r^j$. 
Fix any $i \in k^{-1}(\{j\})$; by (ii)
there exists some point $y_i \in \pi^{-1}(\{\bar y_j\}) \cap Y_r^i$ and the remark \ref{banal} (together with the fact that $\pi (Y_{r}^i) = 
\overline{Y}_{r}^{j}$ by (ii)) then implies that
$$\sys^\diamond_\Gamma (y_i) = \overline{\sys}^\diamond_\Gamma (\bar y_j) = \inf_{\bar y \in \overline Y_r^j} \overline{\sys}^\diamond_\Gamma (\bar y) 
= \inf_{y \in Y_r^i} \overline{\sys}^\diamond_\Gamma \circ \pi (y) = \inf_{y \in Y_r^i} \sys^\diamond_\Gamma (y)\ .$$
From this, from the properness of the action we infer the existence of a torsion-free element $\sigma \in \widehat{\Sigma}_r (y_i)$ satisfying 
$d( y_i , \sigma \, y_i ) = \sys^\diamond_\Gamma (y_i) = \Min_{y \in Y_r^i} \big( \sys^\diamond_\Gamma (y) \big)$, and property (iv) implies that 
$\sigma \in \widehat{\Gamma}_r^i$.
Denote by $[ y_i , \sigma \, y_i ]$ any (length-minimizing) geodesic between the points $y_i$ and $ \sigma\, y_i$, as $\sigma$ is torsion-free, 
every point $u \in [ y_i , \sigma \, y_i ] $ verifies 
$$  \sys^\diamond_\Gamma (u) \le d( u , \sigma \,u) \le d(u, \sigma \, y_i) + d(\sigma \, y_i , \sigma u) = d(u, y_i) + d(u , 
\sigma \, y_i ) = d(y_i, \sigma \, y_i) = \sys^\diamond_\Gamma (y_i) < r\ ;$$
a first consequence is that $ u \in  Y_r^i$, thus that $ \sys^\diamond_\Gamma (u) \ge \sys^\diamond_\Gamma (y_i)$ a second consequence is that
$ \sys^\diamond_\Gamma (u) \le \sys^\diamond_\Gamma (y_i)$; hence $\sys_\Gamma (u) = \sys_\Gamma (y_i)$. It follows that 
all the above inequalities are equalities and thus the union of the two segments $ [u, \sigma \, y_i]$ and $ [\sigma \, y_i , \sigma u]$ 
is a (length-minimizing) geodesic between the points $u$ and $\sigma \, u$.
Denoting by $[\sigma^k y_i , \sigma^{k+1} y_i ]$ the image by $\sigma^k$ of 
$ [ y_i , \sigma \, y_i ] $, it follows that the union $ [ y_i , \sigma \, y_i ] \cup [\sigma \,y_i , \sigma^{2} y_i ]$ is a local geodesic, which is length-minimizing 
on any subsegment of length $T := \sys_\Gamma (y_i)$. In the same way, we prove that 
$$  c := \ldots \cup [\sigma^{-p} y_i , \sigma^{-p+1} y_i ] \cup \ldots \cup[ y_i , \sigma \, y_i ]\cup \ldots \cup [\sigma^{p-1} y_i , \sigma^{p} y_i ] 
\cup [\sigma^{p} y_i , \sigma^{p+1} y_i ]\cup \ldots $$
is a local geodesic entirely contained in $Y_r^i$, which is length-minimizing on any subsegment of length $T $ and verifies $ c(t+ T) = 
\sigma \big(c(t) \big)$ when parametrized by length. As $\pi$ is a locally isometric covering (by (ii)), the image $\pi \circ c$ of
the local geodesic $c$ is a local geodesic of $\overline{Y}$, which is $T$-periodic (for $ \pi \circ c (t + T) = \pi \big( \sigma \circ c(t)\big) = \pi \circ c (t)$),
and which is entirely included in $\pi ( Y_r^i) = \overline Y_r^j$ (equality proved by (ii)).
Moreover the period $T$ of this closed geodesic  coincides with $\sys^\diamond_\Gamma (y_i) = \overline{\sys^\diamond}_\Gamma (\bar y_j)$, which is actually the minimal 
value on $\overline{Y}_{r}^{j}$ of the function $ \bar y \mapsto \overline{\sys}^\diamond_\Gamma (\bar y)$.
\end{proof}

Given any parameters $\delta_0, H_0, D_0, \e'_0> 0$, we recall the universal constants $\alpha_0 := \alpha_0 (\delta_0, H_0, D_0)$,  $\alpha'_0 := \alpha'_0 (\delta_0, H_0)$ and $r_0 :=r_0 (\delta_0, \e'_0)$ defined at \eqref{universalconstants1}, \eqref{universalconstants2} and \eqref{universalconstants3} respectively. We then get the

\begin{theorem}\label{topotubes} \emph{(Thin subset have simple topology)}
Given any parameters $\delta_0, H_0, D_0, \e'_0 , H> 0$, for every element $\Gamma$ of $\text{\rm Hyp}_{\rm action} (\delta_0, H_0, D_0)$ (resp. of 
$\text{\rm Hyp}_{\rm sub} (\delta_0, H_0)$, resp. of $\text{\rm Hyp}_{\rm thick}  (\delta_0, \e'_0)$), for every $r>0$ such that 
$r \le \frac{\alpha_0}{H}$ (resp. such that $r \le  \frac{\alpha'_0}{H}$, resp. such that $r \le  \frac{r_0}{H}$), any proper action (without fixed point, 
by isometries preserving the measure) of 
$\Gamma$ on a metric measured space $ (Y , d, \mu) $ whose entropy is bounded from above by $H$ verifies,
\begin{itemize}
\item[(i)] for every $ i \in I$, $\widehat{\Gamma}_r^i := \{ \g \in \Gamma : \g (Y^i_{r}) = Y^i_{r}\}$ is a virtually cyclic subgroup of $\Gamma$; if moreover 
$(Y,d)$ is a length space, then $\widehat{\Gamma}_r^i $ contains the group generated by the torsion-free elements of 
$\ \bigcup_{y \in Y^i_{r}} \ \widehat \Sigma_{r} (y)$.
 \item[(ii)] If $Y$ is locally path-connected and simply connected then, for every connected component $ Y_r^i$ of $ Y^\diamond_r$ and every $ y \in Y_r^i$, 
if we denote by $ i_*\,: \pi_1 \left(\overline{Y}_{r}^{k(i)} ,\pi (y) \right) \to \pi_1 \big(\overline Y , \pi (y) \big)$ the morphism induced by the inclusion 
$ i : \overline{Y}_{r}^{k(i)} \hookrightarrow \overline Y$, the subgroup $\widehat\Gamma^i_r$ is identified with 
$i_* \pi_1 \left(\overline{Y}_{r}^{k(i)} ,\pi (y) \right) $ by the canonical isomorphism $\Gamma \f  \pi_1 \big(\overline{Y} , \pi (y) \big) $ and 
$i_* \pi_1 \left(\overline{Y}_{r}^{k(i)} ,\pi (y) \right) $ is a virtually cyclic  group.
\end{itemize}
\end{theorem}

\begin{proof}[Proof of (i)] As $\Gamma$ is non virtually cyclic by definition of $\text{\rm Hyp}_{\rm action} (\delta_0, H_0, D_0)$ 
(resp. of $\text{\rm Hyp}_{\rm sub} (\delta_0, H_0)$, resp. of $\text{\rm Hyp}_{\rm thick}  (\delta_0, \e'_0)$), as it admits a 
proper action (by isometries) on some $\delta_0$-hyperbolic space $ (X, d_0)$ such that no element acts as a parabolic isometry (by Lemma \ref{torsionfreehyperbolic}), as $\Gamma_{r} (y)$ is virtually cyclic at every point of $Y^\diamond_r$ by the conclusion (i) (resp. (ii), resp. (iii)) of 
Theorem \ref{transportnil}, we may apply Proposition \ref{transsystprep} (1), which proves the existence of some maximal virtually cyclic group $G^i$ 
which contains all the sets $\Sigma_{r} (y)$ for all the $y \in Y_r^i$.
Let us now show that $ \widehat{\Gamma}_r^i \subset G^i$: for every $\g \in \widehat{\Gamma}_r^i $ and any $y \in Y_r^i$, $y$ and $\g y$ are located 
in the same connected component $Y_r^i$ of $Y^\diamond_r$, 
thus $G^i$ contains $\Sigma_{r} (y)$ and $\Sigma_{r} (\g y) = \g\,\Sigma_{r} (y)\, \g^{-1}$.
Let $\sigma$ be any torsion-free element of $\Sigma_{r} (y)$ (such an element exists for $\sys^\diamond_\Gamma (y) < r$), then $\sigma$ and
$\g \sigma \g^{-1}$ both belong to $G^i$ and the subgroup $\langle \sigma , \g \sigma \g^{-1} \rangle$ is thus virtually cyclic; as $\ell (\sigma) > 0$ by 
Lemma \ref{torsionfreehyperbolic} (ii),  it follows from this and from Corollary \ref{elementaryaction} (i) that $\langle \sigma , \g  \rangle$ is virtually 
cyclic. From this and from the fact that $G^i$ is the maximal virtually cyclic subgroup which contains $\sigma$, we infer that $\langle \sigma , \g  \rangle
\subset G^i$, thus that $\g \in G^i$ for every $\g \in \widehat{\Gamma}_r^i$; this proves that $\widehat{\Gamma}_r^i \subset G^i$.
We conclude that $\widehat{\Gamma}_r^i$ is virtually cyclic.\\
If now $(Y,d)$ is a length space, the fact that $\widehat{\Gamma}_r^i $ contains the group generated by the torsion-free elements of 
$\ \bigcup_{y \in Y^i_{r}} \ \widehat \Sigma_{r} (y)$ is a direct consequence of Proposition \ref{topotubes1} (iv).

\smallskip
\emph{Proof of (ii).} We apply Lemma \ref{lemmetopologique} (v), where we replace the open set $\overline V \subset \overline Y $ of Lemma 
\ref{lemmetopologique} by any connected component $\overline{Y}_{r}^{j}$ of $\overline{Y}_r$ and the set $\left\{ V ^i  : i \in I \right\}$ of 
connected components of $\pi^{-1} (\overline V)$ by the set $\{Y_r^i : i \in k^{-1}(\{j\})\}$ of connected components of $\pi ^{-1}(\overline{Y}_{r}^j)$;
for each $i \in k^{-1}(\{j\})$, the set $\widehat{\Gamma}_V^i$ of elements of $\Gamma $ which 
stabilize $V_i$ is thus replaced by the set $\widehat{\Gamma}_r^i$ of elements of $\Gamma $ which stabilize $Y_r^i$.
Property (ii) is then a direct consequence of Lemma \ref{lemmetopologique} (v).\\
Now property (i) implies that $i_* \pi_1 \left(\overline{Y}_{r}^{k(i)} ,\pi (y) \right) \simeq \widehat{\Gamma}_r^i$ is virtually cyclic.
\end{proof}

\subsubsection{Examples and counter-examples}\label{examplescounter}

In the following examples we refer to the notations of subsection \ref{topominces}: given any $a, b$ such that $0 < a < b$, and defining 
$\delta_0 = \frac{\ln 3}{a}$ and $H_0 = (n-1) b$, every group $\Gamma$ under consideration in these examples is an element of the set of 
groups $\text{\rm Hyp}_{\rm thick} (\delta_0, \e'_0)$ (resp. of $\text{\rm Hyp}_{\rm action} (\delta_0, H_0, D_0)$) introduced in 
Definition \ref{systoleaction} (resp. in Definition \ref{hyperbolicaction}), the values of the other parameters $ \e'_0,  D_0$ being specified further.
In all the following examples (except the fourth one), the group $\Gamma$ is torsion-free, because it is isomorphic to a subgroup of the fundamental group of 
a complete $n$-dimensional Riemannian manifold whose sectional curvature is $\le -a^2$ and whose injectivity radius is $\ge \e'_0/2$.
(resp. to a subgroup of the fundamental group of a compact $n$-dimensional
Riemannian manifold whose sectional curvature $\sigma$ verifies\footnote{Notice that the assumption $- b^2 \le \sigma \le - a^2 $ and the comparison 
Theorems prove that the Riemannian universal covering of this manifold is a $\delta_0$-hyperbolic space, with Entropy $\le H_0$, where $\delta_0
= \frac{\ln 3}{a}$ by Proposition 1.4.3 page 12 of \cite{CDP} and $H_0 = (n-1) b$ by Bishop-Gromov Comparison Theorem.} $- b^2 \le \sigma \le - a^2 $ and whose diameter is bounded from above by $D_0$).\\
We shall consider (in these examples) proper actions (by isometries without fixed point) of the same group $\Gamma$ on another simply connected 
Riemannian manifold $(Y, \tilde g)$ whose entropy is bounded from above by some constant $H$. 
Except in Example (4), as $\Gamma$ is torsion-free, the two functions $\text{\rm sys}^{\diamond}_\Gamma (\cdot )$ and $\text{\rm sys}_\Gamma (\cdot ) $ coincide on the whole 
of $Y$; for this reason, for every positive $r \le \frac{r_0}{H}$ (resp. $r \le \frac{\alpha_0}{H}$) where $r_0 :=r_0 (\delta_0, \e'_0)$ is defined at \eqref{universalconstants3} (resp. where $\alpha_0 := \alpha_0 (\delta_0, H_0, D_0)$ is defined at \eqref{universalconstants1}), the $ r$-thin subset 
$Y_{r}$ and the torsion-free $ r$-thin subset $Y^{\diamond}_{r}$ coincide with $\{ y \in Y : \text{\rm sys}_\Gamma (y)  < r\}$ and their images by 
$\pi : (Y,g) \f  \overline Y = \Gamma \backslash Y$ also coincide, i. e. $\overline Y_r^{\diamond} = \overline Y_r = \{ \bar y \in \overline Y : 
\overline{\text{\rm sys}}_\Gamma (y)  < r\}$.
\normalsize

\begin{itemize}

\item[(1)] When $(M, g_0)$ is a complete $n$-dimensional manifold with sectional curvature $ \sigma \equiv -1$ and 
with injectivity radius
$\ge \frac{\e'_0}{2} > 0$, its universal cover $(\widetilde M , \tilde g_0) \simeq (\mathbb H^n , can.)$ is $\ln 3$-hyperbolic by Proposition 1.4.3 page 12 
of \cite{CDP} and the systole of the action (by deck-transformations) of its fundamental group $G$ on its universal cover $(\widetilde M , \tilde g_0) $ is
bounded from below by $\e'_0$. It follows that $G$ and every non cyclic\footnote{As $G$ acts on $\mathbb H^n $ without fixed point, it is 
torsion-free by Lemma \ref{fixedpoint}, hence every virtually cyclic subgroup of $G$ is cyclic by Lemma \ref{virtuelcyclique}.} subgroup $\Gamma$ of $G$ belong to $\text{\rm Hyp}_{\rm thick} (\ln 3, \e'_0)$.

Considering now the induced action of $\Gamma$ on $(Y,d, \mu) := (\widetilde M , \tilde g_0 , dv_{\tilde g_0})$, the following results are classical in this 
case (see \cite{Ba-Gr-Sc}): each connected component 
$Y^i_{r}$ (resp. $\overline Y^k_{r}$) of the $ r$-thin subset $Y^{\diamond}_{r} = Y_{r} \subset Y$ (resp. of the $ r$-thin subset 
$ \overline Y^{\diamond}_r = \overline Y_r \subset \overline Y$) is 
then diffeomorphic to $ \R \times \mathbb B^{n-1}$ (resp. to $\mathbb S^1 \times \mathbb B^{n-1}$), the groups $\widehat{\Gamma}_r^i := 
\{ \g \in \Gamma : \g (Y^i_{r}) = Y^i_{r}\}$ and $ \pi_1 \left(\overline Y^k_{r} \right)$ being isomorphic to $\Z$ (compare with Theorem \ref{topotubes}
and Theorem \ref{tubelong} (i) and (v)).
Moreover, as proved in Proposition \ref{topotubes1} (v), each bounded connected component $\overline Y^k_{r}$ contains some minimizing closed periodic geodesic 
$\bar c_k$ (of length $\e_k$) and a tubular neighbourhood of this geodesic of radius $ C(n) \,\ln \left( \frac{r}{\e_k}\right)$, and this is consistent with the 
general estimate of the radius of this tubular neighbourhood given by Theorem \ref{tubelong} (iii).\\
The generalization of these results to  complete Riemannian manifolds whose sectional curvature verifies $- b^2 \le \sigma \le - a^2 $ is also classical (see \cite{Ba-Gr-Sc}).
\item[(2)] One application of the results of the present section is the following

\begin{corollary} For any differentiable manifold $X$ which admits a complete Riemannian metric $g_0$ with sectional curvature $ \sigma \le -a^2 < 0$ and injectivity
radius $\ge i_0 > 0$, for every other complete Riemannian metric $g$ on $X$ with bounded Entropy\footnote{We recall that the Entropy of $(X,g)$ is the entropy of its universal cover $ \big(\widetilde X, d_{\tilde g},\, dv_{\tilde g} \big)$.}, all the Margulis' properties proved in the present section are still 
valid for $(X,g)$, viewed as the quotient of the measured metric space  $ \big(\widetilde X, d_{\tilde g},\, dv_{\tilde g} \big)$ by the action of the fundamental group $\Gamma$ of $X$.
\end {corollary}

\begin{proof} Let us first consider the action of $\Gamma$ on the metric space $\big(\widetilde X, d_{\tilde g_0} \big)$: as mentioned above, 
$\big(\widetilde X , d_{\tilde g_0} \big)$ is $\delta_0$-hyperbolic with $\delta_0 := \frac{1}{a} \ln 3$ and, the global  systole 
of the action of $\Gamma$ on this space being bounded from below by $ \e'_0 := 2 i_0$, every $\g \in \Gamma^*$ verifies $\ell (\g) \ge 2 i_0$. This implies that
$ \Gamma$ is a torsion-free element of $\text{\rm Hyp}_{\rm thick} (\delta_0, \e'_0)$.
As the Entropy of the measured metric space $ \big(\widetilde X , d_{\tilde g} \,,\, dv_{\tilde g} \big)$ is bounded above by a given constant (say $H$),
we may apply all the results of the present section \ref{transport} to the action of $\Gamma$ on $ \big(\widetilde X , d_{\tilde g} \,,\, dv_{\tilde g} \big)$.\end{proof}
\end{itemize}

In the general case of subgroups of hyperbolic groups (or of groups which admit an action on some Gromov-hyperbolic space 
$(X,d_0)$), considering their actions on any measured metric space $(Y,d,\mu)$, the results of subsections \ref{topominces} and \ref{tubessimples} aim to generalize the properties of hyperbolic or negatively curved 
manifolds mentioned in the reference-example (1). 
In comparison, the results of subsections \ref{topominces} and \ref{tubessimples} seem weaker: in fact, the strong version (valid in the negatively curved 
case, see Example (1) above) cannot be expected in this general setting, as proved by the following examples:

\begin{itemize}

\item[(3)] {\em  Comparing to the reference-example (1), in the general case, the connected components $\overline Y^k_{r}$ of the $ r$-thin subset $\overline Y^{\diamond}_{r}$ 
are no longer homeomorphic to à 
$\mathbb S^1 \times \mathbb B^{n-1}$, moreover the fundamental group of $\overline Y^k_{r}$ is generally not virtually cyclic: one cannot expect
a better result than the virtual cyclicity (proved in Theorem \ref{topotubes}) of the image of the fundamental group of $\overline Y^k_{r}$ by the homomorphism $i_*$ between fundamental
groups induced by the inclusion $ i : \overline{Y}_{r}^{k} \hookrightarrow \overline Y$:}

\smallskip
Indeed, let $\Sigma$ be a compact surface with genus $\ge 2$; let $g_1$  be a first hyperbolic Riemannian metric on $\Sigma$ whose injectivity radius is bounded 
from below by a constant $\e'_1 > 0$; the fundamental group $\Gamma$ of $\Sigma$ is then a torsion-free element of $\text{\rm Hyp}_{\rm thick} 
(\delta_0, \e'_0)$, where $\delta_0 = \ln 3$ and $\e'_0 = 2 \e'_1$. Then the corresponding universal constant $r_0 = r_0 (H_0 , \e'_0)$ defined in \eqref{universalconstants3} is equal to $\dfrac{2 \e'_1 \ln 2}{13 \ln 3 + 8 \e'_1} < 1$. 

Fix any positive value $r \le r_0 $. Let $g_0$  be another hyperbolic Riemannian metric on $\Sigma$ with big diameter and which admits a small closed geodesic $c$ of length $\e << r$;  
the connected component of the $r$-thin subset which contains $c$ is then a tubular neighbourhood $U$ of $c$ with radius $L \ge \ln (r/\e)$. 
Denote by $\varrho (\cdot ) = d_{g_0}(\cdot ,c)$ the function distance to this geodesic, choose points $\bar y_1 , \ldots , \bar y_N \in U$ such that 
$\inj_{g_0} (\bar y_i) \le r/4$ and  $|\varrho (\bar y_i)-\varrho (\bar y_j)|\ge 5\varepsilon ^2$ when $i\ne j$; define $V = 
\bigcup_{i=1}^N  B_{g_0}(\bar y_i,\varepsilon ^2)$ and $W = \bigcup_{i=1}^N  B_{g_0} (\bar y_i,2\varepsilon ^2)$ 
and construct a $C^\infty$ function $f : \Sigma \to [1,+\infty[$ such that $f = 1/\varepsilon ^3$ on $V$ and $f=1$ on $\Sigma \setminus W$.\\ 
Define the new metric $g = f^2 g_0$ on $\Sigma$ and the pulled-back metrics $\tilde g$ and $\tilde g_0$ on the universal covering $Y := \widetilde \Sigma$ 
of $\Sigma$; as $g\ge g_0$, one gets $\tilde g \ge \tilde g_0$ and thus $B_{(Y,\tilde g)}(y,R) \subset B_{(Y,\tilde g_0)}(y,R)$ for every $y \in Y$; this implies 
that $\Ent (Y , \tilde g) \le \Ent (Y , \tilde g_0) =1$.\\
Consider now the action (as fundamental group) of $\Gamma$ on $\big(Y\,,\, d_{\tilde g} \,,\, dv_{\tilde g} \big)$ (endowed with the Riemannian distance and
measure associated to $\tilde g$) and the quotient $(\overline Y , \bar d) := (\Gamma \backslash Y , \bar d ) = (\Sigma , d_g)$. Then all the results of the
present section \ref{transport} applies to the action of $\Gamma$ on $\big(Y\,,\, d_{\tilde g} \,,\, dv_{\tilde g} \big)$; in particular Theorems \ref{topotubes}
and \ref{tubelong} are valid here. For every
$i \in \{1 , 2 , \ldots , N\}$, every non homotopically trivial loop of $ \overline Y = \Sigma$ (with base-point $\bar y_i$) must go outside 
$B_{g_0}(\bar y_i,\varepsilon ^2)$ (because $B_{g_0}(\bar y_i,\varepsilon ^2)$ is diffeomorphic to $\mathbb B^n$), and thus its length (with respect to 
the new metric $g$) is at least $2/\e$. With the notations of Theorem \ref{topotubes}, let $\overline Y^k_{r}$ be any connected component of the $ r$-thin 
subset $\overline Y^{\diamond}_{r} = \overline Y_{r} =\{ \bar y \in \overline Y : \overline{\sys}_\Gamma (\bar y) < r  \}$, 
where the systole is measured with respect to the new metric $g$;
as $\overline{\sys}_\Gamma (\bar y_i) \ge 2/\e > r$, then $\overline Y^k_{r} = U \setminus \Bigl(\bigcup_{i=1}^N V_i\Bigr)$, each $V_i$ being a closed 
neighbourhood of $\bar y_i$ contained in $B_{g_0}( \bar y_i,2\varepsilon ^2)$. Hence, though all the hypotheses of Theorem \ref{topotubes} are verified, 
the fundamental group of $\overline{Y}_{r}^{k}$ is not virtually 
cyclic (it contains a free group with $N$ generators); however, as predicted by Theorem \ref{topotubes}, its image by the homomorphism 
$ i_*$ induced by the inclusion $ i : \overline{Y}_{r}^{k} \hookrightarrow \overline Y$ is a virtually cyclic subgroup of the fundamental group of $\overline{Y}$.

\item[(4)] {\em In Theorem \ref{topotubes} (ii), the conclusion \lq \lq the image by $ i_*$ of the fundamental group of $\overline{Y}_{r}^{k}$ is virtually
cyclic" cannot be replaced by the same conclusion with \lq \lq cyclic" instead of \lq \lq virtually cyclic" :}

\smallskip
In fact, let us consider the Riemannian surface $ (\Sigma , g_0)$ of the previous example and choose $ (\overline Y, g) := (\Sigma \times \R \text{P}^n ,\, g_0\times \e^2 \cdot can)$, the total space of its Riemannian universal covering is then 
$ (Y , \tilde g) = (\widetilde \Sigma \times \mathbb S^n, \,\tilde g_0 \times \e^2 \cdot can)$ 
and its fundamental group is $\Gamma = \Gamma_0 \times \Z/2\Z$, where $\Gamma_0$ is the fundamental group of $\Sigma$. With the same  parameters 
(and values of these parameters) as in the previous example, $(Y, \tilde g)$ is $(\ln 3 + \pi \e)$-hyperbolic, hence, as all the elements
$(g, s) \in  \Gamma_0 \times \Z/2\Z $ such that $g\ne e$ are torsion-free and verify $\ell(g,s) = \ell(g) \ge \e'_0$, it comes that 
$\Gamma \in \text{\rm Hyp}_{\rm thick} (\delta_0, \e'_0)$, where $\delta_0 = \ln 3 + \pi \e$ and $\e'_0 = 2 \e'_1$.
As $\e << r$, it follows from example (1) that the connected component $\overline Y^k_{r}$ of the $r$-thin subset which contains $c$ is in this case
diffeomorphic to $U\times \R \text{P}^n $, where $U$ is a tubular neighbourhood of $c$ with radius $L \ge \ln (r/\e)$ in the factor $\Sigma$ (see example (3)). 
Thus the image of the fundamental group of $\overline{Y}_{r}^{k}$ by the homomorphism $ i_*$ induced by the inclusion $ i : \overline{Y}_{r}^{k} \hookrightarrow \overline Y$ is a subgroup of the fundamental group of $\overline{Y}$ isomorphic to $\Z \times \Z/2\Z$; thus it is virtually 
cyclic but not cyclic.

\item[(5)] {\em Comparing to the reference-example (1), in the general case, each bounded connected component $\overline Y^k_{r}$ of the $ r$-thin subset 
$\overline Y^{\diamond}_{r}$ still contains some minimizing closed periodic geodesic $\bar c_k$ (of length $\e_k <r $), but it no longer contains a big tubular neighbourhood of this geodesic of radius $ C(n) \,\ln \left( \frac{r}{\e_k}\right)$:}

\smallskip
In fact, let us consider the hyperbolic surface $ (\Sigma , g_0)$ with small injectivity radius of the example (3), we proved in example (3) that 
the fundamental group $\Gamma$ of $\Sigma$ is a torsion-free element of $\text{\rm Hyp}_{\rm thick} (\delta_0, \e'_0)$, where $\delta_0 = \ln 3$ and $\e'_0 = 2 \e'_1$. Then the corresponding universal constant $r_0 = r_0 (H_0 , \e'_0)$ defined in \eqref{universalconstants3} is equal to 
$\dfrac{2 \e'_1 \ln 2}{13 \ln 3 + 8 \e'_1} < 1$. 

Fix any positive value $r \le r_0 $. Let $c$ be a small closed geodesic of $(\Sigma , g_0)$, of length $\e$ ($\e << r$), normally parametrized; the connected component of the $r$-thin subset which contains $c$ is then a tubular neighbourhood $U$ of $c$ with radius $L \ge \ln (r/\e)$ (see example (1)). 
Still denote by $\varrho (\cdot ) = d_{g_0}(\cdot ,c)$ the function distance 
to the geodesic $c$; parametrize $U$ by sending each point $y \in U$ on the couple $\big(t , c(s) \big) \in \, ] - L,L[ \times \text{Im} (c)$, where the 
image $\text{Im} (c)$ of $c$ is isometric to a circle of length $\e$, where $c(s)$ is the $g_0$-orthogonal projection of $y$ on the geodesic circle $c$ and 
where $t := \varrho (y)$. We then have $g_0= dt^2 +  (\cosh t)^2 ds^2$ at the point $\big(t , c(s) \big) $.
Define $ N := \left[\ln (r/\e)\right] -1$ and a $C^\infty$ function $u : [- L , L] \f [1,+\infty[$ which verifies
\begin{equation*}
\left\lbrace \begin{array}{l} u = 1 \text{ on } [- L , - L+ \e^2], \text{ on } [ L- \e^2 ,  L ] \text{ and on }\bigcup_{k= -N}^{k = N}\, [k-\varepsilon ^2,k+\varepsilon ^2],\\
 u = \frac{1}{\varepsilon} \text{ on } \bigcup_{k= -N}^{k = N-1} [k+2\varepsilon ^2 , k+1-2\varepsilon ^2]
\end{array}\right. 
\end{equation*}
and consider a new $C^\infty$ metric $g$ which coincides with $g_0$ on the outside of $U$ and is equal to \break
$dr^2 +  (u(t) \cosh t)^2 ds^2$ on $U$; as $u \ge 1$, the pulled-back metrics $\tilde g$ and $\tilde g_0$ on the universal covering $Y := \widetilde \Sigma$ 
of $\Sigma$ verify $\tilde g \ge \tilde g_0$ and thus $B_{(Y,\tilde g)}(y,R) \subset B_{(Y,\tilde g_0)}(y,R)$ for every $y \in Y$; this implies 
that $\Ent (Y , \tilde g) \le \Ent (Y , \tilde g_0) =1$.\\
Consider now the action (as fundamental group) of $\Gamma$ on $\big(Y\,,\, d_{\tilde g} \,,\, dv_{\tilde g} \big)$ (endowed with the Riemannian distance and
measure associated to $\tilde g$) and the quotient $(\overline Y , \bar d) := (\Gamma \backslash Y , \bar d ) = (\Sigma , d_g)$. Then all the results of the
present section \ref{transport} applies to the action of $\Gamma$ on $\big(Y\,,\, d_{\tilde g} \,,\, dv_{\tilde g} \big)$; in particular Proposition 
\ref{topotubes1} and Theorems \ref{topotubes} and \ref{tubelong} are valid here. 
For every $k \in \{-N , \ldots , 0, 1 , \ldots , N-1\}$ and every $y\in \varrho^{-1}(\{k+ \frac{1}{2}\})$, every loop with base-point $y$ either goes out of
$\varrho ^{-1}\big([k+2\varepsilon^2,k+1-2\varepsilon^2]\big)$ (and then its length with respect to the metric $g$ is at least 
$2(\frac{1 }{2} - 2\varepsilon ^2) = 1-4\varepsilon ^2$), or it remains inside $\varrho^{-1}([k+2\varepsilon ^2,k+1-2\varepsilon ^2])$, and then 
(if it is not homotopically trivial), its projection on $\text{Im} (c)$ is surjective, and this means that the loop writes (in coordinates) $\gamma : \tau \mapsto
\big( t(\tau), c(s(\tau)) \big)$ where $s(\tau)$ varies from $0$ to $\e$ (at least). As $g(\dot{\g}(\tau), \dot{\g}(\tau))^{1/2} \ge  u(t(\tau)) \cosh (t(\tau)) \, |s'(\tau)| \ge 
\frac{1}{\e} |s'(\tau)|$, we get that the $g$-length of $\g$ is at least $1$. One thus obtain that $\sys^\diamond_\Gamma (y) = \sys_\Gamma (y) \ge 
1-4\varepsilon ^2 > r $ at any point $y\in \varrho ^{-1}(\{k+ \frac{1}{2}\})$. On the other hand, for every $k \in \{-N , \ldots , 0, 1 , \ldots , N\}$ 
and  every $y \in \varrho^{-1}(\{k\})$, the non homotopically trivial loop $\g : \tau \mapsto (k, c(\tau))$ from $[0, \e]$ to $U$ (with base-point $y$) verifies
$g(\dot{\g}(\tau), \dot{\g}(\tau))^{1/2} = g_0(\dot{\g}(\tau), \dot{\g}(\tau))^{1/2} = \cosh k $, and its $g$-length is thus $\e \cosh k < r$. It follows that
$\sys^\diamond_\Gamma (y) < r $ for every $y \in \varrho^{-1}(\{k\})$ and thus that (with respect to the action of $\Gamma$ on 
$\big(Y\,,\, d_{\tilde g} \,,\, dv_{\tilde g} \big)$) $U$ contains 
 $2N + 1$ different connected components of $\overline Y_{r} $, denoted by 
$Y^{-N}_r, \ldots , Y^{0}_r , \ldots Y^{N}_r$, where each $Y^{k}_r$ is 
included in the portion $U_k = \varrho ^{-1} (]k-\frac{1}{2} , k+ \frac{1}{2}[)$ of the cylinder.
Applying Proposition \ref{topotubes1}, each of these components $Y^{k}_r$ contains at least one periodic geodesic. Moreover the length of each of these
components is smaller than $1$, thus (in contrast with Example (1)) it does not go to $+\infty$ when $\e \f  0$.

\item[(6)] {\em Comparing to the reference-example (1), in the general case, each connected component $\overline Y^k_{r}$ of the $ r$-thin subset 
$\overline Y^{\diamond}_{r}$ may contain several minimizing closed periodic geodesics:}

\smallskip
Let us revisit the example (5), redefining the function $u : [-L , L] \f [ 0 , 2 ]$ as the following one:
\begin{equation*}
\left\lbrace \begin{array}{l} u = 1 \text{ on } [- L , - L+ \e^2], \text{ on } [ L- \e^2 ,  L ] \text{ and on }\bigcup_{k= -N}^{k = N}\, [k-\varepsilon ^2,k+\varepsilon ^2],\\
 u = 2 \text{ on } \bigcup_{k= -N}^{k = N-1} [k+2\varepsilon ^2 , k+1-2\varepsilon ^2],
\end{array}\right. 
\end{equation*}
consider a new $C^\infty$ metric $g$ which coincides with $g_0$ on the outside of $U$ and is equal to \break
$dt^2 +  (u(t) \cosh t)^2 ds^2$ on $U$; as $u \ge 1$, the pulled-back metrics $\tilde g$ and $\tilde g_0$ on the universal covering $Y := \widetilde \Sigma$ 
of $\Sigma$ verify $\tilde g \ge \tilde g_0$ and thus $B_{(Y,\tilde g)}(y,R) \subset B_{(Y,\tilde g_0)}(y,R)$ for every $y \in Y$; this implies 
that $\Ent (Y , \tilde g) \le \Ent (Y , \tilde g_0) =1$.\\
Consider now the action (as fundamental group) of $\Gamma$ on $\big(Y\,,\, d_{\tilde g} \,,\, dv_{\tilde g} \big)$ and the quotient $(\overline Y , \bar d) := (\Gamma \backslash Y , \bar d ) = (\Sigma , d_g)$ (where $d_g$ is the Riemannian distance associated to the new metric $g$), then all the results of the
present section \ref{transport} applies to this action; in particular Proposition \ref{topotubes1} and Theorems \ref{topotubes} and \ref{tubelong} are valid 
here. As $\sys^\diamond_\Gamma (y) \le 2 \e \cosh \varrho (y) <r$ on the set $\{y \in U : \varrho (y) \le \ln \left( \frac{r}{2 \e}\right)\}$, the connected 
component of $\overline Y_r$ which contains the closed geodesic $c$ also contains the long tubular neighbourhood $\varrho ^{-1} ([- L' , L'])$ (where
$L' := \ln \left( \frac{r}{2 \e}\right)$) of this geodesic. Moreover, on each interval $]k- \frac{1}{2},k+\frac{1}{2}[$, the function $b(t) := u(t) \cosh(t)$ 
attains a local minimum at a point $t_k$ and, as every curve $ \g : s \mapsto \big( t(s) , c (s)\big)$ lying in a sufficiently small neighbourhood of 
the curve $ c_k : s \mapsto (t_k , c(s))$ verifies
$$ \sqrt{g \big( \dot{\g} (s) ,  \dot{\g} (s)\big)} = \sqrt{t'(s)^2 + b \big(t(s)\big)^2} \ge b \big(t(s)\big) \ge b(t_k) \,,$$
one has (by integration) ${\rm length} (\g) \ge {\rm length} (c_k)$ and, for every integer $k$, $c_k$ is a (locally) minimizing closed geodesic. Observe that the $c_k$'s are all homotopic to $c$ and this is consistent with Theorem \ref{tubelong} (vi) and (viii).

\end{itemize}

The example (3) shows that the connected components $\overline Y^{k}_{r}$  of the $ r$-thin subset $\overline Y^{\diamond}_{r}$ generally have a topology which is not that simple. Each of these components admits a \lq \lq soul" which is a closed geodesic but, when this geodesic is small, the corresponding 
component $\overline Y^{k}_{r}$ generally does not contain a big tubular neighbourhood of this closed geodesic (see example (5)). These differences with
the negatively curved Riemannian manifolds may be solved by considering tubes around small closed geodesics instead of $ r$-thin subsets: these tubes
have a great radius (see Theorem \ref{tubelong} (iii)) and a rather simple topology, because they are metric balls in a quotient of $Y$ (in most of the applications $Y$ is a simply connected space) by a virtually cyclic group (a cyclic group when $\Gamma$ is torsion-free, by Remark \ref{virtuelcyclique}), see Theorem \ref{tubelong} (ii).
Developing this viewpoint is the aim of the following subsection.

\subsubsection{The topology of Margulis' tubes is almost trivial}\label{tubessimples}

\small
Let $\Gamma$ be an element of one of the sets of groups  described in the beginning of Section \ref{classesalg} and $ (Y , d, \mu) $ any measured length space 
whose entropy is bounded from above, then the systole of any proper co-compact action (by isometries without fixed point preserving the measure) 
of $\Gamma$ on $ (Y , d, \mu) $  is generally not bounded from below when the diameter of the quotient-space $ \overline{Y}:= \Gamma \backslash Y$ is 
not bounded (compare with Theorem \ref{minorsystglobale}).\\
What can we say when the diameter of $\Gamma \backslash Y$ is not bounded ? We shall prove that, in this case, around each point where the systole is 
small, there exists a ball with big radius which is similar to a \lq \lq tube", i. e. which is isometric to a ball of big radius in the quotient of $Y$ by a virtually 
cyclic subgroup of $\Gamma$.

\normalsize

\smallskip
Let us first notice that $(Y,d)$, being a length space (whose distance is never infinite), is path-connected, for (by definition)
$$\forall y , y' \in Y \ \ \ \ d(y,y') := \inf_{c \text{ path from } y \text{ to } y'} \text{length of } c < + \infty \ .$$
For any $ y \in Y$ and any $ R > \frac{1}{2}\, \sys_\Gamma (y)$, we define $ \widehat U_{R} :=  
\bigcup_{\g \in \widehat{\Gamma}_{2 R} (y) } B_Y ( \gamma \,y , R)$ and, in order to apply the results of subsection \ref{topogeneral} in 
this case, we shall first prove that
\begin{equation}\label{connectedcomp}
\widehat U_{R} \text{ is the connected component of }\pi^{-1} \left( B_{\overline Y} (\bar y, R) \right) \text{ which contains } y \ ,
\end{equation}
where $ \bar y := \pi (y)$).
Indeed, denote by $  U_{R} (y)$ the connected component of $\pi ^{-1} \big( B_{\overline{Y}} (\bar y, R) \big)$ which contains $y$
and by $G_y$ the set $\{\g \in \Gamma : \g ( U_{R }(y)) = U_{R}(y)\}$. In the length space $(Y,d)$, every 
open ball is path-connected\footnote{In fact, in a length space, for every open ball $B_Y(x, r)$ 
and every point $z$ of this ball, there exists a continuous path $c$ between $x$ and $z$ such that every point $z'$ of this path verifies $ d(x, z') + d(z', z) 
< d(x,z) + \eta$ (for any $\eta$ satisfying $ 0 < \eta < r-  d(x,z) $), which implies that $ d(x, z') < r$ and that $c$ is contained in $B_Y(x, r)$.}; from this
and from the fact that $ \g \big( B_Y (   y ,  R) \big) =  B_Y (  \g y ,  R)$, as every connected subset $A$ of $ \pi^{-1} \left( B_{\overline Y} 
(\bar y, R ) \right)$  such that $ A \cap U_{R} (y) \ne \emptyset $ is included in $U_{R} (y)$, we 
deduce (by a closed chain of implications) that, for every $\g \in \Gamma$,
$$\ \g y \in   U_{R } (y) \iff  \g^{\pm 1} \big(   U_{R} (y) \big)\,\cap \, U_{R} (y)\ne \emptyset  
\iff \g^{\pm 1} \big(   U_{R} (y) \big) \subset U_{R} (y) \iff \g ( U_{R }(y)) = U_{R} (y) $$
\begin{equation}\label{equivalences}
\iff \g \in G_y  \iff B_Y (  \g y ,  R) \cap U_{R} (y) \ne \emptyset \iff B_Y (  \g y ,  R) \subset  
  U_{R} (y)\ ,
\end{equation}
where, in the first implication, we have used the fact that $y \in U_R(y)$ by definition.\\
As, by definition of the distance $\bar d$ on $ \Gamma \backslash Y$ (see Lemma \ref{autofidele} (i)), $\pi ^{-1} \big(B_{\overline{Y}} (\bar y, R) \big) 
= \bigcup_{\g \in \Gamma } B_Y ( \gamma \,y , R)$, equivalences \eqref{equivalences}  imply that
\begin{equation}\label{unionboules}
 U_{R} (y) =  \bigcup_{\g \in \Gamma} \big(B_Y ( \gamma \,y , R)\cap 
  U_{R} (y) \big) =  \bigcup_{\g \in G_y}  B_Y ( \gamma \,y , R) \ .
\end{equation}

Let us prove that $ \widehat{\Gamma}_{2 R} (y) \subset G_y$: in fact, for every $ \sigma  \in \widehat{\Sigma}_{2 R} (y)$, one has 
$d(y, \sigma \, y) < 2 R$; on the other hand, as $(Y,d)$ is a length space, there exists a \lq \lq quasi middle point" between points $y$ and $\sigma  y$,
i. e. a point $u$ which satisfies 
$$ d(y , u) = d( \sigma  y , u ) < \frac{1}{2} \big( d(y, \sigma  y) + \eta \big) < R\  \text{ for every } \eta \ \text{ such that } 
\ 0 < \eta < 2\, R -  d(y, \sigma  y)\ ;$$
it follows that $ B_Y (  y ,  R) \cap B_Y (  \sigma  y ,  R) \ne \emptyset $, thus (as $B_Y (  y ,  R) \subset U_{R} (y) $ by 
\eqref{equivalences}) that $ B_Y (  \sigma  y ,  R) \cap   U_{R} (y) \ne \emptyset $, and the equivalences \eqref{equivalences} imply
that $ \sigma  \in G_y$. This proves that $ \widehat{\Sigma}_{2 R} (y) \subset G_y$, hence that $ \widehat{\Gamma}_{2 R} (y) \subset G_y$.\\
We are now going to prove that $ G_y \subset \widehat{\Gamma}_{2 R} (y) $: in fact, as the balls of $(Y,d)$ are path-connected, 
$ U_{R } (y)$ is both connected and locally path-connected, thus it is path-connected; as a consequence, for every $ \g \in G_y$, there
exists a continuous path between $y$ and $\g y$ which is contained in $ U_{R} (y)$, and we have $ d( \Gamma\cdot y , u) = 
\bar d (\bar y , \pi (u)) < R$ for every point $u$ of this path (because $\pi (u) \in B_{\overline Y} (\bar y, R ) $).
By the compactness of this path and by continuity of the (strictly positive) function $u \mapsto R - \bar d \big(\bar y, \pi (u) \big)$, there 
exists some $\eta > 0$ such that $ d( \Gamma\cdot y , u) = \bar d \big(\bar y, \pi (u) \big) \le R - \eta$ in every point $u$ of this path.\\
Let us now consider a subdivision $ y=y_0 , \, y_1 ,\, \ldots ,\, y_{n-1} , \, y_n = \g y$ of this path such that $ d(y_{i-1}, y_i) < \eta$ for every  
$i \in \{1, \ldots ,n\}$; as $y_i$ is on the path for every $i \in \{0, \ldots ,n\}$, there exists $ \g_i \in \Gamma$ such that $ d(y_i , \g_i  y) \le R - \eta$ (we make the canonical choices $\g_0 = e$ and $\g_n = \g$). 
For every $i \in \{1, \ldots ,n\}$, we define $\sigma_i = \g_{i-1}^{-1} \g_i$ and obtain
$$ d( y , \sigma_i y) = d( \g_{i-1} y , \g_{i} y ) \le d( y_{i-1} , \g_{i-1} y  ) +  d(y_{i-1}, y_i) +d(y_i , \g_i  y) < 2\, R - \eta \ ;$$
this proves that, for every $i \in \{1, \ldots ,n\}$, $\sigma_i \in \widehat{\Sigma}_{ 2  R } ( y )$; as $\g = \sigma_1 \cdot \sigma_2 \cdot \ldots 
\cdot \sigma_n$, it implies that $\g \in \widehat{\Gamma}_{2 R} (y) $. We conclude that $ G_y \subset \widehat{\Gamma}_{2 R} (y) $ and thus that
$G_y = \widehat{\Gamma}_{2 R} (y) $.
Using this last equality, the writing \eqref{unionboules} of $U_{R} (y)$ and the definition of $\widehat U_{R}$, we get the equality 
$\widehat U_{R} = U_{R} (y)$, which ends the proof of the property announced in \eqref{connectedcomp} and thus proves that
\begin{equation}\label{groupinterpret}
\widehat{\Gamma}_{2 R} (y) = G_y = \left\{\g \in \Gamma : \g ( U_{R }(y)) := U_{R}(y) \right\} = \left\{\g \in \Gamma : 
\g ( \widehat U_{R }) = \widehat U_{R} \right \}
\end{equation}
As a consequence, we can apply the results of subsection \ref{topogeneral}, where we replace $\overline V $ by $ B_{\overline Y} (\bar y, R ) $,
which is open and connected because $(Y,d)$ is a length space; we also replace
$V $ by $ \pi^{-1} \left( B_{\overline Y} (\bar y, R ) \right) = \bigcup_{\g \in \Gamma} B_Y (\g y , R )  $, the connected
component $V^i$ of $V$ by the connected component $\widehat U_{R} = U_{R} (y)$ of $\pi^{-1} \left( B_{\overline Y} 
(\bar y, R ) \right)$, and the subgroup $ \widehat\Gamma^i_V $ of those $\g \in \Gamma$ such that $\g (V^i) = V^i $ by 
$\widehat{\Gamma}_{2 R} (y) $ which (by \eqref{groupinterpret}) is the subgroup of those $\g \in \Gamma$ such that $\g ( \widehat U_{R }) 
= \widehat U_{R}$.\\
Mimicking the arguments and results of subsection \ref{topogeneral}, this allows to define the quotient spaces 
$\widehat Y =   \widehat\Gamma_{2R} (y) \backslash Y $ and 
$\widehat\Gamma_{2R} (y) \backslash \widehat U_{R} $ and the quotient mapping $\widehat \pi : Y \f  \widehat\Gamma_{2R} (y) 
\backslash Y $; the map $\pi$ then goes down to the quotient and provides the mapping $ \bar \pi \,: \, \widehat\Gamma_{2R} (y) \backslash Y   
\f  \Gamma \backslash Y$, which satisfies $ \pi = \bar\pi  \circ  \widehat \pi $. As the inclusion mapping $ j : \widehat U_{R} \hookrightarrow Y $ trivially 
commutes with the two actions of $ \widehat\Gamma_{2R} (y) $ on $ \widehat U_{R} $ and on $Y$, it gives (going down to the quotients) 
the canonical inclusion mapping $ j' : \widehat\Gamma_{2R} (y)\backslash \widehat U_{R}   
\hookrightarrow \widehat\Gamma_{2R} (y) \backslash Y  $; we then define $ p$ (resp. $\widehat p$) as the restriction of $\pi$ (resp. of 
$\widehat \pi $) to $ \widehat U_{R} $ at the origin and to $B_{\overline Y} (\bar y, R )$ (resp. to $\widehat\Gamma_{2R} (y) 
\backslash \widehat U_{R} $) at the aim, equivalently $ \bar j \circ p = \pi \circ j $ (resp. $ j' \circ \widehat p = \widehat \pi \circ j$) where $\bar j$ is
the canonical inclusion map $B_{\overline Y} (\bar y, R ) \hookrightarrow  \bar Y$.
Consider now the map $\bar\pi \circ j' : \widehat\Gamma_{2R} (y)\backslash \widehat U_{R}    \f \Gamma \backslash Y = \overline Y$, 
as its image is included in the image of $\pi \circ j$, thus in $B_{\overline Y} (\bar y, R ) $ (because $\pi \circ j =  \bar \pi  \circ  \widehat \pi \circ j 
= \bar \pi  \circ j' \circ \widehat p$ and $\widehat p $ is surjective), it gives (by restriction at the aim) a map 
$\bar p :  \widehat\Gamma_{2R} (y)\backslash \widehat U_{R}  \rightarrow B_{\overline Y} (\bar y, R ) $ 
such that $\bar\pi  \circ j'  =\bar j \circ \bar p $. Moreover, one gets $ \bar p  \circ \widehat p  = p$ because
$$\bar j \circ (\bar p  \circ \widehat p) = \bar\pi \circ ( j' \circ \widehat p) =  \bar\pi \circ (\widehat \pi \circ j)
= \pi \circ j = \bar j \circ p \ .$$
All these results are summarized in the following diagram:

 \large
\begin{equation}\label{diagramme3} 
\xymatrix{
   \Gamma \curvearrowright     \hspace{-20mm}        & \hspace{2mm} Y \hspace{2mm}  \ar[r]^{\hspace{-10mm} \widehat \pi}   \ar@/^2pc/[rr]^{\pi}                               & \hspace{4mm}   \widehat Y =   \widehat\Gamma_{2R} (y) \backslash Y  \hspace{1mm}  \ar[r]^{\bar \pi} 
                                 &  \hspace{5mm} \overline{Y}  = \Gamma \backslash Y  \hspace{2mm}  \\
\widehat\Gamma_{2R} (y)  \curvearrowright  \hspace{-10mm}    &  \hspace{2mm} \widehat U_{R}  \hspace{1mm}  \ar@{^{(}->}[u]  \ar[r]^{\hspace{-15mm}  \widehat p}     \ar@/_2pc/[rr]^{p}        & \hspace{2mm} B_{\widehat Y} (\widehat y, R)   =   \widehat\Gamma_{2R} (y)\backslash \widehat U_{R}    \hspace{1mm}  \ar[r]^{\hspace{5mm} \bar p}  \ar@{^{(}->}[u]    
                                 &    \hspace{4mm} B_{\overline{Y}} (\bar y, R)      \hspace{5mm} \ar@{^{(}->}[u]  
 }
\end{equation}

\normalsize

In this diagram, the equality $B_{\widehat Y} (\widehat y, R)   =   \widehat\Gamma_{2R} 
(y)\backslash \widehat U_{R}  $ (where $\widehat y := \widehat p (y)$) is validated by the definitions of $\widehat U_{R} $ and of the 
quotient distance $\widehat d$ on $\widehat Y := \widehat\Gamma_{2R} (y) \backslash Y $, which imply that
$\widehat U_{R} = \{z \, : \,d \left(\widehat\Gamma_{2R} (y) \cdot y, z \right) < R\} = \{z \, : \,\widehat d 
\left(\widehat \pi (y) , \widehat \pi (z) \right) < R\} $, and thus that $ \widehat \pi \left(\widehat U_{R}\right) = 
\{ \widehat \pi (z) \, :\, \widehat d \left(\widehat y, \widehat \pi (z) \right) < R\} = B_{\widehat Y} (\widehat y, R) $.\\

\begin{lemma}\label{topotubes2} \emph{(Structure Lemma)}\label{lemmestructure}
For every proper action (by isometries without fixed point) of any group $\Gamma$ on any length space $ (Y , d) $, for any $y \in Y$, let us consider the 
open set $\widehat U_{R} :=  \bigcup_{\g \in \widehat{\Gamma}_{2 R}(y)} B_Y ( \gamma \,y , R)$, with the above notations, we have
\begin{itemize}
\item[(i)] $\widehat U_{R} $ is the connected component of $\pi^{-1} \left( B_{\overline Y} (\bar y, R ) \right)$ which contains $y$ and \break
$\widehat{\Gamma}_{2 R} (y)  = \left\{\g \in \Gamma : \g ( \widehat U_{R }) = \widehat U_{R} \right \}$;
\item[(ii)] the maps $\pi$, $ \widehat \pi$ are locally isometric coverings, and the maps $p$ (resp. $\widehat{p}$) are locally isometric 
coverings from $\widehat U_{R}$ onto $B_{\overline{Y}} (\bar y, R)$ (resp. onto $ B_{\widehat Y} (\widehat y, R)$);
\item[(iii)] $\bar p$ is a locally isometric homeomorphism from $  B_{\widehat Y} (\widehat y, R)$ onto $ B_{\overline{Y}} (\bar y, R)$.
\end{itemize}
\end{lemma}

Notice that, if $ R \le \frac{1}{2}\, \sys_\Gamma (y)$, then $\widehat{\Gamma}_{2 R} (y)$ is trivial, $\widehat U_{R} = B_Y (y , R)$, $\widehat Y= Y$,
$\widehat \pi = \id_Y$, $\widehat p = \id_{B_Y (y , R)}$ and $ B_{\widehat Y} (\widehat y, R) = B_Y (y , R)$, thus this Theorem reduces to classical 
properties of $\pi$.

\begin{proof} The point (i) has been proved in the beginning of the present subsection \ref{tubessimples}. Property (i), the connectedness of the
balls in the length space $(\overline Y , \bar d)$, and the fact that the action of $\Gamma$ is proper and by isometries without fixed point imply that 
the hypotheses of Lemma \ref{lemmetopologique} are verified if we replace the open set $\overline V$ of Lemma \ref{lemmetopologique} by
$ B_{\overline{Y}} (\bar y, R)$ and the connected component $V^i$ of $\pi ^{-1}(\overline V)$ by the connected component $\widehat U_{R} $ 
of $\pi^{-1} \left( B_{\overline Y} (\bar y, R ) \right)$ which contains $y$. The point (ii) is then a direct consequence of Lemma \ref{lemmetopologique} 
(i) and (ii). Lemma \ref{lemmetopologique} (iv) implies that the mapping $\bar p$ is a locally isometric homeomorphism between the two balls 
$ B_{\widehat Y} (\widehat y, R) = \widehat\Gamma_{2R} (y)\backslash \widehat U_{R} $ and $ B_{\overline Y} (\bar y, R ) $. This ends the proof of
(iii).
\end{proof}

For the sake of simplicity, in the quotient $\overline Y := \Gamma \backslash Y$, a loop $\bar c$ in $\overline Y$ will be said to be
(\emph{homotopically}) \emph{torsion-free} if its iterates $k\cdot \bar c$ (for all $k \in \Z^*$) are not homotopically trivial. A pair of loops $\bar c_1$ and $\bar c_2 $ will be said \emph{independent} if they are both homotopically torsion-free and if, for every $p, q \in \Z^*$, $p \cdot \bar c_1$ and 
$q \cdot \bar c_2 $ are not freely homotopic.

Given any parameters $\delta_0, H_0, D_0, \e'_0 > 0$, recall that $ s_0 := s_0 (\delta_0 , H_0, D_0) $, $\alpha_0 := \alpha_0 (\delta_0 , H_0, D_0) $, 
and $r_0 := r_0 (\delta_0, \e'_0)$ are the universal constants respectively defined at \eqref{defminorant}, at \eqref{universalconstants1} and at \eqref{universalconstants3}; define the function $ N' : \R_+^* \times \R_+^* \f  \N^*$ by 
$N' ( \delta_0, \e) := \Max \left( \left[\dfrac{13 \delta_0+  \e}{\e}\right], 2\right)$.
The following Theorem concerns the sets of groups $\text{\rm Hyp}_{\rm convex} 
(\delta_0, H_0, D_0)$ and $\text{\rm Hyp}_{\rm thick}  (\delta_0, \e'_0)$ introduced in Definitions \ref{Busemannaction} and \ref{systoleaction} 
respectively.

\begin{theorem}\label{tubelong}
Given any parameters $\delta_0,\, H_0, \, D_0 ,\, \e'_0,\,  H > 0$, for every element $\Gamma$ of \linebreak
$\text{\rm Hyp}_{\rm convex} (\delta_0, H_0, D_0)$ (resp. of $\text{\rm Hyp}_{\rm thick} (\delta_0, \e'_0)$) and for every positive 
$\e \le \frac{1}{3 n'_0 H}$, where $ n'_0 := N' ( \delta_0 , s_0)$ (resp. $ n'_0 := N' ( \delta_0 , \e'_0)$), if we introduce the universal constant 
$R_{\e} := \frac{1}{4 H}\ \ln  \left(\dfrac{1 }{n'_0 H \e}\right) - \frac{1}{4} n'_0  \e$, then any proper action of $\Gamma$ (by 
isometries without fixed point preserving the measure) on a connected  measured length space $ (Y , d, \mu) $ whose entropy is bounded from above by $H$
enjoys the following properties for any point $y \in Y$ such that $\sys_\Gamma^\diamond (y) \le \e$, for $\bar y := \pi (y) \in \overline Y := 
\Gamma \backslash Y$ and for $ \widehat y := \widehat \pi (y) \in \widehat Y := \widehat \Gamma_{2R_\varepsilon }(y) \backslash Y$ (here we use the 
same notations as above and refer to the diagram \eqref{diagramme3}, where we replace $R$ by $R_\varepsilon$):

\begin{itemize}
\item [(i)] the subgroups $\widehat{\Gamma}_{2R_\varepsilon }(y)$ and $\Gamma_{2R_\varepsilon }(y)$ are virtually cyclic and contain $\Gamma_{r'}(y')$ 
for any $r' <  2 R_{\e}$ and any $y' \in B_Y(y , R_{\e} - \frac{1}{2}\,r')$;

\item [(ii)] $ B_{\overline{Y}} (\bar y, \frac{1}{2}\,R_{\e})$ is isometric to a ball of radius $\frac{1}{2}\, R_{\e}$ in the quotient of $(Y,d)$ by a virtually 
cyclic subgroup of $\Gamma$, more 
precisely this isometry is the map $ \bar p$ from $ B_{\widehat Y} (\widehat y, \frac{1}{2}\,R_{\e})$ to $ B_{\overline{Y}} (\bar y, \frac{1}{2}\,R_{\e})$; 

\item [(iii)] if $\sys_\Gamma^\diamond (y) < \e$, there exists a homotopically torsion-free loop $\bar c$ (resp. $\widehat c$) with base-point $\bar y$ (resp. $\widehat y$) in 
$\overline Y$ (resp. in $\widehat Y$) such that
$\text{length of } \bar c = \text{length of } \widehat c < \e $ and such that the tubular neighbourhood of $\bar c $ in $\overline Y$ of radius 
$\frac{1}{2}\,(R_{\e} - \e)$ is isometric to the corresponding tubular neighbourhood of $\widehat c $ in the quotient $\widehat Y$ of $(Y,d)$ by a 
virtually cyclic subgroup of $\Gamma$.

\item [(iv)] Considering the $r$-thin subset $Y_r^\diamond$ corresponding to positive values $r \le \frac{\alpha_0}{H}$ 
(resp. $r \le \frac{r_0}{H}$), every
of its connected components $Y_r^i $ which intersects $B_Y(y , R_{\e} - \frac{1}{2}\,r)$ has the following property: every torsion-free 
element\footnote{Recall that $\widehat{\Gamma}^i_r $ is the subgroup of the elements $\g \in \Gamma$ which verifies $\g (Y_{r}^i) = Y_{r}^i$, see
subsection \ref{topominces}.} $\g \in \widehat\Gamma^i_r $ admits a power $ \g^k $ which is a non trivial element of 
$\widehat{\Gamma}_{2R_\varepsilon }(y)$.
Moreover, for every pair $Y_r^i $, $Y_r^j $ of connected components of $Y_r^\diamond$ which both intersect $B_Y(y , R_{\e} - \frac{1}{2}\,r)$,
for every torsion-free elements $\g \in \widehat\Gamma^i_r $ and $g \in \widehat\Gamma^j_r $, there exist $p, q \in \Z^*$ such that $ \g^p = g^q$.
\end{itemize}

If moreover $Y$ is simply connected then

\begin{itemize}
\item [(v)] if we denote by $ i_*\,: \pi_1 \left(B_{\overline{Y}} (\bar y, R_{\e}) ,\bar y\right) \to \pi_1 
\big(\overline Y , \bar y \big)$ the morphism induced by the inclusion $ i : B_{\overline{Y}} (\bar y, R_{\e}) \hookrightarrow \overline Y$, 
the subgroup $\, \widehat{\Gamma}_{2R_\varepsilon }(y)$ is identified with 
$i_* \pi_1 \left(B_{\overline{Y}} (\bar y, R_{\e}) ,\bar y \right) $ by the canonical isomorphism $\Gamma \f  \pi_1 \big(\overline{Y} , \bar y \big) $.\\
As a consequence $i_* \pi_1 \left(B_{\overline{Y}} (\bar y, R_{\e}),\bar y \right) $ is a virtually cyclic subgroup of $\pi_1 \big(\overline Y , \bar y \big)$;

\item [(vi)]  for every $r' <  2 R_{\e}$, for any homotopically torsion-free loops $\bar c_1$ and $\bar c_2 $ of lengths $\le r'$, whose images intersect $B_{\overline{Y}} (\bar y, R_{\e}- \frac{r'}{2}) $, 
there exist $p , q \in \Z^*$ and a (free) homotopy between the loops $p\cdot\bar c_1$ and $q \cdot\bar c_2$;

\item [(vii)] In $(\overline Y,\bar d)$, every pair of independent loops $\bar c_1$ and $\bar c_2 $ of lengths $\le \e$ have disjoint tubular 
neighbourhoods of radius $ \left(\frac{1}{2} \,R_{\e} - \frac{1}{4} \, \e\right)$;

\item [(viii)] Considering the $r$-thin subset $\overline Y^\diamond_r $ of $\overline Y := \Gamma \backslash Y$ corresponding to positive values 
$r \le \frac{\alpha_0}{H}$ (resp. $r \le \frac{r_0}{H}$), in every pair $\overline Y_r^i,\ \overline Y_r^j$ of its connected components which 
intersect $B_{\overline{Y}} (\bar y, R_{\e} - \frac{1}{2} \,r )$, loops $\bar c_i$ and $\bar c_j$ of lengths $< r$ (with base-points in 
$\overline Y_r^i$ and $\overline Y_r^j$ respectively) cannot be independent.
\end{itemize}
\end{theorem}

\begin{remarks}
  \begin{itemize}
    \item About Properties (ii) and (v): When $Y$ is simply connected, as $\widehat Y$ is a quotient of a simply connected space by a virtually cyclic group, its 
fundamental group is virtually cyclic. However example (3) of Subsection \ref{examplescounter} proves that the topology of 
$ B_{\widehat Y} (\widehat y, \frac{1}{2}\,R_{\e})$ is not that simple because, in this example,  this ball is a cylinder with many holes.
    \item About Properties (iv) and (viii): example (5) of Subsection \ref{examplescounter} shows that the \lq \lq long tube" $B_{\overline{Y}} (\bar y, R_{\e})$
may intersect many of the connected components of the $r$-thin subset $\overline Y_r^\diamond$. One can verify on this example
that the closed geodesics which are the \lq \lq souls" of these components admit multiples which are homotopically equivalent.
  \end{itemize}
\end{remarks}

\begin{proof}[Proof of Theorem \ref{tubelong}]
Notice that, as $\dfrac{1}{n'_0 H \e} \ge 3$, then $\dfrac{1}{n'_0 H \e} \ln \left( \dfrac{1}{n'_0 H \e}\right) \ge  3 \ln 3 > 3 \ge \dfrac{n'_0 + 4}{n'_0}$, 
and thus $ R_{\e} > \e$; remark also
that every open ball of $(Y,d)$ is path-connected (this is a classical property of length spaces, see the footnote in the beginning of the proof of 
property \eqref{connectedcomp}), thus every ball in the quotient metric space $(\overline Y , \bar d)$ is also path-connected.
It follows from these remarks, from \eqref{connectedcomp}, and from the fact that the action of $\Gamma$ is proper and without 
fixed point, that the hypotheses of Lemmas \ref{lemmetopologique} and \ref{lemmestructure} are verified if we replace the radius $R$ of the balls
involved in Lemma \ref{lemmestructure} by $R_{\e}$, the open set $\overline V$ of Lemma \ref{lemmetopologique} by $ B_{\overline{Y}} 
(\bar y, R_{\e})$ and the connected component $V^i$ of $\pi ^{-1}(\overline V)$ (involved in Lemma \ref{lemmetopologique}) by the 
connected component $\widehat U_{R_{\e}} $ of $\pi^{-1} \left( B_{\overline Y} (\bar y, R_{\e} ) \right)$ which contains $y$, and the 
subgroup $\widehat\Gamma^i_V $ by $ \widehat{\Gamma}_{2R_\varepsilon }(y)$.
\begin{itemize}

\item The first part of (i) is an immediate consequence of Theorem \ref{minorsystglobale} (i), which proves that, under the hypotheses of 
Theorem \ref{tubelong}, $\Gamma_R (y)$ is virtually cyclic for every $R \le 2 R_{\e}$.\\
Now, for every $y' \in B_Y(y , R_{\e} - \frac{1}{2}\,r')$ and any $\sigma \in \Sigma_{r'} (y')$, the triangle inequality gives:
$d( y , \sigma \, y) \le 2\, d( y, y') + d(y' , \sigma \, y' )< 2 \, R_{\e}$ and thus $ \sigma \in \widehat{\Gamma}_{2R_\varepsilon }(y)$;
this proves that $\Sigma_{r'} (y') \subset \widehat{\Gamma}_{2R_\varepsilon }(y)$ and thus that $\Gamma_{r'}  (y') \subset 
\widehat{\Gamma}_{2R_\varepsilon }(y)$.

\item \emph{Proof of (ii)} : As $\pi \big( B_{Y} (y , \frac{1}{2}\,R_{\e}) \big) = B_{\overline{Y}} \big(\bar y ,  \frac{1}{2}\,R_{\e}\big) $ and 
$\widehat \pi \big( B_{Y} (y , \frac{1}{2}\,R_{\e})\big) = B_{\widehat Y} \big(\widehat y , \frac{1}{2}\,R_{\e}\big)$ (by definition of the 
quotient-distances $\bar d$ and $\widehat d$, see above), we have 
$$B_{\overline{Y}} \big(\bar y ,  R_{\e}/2\big) = \overline \pi \circ \widehat \pi \big( B_{Y} (y , R_{\e}/2)\big) = \overline \pi \left( B_{\widehat Y} 
\big(\widehat y , R_{\e}/2\big)\right) = \bar p \left( B_{\widehat Y} \big(\widehat y , R_{\e}/2 \big)\right) ;$$
Lemma \ref{lemmestructure} (iii) then implies that the map $ \bar p$ may be restricted as a locally isometric homeomorphism from 
$B_{\widehat Y} \big(\widehat y , \frac{1}{2}\,R_{\e}\big)$ onto $B_{\overline{Y}} \big(\bar y , 
\frac{1}{2}\,R_{\e}\big)$. In order to show that this restriction of $ \bar p$ is an isometry, it is sufficient to prove that, for every pair $ z , \, z' \in  
B_Y (y , \frac{1}{2}\,R_{\e}) $, one has $\widehat d \big(\widehat \pi (z) ,\widehat \pi (z') \big) = \bar d \big( \pi (z) , \pi (z') \big)$: 
indeed $\widehat{\Gamma}_{2 R_{\e}} (y) \subset \Gamma$ and every $ \g \in \Gamma \setminus \widehat{\Gamma}_{2 R_{\e}} (y) $ belongs to $\Gamma \setminus \widehat{\Sigma}_{2 R_{\e}} (y)$, 
and thus verifies $d(y , \g \, y) \ge 2 R_{\e} $; the triangle inequality then gives:
$$ d( \g\, z , z') \ge d( \g \, y, y) - d(y, z') - d( \g \,y, \g\, z) \ge  2 R_{\e} - d(y, z') - d( y, z)>  R_{\e} > d(z,z') \ ,$$
and consequently
$$\bar d \big( \pi (z) , \pi (z') \big) = \inf_{\g \in \Gamma} d(\g z, z') = \inf_{\g \in \widehat{\Gamma}_{2 R_{\e}} (y)} d(\g z, z') = 
\widehat d \big(\widehat \pi (z) ,\widehat \pi (z') \big) \ .$$
We conclude by the remark that $B_{\widehat Y} \big(\widehat y , \frac{1}{2}\,R_{\e}\big)$ is a ball in the quotient $(\widehat Y, \widehat d )$ of $(Y,d)$ 
by the virtually cyclic group $\widehat{\Gamma}_{2 R_{\e}} (y)$.

\item \emph{Proof of (iii)} : As $ \sys^\diamond_\Gamma (y)  < \e$, there exists a torsion-free element $\sigma \in \Gamma^*$ such that 
$d(y, \sigma y) < \e$, then $\sigma \in \widehat{\Gamma}_{\e}(y) \subset \widehat{\Gamma}_{2R_\varepsilon }(y)$ for $\e < 2 R_{\e}$.
As $(Y,d)$ is a length space, there exists a path $c : [0,1 ] \to Y$ such that $c(0) = y$ and $ c(1) =\sigma y$ satisfying $\text{length of } c < \e$ and, as $R_{\e} > \e $, this path is included in $\widehat{U}_{R_\varepsilon }$. We define $\bar c := \pi \circ c = p \circ c$ and $\widehat c  := \widehat \pi \circ c
= \widehat p \circ c$ and get $\bar c = \bar p \circ \widehat c$; $\bar c$ and $\widehat c$ are loops in $\overline Y$ and $\widehat Y$ respectively,
because $\pi (\sigma y) = \pi (y)$ and  $ \widehat \pi (\sigma y) = \widehat \pi  (y)$ for $\sigma \in \widehat{\Gamma}_{2R_\varepsilon }(y)$; for every $k \in \Z^*$, as $\sigma^k \ne e$, then $\sigma^k y \ne y$ and 
the loops $k\cdot \widehat c $ and $k\cdot \bar c $ are not homotopically trivial.
The tubular neighbourhood of $\bar c $ (resp. of $\widehat c $) in $\overline Y$ (resp. in $\widehat Y$) of radius $\frac{1}{2}\,(R_{\e} - \e)$ being included
in $ B_{\overline{Y}} (\bar y, \frac{1}{2}\,R_{\e})$ (resp. in $B_{\widehat Y} (\widehat y, \frac{1}{2}\,R_{\e})$) by the triangle inequality, (ii) proves that
$\bar p$ is an isometry from the tubular neighbourhood of $\widehat c $ onto the tubular neighbourhood of $\bar c $.

\item \emph{Proof of (iv)} : Let $y_i$ be a point of $Y^i_r \cap B_Y(y , R_{\e} - \frac{1}{2}\,r)$, as $\sys^\diamond_\Gamma (y_i) < r$, there exists some torsion-free
$\sigma \in \widehat{\Sigma}_r (y_i)$; as (by Theorem \ref{topotubes} (i)) $\sigma $ is a torsion-free element of the virtually cyclic group 
$\widehat{\Gamma}^i_r $ then, for every torsion-free $\g \in \widehat{\Gamma}^i_r$, there exist $p, q \in \Z^*$ such that $ \g^p = \sigma^q$.
On the other hand (i) implies that $\sigma \in \widehat{\Gamma}_{2R_\varepsilon }(y)$, thus that $ \g^p = \sigma^q \in 
\widehat{\Gamma}_{2R_\varepsilon }(y)$.\\
Now, for every pair $Y_r^i $, $Y_r^j $ of connected components of $Y_r^\diamond$ which both intersect $B_Y(y , R_{\e} - \frac{1}{2}\,r)$,
for every torsion-free elements $\g \in \widehat\Gamma^i_r $ and $g \in \widehat\Gamma^j_r $, there exist $k, s \in \Z^*$ such that $ \g^k $ and $g^s$
are torsion-free elements of $\widehat{\Gamma}_{2R_\varepsilon }(y)$ and, as $\widehat{\Gamma}_{2R_\varepsilon }(y)$ is virtually cyclic by (i), 
there exist $p, q \in \Z^*$ such that $ \g^p = g^q$.

\item (v) is then an immediate corollary of Lemma \ref{lemmetopologique} (v) and of point (i), which proves that $\widehat{\Gamma}_{2R_\varepsilon }(y)$ is 
virtually cyclic.

\item \emph{Proof of (vi)}:
Let us choose a point $\bar y_1$ (resp. $\bar y_2$) in the intersection of $B_{\overline{Y}} (\bar y, R_{\e}- \frac{r'}{2}) $ with the image of $\bar c_1$ (resp. 
of $\bar c_2 $). As $\pi$ is surjective from $B_{Y} (y, R_{\e}- \frac{r'}{2}) $ to $B_{\overline{Y}} (\bar y, R_{\e}- \frac{r'}{2}) $ (by definition of the 
quotient-distance $\bar d$), one can choose a point $ y_1 \in \pi^{-1} (\{\bar y_1\})$ (resp. $ y_2 \in \pi^{-1} (\{\bar y_2 \})$) such that
$y_1, y_2 \in B_{Y} (y, R_{\e}- \frac{r'}{2})$. Let us consider the homotopy class $\left[\bar c_1 \right] \in \pi_1 \big(\overline{Y} , \bar y_1 \big) $ 
(resp. $\left[\bar c_2 \right] \in \pi_1 \big(\overline{Y} , \bar y_2 \big) $) and its image $\sigma_1$ (resp. $\sigma_2$) by the canonical isomorphism
$ \psi_{y_1} : \left(\pi_1 \big(\overline{Y} , \bar y_1 \big) , + \right)\f  \left(\Gamma , \cdot \right)$ (resp. $\psi_{y_2} : 
\left(\pi_1 \big(\overline{Y} , \bar y_2 \big) ,+ \right) \f  \left(\Gamma , \cdot \right)$) associated to 
the choice of $ y_1 \in \pi^{-1} (\{\bar y_1\})$ (resp. of $ y_2 \in \pi^{-1} (\{\bar y_2 \})$), namely, for $i= 1, 2$, denoting by $c_i $ the continuous path 
of $Y$ satisfying $c_i (0) = y_i$ and $\pi \circ c_i = \bar c_i$, $\sigma_i :=  \psi_{y_i} \left( \left[\bar c_i \right]\right)$ is the only element of $\Gamma$
such that $\sigma_i (y_i)$ is the endpoint of the path $c_i$. By means of this isomorphism the hypothesis $ \forall k \in \Z^* \ \left[ k \cdot \bar c_i \right] \ne 0$ 
is traduced in $ \forall k \in \Z^* \ \sigma_i^k \ne e$, thus each $\sigma_i$ is torsion-free and verifies $d(y_i , \sigma_i\, y_i) \le \text{ length of }  c_i =
\text{ length of }  \bar c_i \le r'$; it follows that $\sigma_i \in \Sigma_{r'}(y_i) \subset \widehat{\Gamma}_{2R_\varepsilon }(y)$, where the inclusion is
a consequence of (i). As $\widehat{\Gamma}_{2R_\varepsilon }(y)$ is virtually cyclic (by (i)) and, as $\sigma_1$ and $\sigma_2$ are torsion-free elements
of $\widehat{\Gamma}_{2R_\varepsilon }(y)$, there exist $p, q \in \Z^*$ such that $\sigma_1^p = \sigma_2^q$ and this implies that the two loops
$p\cdot \bar c_1$ and $q\cdot \bar c_2$ are freely homotopic\footnote{This is a corollary of the following result: \emph{Given $g \in \Gamma$ and any path
$\g_1$ from $y_1$ to $g y_1$ (resp. $\g_2$ from $y_2$ to $g y_2$), the loops $\bar \g_1 := \pi \circ \g_1$ and $\bar \g_2 := \pi \circ \g_2$ are freely
homotopic}: in fact, if $\alpha$ is a path from $y_1$ to $y_2$, the continuous path $ c := \alpha \cdot \g_2 \cdot (g \circ \alpha ^{-1})$ obtained by 
concatenation of the paths $\alpha$, $\g_2 $ and $g \circ \alpha ^{-1}$ has origin at $y_1$ and endpoint at $g y_1$, thus the loops $\pi \circ c$ and 
$\bar \g_1 =\pi \circ \g_1$ are homotopic. On the other hand, as $\bar \alpha := \pi \circ \alpha$ is a path from $\bar y_1$ to $\bar y_2$,
the loops $\pi \circ c = \bar \alpha \cdot \bar \g_2 \cdot \bar \alpha ^{-1}$ and $\bar \g_2$ are 
freely homotopic, thus $\bar \g_1 $ and $\bar \g_2$ are freely homotopic.}.

\item \emph{Proof of (vii)}: arguing by contradiction, suppose that the intersection of the tubular neighbourhoods of radius 
$ \left(\frac{1}{2} \,R_{\e} - \frac{1}{4} \, \e\right)$ of $\bar c_1$ and $\bar c_2 $ is non empty, then there would exist points 
$\bar y_1$ and $\bar y_2$ on the images of $\bar c_1$ and $\bar c_2 $ (respectively) such that $\bar d (\bar y_1 , \bar y_2) < R_{\e}- \frac{\e}{2}$.
As the images of $\bar c_1$ and $\bar c_2 $ would then intersect $B_{\overline{Y}} (\bar y_1, R_{\e}- \frac{\e}{2}) $, and as 
$\overline \sys_\Gamma^\diamond (\bar y_1) \le \e$, using the point (vi) 
(and replacing $r'$ by $\e$ in it), this would imply that $\bar c_1$ and $\bar c_2 $ are not independent, in contradiction with the hypothesis. We conclude that
the two tubular neighbourhoods are disjoint.

\item \emph{Proof of (viii)}: suppose that $\bar c_i$ and $\bar c_j$ are both homotopically torsion-free and denote by $\bar y_i$ and $\bar y_j$ their 
respective base-points; as $\overline Y_r^i$ and $\overline Y_r^j$  both intersect $B_{\overline{Y}} (\bar y, R_{\e} - \frac{1}{2} \,r )$, 
then $\pi^{-1} (\overline Y_r^i)$ and $\pi^{-1} (\overline Y_r^j)$ both intersect $B_{Y} (y, R_{\e} - \frac{1}{2} \,r )$, and there exist two connected components $Y_r^i$ and $Y_r^j$ of  $\pi^{-1} (\overline Y_r^i)$ and $\pi^{-1} (\overline Y_r^j)$ (respectively) which intersect 
$B_{Y} (y, R_{\e} - \frac{1}{2} r )$. As $\pi (Y_r^i) = \overline Y_r^i $ and $\pi (Y_r^j) = \overline Y_r^j$ by Proposition \ref{topotubes1} (o) and (ii), 
there exist $y_i \in Y_r^i$ and $y_j \in Y_r^j$ such that $\pi(y_i) = \bar y_i$ and $\pi(y_j) = \bar y_j$.\\
Using the same notations and arguments as in the proof of (vi), we denote by $\sigma_i$ and $\sigma_j$ the respective images of the homotopy classes 
$\left[\bar c_i \right]$ and $\left[\bar c_j \right]$ by the isomorphisms $ \psi_{y_i} : \pi_1 \big(\overline{Y} , \bar y_i \big) \f  \Gamma $ and $ \psi_{y_j} : 
\pi_1 \big(\overline{Y} , \bar y_j \big) \f  \Gamma$ (these isomorphisms being associated to the aforementioned choices of $y_i \in \pi^{-1}(\{\bar y_i\}) $ 
and of $y_j \in \pi^{-1}(\{\bar y_j\}) $, then $\sigma_i$ and $\sigma_j$ are torsion-free, because these isomorphisms traduce the assumption \lq \lq $\bar c_i$ 
and $\bar c_j$ homotopically torsion-free" by \lq \lq $\sigma_i$ and $\sigma_j$ torsion-free". Furthermore $\sigma_i$ and $\sigma_j$ belong to 
$\widehat\Sigma_r (y_i)$ and $\widehat\Sigma_r (y_j)$ respectively, because (as in the proof of (vi)) $d(y_i , \sigma_i\, y_i) \le \text{ length of }  \bar c_i < r$ 
and $d(y_j , \sigma_j\, y_j) \le \text{ length of } \bar c_j < r$. A consequence is that $\sigma_i \in \widehat\Gamma^i_r $ and 
$\sigma_j \in \widehat\Gamma^j_r $ 
by Theorem \ref{topotubes} (i); applying point (iv), this proves that
there exist $p , q \in \Z^*$ such that $\sigma_i^p = \sigma_j^q$. Using the same argument as in the proof of (vi), this implies that $p\cdot \bar c_i$ and 
$q \cdot \bar c_j$ are freely homotopic.
\end{itemize}
\end{proof}


\section{Applications}

\emph{In this section we consider closed Riemannian manifolds $(M,g)$ or compact polyhedrons $(X,g)$ endowed with a piecewise smooth Riemannian metric. In both cases we consider the entropy of the universal cover $(\widetilde M, d_{\tilde g}, dv_{\tilde g})$ or $(\widetilde X, d_{\tilde g}, dv_{\tilde g})$. For the sake of simplicity we denote it by $\Ent (M, g)$ or $\Ent (X, g)$.}

\subsection{Application to Manifolds}\label{appli:manifolds}
Let $M$ be a closed $n$-dimensional manifold. Any two metrics $g_1$, $g_2$ on $M$ are Lipschitz equivalent so that 
 $\rm {Ent} (M, g_1) >0$  is equivalent to $ \rm {Ent} (M, g_2)>0$. Similarly, for any two finite generating sets $S_1$ and $S_2$
 of a group $\Gamma$, we have that ${\rm Ent} (\Gamma, S_1) >0$ is equivalent to ${\rm Ent} (\Gamma, S_2) >0$. Moreover, the so-called
${\rm \check{S}}$varc-Milnor's Lemma asserts that  
$\rm {Ent} (M, g) >0$ for one (and thus every) metric on $M$ is equivalent to ${\rm Ent} (\Gamma, S) >0$ for one (and thus every) generating set $S$ of $\Gamma$, \cite{Sva, Mil}.
In order to make the above quantities independent of the choices of the metric and of the generating set, we define the minimal entropy ${\rm Minent} (M)$ of a manifold $M$ 
and the algebraic entropy ${\rm Ent}\, (\Gamma )$ of a finitely generated group $\Gamma$
as follows.
 
\begin{equation}
{\rm Minent}(M) := \inf \{ {\rm Ent} ( M, g) {\rm Vol} (M, g)^{1\over n}\,| \, {g\hbox{ Riemannian metric on}} \, M \}\label{minent}
\end{equation}
\begin{equation}
{\rm Ent}(\Gamma ) := \inf \{ {\rm Ent} (\Gamma, S) \,| \, S\, {\rm \hbox{finite generating set of }}\Gamma \}\label{algent}
\end{equation}

In the sequel $M^n$ is a closed $n$-dimensional manifold with fundamental group $\Gamma$. We denote by $B\Gamma = K(\Gamma , 1)$ an Eilenberg McLane space of $\Gamma$.
The class of closed manifolds $M$ with positive minimal entropy contains those having 
a hyperbolic fundamental group thanks to the two following Theorems due to M. Gromov and I. Mineyev.

\begin{theorem}[\cite{Gr2}] Let $M$ be a closed manifold of dimension $n$.
Then, for every Riemannian metric $g$ on $M$,
$${\rm Ent}(M, g) ^n {\rm Vol}(M, g)\geq C(n) \|M\|,$$
where $C(n)$ depends on the dimension $n$.
\end{theorem}

\begin{defi}\label{defessen}
$M^n$ is said to be an essential manifold if the natural homomorphism $f: H_n (M, \mathbf A)  \rightarrow H_n (B\Gamma, \mathbf A)$ is non zero, where $\mathbf A =\mathbf Z$ when $M$ is orientable and $\mathbf A =\mathbf Z/2\mathbf Z$ when $M$ is non orientable. It is said to be essential on $\textbf R$ if
$f: H_n (M, \mathbf R)  \rightarrow H_n (B\Gamma, \mathbf R)$ is not trivial hence injective since $\dim  H_n (M, \mathbf R)=1$.
\end{defi}

\begin{theorem}[\cite{Min}] \label{theo:Min} Let $M$ be a n-dimensional closed manifold, $n\geq 2$. Assume that the fundamental group of $M$ is a  Gromov hyperbolic group and that $M$ is essential on $\textbf R$, then the simplicial volume of $M$ is positive, \sl{i.e.} $\|M\| >0$.
\end{theorem}

\begin{proof}
This is an immediate consequence of the results proved in \cite{Min} which we briefly describe. One has that the bounded cohomology of $M$, denoted by $H^n_b(M, \mathbf R )$, is isomorphic to $H^n_b(\Gamma, \mathbf R ):= H^n_b(B\Gamma, \mathbf R )$ (see \cite{Gr2}, section 3.1, corollary A, p.40), the bounded cohomology of its fundamental group. In \cite{Min} it is proved that under the hypothesis of $\Gamma$ being an hyperbolic group the natural map  $H^n_b(\Gamma, \mathbf R )\to  H^n(\Gamma, \mathbf R )$ is surjective, for $n\geq 2$. Now $H^n(\Gamma, \mathbf R )\simeq H^n(B\Gamma, \mathbf R )\simeq H_n(B\Gamma, \mathbf R )$, where we consider the simplicial cohomology of $B\Gamma$. Let us recall that, for a torsion-free hyperbolic group, $B\Gamma$ can be taken to be the Rips complex which is a contractible simplicial complex on which $\Gamma$ acts freely by isometries with compact quotient; this is a finite simplicial complex. If $\Gamma$ has torsion then its classifying space is of infinite dimension but has finitely many cells in any dimension (see \cite{BH}, p. 468--471 and \cite{GH}, chap. 4). Since $M$ is essential on $\textbf R$ the natural map $f: H_n (M, \mathbf R)  \rightarrow H_n (B\Gamma, \mathbf R)$ is injective and the fundamental class of $M$ is sent to a non zero element in $H_n(B\Gamma, \mathbf R )\simeq H^n(\Gamma, \mathbf R )$ which is represented by a bounded cycle; this implies that the simplicial volume is non trivial. 
\end{proof}

\begin{corollary}\label{varhyp}
Let $M$ be a closed manifold of dimension $n\geq 2$ which is essential on $\textbf R$. We assume that the fundamental group of $M$ is Gromov hyperbolic, then ${\rm Minent} (M)  >0.$
\end{corollary}

On the group side, there are two distinguished classes of finitely generated groups with positive algebraic entropy which are worth mentioning, the non virtually nilpotent solvable groups
and the hyperbolic groups.

\begin{theorem}[\cite{Osin}]\label{theo:Osin} Let $\Gamma$ be a finitely generated solvable group. If $\Gamma$ is not virtually nilpotent,
then ${\rm Ent}(\Gamma ) >0$.
\end{theorem}

\begin{theorem}[\cite{Kou}]\label{theo:Kou} Let $\Gamma$ be a finitely generated group. Assume $\Gamma$ is Gromov hyperbolic  and non virtually cyclic,
then ${\rm Ent}(\Gamma )>0$.
\end{theorem}

When $\Gamma$ is the fundamental group of a closed manifold $M$, it would be natural to have a relation between 
${\rm Minent}( M)$ and ${\rm Ent}(\Gamma )$. The following examples show that 
there exist manifolds $M$ with vanishing minimal entropy but with fundamental group $\Gamma$ such that ${\rm Ent}(\Gamma )>0$
so that a relation such as ${\rm Minent} (M)  \geq C(n){\rm Ent}(\Gamma )$ does not hold in full generality. 

In the example \ref{exaprod}, the fundamental group of $M$ is hyperbolic but $M$ is not essential (see the definition \ref{defessen}).

\begin{exa}\label{exaprod}
Let $N^n$ be any closed manifold with Gromov hyperbolic fundamental group $\Gamma$ and consider $M^{n+k} = N^n \times S^k$ the cartesian
product of $N^n$ with the $k$-sphere, for $k\geq 1$. Then, the fundamental group of $M$ coincides with $\Gamma$ and therefore ${\rm Ent}(\Gamma )>0$
by Theorem \ref{theo:Kou}. On the other hand, considering the product metrics $g_\epsilon := g \times \epsilon ^2 h$ on $M$, where
g is any metric on $N$ and $h$ any metric on $S^k$, we have 
$${\rm Ent}(M, g_\epsilon)= {\rm Ent}(N, g) {\rm\hbox{ and  }}  {\rm Vol}(M, g_\epsilon)= \epsilon ^k {\rm Vol}(N, g) {\rm Vol}(S^k, h)\,,$$
so that by taking $\epsilon$ arbitrarily small, we get ${\rm Minent} (M) =0$. 
\end{exa}

In the example \ref{exasol}, $M$ is essential but its fundamental group is not Gromov hyperbolic.

\begin{exa}\label{exasol}
Let $M_A$ be the mapping torus $M_A := \mathbb{T} ^2\times [0,1] / \sim$, where $(x,0) \sim (Ax, 1)$ in $\mathbb{T} ^2\times [0,1] $ with
$A = \begin{pmatrix}
2 & 1 \\
1 & 1
\end{pmatrix}$.
The fundamental group of $M_A$ is the solvable group $\Gamma _A := \mathbf Z ^2 \rtimes _A \mathbf Z$, where the action of $\mathbf Z$ on $\mathbf Z ^2$
is defined for $k\in \mathbf Z$ and $x\in \mathbf Z ^2$ by $k.x = A^k x$. We also recall that the action of $\Gamma _A$ on $\widetilde M_A = \mathbf R^2 \times \mathbf R$ is given by 
$$(p,k)\cdot (u,t) = (p + A^k u , k + t)\,,$$  
where $(p,k)\in \Gamma _A := \mathbf Z ^2 \rtimes _A \mathbf Z$ and $(u,t)\in \widetilde M_A = \mathbf R^2 \times \mathbf R$.

The group $\Gamma _A$ is non virtually nilpotent, therefore ${\rm Ent}\, \Gamma _A >0$ by Theorem \ref{theo:Osin}.
On the other hand, $M_A$ is a torus bundle over the circle so that there exists a family of collapsing Riemannian metrics $g_\epsilon$ on $M_A$ with bounded sectional
curvature and $ \lim _{\epsilon \to 0} {\rm Vol} (M_A, g_\epsilon) =0$ (This is part of Cheeger-Gromov Theory, see \cite{Che-Gro}). Since the sectional curvature of $g_\epsilon$ is bounded below, ${\rm Ent} (M_A, g_\epsilon) \leq C$ and therefore
${\rm Minent}(M_A )=0$.
\end{exa} 

Finally, in the next example \ref{exa:torsion}, the manifold $M$ is essential as a product of essential manifolds and its fundamental group is Gromov hyperbolic. However it has torsion. Notice that $M$ is not rationally essential, {\sl i.e.} not essential on $\textbf{Q}$.

\begin{exa}\label{exa:torsion} 
The same construction as the one done in Example \ref{exaprod} can be made with $S^k$ replaced by a real projective space, say, for example, $\textbf{R}P^3$ and the same conclusion holds.
\end{exa}

Our aim in what follows is to give an explicit lower bound of the minimal entropy of essential manifolds whose fundamental group is torsion-free and Gromov hyperbolic.

\begin{theorem}\label{theo} Let $M$ be a n-dimensional essential closed manifold, $n\geq 2$. Assume that the fundamental group $\Gamma$ of $M$ is  non cyclic, torsion-free and is a subgroup of finite index in some group $G$ belonging\footnote{We recall that 
$\text{\rm Hyp} (\delta , H)$ is the set of non virtually cyclic groups $G$ which admit a finite system of generators $S$ such that $G$  is $\delta$-hyperbolic 
(with respect to the associated algebraic distance $d_{S}$) and satisfies $\Ent (G, S) \le H $.} to $\text{\rm Hyp} (\delta , H)$ then, for every Riemannian metric $g$ on $M$,
$${\rm Ent}( M, g) ^n {\rm Vol}( M, g)\geq C(n,\delta, H) >0.$$
\end{theorem}

The case $n=1$ is not adressed since the only closed $1$-dimensional manifold is the circle and its fundamental group is elementary. Notice that the quantity $C(n,\delta, H)$ depends on the dimension $n$ and on the fundamental group only. In particular for two essential closed manifolds of the same dimension and having the same fundamental group the above lower bound of the minimal entropy is the same.

The proof of this result relies on two theorems that we state below. One is a theorem proved by S. Sabourau which  states a universal relation between the Margulis invariant
and the Filling radius of a Riemannian manifold $(M, g)$ such that $M$ is orientable and essential. We now give the necessary definitions.

Let $(M, g)$ be a closed Riemannian $n$-manifold.  Mimicking Definition \ref{Margulisconstant}, we define

\begin{defi}
The Margulis invariant of $(M, g)$, denoted by $\text{\rm Marg} (M, g) $ is defined by
$$ \text{\rm Marg} (M, g) := \sup \{r \, | \,  \Gamma _r (x) {\rm\hbox{  is cyclic for  every x}\in M} \}\, ,$$
where $\Gamma _r (x)$ is the subgroup of $\pi _1 (M^n)$ generated by all loops at $x\in M^n$ of length $ \le r$.
\end{defi}

As $\Gamma$ is supposed to be torsion-free, by Remark \ref{Margulisconstantbis}, $\text{\rm Marg} (M, g)$ coincides with the Margulis invariant of the action of $\Gamma$ 
on the Riemannian universal cover $(\widetilde M , \tilde g)$ (see Definition \ref{Margulisconstant}).

We now denote by $d_g$ the distance on $M$ induced by the Riemannian metric $g$. The map
\begin{equation}\label{kura}
\iota : (M, d_g) \rightarrow (L^\infty (M), \| . \|)
\end{equation}
defined by $\iota (x) := d_g (x, .)$ is an isometric embedding called the Kuratowski embedding. For $\epsilon >0$, we denote 
$\iota _\epsilon :M \rightarrow \mathcal U _\epsilon (M)$ the inclusion of $M$ into the $\epsilon$-neighbourhood of $\iota (M)$.

\begin{defi}
The Filling radius of $(M, g)$ is defined as
\begin{equation}\label{fillrad}
{\rm Fillrad} (M, g):= \inf \{\epsilon >0 \, | \,  (\iota _\epsilon)_* ([M])  = 0 \  \text{\rm  in } H_n( \mathcal U _\epsilon (M), \mathbb A ) \},
\end{equation}
where $\mathbb A = \mathbb Z$ or $\mathbb Z_2$ according to whether $M$ is orientable or not.
\end{defi}

The Filling radius has been introduced by M. Gromov in \cite{Gro2} where are related the volume and the systole of an essential Riemannian manifold. 
One step toward this is the following statement, which we will need later.
\begin{theorem}[\cite{Gro2}]\label{gromov}
Let $(M, g)$ be a closed Riemannian manifold of dimension $n$,
then
$${\rm Vol}(M, g) \geq C(n) {\rm Fillrad} (M, g) ^n\,,$$
where $C(n)$ depends on the dimension $n$ only.
\end{theorem}

The following is Theorem 4.5 in \cite{Sab}.
\begin{theorem}[\cite{Sab}]\label{Sab}
Let $(M, g)$ be a n-dimensional essential closed manifold, $n\geq 2$, whose fundamental group is torsion-free and Gromov hyperbolic. Then, $\text{\rm Marg} (M, g) \leq 8 \,{\rm Fillrad} (M, g)$. In particular, 
$${\rm Vol} (M, g) \geq C'(n) \text{\rm Marg} (M, g)^n \,,$$
where $C'(n)$ depends on the dimension $n$ only.
\end{theorem}

\begin{proof}[Proof of Theorem \ref{theo}]
 As $\Gamma$ is non cyclic and torsion-free, it is non virtually cyclic by Remark \ref{virtuelcyclique}, hence $\Gamma$ 
belongs to $\text{\rm Hyp}_{\rm sub} (\delta , H)$ (see Definition \ref{hyperbolicgroup}). We can thus apply Theorem \ref{transportnil} 
(ii) to the action of $\Gamma$ on the Riemannian universal cover $(\widetilde M , \tilde g)$ of $M$, endowed with the Riemannian distance and with the 
Riemannian measure, we get

\begin{equation}\label{theoMarg}
\Ent (M,g) \text{\rm Marg} (M, g) = \Ent (\widetilde M,d_{\tilde g}, dv_{\tilde g} ) \text{\rm Marg}_\Gamma (\widetilde M, d_{\tilde g}) \ge \alpha'_0 (\delta , H)
\end{equation}
(see \eqref{universalconstants2} for the definition of  $\alpha'_0 (.\,, .)$). Now, being a subgroup of finite index in a Gromov-hyperbolic group, 
$\Gamma$ is also Gromov-hyperbolic and, from Theorem \ref{Sab} and from the last inequality, we deduce
$$\Ent (M, g)^n \Vol (M, g) \ge C'(n) \text{\rm Marg} (M, g)^n \Ent(M,g)^n \ge C'(n)  \alpha'_0 (\delta , H)^n\,.$$
This finishes the proof of Theorem \ref{theo}. 
\end{proof}

\begin{exa}[Main example]
Let $N_1$ and $N_2$ be two closed manifolds of the same dimension $n$ and let us assume that each of them carries a Riemannian metric of negative curvature. We consider the connected sum $M=N_1\# N_2$. 

If  $n\geq 3$ the fundamental group $\Gamma$ of $M$ is the free product of the fundamental groups $\Gamma_1$ of $N_1$ and $\Gamma_2$ of $N_2$. The groups $\Gamma_1$ and $\Gamma_2$ are non elementary torsion-free hyperbolic groups hence $\Gamma = \Gamma_1 * \Gamma_2$ is also a 
non elementary torsion-free hyperbolic group (see \cite{GH}, exercise 34, p.19).

Now, the connected sum of two essential manifolds is essential (see \cite{Gro2}, p.3) hence Theorem \ref{theo} applies to $M$, for any 
finite system $S$ of generators of $\Gamma$, we get a positive lower bound of the minimal entropy of $N_1\# N_2$ in terms of the hyperbolicity, of an upper bound, of the Entropy of $(\Gamma , S)$
and of the dimension $n$.
\end{exa}

\subsection{Application to Polyhedrons}\label{subsect:poly}

In this section we extend Theorem \ref{theo} to the case of finite simplicial complexes endowed with a Riemannian (or Finsler) metric. We follow the definition of $\Delta$-complexes in \cite{Hat}, chapter 2. The spaces $X$ that we are considering can be built from collections of disjoint simplices by identifying various subsimplices, where the identifications are performed using linear homeomorphisms; in other words it is a PL-structure on $X$. We then endow each closed simplex with a smooth Riemannian (or Finsler) metric in such a way that it coincides on the intersections; it is sometimes called a piecewise Riemannian (or Finsler) metric and for the sake of simplicity we shall simply call it a Riemannian (or Finsler) metric. Now $X$, endowed with such a  metric, is called a Riemannian polyhedron or simply a  {\sl polyhedron} and the metric will be denoted by $g$ (notice that $C^\infty$ smoothness of the metric on the simplices may not be completely necessary in the results that follow, we could even consider some geodesic measured metric structures, see \cite{Gro2} p. 264, remark (b)). 
Furthermore we notice that, when $X$ is compact, it has a dimension $n$ which is the maximal dimension of a simplex in $X$. We can define its volume, diameter and  entropy since it has a universal cover.  More informations on simplicial complexes can be found in \cite{BH} for the topological aspects and in \cite{Bab} for the Riemannian aspects.
Finally, working with polyhedrons allows to define all the Riemannian invariants that we need, however it is worth noticing that a polyhedron is a special kind of CW-complex since the $k$-simplices are $k$-cells. This allows us to use the fact that continuous maps between cell-complexes are homotopic to a cellular map (see below) and to use the more flexible cellular homology as in Lemma \ref{lem:homology}.

Now, let $K=B\Gamma$ be the classifying space of the fundamental group $\Gamma$ of X. The space $K$ is an aspherical (possibly infinite) simplicial complex.

\begin{defi}
A polyhedron $X$ of dimension $n$ is said to be $n$-{\rm essential} if there exists a continuous map $F: X\to K$ which does not contract to the $(n-1)$-skeleton of K.
\end{defi}
This definition is different from the one used in \cite{Gro2} p.139 and \cite{Gr1} p.260. Indeed, in these references being essential means that there exists an aspherical simplicial complex with the above property, which might not be the classifying space of $\Gamma$. For our purpose we need properties of the free loop space of $K$ described in \cite{Sab}, Proposition 3.6 and Remark 3.8, which are satified by the classifying space $B\Gamma$. Any essential manifold $M^n$ is $n$-essential, however it is not clear to us whether the converse is true or not. 

In the sequel we use the following consequence of $n$-essentiality. For a CW-complexe $X$ we denote by $X^{(k)}$ the $k$-skeleton, that is the closure of the union of the $k$-cells. All homology groups below are with coefficients in $\textbf{k}=\textbf{Z}$ or $\textbf{k}=\textbf{Z}_2$.

\begin{lemma}\label{lem:homology}
If $X$ is a $n$-dimensional $n$-essential polyhedron, the induced map between the relative cellular homologies of $X$ and $K$ respectively, 
$$F_*:H_n(X^{(n)}, X^{(n-1)})\longrightarrow H_n(K^{(n)},K^{(n-1)})\,,$$
is non identically zero.
\end{lemma}
\begin{proof}[Proof of Lemma \ref{lem:homology}]
As mentioned above we can assume that $F$ is a cellular map between $X$ and $K$ and hence that each $p$-cell of $X^{(p)}$ is sent to $K^{(p)}$ for all $0\leq p\leq n$. We consider $X^{(n)}/X^{(n-1)}$, that is the space obtained by identifying  $X^{(n-1)}$ to a point. 
This space is a bouquet of (topological) $n$-spheres, each of them corresponding to a $n$-cell in $X^{(n)}$, similarly we define $K^{(n)}/K^{(n-1)}$. Now, let us take a topological $n$-sphere $S:=S^n$ in $K^{(n)}/K^{(n-1)}$, we can send $K^{(n)}/K^{(n-1)}$ onto $S$ by contracting the other spheres in $K^{(n)}/K^{(n-1)}$ to a point (the vertex of the bouquet). For the details the reader is referred to \cite{span}, Chapter 4, Section 4. Then we consider the composite map sending  $X^{(n)}/X^{(n-1)}$ to $S$ which we call $F_{S}$. If $F_*$ is zero then each of these maps has zero degree for each $S$ and hence the map $F_{S}$ can be deformed to a constant map. This shows that if $F_*$ is identically zero then $F$ can be 
continuously deformed to a map into $K^{(n-1)}$.
\end{proof}

We now state our main result which is a version of Theorem \ref{theo} for polyhedrons.

\begin{theorem}\label{theo:poly} Let $X$ be a n-dimensional $n$-essential compact polyhedron, $n>0$. Assume that the fundamental group $\Gamma$ of $X$ is non cyclic, torsion-free and is a subgroup of finite index in some group $G$ belonging\footnote{We recall that 
$\text{\rm Hyp} (\delta , H)$ is the set of non virtually cyclic groups $G$ which admit a finite system of generators $S$ such that $G$  is $\delta$-hyperbolic 
(with respect to the associated algebraic distance $d_{S}$) and satisfies $\Ent (G, S) \le H $.} to $\text{\rm Hyp} (\delta , H)$ then, for every Riemannian 
metric $g$ on $X$,
$${\rm Ent}( X, g) ^n {\rm Vol}( X, g)\geq C(n,\delta, H) >0.$$
\end{theorem}
Notice that contrary to Theorem \ref{theo}, the case $n=1$ is included in the above result, and examples are given below.

The proof is similar to the one of Theorem \ref{theo}: we however have to replace the filling radius by another invariant introduced by M.~Gromov in \cite{Gro2} p.138, namely the $(n-1)$ contractibility radius. More precisely, we shall again consider the Kuratowsky embedding

\begin{equation}\label{kura2}
\iota : (X, d) \rightarrow (L^\infty (X), \| . \|)
\end{equation}
defined by $\iota (x) := d (x, .)$ which is an isometric embedding, and, as before, for $\epsilon >0$, we denote by
$\iota _\epsilon :X \rightarrow \mathcal U _\epsilon (X)$ the inclusion of $X$ into the $\epsilon$-neighbourhood of $\iota (X)$.

\begin{defi}\label{definition:contrad}
Let $X$ be compact $n$-dimensional polyhedron.
  \begin{itemize}
    \item for $\epsilon >0$, the map $\iota_\epsilon$ is said to be $(n-1)$-contractible if it is homotopic to a map which factors through a $(n-1)$-dimensional polyhedron. That is, $\iota_\epsilon$ is homotopic to a continuous map $f=f''\circ f'$ with,
    $$X\stackrel{f'}\longrightarrow Y \stackrel{f''}\longrightarrow L^\infty (X)\,,$$
    where $Y$ is a $(n-1)$-dimensional polyhedron.
    \item The $(n-1)$ contractibility radius is defined by,
    $$\mathrm{ContRad}(X)=\inf \{ \epsilon >0;\, \iota_\epsilon\, is\,  (n-1){\textrm -contractible}\}\,.$$
  \end{itemize}
\end{defi}
In \cite{Gro2} p.138 a whole bunch of contractibility radii are defined indexed by a non-negative integer not greater than $n$. Here we only consider the possibility of contracting $\iota_\epsilon (X)$ to a $(n-1)$-dimensional polyhedron. 

The proof of Theorem \ref{theo:poly} is the same as the proof of Theorem \ref{theo}: the inequality \eqref{theoMarg} being valid for the action of 
$\Gamma$ on any metric measured space (see Theorem \ref{transportnil} (ii)), we get that ${\rm Ent}( X, g) \text{\rm Marg} (X, g) \ge \alpha'_0 (\delta , H)$. 
In order to prove the inequality ${\rm Vol} (X, g) \geq C'(n) \text{\rm Marg} (X, g)^n$ which concludes, we need the analogous of Theorems \ref{gromov} and \ref{Sab}, where we replace the Filling radius by the Contractibility radius. The analogous of Theorem \ref{gromov} is the following
\begin{theorem}[\cite{Gro2}, p.138]\label{theo:gromov2}
Let $(X,g)$ be an arbitrary $n$-dimensional compact polyhedron endowed with a Riemannian (or Finsler) metric, then
$${\rm Vol}(X, g) \geq C(n) {\rm ContRad} (X, g) ^n\,,$$
where $C(n)$ depends on the dimension $n$ only.
\end{theorem}

In order to obtain the analogous of Theorem \ref{Sab}, we had to adapt the proof of \cite{Sab} to this new context. While writing this text we got the article \cite{Sab-Bab} which gives a proof of the following theorem along the exact same lines and for the sake of completeness we shall only give the basic construction in the proof that follows.

\begin{theorem}[see \cite{Sab-Bab}]\label{theo:sab2}
Let $(X, g)$ be a n-dimensional compact $n$-essential polyhedron endowed with a Riemannian metric, $n\geq 1$, whose fundamental group is torsion-free and Gromov hyperbolic. Then, 

If $n\geq 2$, $\text{\rm Marg} (X, g) \leq 8 \,{\rm ContRad} (X, g)$. In particular, 
$${\rm Vol} (X, g) \geq C'(n) \text{\rm Marg} (X, g)^n \,,$$
where $C'(n)$ depends on the dimension $n$ only. 

Furthermore, if $n=1$ then we have,
$${\rm Vol} (X, g) \geq \text{\rm Marg} (X, g) \,.$$
\end{theorem}

The Margulis invariant being defined as before. 
\begin{proof}[Proof of Theorem \ref{theo:sab2}]
The proof follows exactly the scheme used by S.~Sabourau in the proof of Theorem 4.5 of \cite{Sab} and of Theorem 4.15 in the recently posted \cite{Sab-Bab}. We will not give all details and the reader is referred to these articles.

We proceed by contradiction and assume that ${\rm ContRad} (X, g) < {1\over 8}\mu (g)$. Let us choose $\rho$ such that  ${\rm ContRad} (X, g) < \rho < {1\over 8}\mu (g)$. Since $\rho > {\rm ContRad} (X, g)$ the map $\iota_\rho$ is homotopic to a continuous map $f: X \rightarrow\mathcal{U} _\rho (X)\subset L^\infty (X)$ which factors through a $(n-1)$-dimensional polyhedron $Y$, namely $f=f''\circ f'$, with $f': X\longrightarrow Y$ and $f'': Y \longrightarrow \mathcal{U}_\rho (X)\subset L^\infty (X)$. Since $X$ and $Y$ are simplicial complexes $f'$ can be approximated, and hence replaced, by a simplicial map still called $f'$. By replacing $Y$ by the sub-complex $f'(X)$ we also may assume that $f'$ is surjective.

Since X is $n$-essential there exists a continuous map $F : X\longrightarrow K$ which does not contract to the $(n-1)$-skeleton of $K$. We may identify $X$ with $\iota_\rho (X)=\iota (X)\subset L^\infty (X)$ since $\iota_\rho$ is an isometry hence an homeomorphism. 

The homotopy between $\iota_\rho$ and $f$ is a continuous map $h_0$,
$$h_0: X\times [0,1] \longrightarrow \mathcal{U}_\rho (X)\subset L^\infty (X)\,,$$  
with $h_{0\vert X\times\{ 0 \}}=\iota_\rho$ and $h_{0\vert X\times\{ 1 \}}=f$. Notice that $X\times [0,1]$ can be given the structure of a  $(n+1)$-dimensional polyhedron (after subdivisions). 

We now define the {\sl mapping cylinder} $P$ of the map $f'$ as the quotient of the disjoint union $(X\times [0,1])\sqcup Y$ obtained by identifying $(x,1)\in X\times [0,1]$ with $f'(x)\in Y$. The homotopy $h_0$ factors through this identification and yields a continuous map, which we denote by $h$,
$$h: P\longrightarrow \mathcal{U}_\rho (X)\subset L^\infty (X)\,,$$
which coincides with $\iota_\rho$ on $X\times \{ 0\}$ and with $f''$ on $(X\times \{ 1\})\stackrel{f'}{\sim} Y$.
The mapping cylinder of a simplicial map can be endowed with a simplical structure extending that of $(X\times \{0\})$ (see \cite{Sak}, Chapter 4, p.218).

The idea is now to try to extend $F:X\times \{ 0\}\simeq X\longrightarrow K$ to a continuous map $\bar F: P\longrightarrow K$. If so, $\bar F_{\vert (X\times \{ t\})}$ is a continuous deformation between  $F_{\vert (X\times \{ 0\})}$ and $\bar F_{\vert Y}$ and, furthermore, $P$ retracts onto its subcomplex $Y$ (see \cite{span}, Chapter 1, Section 5, p.~32). Then, since $Y$ and $K$ are CW-complexes, $\bar F_{\vert Y}$ is homotopic to a map sending $Y$ into the $(n-1)$-skeleton on $K$ and this would be in contradiction with the essentiality of $X$. Unfortunately, as in \cite{Sab} and \cite{Sab-Bab}, this scheme of proof is not possible in this context but instead we can extend $F$ to a map with value in another space $Z$. Some details are described below. 

Since $\iota_\rho$ is an isometric embedding the map $F$ can then be written as $F=G\circ\iota_\rho$ with  $G: \iota_\rho (X)\longrightarrow K$. Up to taking subdivisions in $P$, we may assume that the images by $h$ of its  have a diameter less than $\epsilon$. We choose this number so that  $0<\epsilon <{1\over 4}\text{\rm Marg} (X, g)-2\rho$.

We recall that we denote by $P^{(k)}$ the $k$-skeleton of $P$. We now define $\bar{F}$ on $P^{(0)}\cup (X\times \{ 0\})$ so that it coincides with $\iota_\rho$ on $X\times \{ 0\}$. Take $p\in P^{(0)}\setminus (X\times \{ 0\})$ and send it to any choice $\iota_\rho (x)$, $x\in X$, of a point in $\iota _\rho (X)$ nearest to $h(p)$; we then set $\bar F (p)= G\circ \iota_\rho (x)$. For adjacent vertices $p_i$ and $p_j$ of $P^{(0)}$, that are related by an edge, we consider the corresponding chosen points $x_i$ and $x_j$ in $X$ 
 and a shortest path between them (any choice). We now send the edge of $P^{(1)}$ between $p_i$ and $p_j$, to this chosen curve; again, composing with $G$ we get the extension $\bar F$ on $P^{(1)}$. We then have,
$$
\begin{array}{rcl}
d_X(x_i, x_j)=d_{L^\infty}(\iota_\rho (x_i), \iota_\rho (x_j)) & \leq & d_{L^\infty}(\iota_\rho (x_i), h(p_i))+d_{L^\infty}(h(p_i), h(p_j))+ d_{L^\infty}(h(p_j), \iota_\rho (x_j))\\
                     & \leq  & 2\rho + \epsilon := r< {1\over 4}\text{\rm Marg} (X, g)\,.
\end{array}
$$
For the details the reader is referred to \cite{Sab}, proof of Theorem 4.5 and \cite{Sab-Bab}, proof of Theorem 4.15.

This extension sends the boundary of every $2$-cell $\Delta $ of $P$ to the image by $\iota_\rho$ of a curve in $X$ whose length is at most $3r$. If this number were smaller than the systole of $X$, all these curves would be contractible and, as $K$ is aspherical, the obstruction theory would give the desired extension $\bar G$. Unfortunately, $3r$ may be too large and some of the boundary closed curves of some $2$-simplices might be longer than the systole. Let $\Delta$ be a $2$-cell of $P$, if $\bar F (\partial\Delta )$ is not contractible, the properties of the loop space of $K$ given in \cite{Sab}, Proposition 3.6 and Remark 3.8 produce a homotopy from $\bar F (\partial\Delta )$ to a closed curve depending only on the non trivial homotopy class of $\bar F (\partial\Delta )$. This homotopy class is included in a unique maximal infinite cyclic subgroup generated by a primitive loop $\gamma_\Delta$. Following \cite{Sab} we then define,
$$Z=K\cup \big(\bigcup_{\Delta\in\mathcal{C}}D_\Delta \big)\,$$
where we glue a $2$-cell $D_\Delta$ along $\gamma_\Delta$ for each $2$-cell $\Delta$ of $P$ such that $\bar F (\partial \Delta )$ is not contractible in $K$, the family of these $2$-cells being denoted by $\mathcal C$. The space $Z$ is not aspherical, hence it may be difficult to extend the above map to $3$-simplices. However S.~Sabourau overcame this difficulty in \cite{Sab}, by induction, and extend $F$ to a map from $P$ to $Z$, which we call $\bar F$ (see also \cite{Sab-Bab}, proof of Theorem 4.15). Notice that $\bar F$ sends $P^{(1)}\cup (X\times \{ 0\})$ into $K$.

To summarise, the map
$$X\times\{0\}\stackrel{\iota_\rho}\longrightarrow \iota_\rho (X)\stackrel{G}\longrightarrow K\hookrightarrow Z\,,$$
extends to
$$\bar F:P\longrightarrow Z\,,$$
such that $\bar F (P^{(1)}\cup X\times \{0\})= F(X)\subset K\subset Z$.

As mentioned before the maps $\bar F_{\vert (X\times \{ t\})}$ realize a homotopy between $F_{\vert X}$ and $\bar F_{\vert Y}$, where $Y$ is by construction a sub-complex of $P$. The map $\bar F_{\vert Y}$ is homotopic to a cellular map sending $Y$ into the $(n-1)$-skeleton of $Z$ (see \cite{span}, Chapter 7, Section 6, Corollary 18 p.~404). Hence, $\bar F (X\times \{0\})$ retracts in $Z$ to its $(n-1)$-skeleton. On the other hand $F(X)=\bar F (X\times \{0\})$ does not retract to the $(n-1)$-skeleton of $K$ by the definition of $n$-essentiality.

For $n\geq 3$ this is in contradiction with Lemma \ref{lem:homology}. Indeed, the injection $K\hookrightarrow Z$ induces an isomorphism between $H_n(K^{(n)},K^{(n-1)})$ and $H_n(Z^{(n)},Z^{(n-1)})$, if $n\geq 3$, since we just added a bunch of $2$-cells. Then $\bar F (X\times \{0\})$ cannot be identically zero in $H_n(K^{(n)},K^{(n-1)})\simeq H_n(Z^{(n)},Z^{(n-1)})$, however  it is homotopic to $\bar F (X\times \{1\})$ whose image is in the $(n-1)$-skeleton of $Z$ and hence is identically zero in $H_n(Z^{(n)},Z^{(n-1)})$. A contradiction.

The same proof works for $n=2$, that is, the map $F_*$ cannot be identically zero. In this case though the injection $K\hookrightarrow Z$ does not induce an isomorphism since we added $2$-cells in order to get $Z$. However, let us choose $B\in K\subset Z$ a $2$-cell such that the map $X^{(2)}\longrightarrow K^{(2)}\hookrightarrow Z^{(2)}\longrightarrow B/\partial B = S^2$ is not homotopic to a constant, we get a non trivial map $H_2(X,\textbf{k})\longrightarrow H_2(S^2, \textbf{k})$ yielding a contradiction as above.

Finally for $n=1$, $X$ is a graph with possibly edges which start and end at the same vertex. Its universal cover is a simply connected graph hence a tree (see \cite{span}, Chapter 3, Section 7, p.139). Its fundamental group is then a free group generated by primitive cycles and the condition that it is non elementary rules out graphs with only one primitive cycle. The systole of $X$ is the length of the shortest primitive cycle whereas the Margulis invariant is the length of the second shortest primitive loop. On the other hand the volume of $X$ is the sum of the length of all the edges, hence the inequality $\Vol (X)\geq \text{\rm Marg} (X, g)$ is obvious.
\end{proof}

\begin{remark}
It is worth noticing that Theorems \ref{theo:poly} and \ref{theo} could have a different form if we assumed that the fundamental group $\Gamma$ of $M$ or $X$ is a torsion-free and non cyclic subgroup of finite index in some group $G$ belonging\footnote{We recall that 
$\text{\rm Hyp}^{\star}_{\rm action} (\delta, H, D)$ is the set of groups $G$ which admit a proper action by isometries on some connected, non elementary, 
$\delta$-hyperbolic metric space whose entropy and co-diameter are bounded from above by $H$ and $D$ respectively.} to 
$\text{\rm Hyp}^{\star}_{\rm action} (\delta, H, D)$.
Indeed, it suffices, in the above proof, to replace  Theorem \ref{transportnil} (ii) by Theorem \ref{transportnil} (i). Theorems \ref{transportnil} (i) 
and \ref{theo:sab2} then prove that $\Ent (M,g)^n \Vol(M,g)$ (resp. $\Ent (X,g)^n \Vol(X,g)$) is bounded from below by $C'(n)\, \alpha_0 (\delta, H, D)^n$, 
where $\alpha_0 (\delta, H , D)$ is the universal constant 
defined in \eqref{universalconstants1}.
\end{remark}

\subsubsection{Some examples}

{\bf Graphs } 

Let $X$ be a finite connected polyhedron of dimension $1$, in other word, $X$ is a finite connected graph.
We endow $X$ with a polyhedron metric $g$, that is, every $1$-simplex is isometric to an interval of $\textbf R$ whose length is the only parameter. 
The corresponding volume is defined by $$\Vol (X,g) := \sum _{\sigma} l(\sigma),$$ where the $\sigma$'s are the $1$-simplices of $X$
and $l(\sigma)$ is the length of $\sigma$.
 We assume that
the fundamental group $\Gamma$ of $X$ is non elementary, thus $\Gamma$ is a non abelian free group on $r$ generators with $r\geq 2$.
Notice that if $S=(s_1, ...,s_r)$ is a canonical set of generators of $\Gamma$, then $(\Gamma,S)$ is $0$-hyperbolic and $\Ent (\Gamma, S) = \ln (2r-1)$.
Moreover, $X$ is $1$-essential and by Theorem \ref{theo:poly} above we have for every metric $g$ on $X$,
\begin{equation}\label{Ex1-1}
\Ent(X,g) \Vol (X,g) \geq C(r).
\end{equation}

It is known that, for a graph $X$, the rank $r$ of the fundamental group of $X$ is given by 
$$r=1+{1\over 2}\big ( \sum _{x\in V(X)} (k_x-1)\big )\,,$$
where $V(X)$ is the set of vertices of $X$ and $k_x +1$ denotes the valency of $X$ at $x$.
In the case when $k_x\geq 2$ for every $x\in X$, the exact value of the minimal entropy, as well as the minimizing metric, has been computed by S. Lim, \cite{Lim}, see also \cite{KN}, who showed that
\begin{equation}\label{Ex1-2}
\Ent(X,g) \Vol (X,g) \geq \frac{1}{2} \sum _{x\in V(X)} (k_x +1) \ln (k_x).
\end{equation}
Moreover,
the equality in (\ref{Ex1-2}) is achieved if and only if the metric $g$ is given, up to scaling, by 
\begin{equation}\label{Ex1-3}
l(e) := \frac{\ln \left( k_{i(e)}k_{t(e)}\right)}{\sum _{x\in V(X)} (k_x +1) \ln (k_x)},
\end{equation}
for every edge $e$ of $X$, where $i(e)$ and $t(e)$ are the initial and terminal point of $e$.

Notice that the exact value of the minimal entropy of $X$ depends on the simplicial structure of $X$ and not only
on its homotopy type. Indeed, for example, if $X_1$ is a bouquet of $2$ circles, its minimal entropy 
equals $2 \ln 3$, but if $X_2$ is a union of $2$ circles attached at the $2$ ends of a segment, then its minimal entropy
is $3 \ln 2 < 2 \ln 3$.

It raises the question of the dependance of the minimal entropy on the polyhedral decomposition. Notice that the minimal entropy decreases 
while taking subdivision, it might be thus possible that no best decomposition exists.

Note also that when $X$ is a finite graph, the minimal entropy of $X$ tends to infinity when the rank of the fundamental group of $X$ tends to infinity,
while our lower bound $C(r)$ of the minimal entropy given in (\ref{Ex1-1}) is tending to $0$.

\bigskip
{\bf Surfaces attached along a geodesic} 

Let $S_1$ and $S_2$ be two hyperbolic surfaces of genus two. Assume that each of them have a simple 
closed geodesic $\gamma _1\in S_1$ and $\gamma _2\in S_2$ with the same length. We denote by $X$ the polyhedron obtained from $S_1$ and $S_2$ identifying $\gamma _1$ and $\gamma _2$ by an isometry.
The universal cover of $X$ can be described inductively as follows. The universal cover $\widetilde S_i$ of $S_i$, $i=1,2$ are two copies of the Poincar\'e disk. In each $\tilde S _i$, 
denote $\tilde \gamma _i$ any lift of $\gamma_i$. Consider first $\widetilde S_1$ with the collection of all the lifts of $\gamma _1$ and attach to each of these a copy 
of $\widetilde S_2$ identifying $\tilde \gamma _1$ isometrically with some lift $\tilde \gamma _2$ of $\gamma _2$. Then, on every such attached copy of $\widetilde S_2$ consider 
each unidentified lift $\tilde \gamma _2$ and attach to it a copy of $\widetilde S_1$ identifying $\tilde \gamma _2$ with some $\tilde \gamma _1$. 
We keep doing these gluing and get the universal cover of $X$,
which looks like an hyperbolic building. In particular, it is $\CAT(-1)$. The fundamental group $\Gamma$ of $X$ is isomorphic to the amalgamated product $\Gamma _1 \underset{\mathbb Z}{\ast} \Gamma _2$
where $\Gamma _i$ are the fundamental groups of $S_i$, $i = 1,2$, hence $\Gamma$ is Gromov-hyperbolic. Note also that $X$ is $2$-essential.
By Theorem \ref{theo:poly}, we then have, for every metric on $X$, 
$$
\Ent (X,g) ^2 \Vol (X,g) \geq C(\delta, H),
$$
where $\delta\geq 0$ is the hyperbolicity constant of $\Gamma$ and $H$ satisfies $\Ent (\Gamma , \Sigma ) \leq H$ for some generating set $\Sigma$.

\begin{remark}
When the metric on $X$ comes from two hyperbolic metrics on $S_i$ such that the common length of the geodesics $\gamma _i$ is a number $\epsilon$ going to $0$, then, denoting
by $g_\epsilon$ the resulting metric on $X$, we see that $\Ent (X,g_\epsilon )$ tends to $1$ (cf. the remark on p. 4 in  \cite{Bou}) and 
$$\lim _{\epsilon \rightarrow 0} \Ent (X, g_\epsilon ) ^2 \Vol (X, g_\epsilon ) = \Vol (S_1) + \Vol (S_2)\,,$$ where $\Vol (S_i)$ denotes 
the hyperbolic volume of $S_i$.

We then could conjecture that $$\inf _g\{ \Ent (X,g) ^2 \Vol (X,g)\} = 2 \pi (|\chi (S_1)| + |\chi (S_2)|).$$
\end{remark}

Notice that this construction can be generalized considering two closed hyperbolic $n$-manifolds $M_1$, $M_2$, each containing a closed totally geodesic 
codimension one hypersurface $N_1$ and $N_2$ such that $N_1$ is isometric to $N_2$. The complex $X$ obtained by gluing the $M_i$'s 
along $N_i$ satisfies the same conclusion as above, that is
$$
\Ent (X,g) ^n \Vol (X,g) \geq C(\delta, H).
$$

\bigskip
{\bf Surfaces with an attached circle} 

Let $S$ be a closed hyperbolic surface and $C$ be a circle. Given points $y\in S$ and $z \in C$, let $X$ be the complex defined as the disjoint
union of $S$ and $C$ with $y$ and $z$ being identified. We denote by $X = S \vee C$ and $\widetilde S$, $\widetilde C$, the universal cover of $S$, $C$ respectively. The universal cover of $X$ 
is a "tree of hyperbolic plane" and can be described inductively as follows.
To each lift $\tilde y \in \widetilde S$ of $y$, we attach a copy of $\widetilde C$ identifying $\tilde y $ with some $\tilde z \in \widetilde C$ any lift of $z$. Next, to each unidentified $\tilde z$ on some previously attached $\widetilde C$, we attach 
a copy of $\widetilde S$ at some free $\tilde y $ and we repeat these constructions, attaching to each new $\tilde y $ a copy of $\widetilde C$ at some $\tilde z$, etc....
The fundamental group $G:= \pi _1 (X)$ of $X$ is the free product $G = \Gamma \ast \mathbb Z$ where $\Gamma$ is the fundamental group of $S$. If $(s_1, \cdots, s_k)$ is a generating set 
of $\Gamma$ and $t$ a generator of $\mathbb Z$ and if $\Gamma$ is $\delta$-hyperbolic with respect to $(s_i)_{i=1\dots, k}$, then $G$ is $\delta$-hyperbolic with respect
to the generating set $(s_i, t)$. 

We consider a metric $g$ on $X$, defined by a metric $h$ on $S$ and a metric $h'$ on $C$, the latter being characterised only by the length of the circle  $C$. 
The volume $\Vol (X,g) = \Vol (S,h)$ does not depend on $(C,h')$ but only on the sub-complex  $S$ of $X$ of maximal dimension.

One sees that $X$ is $2$-essential, thus, by Theorem \ref{theo:poly} , we have $$\Ent (X,g) ^2 \Vol (X,g) \geq C( \delta, H)$$
where $H = \Ent (G, \Sigma)$ with respect to the generating set $\Sigma = \{s_1\dots, s_k, t\}$.

One weakness of the above lower bound comes from the fact that the lower dimensional simplex, namely $C$,
can be responsible for a large entropy of $X$ while it does not affect the volume of $X$. Indeed, one has, in particular, 
$$\lim_{l(C)\to 0} \Ent (X,g) ^2 \Vol (X,g) = + \infty$$ 
where $l(C)$ is the length of $C$.

On the other hand, we have 
$$\lim_{l(C)\to +\infty} \Ent (X,g) ^2 \Vol (X,g) = \Ent (S,h) ^2 \Vol (S,h)\,.$$

Notice that since $\Vol (X,g) = \Vol (S,h)$ and $\Ent (X,g) \geq \Ent (S,h)$, one has
$$\Ent (X,g) ^2 \Vol (X,g) \geq \Ent(S,h) ^2 \Vol (S,h)\,,$$
so that 
$$ \inf _g \{\Ent (X,g) ^2 \Vol (X,g)\}  =  \inf _h \{\Ent (S,h) ^2 \Vol (S,h) \}= \Ent (S,\mathrm{hyp}) ^2 \Vol (S, \mathrm{hyp})$$
where $\mathrm{hyp}$ denotes any hyperbolic metric on $S$.

We remark that for any simplicial complex $M$ of dimension $m$ and any simplicial complex $N$ of dimension $n<m$, 
then, $X = M \vee N$ satisfies
$$ \inf _g \{\Ent (X,g) ^m \Vol (X,g)\}  = \inf _h \{\Ent (M,h) ^m \Vol (M,h)\}\,.$$

\bigskip
{\bf A $3$-dimensional polyhedron}

We now describe an example of a finite $3$-essential polyhedron of dimension $3$ with torsion-free non elementary Gromov-hyperbolic fundamental group. It will be clear from the construction that this polyhedron 
is not a manifold,

This example comes from a general construction of Gromov hyperbolic Coxeter groups which associates to each hyperbolic $\CAT(0)$-cube complex another hyperbolic $\CAT(0)$-cube complex with higher homological dimension (see, \cite{Os}). 

The fundamental result underlying this construction is the following Theorem \ref{Os} of D. Osajda.
In order to state it, we need the following classical fact. 
\begin{prop}[\cite{Da}\label{Dav}, Section 1.2, Proposition 1.2.3 and Appendix I, Proposition I.6.8]\label{Dav}
Let $X$ be a finite flag simplicial complex of dimension $d$, then there exists a $\CAT(0)$-cube complex $\widetilde Y$ of dimension $d+1$ whose vertices have link $X$ and a right-angled Coxeter group $G$ acting geometrically on $\widetilde Y$. 
When $X$ is moreover assumed to be $5$-large, then $\widetilde Y$ is Gromov hyperbolic.
\end{prop}
Recall that for a simplicial complex, being flag means that any finite set of vertices which are pairwise connected by edges ({\sl i.e.} generating a complete graph) spans a simplex. 
A flag simplicial complex $X$ is $5$-large if 
every $k$-cycle in $X$ has length $k$ larger than or equal to $5$, where a $k$-cycle is a closed combinatorial loop made of $k$ vertices and $k$ edges which is {\sl full}
as a subcomplex of $X$, \cite{Os}, 2.1 p. 355. Here, a subcomplex $Y$ of a simplicial complex $X$ is said to be {\sl full} if every subset of vertices of $Y$ contained in a simplex of $X$ is contained in
a simplex of $Y$. A cube complex is said to be locally $5$-large if the link of each of its vertices is $5$-large (see, \cite{Os}, 2.1 p. 355). 

\begin{theorem}[\cite{Os}, Main Theorem, p. 354]\label{main-theorem}\label{Os}
Let $X$ be a finite $5$-large simplicial complex such that $H^n (X, \mathbf Q) \neq 0$. Consider $G$ and $\widetilde Y$ the right-angled Coxeter group and the $\CAT(0)$-cube complex given in Proposition \ref{Dav}.
Then, $G$ is hyperbolic and there exists a finite index, torsion-free subgroup, $\Gamma \subset G$, such that the quotient $Y:= \widetilde Y \slash \Gamma$ is a compact cube complex satisfying $H^{n+1} (Y, \mathbf Q) \neq 0$.
\end{theorem}
We now briefly describe the different steps of the construction of our example before giving the details.
We start from a compact surface $S$ with a square complexe structure. We then \lq\lq thicken"  $S$ into a $3$-dimensional simplicial complex $X$. The compact surface
has to be chosen in such a way that this \lq\lq thickening" is $5$-large. Theorem \ref{Os} then produces
a $\CAT(0)$-cube complex $\widetilde Y$ of dimension $4$ and a torsion-free hyperbolic group $\Gamma$ acting properly discontinuously on $\widetilde Y$ with compact quotient $Y = \widetilde Y \slash \Gamma$ 
such that $H^3 (Y, \mathbf Q) \neq 0$. To end up, $Y$ can be retracted
onto a simplicial complex $M$ of dimension $3$. This $M$ will be our example, in particular it turns out to be $3$-essential.

In order to define the surface $S$, we consider $\Gamma _0 \subset PSL(2, \mathbf R)$ the subgroup generated
by the reflexions on the sides of a right-angled pentagon of the hyperbolic plane $\mathbb H ^2$. 
By Poincar\'e's Theorem, $\Gamma _0$ is a discrete group with fundamental domain a right-angled  pentagon  $P$ of $\mathbb H^2$ whose $\Gamma _0$-translates generate a tiling of $\mathbb H ^2$.
Its dual tiling is  by hyperbolic squares of angle $2 \pi / 5$ and is also $\Gamma _0$-invariant. Notice that this tiling of $\mathbb H ^2$ by squares defines a $5$-large square complex since the link of each vertex is a pentagon. 
We then pick a torsion-free subgroup $\Gamma _1$ of $\Gamma _0$ in such a way that $\mathbb H ^2 /\Gamma _1$ is a closed surface with a $\CAT(-1)$-square complex structure 
such that every homotopically non trivial edge-loop has length at least $5$.
We denote by $S :=\mathbb H ^2 /\Gamma _1$ this squared surface. The underlying square complex structure on $S = \mathbb H ^2 /\Gamma _1$ is therefore locally $5$-large since, by the latter condition, the closure of the union 
of the $5$ squares adjacent to a vertex in $\mathbb H ^2$ embeds in $S$ so that the link of every vertex of $S$ is also a pentagon.

We now \lq\lq thicken" the surface $S$ as defined by D. Osajda in \cite{Os}, definition 3.1, p.~357.
 The {\sl thickening} $X$ of the square complex $S$ consists in associating to each square of $S$ a standard $3$-simplex keeping track of the vertices and edges of the square. Formally, two more edges appear on the $3$-simplex which correspond to the diagonals of the square. When two adjacent squares meet along an edge we glue the two corresponding $3$-simplices
 along the corresponding edge. The new simplicial complex $X$ thus obtained is clearly of dimension $3$
 and is homotopy equivalent to $S$ (cf. \cite{Os}, Lemma 3.5), therefore $H^2(X, \mathbf Q) \neq 0$. Moreover
 $X$ is a $5$-large simplicial complex by Proposition 3.4 of \cite{Os}. 
 
 To conclude, we have obtained a $5$-large simplicial complex $X$ of dimension $3$ such that $H^2(X, \mathbf Q) \neq 0$.
We now apply Theorem \ref{Os} to it in order to get a $4$-dimensional $\CAT(0)$-cube complex $\widetilde Y$ 
with a proper and discontinuous action of a right-angled Coxeter group $G$ and a torsion-free hyperbolic subgroup $\Gamma \subset G$
 acting properly discontinuously on $\widetilde Y$ with compact quotient $Y = \widetilde Y\slash \Gamma$ 
such that $H^3 (Y, \mathbf Q) \neq 0$. 

The $4$-dimensional cube complex $Y$ retracts on a $3$-dimensional complex $M$. Indeed, recall that 
the thickening $X$ of the square complex $S$ has been constructed replacing each square by the $3$-simplex spanned on the vertices of the square. Therefore, the $3$-simplices of $X$ meet only along edges as the corresponding squares do.  The neighbourhood of every vertex of the $4$-cube $Y$ is a cone over $X$ and therefore  
two adjacent $4$-cubes meet along $2$-cubes, namely cubical cones over edges. The cube complex $\widetilde Y$ is thus a $4$-dimensional cube complex such that cubes of dimension $4$ are attached along $2$-dimensional faces.
As a consequence, the $3$-faces of $\widetilde Y$ are left free and therefore, after taking a barycentric subdivision of  $\widetilde Y$, we can perform a 
$G$-equivariant retraction of $\widetilde Y$ onto a subset $\widetilde M$ of its $3$-skeleton. Let us briefly describe this retraction. Before doing this, we define $\widetilde M$.

We consider a $4$-cube of $\widetilde Y$ and its standard decomposition in cubes of smaller dimension and we describe below a polyhedral subdivision into $4$-simplices. The barycentric subdivision is obtained inductively by 
adding to each $k$-face, ($k=1,2,3,4$), its barycenter and then coning from this new vertex onto the barycentric subdivision of the $(k-1)$-faces. Notice that each $4$-simplex $\sigma$ of this subdivision 
has a unique $3$-face $f(\sigma)$ whose intersection with the boundary of the cube is a $2$-face contained into 
a $2$-face of the standard decomposition of the cube. The union of all these $2$-faces cover the standard $2$-skeleton of the cube. Notice that the barycentre of $C$ is a vertex of $f(\sigma )$, for every $\sigma$. Now for a cube $C$, the union of these $3$-faces $f(\sigma)$, for $\sigma$ running through 
the set of all $4$-simplices of the barycentric subdivision of $C$, is denoted by $F(C)$. The intersection of $F(C)$ with the boundary of $C$ is the union of the $2$-faces of $C$. Since the cubes $C$ of $\widetilde Y$ meet along
their $2$-faces, their $3$-subcomplexes $F(C)$ naturally attach themselves and give rise to a simplicial subcomplex $\widetilde M$ of $\widetilde Y$ of dimension $3$. 

We first perform the retraction inside each $4$-cube of a fundamental domain of the action of $G$ and then extend it. For such a cube $C$ and a $4$-simplex $\sigma$ of the 
barycentric subdivision of $C$, we consider the map that sends all vertices of $\sigma ^{(0)}$ which are  not in $(f(\sigma))^{(0)}$ to the barycentre of $C$ which is a vertex of $f(\sigma )$; we then extend it linearly as a retraction of $\sigma$ onto $f(\sigma )$. We therefore get a simplicial retraction of $C$ onto $F(C)$ which is the identity on the $2$-skeleton of $C$. Since the gluing of the cubes only take place on their $2$-skeleton we can extend it to 
the union of the cubes of a fundamental domain for $G$. We finally extend the retraction of the fundamental domain of $G$ to an equivariant retraction of $\widetilde Y$ onto $\widetilde M$, this is made possible by the fact that
the action of $G$ is by cubical isometries (see \cite{Da} p.11), hence preserving all faces of the cubes and of their barycentric subdivisions.
A nice picture, one dimension less, is given in \cite{Bra-Mcc} on page 2290. The retraction described above, in the context of this article, consists in sending all yellow dots to the centre of the cube and extending it linearly.

Notice that, unlike $\widetilde Y$, the complex $\widetilde M$ may not be $\CAT(0)$, however $M := \widetilde M\slash \Gamma$ is a retract of $Y = \widetilde Y\slash \Gamma$. Note also that the $3$-dimensional simplicial complex $M$ is distinct from $X$.

We now show that $M$ is $3$-essential. Since $M$ is a retract of $Y= \widetilde Y\slash \Gamma$ and $H^3(Y, \mathbf Q) \neq 0$, we have
 $H^3 (M, \mathbf Q) \neq 0$. On the other hand, $\widetilde Y$ is a $\CAT(0)$-space, hence $\widetilde Y\slash \Gamma$ is a classifying space $K(\Gamma,1)$ of $\Gamma$ and so is $M$ since it is a retract of $Y$. 
We deduce from this discussion that the $3$-dimensional complex $M$ is essential since it is a classifying space of $\Gamma$ and $H^3(M, \mathbf Q) \neq 0$, hence $M$
is $3$-essential.

We therefore obtained a $3$-essential compact $3$-dimensional simplicial complex $M$ with hyperbolic fundamental group $\Gamma$. Notice that 
$\Gamma$ cannot be elementary since $H^3(M,\mathbf Q) \neq 0$. We thus can apply Theorem \ref{theo:poly} and for every metric $g$ on $M$ we have
$$
\Ent(M,g) ^3 \Vol (M,g) \geq C(\delta, H),
$$
where $H= \Ent(G, S)$ with $S$ some generating set of $G$ such that $(G,S)$ is $\delta$-hyperbolic.

\begin{remarks}
  \begin{itemize}
    \item[i)] Starting from the $4$-dimensional cube complex $Y$ we could apply the same construction and  get a $5$-dimensional complex for which Theorem \ref{theo:poly} would give a lower bound for its minimal entropy. Inductively, starting from one surface $S$ as above we could produce examples in any odd dimension. We could also change the surface and get another series of spaces. This construction is thus a powerful tool to produce lots of examples.
    \item[ii)] It would be interesting to compute the constant of hyperbolicity $\delta$ and the entropy $H$ of $(G, S)$. 
Notice that 
$\Gamma$ being a subgroup of the Coxeter group $G$, it is a subgroup of $GL(N,\mathbb R)$ for some $N$, and therefore the algebraic entropy of 
non virtually abelian subgroup can also be estimated in term of $N$ by \cite{Bre-Ge}.
     \item[iii)] The simplicial complex $M$ is a finite union of $3$-simplices, which raises the following question: is the minimal entropy achieved by the metric which is 
      the hyperbolic metric of maximal volume on each $3$-simplex ?
  \end{itemize}  
\end{remarks}

\bigskip
{\bf $\mathrm{CAT}(0)$-square complexes}

A square complex is a $2$-dimensional cube complex.
In \cite{KS} A. Kar and M. Sageev have proved the following 
\begin{theorem}\label{square}
Let $G$ be a finitely generated group acting freely on a $\CAT(0)$-square complex. Then either $G$ is virtually abelian or $\Ent (G) \geq 2^{1/10}$.
\end{theorem}
If we furthermore assume that $G$ is Gromov-hyperbolic, we deduce from this result and Theorem \ref{theo:sab2} the following bound on the minimal entropy of compact $\CAT(0)$-square complexes. 
\begin{theorem}\label{ent-square}
Let $X$ be a compact $\CAT(0)$-square complex with non virtually abelian, torsion-free and hyperbolic fundamental group such that $H^2 (X, {\textbf Z}) \neq 0$. Then, for every Riemannian metric $g$ on $X$, we have
$$
\Ent (X,g) ^2 \Vol (X,g) \geq C >0,
$$
where $C$ is a universal constant.

\end{theorem}
\begin{proof}
The universal cover of $X$ is a $\CAT(0)$ metric space, thus it is contractible and since  $H^2 (X, {\textbf Z}) \neq 0$, $X$ is $2$-essential. Moreover, the space $X'$ obtained from $X$ by adding a vertex in the middle of each square and joining it to each vertex of this square is a polyhedron which is $2$-essential. Each metric $g$ on $X$ yields a metric on $X'$, still denoted by $g$, and such that $\Vol (X',g)=\Vol (X, g)$ and $\text{\rm Marg} (X', g)=\text{\rm Marg} (X, g)$. 

Now, applying Theorem \ref{theo:sab2} to $(X',g)$, we get $\Vol (X,g) \geq C' \text{\rm Marg} (X,g)^2$, where $C' = C'(2)$. By the definition of the Margulis invariant $\text{\rm Marg} (X, g)$ and by Remark \ref{Margulisconstantbis}, 
for every $r > \text{\rm Marg} (X, g)$, there exists $x \in X$ such that $ \Sigma_r (x) $ generates a non virtually nilpotent subgroup, thus 
a non virtually abelian subgroup $ \Gamma_r (x) $ of the fundamental group $\Gamma$ of $X$. Theorem \ref{square} then shows that $\Ent \big(\Gamma_r (x),\Sigma_r (x) \big) \ge 
2^{1/10}$ hence, by Lemma \ref{comparentropi}, that
$$ r \Ent (X,g) \ge \big ( \Max_{\sigma \in \Sigma_r (x)} d(x, \sigma   x )\big)\,\Ent (X,g) \ge \Ent \big(\Gamma_r (x),\Sigma_r (x) \big) \ge 2^{1/10}  . $$
When $r$ goes to $\text{\rm Marg} (X, g)$, this proves that $\text{\rm Marg} (X, g) \Ent (X,g)  \ge  2^{1/10} $ and we deduce that
$\Ent (X,g) ^2 \Vol (X,g) \ge C' \,\text{\rm Marg} (X, g)^2 \Ent (X,g)^2 \ge 2^{2/10} C'$.
\end{proof}
It is important to notice that the bound given in this result neither depends on $H$ nor on $\delta$ as in the previous examples. The bound also neither depends on the square complex $X$ nor on the group $G$.

There are a lot of compact $\CAT(0)$-square complexes $Y$ with hyperbolic fundamental group and $H^2 (Y, {\textbf Z}) \neq 0$ and we describe some below. 

For example, consider any finite graph $X$ with girth larger than or equal to $5$ and without end point. From the assumption on the girth, the graph $X$ is $5$-large and, since $X$ has no end point, 
$H^1(X, \mathbf Q) \neq 0$. Therefore,  by Theorem \ref{Os} (relying on Proposition \ref{Dav}), there exists a Gromov hyperbolic $\CAT(0)$-square complex $\widetilde Y$ and a torsion-free group $\Gamma$ acting geometrically 
on $\widetilde Y$ such that $Y:= \widetilde Y\slash \Gamma$ is a compact locally $\CAT(0)$-square complex with torsion-free hyperbolic fundamental group $\Gamma$ and $H^2 (Y, \mathbf Q) \neq 0$.
Graphs such as above can be chosen, for example, among Ramanujan graphs $X=X^{p,q}$, for $p$, $q$ distinct prime numbers congruent to $1$ mod $4$,
(see, \cite{Lu-Ph-Sa}). By \cite{Lu-Ph-Sa}, p.262-263,  these graphs are $p+1$ regular Cayley graphs of $PSL(2, \mathbf Z \slash q \mathbf Z)$ or 
$PGL(2, \mathbf Z \slash q \mathbf Z)$.
By \cite{Lu-Ph-Sa} p. 263, when the Legendre symbol $\left(\frac{q}{p}\right) = -1$, these graphs $X^{p,q}$ are regular graphs of degre $p+1$ with $v=q(q^2 -1)$ vertices 
and whose girth satisfies $\mathrm{girth}(X^{p,q}) \geq 4\log _p ( \frac{q(q^2 -1)}{3} ) \geq 5$ for $q$ large enough.
Therefore, for every fixed $p\geq 2$ and $q$ large enough, 
the graph $X^{p,q}$ provide an infinite family of examples of compact square complexes as above. Notice that the same conditions also hold with different constants in the case when the Legendre symbol 
$\left(\frac{q}{p}\right) = 1$.

Hyperbolic buildings also provide another class of examples of compact $\CAT(0)$-square complexes $X$ with hyperbolic fundamental group and $H^2(X, \mathbf Q) \neq 0$. 
In \cite{Bo} M. Bourdon has considered a family of $2$-dimensional hyperbolic buildings $I_{pq}$ defined for integers $p$, $q$ as follows. Let $R$ be a regular hyperbolic $p$-gon with angle 
$\pi /2$. We associate now to the $p$-gon the complex of groups where for each edge $e$ of $R$ the corresponding group is $G_e = {\textbf Z}  / q {\textbf Z}$ and for each vertex $s$ of $R$, 
the group is $G_s = {\textbf Z}  / q {\textbf Z} \times {\textbf Z}  / q {\textbf Z}$, and where, for each pair of edges $e$ and $e'$ adjacent to a vertex $s$, there is a natural injective morphisms from $G_e$ and
$G_e'$ on each of the factors of $G_s$. The fundamental group $\Gamma _{pq}$ of this complex of groups has the following presentation, 
$\Gamma _{pq} = < s_1, .., s_p, \,| \, s_i ^q = \, [s_i , s_{i+1}] =1\, >$. It acts on a hyperbolic building $X$ whose chambers are the hyperbolic right-angled $p$-gons, the apartments the hyperbolic plane $\mathbb H ^2$
and the link at each vertex is the complete bipartite $q$-graph. Each $p$-gon in $I_{pq}$ can be subdivised by squares in the following way: join the center of each $p$-gon to the middle of each 
of its sides. This subdivision gives rise to a $\Gamma _{pq}$-invariant $\CAT(-1)$-square complex structure on $X=I_{pq}$ where each square has three right angles and one angle equal to $2\pi /p$.}

Let us consider a finite index torsion-free subgroup $\Gamma \subset \Gamma _{pq}$. Such subgroups do exist since the $\Gamma _{pq}$ are $\mathbf R$-linear (see \cite{Hag}, Corollary 1.2). As an example, the kernel of the surjective morphism from $\Gamma _{pq}$ into $(\mathbf Z \slash q \mathbf Z )^p$ is a subgroup of $\Gamma _{pq}$ of index $q^p$.

The hyperbolic building $X$ is  a 
$\CAT(-1)$ space hence $M:= X \slash \Gamma$ is a classifying space of $\Gamma$. In the sequel, we assume $q$ even and claim that $M$ is essential and in particular that it is $2$-essential.
Indeed, let us show that 
$H^2 (M, \mathbf Z) \neq 0$. Let $W := \Gamma _{p2}$ be the Coxeter group of the right-angled pentagon $P\subset \mathbb H^2$. Notice that $W \subset \Gamma_{pq}$ 
since $q$ is even and that $W$ 
stabilizes the totally geodesic chamber $I_{p2} = \mathbb H^2 \subset I_{pq}$. Since $\Gamma$ has finite index in $\Gamma _{pq}$, we deduce that $W\cap \Gamma$ has finite index in $W$,
hence $I_{p2} \slash (W\cap \Gamma)$ is a compact surface $\Sigma$ which is immersed in $X\slash \Gamma$. Since the building has dimension $2$, the surface $\Sigma$ defines a non trivial element 
in $H^2(M, \mathbf Z)$.

To sum up, $M$ is a compact locally $\CAT(0)$-square complex with torsion-free hyperbolic fundamental group and therefore we can apply Theorem \ref{square} 
which yields that, for every metric $g$ on $M$, we have
$$\Ent (M,g) ^2 \Vol (M,g) \geq C\,.$$

As in the previous examples of the $3$-dimensional polyhedron, we could conjecture that the hyperbolic metric achieves the minimal entropy of these hyperbolic buildings. Notice that, by an argument of A.~Katok (cf. \cite{Ka}), the hyperbolic metric realises the minimal entropy in its conformal class.

\subsection{Finiteness and Compactness Results}

\subsubsection{General Definitions and Results}

Let us recall that, on a fixed closed manifold $M$ of dimension $n\geq 3$, Einstein Riemannian metrics are the critical points of the functional,
$$S: g\longmapsto \frac{1}{\Vol (M,g)^{n-2\over n}}\int_M \scal (g)dv_g\,,$$
where $\scal (g)$ is the scalar curvature of the Riemannian metric $g$. Notice that $S$ is scale-invariant. In dimension $2$ and $3$ Einstein metrics have constant sectional curvature, hence the question of describing them makes sense in dimension greater than $3$.

We now call {\sl Einstein structure} an equivalence class of Einstein metrics  for the following relation: $(Mg)$ and $(N,h)$ are equivalent if there exists a diffeomorphism $\phi : M\longmapsto N$ which is an homothety, that is $\phi^*(h)=\lambda g$ for some positive real number $\lambda$. 

We are interested in the critical points of $S$ as well as its critical values. As an example let us mention that there exist examples for which the set of critical values is infinite and has an accumulation point at zero (see [Be], p.471-472, Add. 3). 

We recall that, when a group $\Gamma$ acts properly on a metric space $(X,d)$, the associated entropy of $(X,d)$ is computed with respect to any measure 
$\mu$ invariant by this action of $\Gamma$: indeed, when $\Gamma \backslash  X$ is compact, the entropy of $(X,d)$ does not depend on the 
choice of this measure (see subsection \ref{entropies}).

Let us recall that, to any pair $\delta_0, \e'_0 > 0$, corresponds the class $\text{\rm Hyp}_{\rm thick} (\delta_0, \e'_0)$ of non virtually cyclic 
groups which admit a proper action by isometries on some $\delta_0$-hyperbolic metric space $ (X, d_0)$ such that every 
torsion-free $\g \in \Gamma^*$ verifies $\ell(\g) \ge \e'_0$ (see Definition \ref{systoleaction}) and that $\text{\rm Hyp}_{\rm thick}$ is the union
of all the classes $\text{\rm Hyp}_{\rm thick} (\delta_0, \e'_0)$ for all pairs $\delta_0, \e'_0 > 0$ (see Definition \ref{systoleaction0}).\\
To the class of groups $\text{\rm Hyp}_{\rm thick} (\delta_0, \e'_0)$, one associates the positive universal constant $ r_0 = r_0 (\delta_0 , \e'_0)$ defined 
in \eqref{universalconstants3}.

\subsubsection{Finiteness and compactness results when the Ricci curvature is bounded from below}\label{appli:finiteness}

We now prove the main finiteness and compactness results.

\begin{defi}\label{Riemmetrics}
Given $n \in \N$ ($n \ge 2$), $D, K , i_0 > 0$ and $\delta_0 , \e'_0 > 0$, let ${\cal R}_{\rm univers}^{\delta_0, \e'_0} (n, K , D , i_0)$  (resp. ${\cal R}_{\rm univers}^{\infty} (n, K , D , i_0)$) be the set of Riemannian $n$-dimensional manifolds $(M,g)$ which verify the
hypotheses:
\begin{itemize}
\item[(i)] the fundamental group $\Gamma_M$ of $M$ is torsion-free and belongs to $\text{\rm Hyp}_{\rm thick} (\delta_0, \e'_0)$ (resp. to 
$ \text{\rm Hyp}_{\rm thick}$),
\item[(ii)] $\Ric_g \ge -(n-1) K^2 \cdot g$ and $\diam (M,g) \le D$,
\item[(iii)] the injectivity radius \emph{of its Riemannian universal cover} $(\widetilde M , \tilde g)$ is bounded from below by $i_0$,
\end{itemize}
\end{defi}

Notice that, in this definition, the property (iii) is verified  in particular when the geodesics of $(M,g)$ have no conjugate points.

In the sequel, for any integer $n \ge 2$ and every $\delta_0 , \e'_0, K , D > 0$, we define 
\begin{equation}\label{defS0}
S_0 = S_0  (\delta_0 , \e'_0, n, K , D) :=  \dfrac{\e'_0}{ 13 \delta_0 + \e'_0}\cdot \dfrac{1}{2 (n-1) K}\, e^{- 4 (n-1) K D}\,.
\end{equation}

\begin{theorem}\label{fini0}
For every integer $n \ge 2$, every $D, K , i_0 > 0$ and every $\delta_0, \e'_0 > 0$, there are only finitely many differentiable structures in ${\cal R}_{\rm univers}^{\delta_0, \e'_0} (n, K , D, i_0)$.
\end{theorem}

We now call {\sl Riemannian structure} an equivalence class of smooth ({\sl i.e.} $C^\infty$) Riemannian metrics for the following relation: $(Mg)$ and $(N,h)$ are equivalent if there exists a diffeomorphism $\phi : M\longmapsto N$ which is an isometry, that is $\phi^*(h)= g$. The space of Riemannian structures on a 
given manifold $M$ can thus be viewed as the quotient of the space ${\cal M} (M)$ of Riemannian metrics on $M$ by the action of the group ${\rm Diff} (M)$
of diffeomorphisms of $M$ (every $\varphi \in {\rm Diff} (M)$ acts as $ g \mapsto \varphi^*g$).
A convergence of a sequence of Riemannian structures is thus a convergence of a sequence of corresponding Riemannian metrics $(g_i)_{i \in \N}$ modulo 
re-parametrization by some sequence $(\varphi_i)_{i \in \N}$ of diffeomorphisms. However, in each case, we have to precise the chosen topology of 
${\cal M} (M)$, and the prescribed regularity of the diffeomorphisms involved and of the limit-metric.\\
More precisely we define:

\begin{defi}[(see \cite{AC}, precised in \cite{HH})]\label{Cconvergence} 
Given $n \in \N$ ($n \ge 2$) and $s \in (0 , 1) $, for any sequence $(M_i , g_i)_{i \in \N}$ of smooth compact $n$-dimensional Riemannian manifolds, the sequence of underlying Riemannian structures is said to 
converge in the $C^{0,s}$-topology  if there exists a smooth compact $n$-dimensional differentiable manifold $M$, a $C^{0,s}$ metric $g$ on $M$ and (for 
every sufficienly large integer $i$) $C^{1,s}$-diffeomorphisms  $\varphi_i : M \to  M_i $ such that $\varphi ^*_i g_i$ converges to $g$  in the
$C^{0,s}$ topology when $i \to + \infty$. This means that (for every $s' < s$) there exists a sub-atlas of the $C^{\infty}$ complete atlas of $M$ such that, 
in each chart $(U_k , \psi_k)$ of this sub-atlas, each coordinate of $\varphi ^*_i g_i$  (viewed as a function on the open subset $\psi_k (U_k)$ of 
$\R^n$) converges (when $i \to + \infty$) to the corresponding coordinate of $g$ with respect to the $C^{0,s'}$ norm for functions on $\R^n$, this 
norm being defined by 
$$ \Vert f \Vert_{C^{0,s'}} := \sup_{y \ne x}\frac{|f(y) - f(x)|}{|y-x|^{s'}}\, .$$ 
\end{defi}

In \cite{AC}, it is remarked that, transporting by $\varphi_i $ the $C^\infty$ differentiable structure of $M_i$ creates a $C^\infty$ differentiable structure on
the limit $M$ which a priori depends on $i$ and that we denote by ${\rm Diff}_{\varphi_i} (M)$. However, fixing some $j_0 \in \N$ sufficiently large and the 
corresponding $C^\infty$ differentiable structure ${\rm Diff}_{\varphi_{j_0}} (M)$ on $M$, approximating each $C^{1,s}$-diffeomorphisms  $\varphi_i \circ \varphi_{j_0}^{-1}  : M_{j_0}\to  M_{i}$ (with respect to the $C^{1,s}$-norm in the charts) by a $C^{\infty}$-diffeomorphism $\phi_{j_0 , i} : M_{j_0}\to  M_{i}$ and 
replacing the $C^{1,s}$-diffeomorphisms $\varphi_i : M \to  M_i $ by $\widetilde \varphi_i := \phi_{j_0 , i} \circ \varphi_{j_0}$, one obtains 
$C^{\infty}$-diffeomorphisms $\widetilde \varphi_i$ from $M$ (endowed with the differentiable structure ${\rm Diff}_{\varphi_{j_0}} (M)$) onto $M_i$
such that $\widetilde \varphi_i^* g_i$ still converges to the metric $g$ in the $C^{0,s}$-topology. Notice that all the differentiable structures on $M$ 
obtained by transporting by $\widetilde \varphi_i $ the $C^\infty$ differentiable structure of $M_i$ coincide with ${\rm Diff}_{\varphi_{j_0}} (M)$ and this
fixes the $C^\infty$ differentiable structure of $M$.

\smallskip
Denoting by \lq \lq can" the canonical metric of the simply connected $n$-space with constant sectional curvature $-K^2$ and by $\mathbb B_{K} (R )$ one of 
its balls of radius $R$ (they are all isometric), we have the

\begin{theorem}\label{compacite0}
For every $n \in \N$ ($n \ge 2$), every $s \in ( 0 , 1 ) $ and every $D, K , i_0 > 0$,
on each compact $n$-dimensional manifold $M$ whose fundamental group belongs to $\text{\rm Hyp}_{\rm thick}$, the set of 
Riemannian structures corresponding to ($C^\infty$) Riemannian metrics $g$ such that $(M, g) \in {\cal R}_{\rm univers}^{\infty} (n, K , D , i_0)$ 
(if non empty) has compact closure in the set of $C^{0,s}$ Riemannian structures on $M$. More precisely, every sequence $(g_i)_{i\in \N}$ of smooth Riemannian 
metrics on $M$ such that $(M, g_i) \in {\cal R}_{\rm univers}^{\infty} (n, K , D , i_0)$ admits a subsequence $(g_j)_{j}$ such that the underlying subsequence of Riemannian structures converges in the $C^{0,s}$-topology (see Definition \ref{Cconvergence}).
The limit space is the smooth manifold $M$, endowed with the $C^{0,s}$ Riemannian structure represented by a $C^{0,s}$ metric denoted by 
$g_\infty$.\\
Moreover, the limit-metric $g_\infty$ satisfies the following bounds :
\begin{itemize}
   \item[i)] $\diam (M, g_\infty)\leq D$, 
   \item[ii)] $(M, g_\infty)$ verifies the same Bishop-Gromov inequality as the $g_i$'s, namely, for $r \le R$,
$$\, \dfrac{\Vol_{g_\infty}  B_M(x , R )}{\Vol_{g_\infty}  B_M(x , r )} \le \dfrac{\Vol_{\rm can}  \mathbb B_{K} (R )}{\Vol_{\rm can}  \mathbb B_{K} ( r )}\,,$$
   \item[iii)] for every $\tilde x \in \widetilde M$ and $s' < s$, the function $d_{\tilde g_\infty} (\tilde x, \bullet)^2$ is $C^{1,s'}$ on every ball of 
$(\widetilde M , \tilde g_\infty)$, centered at $\tilde x$ and of radius $< i_0$; in this sense, the injectivity radius of $(\widetilde M , \tilde g_\infty)$ is $\ge i_0$ (see Remark \ref{analogues} (2)).
\end{itemize}
\end{theorem}

\begin{remarks}\label{analogues}
\emph{
\begin{itemize}
\item[(1)] Our proof uses the work \cite{AC}. However, the main difference between Theorems \ref{fini0} and \ref{compacite0} and the results in 
\cite{AC} is the fact that \cite{AC} assumes the injectivity radii of the compact Riemannian manifolds $(M,g)$ under consideration to be  
bounded from below by the universal constant $i_0$ while we only assume here that the injectivity radii of their Riemannian universal cover $(\widetilde M , \tilde g)$ is bounded from below by $i_0$. This last hypothesis is much weaker than the previous one. Indeed, for example, it is verified
by every compact Riemannian manifold $(M,g)$ whose sectional curvature is non-positive and, more generally, by every compact Riemannian manifold $(M,g)$ without 
conjugate points. In these two cases, the Riemannian universal cover $(\widetilde M , \tilde{g})$ of $(M,g)$ has infinite injectivity radius,
while the injectivity radius of $(M,g)$ itself could be arbitrarily small.
\item[(2)] For any $C^{\infty}$-metric on $M$ the injectivity radius $\inj (\tilde x)$ of $(\widetilde M , \tilde{g})$ at the point $\tilde x$ is defined as the 
supremum of the positive values of $r$ such that the geodesic ball $B_{(\widetilde M , \tilde{g})} (\tilde x , r)$ verifies one of the two following equivalent conditions:
\begin{itemize}
\item[(a)] the function $d_{\tilde{g}} (\tilde x, . )^2$ is smooth when restricted to $B_{(\widetilde M , \tilde{g})} (\tilde x , r)$,
\item[(b)] the exponential map is injective from the ball of radius $r$ in the tangent space $T_{\tilde x} \widetilde M$ onto $B_{(\widetilde M , \tilde{g})} (\tilde x , r)$.
\end{itemize}
In the case where the limit-metric is $C^{0,s}$, only definition (a) makes sense (where smoothness is replaced by $C^{1,s'}$  for every $s' < s$), because the exponential map is not
correctly defined; indeed, two geodesics which coincide on any interval $[0, t]$ may branch at time $t$ (an example is given
by the piecewise $C^1$ metric obtained by gluing two copies of $\R^2\setminus \mathbb B^2$ on their boundary $S^1 = \partial B^2$).
\item[(3)] As the limit-metric $g_{\infty}$ is only $C^{0,s}$, the condition \lq \lq Ricci curvature bounded from below" makes no sense  at the limit, we must
thus find an analogous condition. Following the viewpoint of J. Cheeger and T. Colding described in their works about the structure of Riemannian manifolds with Ricci-curvature bounded below, the analogous assumption that we use is the Bishop-Gromov inequality. Hence,
 in Theorem \ref{compacite0}, we show that the limit spaces $(M,g_\infty)$ all verify this Bishop-Gromov inequality.
\end{itemize}
}
\end{remarks}

\begin{proof}[Proof of Theorems \ref{fini0} and \ref{compacite0}]
Consider any $(M, g) \in {\cal R}_{\rm univers}^{\delta_0, \e'_0} (n, K , D , i_0)$ and its Riemannian universal covering 
$\pi : (\widetilde M , \tilde g) \to (M,g)$ and (for sake of simplicity) denote by $\Gamma$ the fundamental group of $M$, viewed as the group of 
deck-transformations of this Riemannian universal covering. Let us first prove that 
\begin{equation}\label{minorinj}
\forall (M, g) \in {\cal R}_{\rm univers}^{\delta_0, \e'_0} (n, K , D , i_0) \ \ \ \ \ \ \inj (M,g) \ge \Min \left( i_0 ; S_0  ( \delta_0 , \e'_0, n ,K , D)\right) \, .
\end{equation}
Indeed, let $S_0 = S_0  ( \delta_0 , \e'_0,n, K , D)$ for the sake of simplicity, for every $\tilde x \in \widetilde M$ and every positive 
$r \le \frac{S_0}{2}$, if we set $ x := \pi (\tilde x )$,
$\pi$ maps $ B_{\tilde g}(\tilde x , r)$ into $B_g (x,r) $. Moreover $\pi$ is surjective from $ B_{\tilde g}(\tilde x , r)$ onto 
$B_g (x,r) $ because, for every $y \in B_g (x,r)$, as $\inf_{\tilde y \in \pi^{-1} (y)} d_{\tilde g} (\tilde x , \tilde y)  = d_g  \big(\pi (\tilde x ), y ) \big)<r\, $,
there exists $\tilde y \in \pi^{-1} (y)$ such that $d_{\tilde g} (\tilde x , \tilde y) < r$. Let $\tilde y, \tilde z$ be any two points of $ B_{\tilde g}(\tilde x , r)$
such that $\pi(\tilde z) = \pi (\tilde y)$, there then exists some $\g \in \Gamma$ such that $\tilde z = \g \tilde y$. Theorem \ref{minorsystglobale} (ii) (and 
the fact that $\Ent (M , g) \le (n-1) K $ by Bishop-Gromov's comparison Theorem) implies that 
$\inf_{p \in M} \text{\rm sys}_{\Gamma}^{\diamond} (p) > S_0  ( \delta_0 , \e'_0, n , K , D) $, hence, if $\g $ is non trivial, that
$d_{\tilde g} (\tilde y , \tilde z) = d_{\tilde g} ( \tilde y, \g \tilde y) > S_0 \ge 2\, r$, in contradiction with the fact that 
$\tilde y, \tilde z \in B_{\tilde g}(\tilde x , r)$. Consequently $\g$ is trivial and $\tilde z = \tilde y$. It follows that $\pi$ is bijective from $ B_{\tilde g}(\tilde x , r)$ onto 
$B_g (x,r) $.\\
Every pair $c_1 , c_2$ of normal (locally minimizing) geodesics of $(M,g)$ issued from $x $ with distinct initial unit speeds can be lifted as a pair of normal (locally minimizing) geodesics $\widetilde{c_1} , \widetilde{c_2}$ of 
$ (\widetilde M , \tilde g)$ issued from 
$\tilde x$ with distinct initial unit speeds. As the injectivity radius of $(\widetilde M , \tilde g)$ is bounded from below by $i_0$, we have 
$\widetilde{c_1}(t) \ne \widetilde{c_2} (s)$ for every $t,s \in \, ]0 , i_0[$, thus for every $t,s \in \, ]0 , \e'_0[$, where 
$\e'_0 := \Min \left( i_0 , S_0 \right) $; as $\pi$ is bijective from $ B_{\tilde g}(\tilde x , \e'_0)$ onto $B_g (x,\e'_0) $, we have
$\pi \circ \widetilde{c_1}(t) \ne \pi \circ \widetilde{c_2} (s)$, and thus $c_1 (t) \ne c_2 (s)$ for every $t,s \in \, ]0 , \e'_0[$. This proves \eqref{minorinj}.

\smallskip
Now, by the Bishop-Gromov's comparison Theorem, we have
\begin{equation}\label{majorvol}
\forall (M, g) \in {\cal R}_{\rm univers}^{\delta_0, \e'_0} (n, K , D , i_0) \ \ \ \ \ \ \Vol (B_g( x , D))\le \Vol (B_{\tilde g} (\tilde x, D))\le 
\Vol_{\rm can}  \mathbb B_{K} (D) \, .
\end{equation}
As all the Riemannian $n$-manifolds $(M, g) \in {\cal R}_{\rm univers}^{\delta_0, \e'_0} (n, K , D , i_0)$ are compact, verify $\Ric_g \ge -(n-1) K^2 \cdot g$,
have volume uniformly bounded from above (by \eqref{majorvol}) and injectivity radius uniformly bounded from below (by \eqref{minorinj}), they all satisfy 
the hypotheses of the compactness theorem proved in \cite{AC} (Theorem~0.2).\\ 
This result proves that ${\cal R}_{\rm univers}^{\delta_0, \e'_0} (n, K , D , i_0)$ contains a finite number of differentiable structures, and hence 
proves Theorem \ref{fini0}.\\
On any fixed compact $n$-dimensional manifold $M$ the set of Riemannian structures corresponding to metrics $g$ such that $(M,g) \in {\cal R}_{\rm univers}^{\infty} (n, K , D , i_0)$ is 
either empty or precompact: indeed, as the fundamental group $\Gamma$ of $M$ is fixed, if this group belongs to $\text{\rm Hyp}_{\rm thick}$, it belongs to
$\text{\rm Hyp}_{\rm thick} (\delta_0, \e'_0)$ for some fixed $\delta_0\geq 0$ and $\e'_0 >0$. It follows that the set of metrics 
$g$ such that $(M,g) \in {\cal R}_{\rm univers}^{\infty} (n, K , D , i_0)$ is included in the set of metrics $g$ such that 
$(M,g) \in {\cal R}_{\rm univers}^{\delta_0, \e'_0} (n, K , D , i_0)$. A consequence is that all the metrics $g$ such that 
$(M,g) \in {\cal R}_{\rm univers}^{\infty} (n, K , D , i_0)$ verify $\Ric_g \ge -(n-1) K^2 \cdot g$, have volume uniformly bounded from above 
(by \eqref{majorvol}) and injectivity radius uniformly bounded from below (by \eqref{minorinj}). We can thus again apply the compactness theorem in \cite{AC} (Theorem~0.2), which proves that the set of Riemannian structures corresponding to metrics $g$ on $M$ such that 
$(M,g) \in {\cal R}_{\rm univers}^{\infty} (n, K , D , i_0)$ (if non empty) is precompact in the $C^{0,s}$ topology for every $s \in (0 , 1)$: this means
that, for every sequence of metrics $(g_i)_{i \in \N}$ on $M$ such that $(M,g_i) \in {\cal R}_{\rm univers}^{\infty} (n, K , D , i_0)$, there exists a
subsequence $(g_j)_{j}$ and diffeomorphisms $\varphi _j : M \to M$ such that the sequence of metrics $\varphi ^*_j g_j$ converges (in the $C^{0,s'}$ topology,
for every $s' <s$) to some $C^{0,s}$ Riemannian metric $g_{\infty}$.

\smallskip
Let us now prove that the Riemannian manifold $ (M, g_{\infty})$ verifies properties analogous to the ones satisfied by the Riemannian manifolds $ (M, g)$
lying in ${\cal R}_{\rm univers}^{\infty} (n, K , D , i_0)$. As the fundamental group $\Gamma$ of $M$ is torsion-free and belongs to 
$ \text{\rm Hyp}_{\rm thick}$, we only have to prove that the diameter is bounded from above, that the injectivity radius is bounded from below
(by constants comparable to $D$ and $i_0$ respectively) and that $ (M, g_{\infty})$ satisfies the Bishop-Gromov inequality. 
The convergence of the metrics $(\varphi ^*_j g_j)_{j \in \N}$ being $C^{0,s}$, it is a fortiori $C^{0}$, and there thus exists a strictly positive sequence 
$(\e_j)_{j \in \N}$ (such that $\lim_{j \to +\infty} \e_j = 0$) which verifies
\begin{equation}\label{majormetric}
(1-\varepsilon _j)^2 g_\infty \le \varphi ^*_j g_j \le (1+\varepsilon _j)^2 g_\infty \, .
\end{equation}
A consequence of \eqref{majormetric} is that $(1-\varepsilon _j)d_{g_\infty} \le d_{\varphi ^*_jg_j} \le (1+\varepsilon _j) d_{g_\infty}$, which implies
firstly that 
$$\diam (g_\infty) = \lim_{j\to +\infty} \diam (\varphi ^*_j g_j) \le D~,$$
and secondly that (for every $p \in M$)
$$ B_{\varphi ^*_jg_j } \big( p , (1-\varepsilon _j) R \big)  \subset  B_{g_\infty} (p , R) \subset B_{\varphi ^*_jg_j } \big( p , (1+\varepsilon _j) R \big) \, .$$
These inclusions and the comparison $ (1-\varepsilon _j)^n \le \dfrac{dv_{\varphi ^*_j g_j }}{dv_{g_\infty}} \le (1+\varepsilon _j)^n$ between the 
Riemannian measures (coming from \eqref{majormetric}) yield
\begin{equation*}
\dfrac{\Vol_{\varphi ^*_j g_j } \left( B_{\varphi ^*_jg_j } \big( p , (1-\varepsilon _j) R \big) \right)}{(1+\varepsilon _j)^{n}}  \le  \Vol_{g_\infty} \big( B_{g_\infty} (p , R)\big) \le  \dfrac{\Vol_{\varphi ^*_j g_j } \left( B_{\varphi ^*_jg_j } \big( p , (1+\varepsilon _j) R \big) \right)}{(1-\varepsilon _j)^{n}} \, ,
\end{equation*}
hence that, for every $r, R > 0$ such that $r < R$,
$$ \left(\dfrac{1-\varepsilon _j}{1+\varepsilon _j}\right)^n \dfrac{\Vol_{g_\infty} \big( B_{g_\infty} (p , R)\big)}{\Vol_{g_\infty} \big( B_{g_\infty} (p , r)\big)}
 \le \dfrac{\Vol_{\varphi ^*_j g_j } \left( B_{\varphi ^*_jg_j } \big( p , (1+\varepsilon _j) R \big) \right)}{\Vol_{\varphi ^*_j g_j } \left( B_{\varphi ^*_jg_j } \big( p , (1-\varepsilon _j) r \big) \right)} =  \dfrac{\Vol_{g_j } \left( B_{g_j } \big( \varphi_j ( p) , (1+\varepsilon _j) R \big) \right)}{\Vol_{ g_j } \left( B_{g_j } \big(\varphi_j ( p) , (1-\varepsilon _j) r \big) \right)}$$
$$ \le \dfrac{\Vol_{\rm can}  \mathbb B_{K} \big((1+\varepsilon _j) R \big)}{\Vol_{\rm can}  \mathbb B_{K} ( (1-\varepsilon _j) r )}\, ,$$
where the last inequality is a consequence of the classical Bishop-Gromov inequality (recalling that $\Ric_{g_j } \ge 
- (n-1) K^2 \cdot g_j $). Taking the limit when $j \to +\infty$, we get
$$\dfrac{\Vol_{g_\infty} \big( B_{g_\infty} (p , R)\big)}{\Vol_{g_\infty} \big( B_{g_\infty} (p , r)\big)} \le \dfrac{\Vol_{\rm can}  \mathbb B_{K} (R)}{\Vol_{\rm can}  \mathbb B_{K} (r )}\, ,$$
which is exactly the Bishop-Gromov inequality verified by all the Riemannian manifolds $ (M, g)$ such that $\Ric_g \ge (n-1) K^2 \cdot g$.\\
Let us now define $\varrho_\infty :=  d_{\tilde g_\infty}(\tilde x, .)^2$ and verify that, given $i < i_0$, $\varrho_\infty$ is $C^{1,s'}$ in the ball 
$ B_{(\widetilde M , \tilde g_\infty )} (\tilde x, i)$, that we shall denote by $\widetilde B (\tilde x , i)$ for the sake of simplicity: it is the Property (5) of 
Theorem 0.1 of \cite{AC}, which is proved to be valid for the limit-metric $g_\infty$ in Theorem 0.2 of \cite{AC}. Applying roughly this Theorem seems 
to only prove that $\varrho_\infty$ is $C^{1,s'}$ in every ball of radius $\frac{1}{2} i'_0$, where $i'_0 := \Min \left( i_0 ; S_0 \right)$ is a lower bound of 
the injectivity radii of all the manifolds $(M,g) \in {\cal R}_{\rm univers}^{\infty} (n, K , D , i_0)$ given by \eqref{minorinj}.\\ 
However, the remark (2) p. 267 of \cite{AC} proves that the conclusions of Theorem 0.2 of \cite{AC} are still valid on compact domains of 
$(\widetilde M , \tilde g_\infty)$; this proves\footnote{As the elements of this last proof are distributed in the whole of \cite{AC}, we shall 
summarize them here: from \cite{AC} (Theorem 0.1 and discussion 
pp.266-267), it follows that it is possible to choose the diffeomorphisms $\varphi_j$, a system of $g_\infty$ harmonic charts $\{(U_\ell ,H_\ell)\}_{\ell \in I}$ 
of $\widetilde M$ (with values in the euclidean ball $\mathbb B (r)$ with fixed radius $r$) and, defining $h_j = \varphi ^*_jg_j$, a system of $\tilde h_j$-harmonic charts $\{(U_{j,\ell};H_{j,\ell})\}_{\ell \in I}$ of $\widetilde M$ (with values in $\mathbb B (r)$) such that
$\big\{H^{-1}_\ell [\mathbb B(r/ 2)] \cap H_{j,\ell}^{-1} [\mathbb B(r/ 2)]\big\}_{\ell \in I}$ is still a covering of $\widetilde M$ and moreover that 
$H_{j,\ell} \circ H^{-1}_\ell$ converges (in the $C^{1,s'}$ sense for every $s' <s$) to the canonical injection $\mathbb B(r/ 2)\hookrightarrow \mathbb B(r)$.
Defining $\varrho_j = d_{\tilde h_j} (\tilde x,\cdot )^2$ and $\varrho_\infty  = d_{\tilde g_\infty}(\tilde x, \cdot )^2$ and choosing 
$i_1 \in ]i , i_0[$, the proof of the lemma 1.4 of \cite{AC} establishes that the $\tilde h_j$-Laplacian of $\varrho _j $ is bounded (from above and from below) on the ball 
$B_{(\widetilde M , \tilde g_\infty)} (\tilde x , i_1)$ (whose closure is included in $ B_{(\widetilde M , \tilde g_j)} (\tilde x , i_0)$ for $j$ large enough by \eqref{majormetric}), hence that 
$\varrho_j \circ H^{-1}_{j,\ell}$, whose euclidean Laplacian is bounded (because the coordinates are harmonic), has bounded $H^p_2$-norm (for every $p$) 
on $\mathbb B(r/ 2) \cap H_{j,\ell} \big( \widetilde B (\tilde x , i )\big)$ (see \cite{AC}, formula (0.10)); hence a 
subsequence of the sequence $\big(\varrho_j \circ H^{-1}_{j,\ell}\big)_j$ converges 
(in $C^{1,s'}$-norm on $\mathbb B(r/ 2)$) to a limit-function $f_\ell$, which is $C^{1,s'}$. 
As $\varrho _j $ converges to $\varrho_\infty$ in the $C^0$ sense (on every compact set) by \eqref{majormetric}, recalling that, for any 
$\tilde y \in H^{-1}_\ell [\mathbb B(r/2)]\cap H^{-1}_{j,\ell} [\mathbb B(r/2)] \cap \widetilde B (\tilde x , i)$, the sequences
$ \big(\varrho _j \circ H^{-1}_{j,\ell} [H_{j,\ell}(\tilde y)]\big)_j$ and $\big(f_\ell \circ (H_{j,\ell} \circ H^{-1}_\ell) [H_\ell(\tilde y)]\big)_j$ respectively converge to $ f_\ell [H_{j,\ell} (\tilde y)]$ and $ f_\ell[H_\ell(\tilde y)]$, using the triangle inequality and taking the limit, one obtains that $\varrho_\infty = f_\ell \circ H_\ell$ in
each chart, thus that $\varrho_\infty$ is $C^{1,s'}$ on $\widetilde B (\tilde x , i)$.} that $\varrho_\infty$  
is $C^{1,s'}$ in the ball $\widetilde B (\tilde x , i)$ and this ends the proof.
\end{proof}

\subsubsection{Finiteness and compactness Results for Einstein Structures}\label{Einsteindiscrete}

We now wish to apply the previous results to the study of Einstein structures.

Let us recall that, given any parameters $\delta_0, \e'_0 > 0$, the universal constant $ r_0 = r_0 (\delta_0 , \e'_0)$ is defined in \eqref{universalconstants3}.

\begin{theorem}\label{Einstein1}
Given any $K , \delta_0, \e'_0 > 0$, let us consider the set of compact $n$-dimensional manifolds, whose fundamental group is torsion-free and belongs 
to $\text{\rm Hyp}_{\rm thick}  (\delta_0, \e'_0)$ and which admit an Einstein metric $g$ without conjugate point and such that 
$\,{\rm scal} (g)  \diam (g)^2 \ge - n (n-1) K^2$,
\begin{itemize}
\item[(i)] for every $K > 0$, this set is finite,
\item[(ii)] if $K \le \frac{ r_0 (\delta_0 , \e'_0 )}{2 (n-1)}$, this set is empty.
\end{itemize}
\end{theorem}

Conclusion (ii) remains valid if one does not assume that the metric is without  conjugate point.
Conclusion (i) remains valid if the hypothesis \lq \lq $g$ without  conjugate point" is replaced by the following weaker hypothesis: there exists a constant 
$i_0 > 0$ such that, for all the metrics of the set defined in Theorem \ref{Einstein1}, the injectivity radius $\inj (\tilde g)$ of the pulled-back 
metric $\tilde g := \pi^* g$ on the universal cover verifies $\dfrac{\inj (\tilde g)}{\diam (g)} \ge i_0$. Notice that, when $(M,g)$ has no conjugate points, 
then the injectivity radius of its Riemannian universal covering $(\widetilde M , \tilde g)$ is infinite.

The notion of $C^k$-convergence for sequences of Riemannian metrics on $M$ being defined either in a fixed system of local charts 
(the derivatives being then given by the partial derivatives with respect to the coordinates in the charts) or referring to a given Riemannian metric $g_0$ 
on $M$, the $C^k$-convergence is then defined by the following norm on the space of smooth symmetric $2$-tensors: 
$\Vert h \Vert_{C^{k}(g_0)} := \sum_{i = 0}^k \, \Vert (D^{g_0})^i  h \Vert_{g_0}$ (where $D^{g_0}$ is the Levi-Civita covariant derivative associated to
the Riemannian metric $g_0$).

\begin{defi}[see \cite{HH}]\label{Eisteinconverge} Given any compact smooth manifold $M$ and any $k \in \N$, for any sequence $(g_i)_{i \in \N}$ of 
smooth Riemannian metrics, the sequence of underlying Riemannian structures is said to converge in the $C^{k}$-topology if there exists a smooth 
Riemannian metric $g$ on $M$ and, for every sufficienly large integer $i$, $C^{\infty}$-diffeomorphisms  $\varphi_i : M \to  M$ such that $\varphi^*_i g_i$ $C^{k}$-converges to $g$ when $i \to + \infty$.
\end{defi}

\begin{defi}\label{Einsteinmetrics} On any compact manifold $M$, we denote by ${\rm Einst}_M (K,D, i_0) $ the set of Einstein metrics $g$, verifying $\,{\rm scal} (g) \ge - n (n-1) K^2$ (where $n = \dim M$), $\diam (g) \le D$ and such that the injectivity radius \emph{of its Riemannian universal cover} $(\widetilde M , \tilde g)$ is bounded from below by $i_0$.
\end{defi}

Notice that the set of  Einstein metrics $g$ on $M$ without conjugate point and verifying $\,{\rm scal} (g) \ge - n (n-1) K^2$ and $\diam (g) \le D$ is a subset 
of ${\rm Einst}_M (K,D, i_0) $; hence all properties verified on ${\rm Einst}_M (K,D, i_0) $ are verified on this set of metrics.\\
The set of Riemannian structures corresponding to metrics in ${\rm Einst}_M (K,D, i_0) $ has nice compactness properties, summarized in the following results:

\begin{theorem}\label{Einstein2}
On every compact $n$-dimensional manifold $M$ whose fundamental group is torsion-free and belongs to $\text{\rm Hyp}_{\rm thick}$, and for 
every $K , D , i_0> 0$ and every $k \in \N$, the set of Riemannian structures corresponding to metrics belonging to ${\rm Einst}_M (K,D, i_0) $ is (sequentially) 
compact in the $C^k$-topology for every $k \in \N$, {\sl i. e.} from every sequence $(g_i)_{i \in \N}$ belonging to ${\rm Einst}_M (K,D, i_0) $, one can extract a subsequence $(g_j)_{j \in J_k \subset \N}$ such that 
the underlying sequence of Riemannian structures converges, in the $C^{k}$-topology, to some Riemannian structure represented by a metric $g_\infty$ 
which belongs to ${\rm Einst}_M (K,D, i_0) $.
\end{theorem}

Notice that, if the convergent subsequence $(g_j)_{j \in J_k \subset \N}$ and the corresponding diffeomorphisms
are chosen randomly, then the limit $g_\infty$ depends on the value of $k$. Nevertheless, one gets the

\begin{corollary}\label{Einstein3}
On every compact $n$-dimensional manifold $M$ whose fundamental group is torsion-free and belongs to $ \text{\rm Hyp}_{\rm thick}$, and for 
every $K , D , i_0> 0$, from every sequence $(g_i)_{i \in \N}$ of elements of ${\rm Einst}_M (K,D, i_0) $, one can extract a subsequence 
$(g_j)_{j \in J}$ and find diffeomorphisms $\phi_j$ such that $ \big(\phi_j^* g_j\big)_{j \in J} $ $C^\infty$-converges, when $j$ goes to $+\infty$, to some
metric $g$ which belongs to ${\rm Einst}_M (K,D, i_0) $.
\end{corollary}

\begin{proof}[Proof of Theorem \ref{Einstein1}]
For sake of simplicity, we denote by $\Gamma$ the fundamental group of $M$, notice that it is torsion-free by assumption. Applying Corollary \ref{cor:positiveentropy} (iii) with $(Y,d, \mu)=(\widetilde M, d_{\tilde g}, dv_{\tilde g})$ we get that $\Ent (M, g):=\Ent (\widetilde M, d_{\tilde g}, dv_{\tilde g})>0$. Since the metric $g$ is Einstein its Ricci curvature satifies $\Ric_g = \frac{1}{n} \, {\rm scal} (g) \cdot g $, where the scalar curvature ${\rm scal} (g) $ is constant, hence a real number. Now, Bishop-Gromov's 
Inequality shows that $\Ent (M , g) \le \sqrt{\frac{n-1}{n} \,{\rm scal}_- (g)}$, where  ${\rm scal}_-(g) = 
\Max \left(- {\rm scal} (g)\, , \, 0\right)$. The positiveness of the entropy of $(M, g)$ then implies that ${\rm scal}(g)$ is negative.
We may thus apply Theorem \ref{transyst} (iii) (where we replace $H$ by $\sqrt{\frac{n-1}{n} \,{\rm scal}_- (g)}\,$),
which proves that there exists $\tilde p \in \widetilde M$ such that 
\begin{equation}\label{majorscal}
4 \, {\rm scal}_- (g) \cdot \diam (g)^2 \ge {\rm scal}_- (g) \cdot\text{\rm sys}_\Gamma (\tilde p)^2 = {\rm scal}_- (g) \cdot \sys^{\diamond}_\Gamma 
(\tilde p)^2 \ge \frac{n}{n-1} \, r_0 (\delta_0 , \e'_0)^2 \, .
\end{equation}
This implies  (ii) since it proves that $ n (n-1) K^2 \ge - {\rm scal} (g)  \diam (g)^2 \ge    \frac{n}{4 (n-1)}\, r_0 (\delta_0 , \e'_0 )^2$.

\smallskip
For sake of simplicity, we now let ${\cal M}_K $ be the set of compact $n$-dimensional manifolds, whose fundamental group is torsion-free and belongs 
to $\text{\rm Hyp}_{\rm thick}  (\delta_0, \e'_0)$ and which admit an Einstein metric $h$ without conjugate point and such that 
$\,{\rm scal} (h)  \diam (h)^2 \ge - n (n-1) K^2$. By rescaling, on every manifold $M  \in {\cal M}_K$ there exists 
some Einstein metric $g$, without conjugate point, such that $\diam (g) =1$ and $\,{\rm scal} (g) \ge - n (n-1) K^2$, and this yields $\Ric_g \ge - (n-1) 
K^2 \cdot g $. It follows that $(M,g) \in {\cal R}_{\rm univers}^{\delta_0, \e'_0} (n, K , 1 , +\infty)$ and thus, using Theorem \ref{fini0}, that there are only 
a finite number of differentiable structures in ${\cal M}_K $.
\end{proof}

\begin{proof}[Proof of Theorem \ref{Einstein2}]
For the sake of simplicity, let us denote by $\Gamma$ the fundamental group of $M$; as this group is fixed and belongs to $\text{\rm Hyp}_{\rm thick}$, 
it belongs to $\text{\rm Hyp}_{\rm thick} (\delta_0, \e'_0)$ for some fixed $\delta_0\geq 0$ and $\e'_0 >0$. By Definitions \ref{Einsteinmetrics} and \ref{Riemmetrics}, ${\rm Einst}_M (K,D, i_0) \subset 
{\cal R}_{\rm univers}^{\delta_0, \e'_0} (n, K , D , i_0)$ and it follows from Theorem \ref{compacite0} that the set of Riemannian structures corresponding 
to Riemannian metrics $g \in{\rm Einst}_M (K,D, i_0)$ (if not empty) has compact closure in the set of $C^{0,s}$ Riemannian structures on $M$.

\smallskip
On the other hand, for every $ g \in {\rm Einst}_M (K,D, i_0)$, as $(M,g) \in {\cal R}_{\rm univers}^{\delta_0, \e'_0} (n, K , D , i_0)$, we may 
apply the inequalities \eqref{minorinj} and \eqref{majorvol}, which yield
\begin{equation}\label{minorinjvol}
\inj (g) \ge \Min \left( i_0; \,S_0  ( \delta_0 , \e'_0, n , K , D)\right) \quad \text{and} \quad \Vol (g) \le \Vol (B_g( x , D))
\le \Vol_{\rm can}  \mathbb B_{K} (D) \, .
\end{equation}
As $\Ric_g \ge - (n-1) K^2 \cdot g$ and as $D^k \Ric_g = 0$ for every $k \in \N^*$, the hypotheses of the main theorem of \cite{HH} are all satisfied 
and it implies that, for every $k \in \N$, from every sequence $(g_i)_{i \in \N}$ of metrics belonging to ${\rm Einst}_M (K,D, i_0)$, one can extract a
subsequence $(g_j)_{j}$ and a sequence of $C^{k+1}$-diffeomorphisms $\varphi _j : M \to M$ such that the sequence of metrics $\varphi ^*_j g_j$ converges (in the $C^{k}$ topology) to some $C^{k}$-metric $g_\infty$. For the sake of simplicity, let us denote by $h_j$ the metric $\varphi ^*_j g_j$. As $k$ may be chosen $\ge 2$ (and as the Ricci curvature may be written in 
terms of the first and second derivatives of the coordinates of the metric), the equalities $\Ric_{h_j} = \frac{1}{n} \,{\rm scal} (h_j) \cdot h_j$ 
give, at the limit, $\Ric_{g_\infty} = \frac{1}{n} \,{\rm scal} (g_\infty)\cdot g_\infty$ and $g_\infty$ is Einstein, thus
it is $C^{\infty}$. Going to the limit, we also get $\,{\rm scal} (g_\infty) \ge - n (n-1) K^2$ and $\diam (g_\infty) \le D$. In order to end the proof, it is thus sufficient
to show that, for every $\tilde x \in \widetilde M$, $\inj_{\tilde g_\infty} (\tilde x) \ge \limsup_{j\to +\infty} \inj_{\tilde h_j} (\tilde x)$, where $(\widetilde M , \tilde g_\infty)$ and $(\widetilde M , \tilde h_j)$ are the riemannian covers of $(M,  g_\infty)$ and $(M , h_j)$ respectively. This is proved by the following Lemma, whose consequence is that every sequence of Riemannian structures corresponding to metrics in ${\rm Einst}_M (K,D, i_0) $ admits a converging subsequence (in the $C^k$-topology) whose limit is a Riemannian structure corresponding to some $g_\infty \in {\rm Einst}_M (K,D, i_0) $.
\end{proof}

\begin{lemma}\label{lemma:limitinj}
For any $k\ge 0$, on any smooth complete Riemannian manifold $(X,g)$, for any sequence of $C^\infty$ complete metrics $(g_i)_{i \in \N}$ which 
$C^k$-converges to $g$, one has $\inj_{g} (x) \ge \limsup_{i\to +\infty} \inj_{g_i} (x)$ for every $x \in X$.
\end{lemma}

\begin{proof}[Proof of Lemma \ref{lemma:limitinj}]
Arguing by contradiction, suppose that $\inj_{g} (x) < \limsup_{i\to +\infty} \inj_{g_i} (x)$ for some $x \in X$, there then exists a subsequence
$(g_j)_{j}$ and a sequence $(\e_j)_j$ of real numbers going to zero such that $\lim_{j\to +\infty} \inj_{g_j} (x) = \limsup_{i\to +\infty} \inj_{g_i} (x)$, thus
such that $\inj_{g} (x) < \lim_{j\to +\infty} \inj_{g_j} (x) = \inj_{g_j} (x) + \e_j$. Therefore, there exits some $g$-geodesic $c : [0, +\infty [ \to X$ and some 
real number $t_0 $ 
($\inj_{g} (x) \le t_0<  \lim_{j\to +\infty} \inj_{g_j} (x)$) such that $c$ is minimizing between $c(0)$ and $c(t_0)$ and is not minimizing between $x=c(0)$ and 
$c(t)$ (for every $t > t_0$). Let us denote by $c_j$ the $g_j$-minimizing geodesic between $c(0)$ and $c(t_0)$, $c_j$ still minimizes between $c_j(0) = c(0)$ 
and $c_j (t_j)$, where $t_j = \Min \big(\inj_{g_j} (x) , t_0 + 1) $, we thus have $t_0 < t_j $ and thus $t'_j :=  d_{g_j}\big(x, c (t_0) \big) < t_j$ when $j$ is large enough. By compactness of the closed ball 
of radius $t_0 + 1$, there exists a subsequence of the sequence $\big(c_j (t_j)\big)_{j}$ (still denoted by $\big(c_j (t_j)\big)_{j}$ for the sake of simplicity) 
such that $c_j (t_j)$ converges (with respect to the Riemannian distance $d_g$) to some point $y \in X$; as $ c_j (t'_j) = c(t_0)$, this yields
$$ d_{g_j}\big(x, c (t_0) \big) + d_{g_j}\big(c (t_0) , c_j(t_j) \big) = t'_j + |t_j - t'_j| = t_j = d_{g_j}\big(x, c_j(t_j) \big)\ \ \ \ \text{(for $t'_j 
< t_j$)}\, ;$$
taking the limits of both sides when $j \to +\infty$ yields
$$d_{g}\big(x, c (t_0) \big) + d_{g}\big(c (t_0) , y \big) =  d_{g}\big(x, y \big) = \lim_{j\to +\infty} t_j = 
\Min \big(\lim_{j\to +\infty} \inj_{g_j} (x) , t_0 + 1 \big) > t_0 \, .$$
It follows that $ d_{g}\big(c (t_0) , y \big) > 0$ and that, if $\g$ is a minimizing geodesic from $c(t_0)$ to $y$, the path obtained by concatenation of $c$ and $\g$
is minimizing, and thus that it is a minimizing geodesic which coincides with the geodesic $c$.  Consequently, there exists $t > t_0$ such that $c$ minimizes between 
$c(0)$ and $c(t)$, in contradiction with the hypotheses. It is thus impossible that $\inj_{g} (x) < \limsup_{i\to +\infty} \inj_{g_i} (x)$ for some $x \in X$,
this ends the proof.
\end{proof}

\begin{proof}[Proof of Corollary \ref{Einstein3}]
By Theorem \ref{Einstein2}, for every $k \in \N$, from any sequence $(g_i)_{i \in \N}$ of elements of ${\rm Einst}_M (K,D, i_0) $, one can extract a 
subsequence $(g_j)_{j \in J_k \subset \N}$ and find $C^{\infty}$-diffeomorphisms $\varphi_j^k$ such that $ \big((\varphi_j^k)^* g_j \big)_{j \in J_k}$ 
$C^k$-converges (when $j \to +\infty$ in $J_k$) to some metric $g^k$ which also belongs to ${\rm Einst}_M (K,D, i_0) $. Constructing the infinite subset
$J_{k+1}$ by extraction from the infinite subset $J_k$, we obtain that $J_p \subset J_k$ for every $p \ge k$. This implies that, for the Lipschitz distance
$d_L$ between the corresponding Riemannian structures (modulo isometries), one has: $d_L (g_j , g^p) = d_L ((\varphi_j^p)^* g_j  , g^p) \f 0$ and $d_L (g_j , g^k) = d_L ((\varphi_j^k)^* g_j  , g^k) \f 0$ when $j$ goes to $+\infty$ 
in $J_p \subset J_k$; a consequence is that $d_L (g^k , g^p) = 0$, thus that $g^p$ is isometric to $g^k$. There therefore exists a fixed Riemannian metric $g$
and a sequence of diffeomorphisms $(\psi_k)_{k \in \N}$ such that $g^k = (\psi_k)^*g$ for every $k \in \N$, hence $g$ also belongs to 
${\rm Einst}_M (K,D, i_0) $.\\
As, for every smooth diffeomorphism $\varphi$, the map $h \mapsto \varphi^* h$, from the set of symmetric 2-tensors into itself, is continuous with respect to 
the $C^k (g_0)$-norm (see before Definition \ref{Eisteinconverge}), by defining $\phi^k_j := \varphi_j^k  \circ \psi_k^{-1}$, we obtain that (for every $k \in \N$)
$\big( (\phi_j^k )^* g_j \big)_{j \in J_k}= \left(\big(\psi_k^{-1}\big)^* \big((\varphi_j^k )^* g_j \big)\right)_{j \in J_k}$ $C^k$-converges to 
$\big(\psi_k^{-1}\big)^* g^k = g$ when $j$ goes to $+\infty$ in $J_k$. Let us now define, for every $p \in  \N$, 
$$J'_p = \{ j \in J_p : \Vert g -  (\phi_j^p )^* g_j  \Vert_{C^{p}(g_0)}  \le 2^{- p} \} \, ;$$
by definition of a $C^{p}(g_0)$-converging sequence, $ J_p \setminus J'_p$ is finite. Choose one integer $j_p$ in each subset $J'_p$ such that $j_p \ge p$, 
this yields, for every $k \in \N$ and every $p \ge k$,
$$\Vert g -  (\phi_{j_p}^p )^* g_{j_p} \Vert_{C^{k}(g_0)}  \le \Vert g -  (\phi_{j_p}^p )^* g_{j_p} \Vert_{C^{p}(g_0)}  \le 2^{- p} \, .$$
Hence $ (g_{j_p} )_{p \in \N}$ is a subsequence of the initial sequence of metrics which verifies the following property: there exists a 
sequence of diffeomorphisms $\big( \phi_{j_p}^p \big)_{p \in \N} $ such that $\big((\phi_{j_p}^p )^* g_{j_p}\big)_{p \in \N} $ $C^k$-converges to $g$
for every $k \in \N$. By definition, this means that $\big((\phi_{j_p}^p )^* g_{j_p}\big)_{p \in \N} $ $C^\infty$-converges to $g$.
\end{proof}

A list of open questions were asked by A.L.~Besse (\cite{Be}, pp.~354--355). We address some of them below.

\begin{itemize}
   \item[1)] Is it possible, for some manifolds $M^n$, to show and compute a number $\varepsilon = \varepsilon (M) >0$  such that the functional $S$ defined above has no critical values in the interval 
   $(-\varepsilon (M),\varepsilon (M))$~? (\cite{Be}, question 12.63).
   \item[2)] Is it possible, for some manifolds $M^n$, to show that the set of critical values of $S$ ({\sl i.e} the values $S(g)$ for $g$ an Einstein metric) is a closed and discrete set (for example in negative curvature) ?(\cite{Be}, question 12.63).
    \item[3)] For which manifolds is this set finite ? (\cite{Be}, question 12.62).
\end{itemize}

Answers to question 1) have already been given. In particular M.~Gromov in \cite{Gr2} shows that, if the simplicial volume of $M^n$ is non-zero (see \cite{Gr2} for the definition), then there exists a constant $C_n$ such that any 
Einstein metric satisfies:
$$S(g) \le -C_n \hbox{~~(SimplicialVolume ($M$))}^{2/n}.$$

We precise this inequality and extend it to some manifolds whose simplicial volume vanishes. We also answer positively to Question 2) (in any dimension and in negative curvature) and to Question 3) (in dimension 4). We furthermore show, in negative curvature, that to each critical value of $S$, corresponds a finite number of critical points.

\begin{corollary}\label{cor:S(Einstein)}
Let $M$ be a closed differentiable manifold of dimension $n\ge 4$, then:
  \begin{itemize}
    \item[(i)] the image by $S$ of the set of negatively curved Einstein metrics is a closed and discrete subset of 
$(- \infty ,\, 0\,)$.
    \item[(ii)] For all $K\in \R$, the set of negatively curved Einstein structures $[g]$ on $M$ which satisfy $S([g]) \ge -K^2$, is 
finite.
  \end{itemize}
\end{corollary}

\begin{proof}
Up to rescaling we may consider the set ${\cal E}_K$ of Einstein metrics  $g$ on $M$ which satisfy  $\Vol(g) =1$, $\scal (g) \ge -K^2$ and $\sigma (g) <0$, where
$\sigma (g)$ is the sectional curvature of $g$. If ${\cal E}_K \ne \emptyset$, let us fix a metric $g_0 \in {\cal E}_K$.

The fundamental group  $\Gamma$ of $M$ acts co-compactly on the simply connected space $(\widetilde M , \tilde g_0)$ where $\tilde g_0$ has pinched negative curvature (and hence is Gromov-hyperbolic). 
This action is furthermore without fixed point, hence $\Gamma$ is torsion-free by Lemma \ref{fixedpoint}, thus all the elements $\g \in \Gamma^*$ are hyperbolic isometries
and verify $\inf_{\tilde x \in \widetilde M} d_{\tilde g_0} (\tilde x , \g \tilde x) = \ell (\g)$ by Lemma \ref{invariantgeod}; this implies that, for every 
$\g \in \Gamma^*$, $ \ell (\g) \ge \Sys_\Gamma (\widetilde M , \tilde g_0) > 0 $, and hence that $\Gamma \in \text{\rm Hyp}_{\rm thick}$. 
On the other hand, every metric $g \in {\cal E}_K$  verifies $-K^2 \leq \scal (g) < \sigma _g < 0$, the main theorem of \cite{Gr5} gives an upper bound of the diameter of  $(M,g)$ in terms of  $\Vol(g)$, hence by a constant 
$C(n)  K^{nr-1}$ (where $r = 1$ if $n \ge 8$ and $r= 3$ if $4 \le n \le 7 $). The injectivity radius of the universal cover $(\widetilde M , \tilde g)$
of $(M,g)$ being infinite, we deduce that ${\cal E}_K \subset {\rm Einst}_M (K' ,D, +\infty) $, where $K' := \dfrac{K}{\sqrt{n (n-1)}}$ and 
$D := C(n)  K^{nr-1}$. We can then apply Corollary \ref{Einstein3}: from every sequence $(g_i)_{i \in \N}$ of elements of ${\cal E}_K$ we can 
extract a subsequence $(g_j)_{j \in J}$ and find diffeomorphisms $\phi_j$ such that $ \big(\phi_j^* g_j\big)_{j \in J} $ $C^\infty$-converges to some
Einstein metric $g^*$ which satisfies $\Vol (g^*) =1$, $\scal (g^*) \ge -K^2$, $\sigma (g^*) \le 0$ and $\diam (g^*) \le C(n)  K^{nr-1}$.

Let us denote by ${\cal E}$ the set of all Einstein metrics on $M$ of volume equal to $1$. A theorem of N.~Koiso (see, for example, \cite{Be}, p.~351, theorem 12.49) shows the existence of $C^\infty$-neighbourhoods of 
$g^*$ in the space of Riemannian metrics, denoted by $U$ and $U'$, of a real analytic submanifold $Z$ of finite dimension in the same space, such that any element of ${\cal E} \cap  U'$ is isometric to an element of 
${\cal E} \cap Z \cap U$, and such that ${\cal E} \cap Z \cap U$ is a real analytic and path-connected subset of $Z\cap U$ (see \cite{Be}, p.~352, Corollary 12.52). It follows that there exists a sequence of diffeomorphisms $\varphi _j : M \to M$ such that, 
for any $j \ge j_0$, $\varphi ^*_jg_j \in {\cal E} \cap Z \cap U$, hence $\varphi ^*_jg_j $
is connected to $g^*$ by a path of Einstein metrics $t \mapsto g_t$ such that $g_0 = \varphi ^*_{j} g_{j}$ and $g_1 = g^*$. Let $I = 
\{t : g_t \text{ is isometric to } g_0\}$, then $I$ is open because, for every $t_0 \in I$, as $g_{t_0}$ is isometric to $g_0$, $g_{t_0}$ is Einstein with negative
sectional curvature and (see \cite{Be}, p.~357, corollary 12.73) the Einstein structure of $g_{t_0}$ is isolated in ${\cal E}$, thus 
$g_t$ is isometric to $g_{t_0}$ (and to $g_0$) for every $t$ in a neighbourhood of $t_0$, hence $t \in I$ for every $t$ in a neighbourhood of $t_0$. 
The set $I$ is also closed for, if $(t_p)_{p \in \N}$ is any 
sequence of elements of $I$ which converges to $t \in [0 , 1]$ then the Lipschitz distance $d_L$ between Riemannian structures verifies
$$d_L \big(\lceil g_{t_p}\rceil , \lceil g_0 \rceil \big) = 0 \ \ \ \text{ and }\ \ \ d_L \big(\lceil g_{t_p}\rceil   , \lceil g_{t}\rceil  \big) \f  0  \ \ \ \text{ when }
\ \ \ t_p \to t \, ,$$
where $\lceil g \rceil $ is the Riemannian structure corresponding to the metric $g$. An immediate consequence is that $d_L \big(\lceil g_{t}\rceil  , \lceil g_0 
\rceil  \big) 
= 0$, thus that $g_t$ is isometric to $g_0$, which proves that $t \in I$. It follows that $I = [0 , 1]$ and that any sequence of Einstein structures of strictly negative sectional curvature and on which the functional $S$ is bounded below has a stationary subsequence, which implies {\it (ii)}. Finally,  {\it (i)} follows immediately from {\it (ii)}.
\end{proof}

In dimension $4$ the Allendoefer-Chern-Weil theory (see, for example, \cite{Be}, Sections 6.31 and 6.34, p.~161) gives the following formulas for the Euler characteristic $\chi (M)$ and the signature $\tau (M)$ of a closed Riemannian manifold $(M, g)$:
$$
\begin{array}{rcl}
   8\pi^2\chi (M)   & = & \displaystyle\int_M (\lVert W^+(g)\rVert^2 + \lVert W^-(g)\rVert^2  - \lVert Z (g)\rVert^2 +  \lVert U(g)\rVert^2 )dv_g\,,\\
   12\pi^2\tau (M) & = & \displaystyle\int_M (\lVert W^+(g)\rVert^2 - \lVert W^-(g)\rVert^2 )dv_g\,.
\end{array}
$$
Where $W^+(g)$, $W^-(g)$, $Z (g)$ and $U (g)$ are the irreducible components of the curvature tensor of $g$ under the action of $SO(4)$ on $T_xM$ for each $x\in M$. For an Einstein metric, we have, by definition,
$$Z(g)=0\quad \textrm{ and }\quad  \lVert U (g)\rVert^2= {1\over 24} \scal (g) ^2\,,$$
from which we deduce that,
$$S(g)\geq -4\sqrt{6}\pi \sqrt{2\chi (M)- 3\lvert\tau (M)\rvert}\,.$$
With this lower bound on $S(g)$ we can apply Corollary \ref{cor:S(Einstein)} and get,

\begin{corollary}\label{cor:Einstein4}
Let $M$ be a closed differentiable manifold of dimension $4$, the set of Einstein structures with negative sectional curvature on $M$ is finite.
\end{corollary}
Concerning the question 1) above, about the possibility of finding an interval around $0$ without critical values of $S$, we get the following proposition.

\begin{prop}\label{prop:S(hyp)}
Let $M$ be a closed differentiable manifold of dimension $n\geq 4$, and let us assume that $M$ admits a non-zero degree continuous map $f : M\to X$, where $X$ is a closed differentiable manifold which carries a metric $g_0$ with negative sectional curvature. Then, for any Einstein metric $g$ on $M$, we have
$$S(g)\leq -n(n-1) \lvert \mathrm{deg}(f)\rvert^{2/n}\lvert \max (\sigma (g_0)) \rvert \Vol (g_0)^{2/n} <0\,.$$
The equality occurs if and only if both metrics have constant sectional curvature and if there exists $\lambda > 0$ such that the map $f$ is homotopic to a Riemannian covering from $(M, \lambda \cdot g)$ onto $(X, g_0)$.
\end{prop}

\begin{proof}
With the same proof than that of  Corollary 1.4 (p.~156-157) from \cite{BCG-Schwarz}, one can easily show that,
$$\Ent (M, g)^n\Vol (g)\geq \lvert \mathrm{deg}(f)\rvert (n-1)^n \lvert \max (\sigma (g_0))\rvert^{n/2}\Vol (g_0)\,,$$
with equality when both metrics have constant curvature and $f$ is homotopic to a locally homothetic covering. We finish by applying Bishop's comparison Theorem which yields:
$$\Ent (M, g)^2\Vol (g)^{2/n}\leq -(n-1)({\scal (g)\over n})\Vol (g)^{2/n}=-\big ({n-1\over n}\big )S(g)\,.$$
\end{proof}
\begin{remarks}\label{rem:S(hyp)}
There is one interesting particular case. When $g_0$ is a locally real hyperbolic metric (with constant curvature equal to $-1$) on $M$, $g$ another Einstein metric on $M$ and $f=id$, the previous inequality then becomes,
$$S(g)\leq S(g_0)< 0\,,$$
and equality if and only if $g$ and $g_0$ define the same Einstein structure.

In fact, in \cite{BCG-GAFA}, Theorem 9.6, p.~774 it is proved that, in dimension $4$, if $(M,g_0)$ is locally real hyperbolic then $g_0$ is the only Einstein metric up to homothety.
\end{remarks}
The next three results go further in the finiteness results and isolation phenomenons.

\begin{prop}\label{prop:Einsteinneg}
For $K>0$ and $n\geq 4$, there is only a finite number of closed manifolds of dimension $n$  (modulo diffeomorphisms) which admit an Einstein metric of negative sectional curvature satisfying, 
$$S(g)\geq -K^2\,.$$
In dimension $4$, among the closed manifolds satisfying $\chi (M)-{3\over 2}\lvert \tau (M)\rvert\leq K$, only a finite number of them admit Einstein metrics of negative curvature.
\end{prop}

\begin{proof}
Let $(M,g)$ be a closed Einstein manifold, of sectional curvature $\sigma (g)<0$ and satisfying $S(g)\geq -K^2$. Up to rescaling we may assume that $\Vol (g)=1$ which implies that $\scal (g)\geq -K^2$. Following the proof of Corollary \ref{cor:S(Einstein)} we deduce that $-K^2\leq \sigma (g)<0$ and that $\diam (g)$ is bounded above by a constant of the form $C(n)  K^{nr-1}$. J.~Cheeger's Comparison Theorem (\cite{Ch}) then shows that,
$$\Vol (g)\leq \int_0^{\diam (g)}l(c)\cosh (Kt) \big ( {\sinh (Kt)\over K} \big )^{n-2}dt\,,$$
for every periodic geodesic $c$, where $l(c)$ is its length. Consequently one has,
$$\inj (g)=\inf_c(l(c))\geq K^{n-1}\Big (\int_0^{C(n)K^{nr}} \cosh (u)(\sinh (u))^{n-2}du  \Big )^{-1}\,.$$
Hence, the Einstein metrics under consideration have universal bounds on their sectional curvature, diameter and injectivity radius. Cheeger's Finiteness Theorem then shows that there are only finitely many possible differentiable manifolds (modulo diffeomorphisms).

In dimension $4$, the proof of Corollary \ref{cor:Einstein4} shows that,
$$S(g)\geq -4\sqrt{6}\pi\sqrt{2\chi (M)-3\lvert \tau (M)\rvert }\geq -8 \sqrt{3}\pi\sqrt{K}\,,$$
if $\chi (M)-{3\over 2}\lvert \tau (M)\rvert\leq K$; from the preceeding argument the set of differentiable manifolds admitting such an Einstein metric is finite.
\end{proof}

For a Riemannian manifold $(M,g)$ let us now define $\textrm{sgl}(\widetilde M, \tilde g)$ to be the length of the smallest geodesic loop in the universal cover $(\widetilde M, \tilde g)$. If there is no such loop, in particular if the curvature of $(M,g)$ is non positive or more generally if 
$(\widetilde M, \tilde g)$ is a Busemann space, we have $\textrm{sgl}(\widetilde M, \tilde g)=+\infty$ by Remark \ref{unicite}.

\begin{prop}\label{prop:gap}
Let $\delta_0 \geq 0$ and $\e'_0 > 0$, if $M$ is a closed differentiable manifold of dimension $n\geq 4$ and whose fundamental group is torsion-free  and belongs to $\text{\rm Hyp}_{\rm thick} (\delta_0, \e'_0)$ then, for any Einstein metric $g$ on $M$, one has
$$S(g)\leq -C(n) \min\{\mathrm{sgl}(\widetilde M, \tilde g)^2\lvert\scal (g)\rvert ;\, r_0(\delta_0, \e'_0)^2 \}\,,$$
where $C(n)$ is a constant depending on $n$ only.
In particular, if $g$ is non positively curved or more generally if 
$(\widetilde M, \tilde g)$ is a Busemann space, we get
$$S(g)\leq -C(n) r_0(\delta_0, \e'_0)^2 \,.$$
\end{prop}
Notice that the upper bound of $S(g)$ given by this last inequality only depends on the fundamental group modulo isomorphisms.
\begin{proof}
For the sake of simplicity we call $\Gamma$ the fundamental group of $M$ and we denote by $H$ the entropy of $(M,g)$. From Theorem \ref{transyst} and the fact that $\Gamma$ is torsion-free, there exists a point $\tilde m\in \widetilde M$ such that $ \sys_\Gamma (\tilde m) \ge \frac{r_0}{H} $. If $m$ is the projection of $\tilde m$ on $M$ by the covering map this implies that the ball $B_{\widetilde M}(\tilde m, r)\subset \widetilde M$ projects isometrically onto the ball $B_M(m, r)\subset M$ for any radius $r< {r_0\over 2 H}$.

On the other hand, the main theorem of \cite{Sa} asserts that, for $r<\mathrm{sgl}(\widetilde M, \tilde g)/2$, we have 
$$\Vol (B_{\widetilde M}(\tilde m, r))\geq C_1(n) r^n\,,$$
for some number $C_1(n)$ depending on $n$ only. Then, if $R={1\over 2}\min\{\mathrm{sgl}(\widetilde M, \tilde g), {r_0\over H}\}$ we get,
$$\Vol (M,g)^{1/n}\geq \Vol (B_M(m, R))^{1/n}=\Vol (B_{\widetilde M }(\tilde m, R))^{1/n}\geq {1\over 2}C_1(n)^{1/n}\min\Big\{\mathrm{sgl}(\widetilde M, \tilde g), {r_0\over H}\Big\}\,,$$
and Bishop's Comparison Theorem gives,
$$H^2=\Ent (M,g)^2\leq -{n-1\over n}\scal (g)\,.$$ 
Combining these two inequalities and noticing that ${\lvert \scal (g)\rvert\over H^2}>1$ yield the desired results for $C(n)={1\over 4}C_1(n)^{2/n}$.
\end{proof}

\begin{remark}
We wish to  emphasise that among the manifolds which verify the hypotheses of Proposition \ref{prop:gap} there are some with zero simplicial volumes. Indeed, examples are given by products of a sphere and a  negatively curved closed manifold and connected sums of these products with a simply connected manifold. Proposition \ref{prop:gap} thus improves M.~Gromov's result mentioned above.
\end{remark}

The last result of this section is the

\begin{prop}\label{hypisole}
Let $V$ be a positive real number and $n\geq 4$ be an integer. There exists $\varepsilon = \varepsilon (n, V) >0$ such that every closed Einstein manifold of dimension $n$ and of negative curvature, $(M^n,g)$, which admits a continuous map of non-zero degree onto a closed locally real hyperbolic manifold $(X^n,g_0)$  (\sl{i.e.}, $\sigma (g_0)\equiv -1$) with $\Vol (X^n,g_0)\leq V$, and whose scalar curvature satisfies $S(g)> S(g_0)-\varepsilon$, is homothetic to $(X^n,g_0)$. 
\end{prop}
We recall that, by homothetic, we mean that there exists $\lambda >0$ such that $(M,\lambda g)$ is isometric to $(X,g_0)$. This result applies, in particular, when $M=X$ and the map is the identity and as such can be viewed as a gap theorem.

\begin{proof}
From Proposition \ref{prop:Einsteinneg} there is only finitely many closed manifolds, $M_1, \dots, M_p$, such that the set,
$$\Sigma_{M_i}=\Big\{ S(g) : g \mathrm{\hbox{ Einstein metric on }} M_i \mathrm{\hbox{ such that }} \sigma (g)<0 \mathrm{\hbox{ and }} S(g)\geq -n(n-1)(V^{2/n}+1) \Big\}\,,$$
is non empty. The hyperbolic manifolds that we are considering, that is with $\Vol (X^n,g_0)\leq V$, are members of the set $\{ M_1, \dots, M_p\}$. Let us define $\Sigma=\bigcup_{i=1}^p\Sigma_{M_i}$; from Corollary \ref{cor:S(Einstein)}, $\Sigma$ is a finite set whose elements will be denoted by $s_1,\dots ,s_N$ and we set,
$$\varepsilon = \min\Big\{ \min_{1\leq i\ne j\leq N}(\lvert s_i-s_j\rvert),\, n(n-1)\Big\}\,.$$
The number $\varepsilon$ only depends on $n$ and $V$. Let $(X^n,g_0)$ be a locally real hyperbolic manifold of volume less than or equal to $V$, let $M^n$ admit a non-zero degree map $f$ onto $X^n$ and $g$ be an Einstein metric on $M^n$ with negative curvature and satisfying $S(g)> S(g_0)-\varepsilon$. Remark \ref{rem:S(hyp)} implies that $S(g_0)-\varepsilon<S(g)\leq S(g_0)$, hence that $\lvert S(g)-S(g_0)\rvert<\varepsilon$. As $S(g)$ and $S(g_0)$ belong to $\Sigma$, the choice of $\varepsilon$ shows that $S(g)=S(g_0)$. The inequality stated in Proposition \ref{prop:S(hyp)} proves that $\lvert \mathrm{deg}(f)\rvert = 1$; from this and from the equality case of Proposition \ref{prop:S(hyp)} (see also
Remark \ref{rem:S(hyp)}) we deduce that $f$ is homotopic to a homothety.
\end{proof}

\section{Appendix : Basic results on Gromov-hyperbolic spaces}\label{outils}

\subsection{Why is there quantitative gaps between basic sources, and how to precise coherently the constants:}\label{coherence}

\small
These quantitative gaps are mainly due to the fact that each of the basic references (for example \cite{GH}, \cite{CDP}, \cite{BH}) starts from a different 
definition of Gromov $\delta$-hyperbolic spaces and use different methods in order to prove the basic results. A consequence is that, depending on the 
source, all these definitions and results (though qualitatively equivalent) may differ by a multiplicative factor: these differences are not worries when one is
working modulo quasi-isometries, but they become a real problem when one wants to pick results from different sources and to chain them together and when
the final value is important, by example when one aims to determine a precise explicit universal bound for some invariant (especially if one wants this bound to be as independent as possible of the $\delta$-hyperbolic space under consideration).
We thus propose in this subsection a journey in order to raise these ambiguities and to control all these constants in a coherent way. As these results 
are (modulo multiplicative constants) classical, the proofs of the most classical ones will be left to the reader.

\normalsize

\smallskip
\emph{Given any three nonnegative numbers $\alpha , \beta, \gamma$, we define the tripod $T := T(\alpha , \beta, \gamma)$ as the metric simplicial tree 
with $3$ vertices
$x', \,y' ,\, z'$ of valence $1$ (the \lq \lq endpoints"), one vertex $c$ of valence $3$ (the \lq \lq branching point"), and $3$ edges  $[cx'] , \, [cy'] , \, [cz']$ 
of respective lengths $\alpha , \beta, \gamma$ (the \lq \lq branches"). We denote by $d_T (u,v) $ the distance on this tree between two points $u,v \in 
T$, i. e. the minimal length of a path contained in $T$ and joining $u$ to $v$).}

\smallskip
For sake of simplicity, we only consider geodesic metric spaces (see Definition in section \ref{notations}). In such a space a geodesic triangle 
$\Delta = [x , y , z]$ is the union of three geodesics $[x , y ]$, $[y , z ]$ and $[z , x ]$. Given any three points $x,y,z$ in any geodesic metric space, there
exists at least one geodesic triangle $\Delta = [x , y , z]$ whose sides have respective lengths $d(x,y),\, d(y,z)$ and $d(x,z)$. 

\begin{lemma}\label{prodist} 
To any geodesic triangle $\Delta$ corresponds a metric tripod $(T_\Delta , d_T)$ and a surjective map $f_\Delta : \Delta \f  T_\Delta$ (called the 
\lq \lq approximation of $\Delta$ by a tripod") such that, in restriction to each side of $\Delta$, $f_\Delta$ is an isometry,
\end{lemma}

Indeed, $T_\Delta$ is constructed as the tripod $ T(\alpha , \beta, \gamma)$, where (by the triangle inequality) $(\alpha , \beta, \gamma)$ is the 
unique element of $[0 , +\infty[^3$ such that $d(x, y) = \alpha + \beta$, $d(x, z) = \alpha + \gamma$ and $d(y, z) = \beta + \gamma$. This choice of 
$(\alpha , \beta, \gamma)$ implies the existence of the map $f_\Delta : \Delta \f  T_\Delta$ as asserted in Lemma \ref{prodist}.

\begin{defis}\label{hypdefinition}
A geodesic triangle $\Delta$ of $(X,d)$ is said to be $\delta$-thin if, for every $u \in T$ and every $x,  y \in f_\Delta^{-1} (\{u\})$, one has $d(x,y) \le 
\delta$.\\
In the whole of this paper, a metric space is said to be $\delta$-hyperbolic if it is geodesic, proper, and if all its geodesic triangles are $\delta$-thin.
\end{defis}

The following results are well known (and often taken as definitions of $\delta$-hyperbolicity); their proof may be found (with variable constants, depending 
on the choice of the definition) in any classical source (see for example \cite{GH}, \cite{CDP}, \cite{BH}).

\begin{lemma}\label{proprietes} For every $\delta$-hyperbolic space $(X,d)$, one has:

\begin{itemize}

\item[(i)] For every geodesic triangle $\Delta $, its approximation $f_\Delta : \Delta \to (T_\Delta, d_T)$ by a tripod verifies
 $$d(u,v) - \delta  \le d_T \big(f_\Delta (u) , f_\Delta (v) \big) \le d(u,v) \ .$$



\item[(ii)] (Quadrangle Lemma)  For every four points $  x, y , z, w \in X$, one has
$$d(x,z) + d(y,w) \le \Max \big( d(x,y) + d(z,w) \, ;\, d(x,w) + d(y,z) \big) + 2\,\delta \ , ;$$
if moreover if $ [x , y ]$, $ [y , z]$, $ [z , w ]$ and $ [x , w]$ are geodesic segments between the vertices of the quadangle $[x , y , z , w ]$ then, for every 
$v \in [x , w ]$, one has
$$\  d(v\, ,\, [x , y ]\cup [y , z ]\cup [z , w ]) \le 2\, \delta\ .$$
\end{itemize}
\end{lemma}

We often use the following sharper version of the above Quadrangle Lemma:

\begin{lemma}\label{rectangletriangle} 
On any geodesic triangle $\Delta = [x , y , z]$, for every point $u \in [y , z]$, one has
$$ d(x, u) + d(y , z) \le \Max \big( d(x,y) + d(u,z) \, ;\, d(x,z) + d(y , u) \big) + \delta$$
\end{lemma}

\begin{proof} Using the approximation $f_\Delta : \Delta \to (T_\Delta, d_T)$ by the associated tripod, as $f_\Delta (u)$ is on the union of the branches 
$[c , f_\Delta (y)]$ and $[c , f_\Delta (z)]$ of this tripod (where $c$ is the branching point of $T$), we easily verify that
$$d_T(f_\Delta (x) , f_\Delta (u)) + d_T(f_\Delta (y) , f_\Delta (z)) = $$
$$ = \Max \big( d_T(f_\Delta (x),f_\Delta (y)) + d_T( f_\Delta (u),f_\Delta (z)) \, ;\, d_T(f_\Delta (x),f_\Delta (z)) + d_T(f_\Delta (y) ,  f_\Delta (u)) \big)  .$$
We now use Lemma \ref{prodist}, which proves that  $d_T(f_\Delta (y) , f_\Delta (z)) $, $d_T(f_\Delta (x),f_\Delta (y))$, $d_T(f_\Delta (x) ,f_\Delta (z))$, $d_T( f_\Delta (u) ,f_\Delta (z))$ and $d_T(f_\Delta (y) ,  f_\Delta (u))$ are respectively equal to
$d(y , z) $, $d(x,y)$, $d(x,z)$, $d(u,z)$ and $d(y , u)$ 
and Lemma \ref{proprietes} (i), which proves that $d(x,u) \le d_T \big(f_\Delta (x) , f_\Delta (u) \big) + \delta $ and ends the proof.
\end{proof}

\subsection{Projections and quasi equality in the triangle inequality}\label{sectionprojections}

\begin{defi}\label{projete}
For every closed subset $F$ of any metric space $(X,d)$ and every point $x \in X$, a \lq \lq projection of $x$ on $F$" is any point $\bar x \in F$ such that
$ d(x,\bar x) = d(x,F) := \inf_{z \in F} d(x,z)$.
\end{defi}
When it exists, a projection of $x$ on $F$ is generally not unique.
By continuity of the distance there always exist a projection of $x$ on $F$ when $F$ is compact, and when $F$ is closed if the metric space is proper; 
however the assumption \lq \lq proper space" is not necessary when $F$ is the image of a geodesic, as proved by the following
\begin{lemma}\label{projectiongeod}
If $c$ is a geodesic line (or radius or segment) in a metric space $(X,d)$, every point $x \in X$ admits a projection on the image ${\rm Im}(c)$ of $c$
and the map $ x \mapsto d(x , {\rm Im}(c) )$ is Lipschitz with Lipschitz constant $1$.
\end{lemma}
The proof, elementary, is left to the reader.
\begin{lemma}\label{minorangle} In a $\delta$-hyperbolic space (thus geodesic and proper) $(X,d)$, for any three points $x , y , z$ and any geodesic
segment $[y , z ]$ joining $y$ to $z$, if $ d(x,y) \le d \big( x , [y,z]\big) + \eta $, then $ d(x,z) \ge d(x, y) + d(y, z) - 2\, ( \eta + \delta) $.
\end{lemma}

\begin{proof}
Let $\Delta$ be any geodesic triangle with vertices $x , y , z$ whose third side $[y , z ]$ is the geodesic segment under consideration, let 
$(T_\Delta , d_T)$ be the tripod associated to $\Delta$ and $f_\Delta : \Delta \to T_\Delta$ the approximation of $\Delta$ by this tripod. If 
$x',y',z'$ are the three endpoints $f_\Delta (x) , f_\Delta (y) , f_\Delta (z)$ of the tripod, the branching point $c$ satisfies
$ f_\Delta^{-1} (\{c\}) = \{c_x , c_y, c_z\}$, where $c_x$, $c_y$ and $c_z$ respectively lie on $[y,z]$, $[x,z]$ and $[x,y]$).
One thus gets:
$$d_T (x',y') = d(x, y) \le d \big( x , [y,z]\big) + \eta \le d(x, c_x) + \eta \le d_T ( f_\Delta (x) , f_\Delta (c_x))+ \delta + \eta = d_T ( x', c)+ \delta + \eta\ ,$$
where the first, second, third, fourth and fifth equalities or inequalities respectively follow from Lemma \ref{prodist}, from the hypothesis, from the fact that $c_x \in [y,z]$, from Lemma \ref{proprietes} (i) and from the fact that $f_\Delta (c_x) = c$. From this we deduce that 
$ d_T ( y', c) = d_T (x',y') -  d_T ( x', c) \le \delta + \eta$; a consequence of this last inequality, of Lemma \ref{prodist} and of the definition of the distance
on the tripod is that
$$ d(x,z) = d_T (x',z') = d_T ( x', y') + d_T ( y', z') - 2 \, d_T (y', c) \ge d(x, y) + d ( y , z ) - 2\, (\delta + \eta )\ .$$
\end{proof}
Replacing $\eta$ with $0$ in Lemma \ref{minorangle}, we get the following immediate corollary, where the geodesic $c$ under consideration is either a geodesic
line, or a geodesic ray, or a geodesic segment.

\begin{lemma}\label{projection} \emph{(quasi equality in the triangle inequality)}
In any $\delta$-hyperbolic space $(X,d)$, for any geodesic $c$ of $(X,d)$ and every point $y$ on this geodesic, for every point $x \in X$, any of its projections
$\bar x$ on the geodesic verifies $ d(x,y) \ge d(x,\bar x) + d(\bar x , y) - 2\, \delta$.
\end{lemma}

\begin{lemma}\label{ecartement}
In any $\delta$-hyperbolic space $(X,d)$, for any four points $x,\, y, \,\bar x ,\,\bar y $ such that $\bar x $ and $\bar y $ are projections of $x$ and $y$ 
(respectively) on some geodesic segment $ [\bar x , \bar y]$ between $\bar x$ and $\bar y$, then
$$ d(\bar x , \bar y) > 3\, \delta  \implies  d(x,y) \ge d(x,\bar x) + d(\bar x , \bar y) + d( \bar y , y) - 6\, \delta \ .$$
\end{lemma}

\begin{proof}
The lemma \ref{projection} gives $ d (x, \bar y) \ge d(x,\bar x) +  d(\bar x , \bar y) - 2\, \delta$ and $ d (y, \bar x) \ge d(y,\bar y) +  d(\bar y , \bar x) - 
2\, \delta$. From this and from the hypothesis $d(\bar x , \bar y) > 3\, \delta $, we get:
$$d(x , \bar y) + d(y , \bar x) - 2\, \delta \ge  d(x,\bar x) + 2\, d(\bar x , \bar y) +  d(y,\bar y)  - 6\, \delta > d(x,\bar x)  +  d(y,\bar y) \ ,$$
and, using this inequality and the quadrangle Lemma \ref{proprietes} (ii), we obtain
$$ d(x,\bar x) + 2\, d(\bar x , \bar y) +  d(y,\bar y)  - 6\, \delta \le d(x , \bar y) + d(y , \bar x) - 2\, \delta \le \Max \left( d(x,\bar x) + d(y, \bar y) \, ;
\, d(x,y) + d(\bar x , \bar y) \right) $$
$$ \le d(x,y) + d(\bar x , \bar y) \ ,$$
and this concludes.
\end{proof}

 \subsection{Asymptotic and diverging geodesic lines}

\small
In the following Proposition, points (i) and (iii) are classical: indeed, for any pair of geodesic rays $ \g_1, \,\g_2 : [0 , +\infty[ \f X$, the definition of the 
ideal boundary of a Gromov-hyperbolic space $(X,d)$ (see Definition  7.1, p. 117 of \cite{GH} and the Definition given in chapter 2 of \cite{CDP}, these two 
definitions being equivalent by  Proposition 7.4 p. 120 of \cite{GH}) implies the equivalence between the equality $ \gamma_1(+\infty) = \gamma_2(+\infty)$ 
and the finiteness of the Hausdorff distance between the images of $ \g_1 $ and $\g_2 $. In the case where $ \g_1 $ and $\g_2 $ are geodesic lines (as in (i)), 
it is sufficient to consider each $\g_i$ as the union of two geodesic rays. The proofs of properties (i) and (ii) being classical, we shall leave them to the reader
and only precise what are the constants which occur according to our definition of $\delta$-hyperbolic spaces.\\
On the contrary, we shall give a proof of the point (iii) of the same Proposition, though it is similar to a classical result, which says that, given three
geodesic rays $\g,\, \g_1, \,\g_2 : [0 , +\infty[ \f X$ such that $ \g_1 (+ \infty) = \g_2 (+ \infty) =\g (+ \infty)$, the Busemann function $b_\gamma$ 
associated to $\g$ verifies $\sup_{t \in [0 , +\infty[ }| b_\gamma \big(\g_1 (t) \big) - b_\gamma \big(\g_2 (t) \big)| \le C(\delta) $, where the bound 
$C(\delta) $ only depends on the hyperbolicity constant $\delta$ for a good choice of the origins of the rays. In fact, in the point (iii) of next Proposition, we need to precise quantitatively the choice of these origins and the value of $C(\delta) $.

\normalsize

\begin{prop}\label{geodasympt}
In any $\delta$-hyperbolic (thus geodesic and proper) space $(X,d)$ 

\begin{itemize}
\item [(i)] if $\gamma_1, \, \gamma_2 : \R \to X$ are two geodesic lines verifying $ \gamma_1(+\infty) = \gamma_2(+\infty)$ and 
$ \gamma_1(-\infty) = \gamma_2(-\infty)$ then, for every $t \in \R$, there exists $s = s(t) \in \R$ such that $ d \big(\gamma_1 (t) , \gamma_2 (s) \big) 
\le 2\, \delta$.\\
Moreover there exists choices of the origins $t_0$ and $s_0$ of $\gamma_1$ and $\gamma_2$ (respectively) such that $d \big(\g_1 ( t_0 + t) , 
\g_2 ( s_0 + t) \big) \le  4 \, \delta$ for every $t \in \R $.
\end{itemize}
On the other hand, if $\gamma_1, \,\gamma_2 : [ 0 , +\infty [ \to X $ are two geodesic rays satisfying $ \gamma_1(+\infty) = \gamma_2(+\infty)$,
\begin{itemize}
\item [(ii)] for every $t \in  [ d(\gamma_1 (0) , \gamma_2 (0) ) \, , \, +\infty [$, there exists $s \in  [0 , +\infty [ $ such that 
$ d(\gamma_1 (t) , \gamma_2 (s)) \le 2\, \delta$,

\item [(iii)] there exist $t_1 , \, t_2 \ge 0$, verifying $ t_1 + t_2 = d \big( \gamma_1 (0) ,  \gamma_2 (0)\big)$, such that 
$ d \left(\gamma_1 (t_1 + t) , \gamma_2 (t_2 + t) \right) \le 2\, \delta$ for every $ t \in \R^+$.

\end{itemize}
\end{prop}

\begin{proof}[Proof of (iii)]
Using (ii), if $t'$ is great enough (this implies that $t' >  d\big(\gamma_1 (0) , \gamma_2 (0) \big) + 2\, \delta$), there exists $s'$ such that 
$ d\big(\gamma_1 (t') , \gamma_2 (s')\big) \le 2\, \delta$. For sake of simplicity, make $x_0 = \gamma_1 (0)$, $y_0 = \gamma_2 (0)$, 
$x = \gamma_1 (t')$, $y= \gamma_2 (s')$ and consider two geodesic triangles $\Delta = [ x_0 , x , y_0]$ and $\Delta' = [ x , y , y_0]$ (with common 
side $[x ,  y_0]$) and their approximations by tripods $ f_\Delta : (\Delta , d) \to  (T_\Delta , d_T)$ and $ f_{\Delta'} : (\Delta' , d) \to (T_{\Delta'} , d_{T'})$
(as in Lemma \ref{prodist}). Denoting by $c$ (resp. by $c'$) the branching point of $T_{\Delta}$ (resp of $T_{\Delta'}$), we have 
$ f_\Delta^{-1} \big( \{ c\}\big) = \{ c_0 , c_1 , c_2\}$, where $c_0 $, $ c_1 $ and $c_2$ respectively lie on the sides $[ x_0 , y_0]$, $[ x_0 , x]$ and 
$[ y_0 , x]$ of $\Delta$. Define
\begin{equation*} 
t_1 = d(x_0 , c_1) = d(x_0, c_0) = d_{T} \left( f_\Delta (x_0) , c\right) \ \ \ \   , \ \  \ \  t_2 = d(y_0 , c_2) = d(y_0, c_0) =
d_T \left( f_\Delta (y_0) , c\right) \  ,
\end{equation*}
where $ t_1 + t_2 = d(x_0 , y_0 ) $ and where (by the last equations and by Lemma \ref{prodist})
$$t_1 - t_2 = d_{T} \left( f_\Delta (x_0) , c\right) - d_T \left( f_\Delta (y_0) , c\right) = 
d_{T} \left( f_\Delta (x_0) , f_\Delta (x)\right) - d_T \left( f_\Delta (y_0) ,  f_\Delta (x)\right) $$
\begin{equation}\label{busetripode}
= d ( \g_1(t'),x_0) -  d (\g_1(t'), y_0)\ ;
\end{equation}
as $\g_1 (t_1) = c_1$, Lemma \ref{proprietes} (i) implies that $d \left( \g_1 (t_1) , c_2\right) \le \delta$ and, for every $t \in [0, t' - 2\delta -t_1]$, the point 
$\g_1 (t_1 + t)$ relies on $[ c_1 , x]$, thus the point $ f_\Delta \big( \g_1 (t_1 + t) \big)$ is located on the branch $ [ c , f_{\Delta} (x) \,]$ of the tripod 
$T_\Delta$. Denote by $u$ the point of the geodesic segment $[c_2 , x] \subset [y_0 , x]$ such that $  f_\Delta (u) = f_\Delta \big(\g_1 (t_1 + t) \big) $,
thus verifying $ d(c_2, u) = d \left(c_1 , \g_1 (t_1 + t)\right) = t $ (by Lemma \ref{prodist}); hence $u$ satisfies the following three properties:
\begin{equation}\label{triprop}
 d(y_0, u) = t_2 +  t \ \ , \ \  d(x,u) = t' - t_1 -t \ge 2\,\delta \ \ , \ \   d \left(\g_1 (t_1 + t) , u \right) \le \delta \ ,
\end{equation}
where the first equality follows from the properties  $d(c_2 , u) = t$, $ d(y_0 , c_2) = t_2$ and $  c_2 \in [y_0 , u]$ and where the two last inequalities
result from the equality $  f_\Delta (u) = f_\Delta \big(\g_1 (t_1 + t) \big) $, from Lemma \ref{proprietes} (i) (which proves that $d \left(\g_1 (t_1 + t) , u \right) 
\le \delta$) and from the fact that $ d(x,u) = d_T \big( f_{\Delta} (x) ,  f_{\Delta} (u) \big) = d \big(x , \g_1 (t_1 + t) \big) = t' -t_1 -t  $ by Lemma 
\ref{prodist}).

\smallskip
As $u \in [ x ,  y_0]$, by construction of the tripod $(T_{\Delta'} , d_{T'})$ (see before Lemma \ref{prodist}), $f_{\Delta'} (u)$ lies in the union of the two 
branches $ [c' , f_{\Delta'} (x)] $ and  $ [c' , f_{\Delta'} (y_0)]$ of $T_{\Delta'}$. As $ d(x,y) = d\big(\gamma_1 (t') , \gamma_2 (s')\big) \le 2\, \delta$,
Lemma \ref{prodist}, the second of the properties \eqref{triprop} and the definition of $t$ give:
$$ d_{T'} \left( f_{\Delta'} (x) ,  c' \right) \le d_{T'} \left( f_{\Delta'} (x) ,  f_{\Delta'} (y)\right) = d(x,y)
\le  2\, \delta \le  d(x,u) = d_{T'} \left( f_{\Delta'} (x) ,   f_{\Delta'} (u)\right) \ ; $$
thus, as an immediate consequence, $ f_{\Delta'} (u) \notin \, ] c' , f_{\Delta'} (x)]$ and $f_{\Delta'} (u) $ lies on the branch $[ c' , f_{\Delta'} (y_0)]$ 
of $T_{\Delta'}$. Hence, there exists $u' = \g_2 (s) \in [y_0 , y] = \g_2 \big([0 , s']\big)$ such that $   f_{\Delta'} (u') = f_{\Delta'} (u)$, thus (using Lemma
\ref{proprietes} (i)) such that $d \left(u , u' \right) \le \delta$; in addition, Lemma \ref{prodist}, the equality $ f_{\Delta'} (u') = f_{\Delta'} (u)$ and the first of 
properties \eqref{triprop} imply that $ s = d(y_0 , u')  =  d_{T'} \left(  f_{\Delta'} (y_0) , f_{\Delta'} (u) \right) = d( y_0 , u ) =  t_2 +  t  $, thus that 
$ u' =  \g_2 (t_2 + t) $; the triangle inequality and the last of properties \eqref{triprop} then give: $ d \left( \g_1 (t_1 + t) , \g_2 (t_2 + t) \right) \le  
d \left( \g_1 (t_1 + t) , u\right) + d(u,u') \le 2\, \delta$ for every $t \in [0, t' - 2\delta -t_1]$, thus (when $ t' \f +\infty$) for every $ t \in [0 , +\infty [$.
\end{proof}

For every geodesic $c$ of $(X,d)$ and every $x \in X$, call $\pi_c(x)$ the set of all the projections of the point $x$ on the image ${\rm Im}(c)$ of $c$; Lemma 
\ref{projectiongeod} guarantees that, if $c$ is a geodesic segment, or ray, or line, then $\pi_c(x)$ is never empty. The image $\pi_c(E)$ of a subset
$E \subset X$ by the projection on $c$ is defined as $\cup_{x \in E}\, \pi_c(x)$. We have the

\begin{prop}\label{projgeod}
If $c_1,\, c_2$ are any two geodesic lines in any $\delta$-hyperbolic space $(X,d)$ verifying $ \{ c_1  (-\infty) \, ,\, c_1  (+\infty)\} \cap 
\{ c_2  (-\infty) \, ,\, c_2  (+\infty)\} = \emptyset $ (in the ideal boundary), then
\begin{itemize}
\item [(i)] the subset $\Pi$ of $\R$ such that $\g_1 (\Pi) = \pi_{c_1}({\rm Im} (c_2))$ is bounded from below and from above,

\item [(ii)] there exists points $x_0 = c_1 (s_0)$ and $x'_0 = c_1 (s_1)$ on the geodesic line $c_1$ such that, for every $s \in \R$, the gap between the 
two sides of the triangle inequality is bounded, precisely:
$$\left\lbrace \begin{array}{l} 
\limsup_{t \to +\infty} \left[ d\left(c_1(s) , x_0\right) + d \left( x_0 , c_2 (t )\right) - d\left(c_1(s) , c_2 (t)\right)  \right] \le 5 \delta \\
\limsup_{t \to -\infty} \left[ d\left(c_1(s) , x'_0\right) + d \left( x'_0 , c_2 (t )\right) - d\left(c_1(s) , c_2 (t)\right)  \right] \le 5 \delta
\end{array}\right. $$
\end{itemize}
\end{prop}

\begin{proof}[Proof of (i)]
We shall only prove that $\Pi$ admits a bound from above (the proof of the existence of a bound from below is the same, while changing the orientation of 
$c_1$). By contradiction, suppose the existence of a sequence $(s_n)_{n \in \mathbb N}$ in $\Pi$ which goes to $+\infty$, there then exists a real sequence
$(t_n)_{n \in \mathbb N}$ such that $ c_1 (s_n) \in \pi_{c_1}\left(c_2 (t_n)\right)$, and then
\begin{itemize}
\item either $(t_n)_{n \in \mathbb N}$ is bounded, but then $d \left(c_1 (0) , c_2 (t_n)\right)$ is bounded by some constant $ D> 0$ and 
$d\left(c_2 (t_n) , c_1 (s_n)\right) = d\left(c_2 (t_n) , {\rm Im} (c_1)\right) \le d \left(c_1 (0) , c_2 (t_n)\right) \le D $, and this implies (by the triangle 
inequality) that $ s_n = d\left(c_1 (0) , c_1 (s_n)\right) \le 2\,D$;

\item or there exists a subsequence of $(t_n)_{n \in \mathbb N}$ (still denoted by $(t_n)_{n \in \mathbb N}$) which goes to $\pm \infty$. When 
$n \to +\infty$, as the limits (on the ideal boundary $\partial X$) of $ c_1 (s_n) $ and  $ c_2 (t_n) $ do not coincide, Gromov's product\footnote{Gromov's product $\big( x|y \big)$ with respect to a chosen origin $ x_0$ is defined as $\big( x|y \big) := \frac{1}{2}\,\big( d(x_0 , x) + d(x_0 , y) - d(x,y)\big)$. 
The fact that $\big(c_1 (s_n) | c_2 (t_n) \big)$  is bounded when $c_1 (+ \infty) \ne c_2 (+ \infty)$ is a direct consequence of the definition of the 
ideal boundary given (for instance) in \cite{CDP}, section 2.1 p. 16, see Proposition 2.1.2 p. 18 of \cite{CDP}.} 
$ \big(c_1 (s_n) | c_2 (t_n) \big)$ (with respect to the origin $ c_1 (0)$) is bounded by some constant $D'$, and (using also the inequality
$d\left(c_2 (t_n) , c_1 (s_n)\right) = d\left(c_2 (t_n) , {\rm Im} (c_1)\right) \le d\left(c_1 (0) , c_2 (t_n)\right)$) this implies: 
$$ s_n \le  d\left(c_1 (0) , c_1 (s_n)\right) + 
d\left(c_1 (0) , c_2 (t_n)\right) -  d\left(c_2 (t_n) , c_1 (s_n)\right) = 2 \big(c_1 (s_n) | c_2 (t_n) \big)\le 2\, D' \  .$$
\end{itemize}
This contradiction proves the non-existence of a sequence $(s_n)_{n \in \mathbb N}$ in $\Pi$ which goes to $+\infty$.
\end{proof}

\begin{proof}[Proof of (ii)]
Denote by $ K_t $ the closure of $\pi_{c_1} \left(c_2 ([t , + \infty [)\right)$; for every sequence $ \left(t_n\right)_{n \in \mathbb N}$ going to $+\infty$ define 
$\cap_{n\in \mathbb N} \  K_{t_n} = K $ (verify that $K$ does not depend on the choice of the sequence). Property (i) implies that each $K_t$ is bounded 
and thus compact (because $(X,d)$ is proper), it results that $K$ (as a decreasing intersection of the compact sets $K_n$) is a non empty compact set.
Let us now prove that
\begin{equation}\label{lengthK}
K \text{ is included in a segment } [c_1(a), c_1 (b)] \text{ of length} \le 3 \delta \text{ of the geodesic } c_1 \ .
\end{equation}
Indeed, for every $n \in \N$, $K \subset K_n \subset \text{\rm Im} (c_1)$. Furthermore,
every $ x , x' \in K$ are (respectively) limits of sequences $ \left(x_n\right)_{n \in \mathbb N}$ and $ \left(x_n'\right)_{n \in \mathbb N}$ of 
elements of $\pi_{c_1} \left(c_2 ([ n , + \infty [)\right)$: in fact, as $x$ (resp. $x'$) lies in the closure $K_n$ of $\pi_{c_1} \big(c_2 ([n , + \infty [)\big)$,
we only have to choose $x_n$ (resp. $x'_n$) in the (non empty) intersection of $\pi_{c_1} \left(c_2 ([n , + \infty [)\right)$ with 
$ B_X \left( x , \frac{1}{n}\right)$ (resp. with $ B_X \left( x' , \frac{1}{n}\right)$). There thus exist $ t_n , \, t_n' \in [ n , + \infty [$ such that $x_n$ 
(resp. $x_n'$) is a projection on $c_1$ of $c_2 (t_n)$ (resp. of $c_2 (t_n')$). Using two times the lemma \ref{projection}, and afterwards the quadrangle 
Lemma \ref{proprietes} (ii), we get:
$$d\left(x'_n , c_2 (t_n' )\right) + 2 \, d\left(x_n , x_n' \right) + d\left(x_n , c_2 (t_n )\right) \le d\left(x_n , c_2 (t_n' )\right) + d\left(x_n' , c_2 (t_n )\right) +
4 \delta$$
\begin{equation}\label{projgeod1}
\le \Max \left[ d\left(x_n , x_n' \right) + d\left( c_2 (t_n ) , c_2 (t_n' )\right)  \,;\, d\left(x_n , c_2 (t_n )\right) + d\left(x_n' , c_2 (t_n' )\right)\right] +
6 \delta\ \ .
\end{equation}
As $ d\left(x_n , c_1 (0 )\right) $ and $ d\left(x_n' , c_1 (0 )\right) $ are bounded by property (i), then $d\left(x_n , c_2 (0 )\right)$ and 
$ d\left(x_n' , c_2 (0 )\right)$ are bounded by some constant (say $C$); from this and the triangle inequality, when $n$ is great enough, follows that:
$$ d\left(x_n , c_2 (t_n )\right) \ge t_n - d\left(x_n , c_2 (0 )\right)  \ge t_n - C \ \ , \ \  d\left(x_n' , c_2 (t_n' )\right) \ge t_n' - d\left(x_n' , c_2 (0 )\right)  
\ge t_n' - C \ \ ,$$
As $ t_n + t_n' - | t_n - t_n' | = 2\,\Min ( t_n , t_n') \ge 2\,n$, when $n$ is great enough it comes that
$$ d\left(x_n , c_2 (t_n )\right) + d\left(x_n' , c_2 (t_n' )\right) \ge  t_n + t_n' - 2\,C \ge | t_n - t_n' | + 2\,n - 2\,C > d\left( c_2 (t_n ) , c_2 (t_n' )\right) 
+  d\left(x_n , x_n' \right)\ ,$$
where the last inequality (which is valid as soon as $n > 2 C$) is a consequence of the fact that $d(x,x') \le 2 C$ by the point (i) and the triangle inequality.
Plugging this last series of inequalities in \eqref{projgeod1}, we obtain:
$$d\left(x'_n , c_2 (t_n' )\right) + 2 \, d\left(x_n , x_n' \right) + d\left(x_n , c_2 (t_n )\right) - 6\,\delta \le d\left(x_n , c_2 (t_n )\right) + 
d\left(x_n' , c_2 (t_n')\right) \ \ .$$
This implies that $d\left(x_n , x_n' \right) \le 3\,\delta$ and (going to the limit) that $d\left(x , x' \right) \le 3\,\delta$, this proves that $\text{diameter of } K
\le 3\,\delta$ and, as $K \subset \text{\rm Im} (c_1)$, this achieves the proof of \eqref{lengthK}.

\smallskip
As $c_1 : \R \f X$ is an isometric (proper) injection, the compact sets $K'_t := c_1^{-1} (K_t)$ and $K' := c_1^{-1} (K)$ are compact subset of $\R$
verifying $K' \subset [a , b]$ (by \eqref{lengthK}) and such that $K'$ is the decreasing intersection of the compact sets $K'_{t}$ when $t \to + \infty$;
it follows that, for every $\e > 0$, there exists $A_{\e} \in \R^+$ such that $K'_{t} \subset [a-\e , b+\e]$ for every $t \ge A_{\e}$. 
We choose $x_0 = c_1 \left( \frac{a+b}{2}\right)$ (the middle-point of $[c_1(a), c_1 (b)]$).
For every $t \in \mathbb R$, if $ z= c_1 (\tau)$ is a projection of 
$ c_2 (t)$ on the geodesic line $c_1$, one has $ z\in \pi_{c_1} \left(c_2 ([t , + \infty [)\right) = K_{t} = c_1 (K'_{t})$; hence, for every $t \ge A_{\e}$,
one gets that $\tau \in [a-\e , b+\e]$ and thus $ d(x_0 , z) \le \frac{3}{2}\,\delta + \e$. Using Lemma \ref{projection} and the triangle inequality, we 
conclude that, for every $s \in \R$ and every $t \ge A_{\e}$, 
$$d\left(c_1(s) , c_2 (t)\right) \ge  d\left(c_2 (t) , z\right) + d\left(z , c_1(s)\right) - 2 \, \delta
\ge d\left(c_2 (t) , x_0\right) +  d\left(x_0 , c_1(s)\right) - 5 \, \delta - 2 \e \ ;$$
and, when $\e \to 0$, this achieves the proof of the first inequality of (ii).

\smallskip
Now, defining $\tilde c_2 (t) = c_2 (-t)$ and applying the first inequality of (ii) to the pair of geodesic lines $c_1 , \tilde c_2$, we prove that there exists
a point $x'_0$ on the geodesic line $c_1$ such that, for every $s \in \R$,
$$\limsup_{t \to +\infty} \left[ d\left(c_1(s) , x_0\right) + d \left( x_0 , c_2 (-t )\right) - d\left(c_1(s) , c_2 (-t)\right)  \right] \le 5 \delta $$
and this ends the proof of the second inequality of (ii)
\end{proof}

\subsection{Isometries and displacements}

\subsubsection{Discrete subgroups of the group of isometries}\label{discretude}

\small 
Considering the action (by isometries) of any group $\Gamma$ on any metric space, we recall that, for every $R> 0$, and every $x \in X$, we have defined
$\Sigma_R (x) := \{\g \in \Gamma : \g\,x \in \overline B_X(x, R)\}$.

\normalsize

\medskip
The following Proposition precises the equivalences between a \lq \lq proper" and a \lq \lq discrete" action (see Definitions \ref{proprediscr0} of these two notions)

\begin{prop}\label{discret1}
On a metric space $(X,d)$ every faithful and proper action by isometries is discrete.
Conversely, if $(X,d)$ is a proper metric space, then every faithful and discrete action by isometries is proper.
\end{prop}

\begin{proof} Consider any faithful action (by isometries) of a group $\Gamma$ on any metric space $(X,d)$, this authorizes to view $\Gamma$ as a
subgroup of the group $\text{Isom} (X,d)$ of isometries of $(X,d)$, endowed with the topology of uniform convergence on compact sets.\\
If this action is supposed proper, let $\left( \g_{n}\right)_{n \in \N}$ be any converging sequence of elements of $\Gamma$ and denote by 
$ g \in \text{Isom} (X,d)$ its limit; let $x$ be any point in $X$, there then exists $ N \in \N$ such that, for every $p\in \N$, one has 
$ d \big( x, \g_N^{-1}\g_{N+p}\, x\big) = d\big(\g_N \, x , \g_{N+p}\, x \big) < 1$, 
and the sequence $ p \mapsto \g_N^{-1}\g_{N+p}$, running in the finite set $\Sigma_1 (x)$, is stationary because it converges. This proves that the initial
sequence $ n \mapsto \g_{n}$ is also stationary, thus that $\Gamma$ is a closed discrete subset of $\text{Isom} (X,d)$.\\
Conversely, if $\Gamma$ is a discrete subgroup of $\text{Isom} (X,d)$ and if $(X,d)$ is a proper space, arguing by contradiction, suppose that there exist 
$x_0 \in X$ and $R_0 > 0$ such that $\Sigma_{R_0} (x_0)$ is infinite, Ascoli's Theorem then proves that, for any $k \in \N^*$, $\Sigma_{R_0} (x_0)$ is 
a relatively compact subset of the set of continuous mappings from $ \overline B_X(x_0, k)$ to $X$ (endowed with the topology of uniform 
convergence)\footnote{Indeed, the hypotheses of Ascoli's Theorem are verified for $ \overline B_X (x_0, k)$ is compact and $(X,d)$ (being proper) is 
a complete metric space, $\Sigma_{R_0} (x_0)$, being included in $\text{Isom} (X,d)$, is an equicontinuous family and, for every 
$x  \in \overline B_X(x_0, k)$, $\{\g x  : \g \in \Sigma_{R_0} (x_0)\}$, being included in the compact ball $\overline B_X (x_0, R_0 + k)$, is a 
relatively compact subset of $(X,d)$.}. 
As $\Sigma_{R_0} (x_0)$ is infinite, stable under the mapping $\g \mapsto \g^{-1}$, and relatively compact, it 
contains a sequence of \emph{distinct} elements and, from this sequence, we can extract a subsequence $n \mapsto \g_{n}^k$ such 
that $\left( \g_{n}^k\right)_{n \in \N}$ and $\left( \left(\g_{n}^k\right)^{-1}\right)_{n \in \N}$ uniformly converge on $\overline B_X(x_0, k)$, with 
respective limits continuous mappings $f_k$ and $h_k$ from $\overline B_X(x_0, k)$ to $X$. By successive extractions of subsequences, we can
construct $f_{k+1}$ and $h_{k+1}$ from $\overline B_X(x_0, k +1)$ to $X$ as extensions of $f_k$ and $h_k$ (respectively); this process constructs
continuous mappings $f$ and $h$ from $X$ to $X$, which coincide (for every $k \in \N^*$) with $f_k$ and $h_k$ on $\overline B_X(x_0, k)$.  By a diagonal 
process\footnote{The uniform convergence on $\overline B_X(x_0, k)$ of the two sequences $n \mapsto \g_{n}^{k}$ and 
$n \mapsto \left(\g_{n}^k\right)^{-1}$ imply that, to every 
$k \in \N^*$ corresponds an infinite set $A_k \subset \{ \g_{n}^{k} : n \in \N\}$ whose elements $\g_{n}^{k}$ verify simultaneously 
$ \sup_{x \in \overline B_X(x_0, k)} d \left( f (x) , \g_{n}^{k} (x) \right) < 1/k $ and $\sup_{x \in \overline B_X(x_0, k)} 
d \left( h (x) , \left(\g_{n}^{k}\right)^{-1} (x) \right) 
< 1/k $; we choose $ \g_{n_1}^{1}$ in $A_1$ and, if $\g_{n_1}^{1}, \ldots ,  \g_{n_{k-1}}^{k-1}$ have already been chosen, we choose 
$\g_{n_{k}}^{k}$ as an element of $A_k \setminus \{  \g_{n_1}^{1}, \ldots ,  \g_{n_{k-1}}^{k-1}\}$. In this way, we obtain a sequence $k \mapsto 
\g_{n}^{k}$ of distinct elements such that $\left( \g_{n_k}^{k}\right)_{k \in \N^*}$ and $\left( \left(\g_{n_k}^{k}\right)^{-1}\right)_{k \in \N^*}$ 
uniformly converge to $f$ and $h$ (respectively) on each closed ball.} we can choose, for each $k \in \N^*$, an element $\g^k_{n_k}$ in each sequence 
$\left( \g_{n}^k\right)_{n \in \N}$ such that the two sequences
$k \mapsto \g^k_{n_k}$ and $k \mapsto \left(\g^k_{n_k}\right)^{-1}$ are made of distinct elements and converge respectively to $f$ and $g$ (uniformly
on each compact set).\\
As every $\g_{n_k}^{k}$ and $\left(\g_{n_k}^{k}\right)^{-1}$ are isometries, their uniform limits $f $ and $h$ also preserve distances and also verify
$h \circ f =\lim_{k \to +\infty}  \left(\g_{n_k}^{k}\right)^{-1} \circ \g_{n_k}^{k} = \id_X$ and (by similar arguments) $f \circ h = \id_X$; this proves that
$f$ and $h$ are isometries. Another consequence is that $\lim_{k \to +\infty}  \left(\g_{n_k}^{k}\right)^{-1} \circ \g_{n_{k+1}}^{k+1} = h \circ f = \id_X$
and, as  $\id_X$ is isolated in $\Gamma$, that $\left(\g_{n_k}^{k}\right)^{-1} \circ \g_{n_{k+1}}^{k+1} = \id_X$, thus that $\g_{n_{k+1}}^{k+1} = 
\g_{n_{k}}^{k}$ when $k$ is great enough, in contradiction with the fact that the sequence $\left( \g_{n_k}^{k}\right)_{k \in \N^*}$ is made of distinct 
elements. This contradiction can be avoided only if $ \Sigma_R (x_0)$ is finite for every $R>0$.
\end{proof}

\begin{lemma}\label{autofidele}
Every proper action (by isometries) of a group $\Gamma$ on a metric space $(X,d)$ verifies:
\begin{itemize}
\item[(i)] the quotient space $\Gamma\backslash X$ is a metric space when endowed with the quotient-distance $\bar d$ defined by
$ \bar d (\Gamma \cdot x, \Gamma \cdot y) :=  \inf_{\g \in \Gamma} d( x , \g \, y) $,

\item[(ii)] if $(X,d)$ is a proper space and if $ (\Gamma\backslash X , \bar d)$ has finite diameter, then $\Gamma\backslash X$ is compact,

\item[(iii)] if $\Gamma$ acts without fixed point, then the action is faithful and discrete.

\item[(iv)] if $\Gamma$ is torsion-free, then the action is faithful, discrete and without fixed point
\end{itemize}
\end{lemma}

\begin{proof}
The proof of (i) is classical; the proof of (ii) is straightforward because $(X,d)$ is a proper space. These two proofs are thus left to the reader.\\
\emph{Proof of (iii)}: if the action is fixed-point-free, the stabilizer of every point is trivial, thus the action is faithful and, as it is  proper by hypothesis, it is 
discrete with the help of the first part of Proposition \ref{discret1}.\\
\emph{Proof of (iv)}: The action being proper, the stabilizer of every point is a finite subgroup, which is trivial because $\Gamma$ is torsion-free, this implies 
that the action is faithful and fixed-point-free, thus discrete by the point (iii).
\end{proof}

The following remark is well known by the specialists of group theory

\begin{remark}\label{virtuelcyclique}
Every non trivial, torsion-free, virtually cyclic group is isomorphic to $(\Z , +)$.
\end{remark}

\begin{proof}
As the group $G$ under consideration is non trivial and torsion-free, it is infinite and, being virtually cyclic, it contains a cyclic (infinite) subgroup $H$ of finite 
index in $G$, thus $H$ is isomorphic to $\Z$. For every $g \in G$, define $ \varrho (g) (g' H) = (g g')H $, then $\varrho$ is a morphism from $G$ into the (finite)
group of permutations of $G/H$, whose kernel is a normal subgroup of finite index in $G$ and is contained in $H$. Replacing $H$ by $\Ker \varrho$, we reduce the 
problem to the case where $H$ is a normal subgroup of finite index in $G$ and is isomorphic to $\Z$; we fix one generator of this subgroup, denoted by $\tau$.\\
Denote by $ r : G \f  {\rm Aut}(H)$ the conjugacy morphism defined by $ r(g)  (h) =  g h g^{-1} $ then, for every $g \in G$, $r(g)$ maps $\tau$ onto $\tau$ or
$\tau^{-1}$. Suppose that there exists some $g \in G$ such that $r(g) \ne \id_H$, thus such that $g \tau g^{-1} = \tau^{-1}$, then $g h g^{-1} = h^{-1}$ for 
every $h \in H$; as there exists $n \in \N^*$ such that $g^n \in H$, it follows that $g\, g^n g^{-1} = g^{-n}$, thus that $ g^{2 n} = e$, in contradiction with
the \lq \lq torsion-free" hypothesis. It comes that $ g h g^{-1} = r(g) (h) = h$ for every $(g,h) \in G\times H$, hence $H$ is contained in the center of $G$;
let $ N := \# \big( G / H \big)$, then the transfer homomorphism $ V : G \to H $ verifies $V(h) = h^N$ for every $ h \in H$ and is thus injective when restricted to
$H$, by the \lq \lq torsion-free" hypothesis, and the quotient map $G \f G/H$ induces an injective map $\Ker (V) = \Ker (V) / \big(\Ker (V) \cap H\big) \f  G /H$.
This implies that $\Ker V$ is finite, thus it is trivial by the \lq \lq torsion-free" hypothesis; it follows that $V$ is an isomorphism from $G$ onto some (infinite) subgroup 
of $H$, which is isomorphic to $\Z$.
\end{proof}

\subsubsection{Isometries and displacements in Gromov-hyperbolic spaces}\label{basiques}

Let $(X , d)$  be any $\delta$-hyperbolic (thus geodesic and proper) space, denote by $\partial X $ its ideal boundary\footnote{For two definitions of this ideal 
boundary and of its topology, see Definition 7.1 p. 117 of \cite{GH} and chapter 2 of \cite{CDP}, these two definitions being equivalent by  Proposition 
7.4 p. 120 of \cite{GH}.}. It is well known that every isometry $\gamma$ of $(X , d)$ can be extended as a continuous mapping (moreover Lipschitz for
a convenient metric) from $X \cup \partial X $ to $X \cup \partial X $ (see for example Proposition 11.2.1 p. 134 of \cite{CDP}). An isometry $\g$ of $(X , d)$
is said to be
\begin{itemize}
\item \emph{elliptic} if, for at least one $x \in X$ (thus for every $x \in X$), the sequence $ k \mapsto \gamma^k x$ is bounded,

\item \emph{parabolic} if, for at least one $x \in X$ (thus for every $x \in X$), the sequence $ k \mapsto \gamma^k x$ admits one and only one accumulation 
point, denoted by $\gamma^\infty$, located on the ideal boundary $\partial X $ ($\gamma^\infty$ does not depend on the choice of $x$),

\item \emph{hyperbolic} if, for at least one $x \in X$ (thus for every $x \in X$), the map $ k \f  \gamma^k x$ is a quasi-isometry from $\Z$ to $X$.
\end{itemize}

The following Theorem is classical

\begin{theorem}\label{ellparahyp} \emph{(see for example \cite{CDP}, Th\'eor\`eme 9.2.1 p. 98)}
On any $\delta$-hyperbolic space, every isometry is either elliptic, or parabolic, or hyperbolic.
\end{theorem}

If $\g$ is an hyperbolic (resp. parabolic) isometry, it is a classical result [see for example \cite{CDP}, Proposition 10.6.6 p. 118, (resp. \cite{GH}, 
Th\'eor\`eme 17 in Chapter 8)] that the action (extended as explained above) of $\g$ on $ X \cup \partial X $ admits exactly two (resp. one) fixed points, 
which are the limits $\gamma^+$ and $\gamma^-$ of $ \gamma^p x$ and $ \gamma^{-p} x$ when $ p \to +\infty$ (resp. which is the limit $\gamma^\infty$ 
of $ \gamma^k x$ when $k \f \pm \infty$).\\
The following remark is well known by the specialists; its proof is trivial (and left to the reader) if one notices that, by Lemma \ref{discret1}, every discrete
subgroup of the group $\text{Isom} (X,d)$ of isometries of $(X,d )$ acts properly on $(X,d )$.

\begin{remark}\label{kpointsfixes}
On any $\delta$-hyperbolic space $(X,d )$, if $\Gamma$ is a discrete subgroup of $\text{Isom} (X,d)$, then
\begin{itemize}
\item[(i)] an element of $\Gamma^*$ is elliptic if and only if it has torsion; if $\Gamma$ is torsion-free, every $\g \in \Gamma^*$ is either hyperbolic or parabolic;
\item[(ii)] for every $g \in \Gamma^*$ and every $k \in \Z^*$ such that $g^k \ne id_X$, $g$ is hyperbolic (resp. parabolic, resp. elliptic) if and only if $g^k$ is 
hyperbolic (resp. parabolic, resp. elliptic); moreover, in the cases where $g $ is hyperbolic or parabolic, then $ g^k$ and $g$ have the same set of fixed points.
\end{itemize}
\end{remark}

In case of an hyperbolic isometry $\g$, we introduce (in Definitions \ref{faisceau}) the set ${\cal G} (\g)$ of the oriented (generally not unique) geodesic lines from $ \gamma^-$ to $ \gamma^+$ ($\gamma^- , \gamma^+ \in \partial X$ being the two fixed points of $\g$) and the union $ M(\g) \subset X $ of the 
images of these geodesic lines.

\begin{remark}\label{invariance}
On any $\delta$-hyperbolic space $(X,d )$, for every hyperbolic isometry, if $c \in {\cal G} (\g)$, then $\g^k \circ \, c \in {\cal G} (\g)$  for every $k \in \Z$;
consequently $M(\g)$ is (globally) invariant under $\g^k$.
\end{remark}
The proof, straightforward, is left to the reader.

\begin{defis}\label{deplacements}
To each non trivial isometry $\g$ of a $\delta$-hyperbolic space $(X, d)$, one associates:
\begin{itemize}
\item its  \emph{asymptotic displacement} $\ell(\g)$, i. e. the limit\footnote{By sub-additivity, this limit exists and, by 
the triangle inequality, it does not depend on the point $x$, see for example \cite{CDP}, 
Proposition 10.6.1 page 118).} (when $k \f + \infty$) of $\frac{1}{k} \  d (x, \g^k x) $,

\item its \emph{minimal displacement} $s(\g) :=\inf_{x \in X} d(x, \g \,x)$.
\end{itemize}
\end{defis}

Notice that $\ell(\g^k ) = |k| \, \ell(\g)$ for every $ k \in \Z$ and that \footnote{Indeed, for every $x \in X$ and every $k \in \N$, the triangle inequality gives 
$\dfrac{1}{k} \,d (x, \g^k x) \le d (x, \g \, x)$.\\} $ \ell(\g) \le s(\g)$.

The following Lemma is classical (see for instance \cite{CDP}, Proposition 10.6.3, p. 118): 
\begin{lemma}\label{ellpositive}
On any $\delta$-hyperbolic space $(X,d )$, for every isometry $\g$, $\g$ is hyperbolic if and only if $\ell(\gamma) > 0$.
\end{lemma}

\begin{lemma}\label{puissances} 
On any $\delta$-hyperbolic space $(X,d )$, for every isometry $\g$ and every $x \in X$,
\begin{itemize}
\item[(i)] $ d(x , g^2 x) \le d(x, g\, x) + \ell(g) + 2\,\delta\ $,
\item[(ii)] $\forall p \in \N \ \ \ \ \ d(x , g^{2^{p}} x) \le  d(x, g\, x) + (2^{p}-1)\, \ell(g) + 2\,p\,\delta\ $,
\item[(iii)] $\forall n \in \N^* \ \ \ \ \  d(x , g^n x) < d(x, g\, x) + (n-1)\, \ell(g) + 4\ \delta\, \dfrac{\ln (n)}{\ln (2)} \ $.
\end{itemize}
\end{lemma}

\begin{proof}[Proof of (i)]
The proof of (i) is a simple rearrangement of the arguments used in the proof of the Lemme 9.2.2 pp 98-99 of \cite{CDP}: define $ a_n := d(x , g^n x)$ and
apply Lemma \ref{proprietes} (ii) to the quadrangle with vertices $ x,\ g\,x,\ g^2 x $ and $g^n x$, this gives:
\begin{equation}\label{recurrence}
a_2 + a_{n-1} \le \Max ( a_1 + a_{n-2} \,,\, a_1 + a_{n} ) + 2 \, \delta
\end{equation}
for every $ n \ge 3$, and thus for every $ n \ge 2$.\\
By contradiction, suppose that there exists a point $x$ such that $ d(x , g^2 x) > d(x, g\, x) + \ell(g) + 2\,\delta$, which means that there exists 
$\alpha > \ell(g)$ which verifies
\begin{equation}\label{recurrence1}
a_2 \ge a_1 + \alpha + 2\, \delta \ \ ,
\end{equation}
the inequality \eqref{recurrence} then implies that
\begin{equation}\label{recurrence2}
\forall n \ge 2 \ ,\ \ a_{n-1} + \alpha \le \Max (a_{n-2}\,,\, a_{n})\ .
\end{equation}
Let us prove the property 
\begin{equation}\label{recurrence3}
\forall n \in \mathbb N \ \ \  a_{n+1} \ge  a_{n} + \alpha \ \ :
\end{equation}
the proof of this property is by iteration: indeed the inequality $a_{n+1} \ge  a_{n} + \alpha$ is verified when $n = 1$ (by the inequality \eqref{recurrence1}) 
and also when $n = 0$ because, by \eqref{recurrence1} and the triangle inequality, $a_0 + \alpha + a_1 + 2\,\delta =  \alpha + a_1 + 2\,\delta \le a_2 \le 2\, a_1$ (for $a_0 = 0$). Following the iteration process, assuming the inequality $a_{n+1} \ge  a_{n} + \alpha$, we now prove that 
$a_{n+2} \ge a_{n+1} + \alpha$: indeed, from \eqref{recurrence2}, we get $\Max ( a_{n+2} \, , \,  a_{n}) \ge  a_{n+1} + \alpha \ge a_n + 2\, \alpha$, 
thus that $a_{n+2} \ge  a_{n+1} + \alpha $, and (by iteration) this ends the proof of (\ref{recurrence3}).\\
A consequence of \eqref{recurrence3} is that $a_n \ge n\,\alpha$, thus that $\lim_{n \to + \infty} \left( \frac{1}{n}\,d(x, g^n x) \right) \ge \alpha > \ell(g)$,
in contradiction with the Definition \ref{deplacements} of $\ell(g)$. The only solution to avoid this contradiction is thus that, for every $x \in X$, one has
$d(x , g^2 x) \le d(x, g\, x) + \ell(g) + 2\,\delta $.
\end{proof}
\begin{proof}[Proof of (ii)]
Let us denote by $({\cal Q}_p)$ the property \lq \lq $ \,d(x , g^{2^{p}} x) \le  d(x, g\, x) + (2^{p}-1)\, \ell(g) + 2\,p\,\delta\,$", this property is trivially
verified when $p=0$ and, by (i), when $p=1$; if $({\cal Q}_p)$ is verified, then $({\cal Q}_{p+1})$ is satisfied also because, using (i) and the fact that
$\ell(g^k) = |k| \,\ell(g)$, one gets
$$ d(x , g^{2^{p+1}} x) \le d(x, g^{2^{p}}\, x) + \ell(g^{2^{p}}) + 2\,\delta \le  d(x, g\, x) + (2^{p+1}-1)\, \ell(g) + 2\,(p+1)\,\delta\ ,$$
which proves $({\cal Q}_p)$ for every $p \in \N$, and thus proves (ii).
\end{proof}
\begin{proof}[Proof of (iii)]
Let us define the property 
$$ \big({\cal P}_p \big) \ \  : \ \ \  \forall n \in \N^* \  \text{ such that }\ 2^{p} < n < 2^{p+1}\ ,\ \ \ \ d(x , g^n x) \le d(x, g\, x) + (n-1)\, \ell(g) + 4\,\delta \, p \ .$$
As (ii) solve the case where  $n = 2^{p}$, in order to prove (iii), it is sufficient to verify property $\big({\cal P}_p \big)$ for every $p \in \N^*$.
For every $p \in \N^* $ and every $ n $ such that $ 2^{p} < n < 2^{p +1}$, the quadrangle Lemma \ref{proprietes} (ii) implies that
$$ d(x , \, g^{n} x)  + d( g^{2^{p}} x ,\, g^{2^{p + 1}} x) - 2 \, \delta \le \Max \left[ d(x,\, g^{2^{p}} x) + 
d(g^{n} x,\, g^{2^{p+1}} x)  \, ; \, d(x,\, g^{2^{p+1}} x) + d(g^{n} x,\, g^{2^{p}} x)  \right]$$
$$\le \Max \left[ d(x,\, g^{2^{p}} x) + d(g^{n} x,\, g^{2^{p+1}} x)  \, ; \, d(x, g^{2^{p}}\, x) + 
2^{p}\, \ell(g) + 2\,\delta+ d(g^{n} x,\, g^{2^{p}} x)  \right] \ ,$$
where the second inequality comes from point (i). Simplifying, we obtain:
\begin{equation}\label{powern}
d(x , \, g^{n} x)  \le \Max \left[  d( x,\, g^{2^{p+1} - n} x)  \, ; \, 
2^{p}\, \ell(g) + 2\,\delta+ d( x ,\, g^{n -2^{p}} x)  \right] + 2 \, \delta \ .
\end{equation}
Applying the inequality \eqref{powern} to the case $ p = 1$ and $n = 3$, it comes:
$$d(x , \, g^{3} x) \le 2 \, \delta + \Max \left[  d( x,\, g\,  x)  \, ; \, 
2\, \ell(g) + 2\,\delta+ d( x ,\, g\, x)  \right] = d( x ,\, g\, x) + 2\, \ell(g) + 4\,\delta \ ,$$
and this proves property $\big({\cal P}_1 \big)$. We now show that, for every  $p \ge 2$, if $\big({\cal P}_1 \big) , \ldots ,\big({\cal P}_{p-1} \big)$ are all 
verified, then $\big({\cal P}_p \big)$ is verified too: thus, assuming that $\big({\cal P}_1 \big) , \ldots ,\big({\cal P}_{p-1} \big)$ are verified and property (ii),
the inequality $d(x , g^i x) \le d(x, g\, x) + (i-1)\, \ell(g) +  4\,\delta \, (p-1)$ is verified for every $i \in \N$ such that $1 \le i \le 2^p$, and this implies that, 
for every $ n \in \N^* $ such that $ 2^{p} < n < 2^{p+1}$, then $ 1 \le 2^{p+1} - n , n -2^{p} \le 2^p -1$ and we get that
$$ d( x,\, g^{2^{p+1} - n} x) \le  d( x,\, g\, x) + \left( 2^{p+1} - n - 1\right) \ell(g) + 4\,  (p-1) \,\delta \ ,$$
$$ d( x,\, g^{n - 2^{p} } x) \le  d( x,\, g\, x) + \left( n - 2^{p} - 1\right) \ell(g) + 4\,  (p-1) \,\delta \ ,$$
Plugging these two last inequalities in \eqref{powern} and noticing that $ n-1 \ge 2^{p+1} - n$, it follows that
$$d(x , \, g^{n} x)  \le  d( x,\, g\, x) + 4\, p\,\delta + \Max \left[ \left( 2^{p+1} - n - 1\right) \ell(g) - 2\, \delta \, ; \, 
 \left( n - 1\right) \ell(g) \right] \le d( x,\, g\, x) + ( n  - 1 ) \, \ell(g)+ 4\, p\,\delta\,.$$
This proves that $\left[\big({\cal P}_1 \big) \text{ and } \ldots  \text{ and } \big({\cal P}_{p-1} \big)\right] \implies  \big({\cal P}_{p} \big )$.
This proves property $\big({\cal P}_p \big)$ for every $ p \in \N^*$ and ends the proof as announced above.
\end{proof}

\begin{lemma}\label{milieux} 
On any $\delta$-hyperbolic space $(X,d )$, for every isometry $g$ and every $x \in X$, if $m$ is the middle-point of some geodesic $[ x , g \,x ]$ from $x$ 
to $gx$, then
$$ \Max \left[ 0 \,,\,  d(x , g^2 x) - d(x, g\, x) \right] \le  d(m, g\, m) \le \Max \left[ 0 \,,\,  d(x , g^2 x) - d(x, g\, x) \right] + \delta \le  \ell(g) + 3\,\delta \ .$$
\end{lemma}

\begin{proof}
The proof consists in updating the proof of Lemma 9.3.1 p. 101 of \cite{CDP}. After having chosen two geodesic segments $[x , g x ]$ and $[x , g^2  x ]$, 
we fix the geodesic segment $[g  x , g^2 x ]$ as the image by $g$ of $[x , g  x ]$, then $g m$ is the middle-point of $[ g\, x , g^2 x ]$.
Apply the approximation Lemma \ref{prodist} to the geodesic triangle $\Delta = [x\, , g\,  x \,,  g^{2} x]$ which is the union of these three geodesic segments
and to the approximation map $ f_\Delta$ from $\Delta $ to the metric tripod $(T_\Delta , d_T)$; denote by $\alpha ,\ \beta$ and $\gamma$ the lengths of
the branches of $T_\Delta$ with respective endpoints $ f_\Delta (x) , \ f_\Delta ( g \, x )$ and $ f_\Delta ( g^{2} x)$, and by $c$ the branching point of this
tripod. Lemma \ref{prodist} guarantees that
$$ \alpha + \beta = d_T\big(  f_\Delta (x) ,  f_\Delta (g\,x) \big) = d(x , g\,x) = d( g\, x , g^2 x) =
d_T\big(  f_\Delta (g\,x) ,  f_\Delta (g^2x) \big) = \beta + \gamma \ ,$$
thus that $\gamma = \alpha$; Lemma \ref{prodist} also implies that
\begin{equation}\label{milieu}
d_T \big( f_\Delta (m) ,  f_\Delta (x) \big)  = d_T \big( f_\Delta (m) ,  f_\Delta (g\, x) \big) =
\frac{\alpha + \beta}{2}  = d_T \big(f_\Delta ( g \,m) ,  f_\Delta (g \,x) \big) \ .
\end{equation}
\begin{itemize}
\item If  $\beta \ge \alpha$, it follows from \eqref{milieu} that $f_\Delta (m)$ and $f_\Delta (g\, m)$ both lie on the 
same branch $ [c \, , f_\Delta (g\,x)  ]$ of the tripod and that $f_\Delta (m) = f_\Delta (g\, m)$ because $d_T \big(f_\Delta ( g \,m) ,  f_\Delta (g \,x) \big) =
d_T \big( f_\Delta (m) ,  f_\Delta (g\, x) \big) $; by Lemma \ref{proprietes} (i), this implies that $ d(m, g\, m) \le \delta$

\item  If  $\beta < \alpha$ (i. e. if $ d(x, g^2x ) >  d(x, g\,x )$),  it follows from \eqref{milieu} that $f_\Delta (m)$ and $f_\Delta (g\, m)$ respectively lie on 
the branches $ [c \, , f_\Delta (x)  ]$ and $ [c \, , f_\Delta (g^2 x)  ]$ of the tripod; from this and from Lemmas \ref{prodist} and \ref{proprietes} (i), we 
deduce that
$$ d(x, g^2x ) - d(x, g\,x ) =  \gamma - \beta  = d_T \big(f_\Delta ( x) ,  f_\Delta (g^2x) \big) - d_T \big(f_\Delta ( m) ,  f_\Delta (x) \big) - d_T \big(f_\Delta ( g \,m) ,  f_\Delta (g^2 x) \big) $$
$$ =d_T \big( f_\Delta ( m)  , f_\Delta ( g \,m) \big) \le d (m,g \, m) \ .$$
By the approximation Lemma \ref{proprietes} (i), it then comes that
$$ d(x, g^2x ) - d(x, g\,x ) \le d (m,g \, m) \le d_T \big( f_\Delta ( m)  , f_\Delta ( g \,m) \big) + \delta \le d(x, g^2x ) - d(x, g\,x ) + \delta\ ,$$
\end{itemize}
This proves the first and second inequalities of Lemma \ref{milieux}, the third one then follows from Lemma \ref{puissances} (i)
\end{proof}

\begin{lemma}\label{Busemann}
On a $\delta$-hyperbolic space $(X , d)$, every parabolic isometry $\g$ satisfies the following property: for every geodesic ray $c$ such that $c (+\infty)
\in \partial X$ 
is the single fixed point $\g^\infty $ of $\g$, there exists $T = T(\g) \le d \big(c(0) ,  \g \circ\, c(0) \big)$ such that 
$d\left( c(t), \g\big(c(t)\big) \right) \le 7\delta $ for every $ t \in [T , + \infty[$.
\end{lemma}

\begin{proof}
As $\g\circ\, c(+\infty)=\g( \g^\infty)=\g^\infty = c(+\infty)$, we may apply Proposition \ref{geodasympt} (iii) to the geodesics $c$ and 
$ \g \circ \, c$ and conclude that there exist $t_1 , \, t_2 \ge 0$, verifying $ t_1 + t_2 = d \big( c (0) ,  \gamma \circ c (0)\big)$, 
such that $  d \left(c (t_1 + s) ,  \g \circ \, c \, (t_2 + s) \right) \le 2\, \delta$ for every $ s \ge 0 $. Let us define $ T := \max (t_1,t_2)$, then two cases are possible:
\begin{itemize}
 \item If $| t_2 - t_1| \le 3\, \delta $ then, for every $ s \in [0 , + \infty[$, the last inequality and the triangle one give:
$$  d \big( c ( t_2 + s) ,   \g \circ \, c ( t_2 + s) \big) \le d \big( c ( t_2 + s) ,  c (t_1 + s) \big) + 
d \big( c (t_1 + s) ,   \g \circ \, c ( t_2 + s) \big)\le | t_2 - t_1|  + 2\,\delta  \le 5 \delta \ ,$$ 
and, making $s := t - t_2$, which is nonnegative when $t \ge T$, this ends the proof in this case.

\item If $| t_2 - t_1| > 3\, \delta $, for every $ s \in [0 , + \infty[$, let us consider a geodesic triangle \linebreak
$\Delta = \big[  c (t_1 + s) \,,\, \g \circ \, c (t_2 + s)\, ,\, c(t_2 + s) \big]$ and its approximation by the associated tripod $ f_\Delta : (\Delta , d) \to  (T_\Delta , d_T)$, the branching point of this tripod being denoted by 
$e$ and the lengths of the three branches with respective ending points  $  f_\Delta \big( c (t_1 + s) \big)$, $ f_\Delta ( \g \circ \, c(t_2 + s) )$ and 
$  f_\Delta ( c(t_2 + s) )$ being denoted by $ \alpha$, $ \beta$ and $ \gamma$.\\
Lemma \ref{prodist} and the previous inequality $ d \left(c (t_1 + s) ,  \g \circ \, c \, (t_2 + s) \right) \le 2\, \delta$ yield
$\alpha + \gamma = | t_2 - t_1 |$, $\alpha + \beta \le 2\,\delta$, $\gamma + \beta = d \big(c(t_2 + s),  \g \circ \, c (t_2 + s) \big)$, and consequently 
$\gamma - \beta = (\alpha + \gamma)-(\alpha + \beta) \ge | t_2 - t_1 |-  2\,\delta > \delta > 0$.
This last inequality implies that, if $m_s$ is the middle point of the geodesic side $ \left[ c (t_2 + s) ,   \g \circ \, c \, (t_2 + s) \right]$ of $\Delta$, the point 
$ f_\Delta ( m_s)$ lies on the branch $ \big[ e , f_\Delta ( c(t_2 + s) ) \big]$ of the tripod $T_\Delta$; there thus exists a point $m'_s = c \big(\tau (s) \big) $ of the geodesic side 
$ \big[  c(t_2 + s) , c (t_1 + s) \big]$ of $\Delta$ such that $ f_\Delta ( m'_s) = f_\Delta ( m_s)$.
Applying Lemma \ref{proprietes} (i), we obtain that 
$d(m_s, m'_s) \le \delta$. From this, from the triangle inequality and from Lemma \ref{milieux}, for every $s \ge 0$, we deduce:
\begin{equation}\label{controletau}
 d \big(c \big(\tau (s) \big) , \g\, c \big(\tau (s) \big) \big)  = d(m'_s, \g\, m'_s) \le d(m_s, \g\, m_s) + 2 \, \delta \le  \ell(\g) + 5\,\delta = 5\,\delta\, ,
\end{equation}
where this last equality comes from the fact that $\ell (\g) = 0$ for $\g$ is parabolic.\\
As $ s + \min (t_1, t_2)\le \tau (s) \le s + \max (t_1, t_2)$, it follows first that $\tau (s) \le T+s$, on the other hand, defining $s' := s + |t_2 - t_1|$
one has 
$$\tau (s') = \tau (s + |t_2 - t_1|) \ge s + |t_2 - t_1| + \min (t_1 , t_2) = s +  \max (t_1, t_2) = T+ s \, .$$
As $\tau (s) \le T+s \le \tau (s') $, from \eqref{controletau} and from the convexity (modulo an error $\le 2 \delta$) of the distance between geodesics in a 
$\delta$-hyperbolic space, we draw, for every $s\ge 0$,
$$d \big( c(T+s) , \g \circ \, c (T+s) \big) \le \Max \left[ d \big( c\big(\tau (s)\big) , \g \circ \, c\big(\tau (s) \big) \big) \,;\, d \big( c\big(\tau (s')\big) , 
\g \circ\,  c\big(\tau (s' )\big) \big) \right] + 2\, \delta\le 7\,\delta\,, $$
and, making $s := t - T$, which is nonnegative when $t \ge T$, this ends the proof in this case.
\end{itemize}
\end{proof}

\begin{lemma}\label{quasigeod}
Let $g$ be any non trivial isometry of any $\delta$-hyperbolic space $(X, d)$, then
\begin{itemize}
\item [(i)] $\ell(g) \le s(g) \le \ell(g) + \delta$,
\item [(ii)] if $s (g) > \delta$, for every $p \in \N^*$, one has the property:
\begin{equation*}
\big({\cal Q}_p\big)\ \ \ \ : \ \ \ \  \forall x \in X\ \ \ \ d \big(x, g^{2^p} x \big)  \ge d (x, g\, x) + \left( 2^{p} -1\right)  \left(s(g) - \delta \right) 
\end{equation*}
\item [(iii)] if $s (g) > 3 \,\delta$, for every $x \in X$, the sequence $ \big( d( x , g^n x)\big)_{n \in \N^*}$ is strictly increasing.
\end{itemize}
\end{lemma}

\begin{proof}
Denote by $m_x$ the middle-point of any geodesic $[x , g\,x ]$ from $g$ to $g\,x $. If $s (g) > \delta$, one has $d(m_x, g\, m_x) > \delta$ and Lemma
\ref{milieux} then proves that
\begin{equation}\label{quasigeod1} 
\forall x \in X \ \ \ \ s(g) -  \delta \le d(m_x, g\, m_x) -  \delta \le d(x , g^2 x) - d(x, g\, x) \  .
\end{equation}
We are going to prove (by iteration) property $\big({\cal Q}_p\big)$ for every $p$, an immediate corollary being:
\begin{equation}\label{quasigeod2}
s \big( g^{2^p} \big)  - \delta \ge  2^{p}   \left( s(g) - \delta \right) \ .
\end{equation}
Property $\big({\cal Q}_1\big)$ is verified by \eqref{quasigeod1}; suppose $\big({\cal Q}_p\big)$ to be verified, then its corollary \eqref{quasigeod2} is
also verified, in particular one has $s \big( g^{2^p} \big) > \delta$ and one may apply inequality \eqref{quasigeod1}, and afterwards $\big({\cal Q}_p\big)$ 
and its corollary, which give:
$$d \big(x, g^{2^{p+1}} x \big)  \ge d \big(x, g^{2^p} x \big) + s \big( g^{2^p} \big) - \delta \ge d (x, g\, x) + \left( 2^{p+1} -1\right)  \left(s(g) - \delta 
\right) \ ;$$
a consequence is that $\big({\cal Q}_p\big) \implies \big({\cal Q}_{p+1}\big)$ for every $p \in \N^*$, and thus that property $\big({\cal Q}_p\big)$ is verified 
for every $p \in \N^*$. This proves (ii).

\smallskip
The proof of the inequality $ \ell(g) \le s(g)$ has been given after Definitions \ref{deplacements}. Let us now prove that $ s(g) \le \ell(g) + \delta$: if
$s(g) \le \delta$ this inequality is trivially verified, let us thus suppose that $ s(g) > \delta$; this allows to apply property $\big({\cal Q}_p\big)$, which gives:
$ \ell(g) = \lim_{p \to +\infty} 2^{-p}\, d \big(x, g^{2^p} x \big) \ge s(g) - \delta$. This proves (i).

\smallskip
By (i), the hypothesis $s(g)> 3\, \delta $ implies that $\ell(g) > 2 \, \delta$, thus $g$ is an hyperbolic isometry by Lemma \ref{ellpositive}.
Applying \eqref{quasigeod1}, and afterwards the quadrangle Lemma \ref{proprietes} (ii), we obtain:
$$ d \big(x \,, g^{n+1} x \big) + d(x, g\,x ) + s(g) - \delta \le d \big(x \,, g^{n+1} x \big) + d \big(x \,, g^{2} x \big) =
d \big(x \,, g^{n+1} x \big) + d \big(g^{n} x \,, g^{n+2} x \big) $$
$$\le \Max \left[d \big(x \,, g^{n} x \big) + d \big(g^{n+1} x \,, g^{n+2} x \big) \, ,\, 
d \big(x \,, g^{n+2} x \big) + d \big(g^{n} x \,, g^{n+1} x \big)\right]+ 2\,\delta \ .$$
From this, defining $d_n = d \big(x \,, g^{n} x \big) $, we deduce that $ d_{n+1} \le \Max \left[ d_{n} , d_{n+2}\right] + 3\, \delta - s(g) < 
\Max \left[ d_{n} , d_{n+2}\right]$ for $s(g)> 3\, \delta $; as $ d_1 < d_2$ by \eqref{quasigeod1}, this implies (arguing by iterations) that the sequence 
$ \big( d_n\big)_{n \in \N^*}$ is strictly increasing. This proves (iii).
\end{proof}

The following Lemma is an explicit and quantified version of the lemma 9.2.3 p. 98 of \cite{CDP}:

\begin{lemma}\label{minorell}
Let $(X,d)$ be a $\delta$-hyperbolic space, $\alpha$ be any strictly positive real value, 
and $\{g,h\}$ be any pair of isometries of $(X,d)$ which satisfies the property:
\begin{equation}\label{minorell1}
\exists x \in X \ \ \  s. t.\ \ \ d(g x , h x) \ge \Max \,[d(x,g x)\,;\, d(x, h x) ] + 5\,\delta + \alpha \ ,
\end{equation}
then either $\ell (g) \ge 3\delta + \alpha$, or $\ell (h) \ge 3\delta + \alpha$, or $\ell (g h) = \ell (h g) \ge 2\, \alpha$.
\end{lemma}

\begin{proof} The proof mimics the one of Lemma 9.2.3 of \cite{CDP} (pp. 98-101).\\
Let us suppose that $\ell (g),\, \ell(h) < 3\delta + \alpha$ (if not, the lemma is automatically verified); applying the
quadrangle Lemma \ref{proprietes} (ii) to the points $x,\,g x,\, g^2 x,\, g h x$, we get:
$$ d(x,g x) + d(g x , h x) - 2\,\delta = d(x,g x) + d(g^2 x , g h x) - 2\,\delta $$
$$\le \Max \left[ d(x, g^2 x) +  d(g x , g h x ) \, ;\, d(g x, g^2 x) +  d( x , g h x ) \right]$$
$$\le \Max \left[ d(x, g x)  + \ell(g) + 2 \, \delta +  d(x , h x ) \, ;\, d(x, g x) +  d( x , g h x ) \right] \ ,$$
where the last inequality is derived from Lemma \ref{puissances} (i).
As, by assumption, 
$$d(g x , h x) \ge d(x, h x) + 5\,\delta + \alpha > d(x, h x) + \ell(g) + 2\delta \ ,$$
the above inequality gives: $ d(g x , h x) \le  d( x , g h x ) + 2 \, \delta$; by a similar proof (exchanging the names of
$g$ and $h$) we get $ d(g x , h x) \le  d( x , h  g x ) + 2 \, \delta$ and thus, using the assumption \eqref{minorell1}:
\begin{equation*}
\Min[ d(x, g h x) \, ; \, d(x, h g x) ] \ge \Max \,[d(x,g x)\,;\, d(x, h x) ] + 3\,\delta + \alpha \ .
\end{equation*}
This and the quadrangle Lemma \ref{proprietes} (ii) (applied to points $x,\,g x,\, g h x , \, g h g x$) imply:
$$2\, \Max \,[d(x,g x)\,;\, d(x, h x) ] + 4\,\delta + 2\, \alpha \le d(x, g h x) + d(x, h g x) - 2\,\delta$$ 
$$= d(x, g h x) + d (g x,  gh g x)- 2\,\delta 
\le  \Max \left[ d(x, g x) +  d(g h x , g h g x ) \, ;\, d(x, g h g x) +  d( g x , g h x ) \right]  $$
$$=  \Max \left[ 2 \,d(x, g x)\, ;\, d(x, g h g x) +  d( x ,  h x ) \right] = d(x, g h g x) +  d( x ,  h x ) \ ,$$
which implies that
\begin{equation*}
 d(x, g h g x) \ge \Max \,[d(x,g x)\,;\, d(x, h x) ] + 4\,\delta + 2 \,\alpha \ .
\end{equation*}
This and the quadrangle Lemma \ref{proprietes} (ii) (applied to the points $x,\, g h x , \, g h g x , \, g h g h x$) yields:
$$ d(x, g h x) + \Max \,[d(x,g x)\,;\, d(x, h x) ] + 2\,\delta + 2 \,\alpha \le d(x, g h x) + d(x, g h g x) - 2\, \delta$$
$$= d \big(g h  x , (g h)^2 x \big) +  d \big(x, g h g x \big) - 2\, \delta \le 
\Max \left[ d(x, g h x) +  d \big(g h g x , (g h)^2 x \big) \, ;\, d \big( x , (g h)^2 x \big) +  d( g h x , g h g x ) \right] $$
$$ = \Max \left[ d(x, g h x) +  d \big(x , h x \big) \, ;\, d \big( x , (g h)^2 x \big) +  d(  x ,  g x ) \right] 
= d \big( x , (g h)^2 x \big) +  d(  x ,  g x ) \ . $$
This last inequality and the lemma \ref{puissances} (i) imply that:
$$ d(x, g h x) + 2\,\delta + 2 \,\alpha \le d \big( x , (g h)^2 x \big)  \le  d(x, g h x) + 2\,\delta + \ell (gh) \ ,$$
and we conclude that $ \ell (hg) = \ell (h (gh) h^{-1}) = \ell (gh) \ge 2 \,\alpha$.
\end{proof}
\subsection{Margulis' domains of an hyperbolic isometry}

\begin{defis}\label{faisceau}
On any $\delta$-hyperbolic space $(X,d )$, for every hyperbolic isometry $\g$, we denote by $ \gamma^-$ and $ \gamma^+$ the fixed points of $\g$,
by ${\cal G} (\g)$ the set of the geodesic lines such that $ c(-\infty) = \gamma^-$ and $ c(+ \infty) = \gamma^+$ and by $M(\g)$ the subset of 
$X$ obtained as union of these geodesic lines.\\
On the other hand, $M_{\rm{min}} (\g)$ is the name of the (non empty\footnote{We shall prove in Lemma \ref{tube} (iv) that this set is non empty and, 
by continuity, closed.}) set of points of $X$ where the function $ x \f d(x, \g\, x)$ attains its infimum $ s(\gamma)$.
\end{defis}

Though they generally do not coincide, the sets $M(\g)$ and $M_{\rm{min}} (\g)$ are equal when the distance is convex (see Lemma \ref{invariantgeod}).
Furthermore, in the case of Hadamard spaces or of $\CAT (-1)$ spaces, the sets $M(\g)$ and $M_{\rm{min}} (\g)$ coincide with the unique geodesic line
which joins points $ \gamma^-$ and $\gamma^+$.

\begin{defi}\label{rayondeplact}
To each isometry $\g$ of $(X, d)$ and to each point $x \in X$, one associates its \emph{displacement radius at $x$},
defined as $R_{\g} (x) = \inf_{n \in \N^*} d(x, \g^n \,x) $. This infimum is attained when $\g$ is non elliptic.
\end{defi}

Notice that $R_{\g} \equiv 0$ when $ \g$ has torsion, we shall always exclude this case in the sequel.\\ 
The proof of the fact that the infimum is attained is as follows: when $\g$ is non elliptic, it is parabolic or hyperbolic
by Theorem \ref{ellparahyp}, and then $g^k x$ goes to infinity when $k \f \pm \infty$; this proves that $\{ k \in \Z^* : d(x, \g^k \,x)  \le R_{\g} (x) +1\}$ 
is finite; there thus exists $n\in \N^*$ such that $d(x , \g^n \, x ) = \inf_{n \in \N^*} d(x, \g^n \,x) =  R_\gamma (x) $.

\begin{lemma}\label{geodmin}
For every $ x \in M_{\rm{min}} (\g)$ and any choice of a geodesic segment $ \left[  x , \g  x \right] $ from $x$ to $\g x$, the union $\cup_{p \in \Z}\g^p 
\left(\left[  x , \g \, x \right] \right)$ is a $\g$-invariant local geodesic included in $M_{\rm{min}} (\g)$.
\end{lemma}

\begin{proof}
Denote by $\left[ \g^p x \,, \g^{p+1} x \right]$ the geodesic $ \g^p \left(\left[  x , \g \, x \right] \right)$ from $\g^p x $ to $\g^{p+1} x $; for every 
$:u \in \left[ \g^p x \,, \g^{p+1} x \right]$, one has
$$d(u, \g u) \le d(u, \g^{p+1} x ) + d(\g^{p+1} x , \g u) =  d(\g^{p} x , u) + d(u, \g^{p+1} x ) = 
 d(\g^{p} x , \g^{p+1} x ) = d(x, \g x) = s(\gamma)\, .$$
This implies first that $d(u, \g u) = s(\gamma)$, thus that $ u \in M_{\rm{min}} (\g)$ and, in a second time, that the above inequalities are all equalities, thus 
that $d(u, \g u) =  d(u, \g^{p+1} x ) + d(\g^{p+1} x , \g u) $, which means that $ [ u, \g^{p+1} x ] \cup [ \g^{p+1} x , \g u ]$ is a (minimizing) geodesic 
from $u$ to $\g \, u$. It follows that  $\cup_{p \in \Z}\g^p \left(\left[  x , \g \, x \right] \right)$ is a continuous path which is minimizing on any sub-path of 
length $ s(\gamma)$, and consequently it is a local geodesic.
\end{proof}

\begin{lemma}\label{distgeod}
If $  \ell(\g) > 3\, \delta$, one has $d \big(x, M(\g) \big) \le \frac{1}{2} \big(d(x, \g \,x) - \ell(\g)\big) + 3\,\delta $ for every $x \in X$.
\end{lemma}

\begin{proof}
Let us consider any geodesic $ c_\varepsilon \in \cal{G} (\g)$ such that $d (x, c_\varepsilon ) \le d \left(x,  M(\g) \right) + \frac{1}{2}\,\e$. For every 
$ k \in \Z$, denote by $c_\varepsilon (t_k)$ a projection of the point $\g^k x$ onto the image of the geodesic $c_\varepsilon$. We now prove the property 
\begin{equation}\label{shift}
\exists p \in \N  \  \text{tels que} \ d \big(c_\varepsilon (t_p) , c_\varepsilon (t_{p+1}) \big) > \ell (\g) -  \frac{1}{2}\,\e \ .
\end{equation}
Indeed, let $c_\varepsilon (t_{k}')$ be a projection of $\g^k (  c_\varepsilon (t_{0} )  )$ on the image of the geodesic $c_\varepsilon$, Proposition \ref{geodasympt} (i) (and the fact that $c_\varepsilon$ and $ \g^k \circ \, c_\varepsilon$ are two geodesics from $\gamma^-$ to $\gamma^+$) 
guarantees that $d \left( \g^k c_\varepsilon (t_{0} ) \,,\, c_\varepsilon (t_{k}')\right) \le 2\, \delta$, and this yields
$$ d \left( \g^k x , c_\varepsilon (t_{k})\right) \le  d \left( \g^k x , c_\varepsilon (t_{k}')\right) \le d \left( \g^k x , \g^k c_\varepsilon (t_{0} ) \right) + 
d \left( \g^k c_\varepsilon (t_{0} ) \,,\, c_\varepsilon (t_{k}')\right) \le d (x, c_\varepsilon ) + 2 \delta \ ;$$
we deduce that $d \big(  c_\varepsilon (t_{0})  ,  c_\varepsilon (t_{k}) \big) \ge d \big( x , \g^k x \big) - d \big( x , c_\varepsilon (t_{0})\big) - 
d \left( \g^k x , c_\varepsilon (t_{k})\right) \ge  d \big( x , \g^k x \big) - 2\, d (x, c_\varepsilon ) - 2\, \delta $, thus that
$$ \lim_{k \f +\infty} \left(\frac{1}{k}\, \sum_{i = 0}^{k - 1} d \big(  c_\varepsilon (t_{i})  ,  c_\varepsilon (t_{i+1}) 
\big) \right) \ge \lim_{k \f +\infty} \left( \frac{1}{k}\,  d \big(  c_\varepsilon (t_{0})  ,  c_\varepsilon (t_{k}) \big) \right)
= \lim_{k \f +\infty} \left(  \frac{1}{k} \, d \big( x , \g^k x \big)\right) = \ell(\g) \,;$$
a consequence is that $ \sup_{p \in \N} \  d \big(  c_\varepsilon (t_{p})  ,  c_\varepsilon (t_{p+1}) \big) \ge  \ell(\g)$, and this proves property \eqref{shift}.

\smallskip
When $ \ell (\g) > 3 \,\delta $, choosing $\e$ small enough, property \eqref{shift} implies the existence of some $p\in \N$ such that 
$d \big(c_\varepsilon (t_p) , c_\varepsilon (t_{p+1}) \big) > \ell (\g) -  \frac{1}{2}\,\e > 3 \,\delta $, and it then follows from Lemma \ref{ecartement} that
$$ d \left( x ,  \g \, x\right) = d \left( \g^{p} x ,  \g^{p +1} x\right) \ge d \left( \g^{p} x ,  c_\varepsilon (t_p)\right) + 
d \big(c_\varepsilon (t_p) , c_\varepsilon( t_{p +1})\big)  + d \left( c_\varepsilon( t_{p +1}) , \g^{p +1} x \right) - 6
\, \delta$$
\begin{equation}\label{distgeod1}
\ge 2\, d \big( x,  M(\g) \big) + \ell(\g) - \frac{1}{2}\,\e  - \,6\, \delta \ ,
\end{equation}
where the last inequality is a consequence of the inequality $d \left( \g^{p} x ,  c_\varepsilon (t_p)\right) = d \left( x , \g^{-p}\circ \, c_\varepsilon \right) \ge 
d \big( x,  M(\g) \big) $ (for the geodesic $\g^{-p}\circ \, c_\varepsilon $ belongs to $\cal{G} (\g)$). We end the proof by making $\e$ go to zero in the
inequality \eqref{distgeod1}.
\end{proof}

\begin{defi}\label{Margudomain}
In a $\delta$-hyperbolic space $(X,d)$, for every isometry $\g \ne \id$, every $R \in \, ] 0 , +\infty[$ and every $p \in \Z^* $, one defines
\begin{itemize}
\item the $p$-th Margulis domain $M^p_R (\g) := \{ x \in X : d(x , \g^p \, x ) \le R\}$,
\item the Margulis domain $M_R (\g) := \cup_{k \in \mathbb Z^*} M^k_R (\g)$.
\end{itemize}
\end{defi}

The following remark is an immediate consequence of the definition of $ R_\gamma (x) $ as $\inf_{n \in \N^*} d(x, \g^n \,x) $ and of the fact that this
infimum is attained (see Definition \ref{rayondeplact}).

\begin{remark}\label{Margudomain1} 
In a $\delta$-hyperbolic space $(X,d)$, for every non elliptic isometry $\g$ and every $R \in \, ] 0 , +\infty[$, one has 
$M_R (\g) = \{ x : R_\gamma (x) \le R\}$.
\end{remark}

\begin{lemma}\label{MRnontout} In every non elementary $\delta$-hyperbolic space $ (X,d)$, for every hyperbolic isometry $\gamma$ and every real 
number $R> 0$, $ X \setminus M_R (\g) \ne \emptyset$.
\end{lemma}

\begin{proof}
Observe that $\ell(\g) > 0$ by Lemma \ref{ellpositive}. Define $k_0 := \left[\frac{R}{\ell (\g)}\right]$ and notice that $M^k_R (\g) = \emptyset$ 
for every integer $ k \ge  k_0 +1$, because every $x \in X$ verifies $  d(x , \g^k \, x ) \ge \ell (\g^k) \ge (k_0 + 1) \, \ell (\g)> R$ (see the definition 
\ref{deplacements} of $\ell (\g)$ and the properties following this definition).; consequently $M_R (\g) = \cup_{0 < k \le k_0} \, M^k_R (\g)$.
Denote by $\{ \g^-  ,\g^+\}$ the set of fixed points of $\g$, it is included in the ideal boundary $\partial X$ of $(X,d)$ (see the beginning of section \ref{basiques}); choose some point $\theta \in \partial X \setminus \{ \g^-  ,\g^+\}$ and any geodesic ray $c_0$ such that $c_0 (+ \infty) = \theta$.
For every $ k \in \N^*$, the geodesic ray $\g^k \circ \,c_0$ verifies $\g^k \circ \,c_0 (+ \infty) = \g^k \theta \ne \theta $, thus 
$ d \big(c_0(t) \, , \g^k \circ \, c_0 (t) \big) \f + \infty$ when $t \f + \infty$, and this guarantees the existence of some $T_k > 0$ such that $ c_0 (t)
\notin  M_R^k (\g)$ for every $t > T_k $. Define $T := \max_{0 < k \le k_0} (T_k)$, then, for every $ t \in \, ]\, T , +\infty [$, 
$c(t) \notin \cup_{0 < k \le k_0} M_R^k (\g)$, hence $c(t) \notin M_R (\g)$. 
\end{proof}

\begin{lemma}\label{MRferme}
 For every hyperbolic isometry $\gamma $ and every real number $R> 0$ and, as 
\begin{itemize}
\item[(i)] for every $ k \in \Z^*$, $M^k_R (\gamma)$ is empty when $ R <  |k| \,\ell(\gamma)$ and non empty when $ R >  |k| \,\ell(\gamma) + \delta$,
\item[(ii)] $M_R (\gamma)$ is a closed set which is empty when $ R < \ell(\gamma)$ and non empty when $ R >  \ell(\gamma) + \delta$.
\end{itemize}
\end{lemma}

\begin{proof}
Notice that $\ell(\g) > 0$ by Lemma \ref{ellpositive}. If $M^k_R (\gamma)$ is not empty, for every $x \in \,M^k_R (\gamma)$, the basic properties of
functions $s(\cdot )$ and $\ell(\cdot )$ (see comments after Definitions \ref{deplacements}) give $ R \ge  d(x , \gamma^k \, x ) \ge s (\g^k) \ge \ell (\g^k) = 
|k| \,\ell(\gamma)$; this proves that $M^k_R (\gamma) = \emptyset$ when $ R <  |k| \,\ell(\gamma)$. On the other hand, by the Definition \ref{deplacements}
of the function $s(\cdot )$, $M^k_R (\gamma)\ne \emptyset \implies s (\g^k) \le R$ and $s (\g^k) < R \implies M^k_R (\gamma)\ne \emptyset$; as $ s (\g^k) 
\le  \ell (\g^k) + \delta =  |k| \,\ell(\gamma) + \delta$ by Lemma \ref{quasigeod} (i), then $ R > |k| \,\ell(\gamma) + \delta \implies M^k_R (\gamma) \ne 
\emptyset$.\\
By continuity of $x \mapsto d(x , \gamma^k \, x )$, each $M^k_R (\gamma)$ is a closed subset of $(X,d)$. Hence $M_R (\gamma)$ is closed for it is a finite 
union of the closed sets $M^k_R (\gamma)$ (corresponding, by (i), to the $k$'s such that $|k| \le  \frac{R}{\ell(\gamma)}$). If $ R < \ell(\gamma)$, then 
$ R <  |k| \,\ell(\gamma)$ for every $ k \in \Z^*$ and  it follows from (i) that $ \forall k \in \Z^*\  M^k_R (\gamma) = \emptyset$, thus that their union 
$M_R (\gamma)$ is empty. If $ R > \ell(\gamma) + \delta$, (i) guarantees that $M^1_R (\gamma) \ne \emptyset$, thus that $M_R (\gamma) \ne \emptyset$. 
\end{proof}

\begin{lemma}\label{distants}
For every hyperbolic isometry $\gamma $, for every $ r \,, \, R \in \, ]0 , +\infty[$ ($ r < R$) such that $M_r(\g) \ne \emptyset$, every point $x$ of the closure of 
$X \setminus M_R (\gamma)$ verifies $ d( x , M_{r} (\gamma))  \ge \frac{1}{2}\, (R - r)$.
\end{lemma}

\begin{proof} Denote by $\bar x$ any projection of $x$ onto the closed set $ M_{r} (\gamma)$. By Definition \ref{rayondeplact} and Remark 
\ref{Margudomain1}, there exists $k \in \N^*$ such that $d (\bar x , \g^k \bar x) \le r$. As $d(x , \g^k x) \ge R$, we get
$$ 2\; d( x , M_{r} (\gamma)) + r = 2\; d(x, \bar x) + r \ge  d(x, \bar x) + d ( \bar x , \g^k \bar x) + d(\g^k \bar x, \g^k  x) \ge d(x , \g^k x) \ge R \ ,$$
and this ends the proof.
\end{proof}

\begin{lemma}\label{tube}
For every hyperbolic isometry $\gamma $, for every $R > 0$ such that $M_R(\g) \ne \emptyset$,
\begin{itemize}
\item[(i)] every geodesic line $c $ connecting the fixed points $\gamma^-$ and $\gamma^+$ of $\g$ verifies (for every 
$x \in M_R (\gamma) $) : $ \Min_{t \in \R}\, d \big( x , c (t) \big) \le\frac{1}{2} \left(\frac{7 \, \delta}{ \ell(\g)} + 1 \right) R + \frac{7}{2}\, \delta$,

\item[(ii)] given any sequence $ \big( x_n\big)_{n \in \N}$ going to infinity in $M_R (\gamma)$ all limit points (in $\partial X$) of every subsequence are in
$\{ \g^- , \g^+\} $.

\item[(iii)] given any origin $x_0 \in X$ and the corresponding Dirichlet domain $D_{\g} (x_0) := \{ x  : \forall k \in \Z \ \ d(x_0,x) \le d(\g^k x_0 , x)\}$ of the 
action of the group $ \langle \g \rangle$, then $X = \cup_{k \in \Z} \ \g^k \big( D_{\g} (x_0) \big)$ and $ M_R (\gamma) \cap D_{\g} (x_0) $ is compact,

\item[(iv)] the function $x \mapsto d(x , \g x)$ attains its minimum on $X$.
\end{itemize}
\end{lemma}

\begin{proof}[Proof of (i)] Given any $x \in M_R (\gamma)$, by Definitions \ref{rayondeplact} and \ref{Margudomain}, there exists $p \in \N^*$ such that 
$ d(x, \g^p x ) \le R$, thus verifying $R \ge d(x, \g^p x ) \ge s(\g^p) \ge p\, \ell(\g)$, hence $p \le \dfrac{R}{\ell(\g)}$. Considering the integer $k$ such 
that $ (k-1)\,p\,\ell(\g)\le 7 \delta < k\,p\,\ell(\g)$, denote by $ y_k$ (resp. $ y_0$) a projection of the point $\g^{kp} (x)$ (resp. of the point  $ x$) onto the 
geodesic line $c$; denote by $ y'_k$ a projection of $\g^{kp} (y_0)$ onto the geodesic line $c$ and by $z_0$ a projection of $y_0$ onto the geodesic line 
$ \g^{kp} \circ \, c$; as $c $ and $ \g^{kp} \circ \, c $ are two geodesic lines from $\gamma^-$ to $\gamma^+$, the proposition \ref{geodasympt} (i)
implies that
\begin{equation}\label{approche}
d \left( \g^{kp} y_0  , y_k')\right) = d \left( \g^{kp} y_0 , c \right) \le 2\, \delta\ \ \ ,\ \ \ 
d \left( y_0  , z_0)\right) = d \left( y_0 ,  \g^{kp} \circ \,c )\right) \le 2\, \delta
\end{equation}
The quadrangle Lemma \ref{proprietes} (ii) and the equality $  d( \g^{kp} x ,  \g^{kp} y_0) = d(x , y_0)$ imply that
\begin{equation}\label{rectangulaire}
d(y_0 , \g^{kp} x) + d(x,  \g^{kp} y_0) - 2 \, \delta  \le \Max \left[ d(x , \g^{kp} x) + d(y_0,  \g^{kp} y_0) \,,\,
2 \, d(x , y_0)\right] \, ;
\end{equation}
using first the triangle inequality, and afterwards Lemma \ref{projection} and inequalities \eqref{approche}, we obtain
$$  d(x,  \g^{kp}  y_0) \ge  d(x,  y'_k) - d ( y'_k , \g^{kp} y_0) \ge d(x, y_0) + d(y_0 ,  y'_k ) - 4\,\delta
\ge d(x, y_0) + d(y_0 ,  \g^{kp} y_0 ) - 6\,\delta \ \ \ \text{and}$$
$$  d(y_0 ,  \g^{kp} x) \ge  d(z_0 ,  \g^{kp} x) - d(y_0 , z_0) \ge d( \g^{kp} x ,  \g^{kp} y_0 ) + d( \g^{kp} y_0 , z_0) - 4 \, \delta
 \ge d( x ,   y_0 ) + d(y_0 ,  \g^{kp} y_0 ) - 6\,\delta \, ,$$
where we used the fact that $y_0$ and $\g^{kp} y_0$ are respective projections of $x$ and $\g^{kp} x$ onto the geodesic lines $ c$ and $ \g^{kp} \circ \, c$.
These two inequalities yield
$$ d(y_0 , \g^{kp} x) + d(x,  \g^{kp} y_0) - 2\, \delta \ge  2 \, d(x, y_0)  + 2\; d(y_0 ,  \g^{kp} y_0 ) - 14 \delta >  2 \, d(x, y_0) \,,$$
where the last inequality stems from the choice of $k$, which guarantees that $ d(y_0 ,  \g^{kp} y_0 ) \ge \ell(\g^{kp}) = k\,p\, \ell(\g) > 7 \, \delta $.
Carrying forward these two last inequalities in \eqref{rectangulaire}, we get
$$ 2 \, d(x, y_0) +   d(y_0 ,  \g^{kp} y_0 ) - 7 \delta  <  2 \, d(x, y_0) +  2\, d(y_0 ,  \g^{kp} y_0 ) - 14 \delta <  d(y_0 , \g^{kp} x) + d(x,  \g^{kp} y_0) - 
2 \, \delta $$
$$   \le   d(x , \g^{kp} x) + d(y_0,  \g^{kp} y_0) \le k \,R + d(y_0,  \g^{kp} y_0) \le \left(\frac{7 \, \delta}{p \; \ell(\g)} + 1 \right) R + d(y_0,  \g^{kp} y_0)\,,$$
$$ \text{and thus: }\ \ \ \Min_{t \in \R}\, d \big( x , c(t) \big) =  d( x , y_0 )  \le \frac{1}{2} \left(\frac{7 \, \delta}{p \; \ell(\g)} + 1 \right)
 R +  \frac{7}{2}\, \delta \le \frac{1}{2} \left(\frac{7 \, \delta}{ \ell(\g)} + 1 \right) R + \frac{7}{2}\, \delta\ ;$$
and this concludes.
\end{proof}

\begin{proof}[Proof of (ii) and (iii)] Let $c$ be any geodesic line from $\gamma^-$ to $\gamma^+$.\\
Given any sequence $ \big( x_n\big)_{n \in \N}$ going to infinity 
in $M_R (\gamma)$, point (i) guarantees the existence of some $t_n \in \R$ such that $ d \big( x_n , c (t_n) \big) \le \frac{1}{2} \left(\frac{7 \, \delta}{ \ell(\g)} 
+ 1 \right) R + \frac{7}{2}\, \delta$, thus $c(t_n)$ goes to infinity (and thus $|t_n| \f +\infty$) when $n \f + \infty$. Consequently, as $d \big(x_n , c(t_n)\big)$ 
is bounded independently of $n$, when $n \f + \infty$, either $t_n \f +\infty$ and then $c(t_n)$ and $x_n$ both go to $\g^+$, or $t_n \f - \infty$ and then 
$c(t_n)$ and $x_n$ both go to $\g^-$, or $I_+ :=\{n : t_n \ge 0\}$ and $I_- := \{n : t_n < 0\}$ are both infinite and the sub-sequences 
$ \big( c(t_n) \big)_{n \in I+}$ and $ \big( x_n\big)_{n \in I_-}$ are respectively going to $\g^+$ and $\g^-$. This proves (ii).

\smallskip
If $c(t_k)$ is a projection of $\g^k \circ\, c(0) $ on the 
image of $c$, Proposition \ref{geodasympt} (i) and the fact that $c $ and $ \g^{k} \circ \, c $ are two geodesics from $\gamma^-$ to $\gamma^+$
guarantees that $d \big(c(t_k) ,  \g^k \circ\, c(0)\big) \le 2\, \delta$, and thus that, for every $k \in \Z$,
\begin{equation}\label{pas}
|t_{k+1} - t_k| = d \big( c(t_k) , c(t_{k+1})\big) \le d \big(  \g^k c(0) ,   \g^{k+1} c(0) \big) + 4\,\delta
\le d \big( c(0) ,   \g\, c(0) \big) + 4\,\delta \ .
\end{equation}
Fix any $t \in \R$. When $k$ goes to $+\infty$ (resp. to $-\infty$), $ \g^k c(0) $ goes to $\gamma^+$ (resp. to $\gamma^-$), thus $ c(t_k) $ also goes
to $\gamma^+$ (resp. to $\gamma^-$) for $d \big(c(t_k) ,  \g^k \circ\, c(0)\big) \le 2\, \delta$ by Proposition \ref{geodasympt} (i), and consequently
$t_k$ goes to $+\infty$ (resp. to $-\infty$); it follows that there exists $k' \in \Z$ such that $ t_{k'} \le t < t_{k'+1}$. It comes from this and from \eqref{pas}
that there exists $k \in \Z$ such that $ |t - t_k| \le\frac{1}{2} \,d \big( c(0) ,   \g\, c(0) \big) + 2\,\delta $. Hence, for every $t \in \R$, there exists $k \in \Z$ 
such that 
\begin{equation}\label{pas1}
d \big( c(t) ,  \g^k c(0) \big) \le d \big( c(t) , c(t_k) \big) + d \big( c(t_k) , \g^k c(0) \big) \le 
\frac{1}{2} \,d \big( c(0) ,   \g\, c(0) \big) + 4\,\delta\ .
\end{equation}
Consider any $x \in M_R (\gamma)$ and one of its projections $c(t)$ onto the geodesic line $c$; a consequence of \eqref{pas1} and of point (i) is the existence
of one $k \in \Z$ satisfying $ \, d \big( x ,  \g^k c(0) \big) \le d \big( x ,  c(t) \big) + d \big( c(t) ,  \g^k c(0) \big) \le R_1$, where $ R_1 :=  
\dfrac{1}{2} \left(\dfrac{7 \, \delta}{ \ell(\g)} + 1 \right) R + \frac{1}{2} \,d \big( c(0) ,   \g\, c(0) \big) + \dfrac{15}{2}\, \delta $; it follows that
$d \big( x ,  \g^k x_0 \big) \le d \big( x ,  \g^k c(0) \big) + d \big( \g^k c(0) , \g^k x_0 \big) \le d \big(x_0 , c(0)\big) + R_1 $. For every $ x \in 
M_R (\gamma) \cap D_{\g} (x_0) $, one has $d(x,x_0) \le d \big( x ,  \g^k x_0 \big) $ by definition of the Dirichlet domain, and the previous 
inequality gives $d \big( x , x_0 \big) \le d \big(x_0 , c(0)\big) + R_1 $. Thus $ M_R (\gamma) \cap D_{\g} (x_0) $, being a closed bounded subset of the proper 
space $(X,d)$, is compact.
\end{proof}

\begin{proof}[Proof of (iv)] Denote by $f$ the continuous function $ x \mapsto  d(x , \g x)$ and let $ s(\g) := \inf_{x \in X} f(x)$ and $R := s(\g) + 1$; as 
$\{x : f(x) \le R\} = M_R^1 (\gamma) \subset M_R (\gamma) $, it is sufficient to prove that $f$ attains its minimum when restricted to $M_R (\gamma) = 
\cup_{ k \in \Z} \left(\g^k \left(D_{\g} (x_0) \right) \cap M_R (\gamma) \right)$ and, as $ f \circ \g^k = f$ and as $M_R (\gamma) $ is stable by $\g^k$, it is 
sufficient to prove that $f$ attains its minimum when restricted to $ M_R (\gamma) \cap D_{\g} (x_0) $. This is verified for $ M_R (\gamma) \cap D_{\g} (x_0) $
is compact by (iii).
\end{proof}

\begin{lemma}\label{quasiconvexe}
For every hyperbolic isometry $\gamma $ and every $R> 0$, one has
\begin{itemize}
\item[(i)] for every $p \in \N^* $ and every $ x \in M_R^p (\gamma) $, every geodesic segment from $x$ to $\g^{p} x$ or to $\g^{- p} x$ is included in
$M_R^p (\gamma)$,

\item[(ii)] for every $ x , \, y \in M_R^p (\gamma) $, every geodesic from $x$ to $y $ is included in $M_{R + 2\,\delta}^p (\gamma)$.

\item[(iii)]  for any $ x \in M_R (\gamma) $ and every $k \in \Z$, any geodesic from $x$ to $\g^k x$ is contained in $M_{R + 2\,\delta} (\gamma)$.
\end{itemize}
\end{lemma}

\begin{proof}[Proof of (i) and (ii)] We suppose that $ M_R^p (\gamma) \ne \emptyset$ (elsewhere (i) and (ii) are trivially verified). For every $ x \in 
M_R^p (\gamma) $ and any $ q \in \{-p , p\}$, choose any geodesic segment $[ x , \g^q x]$ from $x$ to $\g^q x$; for every $u \in[ x , \g^q x]$, 
the triangle inequality yields
$$ d( u , \g^p u )= d( u , \g^q u ) \le d( u , \g^q x ) + d( \g^q x , \g^q u ) =  d(x, u) + d( u , \g^q x ) = d( x , \g^q x ) \le R \ ,$$
hence $u \in M_R^p (\gamma)$ and this proves (i).\\
For every $ x , \, y \in M_R^p (\gamma) $ and any (normally parametrized) geodesic segment $c$ from $x$ to $y$, if $ L := d(x,y)$, the quasi-convexity 
property of $\delta$-hyperbolic spaces (see Proposition 25 , p. 45 of \cite{GH}) and (afterwards) the fact that $c(0) =x$ and $c(L) = y$ both belong to 
$M_R^p (\gamma) $ give:
$$ d \big( c(t)\, , \, g^p \circ \, c(t) \big) \le \left( 1 - \dfrac{t}{L} \right)  d \big( c(0)\, , \, g^p  \circ \,c(0) \big) + 
 \dfrac{t}{L} \ d \big( c(L) \,,\, g^p  \circ \,c(L) \big) + 2\,\delta \le R + 2\,\delta \ ,$$
and this proves (ii) because $c$ is entirely contained in $M_{R + 2\,\delta}^p (\gamma)$.
\end{proof}

\begin{proof}[Proof of (iii)] By definition of $ M_R (\gamma) $, for any $ x \in M_R (\gamma) $, there exists some $p \in \N^*$ such that $ x \in 
M_R^p (\gamma) $. For every $k \in \Z$, as $\g^k x$ also belongs to $M_R^p (\gamma)$, it comes from (ii) that every geodesic segment from $x$ to $\g^k x$ 
is contained in $M_{R + 2\,\delta}^p (\gamma)$.
\end{proof}

\subsection{Hyperbolic isometries of Busemann Gromov-hyperbolic spaces}

\begin{defi}\label{dconvexe}
A metric space $(X,d)$ is \emph{a Busemann space} if it is geodesic, proper and if its distance $d$ is convex: i. e. if,
for every pair of geodesic segments $c_1\, , c_2 $ (reparametrized in order that $c_1\, , c_2 $ are defined
from $[0,1]$ to $X$), the function $t \mapsto d(c_1(t) , c_2 (t) )$ is convex.
\end{defi}

Examples of Busemann spaces are given by simply connected Riemannian manifolds whose sectional curvature is non-positive and, more generally,
by $\CAT (0)$ spaces.

\begin{remark}\label{unicite} In a Busemann space any two points are connected by a unique geodesic segment.
\end{remark}
The proof, elementary, is left to the reader.

\begin{lemma}\label{convex1}  In a Busemann space, for any pair $c_0\, , c_1 $ of geodesic segments the function $t \mapsto  d \big(c_1 (t) , 
{\rm Im}(c_0)\big)$ is convex.
\end{lemma}

\begin{proof} Denote by $[\,0 , a_0\,]$ and $[\,0 , a_1\,]$ the intervals of definition of $c_0$ and $c_1$ (respectively); for every $ t_1 , \, t_2 \in 
[\,0 , a_1\,]$ ($t_1 < t_2$), denote by $c_0 (s_1) $ and $c_0 (s_2)$ the respective projections of $c_1 (t_1) $ and $c_1 (t_2)$ onto the image of $c_0$.
For every $ t \in [t_1 , t_2]$, we define $ \beta := \dfrac{t - t_1}{t_2 - t_1}$ (i. e. $ t = (1- \beta)\,t_1 + \beta \, t_2  $) and the convexity of the distance
yields:
$$d \big(c_1 (t) , {\rm Im} (c_0)\big)  \le d \left( c_1 \big((1- \beta)\,t_1 + \beta \, t_2\big)\, , 
\,c_0 \big((1- \beta)\,s_1 + \beta \, s_2\big)\right) \le (1- \beta) \, d \big( c_1 (t_1), c_0 (s_1) \big) +$$
\begin{equation*}
+ \beta \;  d \big( c_1 (t_2), c_0 (s_2) \big) =  (1- \beta) \, 
d \big( c_1 (t_1),  {\rm Im} (c_0) \big)+ \beta \;  d \big( c_1 (t_2), {\rm Im} (c_0) \big)\ ,
\end{equation*}
and this ends the proof.
\end{proof}

\begin{lemma}\label{localglobal} 
On a Busemann space $(X,d)$, every local geodesic is a (minimizing) geodesic.
\end{lemma}

\begin{proof} Consider any local geodesic $c : I \to X$ (see definition in section 2) and any pair $t_1 , t_2 \in I$ ($t_1 < t_2$); denote by $[x,y]$ the 
(minimizing) geodesic segment from $ x := c(t_1)$ to $ y := c(t_2)$ and by $K $ the closed nonempty subset of $ [t_1 , t_2]$ where the function 
$f : s \mapsto d( c(s) , [x,y])$ attains its maximum (on $ [t_1 , t_2]$). For every $t \in K \cap \, ]t_1 , t_2[$, there exists an open interval $J$ containing $t$ 
on which the restriction of $c$ is a (minimizing) geodesic and on which $f$ is convex (by Lemma \ref{convex1}), as $f$ attains its maximum at $t$, its convexity 
implies that $f$ is constant on $J$, consequently $ K \cap \,]t_1 , t_2[ $ is open and closed in $ ]t_1 , t_2[$ thus either $ K \cap \, ]t_1 , t_2[ = ]t_1 , t_2[ $ (and then,
for every $t \in \, ]t_1 , t_2[$, $f(t) = f(t_1) = 0$ by continuity) or $ K \cap \, ]t_1 , t_2[ = \emptyset $ (and then $K \subset \{t_1 , t_2\}$, thus $f = 0$ on
$ [t_1 , t_2]$). In every case $c( [t_1 , t_2])$ coincides with the geodesic segment $[c(t_1) , c(t_2)]$ and $c$ is minimizing between any pair of its points.
\end{proof}

\begin{lemma}\label{fixedpoint} 
On every Busemann space $(X,d)$, every isometry $\g \ne \id_X$ of finite order admits some fixed point.
Reciprocally, for every isometry $\g \ne \id_X$ such that the subgroup $\langle \g \rangle $ is discrete, if $\g$ admits a fixed point, then it has finite order.
\end{lemma}

\begin{proof} The converse property being proved in Lemma \ref{autofidele} (iv), it is sufficient to prove that every isometry $\g \ne \id_X$ of finite order 
$N$ admits some fixed point.
Let $G := \{\id_X , \g , \ldots , \g^{N-1}\}$ be the subgroup generated by $\g$. Choosing an initial point $u_0$, there exists $x_0 \in X$ 
and $R> 0$ such that $G \cdot u_0 \subset \overline B_X (x_0 , R)$. We construct by iteration the sequence $(u_n)_{n \in \N}$ by choosing $u_{n+1}$
as the middle-point of the geodesic segment $[u_n , \g u_n]$. Let us prove by iteration the property $G \cdot u_n \subset \overline B_X (x_0 , R)$: this 
property is true when $n=0$; suppose now that it is true at step $n$ and prove it at step $ n+1$: as $\overline B_X (x_0 , R)$ is convex and as 
$G \cdot u_n \subset \overline B_X (x_0 , R)$, for every $k \in \{0 , \ldots , N-1\}$, one has $\g^k ([u_n , \g u_n]) =[\g^k u_n , \g^{k+1} u_n] \subset \overline B_X (x_0 , R)$, thus its middle-point $\g^k u_{n+1}$ lies in $\overline B_X (x_0 , R)$, this proves that $G \cdot u_{n+1} \subset 
\overline B_X (x_0 , R)$ and thus that $G \cdot u_n \subset \overline B_X (x_0 , R)$ for every $n \in \N$.\\
As $\overline B_X (x_0 , R)$ is compact, there exists a subsequence $(u_{n_k})_{k \in \N}$ of the sequence $(u_n)_{n \in \N}$ which converges, denote by
$u$ its limit. As $d(u_{n+1} , \g u_{n+1}) \le  d(u_{n+1} , \g u_{n}) + d(u_{n} , u_{n+1}) = d(u_{n} , \g u_{n}) $, the sequence $n \mapsto 
d(u_{n} , \g u_{n}) $ is decreasing and converges to some value $\alpha \ge 0$, it follows that $\alpha = d(u , \g u)$.\\
If $\alpha > 0$, let $v$ be the middle-point of $[u , \g u]$, there exists a positive sequence $\e_k $ (going to zero when $k \to +\infty$) such that 
$d(u , u_{n_k}) < \e_k$ and $d(\g u ,  \g u_{n_k}) < \e_k$; the convexity of the distance (applied to the middle-points of $[u , \g u]$ and 
$[u_{n_k} , \g u_{n_k}]$) then implies $d(v, u_{n_k + 1}) < \e_k$; this (and the fact that $n \mapsto d(u_{n} , \g u_{n}) $ is decreasing) yields
$$d(u , \g u) = d(\g u , \g v ) + d( v , \g u)  \ge d(v , \g v) \ge d (u_{n_k + 1} , \g u_{n_k+1}) - 2 \e_k \ge d(u , \g u) -  2 \e_k \, .$$
When $\e_k \to 0$, this proves that $d(v , \g v) = d( v , \g u) + d(\g u , \g v )$, thus that, for every $k \in \{0 , \ldots , N-1\}$, the length of
$[\g^k v, \g^{k+1} u ] \cup [\g^{k+1} u , \g^{k+1} v]$ is equal to $d(\g^{k} v , \g^{k+1}  v) $; hence 
$[\g^k v, \g^{k+1} u ] \cup [\g^{k+1} u , \g^{k+1} v]$ is the (minimizing) geodesic from $\g^{k} v $ to $ \g^{k+1} v $. It follows that the union of all the
$[\g^k u, \g^{k+1} u ]$ is a local geodesic, thus a (globally minimizing) geodesic by Lemma \ref{localglobal}, in contradiction with the fact that it is a closed
path.\\
The only possibility is thus that $\alpha = 0$ and this implies that $u$ is a fixed point of $\g$.
\end{proof}

\begin{lemma}\label{invariantgeod} 
On any $\delta$-hyperbolic Busemann space $(X,d)$, for every hyperbolic isometry $\g$ and every $ x \in M_{\rm{min}} (\g)$, the union (or more precisely
the concatenation), for all the $ p \in \Z$,
of the geodesic segments $ \left[ \g^p  x , \g^{p+1}  x \right]$ is a $\g$-invariant geodesic line $c_{\g}$ from $\g^-$ to $\g^+$, contained in 
$M_{\rm{min}} (\g)$, such that $ \g \big( c_{\g} (t) \big) = c_{\g} \big(t + \ell (\g) \big)$. Consequently $s(\g) = d(x,\g x) = \ell (\g)$.
\end{lemma}

\begin{proof} By Lemma \ref{unicite} $ \left[ \g^p  x , \g^{p+1}  x \right] = \g^p \left(\left[  x , \g \, x \right] \right)$ is the unique geodesic segment from
$ \g^{p}  x$ to $ \g^{p+1}  x$. Then Lemma \ref{geodmin} proves that the concatenation $c_\gamma := \cup_{p \in \Z} \left[ \g^p  x , \g^{p+1}  x \right]$
(oriented in the sense of increasing values of $p$) 
is a $\g$-invariant local geodesic included in $M_{\rm{min}} (\g)$ and Lemma \ref{localglobal} guarantees that this local geodesic is a (minimizing) geodesic.
This implies that, for every $k \in \Z$, $d ( x, \g^k x) = |k| \, d(x, \g x)$, thus that $c_\gamma$ is defined on $] - \infty , +\infty [$ and verifies 
$c_{\g} (\pm \infty) = \lim_{p \to +\infty} \left(\g^{\pm p} x\right) = \g^{\pm}$. This also implies that $\ell (\g) = \lim_{p \to +\infty} \frac{1}{p}\, d(x, \g^p x) 
= d(x, \g x) = s(\g)$. Notice that, for every $t \in \R$, there exists $p \in \Z$ such that $c_\gamma (t) \in \left[  \g^{p-1}  x , \g^p  x \right]$, and then
$d\big( c_\gamma (t) ,  \g c_\gamma (t)\big) = d\big(  c_\gamma (t) ,  \g^p  x \big) + d\big( \g^p  x  ,  \g c_\gamma (t)\big) = 
d\big( \g^{p-1}  x  ,  c_\gamma (t)\big) + d\big(  c_\gamma (t) ,  \g^p  x \big) = d\big( \g^{p-1}  x ,  \g^p  x \big) = d(x, \g x) = \ell (\g)$. This proves that,
if $c_\gamma (s) := \g c_\gamma (t)$, then $s-t = d\big( c_\gamma (t) ,  \g c_\gamma (t)\big) =\ell (\g)$ and this achieves the proof.
\end{proof}

Lemma \ref{invariantgeod} does not claim that there exists a unique geodesic from $\g^-$ to $\g^+$; indeed this is generally false: a 
counter-example can be constructed starting from a simply connected surface $S$ with sectional curvature $\le -1$, choosing an hyperbolic isometry $\g$
and its $\g$-invariant geodesic $c_\gamma$ and gluing a flat strip of width $\e$ between the two connected components of $ S \setminus c_\gamma$.

\subsection{Elementary subgroups}\label{nilpotents}

In this section we consider \emph{discrete} subgroups $G$ of the isometry group of any Gromov-hyperbolic space $X$. A classification of these groups has been 
sketched by M. Gromov \cite{Gr4} (see also \cite{SUD}, \cite{CCMT}), in terms of their limit set $LG$ (i. e. the set of accumulation points of any orbit of the 
action of $G$ on $X$); he classified these groups in the following classes\footnote{Here we only consider groups acting discretely; excluding for instance 
\emph{focal} groups whose action is not discrete.}:

\begin{itemize}
\item {\em elliptic} groups (also said {\em bounded}): finite groups whose all orbits are thus bounded;
\item {\em parabolic} (or, according to the original terminology, {\em horocyclic}) groups: infinite groups $G$ such that $\# (LG) = 1$. A parabolic 
group thus only contains parabolic or elliptic elements\footnote{Indeed, if $G$ contains an hyperbolic isometry $g$, its fixed points $g^+$ and $g^-$ are 
accumulation points of the sequence $\left( g^k \, x\right)_{k \in \Z}$, they thus belong to $LG$.}.
\item {\em \lq \lq lineal"} groups: infinite groups $G$ such that $\# (LG) = 2$; a \lq \lq lineal" group only contains hyperbolic or elliptic elements (for a proof, see
for instance \cite{C} section 3.4.2).
\item groups of {\em general type}: groups $G$ such that $\# (LG) \ge 3$ (in this case $LG$ is infinite); this is equivalent to say that $G$ contains (at least) two
hyperbolic elements whose sets of fixed points are disjoint\footnote{If $G$ contains two hyperbolic elements $ g_1$ and $ g_2$ whose sets of fixed points 
$\text{Fix}(g_1)$ and $\text{Fix}(g_2)$ are disjoint, a trivial consequence is that $\# (LG) > 2$, for $\# (LG)$  contains $\text{Fix}(g_1) \cup \text{Fix}(g_2)$
the converse implication is proved for instance in \cite{C}, Lemma 3.7. The fact that, if $ \# (LG)>2$, then $ LG$ is infinite uncountable is announced in 
\cite{Gr4},  section 3.5, Theorem p.194; one can find a complete proof in \cite{SUD}, Proposition 6.2.14.}.
\end{itemize}

In the three first cases, the action of $G$ on $X$ is said to be {\em elementary}. By extension, an hyperbolic group will be said \lq \lq elementary" if the
action of $G$ (by left translations) on $G$ (endowed with the algebraic word-metric) is elementary; as, for this action, $LG$ coincides with the ideal boundary 
$\partial G$, an hyperbolic group is \lq \lq elementary" iff $\# (\partial G) \le 2$.\\
Notice that, if the action of $G$ is elementary, the set $\text{Fix}(g)$ of fixed points of any non elliptic element $ g \in G$ coincides\footnote{Indeed, for any parabolic 
(resp. lineal) group $G$, for every parabolic (resp. hyperbolic) isometry $g \in G$, the fixed points of $g$ are the accumulation points of the sequence
$\left( g^k \,x\right)_{k \in \Z}$, thus $\text{Fix}(g) \subset LG$; furthermore, as $\# \big( \text{Fix}(g)\big) = \# (LG)$, then $\text{Fix}(g) =  LG$.} 
with the limit set $LG$ of $G$.

\medskip
In the three following Propositions, we shall recall basic properties of elementary groups, most of them being immediate corollaries (well known among specialists) 
of the above classification and will thus be recalled without proof, with a simple reference.

\begin{prop}\label{actionelementaire} \emph{(Elementary actions)}
Let $ G$ be a discrete subgroup of the isometry group of any Gromov-hyperbolic space $X$, then:
\begin{itemize}
\item[(i)] if $ \g_1, \g_2 \in G $ are two hyperbolic elements with a common fixed point, then they have the same pair of fixed points, i. e. $ \{ \g_1^-, \g_1^+ \} 
= \{ \g_2^-, \g_2^+ \} $;

\item[(ii)] if $\g \in G$ is an hyperbolic isometry, the subgroup $G_{\g} = \left\{ g \in G \, : \, g\big( \{ \g^-, \g^+ \} \big) = \{ \g^-, \g^+ \}\right\}$ is the 
maximal subgroup among all virtually cyclic subgroups of $G$ which contain $\g$; if moreover $G$ is torsion-free, then $G_{\g} = \left\{ g \in G \, : \, 
g ( \g^-) =  \g^- \text{ and } g ( \g^+) =  \g^+\right\}$;

\item[(iii)] if $G$ est amenable (e. g. virtually nilpotent), then the action of $G$ is elementary; 

\item[(iv)] if $G$ is virtually nilpotent, its non elliptic elements are either all parabolic or all hyperbolic and all have the same set of fixed points.

\item[(v)] if $\g \in G$ is an hyperbolic isometry, for every $g\in G$, the subgroup generated by
$\g$ and $g \g g^{-1}$ is virtually cyclic if and only if the subgroup generated by $\g$ and $g $ is virtually cyclic;

\item[(vi)] if $a$ and $b$ are hyperbolic isometries and if $ \langle a , b  \rangle$ is not virtually cyclic then, for every
$p , q \in \Z^*$, $  a^p$ and $ b^q  $ are hyperbolic and $ \langle a^p, b^q  \rangle$ is not virtually cyclic;

\item[(vii)] if $\g \in G$ is an hyperbolic isometry, any subset $S$ of $G$ such that $ \langle \g, g\rangle$ is 
virtually cyclic for every $g\in S$ generates a virtually cyclic group.
\end{itemize}
\end{prop}

For the proof of points (i) and (ii) of Proposition \ref{actionelementaire}, a reference is \cite{C}, Propositions 3.21  and 3.27; though this last Proposition only 
proves the first part of (ii), the second part of (ii) immediately follows from this first one because, for every $g \in G_\gamma$ which swaps points $ \g^- $ 
and $ \g^+ $,
$g^2$ fixes $ \g^- $ and $ \g^+ $ and, as the torsion-free hypothesis implies the absence of elliptic elements (by Remark \ref{kpointsfixes} (i)), $g^2$ is 
parabolic or hyperbolic, thus hyperbolic (because a parabolic isometry cannot fix two points of the ideal boundary), and then $g$ would fix $ \g^- $ and 
$ \g^+ $ by Remark \ref{kpointsfixes} (ii).\\
If the action of $G$ is non elementary, then $G$ contains a free subgroup (see \cite{D}, \cite{Kou}, \cite{CDP} section 11  Prop. 3.1 and\cite{SUD} Prop. 10.5.4) 
and is thus non amenable; this proves (iii).\\
Point (iv) immediately follows from point (iii) (which proves that the action of $G$ is elementary), from the above classification of elementary groups and from 
the remark which follows this classification, which says that every non elliptic $\g \in G$ verifies $\text{Fix}(g) =  LG$.\\
Point (vii) is a direct consequence of (ii), which implies that $ \langle \g, g\rangle \subset G_\gamma$ 
for every $g\in S$.\\
Let us now prove (v): if the subgroup $ \langle\g , g \rangle $ generated by $\g$ and $g $ 
is virtually cyclic, then its
subgroup $ \langle\g , g \g g^{-1} \rangle $ is also virtually cyclic. Conversely, if the subgroup 
$ \langle\g , g \g g^{-1} \rangle $ is virtually cyclic, then the point (ii) implies that the set $\{ \g^-, \g^+ \}$ of fixed
points of $\g$ is globally invariant by $ g  \g g^{-1}$; this proves that $\{ \g^-, \g^+ \}$ is the set of fixed points
of  the hyperbolic isometry $ g  \g^2 g^{-1}$. As the fixed points of $ g  \g^2 g^{-1}$ are also $ g(\g^-)$ and $ g(\g^+)$,
it follows that $ \{g(\g^-) \,,\, g(\g^+)\} = \{ \g^-, \g^+ \} $ and the point (ii) then implies that $ \langle\g , g \rangle $
is virtually cyclic, as being a subgroup of the subgroup $G_{\g}$ which (globally) stabilizes the set $\{ \g^-, \g^+ \}$.\\
Proof of (vi): By Remark \ref{kpointsfixes} (ii), $ a^p$ and $ b^q  $ are hyperbolic and their sets of fixed points verify $\text{Fix}(a^p) =  
\text{Fix}(a)$ and $\text{Fix}(b^q) =  \text{Fix}(b)$; if $ \langle a^p, b^q  \rangle$ is virtually cyclic, it comes from (iv) that $\text{Fix}(a^p) =  
\text{Fix}(b^q) $, thus that $\text{Fix}(a) = \text{Fix}(b)$; by (ii), this implies that $ \langle a, b  \rangle$ is virtually cyclic.

\begin{corollary}\label{elementaryaction}
Let $\Gamma$ be a group acting properly (but eventually not faithfully) on a Gromov-hyperbolic space $(X,d)$, and let $\g $ be any element of $\Gamma$ 
such that $\ell(\g) > 0$, then
\begin{itemize}
\item[(i)] for every $ g \in \Gamma$, the subgroup generated by $\g$ and $ g \g g^{-1}$ is virtually cyclic 
if and only if the subgroup generated by $\g$ and $ g $ is virtually cyclic,
\item[(ii)] for every $ g \in \Gamma$ satisfying $\ell(g) > 0$, if there exists $p, q \in \Z^*$ such that $ \g^p$ and $g^q$ generate a virtually cyclic 
subgroup, then $\g$ and $g$ generate a virtually cyclic subgroup of $\Gamma$,
\item[(iii)] any subset $S \subset \Gamma^*$ such that $\langle \g , \sigma \rangle$ is virtually cyclic for every $\sigma \in S$ generates a virtually 
cyclic subgroup of $\Gamma$.
\end{itemize}
\end{corollary}

\begin{proof}
Let $\varrho : \Gamma \f \text{Isom} (X,d)$ be the representation corresponding to the action of $\Gamma$ on $ (X,d)$, Lemmas \ref{reductisom} (vi) 
and \ref{ellpositive} prove that $\varrho (\g)$ is an hyperbolic isometry of $(X,d)$.\\
\emph{Proof of (i)}: Proposition \ref{actionelementaire} (v) implies that, for every $g \in \Gamma$, the subgroup generated by 
$\varrho (\g)$ and $ \varrho (g \g g^{-1})$ is virtually cyclic if and only if the subgroup $\varrho (\langle \g , g \rangle )$ generated by $\varrho (\g)$ 
and $ \varrho (g)$ is virtually cyclic, thus (by Lemma \ref{reductisom} (vii)) if and only if the subgroup $\langle \g , g \rangle $ is virtually cyclic.\\
\emph{Proof of (ii)}: For every $g \in \Gamma$ satisfying $\ell(g) > 0$, $\varrho (g)$ is also an hyperbolic isometry of $(X,d)$ by Lemma \ref{ellpositive}; if 
$\langle \g^p, g^q \rangle$ is virtually cyclic, then $\langle \varrho (\g)^p, \varrho (g)^q \rangle$  is virtually cyclic too and Proposition 
\ref{actionelementaire} (vi) then guarantees that $\langle \varrho (\g), \varrho (g) \rangle = \varrho (\langle \g , g \rangle )$  is virtually cyclic, 
it then follows from  Lemma \ref{reductisom} (vii) that $\langle \g , g \rangle $ is virtually cyclic.\\
\emph{Proof of (iii)}: for any $S \subset \Gamma^*$, if $\langle \g , \sigma \rangle$ is virtually cyclic for every $\sigma \in S$,  then 
$\langle \varrho (\g), \varrho (\sigma) \rangle $  is virtually cyclic for every $\varrho (\sigma) \in \varrho (S)$ and it follows from Proposition 
\ref{actionelementaire} (vii) that $\langle \varrho (S) \rangle $ is virtually cyclic, thus (by Lemma \ref{reductisom} (vii)) that $\langle S \rangle $
is virtually cyclic.
\end{proof}

In the co-compact case, one has the following classical results:

\begin{prop}\label{actioncocompacte} 
\emph{(Cocompact actions, see \cite{GH}, \cite{CDP}, \cite{BH})} Let $G$ be a discrete co-compact group of isometries of a Gromov-hyperbolic space $X$, then
\begin{itemize}
\item[(i)] $G$ is a finitely generated Gromov-hyperbolic group;
\item[(ii)] $G$ does not contain any parabolic isometry;
\item[(iii)] $G$ is elementary if and only if the space $X$ is elementary;
\item[(iv)] $G$ is elementary if and only if it is virtually cyclic;
\item[(v)] every virtually nilpotent subgroup of $G$ is virtually cyclic.
\end{itemize}
\end{prop}

For a proof of (i), a reference is \cite{GH}: indeed, in the co-compact case, the space and the group are quasi-isometric, thus the space is Gromov-hyperbolic
if and only if the group is\footnote{The notion of Gromov-hyperbolicity does not depend on the system of generators, however the hyperbolicity constant $\delta$
highly depend on this system. This is a problem that we have to take into account all over this paper.}.\\
For point (ii), it is classical that an hyperbolic group contains no parabolic element (for the action by left-translations, see \cite{GH} section 8, th\'eor\`eme 29 
or \cite{CDP} section 9, th\'eor\`eme 3.4). Proving the same result for any co-compact action of the same group on any Gromov-hyperbolic space is part of the
folklore, the main argument being the quasi-isometry between the group and the space\footnote{More precisely, the characteristics of the action of $g \in G$ on
$(X,d)$ is determined by its action on $\partial X$ and its fixed points, as the ideal boundary $\partial G$ of $G$ (for the algebraic word-metric) is identified to
$\partial X$ by means of a homeomorphism which intertwines the actions of $G$ on $X$ and $G$, the characteristics of the action of $g \in G$ on $(X,d)$
are determined by its action on $\partial G$, it thus cannot be parabolic.}. The same argument proves that $G$ is elementary iff its co-compact action on $X$ 
is elementary.\\
Point (iii) follows from this (i. e. $\partial X \simeq \partial G$) and from the fact that, in the co-compact case, $LG = \partial X $.\\ 
For point (iv), either $G$ is finite, and then (iv) is trivially verified, or $G$ is infinite and Lemma \ref{loxodromique} then implies that $\partial G$ has exactly two
points and that $G$ contains a hyperbolic isometry $\gamma$ (thus $\{ \g^-, \g^+ \} = \partial G$ and the group $G_\gamma $ of the elements $g$ such that
$g (\{ \g^-, \g^+ \}) = \{ \g^-, \g^+ \}$ coincides with $G$). Applying Proposition \ref{actionelementaire} (ii) then proves that $G_\gamma = G $ is virtually 
cyclic.\\
Point (v) is a consequence of the fact that $G$ is an hyperbolic group (by (i)) and of the fact that every infinite nilpotent subgroup $G'$ of $G$ contains 
a cyclic subgroup of finite index (by \cite{GH}, section 8, Th\'eor\`eme 37  p. 157).

\bibliographystyle{alpha}
\bibliography{bibliography-Margulis}

\begin{thebibliography}{CCMT15}
\bibitem[BZ88]{Bu-Za}
Yu.D. {Burago} and V.A. {Zalgaller}.
\newblock {\em {Geometric inequalities. Transl. from the Russian by A. B.
  Sossinsky. Transl. from the Russian by A.B. Sossinsky.}}
\newblock Berlin etc.: Springer-Verlag, transl. from the russian by a.b.
  sossinsky edition, 1988.
\end{thebibliography}

\bigskip
\textsc{G. Besson, Institut Fourier, CNRS et Universit\'e Grenoble Alpes,
    CS 40700, 38058 Grenoble Cedex 09, France}\par\nopagebreak
  \textit{E-mail address}: \texttt{g.besson@univ-grenoble-alpes.fr}

\medskip

\textsc{G. Courtois, Institut de Math\'ematiques de Jussieu-Paris Rive Gauche, CNRS et Universit\'e Pierre et Marie Curie, 4 place Jussieu, 75232 Paris Cedex 09, France}\par\nopagebreak
  \textit{E-mail address}: \texttt{gilles.courtois@imj-prg.fr}
  
\medskip

\textsc{S. Gallot, Institut Fourier, Universit\'e Grenoble Alpes,
    CS 40700, 38058 Grenoble Cedex 09, France}\par\nopagebreak
  \textit{E-mail address} : \texttt{sylvestre.gallot@univ-grenoble-alpes.fr}
  
\medskip

\textsc{A. Sambusetti, Dipartimento di Matematica, Sapienza Universit\`a  di Roma, Piazzale Aldo Moro 5, 00185 Roma, Italy}\par\nopagebreak
   \textit{E-mail address} : \texttt{sambuset@mat.uniroma1.it}

\end{document}